\documentclass[12pt]{article}%twoside
\usepackage{litpreamble}

\begin{document}

\begin{titlepage}
    \begin{center}
%        \vspace*{0.05cm}
            
        \LARGE
        \textbf{Blow-ups of Lie groupoids and Lie algebroids}
            
        \vspace{0.1cm}
        \large
        A comprehensive study of the deformation to the normal cone and blow-up constructions for Lie groupoids and Lie algebroids
        
        \vspace{1cm}
            
        \includegraphics[width=0.9\textwidth]{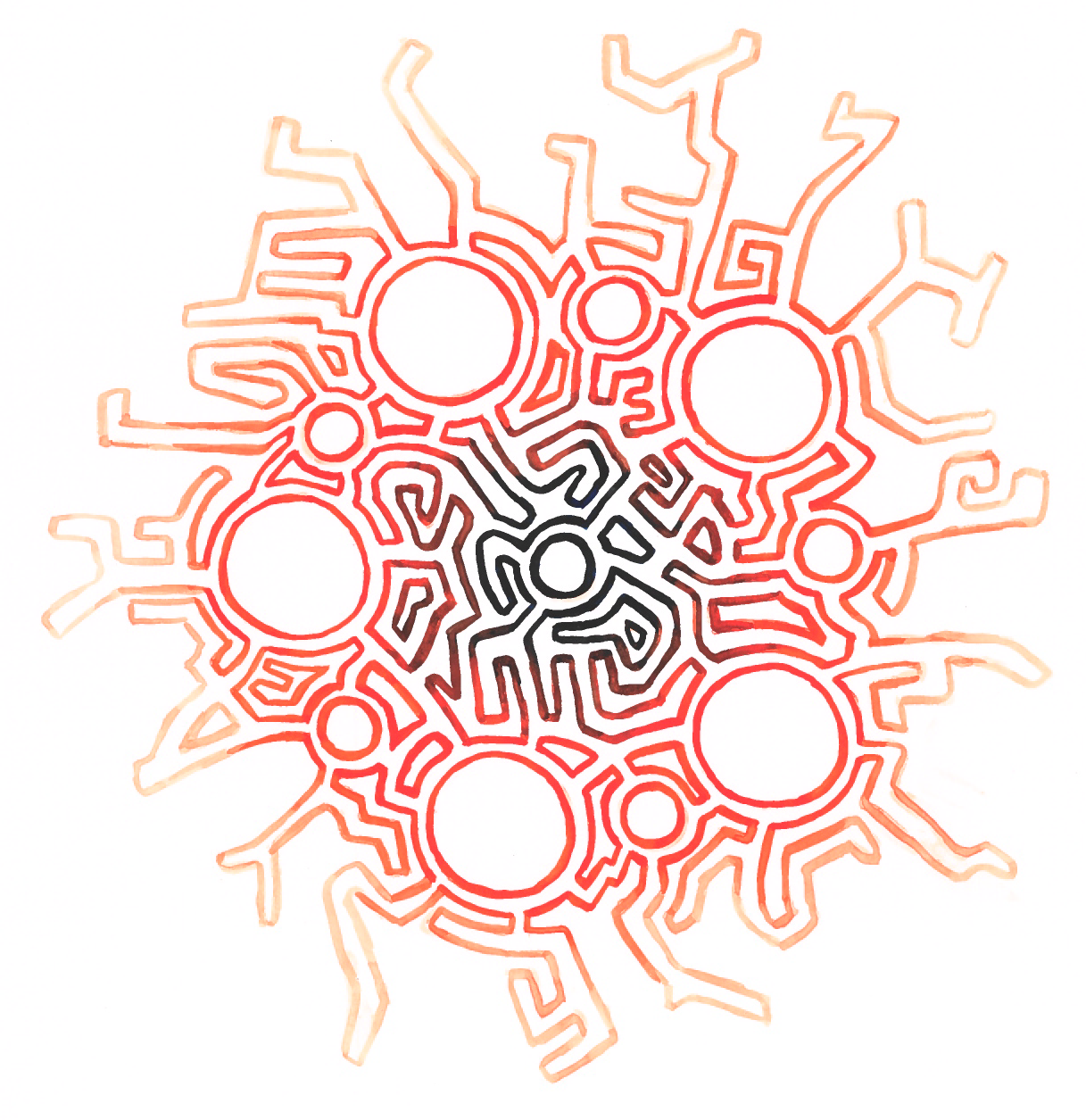}
        
        \vspace{1cm}
        
        \small A thesis presented for the degree of master's of science at \\
        \small Radboud University \\
        \includegraphics[width=0.1\textwidth]{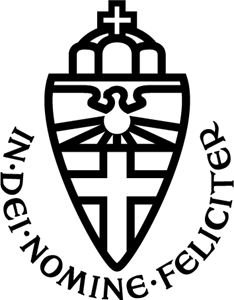} \\
        September 2021 \\
        
        \vspace{0.5cm}
        
        \flushleft{%\large 
        \textit{Author:}  Lennart Obster \hfill \textit{Supervisor:}  dr. Ioan M\v{a}rcu\c{t}
        } \\
        \hfill%\large 
        \textit{Second reader:} prof. dr. Klaas Landsman

    \end{center}
\end{titlepage}

\thispagestyle{myplain}
\tableofcontents\newpage

\setcounter{page}{1}
\section{Introduction}
The theory of Lie groupoids and Lie algebroids is nowadays very popular. The theory goes back to \cite{Pradines1,Pradines2}, where both Lie groupoids and Lie algebroids are already studied. The theory of Lie groupoids flourished once Alain Connes discovered the use of Lie groupoids in noncommutative geometry (see e.g. \cite{Connes,Connesbook}). Once Lie algebroids were recognised as natural objects to study in Poisson Geometry (see e.g. \cite{Weinstein,Karasev}), many mathematicians regained interest in the theory of Lie algebroids. Nowadays, there are numerous applications in Poisson geometry \cite{Cattaneo}, noncommutative geometry \cite{CannasWeinstein}, and geometric quantisation \cite{Kontsevich}, to name a few.

In \cite{gualtieri2012symplectic} and \cite{2017arXiv170509588D}, a particular geometric construction of Lie groupoids is explained which builds from a Lie groupoid and a closed Lie subgroupoid (the natural subobject which is also closed) a new Lie groupoid. One can think about this construction as replacing the closed Lie subgroupoid with a new closed Lie subgroupoid. In its most trivial form, this so-called blow-up construction corresponds to a blow-up construction for manifolds. For manifolds, blow-up constructions have been studied in the 90’s by Richard Melrose, who used spherical blow-ups in his geometric approach to the Atiyah-Patodi-Singer Index theorem \cite{Melrosegreen}. However, (iterated) spherical blow-ups “add a boundary” to the manifold in question, and one is therefore led to study manifolds with corners (that is, locally isomorphic to an open subset of some $\mathbb{R}^k \times \mathbb{R}_{\ge0}^{\ell}$). With our eye on the theory of Lie groupoids and Lie algebroids, we will not be interested in this type of blow-up. Instead, we will be interested in another type of blow-up intimately related to it: the projective blow-up construction. This type of blow-up is well-known, and studied extensively, in the world of algebraic geometry.
In differential geometry, blow-ups have proven useful for many different purposes, albeit in disguise. For example, a special type of (projective) blow-up is the connected sum operation with a real projective space. %These types of connected sums are still considered vital to the study of classifying, say, smooth manifolds in dimensions $2$ and $3$. 
Blow-ups have also been an important tool in “desingularising” certain geometric structures; see e.g. \cite{proper}. A comprehensive study of blow-ups is on its place, and in the thesis we will dive into the theory of blow-ups and view it from many different angles discussed in a variety of literature from the last couple of years.

Different classes of examples come into play when looking at the Lie groupoid and Lie algebroid blow-up constructions. One can purely view these blow-up constructions as yielding new types of examples of both Lie groupoids and Lie algebroids (see e.g. \cite{Songhao}).
But other applications include so-called ``symplectic groupoids'', ``log-tangent bundles'' and ``scattering tangent bundles'' \cite{Melinda}, which are used in the study of Poisson manifolds (and their cohomology theory). In this thesis we will focus on the general theory, and as a result we can apply the blow-up construction to many different settings, such as foliations and differentiable stacks.

What is brought to the attention in \cite{2017arXiv170509588D}, is that one can approach blow-ups via another (common) construction in non-commutative geometry: the deformation to the normal cone construction; yet another construction widely used in algebraic geometry. While needing to explain the theory of deformation to the normal cones, other, formerly hard (or unknown) results regarding blow-ups come to light much easier, or can even be proven immediately. With our mind set on the theory of Lie groupoids and Lie algebroids, we will see that it is a natural step to “factor through” the deformation to the normal cone. %While doing this, we will see that, 
Not only are these types of deformation spaces useful for our understanding of blow-ups, they provide us with a powerful tool to shed new light on many widely used theorems like the Darboux theorem, the Weinstein theorem, and the splitting theorem for Lie algebroids (see e.g. \cite{Bischoff_2020,sadegh2018eulerlike}). In this direction, we will digress into some of the applications.

While some differential geometers are more and more interested in blow-ups, there are only a few sources focused on discussing the basic theory of blow-ups, especially concerning blow-ups of Lie groupoids and Lie algebroids. In fact, there seems to be no geometric approach to the blow-up of Lie algebroids (for the codimension $1$ case, in \cite{Melinda} and 
\cite{gualtieri2012symplectic} there is an algebraic construction of a Lie algebroid blow-up, and in \cite{2017arXiv170509588D} it is mentioned what the Lie algebroid of a Lie groupoid is). In this thesis we will provide background on Lie groupoids and Lie algebroids, we will broadly discuss the deformation to the normal cone construction, and explain how we can construct new Lie groupoids and new Lie algebroids out of it, and do so for blow-ups as well, again, both in the Lie groupoid and Lie algebroid case. In addition, we will digress a number of times, providing the interested reader with more background, or create awareness for related topics. For example, we will go into the theory of foliations and holonomy, the theory of $\VB$-groupoids and $\VB$-algebroids and applications of the theory of deformation to the normal cones; most notably to the splitting theorem for Lie algebroids. More precisely:
\begin{itemize}
    \item Chapter \ref{sec: Lie groupoids} and \ref{sec: Lie algebroids} give a general introduction to the theory of Lie groupoids and Lie algebroids;
    \item Chapter \ref{sec: deformation to the normal cone} gives a thorough introduction to deformation to the normal cones, with digressions into $\VB$-groupoids, double vector bundles, $\VB$-algebroids, Euler-like vector fields and applications of deformation to the normal cones; 
    \item Chapter \ref{sec: Blow-up of a pair of Lie groupoids} is dedicated to blow-ups of manifolds, Lie groupoids, Lie algebroids, and differentiable stacks.
\end{itemize}
It is important to mention that this thesis is far from complete in the sense that there are things that deserve consideration and are not included. The contents of the thesis are such that the main goal is to provide a broad and thorough introduction to blow-ups. Therefore, there is little to no motivation to be found why one would want to learn about blow-ups. However, as already mentioned in this introduction, there are many (hidden) applications of blow-ups already in the literature, and blow-ups seem like a powerful tool that can be of great importance for certain applications. %At the end a discussion can be found in which we explain known results, and new directions one can apply, or mimic, the blow-up constructions to new situations that might be of interest to the reader.

%Lastly, unfortunately, there are no pictures in the thesis, but I can only encourage to draw some of the simple examples for yourself.

I hope you have a good time reading and learning about the things I learned about during the last year.

\section*{Acknowledgements}
I would like to thank Ioan M\v{a}rcu\c{t}, my thesis supervisor, for all the (online) discussions, help, and advice. I am especially very happy that I got the room to learn about the things that I found interesting myself, and that I could always ask for help whenever I got stuck on something. Of course many thanks to my friends and family for the great times I've had in my time as a (master's) student. Also a big thank you to Klaas Landsman, who is so nice to be the second reader of this thesis. Lastly, I would like to express my appreciation for Fenna Sanders who made the very nice title page artwork.
\addtocontents{toc}{\protect\thispagestyle{myplain}}\newpage

\section{Lie Groupoids}\label{sec: Lie groupoids}
\cfoot{\thepage}
In this chapter we will explain the basics of the theory of Lie groupoids. The theory is mainly based on \cite{ruimarius,mackenziebook,meinrenken,moerdijk_mrcun_2003} (in every definition, remark, etc. we will give a reference which fits the content the best). Often we think of a groupoid as a generalisation of a group; instead of allowing just one point of symmetry (the identity element in a group) we allow for more. One way to formalise this is by generalising the notion of a group, thought of as a category with only one object, such that every arrow has an inverse.
%\begin{figure}[h!]
%\centering

%\definecolor{ffffww}{rgb}{1.,1.,0.4}
%\definecolor{ffzzcc}{rgb}{1.,0.6,0.8}
%\begin{tikzpicture}[line cap=round,line join=round,>=triangle 45,x=1.0cm,y=1.0cm]
%\clip(-2.,-0.1) rectangle (2.,3.9);
%\draw [line width=2.pt,color=ffffww] (0.,0.4426877757896084) circle (0.4426877757896084cm);
%\draw [line width=2.pt,color=ffffww] (0.,0.8853755515792168) circle (0.8853755515792168cm);
%\draw [line width=2.pt,color=ffffww] (0.,1.6070346529164505) circle (1.6070346529164505cm);
%\draw [line width=2.pt,color=ffffww] (0.,1.2535722080743734) circle (1.2535722080743734cm);
%\draw [line width=1.2pt,dotted,color=ffffww] (0.,3.314069305832901)-- (0.,3.85);
%\begin{scriptsize}
%\draw [fill=ffzzcc] (0.,0.) circle (3.0pt);
%\draw[color=ffzzcc] (0.15267053632197786,0.43229182152954737) node {$e$};
%\end{scriptsize}
%\end{tikzpicture}

%    \caption{A group visualised as a category. The unique object $e$ can be thought of as the unit element and the circles represent two arrows of the category: "moving left" is inverse to "moving right".}
%    \label{fig: group_as_category}
    
%\end{figure}

\begin{defn}\cite{moerdijk_mrcun_2003}\label{defn: groupoid}
A \textit{groupoid} is a small category such that every arrow has an inverse.
\end{defn}

Since our focus lies on groupoids in the smooth category, i.e. Lie groupoids, it is useful to unravel the above definition in a specific way that allows us to indicate how we should equip a groupoid with extra structure so that it becomes a Lie groupoid. 

\begin{rema}\cite{moerdijk_mrcun_2003}\label{rema: groupoid spelled out definition}
Let $\groupoid$ be a groupoid, where $\group$ denotes the set of arrows and $\base$ the set of objects. We can associate so-called \textit{structure maps} to this groupoid:
\[\group \xrra[\target]{\source} \base, \hspace{0.1cm} \group^{(2)} \xra{\mult} \group, \hspace{0.1cm} \group \xra{\inv} \group, \hspace{0.1cm} \base \xra{\identity} \group,\]
where 
\[\group^{(2)} \coloneqq \group \tensor[_\source]{\times}{_\target} \group = \{(g,h) \in \mathcal{G} \times \mathcal{G} \mid \source(g) = \target(h)\} \subset \group \times \group,\] 
which are defined as follows: 
\begin{itemize}
\item The \textit{source map} $\source$ and \textit{target map} $\target$ map an arrow to its source and target respectively. 
\item The \textit{multiplication map} $\mult$ maps a pair of arrows $(g,h)$ to the composition $gh$. This map is well-defined, because the domain is the set $\group^{(2)}$ consisting of \textit{composable arrows}.
\item The \textit{inverse map} $\inv$ maps an arrow $g$ to its inverse $g^{-1}$.
\item The \textit{unit map} $\identity$ maps an object $x$ to the identity morphism $\id_x$ of that object.
\end{itemize}
It is readily verified that these structure maps satisfy the following relations:
\begin{itemize}
    \item $\source(gh) = \source(h)$ and $\target(gh) = \target(g)$ for all $(g,h) \in \group^{(2)}$,
    \item $\source(g) = \target(g^{-1})$ for all $g \in \group$,
    \item $(gh)k = g(hk)$ for all $(g,h),(h,k) \in \group^{(2)}$, \textnormal{and }
    \item $\id_{\target(g)}g = g = g\id_{\source(g)}$ for all $g \in \group$.
\end{itemize}
Conversely, given a pair of sets $(\group,\base)$ together with structure maps $(\source,\target,\mult,\inv,\identity)$ satisfying these relations, then this defines a groupoid.
\end{rema}

\begin{term}\label{term: groupoids}
Whenever we talk about a groupoid $\groupoid$, the structure maps, as defined above, are written as $(\source_\group,\target_\group,\mult_\group,\inv_\group,\identity_\group)$ (or without lower indices if no confusion is possible). We will call $\base$ the base manifold, and simply call $\group$ the groupoid (if no confusion is possible). In some examples we have a different symbol to denote the groupoid $\groupoid$.
\end{term}

\begin{defn}\cite{moerdijk_mrcun_2003}\label{defn: Lie groupoid}
A Lie groupoid is a groupoid $\groupoid$ for which $\group$ is a, possibly non-Hausdorff, smooth manifold, $\base$ is a (Hausdorff) smooth manifold, all structure maps are smooth, and the source map $\source$ is a submersion.
\end{defn}
In the above definition, we require the source map to be a submersion. The reason for imposing this condition is to ensure that all the structure maps are well-defined smooth maps. Such a strong transversality assumption is unnecessary if this would be the only reason. However, we will use that the source map is a submersion a number of times, as we will soon see.
\begin{prop}\cite{ruimarius}
If $\groupoid$ is a Lie groupoid, then the inverse map $\inv$ is a diffeomorphism, the target map $\target$ is a submersion and the unit map $\identity$ is an embedding.
\end{prop}
\begin{proof}
Since $\inv: \mathcal{G} \ra \mathcal{G}$ is smooth and inverse to itself (i.e. an involution), it is a diffeomorphism. The second statement is now a consequence of the relation $\target = \source \circ \inv$. The last statement follows by the relation $\id_\base = \source \circ \identity$, since it shows that $\identity$ has a global left inverse; it is an injective immersion and a topological embedding, so it is a smooth embedding.
\end{proof}
\begin{rema}\cite{ruimarius}\label{rema: Lie groupoid}
As already hinted upon, there could be an ambiguity in the above definition: for a groupoid $\groupoid$ to be a Lie groupoid, we wrote that the structure maps, in particular the multiplication map
\[\mult: \group^{(2)} \ra \group \textnormal{ given by } (g,h) \mapsto gh,\]
should be smooth; this would not make sense if $\group^{(2)}$ is not a smooth manifold. That it is a smooth manifold follows from $\source$ and $\target$ being submersions. Indeed, suppose $\groupoid$ is a groupoid for which $\group$ is a (possibly non-Hausdorff) smooth manifold, $\base$ is a (Hausdorff) smooth manifold, and the structure maps $\source$ and $\target$ are smooth submersions. Consider the map
\[\source \times \target: \group \times \group \ra \base \times \base.\]
Then $\source \times \target$ is a submersion, so $\group^{(2)} = (\source \times \target)^{-1}(\Delta_\base)$, where $\Delta_\base = \{(x,x) \in \base \times \base\}$ is the diagonal of $\base$, is an embedded submanifold of $\group$. Notice that $\group$ is non-Hausdorff if and only if $\group^{(2)}$ is non-Hausdorff. Indeed, suppose we have a sequence $(g_n)$ in $\group$ that has two limits $g,h \in \group$. Since $\source(g_n) \ra \source(g)$ and $\source(g_n) \ra \source(h)$, we see that $\source(g)=\source(h)=\target(h^{-1})=\target(g^{-1})$, because $\base$ is Hausdorff. This shows that $(g_n,g_n^{-1})$ is a sequence in $\group^{(2)}$ that has the two limits $(g,h^{-1})$ and $(h,g^{-1})$, and these limits are equal if and only if $g=h$.
Even though we allow for $\group$ to be non-Hausdorff, \textbf{we do require all source fibers to be Hausdorff}. It follows that the target fibers are also Hausdorff: if $(g_n)$ is a sequence in $\target^{-1}(x)$ with two limits $g,h \in \target^{-1}(x)$, then $(g_n^{-1})$ is a sequence in $\source^{-1}(x)$, so $g^{-1}=h^{-1}$, and therefore $g=h$. 
%Note that, even though $\group \times \group$ might not be Hausdorff, $\group^{(2)}$ is. To see this, observe first that the diagonal of $\identity(\base) \subset \group$ is closed, as $\base$ is Hausdorff and $\identity$ is an embedding. Moreover, 
%\[\group^{(2)} \tensor[_{\source \times \source}]{\times}{_{\target \times \target}} \group^{(2)} = \{(g,h,g',h') \in \group^{(2)} \times \group^{(2)} \mid \source(g)=\source(g') \textnormal{ and } \target(h) = \target(h')\}\]
%is a closed subset of $\group^{(2)} \times \group^{(2)}$, hence $\Delta_{\group^{(2)}} \subset \group^{(2)} \times \group^{(2)}$?????????????????
%is closed being the pre-image of $\identity(\base) \times \identity(\base)$ of the continuous map 
%\[\group^{(2)} \tensor[_{\source \times \source}]{\times}{_{\target \times \target}} \group^{(2)} \xra{f} \group^{(2)} \times \group^{(2)} \xra{\mult \times \mult} \group \times \group,\]
%where $f$ is the continuous map $(g,h,g',h') \mapsto (g,\inv(g'),\inv(h'),h)$. 
\end{rema} 
Since we allow the total space of a groupoid to be non-Hausdorff, we have to be a bit careful with applying some of the standard theory of smooth manifolds. Typically, a non-Hausdorff manifold does not admit partitions of unity (take for example the line with a double origin). Therefore, in general, Stokes's theorem fails, and we can not make use of the correspondence between manifolds and its algebra of smooth functions. It will be useful to recall some equivalent definitions of a topological space to be Hausdorff. 
\begin{prop}\label{prop: equivalent definitions Hausdorff}
A topological space is $\base$ is Hausdorff if and only if
\begin{enumerate}
    \item every convergent sequence $(x_n)$ in $\base$ has a unique limit;
    \item if $\basetwo$ is a topological space, and $\subbasetwo \subset \basetwo$ is a dense subspace, then a continuous map $\subbasetwo \ra \base$ which extends to a continuous map $\basetwo \ra \base$ extends uniquely.
\end{enumerate}
\end{prop}
%\begin{proof}
%$1.$ If $x,y \in \base$ are distinct and inseparable by open subsets, then $(x_n) \coloneqq (x)$ is a sequence converging to both $x$ and $y$. Conversely, if $(x_n)$ converges to distinct elements $x,y \in \base$, then any two open subsets $U \ni x$ and $V \ni y$ of $\base$ eventually contain $(x_n)$, and therefore have a non-empty intersection. 
%\end{proof}

Here is an equivalent characterisation of a groupoid $\groupoid$ having a Hausdorff total space $\group$.
\begin{prop}\cite{hoyohausdorff}\label{prop: G Hausdorff iff X closed}
Let $\groupoid$ be a Lie groupoid. Then $\group$ is Hausdorff if and only if $\identity(\base) \subset \group$ is a closed subspace.
\end{prop}
\begin{proof}
Suppose $(g_n)$ is a sequence in $\group$ with two distinct limits $g,h \in \group$ (i.e. $\group$ is not Hausdorff). Then $\source(g_n) \ra \source(g)$ and $\source(g_n) \ra \source(h)$, so, since $\base$ is Hausdorff, $\source(g)=\source(h)$%, and similarly $\target(g)=\target(h)$
. The sequence
\[(\id_{\target(g_n)}) = (g_ng_n^{-1})\]
in $\identity(\base)$ has limit $gh^{-1}$ (note: $(g,h) \in \group^{(2)}$ because $\source(g)=\source(h)$), which does not lie in $\identity(\base)$; otherwise $g=h$. This shows that $\identity(\base) \in \group$ is not closed. Conversely, if $\identity(\base) \in \group$ is not closed, then there exists a sequence $(x_n)$ in $\base$ such that $\id_{x_n} \ra g$ for some $g \not\in \identity(\base)$. But the sequence
\[(x_n) = (\target(\id_{x_n}))\] 
also has limit $\target(g)$, so $(\id_{x_n})$ has two distinct limits: $g$ and $\id_{\target(g)}$. This shows that $\group$ is not Hausdorff, as required.
\end{proof}
Of course, there is an appropriate category in which the objects are Lie groupoids and the morphisms are ``Lie groupoid-structure preserving''. 
\begin{defn}\cite{meinrenken}\label{defn: morphisms of groupoids}
Let $\groupoid$ and $\subgroupoid$ be Lie groupoids. A morphism from $\groupoid$ to $\subgroupoid$, denoted simply by $\group \ra \subgroup$, is a functor for which the map between arrows %and the map between objects are 
is smooth.
\end{defn}
\begin{rema}\cite{meinrenken}\label{rema: morphism of groupoids}
Let $\groupoid$ and $\subgroupoid$ be Lie groupoids. Then a morphism $F: \group \ra \subgroup$ comes with a \textit{base map} $f: \base \ra \subbase$, and yields a commutative diagram
\begin{center}
\begin{tikzcd}
    \group \ar[r,"F"] \ar[d,shift left, "\source_\group"] \ar[d,shift right,"\target_\group"'] & \subgroup \ar[d,shift left, "\source_\subgroup"] \ar[d,shift right,"\target_\subgroup"'] \\
    \base \ar[r,"f"] & \subbase.
\end{tikzcd}
\end{center}
More precisely, we obtain two commutative diagrams: one with vertical arrows the target maps, and one with vertical arrows the source maps. Since, for example, $\source_\group$ is a submersion, the base map $\base \ra \subbase$ is smooth. Besides the existence of the above commutative diagram(s), $F$ and $f$ satisfy the relations
\begin{itemize}
    \item $F(gh)=F(g)F(h)$ for all $(g,h) \in \group^{(2)}$, \textnormal{ and}
    \item $F(\id_{x})=\id_{f(x)}$ for all $x \in \base$.
\end{itemize}
Conversely, a smooth pair of maps $(F: \group \ra \subgroup,f: \base \ra \subbase)$ fitting in a commutative diagram as above, and satisfying these relations, defines a morphism of Lie groupoids $\group \ra \subgroup$. However, it suffices to show that $F(gh)=F(g)F(h)$ for all $(g,h) \in \group^{(2)}$, since then 
\[F(\id_x)=F(\id_x\id_x)=F(\id_x)F(\id_x)\]
shows that $F(\id_x)=\id_{\source(F(\id_x))}=\id_{f(x)}$. Notice that a morphism of Lie groupoids $F: \group \ra \subgroup$ is an isomorphism of Lie groupoids if and only if $F$ (that is, the map $\group \ra \subgroup$ on total spaces) is a diffeomorphism. Indeed, if $F$ is a diffeomorphism, then
\begin{align*}
    F^{-1}(gh)&=F^{-1}(F(F^{-1}(g))F(F^{-1}(h))) \\
    &= F^{-1}(F(F^{-1}(g)F^{-1}(h))) = F^{-1}(g)F^{-1}(h),
\end{align*}
which shows that $F^{-1}$ is a morphism of Lie groupoids.
\end{rema}
\begin{defn}\cite{meinrenken}\label{defn: subgroupoid}
A \textit{subgroupoid} $\subgroupoid$ of a groupoid $\groupoid$ is a subcategory such that every arrow has an inverse. If $\groupoid$ is a Lie groupoid, and, in addition, the inclusion functor $\subgroup \hookrightarrow \group$ is a morphism of Lie groupoids, then $\subgroupoid$ is called a \textit{Lie subgroupoid}.
\end{defn}
\begin{rema}\label{rema: subgroupoid}
%\faWarning\hspace{0.1cm} In the above definition, observe that the requirement for the inclusion $\subgroup \hookrightarrow \group$ to be an embedding is not standard, since it is stronger then requiring the inclusion to be a morphism of Lie groupoids. Already in the case of Lie groups (i.e. Lie groupoids over a point), Lie subgroups can be immersed submanifolds; in fact, a Lie subgroup is an embedded submanifold if and only if it is a closed subset. However, in this thesis we are interested in blow-ups of Lie groupoids along closed embedded Lie subgroupoids.
Let $\subgroupoid$ be a Lie subgroupoid of a Lie groupoid $\groupoid$. Then, as manifolds, $\iota: \subgroup \hookrightarrow \group$ is an immersed submanifold. Conversely, a subgroupoid $\subgroup$ of a Lie groupoid $\group$, such that $\subgroup \subset \group$ is an immersed submanifold, is a Lie subgroupoid.
\end{rema}
In the theory of blow-ups, we will only blow up a particular nice class of Lie subgroupoids of a Lie groupoid. We will use the following terminology for this.
\begin{defn}\label{defn: pair of Lie subgroupoids}
A \textit{pair of manifolds} is a tuple $(\base,\subbase)$ such that $\subbase \subset \base$ is a closed embedded submanifold. A \textit{pair of Lie groupoids} is a tuple $(\groupoid,\subgroupoid)$ such that $\subgroupoid$ is a subgroupoid of $\groupoid$ for which $\subgroup \subset \group$ and $\subbase \subset \base$ are closed embedded submanifolds. 
\end{defn}

Later, we will see that, in analogy to Lie groups, to any Lie groupoid we can associate a \textit{Lie algebroid}. This will be discussed later, as will the theory of Lie algebroids.
Next, we will shortly explain some of the basic theory of Lie groupoids, and discuss examples of Lie groupoids.

\begin{term}
From now on we will call a groupoid, without extra structure, a set-theoretic groupoid, and call a Lie groupoid a groupoid (or still a Lie groupoid).
\end{term}

\subsection{Isotropy Lie groups and orbits}\label{sec: Isotropy groups and orbits}

One of the important tools when working with (Lie) groupoids are local bisections. Here we will define the concept of local bisection, and see that every groupoid has many of them. We will apply this theory to show that the target map restricted to a source-fiber has constant rank. From this we can deduce a number of facts, and it allows us to introduce the isotropy groups and orbits associated to a Lie groupoid. This section is based on a treatment in \cite{meinrenken}. We fix a groupoid $\groupoid$.

\begin{defn}\cite{meinrenken}\label{defn: local bisection}
Let $g \in \group$ and let $U \subset \base$ be an open subset. A smooth local section $\sigma: U \ra \group$ of $\source$ through $g$ (i.e. $g \in \sigma(U)$ and $\source \circ \sigma = \id_U$) is called a \textit{local bisection through $g$} if $\target \circ \sigma$ is a diffeomorphism onto an open subset of $\base$. 
%If $\sigma$ is a global section of $\source$, then it is called a (global) bisection.
\end{defn}
We can find a local bisection through every arrow of a groupoid. It rests on the following simple fact from linear algebra.
\begin{lemm}\cite{discretemechanics}\label{lemm: complementary vector space}
Let $V$ be a vector space and let $W,W' \subset V$ be subspaces of equal dimension. Then there is a subspace $C \subset V$ such that $W \oplus C = V = W' \oplus C$.
\end{lemm}
\begin{proof}
Choose a basis $\{b_i\}$ of $W \cap W'$, complete it to bases $\{b_i,w_j\}$ of $W$ and $\{b_i,w_j'\}$ of $W'$ (we use the same index $j$ because $\dim W = \dim W'$), and then complete the basis $\{b_i,w_j,w_j'\}$ of $W + W'$ to a basis $\{b_i,w_j,w_j',v_k\}$ of $V$. The result follows by setting $C \coloneqq \textnormal{span}(w_j + w_j',v_k)$.
\end{proof}

\begin{prop}\cite{discretemechanics}\label{prop: local bisection through g}
Let $g \in \group$. Then, there is a local bisection through $g$.
\end{prop}
\begin{proof}
Since $\source$ and $\target$ are submersions, $\dim \ker d\source(g) = \dim \group - \dim \base = \dim \ker d\target(g)$.
Hence, Lemma \ref{lemm: complementary vector space} applies: there exists a subspace $C \subset T_g\group$ such that 
\[\ker d\source(g) \oplus C = T_g\group = \ker d\target(g) \oplus C.\]
%Take any chart $\varphi: V \ra \mathbb{R}^n$ through $g \in \group$ such that $d\varphi(g)$ maps $C$ isomorphically to $\mathbb{R}^k \times \{0\}$ (by possibly post-composing with a linear endomorphism of $\mathbb{R}^n$). Then the submanifold $\Gamma \coloneqq \varphi^{-1}(\mathbb{R}^k \times \{0\}) \subset V$ of $\group$ 
Take a submanifold $S \ni g$ of $\group$ for which $T_gS
%= d\varphi(g)^{-1}(\mathbb{R}^k \times \{0\}) 
= C$. Now observe that there is an open subset $U \ni g$ of $S$ for which $\target|_{U}: U \ra \target(U)$ and $\source|_{U}: U \ra \source(U)$ are both diffeomorphisms onto open subsets of $\base$ (by the inverse function theorem). The inverse of the latter map %(post-composed with the inclusion $U \hookrightarrow \group$) 
gives the desired local bisection.
\end{proof}
\begin{rema}
Observe that, alternatively, through every $g \in \group$ there exists a smooth local section $\tau: U \ra \group$ of $\target$ such that $\source \circ \tau$ is a diffeomorphism onto its image. 
\end{rema}
Using the notion of local bisection, we can prove the following useful result.
\begin{prop}\cite{meinrenken}\label{prop: t restricted to source-fiber has constant rank}
Let $x \in \base$. Then $\target|_{\source^{-1}(x)}: \source^{-1}(x) \ra \base$ has constant rank. 
\end{prop}
\begin{proof}
Let $g,h \in \source^{-1}(x)$, and choose a bisection $\sigma: U \ra \sigma(U)$ through $gh^{-1}$. Define $V \coloneqq \target(\sigma(U))$ (note that this set is open in $\base$), and consider the ``left-translation by $\sigma$'':
\[f \coloneqq \mult(\sigma \circ \target(\cdot), \cdot)|_{\target^{-1}(U)}: \target^{-1}(U) \ra \target^{-1}(V), \textnormal{ } a \mapsto \sigma \circ \target(a) \cdot a.\]
Observe that $f$ indeed maps into $\target^{-1}(V)$, since $\target(ab) = \target(a)$ for all $(a,b) \in \group^{(2)}$, and that $f$ is well-defined, because $\source \circ \sigma \circ \target = \target$. Moreover, $h \in \target^{-1}(U)$, because $gh^{-1} \in \sigma(U)$, so that $\target(h) = \source(gh^{-1}) \in U$. Now, $f$ maps $h$ to $g$, and we will show that it is even a diffeomorphism with inverse given by the left-translation by $\inv \circ \sigma \circ (\target \circ \sigma)^{-1}$:
\[f' \coloneqq \mult(\inv \circ \sigma \circ (\target \circ \sigma)^{-1} \circ \target(\cdot), \cdot)|_{\target^{-1}(V)}: \target^{-1}(V) \ra \target^{-1}(U).\] 
Observe that $f'$ indeed maps into $\target^{-1}(U)$, since $\target \circ \inv \circ \sigma = \source \circ \sigma = \id_U$, and that $f'$ is well-defined, because $\source \circ \inv = \target$. To see that $f'$ is indeed inverse to $f$, let $b \in \target^{-1}(V)$, and $a \in \target^{-1}(U)$. Then, indeed,
\begin{align*}
    f\circ f'(a) &= \mult\left(\sigma \circ \target\left(\inv \circ \sigma \circ (\target \circ \sigma)^{-1} \circ \target(a)\right),\mult(\inv \circ \sigma \circ (\target \circ \sigma)^{-1} \circ \target(a),a)\right) \\
    &= \mult\left(\sigma \circ (\target \circ \sigma)^{-1} \circ \target(a),\mult(\inv \circ \sigma \circ (\target \circ \sigma)^{-1} \circ \target(a),a)\right) = a, \textnormal{ and}\\
    f' \circ f(b) &= \mult\left(\inv \circ \sigma \circ (\target \circ \sigma)^{-1} \circ \target \circ \sigma \circ \target(b), \mult(\sigma \circ \target(b), b)\right) \\
    &= \mult\left(\inv \circ \sigma \circ \target(b), \mult(\sigma \circ \target(b), b)\right) = b.
\end{align*}
Since $\source \circ f = \source$, we also have that $f|_{\source^{-1}(x)}: \target^{-1}(U) \cap \source^{-1}(x) \ra \target^{-1}(V) \cap \source^{-1}(x)$ is a diffeomorphism. We now have a commutative diagram:
\begin{center}
\begin{tikzcd}
\target^{-1}(U) \cap \source^{-1}(x) \arrow[r, "f|_{\source^{-1}(x)}"] \arrow[d, "\target|_{\target^{-1}(U) \cap \source^{-1}(x)}"] & \target^{-1}(V) \cap \source^{-1}(x) \arrow[d,"\target|_{\target^{-1}(V) \cap \source^{-1}(x)}"] \\
U \arrow[r,"\target \circ \sigma"]& V.
\end{tikzcd}
\end{center}
Since $f|_{\source^{-1}(x)}$ and $\target \circ \sigma$ are diffeomorphisms, we see that $\target|_{\target^{-1}(U) \cap \source^{-1}(x)}$ must have the same rank as $\target|_{\target^{-1}(V) \cap \source^{-1}(x)}$, proving that $\target|_{\source^{-1}(x)}$ has constant rank, as claimed.
\end{proof}

The rest of this section is devoted to treating corollaries of this proposition. We fix $x,y \in \base$.

\begin{coro}\cite{meinrenken}\label{coro: isotropy groups are Lie groups}
The subset $\group(x,y) = \source^{-1}(x) \cap \target^{-1}(y)$ is an embedded submanifold of $\group$. In particular, $\group_x \coloneqq \group(x,x)$ is a Lie group.
\end{coro}
\begin{proof}
Since $\target|_{\source^{-1}(x)}$ has constant rank, it follows that $\group(x,y) = \left(\target|_{\source^{-1}(x)}\right)^{-1}(y)$ is an embedded submanifold of $\source^{-1}(x)$, which in turn is an embedded (Hausdorff) submanifold of $\group$. Clearly, $\group_x$ now is a Lie group with multiplication map $\mult|_{\group_x \times \group_x}$, and inverse map $\inv|_{\group_x}$. 
\end{proof}
\begin{defn}\cite{meinrenken}\label{defn: isotropy group of a Lie groupoid}
The Lie group $\group_x$ of $\groupoid$ is called the \textit{isotropy group} of $x$.
\end{defn}
\begin{coro}\cite{meinrenken}\label{coro: orbits are immersed submanifolds}
The subset $\orbit_x \coloneqq \target(\source^{-1}(x))$ of $\base$ is an immersed submanifold. Moreover, the map $\target|_{\source^{-1}(x)}: \source^{-1}(x) \ra \orbit_x$ is naturally a principal $\group_x$-bundle.
\end{coro}
\begin{proof}
We claim that the action of $\group_x$ on $P_x \coloneqq \source^{-1}(x)$ by right multiplication is free and proper. Using the claim, we obtain a principal $\group_x$-bundle $\pi: P_x \ra P_x/\group_x$. Since $\target|_{P_x}$ is $\group_x$-equivariant, i.e. $\target(gh)=\target(g)$ for all $g \in P_x$ and $h \in \group_x$, it factors through the map
\[\widetilde\target: P_x/\group_x \ra \base \textnormal{ given by } g \cdot \group_x \mapsto \target(g).\]
Notice that $\widetilde\target$ is injective, since if $\target(g) = \target(h)$, then $g=h\cdot (h^{-1}g)$, and $\widetilde\target$ has constant rank, because $\target|_{P_x}$ and $\pi$ have constant rank. Hence, $\widetilde \target$ is an immersion. Since $\widetilde\target$ maps onto $\orbit_x$, it remains to prove the claim. It is almost immediate that the action is free: if $hg = h$, then $g = h^{-1}h = \id_x$. To see that the action is proper, we will prove that, whenever $C,D \subset P_x$ are compact subsets, then the set
\[(\group_x)_{C,D} \coloneqq \{h \in \group_x \mid Ch \cap D \neq \emptyset\}\] 
is a compact subset of $\group_x$. This is clear once we recognise this set as the image of the closed (and hence compact) subset
$C \tensor[_\target]{\times}{_\target} D \coloneqq \{(c,d) \in C \times D \mid \target(c) = \target(d)\}$ of $C \times D$, subject to the continuous map $P_x \tensor[_\target]{\times}{_\target} P_x \ra \group_x$ 
given by $(c,d) \mapsto c^{-1}d$. This proves the statement.
\end{proof}
%\begin{rema}
%From the proof (we will use the same notation from the proof) we can see that $\orbit_x$ is even a regular immersed submanifold, i.e. if $f: X \ra M$ is a smooth map that maps into $\orbit_x$, then the induced map $\widetilde f: X \ra P_x/\group_xx$ is smooth. Indeed, it suffices to prove that $\widetilde f$ is continuous (since $\widetilde \target$ has local left-inverses), which follows from the commutative diagram
%\begin{center}
%\begin{tikzcd}
%& X \ar[d,"\widetilde f"] \ar[rd, "f"] & \\
%P_x \ar[rr,"\target",out=-30,in=210] \ar[r,"\pi"] & P_x/\group_x \ar[r,"\widetilde \target"] & \base
%\end{tikzcd}
%\end{center}
%so that for any open set $U \subset P_x/\group_x$ we have that $\widetilde f^{-1}(U) = f^{-1}\left(\target\left(\pi^{-1}(U)\right)\right)$ is open.
%\end{rema}
\begin{defn}\cite{meinrenken}\label{defn: orbit of a Lie groupoid}
The immersed submanifold $\orbit_x$ of $\base$ is called the \textit{orbit} of $x$. If $\orbit_x = \base$, then $\groupoid$ is called \textit{transitive}.
\end{defn}
Of course, from Corollary \ref{coro: orbits are immersed submanifolds}, we obtain the following result about transitive groupoids.
\begin{coro}\cite{meinrenken}\label{coro: if groupoid is transitive, then target is surjective submersion}
If $\group$ is transitive, then $\target|_{\source^{-1}(x)}: \source^{-1}(x) \ra \base$ is a principal $\group_x$-bundle.
\end{coro}
In the next section we will see that every principal bundle arises in this way; in fact, there is an equivalence of categories between the category of transitive groupoids and the category of principal bundles.
%\begin{proof}
%That $\group$ is transitive, means that $\target|_{\source^{-1}(x)}$ is surjective. In particular, this map is a submersion at some $g \in \source^{-1}(x)$ (e.g. by Sard's theorem), so it must be a submersion everywhere, because it is of constant rank.
%\end{proof}

\subsection{Examples of Lie groupoids}\label{sec: Examples of Lie groupoids}
Here, we discuss some basic examples of Lie groupoids. These introductions to examples are kept short. 
%Later, when we discuss examples of blow-ups, we will naturally go deeper into the theory that we will need. 
Most of these groupoids will merely be a ``rewriting'' of a certain construction that we are familiar with, like actions of Lie groups, path spaces and principal bundles. Later, we will also discuss the basics of foliation theory, so that we can properly introduce the so-called \textit{monodromy and holonomy groupoids}, which produces examples of groupoids that have a non-Hausdorff space of arrows. This section is based on \cite{ruimarius,meinrenken}.

The first example that we will treat is the appropriate ``trivial'' object in the category of groupoids. As such, it is probably one of the most important examples discussed here.

\begin{exam}[Pair groupoids]\cite{ruimarius}\label{exam: pair groupoids}
To any manifold $M$ we can, canonically, associate the \textit{pair groupoid}
\[M \times M \rra M,\]
with $\source(x,y)=y$, $\target(x,y)=x$, and $\mult((x,y),(y,z)) = (x,z)$.
\end{exam}
By viewing a manifold $M$ as a principal $\{e\}$-bundle, where $\{e\}$ is viewed as a trivial group, we can generalise the above example as follows:
\begin{exam}[Gauge groupoids]\cite{ruimarius}\label{exam: Gauge groupoids}
Let $G$ be a Lie group, and let $P$ be a principal $G$-bundle. Consider $P \times P$ with the diagonal action of $G$. We can associate to $P$ the \textit{Gauge groupoid}
\[P \otimes_G P: (P \times P)/G \rra P/G\]
with $\source[p,q] = [q]$, $\target[p,q] = [p]$, and $\mult\left([p,q],[q,r]\right) = [p,r]$. Notice that $P \otimes_G P$ is a transitive groupoid. Moreover, for all $x \in M$, if we pick $p \in \source^{-1}(x)$, then we have an isomorphism $G \xra{\sim} \group_x$ given by $g \mapsto [pg,p]$.
\end{exam}
As already mentioned before, by using Corollary \ref{coro: if groupoid is transitive, then target is surjective submersion}, we can, and will, prove that there is an equivalence of categories between the category of transitive groupoids and the category of principal bundles. 
The reason we include a proof of this equivalence is to stress that the theory of Lie groupoids is rich, and provides us with a framework to study many familiar geometric structures that describe symmetries, such as principal bundles, all at once. Moreover, principal bundles will arise naturally in the theory of blow-ups, and it is good to be aware of this fact.

\begin{rema}\label{rema: transitive groupoids are in one to one correspondece with Gauge groupoids}[Transitive groupoids $\rla$ Gauge groupoids]
We start by extending the assignment ($P \ra M$ principal $G$-bundle) $\mapsto P 
\otimes_G P$ to a functor $\mathcal{F}$. Let $P \ra M$ be a principal $G$-bundle, and let $Q \ra N$ be principal $H$-bundle. A morphism $\Phi_\varphi$ from $P \ra M$ to $Q \ra N$ is given by a $G$-equivariant smooth map 
\[\Phi: P \ra Q\]
with respect to a Lie group morphism
\[\varphi: G \ra H\]
(that is, the $G$-action on $Q$ is given by the action that $\varphi$ induces). Such a morphism induces a groupoid morphism by setting
\[\mathcal{F}(\Phi_{\varphi}): P \otimes_{\group_x} P \ra Q \otimes_{\subgroup_x} Q, \textnormal{ given by } [p,p'] \mapsto [\Phi(p),\Phi(p')]\]
on the level of arrows, and the corresponding base map is the map $P/G \ra Q/H$ that $\Phi$ induces. From the definition of this morphism, it is obvious that $\mathcal{F}$ becomes a functor this way. Moreover, notice that, for $p \in P$, 
\[G \xra{\sim} (P \otimes_{\group_x} P)_{[p]} \textnormal{ by } g \mapsto [pg,p], \textnormal{ and } P \xra{\sim} \source_{P \otimes_{\group_x} P}^{-1}([p]) \textnormal{ by } p' \mapsto [p',p],\] 
and similarly for $q=\Phi(p) \in Q$. In particular, we can extract $\varphi: G \ra H$ and $\Phi_\varphi: P \ra Q$ from the Lie groupoid morphism $\mathcal{F}(\Phi_\varphi)$: using these isomorphisms,
\[\varphi: G \ra H \textnormal{ is given by } g \mapsto \mathcal{F}(\Phi_\varphi)([gp,p]), \textnormal{ and } \Phi: P \ra Q \textnormal{ is given by } p' \mapsto \mathcal{F}(\Phi_\varphi)([p',p]).\]
Similarly, every morphism of groupoids between $P \otimes_{\group_x} P$ and $Q \otimes_{\subgroup_x} Q$ gives rise to a morphism of principal bundles, so this shows that the functor $\mathcal{F}$ is fully faithful. 

To show that the category of transitive groupoids is equivalent to the category of principal bundles via the functor $\mathcal{F}$, it remains to prove that $\mathcal{F}$ is essentially surjective. So, let $\group \rra M$ be a transitive groupoid, and let $x \in M$. We will show that it is isomorphic to the Gauge groupoid associated to the principal $\group_x$-bundle $P_x \xra{\target} M$, where $\group_x$ acts on $P_x \coloneqq \source^{-1}(x)$ by right multiplication (see Corollary \ref{coro: orbits are immersed submanifolds}). Consider the following map
\[P_x \otimes_{\group_x} P_x \ra \group \textnormal{ given by } [p,q] \mapsto pq^{-1}.\]
Then this map is smooth, e.g. because the pre-composition of it with the (surjective) submersion $\pi: P_x \times P_x \ra (P_x \times P_x)/\group_x$ is smooth. Now let $\tau: M \ra P_x$ be any set-theoretical section of $\target|_{P_x}$. Then,
\[\group \ra P_x \otimes_{\group_x} P_x \textnormal{ given by } g \mapsto [g\tau(\source(g)),\tau(\source(g))],\]
is smooth. To see why, observe first that for a different section $\tau': M \ra P_x$ we have 
\begin{align*}
    [g\tau'(\source(g)),\tau'(\source(g))] &= [g\tau (\source(g))\inv\left(\tau(\source(g))\right)\tau'(\source(g)),\tau (\source(g))\inv\left(\tau(\source(g))\right)\tau'(\source(g))] \\
    &= [g\tau (\source(g)),\tau (\source(g))],
\end{align*}
where we used that $\inv\left(\tau(\source(g))\right)\tau'(\source(g)) \in \group_x$. Hence, for all $g \in \group$, and a smooth (!) local section $\sigma: U \ra P_x$ of $\target|_{P_x}$ through $g^{-1}$, we can replace $\tau$ with $\sigma$ once we restrict the map to $\source^{-1}(U) \ni g$, which shows that the map is smooth. These maps can easily be extended to groupoid morphisms, and are readily verified to be inverse to each other. 
\end{rema}
%We also include the following remark about Gauge groupoids here. It will come in useful later when discussing an important example of a blow-up of Lie groupoids. 
%\begin{rema}\label{rema: Gauge groupoid map to pair groupoid}
%Let $P$ be a principal $G$-bundle. Observe that, for $i=1,2$, the map $P \times P \xra{\pr_i} P$ is $G$-equivariant, so it descends to a smooth map $(P \times P)/G \ra P/G$. In particular, we have a canonical smooth map
%\[F_P: (P \times P)/G \ra P/G \times P/G \textnormal{ given by } [p,q] \mapsto ([p],[q])\]
%which is smooth, because post-composing with either projection $P/G \times P/G \ra P/G$ gives one of the former smooth projection maps $(P \times P)/G \ra P/G$ 
%In fact, $F_P$ can be seen as a morphism of groupoids; that is, $F_P$ is even a functor between the Gauge groupoid $P \otimes_G P$ and the pair groupoid $P/G \times P/G \rra P/G$, over the base map $\id_{P/G}$. Indeed,
%\[F_P([p,q] \cdot [q,r]) = ([p],[r]) = ([p],[q]) \cdot ([q],[r]) = F_P([p,q]) \cdot F_P([q,r]), \textnormal{ and } F_P([p,p])=([p],[p]),\]
%where $p,q,r \in P$.
%\end{rema}

If $M \ni x_0$ is a connected manifold, we have an associated universal cover $p: \widetilde M \ra M$; we can take $\widetilde M$ to be the path space consisting of homotopy classes of paths starting at $x_0$ and set $p([\gamma]) = \gamma(1)$. If in Example \ref{exam: Gauge groupoids} we take $G=\pi(M;x_0)$, and $P= \widetilde M$, we get the following class of examples.
\begin{exam}[Fundamental groupoids]\cite{ruimarius}\label{exam: fundamental groupoids}
To any connected manifold $M$ we can associate the \textit{fundamental groupoid}
\[\Pi(M) \rra M,\]
where $\Pi(M)$ is the set of homotopy classes of continuous paths, and $\source([\gamma]) = \gamma(0)$, $\target([\gamma]) = \gamma(1)$, and $\mult([\gamma],[\gamma']) = [\gamma \star \gamma']$, where $\star$ denotes concatenation of paths: $\gamma \star \gamma'|_{[0,\frac{1}{2}]}(t) = \gamma(2t)$ and $\gamma \star \gamma'|_{[\frac{1}{2},1]}(t) = \gamma'(2t-1)$) (see also Remark \ref{rema: Mon to fundamental groupoids}).
\end{exam}

If we have a vector bundle $E \ra M$, we can again infer Example \ref{exam: Gauge groupoids} to obtain another familiar class of examples by taking $G = \textnormal{GL}(k,\mathbb{R})$ and $P=\textnormal{Fr } E$ (the frame bundle of $E$).
\begin{exam}[General linear groupoids]\cite{ruimarius}\label{exam: GL groupoids}
Let $E \ra M$ be a vector bundle. Then there is an associated \textit{general linear groupoid}
\[\textnormal{GL}(E) \rra M\]
(we denoted GL$(E)$ for the set of linear isomorphisms between any two fibres), with $\source(E_x \ra E_y) = x$, $\target(E_x \ra E_y) = y$ and $\mult(E_x \ra E_y, E_y \ra E_z) = E_x \ra E_z$. 
%Note that a representation $E \rightarrow M$ of $\group$ is the same thing as a morphism of Lie groupoids $\group \ra GL(E)$.
\end{exam}

\begin{exam}[Action groupoids]\cite{ruimarius}\label{exam: action groupoids}
Suppose we have a Lie group $G$ acting on a manifold $M$. Then there is an associated \textit{action groupoid}
\[G \ltimes M \coloneqq G \times M \rra M,\]
with $\source(g,x) = x$, $\target(g,x) = g \cdot x$, and $\mult((h,g \cdot x),(g,x)) = (hg,x)$. Observe that in this case $\orbit_x = \target(\source^{-1}(x)) = \{y \in M \mid gx = y$ for some $g \in G\}$ and $\group_x = \{g \in G \mid gx = x\}$, so the usual notions of orbit and isotropy (or stabiliser) group of the group action coincide with the ones defined for the groupoid.  
\end{exam}

%In the last part of this section we go into examples of Lie groupoids that are constructed out of given Lie groupoids. 

\begin{exam}\cite{ruimarius}\label{exam: tangent groupoid}
To any Lie groupoid $\group \rra M$ we can associate the \textit{tangent prolongation groupoid}
\[T\group \rra TM\]
in the obvious way: all structure maps are the differentials of the respective structure maps of $\group \rra M$. Since $d\source$ and $d\target$ are submersions, it is a Lie groupoid.  
\end{exam}

\begin{exam}\cite{meinrenken}\label{exam: pull-back of a groupoid}
Let $\groupoid$ be a set-theoretical groupoid (see Terminology \ref{term: groupoids}) and let $f: M \ra \base$ be a map. Then we can pull the (set-theoretical) groupoid back along $f$ to a (set-theoretical) groupoid
\[f^!\group \rra M,\]
where 
\[f^!\group \coloneqq M \tensor[_{f}]{\times}{_{\source}} \group \tensor[_{\target}]{\times}{_{f}} M = \{(x_1,g,x_2) \in M \times \group \times M \mid \source(g) = f(x_1) \textnormal{ and } \target(g)=f(x_2)\},\]
with $\source(x_1,g,x_2)=x_1$, $\target(x_1,g,x_2)=x_2$, and $\mult((x_1,g,x_2),(y_1,h,x_1)) = (x_2,\mult_\group(g,h),y_1)$ (note: this is well-defined, since $\source(g)=f(x_1)=\target(h)$ and $f(x_2)=\target(g)= \target(gh)$ and $f(y_1)=\source(h)=\source(gh)$). Pay attention to the fact that in case $\groupoid$ is a Lie groupoid and $f$ is a smooth map, the resulting set-theoretical groupoid $f^!\group$ is not necessarily a Lie groupoid (e.g. $M \tensor[_{f}]{\times}{_{\source}} \group \tensor[_{\target}]{\times}{_{f}} M$ does not have to be a smooth manifold). If $f$ is, in addition, a submersion, then $f^!\group$ is a Lie groupoid. For a more general result, we will introduce the notion of \textit{clean intersection} that we will use later as well. See Proposition \ref{prop: existence of pullback Lie groupoid} below for the result.
\end{exam}
So-called clean intersections of maps is a weaker condition than transversality of maps. Often results that hold for transverse maps are easily generalised to clean intersecting maps; the latter notion is often even easier to use in proofs (although there is more to check), and this is the condition we will mainly work with. %In the following definition, replace submanifold with embedded submanifold everywhere, or with immersed submanifold everywhere.
\begin{defn}\cite{meinrenken}\label{defn: clean intersetion}
Let $f: M \ra N$ and $f_i: M_i \ra N$ ($i=1,2$) be smooth maps and let $S,S' \subset M$ be embedded submanifolds.
\begin{enumerate}[label=\Roman*]
    \item The map $f$ is said to have \textit{clean intersection with $S$} if $f^{-1}(S) \subset M$ is an embedded submanifold such that for all $x \in f^{-1}(S)$ we have $T_xf^{-1}(S) = df(x)^{-1}(T_{f(x)}S)$.
    \item If $S' \hookrightarrow M$ has clean intersection with $S$, then $S'$ and $S$ are said to have \textit{clean intersection}.
    \item The maps $f_1$ and $f_2$ have \textit{clean intersection} if $f_1 \times f_2$ has clean intersection with the diagonal $\Delta_N \subset N \times N$. 
\end{enumerate} 
In I, the \textit{excess} $e$ is defined as $e \coloneqq (\dim N - \dim S) - (\dim M - \dim f^{-1}(S))$ (note: this gives rise to a notion of excess in II and III).
\end{defn}
\begin{rema}\cite{meinrenken}\label{rema: clean intersection}
Observe that if we have smooth maps $f_1: M_1 \ra N$ and $f_2: M_2 \ra N$ which are transverse, then they intersect cleanly. Indeed,
\[M_1 \tensor[_{f_1}]{\times}{_{f_2}} M_2 = \{(x_1,x_2) \in M_1 \times M_2 \mid f_1(x_1)=f(x_2)\} = (f \times f)^{-1}(\Delta_N)\]
is a submanifold of $M_1 \times M_2$, and 
\begin{align*}
    T_{(x_1,x_2)}(M_1 \tensor[_{f_1}]{\times}{_{f_2}} M_2) &= \{(\xi_1,\xi_2) \in T_{x_1}M_1 \oplus T_{x_2}M_2 \mid df_1(x_1)\xi_1=df_2(x_2)\xi_2\} \\
    &= d(f_1 \times f_2)(x_1,x_2)^{-1}(T_{f(x_1),f(x_2)}\Delta),
\end{align*}
which shows that, indeed, $f_1$ has clean intersection with $f_2$. Since $\dim M_1 \tensor[_{f_1}]{\times}{_{f_2}} M_2 = \dim M_1 + \dim M_2 - \dim N$, we have $e=0$ in this case.
\end{rema}
\begin{prop}\cite{meinrenken}\label{prop: existence of pullback Lie groupoid}
Let $\groupoid$ be a Lie groupoid, and let $f: M \ra \base$ be a smooth map. If $f \times f$ has a clean intersection with $(\target,\source)$, then $f^!\group$ is a Lie groupoid.
\end{prop}
\begin{proof}
%Observe that we can view $f \times f$ as a Lie groupoid morphism from the pair groupoid of $M$ to the pair groupoid of $N$. 
By assumption,
\[f \times (\source_\group,\target_\group) \times f: M \times \group \times M \ra \base \times \base \times \base \times \base\]
has clean intersection with the diagonal $\Delta \coloneqq \Delta_{X \times X} = \{(x,x,y,y) \in X \times X \times X \times X\}$, i.e.
\[\left(f \times (\source_\group,\target_\group) \times f\right)^{-1}(\Delta) = \{(x,g,y) \in M \times \group \times M \mid \source_\group(g)=f(x) \textnormal{ and } \target_\group(y)=f(y)\} = f^!\group\]
is an embedded submanifold of $M \times \group \times M$ such that for all $(x,g,y) \in f^!\group$ we have 
\begin{align*}
    T_{(x,g,y)}f^!\group &= d(f \times (\source_\group,\target_\group) \times f)(x,g,y)^{-1}(T_{(f(x),\source_\group(g),\target_\group(g),f(y))}\Delta) \\
    &= d(f \times (\source_\group,\target_\group) \times f)(x,g,y)^{-1}(\Delta_{T_{\source_\group(g)}\base} \oplus \Delta_{T_{\target_\group(g)}\base}).
\end{align*}
It remains to show that the structure maps $(\source,\target,\mult,\inv,\identity)$ of $f^!\group$ are smooth and that $\source$ is a submersion. First of all, the source and target maps of $f^!\group$ are given by restriction of the projection maps $\pr_1,\pr_3: M \times \group \times M \ra M$, so they are smooth. To see that $\source$ is a submersion, notice that, for all $(x,g,y) \in f^!\group$, we have
\[d\source(x,g,y)T_{(x,g,y)}f^!\group = T_xM,\]
where we used the above expression for $T_{(x,g,y)}f^!\group$ and that $d\source_\group(g)T_g\group=T_{\source_\group(g)}\base$. The multiplication map is smooth, because it is given by post-composing the map
\begin{align*}
    f^!\group^{(2)} = (M \tensor[_{f}]{\times}{_{\source_\group}} \group \tensor[_{\target_\group}]{\times}{_f} M) \tensor[_{\source}]{\times}{_{\target}} (M \tensor[_{f}]{\times}{_{\source_\group}} \group \tensor[_{\target_\group}]{\times}{_f} M) \xra{(\pr_3,(\pr_2,\pr_5),\pr_4)} M \times \group^{(2)} \times M
\end{align*}
with $\id_M \times \mult_\group \times \id_M$. Similarly, the inverse map $\inv=\id_M \times \inv_\group \times \id_M$ is smooth. The unit map is given by restricting the co-domain of the map $(\id_M,\identity_\group,\id_M): M \ra M \times \group \times M$ to $f^!\group$. Since $f^!\group$ is an embedded submanifold, this map is also smooth. This finishes the proof.
\end{proof}
%\subsection{Representations of Lie groupoids}\label{sec: Representations of groupoids}
As in the case of Lie groups, there is an appropriate notion of an action of a groupoid, and a representation of a groupoid.
\begin{defn}\cite{ruimarius}\label{defn: action of a groupoid}
Let $\groupoid$ be a groupoid and let $\mu: E \ra \base$ be a smooth map (note: we allow $E$ to be non-Hausdorff). Then $E$ is called a \textit{(left) $\group$-space with moment map $\mu$} if there is a multiplication map 
\[\nu: \group \tensor[_\source]{\times}{_\mu} E \xra{\cdot} E,\]
called the \textit{(groupoid) action} of $\group$ on $E$, such that 
\begin{enumerate}[label=\Roman*]
    \item $\mu(g \cdot e) = \target(g)$ (where $(g,e) \in \group \tensor[_\source]{\times}{_\mu} E$),
    \item $g \cdot (h \cdot e) = \mult(g,h) \cdot e$ (where $(g,h) \in \group^{(2)}$ and $(h,e) \in \group \tensor[_\source]{\times}{_\mu} E$) \textnormal{, and}
    \item $\id_{\mu(e)} \cdot e = e$.
\end{enumerate} 
If, in addition, $\mu: E \ra \base$ is a vector bundle, and for all $x,y \in \base$ and $g \in \group(x,y)$, $E_x \xra{g \cdot} E_y$ is a linear isomorphism, then it is called a \textit{representation} of $\group$.
\end{defn}
\begin{rema}\cite{ruimarius}\label{rema: representation equivalent to morphism}
Observe that a representation is equivalent to a morphism $\group \ra \textnormal{GL}(E)$.
\end{rema}
Notice that a $\groupoid$ is a left and right $\group$-space over itself (note: for right $\group$-spaces we swap the roles of $\source$ and $\target$ in Definition \ref{defn: action of a groupoid}). Indeed, it is a left $\group$-space with $\mu \coloneqq \target$ and $\nu \coloneqq \mult$, and it is a right $\group$-space with $\mu \coloneqq \source$ and $\nu \coloneqq \mult$.

We can even generalise action groupoids now to mean the following. 
\begin{exam}[action groupoids: generalisation]\cite{meinrenken}\label{exam: action groupoids generalisation}
Let $\groupoid$ be a groupoid, and let $\mu: E \ra \base$ be a $\group$-space. We can associate the \textit{action groupoid} 
\[\group \tensor[_{\source}]{\times}{_{\mu}} E \rra E,\]
with $\source(g,e) = e$, $\target(g,e)=g \cdot e$, and $\mult((h,g \cdot e),(g,e)) = (hg,e)$.
\end{exam}
Already we are seeing that many group-theoretic definitions and constructions generalise to groupoids in a natural way. In Section \ref{sec: Morita equivalence} we will continue this line of thought.

\subsection{Foliations and the monodromy and holonomy groupoids}\label{sec: Foliations and the monodromy and holonomy groupoids}
In this section we briefly introduce foliations, and quickly go through the construction of the monodromy and holonomy groupoids. These groupoids form an interesting class of examples of groupoids whose space of arrows is in general not Hausdorff. Moreover, foliations are a specific type of Lie algebroid, so they provide a rich source of examples of Lie algebroids. We will explain what a foliation is, discuss some simple examples, and give a detailed overview of the concept of holonomy.  
This section is based on \cite{moerdijk_mrcun_2003}. 

\begin{defn}\cite{moerdijk_mrcun_2003}\label{defn: foliation}
A partition $\mathcal{F}$ of $M$ into immersed submanifolds of codimension $q$ is called a \textit{foliation of $M$ of codimension $q$} if there is an open cover $(U_i)_{i \in I}$ such that 
\[\mathcal{F}|_{U_i} \coloneqq \{\textnormal{connected components of } L \cap U_i \mid L \in \mathcal{F}\}\] 
equals the partition of $U_i$ into the fibers of a submersion $U_i \ra \mathbb{R}^q$.In this case, $(M,\mathcal{F})$ is called a foliated manifold, and we will call the elements of $\mathcal{F}$ the \textit{leaves}.
%If a manifold $M$ satisfies one of the conditions in Proposition \ref{prop: equivalent definitions of foliations}, we call $M$ a \textit{foliated manifold}. We will call a partition $\mathcal{F}$ of $M$ into leaves a \textit{foliation of codimension $q$} and write $T\mathcal{F}$ for its corresponding involutive distribution $E$. Moreover, if $(U,\varphi)$ is a foliation chart, then we call the connected components of $\varphi^{-1}(\mathbb{R}^{n-q} \times \{y\})$ \textit{plaques}.
\end{defn}
In fact, there are many equivalent ways to define foliations. Most importantly, a foliation $\mathcal{F}$ on $M$ induces an involutive distribution $T\mathcal{F}$ of $M$, consisting of vectors which are tangent to the leaves, and, conversely, every involutive distribution induces a foliation on $M$.
\begin{term}\label{term: distribution}
We adopt the convention that, for a manifold $M$, a distribution of $M$ is a vector subbundle of $TM$ over $M$. 
\end{term}
\begin{prop}\cite{moerdijk_mrcun_2003}\label{prop: equivalent definitions of foliations}
Let $M$ be a manifold together with a partition $\mathcal{F}$ of immersed submanifolds of codimension $q$ (with $n \coloneqq \dim M$). The following are equivalent:
\begin{enumerate}[label=\Roman*]
    \item $\mathcal{F}$ is a foliation on $M$,
    \item there is a \textit{foliation atlas} $(\varphi_i: U_i \ra \mathbb{R}^{n-q} \times \mathbb{R}^q)_{i \in I}$, i.e. an atlas for which the transition maps $\varphi_{ij} \coloneqq \varphi_i \circ \varphi_j^{-1}$ are of the form $(x,y) \mapsto (g_{ij}(x,y),h_{ij}(y))$, such that $\mathcal{F}|_{U_i}$ equals the partition of $U_i$ into the so-called \textit{plaques} $\varphi^{-1}(\mathbb{R}^{n-q} \times \{y\})$,
    \item there is an open cover $(U_i)_{i \in I}$ of $M$ together with submersions $s_i:U_i \ra \mathbb{R}^q$ and \textit{Haefliger cocycles} $\gamma_{ij}$ relating them, i.e. diffeomorphisms $\gamma_{ij}: s_j(U_i \cap U_j) \ra s_i(U_i \cap U_j)$ which satisfy $s_i|_{U_i \cap U_j} = \gamma_{ij} \circ s_j|_{U_i \cap U_j}$, such that $\mathcal{F}|_{U_i}$ equals the partition of $U_i$ into the fibers of $s_i$,
    \item there is an \textit{involutive distribution} $E$ of $M$ of rank $n-q$, i.e. $\Gamma(E)$ is closed under taking the usual Lie bracket of vector fields, such that $E$ consists of the vectors which are tangent to the leaves,
    \item there is a \textit{locally trivial differential graded ideal} $J = \textstyle\bigoplus_{k=1}^n J^k$ of rank $q$ in $\Omega(M)$, i.e. every $x \in M$ has an open neighbourhood $U$ for which the ideal $J|_U \subset \Omega(M)|_U$ is generated by $q$ linearly independent $1$-forms (locally trivial) and $dJ \subset J$ (differential), such that $J$ consists of the $k$-forms which vanish when applied to a $k$-tuple of vectors tangent to the leaves.
\end{enumerate}
\end{prop}
\begin{proof} The statements $\RoN{1}$  $\leftrightarrow$ $\RoN{2}$ and $\RoN{1} \leftrightarrow \RoN{3}$ are readily verified. We show $\RoN{3} \ra \RoN{4}$, $\RoN{4} \leftrightarrow \RoN{5}$. The statement $\RoN{4} \ra \RoN{3}$ is known as the \textit{Frobenius theorem}; for a proof, see e.g. \cite{10.2307/2159499}%(here completely integrable is easily seen to imply \RoN{3})
.
%$\RoN{1} \ra \RoN{3}$ Take $s_i \coloneqq \pr_2 \circ \varphi_i$ and $\gamma_{ij} \coloneqq h_{ij}$. \\

\RoN{3} $\ra$ \RoN{4}:  Set $\alpha_i \coloneqq ds_i$ and define
\[E_i \coloneqq \ker \alpha_i.\]
Observe that the diffeomorphisms $\gamma_{ij}$ satisfy the cocycle conditions
\[\gamma_{ik} = \gamma_{ij} \circ \gamma_{jk},\]
so that we obtain a vector bundle $E \ra M$ if we define it locally by $E|_{U_i} \coloneqq E_i$. Then, $E$ is a distribution of $M$. Since for all $\sigma,\tau \in \mathfrak{X}(U_i)$ we have
\[d\alpha(\sigma,\tau) = \alpha([\sigma,\tau]) - \sigma(\alpha(\tau)) + \tau(\alpha(\sigma)),\]
and $\alpha$ is closed, we have $\alpha([\sigma,\tau]) = \sigma(\alpha(\tau)) - \tau(\alpha(\sigma))$. In particular, if $\sigma,\tau \in \Gamma(E|_{U_i})$ it follows that $\alpha([\sigma,\tau]) = 0$, hence $[\sigma,\tau] \in \Gamma(E|_{U_i})$, so $E$ is involutive.

\RoN{4} $\ra$ \RoN{5} Define for all $1 \le k \le n$ the ideal
\[J^k = \{\omega \in \Omega^k(M) \mid \omega(\sigma^1,\dots,\sigma^k) = 0 \textnormal{ for all } \sigma^1,\dots,\sigma^k \in \Gamma(E)\}\]
of $\Omega^k(M)$ and define $J \coloneqq \textstyle\oplus_{k=1}^n J^k$. To see that $J$ is locally trivial of rank $q$, let $x \in M$ and pick an open subset $U$ such that there is a local frame $\sigma^1,\dots,\sigma^n \in \mathfrak{X}(U)$ for which $\sigma^1,\dots,\sigma^{n-q}$ is a frame for $E|_U$. Clearly, $\omega^{n-q+1},\dots,\omega^n$ of the dual frame $\omega^1,\dots,\omega^n$ in $T^*M|_U$ generates $J$. That $J$ is differential, follows (using that $E$ is involutive) from the formula
\begin{align*}
    d\omega(\sigma^1,\dots,\sigma^{k+1}) = &\sum_{i < j} (-1)^{i+j} \omega([\sigma^i,\sigma^j],\dots,\widehat{\sigma^i},\dots,\widehat{\sigma^j},\dots,\sigma^{k+1}) \\
    - &\sum_i (-1)^i \mathcal{L}_{\sigma^i}(\omega(\sigma^1,\dots,\widehat{\sigma^i},\dots,\sigma^{k+1})).
\end{align*}

\RoN{5} $\ra$ \RoN{4} Since $J^1$ is locally generated by a $1$-form, we obtain a distribution $E$ of $M$ by locally defining it to be the kernel of this $1$-form (see \RoN{3} $\ra$ \RoN{4}). Again by the formula of the differential, we see that $E$ is involutive. 
\end{proof}

We now discuss a few simple examples.

\begin{exam}[Foliations by submersions]\cite{moerdijk_mrcun_2003}\label{exam: Stable foliation}
The fibers of a submersion $f: M \ra N$ give rise to a foliation on $M$. This follows by writing $f$ locally as a projection map. In particular, the total space of a groupoid $\group \rra \base$ admits a foliation by its source fibers, and also by its target fibers.
\end{exam}
\begin{exam}\cite{moerdijk_mrcun_2003}\label{exam: quotient foliation}
Let a discrete group $G$ act freely and properly on a manifold $M$, and suppose that $M$ is equipped with a foliation $\mathcal{F}$. If each $g \in G$, viewed as a diffeomorphism $M \ra M$, leaves the leaves of $\mathcal{F}$ invariant (that is, $g$ maps leaves to leaves), then $M/G$ inherits a foliation $\mathcal{F}/G$ from the foliation $\mathcal{F}$ on $M$. Indeed, this follows because the projection map $q: M \twoheadrightarrow M/G$ is locally injective, and therefore we can define a foliation atlas $(q(U_i),\varphi_i \circ q^{-1})$ on $M/G$, where $(U_i,\varphi_i)$ is a foliation atlas on $M$, such that the maps $q|_{U_i}$ are injective.
\end{exam}
The following four simple examples will turn out to be examples for which the monodromy and the holonomy groupoid can be non-Hausdorff, and we will see that (in order of presentation) both can be Hausdorff, one can be Hausdorff while the other one is not, or both can be non-Hausdorff.
\begin{exam}[Kronecker foliation]\cite{moerdijk_mrcun_2003}\label{exam: Kronecker foliation}
The \textit{Kronecker foliation} of the torus $T^2$ is the foliation obtained by fixing an irrational number $r \in \mathbb{R} \setminus \mathbb{Q}$, say $r=\sqrt{2}$, and letting the foliation on $T^2 \cong \mathbb{R}^2/\mathbb{Z}^2$ be the induced foliation from the foliation on $\mathbb{R}^2$ given by lines with slope $r$ (this is a foliation on $\mathbb{R}^2$ which is invariant under the $\mathbb{Z}^2$-action on $\mathbb{R}^2$; see Example \ref{exam: quotient foliation}). Notice that all the leaves are diffeomorphic to $\mathbb{R}$, and are dense in $T^2$.
\end{exam}
\begin{exam}\cite{hoyohausdorff}\label{exam: R^3 setminus 0 foliation}
The foliation induced by the submersion $\pr_3: \mathbb{R}^3 \setminus \{0\} \ra \mathbb{R}$. The leaves of the foliation are the ``horizontal planes'', and the leaves are diffeomorphic to $\mathbb{R}^2$, except for the leaf given by the fiber of $0 \in \mathbb{R}$, which is diffeomorphic to $\mathbb{R}^2 \setminus \{0\}$.
\end{exam}
\begin{exam}\cite{hoyohausdorff}\label{exam: e^-1/t}
Consider the smooth map 
\[f: \mathbb{R} \ra \mathbb{R} \textnormal{ given by }
t \mapsto 
\begin{cases}
e^{-\tfrac{1}{t}} & \text{if}\ t>0 \\
0 & \text{if}\ t \le 0.
\end{cases}\]
This map induces a foliation on $\mathbb{R} \times S^1$ via the vector field
\[\sigma(t,\theta) \coloneqq \frac{\partial}{\partial \theta} + f(t)\frac{\partial}{\partial t};\]
that is, $\sigma$ spans a distribution of $\mathbb{R} \times S^1$, and therefore we obtain a foliation on $\mathbb{R} \times S^1$ (notice that the condition for the distribution to be involutive is trivial in this case). It is not hard to see that elements of $\mathbb{R}_{\le 0} \times S^1$ lie on leaves which are diffeomorphic to $S^1$ (the leaves are given by $\{t\} \times S^1$ where $t<0$), and elements of $\mathbb{R}_{>0} \times S^1$ lie on leaves which are diffeomorphic to $\mathbb{R}$.
\end{exam}
\begin{exam}[Reeb foliation]\cite{moerdijk_mrcun_2003}\label{exam: Reeb foliation}
We will construct a foliation on $S^3$ called the \textit{Reeb foliation}. Before we do this, let $n \ge 1$. Observe that we have a smooth submersion
\[\mathbb{R} \times S^{n-1} \times \mathbb{R} \ra \mathbb{R} \textnormal{ given by } (r,\theta,t) \mapsto (r^2-1)e^t.\]
By considering cylindrical coordinates on the $n$-dimensional disc $D^n$, we may view $(0,1) \times S^{n-1} \subset D^n$. Then, the restriction of the above map to $(0,1) \times S^{n-1} \times \mathbb{R}$ extends smoothly to $\{0\} \times \mathbb{R} \subset D^n \times \mathbb{R}$, so that we obtain a smooth map
\begin{align*}f: D^n \times \mathbb{R} \ra \mathbb{R} \textnormal{ given by } z \mapsto
\begin{cases}
(r^2-1)e^t & \text{if}\ z = (r,\theta,t) \\
-e^t & \text{if}\ z=(0,t)
\end{cases}
\end{align*}
which again is a submersion. We now have an induced stable foliation on the cylinder $D^n \times \mathbb{R}$. By inspecting the map $f$, we have for $s \in \mathbb{R}$ that
\[f^{-1}(s)=\{(r,\theta,t) \in D^n \times \mathbb{R} \mid 1+se^{-t}=r^2\} \cup \{(0,t) \in D^n \times \mathbb{R} \mid 1+se^{-t}=0\},\]
hence, for $d \in \mathbb{Z}$, we see that
\[f^{-1}(se^d) = \{(z,t+d) \in D^n \times \mathbb{R} \mid (z,t) \in f^{-1}(s)\}.\]
So, with respect to the $\mathbb{Z}$-action $(z,t) \cdot d = (z,t+d)$ on $D^n \times \mathbb{R}$, the quotient $(D^n \times \mathbb{R})/\mathbb{Z} \cong D^n \times S^1$ inherits a foliation from $D^n \times \mathbb{R}$. Now, let $n=1$. Then the described foliation is a foliation on the solid torus $D^2 \times S^1$ (notice that $D^2$. By viewing $S^3 \subset \mathbb{C}^2$, we can see that $S^3$ can be realised as two copies of solid tori ``glued'' along their boundary, the torus (surface) $S^1 \times S^1$. To see this, note that 
\[S^3 = \{(z,w) \in \mathbb{C}^2 \mid z^2 +w^2 = 1\} = \{(z,w) \in S^3 \mid |z| \le \tfrac{1}{\sqrt{2}}\} \cup \{(z,w) \in S^3 \mid |w| \le \tfrac{1}{\sqrt{2}}\}.\]
The foliation described before on the solid torus $D^2 \times S^1$ has one leaf equal to $S^1 \times S^1$. Consequently, the foliation induces a foliation on $S^3$, called the Reeb foliation, whose leaves are given by the leaves of the foliation described before on both copies of the solid tori inside $S^3$. So, this foliation partitions $S^3$ into leaves diffeomorphic to $\mathbb{R}^2$, except for one leaf, which is diffeomorphic to the torus $T^2$. 
\end{exam}

We end the section with a useful result regarding the topology of a foliated manifold.
\begin{prop}\label{prop: the collection of leaves intersecting an open subset}
Let $(M,\mathcal{F})$ be a foliated manifold. If $T \subset M$ is a transversal (i.e. a submanifold which is transverse to the leaves), then the collection
\[T_{\textnormal{leaf}} \coloneqq \bigcup_{L \in \mathcal{F} \textnormal{ with } L \cap T \neq \emptyset} L \subset M\]
is an open subset.
\end{prop}
\begin{proof}
%Obviously, $U \subset U_{\textnormal{leaf}}$, so 
Let $y \in L$ with $L \cap T \neq \emptyset$% and $x \not\in T$
, say, $x \in L \cap T$.
Since $L$ is connected, we can find a foliation path $\gamma: [0,1] \ra L$ from $x$ to $y$. We now cover $\gamma([0,1])$ with, say, $\ell+1$ foliation charts $(U_0,\varphi_0),\dots,(U_\ell,\varphi_\ell)$ such that $\gamma([\tfrac{j}{\ell+1},\tfrac{j+1}{\ell+1}]) \subset U_j$. By covering $T$ with smaller transversals (and taking a union of open sets at the end), and shrinking $U_0$, we may assume that $T$ equals the transversal $\varphi_0^{-1}(\{0\} \times \mathbb{R}^q)$ (assuming that $\varphi_0(x) = 0$).
Observe that all the $U_j$ are a union of plaques; the connected components of $\varphi_j^{-1}(\mathbb{R}^{n-q} \times \{z\})$. By shrinking $U_1$, we may assume that all of the plaques in $U_1$ intersect $U_0$. By an inductive argument, we can make sure that all of the plaques in $(U_j,\varphi_j)$ intersect with a plaque in $(U_{j-1},\varphi_{j-1})$ (we possibly need more charts). Again, by an inductive argument, all plaques in $(U_\ell,\varphi_\ell)$ are contained in a leaf that intersects $T$, i.e. $U_\ell \subset T_{\textnormal{leaf}}$. Since $y \in U_\ell$, this proves the statement.
\end{proof}
\subsubsection{Holonomy}\label{sec: Holonomy}
We introduce holonomy in the setting of II of Proposition \ref{prop: equivalent definitions of foliations}. Before we do so, let us start by introducing holonomy groups of a vector bundle by using \textit{parallel transport}. If the reader is familiar with connections on a vector bundle, this is a nice and easy introduction to (linear) holonomy. Later in this section we will see that both notions of holonomy are related. 

Consider a vector bundle $E \ra M$ of rank $k$ equipped with a connection $\nabla$, i.e. an $\mathbb{R}$-linear map $\Gamma(E) \ra \Gamma(T^*M \otimes E)$ satisfying the Leibniz-type rule 
\[\nabla(f \cdot s) = f \cdot \nabla s + df \otimes s \text{ for all } f \in C^{\infty}(M) \text{ and } s \in \Gamma(E).\]
If $\gamma: [0,1] \ra M$ is a path with $\gamma(0)=x$ and $\gamma(1)=y$, then to each $v \in E_x$ we can associate a unique section $\tau_\gamma^v$ of $\gamma^*E$ with the property that 
\[\gamma^*\nabla(\tau_\gamma^v)=0 \text{ subject to the initial condition } \tau_\gamma^v(0)=v\] 
since, locally, the above property becomes a first-order linear initial value problem. Indeed, we can write $\nabla$ locally as a matrix whose entries are $1$-forms $(\omega^j_i)$ (then, $\nabla(e_i)=\textstyle\sum_j\omega^j_i \otimes e_j$ with $(e_j)$ a local frame of $E$). By writing a local section $\tau$ of $\gamma^*E$ as a smooth map  $(\tau^1,\dots,\tau^k): [0,1] \ra \mathbb{R}^k$, the expression $\gamma^*\nabla(\tau)(t)=0$ can be written as
\[\dot\tau^i(t)+\sum_j\tau^j\omega^i_j(\dot\gamma(t))=0,\]
which, indeed, becomes a first-order linear initial value problem $\dot \tau(t)=A\tau(t)$, where $A$ is the matrix $-(\omega^j_i(\dot\gamma(t))$, subject to the initial condition $\tau(0)=v$.
In particular, we obtain a well-defined map
\[\tau_\gamma: E_x \ra E_y \text{ given by } v \mapsto \tau_\gamma^v(1).\]
Now, denote $\gamma^{-1}:[0,1] \ra M$ for the path $t \mapsto \gamma(1-t)$ and write $\star$ for concatenation of paths. It is clear from the definition that $\tau_x$, where $x: [0,1] \ra M$ also denotes the constant path with value $x$, is the identity element of $\GL(E_x)$. It is then also clear that if $\gamma$ is a loop (i.e. $x=\gamma(0)=\gamma(1)=y$), then $\tau_\gamma \in \GL(E_x)$ with inverse $\tau_{\gamma^{-1}}$. Moreover, if $\gamma'$ is another loop based at $x$, then the composition in $\GL(E_x)$ of $\tau_\gamma$ and $\tau_{\gamma'}$ equals $\tau_{\gamma \star \gamma'}$. Summarising this discussion,
\[\{\tau_\gamma \mid \gamma \text{ is a loop based at } x\}\] 
is a subgroup of $\GL(E_x)$. Lastly, notice that if a loop $\gamma$ is contractible, then so is $\gamma' \star \gamma \star (\gamma')^{-1}$ for any other loop $\gamma'$. In particular, the above subgroup of $\GL(E_x)$ contains a normal subgroup consisting of elements $\tau_\gamma$ with $\gamma$ contractible.
We denote the group
\[\{\tau_\gamma \mid \gamma \text{ is a loop based at } x\}/\{\tau_\gamma \mid \gamma \text{ is a contractible loop based at } x\}\]
by $\Hol(\nabla,x)$.
\begin{rema}\label{rema: surjective homomorphism fundamental group to Hol/Hol^0}
With notation as in the above definition, it is clear that there is a surjective group homomorphism $\pi_1(M,x) \ra \Hol(\nabla,x) \text{ given by } \gamma \mapsto \tau_\gamma$.
\end{rema}

We now restrict our attention to foliated manifolds. In this case, there is a distinguished class of paths, namely those paths having the property that they are contained in just one leaf.
\begin{defn}\cite{meinrenken}\label{defn: foliation path}
Let $M$ be a foliated manifold. A path (resp. loop) $\gamma: [0,1] \ra M$ contained in a single leaf $L$ of the foliation (i.e. $\gamma([0,1]) \subset L$) is called a \textit{foliation path} (resp. \textit{foliation loop}).
\end{defn}

Using these special paths, we can define a holonomy group that is similar in nature, but constructed in a different way.
So, let $M$ be a manifold of dimension $n$ equipped with a foliation $\mathcal{F}$ of codimension $q$. Suppose we have a foliation path $\gamma: [0,1] \ra M$ from $x=\gamma(0)$ to $y=\gamma(1)$ contained in a leaf $L$. Divide $\gamma$ into paths $\gamma_1,\dots,\gamma_r$ such that there is a sequence of foliation charts
\[(U_1,\varphi_1),\dots,(U_r,\varphi_r)\]
with the property that, for all $1 \le j \le r$, $\gamma_j([0,1]) \subset U_j$. Now, for all $1 \le j \le r$, let $x_j \coloneqq \gamma_j(0)$ and pick a transversal $T_j \subset U_j$, i.e. a submanifold transverse to the leaves of the foliation, with $x_j \in T_j$, and denote $T \coloneqq T_1$ and $S \coloneqq T_r$. For all $1 \le j \le r$, we can shrink $T_j$, if necessary, to open subsets $T_j' \subset T_j$, containing $x_j$, and intersecting a plaque of $(U_{j+1},\varphi_{j+1})$ only once; that way we obtain a sequence of maps 
\[T_1' \ra \dots \ra T_r',\]
where $T_j' \ra T_{j+1}'$ is defined by mapping a point to the unique point on the same plaque (note that it maps $x_j$ to $x_{j+1})$. These maps are locally given by a translation, so that again, by shrinking the $T_j'$ accordingly, we may assume that the maps in the sequence $T_1' \ra \dots \ra T_r'$ are all diffeomorphisms. The point is this: we will regard the resulting diffeomorphism $T_1' \ra T_r'$ as a germ of a diffeomorphism $T \ra S$ from $x$ to $y$, i.e. regard two diffeomorphisms $U \ra V$ and $U' \ra V'$ (with $U,U' \subset T$ and $V,V' \subset S$ open subsets) mapping $x$ to $y$ as equivalent if they agree on an open subset $W \subset U \cap U'$ around $x$. This way, the above diffeomorphism is a particular choice of representative of a germ of a diffeomorphism $T \ra S$, and is denoted by
\[\hol^{T,S}(\gamma) \in \germ_{x,y}(T,S).\] 
It does matter what transversals $T$ and $S$ we begin with. Note, however, that two transversals $T,T' \ni x$ (resp. two transversals $S,S' \ni y$) give rise to $\hol^{T,T'}(x)$ (resp. $\hol^{S,S'}(y)$), regarding $x$ (resp. $y$) as a constant path, and it is clear from the construction that $\hol^{T,S}(\gamma)$ and $\hol^{T',S'}(\gamma)$ are related by
\[\hol^{T,S}(\gamma) = \hol^{T,T'}(x) \circ \hol^{T',S'}(\gamma) \circ \hol^{S',S}(y).\]
Moreover, by inspecting the construction, we see that probably many choices of $\gamma$ give rise to the same element $\hol^{T,S}(\gamma)$, since we only use a sequence of points that lie on the path. In fact, it only depends on the homotopy class of $\gamma$. Indeed, this fact is readily verified by dividing a homotopy $F$ between $\gamma$ and another path into $d$ smaller homotopies $F_j$ from $\gamma_j$ to $\gamma_{j+1}$ (with $\gamma_1=\gamma$) in such a way that for each $1 \le j \le d$ we can use for $\gamma_j$ and $\gamma_{j+1}$ the same sequence of foliation charts. Restricting to foliation loops, we obtain a map $\hol^{T,T}: \pi_1(L,x) \ra \germ_x(T)$ that, by identifying $\germ_x(T)$ with $\germ_0(\mathbb{R}^q)$, gives rise to a map 
\[\hol: \pi_1(L,x) \ra \germ_0(\mathbb{R}^q).\]
From the construction of $\hol$ it is clear that it is even a group homomorphism. 
\begin{defn}\cite{moerdijk_mrcun_2003}\label{defn: holonomy group of a foliation}
Let $(M,\mathcal{F})$ be a foliated manifold. The holonomy group $\Hol(L,x)$ of $L \in \mathcal{F}$ is the image of $\pi_1(L,x)$ under the group homomorphism $\hol: \pi_1(L,x) \ra \germ_0(\mathbb{R}^q)$. It is determined up to inner automorphism (i.e. up to a conjugation) of $\germ_0(\mathbb{R}^q)$.
\end{defn}
\begin{rema}\cite{moerdijk_mrcun_2003}\label{rema: Hol is well-defined}
Since the holonomy group of $L$ is only determined up to inner automorphism, it is independent of the choice of base-point.
\end{rema}
Observe that the differential $d$ of $\mathbb{R}^q$ induces a homomorphism of groups $\germ_0(\mathbb{R}^q) \ra \GL(\mathbb{R}^q)$ by sending $f$ to $df(0)$. 
\begin{defn}\cite{moerdijk_mrcun_2003}\label{defn: linear Hol of a foliation}
The linear holonomy group $d\Hol(L,x)$ is the image of $\pi_1(L,x)$ under the group homomorphism $d\hol \coloneqq d \circ \hol: \pi_1(L,x) \ra \GL(\mathbb{R}^q)$.
\end{defn}
In fact, the linear holonomy group is the holonomy group of a specific connection that arises naturally on the leaf $L$. We can describe this connection as coming from a certain connection on a vector bundle intimately related to the involutive distribution $T\mathcal{F}$ associated to $\mathcal{F}$, namely the normal bundle of the foliation, which is given by 
\[\normal\mathcal{F} \coloneqq TM/T\mathcal{F} \ra M.\]
However, it is not a connection on this vector bundle in the usual sense. To describe this connection, we will use an alternative description of the normal bundle as a subbundle of $T^*M$ instead of $TM$:
\[T^{\textnormal{ann}}\mathcal{F} \coloneqq \{\alpha \in T^*M \mid \alpha|_{T\mathcal{F}}=0\}.\] 
\begin{defn}\cite{bott}\label{defn: Bott connection}
Let $M$ be a foliated manifold with foliation $\mathcal{F}$ and let $L$ be a leaf. The \textit{Bott-connection} $\nabla^L$ is the $\mathbb{R}$-bilinear map $\Gamma(T^{\textnormal{ann}}\mathcal{F}) \ra \Gamma(T^*\mathcal{F} \otimes T^{\textnormal{ann}}\mathcal{F})$, defined by 
\[\nabla^L(\alpha) \coloneqq \mathcal{L}_\bullet\alpha,\]
where $\mathcal{L}_\sigma$ denotes the Lie derivative along a vector field $\sigma$.
\end{defn}
\begin{rema}\label{rema: Bott connection}
As already mentioned, the Bott-connection is not really a connection on $E^{\textnormal{ann}}$: we only allow vector fields tangent to the foliation. All the other properties of a connection are, however, satisfied. To see that it is well-defined, let $\sigma \in T\mathcal{F}$ and $\alpha \in T^{\textnormal{ann}}\mathcal{F}$. Then $\mathcal{L}_\sigma\alpha=\iota_\sigma\alpha$ by Cartan's magic formula (since $\alpha(X)=0$), so indeed
\[\mathcal{L}_\sigma\alpha(\tau)=(\iota_\sigma\alpha)(\tau)=-\alpha([\sigma,\tau])=0,\]
where we used that $T\mathcal{F}$ is involutive in the last equality. Later, when discussing Lie algebroids, we will see that these types of connections are very natural to look at when working with Lie algebroids. Even though the Bott-connection is not a connection on $T^{\textnormal{ann}}\mathcal{F}$, we can restrict it to a leaf $L$; there, it is a connection.
\end{rema}

Now, the Bott-connection induces a holonomy group on each leaf $L$ of the foliated manifold $M$. In fact, if $\gamma$ is a foliation loop based at $x \in M$, we have $\tau_\gamma$, which we can view as an element of $\GL(\mathbb{R}^q)$, and $d\hol(\gamma)$; but one can show that, up to the correct identifications, these two elements are equal. 

\subsubsection{The monodromy and holonomy groupoids}\label{sec: monodromy and holonomy groupoids}
Here, we will define the monodromy and holonomy groupoid associated to a foliated manifold. Recall from the previous section that we constructed from a foliation path $\gamma$, and two transversals $T \ni \gamma(0) \eqqcolon x$ and $S \ni \gamma(1) \eqqcolon y$, a germ of a diffeomorphism $T \ra S$ from $x$ to $y$ called $\hol^{T,S}(\gamma)$. Moreover, if $\gamma$ is a foliation loop, we view it as the element $\hol(\gamma)$ in $\germ_0(\mathbb{R}^q)$, which is determined up to an inner automorphism of $\germ_0(\mathbb{R}^q)$. In particular, two foliation paths $\gamma$ and $\gamma'$, contained in a leaf $L$, and having the same endpoints $x=\gamma(0)=\gamma'(0)$ and $y=\gamma(1)=\gamma'(1)$, determine the element 
\[\hol(\gamma \star (\gamma')^{-1}) \in \germ_0(\mathbb{R}^q).\]
By requiring that $\hol(\gamma \star (\gamma')^{-1}) = 1$, we define an equivalence relation on the space of paths in $L$ with endpoints $x$ and $y$.
\begin{defn}\cite{moerdijk_mrcun_2003}\label{defn: holonomy class of a foliation path}
Let $\gamma$ be a foliation path. Then the equivalence class of $\gamma$, with respect to the equivalence relation 
\[\gamma' \sim \gamma'' \textnormal{ if and only if } \hol(\gamma' \star (\gamma'')^{-1}) = 1,\]
is called the \textit{holonomy class} of $\gamma$. 
\end{defn}
Now, we can view the holonomy class $[\gamma]$ of a foliation path $\gamma$we can view  as an arrow from $x$ to $y$ with inverse arrow $[\gamma^{-1}]$ (since $\gamma \star \gamma^{-1}$ is homotopic to the constant path $x$). Of course, homotopy (relative to endpoints) also defines an equivalence relation on the space of foliation paths of the leaf $L$ with respect to the endpoints $x$ and $y$, and the homotopy class $[\gamma]$ also can be viewed as an arrow from $x$ to $y$ with inverse arrow $[\gamma^{-1}]$.  
This leads us naturally to the notions of monodromy and holonomy groupoids, which are even, in a natural way, Lie groupoids.
\begin{exam}\cite{moerdijk_mrcun_2003}[Monodromy and holonomy groupoids]\label{exam: monodromy and holonomy groupoids}
Let $M$ be a foliated manifold with foliation $\mathcal{F}$ and let $L$ be a leaf. Then there is an associated \textit{monodromy groupoid} 
\[\Mon(M,\mathcal{F}) \rra M\]
(where we denoted $\Mon(M,\mathcal{F})$ for the set of homotopy classes, relative to endpoints, of foliation paths) with $\source([\gamma])=\gamma(0)$, $\target([\gamma])=\gamma(1)$ and $\mult([\gamma],[\gamma']) = [\gamma \star \gamma']$. 

Similarly, there is an associated \textit{holonomy groupoid}
\[\Hol(M,\mathcal{F}) \rra M\]
(where we denoted $\Hol(M,\mathcal{F})$ for the set of holonomy classes, relative to endpoints, of foliation paths) with $\source([\gamma])=\gamma(0)$, $\target([\gamma])=\gamma(1)$ and $\mult([\gamma],[\gamma']) = [\gamma \star \gamma']$.
\end{exam}
\begin{rema}\cite{moerdijk_mrcun_2003}\label{rema: groupoid morphism Mon -> Hol}
Before we prove that the above groupoids are Lie groupoids, we observe here that two homotopic foliation paths are holonomic to each other, so we have a canonical map 
\[F=F_{\mathcal{F}}: \Mon(M,\mathcal{F}) \ra \Hol(M,\mathcal{F})\]
which, by definition of the (set-theoretic) groupoid structures, obviously extends to a groupoid morphism over $\id_M$. We will use this groupoid morphism in the proof below; it will then also be obvious that, when equipped with the Lie groupoid structures that we describe below, it is a morphism of Lie groupoids.
\end{rema}
\begin{prop}\cite{moerdijk_mrcun_2003}\label{prop: Mon and Hol are Lie groupodis}
The monodromy and holonomy groupoid of a manifold $M$ of dimension $n$, equipped with a foliation $\mathcal{F}$ of codimension $q$, are Lie groupoids of dimension $2n-q$.
\end{prop}
\begin{proof}
Let $\gamma$ be a foliation path with endpoints $x \coloneqq \gamma(0)$ and $y \coloneqq \gamma(1)$ and pick foliation charts $(U_x,\varphi_x)$ and $(U_y,\varphi_y)$ such that $\varphi_x(U_x)=V_x \times W_x \subset \mathbb{R}^{n-q} \times \mathbb{R}^q$ with $V_x$ and $W_x$ connected and simply connected open subsets of $\mathbb{R}^{n-q}$ and $\mathbb{R}^q$ respectively (and similarly $\varphi_y(U_y)=V_y \times W_y$ with $V_y$ and $W_y$ connected and simply connected). Write $\varphi_x(x)=(v_x,w_x)$ and $\varphi_y(y)=(v_y,w_y)$ and shrink $V_x,W_x,V_y$ and $W_y$, if necessary, so that a representative of $\hol^{T,S}(\gamma)$, where $S$ and $T$ are the transversals
\[S \coloneqq \varphi_x^{-1}(\{v_x\} \times W_x) \textnormal{ and } T \coloneqq \varphi_y^{-1}(\{v_y\} \times W_y),\]
is actually a diffeomorphism $S \ra T$. From the construction of $\hol^{T,S}(\gamma)$ (see Section \ref{sec: Holonomy}), it is clear that this diffeomorphism can be realised as a smooth map $H(1,\cdot)$, where $H: [0,1] \times S \ra M$ is a smooth map with the following properties:
\[H(0,\cdot) = \iota_S \textnormal{, } H(1,\cdot)=\iota_T \circ \hol^{T,S}(\gamma) \textnormal{ and for all } s \in S, \textnormal{ } H(\cdot,s) \textnormal{ is a path in some leaf},\]
with $\iota_S: S \hookrightarrow M$ (resp. $\iota_T: T \hookrightarrow M$) the inclusion map.

Now, we will define an injective map
\[\varphi: V_x \times V_y \times W_x \ra \Mon(M,\mathcal{F})\] 
as follows: let $(v,v',w) \in V_x \times V_y \times W_x$. Then $s \coloneqq \varphi_x^{-1}(v,w)$ lies in the same plaque as $s_x \coloneqq \varphi_x^{-1}(v_x,w)$, so we take any path $\delta$ from $s$ to $s_x$ whose image lies entirely in the plaque. Similarly, we define such a path $\delta'$ from $s_y \coloneqq \hol^{T,S}(\gamma)(s_x)$, which equals $\varphi_y^{-1}(v_y,w')$ for some $w' \in W_y$, to $s'\coloneqq \varphi_y^{-1}(v',w')$. We now put
\[\varphi(v,v',w) \coloneqq [\delta' \star H(\cdot,s_x) \star \delta].\]
Notice that the homotopy class of $\delta' \star H(\cdot,s_x) \star \delta$ is independent of the choices we made for $\delta$ and $\delta'$, since we required $V_x,W_x,V_y$ and $W_y$ to be connected and simply connected. Moreover, $\delta' \star H(\cdot,s_x) \star \delta$ is a path from $s = \varphi_x^{-1}(v,w)$ to $s' = \varphi_y^{-1}(v',w')$, so indeed it is an injective map, and if we post-compose $\varphi$ with the groupoid morphism $F$ from Remark \ref{rema: groupoid morphism Mon -> Hol}, it is still injective by the same reasoning. We now equip $\Mon(M,\mathcal{F})$ (resp. $\Hol(M,\mathcal{F})$) with the topology that is generated by the basis given by the images of such maps $\varphi$ (resp. $F \circ \varphi$). To see that the maps $\varphi$ give rise to a smooth structure on $\Mon(M,\mathcal{F})$, notice that the inverse of $\varphi$ is given by 
\[[\gamma] \mapsto \left(\varphi^1_x(\gamma(0)),\varphi_y^1(\gamma(1)),\varphi_x^2(\gamma(0))\right)\]
(where we wrote $\varphi_x=(\varphi^1_x,\varphi^2_x)$ and $\varphi_y=(\varphi^1_y,\varphi^2_y)$), so that if $\psi$ is constructed in the same way as $\varphi$, but from the charts $(U_{x'},\psi_{x'})$ and $(U_{y'},\psi_{y'})$ (that are compatible with the charts $(U_x,\varphi_x)$ and $(U_y,\varphi_y)$), then  
\[\psi^{-1} \circ \varphi(v,v',w) = \left(\psi^1_{x'} \circ \varphi_x^{-1}(v,w), \psi_{y'}^1 \circ \varphi_y^{-1}(v',w'), \psi_{x'}^2 \circ \varphi_x^{-1}(v,w)\right),\]
which is a smooth map by assumption. Similarly, the maps $F \circ \varphi$ give rise to a smooth structure on $\Hol(M,\mathcal{F})$. We also see that the source (resp. target) map of $\Mon(M,\mathcal{F})$ becomes locally a projection map with respect to the charts $(\varphi(V_x \times V_y \times W_x),\varphi^{-1})$ and $(U_x,\varphi_x)$ (resp. $(U_y,\varphi_y)$). So, the source map (and the target map of $\Mon(M,\mathcal{F})$) is a submersion (and, similarly, the source map of $\Hol(M,\mathcal{F})$ is a submersion). Lastly, $\dim V_x + \dim V_y + \dim W_x = 2n-q$, which concludes the proof. 
\end{proof}
We end this section with a few remarks about monodromy and holonomy groupoids. We start with the following simple observation.
\begin{rema}\cite{moerdijk_mrcun_2003}\label{rema: orbits and isotropy groups}
Notice that the orbits of both the monodromy and holonomy groupoid are given by the leaves of the foliation. Moreover, if $L$ is a leaf and $x \in L$, then the isotropy group at $x$ of the monodromy groupoid is given by the fundamental group $\pi_1(L,x)$, and the isotropy group at $x$ of the holonomy groupoid is given by the holonomy group $\Hol(L,x)$. 
\end{rema}
We also obtain a different characterisation of the fundamental groupoid (see Example \ref{exam: fundamental groupoids}).
\begin{rema}\cite{moerdijk_mrcun_2003}\label{rema: Mon to fundamental groupoids}
Since every connected manifold $M$ has a trivial foliation given by $\mathcal{F} \coloneqq \{M\}$, we can give an alternative description of the fundamental groupoid $\Pi(M)$, namely as the monodromy groupoid $\Mon(M,\mathcal{F})$. 
\end{rema}
As already mentioned a few times, monodromy and holonomy groupoids are not necessarily Hausdorff. In fact, there seems to be no relation in general between the Hausdorffness of the monodromy groupoid and the Hausdorffness of the holonomy groupoid. We discuss the monodromy and holonomy groupoids of the four examples (Examples \ref{exam: Kronecker foliation} up to and including \ref{exam: Reeb foliation}) discussed before, and we will see that e.g. the monodromy groupoid can be Hausdorff while the holonomy groupoid is not, and vice versa. To keep it short, we will not go into detail about the structure maps of these groupoids, but focus on whether or not the total spaces are Hausdorff or not.
\begin{rema}\cite{moerdijk_mrcun_2003}\label{rema: monodromy and holonomy groupoids Kronecker foliation}
The monodromy and holonomy groupoids of the Kronecker foliation (see Example \ref{exam: Kronecker foliation}) are readily verified to be equal, since foliation paths are homotopic if and only if they are holonomic; in fact, foliation paths are homotopic/holonomic if and only if they have the same endpoints. To describe this groupoid, notice that the $\mathbb{Z}$-action on $\mathbb{R}^2$ given by $m \cdot (t,x) \coloneqq (t,m+x)$ exhibits $\mathbb{R} \times S^1$ as the quotient space $\mathbb{R}^2/\mathbb{Z}$, and the foliation by lines of slope $r=\sqrt{2}$ on $\mathbb{R}^2$ induces a foliation on $\mathbb{R} \times S^1$. It is clear that, again, the monodromy and holonomy groupoids are equal (following the same reasoning), and is therefore given by $\mathbb{R} \times \mathbb{R} \times S^1 \rra \mathbb{R} \times S^1$% (with source map and target map equal to $\pr_{2,3}$ and $\pr_{1,3}$, respectively)
. The topology on the total space $\mathbb{R} \times \mathbb{R} \times S^1$ is just the product topology. The $\mathbb{Z}$-action on $\mathbb{R} \times S^1$ given by $m \cdot (t,\theta) \coloneqq (m+t,\theta)$, which exhibits $T^2$ as the quotient space $(\mathbb{R} \times S^1)/\mathbb{Z}$ (and induces the Kronecker foliation), defines a free and proper $\mathbb{Z}$-action on the path space of $\mathbb{R} \times S^1$, and therefore on $\mathbb{R} \times \mathbb{R} \times S^1$. It is readily verified that the quotient defines a groupoid $(\mathbb{R} \times \mathbb{R} \times S^1)/\mathbb{Z} \rra T^2$, and that it equals the monodromy groupoid (and holonomy groupoid) of the Kronecker foliation on $T^2$. In particular, the total spaces of the monodromy and holonomy groupoids are Hausdorff in this case.
\end{rema}
\begin{rema}\cite{hoyohausdorff}\label{rema: monodromy and holonomy groupoids R^3 setminus 0}
Consider the foliation by horizontal planes on $\mathbb{R}^3 \setminus \{0\}$ defined by the submersion $\pr_3: \mathbb{R}^3 \setminus \{0\} \ra \mathbb{R}$ (see Example \ref{exam: R^3 setminus 0 foliation}). Notice that the fundamental group of the leaf which is diffeomorphic to $\mathbb{R}^2 \setminus \{0\}$ is isomorphic to $\pi_1(S^1) \cong \mathbb{Z}$. The sequence $\id(1,0,\tfrac{1}{n})$ in $\Mon(\mathbb{R}^3 \setminus \{0\},\mathcal{F})$ converges, but to any point in the isotropy group of $(1,0,0)$ (this follows because, for any $n$, $\id(1,0,\tfrac{1}{n})=[\gamma_{k,n}]$, where $\gamma_{k,n}$ is a loop in $\mathbb{R}^2 \times \{\tfrac{1}{n}\}$ that ``loops'' around the $z$-axis $k$ times). This shows that $\id(\mathbb{R}^3 \setminus \{0\})$ is not closed in the monodromy groupoid, and therefore the total space is not Hausdorff (see Proposition \ref{prop: G Hausdorff iff X closed}). The holonomy of foliation paths of this foliation are readily verified to be trivial, and the holonomy groupoid is given by $(\mathbb{R}^3 \setminus \{0\}) \tensor[_{\pr_3}]{\times}{_{\pr_3}} (\mathbb{R}^3 \setminus \{0\}) \rra \mathbb{R}^3 \setminus \{0\}$, whose total space is Hausdorff. 
\end{rema}
\begin{rema}\cite{hoyohausdorff}\label{rema: monodromy and holonomy groupoids of e^-1/t}
The foliation on $S^1 \times \mathbb{R}$ induced by the vector field $\sigma(\theta,t) \coloneqq \tfrac{\partial}{\partial \theta} + f(t)\tfrac{\partial}{\partial t}$ (see Example \ref{exam: e^-1/t}; also for the map $f$) has a Hausdorff monodromy groupoid, but a non-Hausdorff holonomy groupoid. Indeed, a path from $x$ to $y$, which is contained in a leaf diffeomorphic to $S^1$, is homotopic to a path which is the concatenation of a loop which ``loops'' around $n$ times together with a shortest path from $x$ to $y$, and paths in a leaf which is diffeomorphic to $\mathbb{R}$ are all homotopic to each other relative to endpoints. 
From this it is readily verified that the total space of the monodromy groupoid is given by $\mathbb{R} \times S^1 \times \mathbb{R}$, and it is equipped with the product topology (in fact, the monodromy groupoid is given by the action groupoid $\mathbb{R} \ltimes (S^1 \times \mathbb{R})$, with action given by the flow of $\sigma$). To see that the holonomy groupoid is not Hausdorff, we claim that any holonomy class $[\gamma]$ of a loop $\gamma$ in the leaf $S^1 \times \{0\} \subset S^1 \times \mathbb{R}$, say with $(\theta,0) \coloneqq \gamma(0)$, is the limit of the sequence $[(\theta,-\tfrac{1}{n})]$ (where we view $(\theta,-\tfrac{1}{n})$ as a constant path; see Proposition \ref{prop: G Hausdorff iff X closed}). First of all, $[\gamma]$ and $[\gamma(0)]$ represent different holonomy classes, because a (representative of a) germ of a diffeomorphism $T \ra T$ of a transversal $T$ that $\gamma$ induces is not the identity map when restricted to $T \cap (S^1 \times \mathbb{R}_{>0})$ (by our definition of $f$). However, if $\gamma'$ is a path contained in a leaf $S^1 \times \{t\}$, with $t<0$, then we can take a small transversal $T$ that is contained in $\mathbb{R}_{<0} \times S^1$, and then $\gamma'$ does induce the identity map $\id_T$, so it has trivial holonomy, and the claim follows.
\end{rema}
\begin{rema}\cite{moerdijk_mrcun_2003}\label{rema: monodromy and holonomy groupoids are not Hausdorff in general}
In case of the Reeb foliation of $S^3$ (see Example \ref{exam: Reeb foliation}), the monodromy and holonomy groupoids coincide. In a leaf diffeomorphic to $\mathbb{R}^2$, each path is homotopic (and therefore holonomic) to any other path relative to endpoints. However, in the leaf $T^2 \subset S^3$, the homotopy classes are generated by the generator of the fundamental group of $S^1$; any path from $x$ to $y$ in $T^2$ is homotopic to a path that is the concatenation of a loop, based at $x$, which is ``looped'' around the torus ``vertically'' $n \in \mathbb{Z}$ times and ``horizontally'' $m \in \mathbb{Z}$ times, with a shortest path from $x$ to $y$ (here, we view $T^2$ as an embedded submanifold of $\mathbb{R}^3$ in the standard way).
From this discussion, we can deduce that the total space of the monodromy (= holonomy) groupoid is equal to $S^3 \times \mathbb{R}^2$. The topology on $S^3 \times \mathbb{R}^2$, however, is not the product topology. To see that this topology is not Hausdorff, observe that if we take a loop $\gamma$ that loops around the torus vertically once, then this loop $\gamma$ and the constant path $x \coloneqq \gamma(0)$ are, in the monodromy (= holonomy) groupoid, not separable by open subsets. Indeed, take a basis open subset around $[\gamma]$, say $\varphi_x(V_x \times V_x' \times W_x)$ with $[\gamma]=\varphi(v,v',w)$ (see the proof of Proposition \ref{prop: Mon and Hol are Lie groupodis}). Then, for any point $(v,v',w') \in V_x \times V_x' \times W_x$, with $w'\neq w$, $\varphi(v,v',w')$ is represented by a path in a leaf different from $T^2$, hence contractible. A similar reasoning holds for any basis open subset around $[x]$, so any basis open subset around $[\gamma]$ and any other basis open subset around $[x]$ have non-empty intersection (alternatively, as we did in Remarks \ref{rema: monodromy and holonomy groupoids R^3 setminus 0} and \ref{rema: monodromy and holonomy groupoids of e^-1/t}, one can construct a sequence which has multiple limits).
\end{rema}
Recall that the universal cover of a connected manifold $M$, with $x_0 \in M$, can be realised as $p: \widetilde M \ra M$, where $\widetilde M$ is the path space consisting of homotopy classes of paths starting at $x_0$, and $p([\gamma]) \coloneqq \gamma(1)$. If we compare this with $\Mon(M,\mathcal{F})$ and $\Hol(M,\mathcal{F})$ of a foliated manifold $(M,\mathcal{F})$, we obtain the following result. 
\begin{prop}\cite{moerdijk_mrcun_2003}\label{prop: (universal) cover of Mon and Hol}
If $(M,\mathcal{F})$ is a foliated manifold, and $L \ni x_0$ is a leaf, then the target map of $\Mon(M,\mathcal{F})$ restricts to the map
\[\source^{-1}_{\Mon}(x_0) \xra{\target_{\Mon}%|_{\source_{\Mon}^{-1}(x_0)}
} L\]
which is the universal cover of $L$. Similarly, the target map of $\Hol(M,\mathcal{F})$ restricts to the map
\[\source^{-1}_{\Hol}(x_0) \xra{\target_{\Hol}%|_{\source_{\Hol}^{-1}(x_0)}
} L\]
which is a covering projection of $L$ induced by the kernel of the map $\hol: \pi_1(L,x_0) \ra \Hol(L,x_0)$ (see Definition \ref{defn: holonomy group of a foliation}). Moreover, the morphism of Lie groupoids $F: \Mon(M,\mathcal{F}) \ra \Hol(M,\mathcal{F})$ (see Remark \ref{rema: groupoid morphism Mon -> Hol}) restricts to a covering projection $\source^{-1}_{\Mon}(x_0) \ra \source^{-1}_{\Hol}(x_0)$.
\end{prop}

Later, when introduced to the theory of Lie algebroids, we will also see the notion of singular foliation. This concept generalises the foliations defined here in the sense that one allows partitions into immersed submanifolds varying in dimension% while satisfying a local condition.
. This theory is a bit more subtle, and it will be useful to not use a geometric approach, but instead pass to the theory of sheaves (note: if we regard a foliation as an involutive distribution, then the Serre-Swan theorem leads the way to an approach using sheaf-theory). 

\subsection{Morita equivalence}\label{sec: Morita equivalence}
There is an equivalence relation on the collection of groupoids. If $\groupoid$ is a groupoid, recall the notion of a $\group$-space (see Definition \ref{defn: action of a groupoid}). To define Morita equivalence, we will use the notion of a (right) principal $\group$-bundle. In fact, the theory of principal $\group$-bundles is very similar to the theory of the usual principal bundles in the following sense. We fix a groupoid $\groupoid$ and a right $\group$-space $\mu: P \ra \base$ with multiplication map $\nu$ (note: recall that, for a right $\group$-space, we switch the roles of $\source$ and $\target$ in Definition \ref{defn: action of a groupoid}).
\begin{defn}\cite{gerbes}\label{defn: free and proper groupoid action}
The groupoid action is called \textit{free} if whenever $p \cdot g = p$, then $g=\id_{\mu(p)}$. The groupoid action is called \textit{proper} if $(\pr_1,\nu): P \tensor[_{\mu}]{\times}{_{\target}} \group \ra P \times P$ is a proper map.
\end{defn}
The following result holds also in the groupoid setting. 
\begin{prop}\cite{gerbes}\label{prop: quotient by groupoid action}
If the groupoid action is free and proper, then $P/\group$ is a smooth manifold. 
\end{prop}\vspace{-\baselineskip}
\begin{defn}\cite{gerbes}\label{defn: principal groupoid bundles}
A \textit{(right) principal $\group$-bundle over $M$} is a surjective submersion $\pi_M: P \ra M$ (where $P$ is a right $\group$-space as before), called the \textit{structure map}, such that $\pi_M(p)=\pi_M(p')$ if and only if there is a unique $g \in \group$ such that $p \cdot g = p'$.
\end{defn}
%Let $\subgroupoid$ be another groupoid, and let $\mu_Q: Q \ra N$ be a right $\subgroup$-space. Then a groupoid morphism $\varphi: \group \ra \subgroup$ turns $Q \ra \subbase$ into a right $\group$-bundle in a canonical way. Indeed, we define \[\pi: F \ra \subbase\]
%\begin{defn}\label{defn: principal groupoid bundles morphisms}
%Suppose $E \ra M$ is a principal $\group$-bundle, and let $F \ra N$ be a principal $\subgroup$-bundle. A morphism between them is a groupoid morphism $\varphi: \group \ra \subgroup$, together with a commutative diagram (of smooth maps)
%\begin{center}
%\begin{tikzcd}
%    E \ar[r,"\Phi_\varphi"] \ar[d] & F \ar[d] \\
%    M \ar[r] & N
%\end{tikzcd}
%\end{center}
%where $\Phi_\varphi$ is a $\group$-equivariant map with respect to the $\group$-action on $F$ induced by $\varphi$.
%\end{defn}
That free and proper Lie group actions are in one-to-one correspondence with principal bundles has an analogue for Lie groupoid actions as well:
\begin{prop}\cite{gerbes}
The $\group$-space $P$ is a principal $\group$-bundle over a manifold $M$ if and only if the groupoid action is free and proper, and $P/\group \cong M$.
\end{prop}
Fix another groupoid $\subgroupoid$ and suppose $P$ is, in addition to being a right $\group$-space, a left $\subgroup$-space. Moreover, assume that there are surjective submersions $\pi_\base: P \ra \base$ and $\pi_\subbase: P \ra \subbase$ such that $P$ is a left principal $\subgroup$-bundle over $\base$ and a right principal $\group$-bundle over $\subbase$. That is, $P$ carries simultaneously the structure of a left principal $\subgroup$-space and a right principal $\group$-space.
\begin{defn}\cite{gerbes}\label{defn: Morita equivalence}
If the two principal groupoid actions on $\base \xla{\pi_\base} P \xra{\pi_\subbase} \subbase$ commute, then this structure defines a \textit{Morita equivalence} between $\groupoid$ and $\subgroupoid$, and $\groupoid$ and $\subgroupoid$ are called \textit{Morita equivalent}; we write $\group \sim \subgroup$.
\end{defn}
As the definition suggests, Morita equivalence is an equivalence relation. Calling it an equivalence relation, however, is misleading, since the collection of pairs of groupoids form a category, and not a set. But, clearly, the reflexivity, symmetry, and transitivity conditions still make sense in a category. The objects of interest for us are the equivalence classes of Morita equivalence, which have a special name. 
\begin{defn}\cite{gerbes}\label{defn: Differentiable stack}
A \textit{differentiable stack} $\mathfrak{X}$ is an equivalence class with respect to Morita equivalence of Lie groupoids.
\end{defn}
\begin{rema}\cite{gerbes}\label{rema: Differentiable stack}
A \textit{stack} is a categorical notion, and by equipping a specific kind of stack with extra structure, it becomes a differentiable stack, which will be equivalent to how we define it here. While for some applications this way of introducing differentiable stacks is useful, we will not use it, and use the more geometric definition presented here.
\end{rema}
We end the section with a proof that Morita equivalence is actually an equivalence relation.
\begin{prop}\cite{gerbes}\label{prop: Morita equivalence is an equivalence relation}
Morita equivalence is an equivalence relation. 
\end{prop}
\begin{proof}
We fix groupoids $\groupoid$, $\grouptwoid$, and $\mathcal{P} \rra U$.

(transitivity) Recall that $\groupoid$ is a left and right $\group$-space. Since $\source=\source_\group$ is a surjective submersion and satisfies $\source(g)=\source(g')$ if and only if $g'=(g'g^{-1})g$, and $\target=\target_\group$ is a surjective submersion as well (and satisfies $\target(g)=\target(g')$ if and only if $g'=g(g^{-1}g')$), we see that $\base \xla{\source} \group \xra{\target} \base$ is a Morita equivalence $\group \sim \group$. 
%Since $\subbase \xla{\source_\subgroup} \subgroupoid \xra{\target_\subgroup} \subbase$ is a Morita equivalence following the same reasoning, and $\source_{\subgroup} = \source|_\subgroup$ and $\target_{\subgroup}=\target|_\subgroup$, we see that $(\groupoid,\subgroupoid) \sim (\groupoid,\subgroupoid)$.

(symmetry) A left principal $\group$-bundle $P$, with $\group$-action $\nu_L$ and moment map $\mu$, can be \textit{reversed} to a right principal $\group$-bundle by setting
\[\nu_R(p,g) \coloneqq \nu_L(g^{-1},p) \textnormal{ for all } (p,g) \in P \tensor[_{\mu}]{\times}{_\target} \group\]
Of course, similarly, a right principal $\group$-bundle can be reversed to a left principal $\group$-bundle. Now, if $\base \xla{\pi_\base} P \xra{\pi_\basetwo} \basetwo$ is a Morita equivalence $\group \sim \grouptwo$, then $\basetwo \xla{\pi_\basetwo} P \xra{\pi_\base} \base$, where we reversed the principal bundle structures on both $\base \xla{\pi_\base} P$ and $P \xra{\pi_\basetwo} \basetwo$, is a Morita equivalence $\grouptwo \sim \group$. 

(transitivity) Let $\base \xla{\pi_\base^P} P \xra{\pi_\basetwo^P} \basetwo$ and $\basetwo \xla{\pi_\basetwo^Q} Q \xra{\pi_U^Q} U$ be Morita equivalences. Since $\pi_\basetwo^P$ and $\pi_\basetwo^Q$ are submersions, $P \tensor[_{\pi_\basetwo^P}]{\times}{_{\pi_\basetwo^Q}} Q$ is a smooth manifold. Now, $\grouptwo$ acts on $P \tensor[_{\pi_\basetwo^P}]{\times}{_{\pi_\basetwo^Q}} Q$ by
\[k \cdot (p,q) = (pk^{-1},kq), \textnormal{ where } k \in \grouptwo \textnormal{ and } (p,q) \in P \tensor[_{\pi_\basetwo^P}]{\times}{_{\pi_\basetwo^Q}} Q.\]
Since $\grouptwo$ acts freely and properly on $P$ and $Q$, it is readily verified that $\grouptwo$ also acts freely and properly on $P \tensor[_{\pi_\basetwo^P}]{\times}{_{\pi_\basetwo^Q}} Q$, so $\widetilde{P} \coloneqq P \tensor[_{\pi_\basetwo^P}]{\times}{_{\pi_\basetwo^Q}} Q/\grouptwo$ is a smooth manifold. By definition of $\widetilde{P}$, we have two canonical projections
\[\pi_\base^{\widetilde{P}}: \widetilde{P} \ra \base \textnormal{ given by } [p,q] \mapsto \pi_\base^P(p) \textnormal{ and } \pi_U^{\widetilde{P}}: \widetilde{P} \ra U \textnormal{ given by } [p,q] \mapsto \pi_U^Q(q).\]
Since these maps are both induced from surjective submersions, they are surjective submersions. We equip $\widetilde{P} \ra \base$ with a right $\mathcal{P}$-action by setting
\[[p,q] \cdot h = [p, q \cdot h], \textnormal{ for all } h \in \mathcal{P} \textnormal{ and } [p,q] \in \widetilde{P} \tensor[_{\pi_\base^{\widetilde{P}}}]{\times}{_{\target_{\mathcal{P}}}} \mathcal{P},\]
and we equip $\widetilde{P} \ra U$ with a left $\group$-action by setting
\[g \cdot [p,q] = [g \cdot p, q], \textnormal{ for all } g \in \group \textnormal{ and } [p,q] \in \group \tensor[_{\target_\group}]{\times}{_{\pi^{\widetilde{E}}_U}} \widetilde{E}.\]
It is now readily verified that $\base \xla{\pi_\base^{\widetilde{P}}} \widetilde{P} \xra{\pi_U^{\mathcal{P}}} U$ is a Morita equivalence.
\end{proof}
\addtocontents{toc}{\protect\thispagestyle{myplain}}\newpage

\section{Lie algebroids}\label{sec: Lie algebroids}
In analogy to Lie groups, we have an infinitesimal object of a Lie groupoid, fittingly called the Lie algebroid. Since, roughly speaking, a Lie groupoid describes multiple smooth symmetries on equal footing, while a Lie group describes only one, a Lie algebroid is a vector bundle, while a Lie algebra is a vector space. Recall that to each identity element of the groupoid we have constructed the isotropy Lie group, so we can associate its corresponding Lie algebra. These Lie algebras are related by the Lie algebroid, but they might not all have the same dimension (in particular, they will, in general, not be the fibers of the underlying vector bundle). We will approach Lie algebroids by introducing and discussing the general notion. Later, we will explain how we obtain a Lie algebroid from a Lie groupoid. Much of the theory can be found in \cite{ruimarius,rui,meinrenken,moerdijk_mrcun_2003,mackenziebook,mackenziepaper1,mackenziepaper2}.

\begin{term}\label{term: vector bundle with Lie bracket}
Let $E \ra \base$ be a vector bundle. We call a Lie bracket $[\cdot,\cdot]_E$ on the global sections $\Gamma(E) \coloneqq \Gamma(X,E)$, viewed as an $\mathbb{R}$-vector space, a \textit{Lie bracket (on $E$)}.
\end{term}

\begin{defn}\cite{ruimarius}\label{defn: Lie algebroid}
A Lie algebroid $\algebroid$ is a vector bundle equipped with a Lie bracket $\bracket{\cdot,\cdot}$ and a vector bundle map $\anchor_{\algebr}: \algebr \ra T\base$, called the \textit{anchor map}, satisfying the Leibniz rule:
\[\bracket{\alpha,f\beta} = f\bracket{\alpha,\beta} + \anchor_{\algebr}(\alpha)(f)\beta, \textnormal{ where } \alpha,\beta \in \Gamma(\algebr) \textnormal{ and } f \in C^\infty(\base).\]\vspace{-\baselineskip}
\end{defn}
\begin{rema}\label{rema: Lie algebroid}
To explain the notation: $\anchor_{\algebr}: \algebr \ra T\base$ induces a $C^\infty(\base)$-linear map $\Gamma(\algebr) \ra \mathfrak{X}(\base)$, also called $\anchor_{\algebr}$, defined by $(\anchor_A\alpha)(x) = \anchor_A(x)(\alpha(x))$. Moreover, we wrote $\sigma(f)$ for the Lie derivative $\mathcal{L}_\sigma f$ of $f \in C^\infty(\base)$ along $\sigma \in \mathfrak{X}(\base)$. Later, we will also introduce the notation $\Lie_{\anchor_A}$ for $\anchor_A$ seen as a map of sections $\Gamma(\algebr) \ra \mathfrak{X}(\base)$.
\end{rema}
\begin{term}\label{term: Lie algebroid notation}
Whenever $\algebroid$ is a Lie algebroid, we always write $\bracket{\cdot,\cdot}$ for its Lie bracket and $\anchor_\algebr$ for its anchor map (or just $\anchor$ if no confusion is possible). We reserve $[\cdot,\cdot]$ for the usual Lie bracket on $T\base$. 
\end{term}
Standard examples of Lie algebroids include the tangent bundle of a manifold, a foliation of a manifold (as an involutive distribution; see Proposition \ref{prop: equivalent definitions of foliations}), and the cotangent bundle of a Poisson manifold. Later, we will introduce these, and more, examples in detail. To be completely clear: when discussing the general theory of vector bundles or Lie algebroids, \textbf{we will assume that the total and base spaces of a vector bundle are Hausdorff} (we will discuss the normal bundle and apply the construction to Lie groupoids, for example, so in such cases we do not assume the normal bundle to be Hausdorff, even though it is a vector bundle). For the rest of this section, we fix a Lie algebroid $\algebroid$.
%We note here that if $A$ is a Lie algebroid with bracket $\bracket{\cdot,\cdot}$, then $\anchor$ is uniquely determined by the Leibniz rule.
\begin{prop}\cite{meinrenken}\label{prop: anchor induces Lie algebra morphism on sections}
The anchor map $\anchor: \Gamma(\algebr) \ra \mathfrak{X}(\base)$ is a morphism of Lie algebras.
\end{prop}
\begin{proof}
Let $\alpha,\beta \in \Gamma(\algebr)$. Then, for all $f \in C^\infty(\base)$ and $\gamma \in \Gamma(\algebr)$, we have, by repeated applications of the Leibniz rule, that
\begin{align*}
    \anchor(\alpha)\left(\anchor(\beta)f\right)\gamma &= \bracket{\alpha,\anchor(\beta)(f)\gamma} - \anchor(\beta)(f)\bracket{\alpha,\gamma} \\
    %&= \bracket{\alpha,\bracket{\beta,fX} - f\bracket{\beta,X}} - \bracket{\beta,f\bracket{\alpha,X}} + f\bracket{\beta,\bracket{\alpha,X}} \\
    &= \bracket{\alpha,\bracket{\beta,f\gamma}} - \bracket{\alpha,f\bracket{\beta,\gamma}} - \bracket{\beta,f\bracket{\alpha,\gamma}} + f\bracket{\beta,\bracket{\alpha,\gamma}},
\end{align*}
so since the two middle-terms in the last expression are symmetric in $\alpha$ and $\beta$, we see that
\begin{align*}
    [\anchor(\alpha),\anchor(\beta)](f)\gamma &= \anchor(\alpha)\left(\anchor(\beta)f\right)\gamma - \anchor(\beta)\left(\anchor(\alpha)f\right)\gamma \\
    &= \bracket{\alpha,\bracket{\beta,f\gamma}} - \bracket{\beta,\bracket{\alpha,f\gamma}} + f\left(\bracket{\beta,\bracket{\alpha,\gamma}} - \bracket{\alpha,\bracket{\beta,\gamma}}\right) \\
    &= \bracket{\bracket{\alpha,\beta},f\gamma} - f\bracket{\bracket{\alpha,\beta},\gamma} = \anchor(\bracket{\alpha,\beta})(f)\gamma,
\end{align*}
where in the third equality we used twice the Jacobi identity. This proves the statement.
\end{proof}

In the rest of the section we will describe the local nature of Lie algebroids.
First of all, Lie algebroids induce a Lie algebroid structure on each open subset, and therefore we have a ``sheaf of Lie algebras''. %To prove this, we will use the concept of a ``Lie algebroid connection''.
%\begin{defn}\label{defn: Lie algebroid connection}
%Let $\algebroid$ be a Lie algebroid, and let $E \ra \base$ be a vector bundle. An $\algebr$-connection on $E \ra \base$
%\end{defn}
\begin{prop}\cite{mackenziebook}\label{prop: Lie bracket on A defines a Lie bracket on sheaf of A}
The Lie bracket $\bracket{\cdot,\cdot}$ on $\algebroid$ induces a Lie bracket on the sheaf of sections of $\algebr$; that is, for each open subset $U \subset \base$, $\Gamma(U,\algebr)$ inherits a Lie bracket from $\Gamma(\base,\algebr)$ such that if $V \subset U$ is an open subset, then the restriction map $\Gamma(U,\algebr) \ra \Gamma(V,\algebr)$ is a morphism of Lie algebras.
\end{prop}
\begin{proof}
Let $U \subset \base$ be an open subset. The construction of the Lie bracket on $\Gamma(U,\algebr)$ is as follows: let $\beta,\gamma \in \Gamma(U,\algebr)$. Whenever $x \in U$, then we can find extensions of $\beta|_V$ and $\gamma|_V$, for some open set $V \ni x$ with $V \subset U$, to global sections $\widetilde{\beta}$ and $\widetilde{\gamma}$, respectively% (which exist by a standard argument involving a partition of unity)
. We now define
\begin{equation}\label{eq: local Lie bracket is well-defined}
    [\beta,\gamma]_U \coloneqq \bracket{\widetilde{\beta},\widetilde{\gamma}}|_U.
\end{equation}
To see that this is well-defined, %we have to show that, whenever we have extensions of $\beta|_V$ and $\gamma|_V$, with $V \ni p$ and $V \subset U$, then the difference of the two expressions in \eqref{eq: local Lie bracket is well-defined} vanishes along $U$. First of all, 
it suffices to prove that the expression is well-defined in case $\gamma|_V=0$ for some open set $V \ni x$ with $V \subset U$. Indeed, if $\sigma, \tau \in \Gamma(U,\algebr)$, then for extensions $\widetilde{\sigma}^i$ and $\widetilde{\tau}^i$ ($i=1,2$) of $\sigma|_V$ and $\tau|_V$, respectively, we have
\begin{align*}
    \bracket{\widetilde{\sigma}^1,\widetilde{\tau}^1} - \bracket{\widetilde{\sigma}^2,\widetilde{\tau}^2} &= \bracket{\widetilde{\sigma}^2,\widetilde{\tau}^1} - \bracket{\widetilde{\sigma}^2,\widetilde{\tau}^1} + \bracket{\widetilde{\sigma}^1,\widetilde{\tau}^1} - \bracket{\widetilde{\sigma}^2,\widetilde{\tau}^2} \\
    &= \bracket{\widetilde{\sigma}^1-\widetilde{\sigma}^2, \widetilde{\tau}^1} - \bracket{\widetilde{\sigma}^2,\widetilde{\tau}^1-\widetilde{\tau}^2} = -\bracket{\widetilde{\tau}^1,\widetilde{\sigma}^1-\widetilde{\sigma}^2} - \bracket{\widetilde{\sigma}^2,\widetilde{\tau}^1-\widetilde{\tau}^2},
\end{align*}
so, since $\widetilde{\sigma}^1 - \widetilde{\sigma}^2$ and $\widetilde{\tau}^1 - \widetilde{\tau}^2$ are both extensions of $0=(\sigma-\sigma)|_V=(\tau-\tau)|_V$, the above expression equals to zero if \eqref{eq: local Lie bracket is well-defined} holds for all $\gamma \in \Gamma(U,\algebr)$ such that $\gamma|_V=0$ for some open set $V \ni x$ with $V \subset U$.

We have reduced the problem to showing that: if $\beta,\gamma \in \Gamma(U,\algebr)$, such that $\gamma|_V=0$ for some open subset $V \ni x$ with $V \subset U$, then $\bracket{\beta,\gamma}(y)=0$ for all $y \in V$. This is indeed the case: we can find a smooth map $f: \base \ra \mathbb{R}$ such that $f^{-1}(0)=y$ and $f^{-1}(1)= \base \setminus V$. Then, $\gamma=f\gamma$, so, by the Leibniz rule,
\[\bracket{\beta,\gamma}(y)=\bracket{\beta,f\gamma}(y)=f(y)\bracket{\beta,\gamma}(y) + \anchor_\algebr(\beta)(f)(y)\gamma(y)=0.\]
That the restriction map $\Gamma(U,\algebr) \ra \Gamma(V,\algebr)$ is a morphism of Lie algebras, is obvious by definition of the Lie bracket. This proves the statement.
%In this case, it is useful to assume $U$ is the domain of a chart, which is justified, because we can cover $U$ with charts defined on open subsets $U_j \subset U$ and then glue the sections $\bracket{\widetilde{\beta},\widetilde{\gamma}}|_{U_j}$ (uniquely) to the section $\bracket{\widetilde{\beta},\widetilde{\gamma}}|_{U}$. Now, write $\widetilde{\gamma} \coloneqq \textstyle\sum_j c_j \alpha^j$ for a local frame $\alpha^1,\dots,\alpha^r$ of $\algebr$ and smooth functions $c_j \in C^\infty(\base)$ that, by definition of $\widetilde{\gamma}$, have to vanish when restricted to $U$. Then, for all $x \in U$, we have
%\[\left(\anchor(\widetilde{\beta})(c_j)\widetilde{\gamma}\right)(x) = db_j(x)\left(\anchor(x)\widetilde{\beta}(x)\right) \cdot 0 = 0.\]
%Then, by the Leibniz rule,
%\[\bracket{\widetilde{\beta},\widetilde{\gamma}} = \sum_j c_j\bracket{\widetilde{\beta},\alpha^j} + \anchor(\widetilde{\beta})(c_j)\gamma\]
%equals to zero when restricted to $U$. This proves that \eqref{eq: local Lie bracket is well-defined} is well-defined. By construction, it is obvious that the restriction maps $\Gamma(U,\algebr) \ra \Gamma(V,\algebr)$ are Lie algebra morphisms.
\end{proof}

\begin{term}\label{term: Lie bracket on open subset}
From now on, whenever $\algebroid$ is a Lie algebroid, then, for all open subsets $U \subset \base$, we also write $\bracket{\cdot,\cdot}$ for the Lie bracket on $\Gamma(U,\algebr)$.
\end{term}
Next, we will write down the so-called \textit{structure functions and equations}. These equations encode the defining properties of the Lie algebroid locally.
\begin{rema}\cite{rui}\label{rema: Lie algbroid local nature}
Let $U \subset \base$ be an open subset on which we have a vector bundle chart with local coordinates $(x^1,\dots,x^n): U \ra \mathbb{R}^n$ and local frame $\alpha^1,\dots,\alpha^r \in \Gamma(U,\algebr)$. Since, for all $1 \le j \le r$, $\anchor \alpha^j \in \mathfrak{X}(U)$, we can write
\[\anchor \alpha^j = \sum_i b^{ji} \frac{\partial}{\partial x^i}, \textnormal{ where } b^{ji} \in C^\infty(U),\]
and, similarly, for all $1 \le j,k \le r$, we can write
\[\bracket{\alpha^j,\alpha^k} = \sum_\ell c_\ell^{jk}\alpha^\ell, \textnormal{ where } c_\ell^{jk} \in C^\infty(U).\]
The smooth maps $b^{ji},c^{jk}_\ell \in C^\infty(U)$ are called the \textit{structure functions} of $A$. We can now translate the defining relations for $\bracket{\cdot,\cdot}$ to be a Lie bracket that satisfies the Leibniz rule with respect to $\anchor$ into partial differential equations in terms of the structure functions. The Leibniz rule, in these local coordinates, can be written as
\begin{equation}\label{eq: local expression Leibniz rule}
    \bracket{\alpha^j,f\alpha^k} = f\bracket{\alpha^j,\alpha^k} + \anchor \alpha^j(f)\alpha^k = f\sum_\ell c_\ell^{jk}\alpha^\ell + \sum_i b^{ji} \frac{\partial f}{\partial x^i} \alpha^k, \textnormal{ for all } f \in C^\infty(U).
\end{equation}
The anti-symmetry of the Lie bracket can be expressed as
\[\sum_\ell c_\ell^{jk}\alpha^\ell = [\alpha^j,\alpha^k] = -[\alpha^k,\alpha^j] = \sum_\ell -c_\ell^{kj}\alpha^\ell,\]
i.e. the $c^{jk}_\ell$ satisfy the relations
\begin{equation}\label{eq: structure equation 1}
    c_{\ell}^{jk}=-c_{\ell}^{kj} \textnormal{ for all } 1 \le \ell \le r.
\end{equation}
Before we write down the Jacobi identity locally, we introduce the \textit{cyclic sum} notation $\textstyle\bigodot$: it means that we sum over all cyclic permutations of the subscripts (see Definition \ref{defn: cyclic sum} below). Now,
\begin{align*}
    0 = \bigodot_{i,j,k}\bracket{\alpha^i,\bracket{\alpha^j,\alpha^k}} &= \bigodot_{i,j,k}\sum_\ell \bracket{\alpha^i,c_\ell^{jk}\alpha^\ell} = \bigodot_{i,j,k}\sum_\ell \left(c_\ell^{jk}\sum_s c_s^{i\ell}\alpha^s + \sum_{s'} b^{is'} \frac{\partial c_\ell^{jk}}{\partial x^{s'}} \alpha^\ell\right)
\end{align*}
(where we used the Leibniz rule in the last equality), so, by collecting the terms on the right-hand side, we conclude that the Jacobi identity in these local coordinates becomes
\begin{equation}\label{eq: structure equation 2}
    \bigodot_{i,j,k}\sum_s c_s^{jk}c_\ell^{is} + b^{is} \frac{\partial c_\ell^{jk}}{\partial x^s}=0 \textnormal{ for all } 1 \le i,j,k,\ell \le r.
\end{equation}
The equations (\ref{eq: structure equation 1}) and (\ref{eq: structure equation 2}) are called the \textit{structure equations} of $\algebr$.
\end{rema}
\begin{defn}\label{defn: cyclic sum}
Let $\{x_i\}_{i \in \mathbb{N}}$ be a subset of a set $\base$ and let $f: \base^n \ra \base$ be a map. Then 
\[\bigodot_{i_1,\dots,i_n} f(x_{i_1},\dots,x_{i_n}) \coloneqq f(x_{i_1},\dots,x_{i_n}) + f(x_{i_2},\dots,x_{i_n},x_{i_1}) + \dots + f(x_{i_n},x_{i_1},\dots,x_{i_{n-1}})\]
is called the \textit{cyclic sum} (with respect to $f$ and $i_1,\dots,i_n \in \mathbb{N}$).
\end{defn}
In the next section we will briefly introduce so-called \textit{anchored vector bundles}. These vector bundles are like Lie algebroids without necessarily having a Lie bracket. Some properties of Lie algebroids, of course, already hold for anchored vector bundles and sometimes we want to be precise about this. Moreover, it is often easy to describe an anchored vector bundle structure on a vector bundle. We will see that then it suffices to define Lie brackets on a local frame satisfying the structure equations described in Remark \ref{rema: Lie algbroid local nature}.

\subsection{Anchored vector bundles}\label{defn: anchored vector bundles}
First of all, to be precise, we put here our convention on vector subbundles of a vector bundle.
\begin{term}\label{term: vector subbundle}
Let $\pi_E: E \ra M$ be a vector bundle. By a vector subbundle $\pi_F: F \ra N$ we mean that $F \subset E$ and $N \subset M$ are embedded submanifolds, such that the structure maps of $E$, i.e. $\pi_E$, $+_E: E \tensor[_{\pi_E}]{\times}{_{\pi_E}} E \ra E$, and $\cdot_E: \mathbb{R} \times E \ra E$ restrict to the structure maps of $F$.
\end{term}
In Definition \ref{defn: Lie algebroid}, we defined a Lie algebroid as a vector bundle $\algebroid$ equipped with a vector bundle morphism $\algebr \ra T\base$, called the anchor map, and a Lie bracket such that the bracket and the anchor map are compatible via a Leibniz-type rule. Of course, we can define vector bundles $E \ra \base$ that are equipped only with either a vector bundle morphism $E \ra T\base$, or a Lie bracket. %In fact, giving a definition for the first type of vector bundle will be useful to indicate that some results, or definitions, make sense in this more general setting. 
For us, the latter type of vector bundle does not yield an interesting new concept; if we have a vector bundle $E \ra \base$ equipped with a Lie bracket, then two vector bundle morphisms $E \ra T\base$ that turn $E$ into a Lie algebroid must be equal due to the Leibniz rule. However, it will be useful, mainly for convenience and notational purposes, to introduce the first type of vector bundle.
\begin{defn}\cite{meinrenken}\label{defn: anchored vector bundles}
A vector bundle $E \ra \base$ is called an \textit{anchored vector bundle} if it comes with a vector bundle morphism $\anchor_E: E \ra T\base$, called the \textit{anchor map}.
\end{defn}

So, a Lie algebroid $\algebroid$ with anchor map $\anchor$ is an anchored vector bundle, and the Lie bracket $\bracket{\cdot,\cdot}$, in turn, uniquely determines the anchor map $\anchor$. 
\begin{term}\label{term: anchored vector bundle}
If $E \ra \base$ is an anchored vector bundle, then we denote its anchor map by $\anchor_E$ (or just $\anchor$ if no confusion is possible).
\end{term}

\begin{defn}\cite{meinrenken}\label{defn: anchored vector bundle morphism}
Let $E \ra \base$ and $F \ra \subbase$ be two anchored vector bundles. A \textit{morphism of anchored vector bundles} between them is a vector bundle morphism $E \ra F$ such that the diagram
\begin{center}
\begin{tikzcd}
    E \ar[r] \ar[d,"\anchor_E"] & F \ar[d,"\anchor_F"] \\
    T\base \ar[r] & T\subbase
\end{tikzcd}
\end{center}
commutes, where the lower horizontal map is the differential of the base map of $E \ra F$.
\end{defn}
Observe that we have not yet defined morphisms of Lie algebroids. The reason is that we need to do some preparation for this: it is not just a morphism of anchored vector bundles that intertwines the Lie brackets. Simply put, the reason is that this does not make sense in general: a vector bundle morphism does not necessarily induce a map on sections. A minimal example is the following.
\begin{exam}\label{rema: vector bundle does not induce map on sections}
Consider the vector bundles $E \coloneqq T\mathbb{R} \ra \mathbb{R}$ and $F \coloneqq \mathbb{R} \ra \{\pt\}$. Then $E$ and $F$ are naturally anchored vector bundles and the map
\[g: E \ra F \textnormal{ given by } (x,v) \mapsto (\pt,v)\]
defines an anchored vector bundle morphism. Moreover, $E$ comes with the usual Lie bracket on $T\mathbb{R}$ and $F$ comes with the zero Lie bracket, so they are even Lie algebroids. However, the map $g$ does not induce a map
\[G: \Gamma(E) \ra \Gamma(F),\]
since if $\sigma \in \Gamma(E)$, then $g(0,\sigma(0))=(\pt,\sigma(0))$ and $g(1,\sigma(1))=(\pt,\sigma(1))$ should both be equal to $G(\sigma)(\pt)$, which is not the case if $\sigma(0) \neq \sigma(1)$.
\end{exam}
Nonetheless, vector bundles over the same base do induce a map on sections (in the obvious way). In fact, we will see that if two Lie algebroids have the same base, then a morphism between them is a morphism of anchored vector bundles that intertwines the Lie brackets, but more on this later. We will now describe how to try to turn an anchored vector bundle into a Lie algebroid (with the same anchor map).
\begin{rema}\label{rema: defining a Lie algebroid locally out of an anchored vector bundle}
Let $E \ra \base$ be an anchored vector bundle with anchor map $\anchor: E \ra T\base$. If we want to describe a Lie algebroid structure on $E$ compatible with the anchor map (via the Leibniz rule), then we can do this locally. By this we mean that we can describe a Lie bracket on each $E|_{U_t}$ for an open cover $\{U_t\}$ of $\base$ (which are the domains of vector bundle charts), and impose conditions locally to ensure the Lie bracket glues to a Lie bracket on $E$ compatible with $\anchor$. To do this, let $s^1,\dots,s^r$ be a local frame for $E$ defined on an open set $U \subset \base$. Write, for all $1 \le j \le r$, \[\anchor s^j = \sum_i b^{ji} \frac{\partial}{\partial x^i}, \textnormal{ where } b^{ji} \in C^\infty(U).\] 
We now define, for all $1 \le j,k \le r$, the expressions
\[[s^j,s^k]_E \coloneqq \sum_\ell c^{jk}_\ell s^\ell, \textnormal{ where } c^{jk}_\ell \in C^\infty(U).\]
Then, following the equations \eqref{eq: local expression Leibniz rule} from Remark \ref{rema: Lie algbroid local nature} (the Leibniz rule), there is at most one extension of the above expressions to all of $\Gamma(U,E)$ that could define a Lie bracket. In particular, in order for the definitions of $[s^j,s^k]_E$ to describe a Lie bracket on $E|_U$, it suffices to extend the definition to all of $\Gamma(U,E)$ using the Leibniz rule and check that the structure equations \eqref{eq: structure equation 1} (skew-symmetry) and \eqref{eq: structure equation 2} (Jacobi-identity) of Remark \ref{rema: Lie algbroid local nature} are satisfied. Indeed, observe that the maps
\[[\bullet,\star]_E + [\star,\bullet]_E, \textnormal{ and } \bigodot_{\bullet,\star,\diamond} [[\bullet,\star]_E,\diamond]_E \textnormal{ (the Jacobiator)}\]
are $C^\infty(U)$-linear by the Leibniz rule, so they are satisfied if and only if they are satisfied in case all sections involved are part of the local frame $s^1,\dots,s^r$.
Now, if we describe a Lie algebroid structure on each $E|_{U_t}$ for an open cover $\{U_t\}$ of $\base$ (which are the domains of vector bundle charts), then, if we make sure the defining functions $c^{jk}_\ell$ agree on overlaps, the Lie brackets on $E|_{U_t}$ glue to a Lie bracket $[\cdot,\cdot]_E$ on $E$. Indeed, the $\mathbb{R}$-bilinearity, skew-symmetry, Jacobi-identity, and Leibniz rule can all be phrased as describing equalities of sections, so if these relations are checked locally on the $U_t$, where they hold by assumption, then the relations also hold on all of $\base$ (by uniqueness of gluing sections).
\end{rema}

\subsection{Lie algebroid cohomology}\label{sec: Lie algebroid cohomology}
Many things that are defined for tangent bundles make sense for Lie algebroids as well. Here we introduce the de Rham-complex for Lie algebroids. This concept generalises many cohomology theories, like the usual de Rham-cohomology, but also the Chevalley-Eilenberg cohomology and Poisson cohomology; more on that later. We start by recalling the notion of graded algebras.
\begin{defn}\cite{AtiyahMac}\label{defn: graded derivation}
Let $T \coloneqq \textstyle\bigoplus_{n \ge 0} T_n$ be a graded algebra over a commutative ring $R$ with unit. An element $a \in T_n$ is called \textit{homogeneous of degree $n$}, and we denote $T_{\textnormal{hom}}$ for the set of homogeneous elements, which comes with a canonical map $T_{\textnormal{hom}} \ra \mathbb{Z}_{\ge 0}$ denoted by $|\cdot|$. For $n < 0$, $T_n \coloneqq 0$.
\end{defn}
The main reason we put this notion here is to properly introduce graded derivations, which allows us to ease the notation in this section. The reason these notions will come up, is because we can turn $\Omega^\bullet(E) \coloneqq \Gamma(\textstyle\bigwedge^{\bullet} E^*)$, where $E \ra \base$ is a vector bundle, into a graded algebra over $C^\infty(\base)$ by introducing the wedge product.
\begin{defn}\cite{meinrenken}
An $R$-linear map $D: T \ra T$ is called a \textit{(graded) derivation of degree} $d$ if it restricts to maps $T_k \ra T_{k+d}$ and it satisfies the Leibniz-type rule
\[D(ab) = D(a) \cdot b + (-1)^{|a|  d}a \cdot D(b) \textnormal{ for all } a,b \in T_{\textnormal{hom}}.\]
If $D$ is a graded derivation of degree $1$, and $D^2=0$, then it is called a \textit{(graded) differential}.
\end{defn}
\begin{prop}
Let $D$ and $D'$ be two graded derivations of degrees $d$ and $d'$, respectively, of a graded algebra $T$. Then 
\[[D,D'] \coloneqq DD'-(-1)^{dd'}D'D\] 
is a graded derivation of degree $d+d'$.
\end{prop}
\begin{proof}
Let $a,b \in T_{\textnormal{hom}}$. Then,
\begin{align*}
    DD'(ab) &= D\left(D'(a) \cdot b + (-1)^{|a|  d'}a \cdot D'(b)\right) \\
    &= DD'(a) \cdot b + (-1)^{(|a|+d')  d}D'(a) \cdot D(b) \\
    &+ (-1)^{|a| d'}D(a) \cdot D'(b) + (-1)^{|a|  (d+d')}a \cdot DD'(b).
\end{align*}
The result follows by noting that the middle two terms of this expression vanish against the middle two terms of the same expression for $-(-1)^{dd'}D'D(ab)$.
\end{proof}
In the rest of this section, we fix a Lie algebroid $\algebroid$. We can define a Lie derivative, interior derivative, and differential on $\algebr$ by using the same formulas as in the case $\algebr=T\base$, but by making sure they make sense for Lie algebroids. This suggest the following definition.

\begin{defn}\cite{meinrenken}\label{defn: de Rham differential and Lie derivative}
For all $\alpha \in \Gamma(A)$, $k \ge 0$, the \textit{Lie derivative} $\Lie_{\anchor}(\alpha): \Omega^k(A) \ra \Omega^k(A)$ is given by
\begin{align*}
\left(\Lie_{\anchor}(\alpha)\omega\right)(\alpha_1,\dots,\alpha_k) = &\anchor(\alpha)\left(\omega(\alpha_1,\dots,\alpha_k)\right) \\
- &\sum_i\omega(\alpha_1,\dots,\bracket{\alpha,\alpha_i},\dots,\alpha_k) , \textnormal{ where } \alpha_1,\dots,\alpha_k \in \Gamma(A),
\end{align*}
and the \textit{interior derivative} $\iota_\alpha:  \Omega^k(A) \ra \Omega^{k-1}(A)$ is given by
\begin{align*}
    \left(\iota_\alpha\omega\right)(\alpha_1,\dots,\alpha_k) = \omega(\alpha,\alpha_1,\dots,\alpha_k) \textnormal{ if } k\ge1 \textnormal{ and } \iota_\alpha\omega=0 \textnormal{ if } k=0.
\end{align*}
Moreover, for all $k \ge 0$, the \textit{de Rham-differential} $d_\algebr: \Omega^k(A) \ra \Omega^{k+1}(A)$ is given by
\begin{align*}
    \left(d_A\omega\right)(\alpha_1,\dots,\alpha_{k+1}) = &\sum_{i < j} (-1)^{i+j} \omega(\bracket{\alpha_i,\alpha_j},\alpha_1,\dots,\widehat{\alpha_i},\dots,\widehat{\alpha_j},\dots,\alpha_{k+1}) \\
    - &\sum_i (-1)^i\anchor(\alpha_i)(\omega(\alpha_1,\dots,\widehat{\alpha_i},\dots,\alpha_{k+1})), \textnormal{ where } \alpha_1,\dots,\alpha_{k+1} \in \Gamma(A).
\end{align*}\vspace{-\baselineskip}
\end{defn}
Notice that for all $\alpha \in \Gamma(A)$, $\anchor(\alpha) = \Lie_{\anchor}(\alpha)$ and that for a vector bundle $E \ra \base$ we can define the interior derivative as well (by using the same formula). A tedious calculation shows:
\begin{prop}\cite{meinrenken}\label{prop: d_A, iota and Lie derivative are graded derivations}
The maps $d_A$ assemble into a graded derivation of degree $1$, and for all $\alpha \in \Gamma(A)$, the maps $\iota_\alpha$ and the maps $\Lie_{\anchor}(\alpha)$ assemble into graded derivations of degree $-1$ and $0$, respectively.
\end{prop}
In fact, $d_A$ is a differential. Verifying this is not hard, and we will use the identities of \textit{Cartan Calculus} that hold in $\Omega^\bullet(A)$ to prove this. In fact, conversely, we will see that a vector bundle, equipped with a differential, is a Lie algebroid. We will prove this fact, and for the proof it will be useful to prove the Cartan Calculus identities as well.
\begin{lemm}\cite{meinrenken}\label{lemm: interior product of bracket}
For all $\alpha,\beta \in \Gamma(A)$, we have 
\[[\Lie_{\anchor}(\alpha),\iota_\beta] = \Lie_{\anchor}(\alpha) \circ \iota_\beta - \iota_\beta \circ \Lie_{\anchor}(\alpha) = \iota_{\bracket{\alpha,\beta}}.\]\vspace{-\baselineskip}
\end{lemm}
\begin{proof}
Put $\alpha_1 \coloneqq \beta$, let $\omega \in \Omega^k(\algebr)$ and let $\alpha_2,\dots,\alpha_k \in \Gamma(A)$. Then 
\begin{align*}
    &\iota_{\alpha_1} \circ \left(\Lie_{\anchor}(\alpha)\omega\right)(\alpha_2,\dots,\alpha_k) = \anchor(\alpha)\left(\omega(\alpha_1,\dots,\alpha_k)\right) - \sum_i\omega(\alpha_1,\dots,\bracket{\alpha,\alpha_i},\dots,\alpha_k) \\
    &= \anchor(\alpha)\left(\iota_{\alpha_1}\omega(\alpha_2,\dots,\alpha_k)\right) - \sum_{i>1}\iota_{\alpha_1}\omega(\alpha_2,\dots,\bracket{\alpha,\alpha_i},\dots,\alpha_k) - \iota_{\bracket{\alpha,\alpha_1}}\omega(\alpha_2,\dots,\alpha_k) \\
    &= (\Lie_{\anchor}(\alpha)\iota_{\alpha_1}\omega)(\alpha_2,\dots,\alpha_k) - \iota_{\bracket{\alpha,\alpha_1}}\omega(\alpha_2,\dots,\alpha_k),
\end{align*}
which proves the statement.
\end{proof}

\begin{lemm}\cite{meinrenken}[Magic formula]\label{lemm: Cartan's magic formula}
For all $\alpha \in \Gamma(A)$, we have
\[[d_A,\iota_\alpha] = d_A \circ \iota_\alpha + \iota_\alpha \circ d_A = \Lie_{\anchor}(\alpha).\]\vspace{-\baselineskip}
\end{lemm}
\begin{proof}
Put $\alpha_1 \coloneqq \alpha$, let $\omega \in \Omega^k(\algebr)$ and let $\alpha_2,\dots,\alpha_{k+1} \in \Gamma(A)$. Then
\begin{align}
    d_\algebr(\iota_{\alpha_1}\omega)(\alpha_2,\dots,\alpha_{k+1}) = &\sum_{2 \le i < j \le k+1} (-1)^{i+j} (\iota_{\alpha_1}\omega)(\bracket{\alpha_{i},\alpha_{j}},\alpha_2,\dots,\widehat{\alpha_{i}},\dots,\widehat{\alpha_{j}},\dots,\alpha_{k+1}) \nonumber\\
    -&\sum_{i=2}^{k+1} (-1)^{i-1}\anchor(\alpha_i)((\iota_{\alpha_1}\omega)(\alpha_2,\dots,\widehat{\alpha_{i}},\dots,\alpha_{k+1})) \label{eq: magic formula}
\end{align}
and $\iota_{\alpha_1}(d_\algebr\omega)(\alpha_2,\dots,\alpha_{k+1})$ is simply the expression from Definition \ref{defn: de Rham differential and Lie derivative}. Now, if we add the two expressions, then all terms in \eqref{eq: magic formula} vanish (by using the skew-symmetry of $\bracket{\cdot,\cdot}$), and therefore
\begin{align*}
    (d_\algebr(\iota_{\alpha_1}\omega) + \iota_{\alpha_1}(d_\algebr\omega))(\alpha_2,\dots,\alpha_{k+1}) &=
    %&\sum_{2 \le i < j \le k+1} (-1)^{i+j} (\iota_{\alpha_1}\omega)(\bracket{\alpha_{i},\alpha_{j}},\alpha_2,\dots,\widehat{\alpha_{i}},\dots,\widehat{\alpha_{j}},\dots,\alpha_{k+1}) \\
    %&-\sum_{i=2}^{k+1} (-1)^{i+1}\anchor(\alpha_i)((\iota_{\alpha_1}\omega)(\alpha_2,\dots,\widehat{\alpha_{i}},\dots,\alpha_{k+1})) \\
    %&+ \sum_{1 \le i < j \le k+1} (-1)^{i+j} \omega(\bracket{\alpha_{i},\alpha_{j}},\alpha_1,\dots,\widehat{\alpha_{i}},\dots,\widehat{\alpha_{j}},\dots,\alpha_{k+1}) \\
    %&-\sum_{i=1}^{k+1} (-1)^i\anchor(\alpha_i)(\omega(\alpha_1,\dots,\widehat{\alpha_{i}},\dots,\alpha_{k+1})) \\
    -\sum_{2 \le j \le k+1} (-1)^j \omega(\bracket{\alpha_1,\alpha_{j}},\alpha_2,\dots,\widehat{\alpha_{j}},\dots,\alpha_{k+1}) \\
    &+ \anchor(\alpha_1)(\omega(\alpha_2,\dots,\alpha_{k+1})) \\
    &= -\sum_{2 \le j \le k+1} \omega(\alpha_2,\dots,\bracket{\alpha_1,\alpha_{j}},\dots,\alpha_{k+1}) \\
    &+ \anchor(\alpha_1)(\omega(\alpha_2,\dots,\alpha_{k+1})) = (\Lie_\anchor(\alpha_1)\omega)(\alpha_2,\dots,\alpha_{k+1}).
\end{align*}
This proves the statement.
\end{proof}
The following three identities we will prove using the fact that if for all $\alpha \in \Gamma(A)$, $\iota_\alpha$ commutes with a given derivation, then the derivation must be zero, because $\iota_\alpha f=0$ for all smooth functions $f \in C^\infty(\base)=\Omega^0(A)$, so the derivation is zero by using an inductive argument.
\begin{lemm}\cite{meinrenken}\label{lemm: Lie derivative of bracket}
For all $\alpha,\beta \in \Gamma(A)$, we have \[[\Lie_{\anchor}(\alpha),\Lie_{\anchor}(\beta)] = \Lie_{\anchor}(\alpha) \circ \Lie_{\anchor}(\beta) - \Lie_{\anchor}(\beta) \circ \Lie_{\anchor}(\alpha) = \Lie_{\anchor}(\bracket{\alpha,\beta}).\]\vspace{-\baselineskip}
\end{lemm}
\begin{proof}
Let $\gamma \in \Gamma(A)$. Then
\begin{align*}
&\iota_\gamma \circ (\Lie_{\anchor}(\bracket{\alpha,\beta}) - \Lie_{\anchor}(\alpha) \circ \Lie_{\anchor}(\beta) + \Lie_{\anchor}(\beta) \circ \Lie_{\anchor}(\alpha)) \\
%&= (\Lie_{\anchor}(\bracket{\alpha,\beta}) \iota_\gamma - \iota_{\bracket{\bracket{\alpha,\beta},\gamma}}) - (\Lie_{\anchor}(\alpha)\iota_\gamma - \iota_{\bracket{\alpha,\gamma}}) \circ \Lie_{\anchor}(\beta) + (\Lie_{\anchor}(\beta)\iota_\gamma - \iota_{\bracket{\beta,\gamma}}) \circ \Lie_{\anchor}(\alpha) \\
&=(\Lie_{\anchor}(\bracket{\alpha,\beta}) \circ \iota_\gamma - \iota_{\bracket{\bracket{\alpha,\beta},\gamma}}) - \Lie_{\anchor}(\alpha) \circ (\Lie_{\anchor}(\beta) \circ \iota_\gamma - \iota_{\bracket{\beta,\gamma}}) + (\Lie_{\anchor}(\beta) \circ \iota_{\bracket{\alpha,\gamma}} - \iota_{\bracket{\beta,\bracket{\alpha,\gamma}}}) \\ 
&+ \Lie_{\anchor}(\beta) \circ (\Lie_{\anchor}(\alpha) \circ \iota_\gamma - \iota_{\bracket{\alpha,\gamma}}) - (\Lie_{\anchor}(\alpha) \circ \iota_{\bracket{\beta,\gamma}} - \iota_{\bracket{\alpha,\bracket{\beta,\gamma}}}) \\
&= (\Lie_{\anchor}(\bracket{\alpha,\beta}) - \Lie_{\anchor}(\alpha) \circ \Lie_{\anchor}(\beta) + \Lie_{\anchor}(\beta) \circ \Lie_{\anchor}(\alpha)) \circ \iota_\gamma - \bigodot_{\alpha,\beta,\gamma} \iota_{\bracket{\bracket{\alpha,\beta},\gamma}},
\end{align*}
where we used Lemma \ref{lemm: interior product of bracket} in the first equality and the skew-symmetry of $\bracket{\cdot,\cdot}$ in the second equality. By the Jacobi identity, the result follows.
\end{proof}
\begin{lemm}\cite{meinrenken}\label{lemm: d_A commutes with Lie derivative}
For all $\alpha \in \Gamma(A)$, we have 
\[[\Lie_{\anchor}(\alpha),d_A] = \Lie_{\anchor}(\alpha) \circ d_A - d_A \circ \Lie_{\anchor}(\alpha)=0.\]\vspace{-\baselineskip}
\end{lemm}
\begin{proof}
Let $\beta \in \Gamma(A)$. Then
\begin{align*}
    &\iota_{\beta}\circ (\Lie_{\anchor}(\alpha) \circ d_A - d_A \circ \Lie_{\anchor}(\alpha)) = (\Lie_{\anchor}(\alpha) \circ \iota_{\beta} - \iota_{\bracket{\alpha,\beta}}) \circ d_A - (\Lie_{\anchor}(\beta) - d_A \circ \iota_\beta) \circ \Lie_{\anchor}(\alpha) \\
    &= \Lie_{\anchor}(\alpha) \circ (\Lie_{\anchor}(\beta) - d_A \circ \iota_\beta) - \iota_{\bracket{\alpha,\beta}} \circ d_A - \Lie_{\anchor}(\beta) \circ \Lie_{\anchor}(\alpha) + d_A \circ (\Lie_{\anchor}(\alpha) \circ \iota_\beta - \iota_{\bracket{\alpha,\beta}}) \\
    %&=  -\Lie_{\anchor}(\alpha) \circ d_A \circ \iota_\beta + d_A \circ \iota_{\bracket{\alpha,\beta}} + d_A \circ (\Lie_{\anchor}(\alpha) \circ \iota_\beta - \iota_{\bracket{\alpha,\beta}}) \\
    &=  -(\Lie_{\anchor}(\alpha) \circ d_A - d_A \circ \Lie_{\anchor}(\alpha)) \circ \iota_\beta,
\end{align*}
where we used Lemma \ref{lemm: Cartan's magic formula} and Lemma \ref{lemm: Lie derivative of bracket} in the last equality. This proves the statement.
\end{proof}

\begin{prop}\cite{meinrenken}\label{prop: d_A^2=0}
We have $d_A^2=0$.
\end{prop}
\begin{proof}
Let $\alpha \in \Gamma(A)$. Then
\begin{align*}
    \iota_\alpha \circ d_A^2 &= (\Lie_{\anchor}(\alpha) - d_A \circ \iota_\alpha) \circ d_A \\
    &= \Lie_{\anchor}(\alpha) \circ d_A - d_A \circ (\Lie_{\anchor}(\alpha) - d_A \circ \iota_\alpha) \\
    &= d_A^2 \circ \iota_\alpha,
\end{align*}
where we used Lemma \ref{lemm: d_A commutes with Lie derivative} in the last equality. This proves the statement.
\end{proof}

From this we see that $(\Omega^\bullet(\algebr),d_\algebr)$ is a cochain complex.

\begin{defn}\cite{meinrenken}\label{defn: cohomology}
The cohomology groups $H^\bullet(\algebr)$ corresponding to the cochain complex $(\Omega^\bullet(\algebr),d_\algebr)$ are called the Lie algebroid cohomology groups.
\end{defn}

If $E \ra \base$ is a vector bundle, we mentioned that $\Omega^\bullet(E) \coloneqq \Gamma(\textstyle\bigwedge^\bullet E^*)$ admits a wedge product (analogous to the case $E=T\base$) which turns $\Omega^\bullet(E)$ into a graded algebra over the ring of smooth functions $C^\infty(\base)$. Here is the converse statement to the discussion above:
\begin{prop}\cite{meinrenken}\label{prop: differential determines Lie algebroid}
Let $E \ra \base$ be a vector bundle with a differential $d_E$ on $\Omega^\bullet(E)$ (i.e. a derivation of degree $1$ satisfying $d_E^2=0$). Then there is a unique Lie algebroid structure on $E$ for which $d_E$ is the de Rham-differential.
\end{prop}
\begin{proof}
For all $\alpha \in \Gamma(E)$ and $k \ge 0$ we can define a map $\Lie_{\anchor}(\alpha): \Omega^k(E) \ra \Omega^k(E)$ by 
\[\Lie_{\anchor}(\alpha) \coloneqq d_E \circ \iota_\alpha + \iota_\alpha \circ d_E.\]
Observe that these maps are derivations, because $\iota_\alpha$ and $d_E$ are derivations. Since $\Omega^0(E) = C^\infty(\base) = \Omega^0(T\base)$, the map $\alpha \mapsto (\Lie_{\anchor}(\alpha): C^\infty(\base) \ra C^\infty(\base))$, which is $C^\infty(\base)$-linear (because $\alpha \mapsto \iota_\alpha$ is), determines a morphism of vector bundles $\anchor: E \ra T\base$. We will show that the formula
\[\iota_{[\alpha,\beta]_E} \coloneqq \Lie_{\anchor}(\alpha) \circ \iota_\beta - \iota_\beta \circ \Lie_{\anchor}(\alpha): \Omega^1(E) \ra C^\infty(\base) \textnormal{ for all } \alpha,\beta \in \Gamma(E)\]
defines a Lie bracket on $E$. The bracket is obviously $\mathbb{R}$-bilinear and skew-symmetric, so it remains to prove the Jacobi identity. We will do this by proving the identity $\Lie_{\anchor}([\alpha,\beta]_E) = [\Lie_{\anchor}(\alpha),\Lie_{\anchor}(\beta)]$
%=\Lie_{\anchor}(\alpha)\Lie_{\anchor}(\beta) - \Lie_{\anchor}(\beta)\Lie_{\anchor}(\alpha)$ 
without using the Jacobi identity, so that the Jacobi identity follows by tracing the proof of Lemma \ref{lemm: Lie derivative of bracket}. 

First of all, by induction, our definition of the bracket gives rise to the formula for $\Lie_{\anchor}(\alpha): \Omega^k(E) \ra \Omega^k(E)$ we had before for Lie algebroids, so, in particular, we see that  
\[\iota_{[\alpha,\beta]_E} = \Lie_{\anchor}(\alpha) \circ \iota_\beta - \iota_\beta \circ \Lie_{\anchor}(\alpha)\]
by the same proof as in Lemma \ref{lemm: interior product of bracket}. Now,
\begin{align*}
    \Lie_{\anchor}(\alpha) \circ d_E &= d_E \circ \iota_\alpha \circ d_E = d_E \circ \Lie_{\anchor}(\alpha),
\end{align*}
by using $d_E^2=0$ and the definition of $\Lie_{\anchor}(\alpha)$, and so,
\begin{align*}
    \Lie_{\anchor}([\alpha,\beta]_E) &= d_E \circ \iota_{[\alpha,\beta]_E} + \iota_{[\alpha,\beta]_E} \circ d_E \\
    &= d_E \circ (\Lie_{\anchor}(\alpha) \circ \iota_\beta - \iota_\beta \circ \Lie_{\anchor}(\alpha)) + (\Lie_{\anchor}(\alpha) \circ \iota_\beta - \iota_\beta \circ \Lie_{\anchor}(\alpha)) \circ d_E \\
    &= \Lie_{\anchor}(\alpha) \circ d_E \circ \iota_\beta - d_E \circ \iota_\beta \circ \Lie_{\anchor}(\alpha) + \Lie_{\anchor}(\alpha) \circ \iota_\beta\circ d_E - \iota_\beta \circ d_E \circ\Lie_{\anchor}(\alpha) \\
    &= \Lie_{\anchor}(\alpha) \circ (d_E \circ \iota_\beta + \iota_\beta\circ d_E) - (d_E \circ \iota_\beta  + \iota_\beta \circ d_E) \circ\Lie_{\anchor}(\alpha) = [\Lie_{\anchor}(\alpha),\Lie_{\anchor}(\beta)].
\end{align*}
To prove that $E \ra M$ is a Lie algebroid with the bracket $[\cdot,\cdot]_E$ and anchor map $\anchor$, it remains to prove the Leibniz rule:
\begin{align*}
    \iota_{[\alpha,f\beta]_E} &= \Lie_{\anchor}(\alpha)\circ\iota_{f\beta} - \iota_{f\beta} \circ \Lie_{\anchor}(\alpha) \\
    &= \Lie_{\anchor}(\alpha)\circ(f\iota_\beta) - f\iota_\beta \circ \Lie_{\anchor}(\alpha) \\
    &= f\Lie_{\anchor}(\alpha)\circ\iota_\beta + \Lie_{\anchor}(\alpha)(f)\iota_\beta - f\iota_\beta \circ \Lie_{\anchor}(\alpha) = \iota_{f[\alpha,\beta]_E} + \Lie_{\anchor}(\alpha)(f)\iota_\beta,
\end{align*}
where we used that $\Lie_{\anchor}(\alpha)$ is a derivation in the third equality. By the Cartan Calculus identities that hold, it is readily verified by induction that $d_E$ equals the de Rham-differential. Again by the Cartan Calculus identities, the Lie bracket is uniquely determined by the anchor map, and the anchor map is uniquely determined by the differential, so the Lie algebroid structure is uniquely determined. This proves the statement.
\end{proof}
In fact, a morphism of vector bundles $E \ra F$ between two vector bundles $E \ra \base$ and $F \ra \subbase$ induces a graded algebra morphism $\Omega^\bullet(F) \ra \Omega^\bullet(E)$ 
%$\omega \in \Omega^k(F) \mapsto x \mapsto \omega(f(x))$
in the obvious way. If $E$ and $F$ are both equipped with a differential, i.e. are Lie algebroids, then there is an obvious candidate for calling such a vector bundle morphism a Lie algebroid morphism; namely, by requiring it to be a cochain map.
%\begin{defn}\label{defn: morphism of Lie algebroids using differential}
%A morphism between two Lie algebroids $\algebroid$ and $\subalgebroid$ is a bundle map $\algebr \ra \subalgebr$ such that the induced map $\Omega^\bullet(\subalgebr) \ra \Omega^\bullet(\algebr)$ is a cochain map.
%\end{defn}
That is, we manually extend the statement in Proposition \ref{prop: differential determines Lie algebroid} to mean an isomorphism of categories. In practice it is useful to have a more concrete definition of a morphism of Lie algebroids, but to write this down it will be necessary to introduce some specific examples. So, before doing this, we will go into examples of Lie algebroids, starting with the Lie algebroid associated to a Lie groupoid. Later, we will show that both definitions of morphisms of Lie algebroids coincide. 

\subsection{The Lie algebroid of a Lie groupoid}\label{sec: The Lie algebroid of a Lie groupoid}

Let $\groupoid$ be a groupoid. As $\ker d\source \subset T\group$ has constant rank (since $\source$ is a submersion), it is a vector bundle over $\group$. We will see that the vector bundle $\algebr \coloneqq \ker d\source|_{\identity(\base)}$ (which is a vector bundle over $\base$) can be made into a Lie algebroid. Compare this discussion with the usual way a Lie algebra is constructed out of a Lie group.

First of all, we set $\anchor \coloneqq d\target|_A: \algebr \ra T\base$. Now, we need to define a Lie bracket $\bracket{\cdot,\cdot}$ on $\algebr$ such that, together with this anchor map, the Leibniz rule
\[\bracket{\alpha,f\beta} = f\bracket{\alpha,\beta} + \anchor(\alpha)(f)\beta, \textnormal{ where } \alpha,\beta \in \Gamma(\algebr) \textnormal{ and } f \in C^\infty(\base),\]
is satisfied. To do this, consider the right-invariant vector fields on $\group$:
\[\mathfrak{X}_{R\textnormal{-inv}}(\group) \coloneqq \{\sigma \in \Gamma(\ker d\source) \mid \sigma(gh) = dR_h(g)\sigma(g) \textnormal{ for all } (g,h) \in \group^{(2)}\},\]
or, using the identification of vector fields and $C^\infty(\group)$-derivations, we have
\[\mathfrak{X}_{R\textnormal{-inv}}(\group) = \{\sigma \in \Gamma(\ker d\source) \mid \sigma(f) \circ R_h = \sigma(f \circ R_h) \textnormal{ for all } h \in \group \textnormal{ and } f \in C^\infty(\base)\},\]
where, for all $h \in \group$, $R_h: \source^{-1}(\target(h)) \ra \source^{-1}(\source(h))$ denotes right-multiplication by $h$ (note: if $\group$ is non-Hausdorff, then the identification with $C^\infty(\group)$-derivations fails; this identification is not necessary, but we will still use the identification to ease the computations a bit). Notice that, indeed, for all $g \in \source^{-1}(\target(h))$, $dR_h(g)$ maps $\ker d\source(g)$ into $\ker d\source(gh)$: this follows by the observation that $\source \circ R_h$ maps into $\{\source(h)\} \subset \base$, and a straightforward application of the chain rule. 

\begin{prop}\cite{ruimarius}\label{prop: sections of Lie algebroid of a Lie groupoid are in bijective correspondence with right invariant vector fields}
Let $\groupoid$ be a groupoid and set $\algebr = \ker d\source|_{\identity(\base)}$. Then the right invariant vector fields $\mathfrak{X}_{R\textnormal{-inv}}(\group)$ are in bijective correspondence with $\Gamma(\algebr)$. In particular, we obtain a Lie bracket $\bracket{\cdot,\cdot}$ on $\algebr$.
\end{prop}
\begin{proof}
Let $\alpha \in \Gamma(\algebr)$. Define 
\[\sigma_\alpha \in \Gamma(\ker d\source) \textnormal{ by } g \mapsto dR_g(\id_{\target(g)})\alpha(\target(g)).\]
Then $\sigma_\alpha$ is right-invariant: for all $(g,h) \in \group^{(2)}$ we have
\begin{align*}
    \sigma_\alpha(gh) &= dR_{gh}(\id_{\target(gh)})\alpha(\target(gh)) = dR_h(g) \circ dR_g(\id_{\target(g)})\alpha(\target(g)) = dR_h(g)\sigma_\alpha(g),
\end{align*}
where in the second equality we used that $R_{gh} = R_h \circ R_g$ and that $\target(gh) = \target(g)$.
Conversely, for $\sigma \in \mathfrak{X}_{R\textnormal{-inv}}(\group)$ we define 
\[\alpha_\sigma \in \Gamma(A) \textnormal{ by } x \mapsto \sigma(\id_x).\]
Then the above assignments are readily verified to be inverse to each other. This proves the first statement.

For the last statement, note that for all $\sigma,\tau \in \mathfrak{X}_{R\textnormal{-inv}}(\group)$, all $f \in C^\infty(\group)$ and all $(g,h) \in \group^{(2)}$ we have
\begin{align*}
\sigma(\tau f)(gh) &= d(\tau f)(gh)\sigma(gh) \\
&= d(\tau f)(gh)dR_h(g)\sigma(g) \\
&=  d(\tau f \circ R_h)(g)\sigma(g) \\
&= d(\tau(f \circ R_h))(g)\sigma(g) = \sigma(\tau(f \circ R_h))(g).
\end{align*}
It follows that $\mathfrak{X}_{R\textnormal{-inv}}(\group)$ is a Lie subalgebra of $\mathfrak{X}(\group)$ (equipped with the usual Lie bracket of vector fields). Indeed, for all $\sigma,\tau \in \mathfrak{X}_{R\textnormal{-inv}}(\mathcal{G})$, all $f \in C^\infty(\group)$, and all $(g,h) \in \group^{(2)}$, we have by the above calculation that
\begin{align*}
    [\sigma,\tau](f) \circ R_h(g) &= \sigma(\tau f)(gh) - \tau(\sigma f)(gh) \\
    &= \sigma(\tau(f \circ R_h))(g) - \tau(\sigma(f \circ R_h))(g) = [\sigma,\tau](f \circ R_h)(g).
\end{align*}
This equips $\mathfrak{X}_{R\textnormal{-inv}}(\group)$ with the usual Lie bracket $[\cdot,\cdot]$ of $\mathfrak{X}(\group)$ and so we can transfer this Lie bracket of $\mathfrak{X}_{R\textnormal{-inv}}(\group)$ to $\Gamma(\algebr)$ by using the bijection constructed above.
\end{proof}
\begin{prop}\cite{ruimarius}\label{prop: Lie bracket of Lie algebroid of a Lie groupoid satisfies the Leibniz identity}
Let $\groupoid$ be a Lie groupoid and consider $\algebr \coloneqq \ker d\source|_{\identity(\base)}$ together with the Lie bracket $\bracket{\cdot,\cdot}$ induced by $\mathfrak{X}_{R\textnormal{-inv}}(\group)$ (see Proposition \ref{prop: sections of Lie algebroid of a Lie groupoid are in bijective correspondence with right invariant vector fields}) and anchor map $\anchor = d\target|_\algebr$. Then the Leibniz rule of Definition \ref{defn: Lie algebroid} is satisfied.
\end{prop}
\begin{proof}
Let $\alpha,\beta \in \Gamma(\algebr)$ and $f \in C^\infty(\base)$. First of all, for all $g \in \mathcal{G}$, with notation as in the proof of Proposition \ref{prop: sections of Lie algebroid of a Lie groupoid are in bijective correspondence with right invariant vector fields}, we have
\[\sigma_{f\beta}(g) = dR_g(\id_{\target(g)})(f\beta)(\target(g)) = f(\target(g))dR_g(\id_{\target(g)})\beta(\target(g)) = f\circ \target(g)\sigma_\beta(g),\]
so $\sigma_{f\beta} = (f \circ \target)\sigma_\beta$. Moreover, note that
\begin{align*}
    \sigma_\alpha(f \circ \target)(g) &= d(f \circ \target)(g)\sigma_\alpha(g) \\
    &= df(\target(g)) \circ d\target(g)dR_g(\id_{\target(g)})\alpha(\target(g)) \\
    &= df(\target(g)) \circ d\target(\id_{\target(g)})\alpha(\target(g)) = \anchor(\alpha)(f) \circ \target(g),
\end{align*}
where in the third equality we used that $\target \circ R_g = \target$. It follows that
\begin{align*}
    [\sigma_\alpha,\sigma_{f\beta}] &= [\sigma_\alpha,(f \circ \target)\sigma_\beta] \\
    &= (f \circ \target)[\sigma_\alpha,\sigma_\beta] + \sigma_\alpha(f \circ \target)\sigma_\beta \\
    &= (f \circ \target)[\sigma_\alpha,\sigma_\beta] + (\anchor(\alpha)(f) \circ \target)\sigma_\beta,
\end{align*}
where in the second equality we used the Leibniz rule for $[\cdot,\cdot]$. For all $x \in \algebr$ we now have
\begin{align*}
    \bracket{\alpha,f\beta}(x) &= [\sigma_\alpha,\sigma_{f\beta}](\id_x) \\
    &= (f \circ \target)(\id_x)[\sigma_\alpha,\sigma_\beta](\id_x) + (\anchor(\alpha)(f) \circ \target)(\id_x)\sigma_\beta(\id_x) \\
    &= f(x)[\alpha,\beta](x) + \anchor(\alpha)(f)(x)\beta(x),
\end{align*}
which proves the result.
\end{proof}
\begin{defn}\cite{ruimarius}\label{term: Lie algebroid of a Lie groupoid}
Let $\groupoid$ be a groupoid. The Lie algebroid $\algebr \coloneqq \ker d\source|_{\identity(\base)}$, together with the Lie bracket $\bracket{\cdot,\cdot}$ induced by $\mathfrak{X}_{R\textnormal{-inv}}(\group)$ and anchor map $\anchor \coloneqq d\target|_\algebr$, is called the Lie algebroid of the groupoid $\groupoid$.
\end{defn}

In the case that the base is a point, we recover the usual Lie algebra associated to a Lie group. Conversely, we might wonder if every Lie algebroid integrates to a Lie groupoid, as is the case in the Lie algebra setting. However, this is not the case in general. 

\begin{defn}\cite{ruimarius}\label{defn: integrable Lie algebroid}
A Lie algebroid is called \textit{integrable} if it is isomorphic to the Lie algebroid of a Lie groupoid.
\end{defn}
In \cite{ruimarius} it is explained in detail when a Lie algebroid is integrable, and when it is not. In the case of Lie algebras, we can integrate a Lie algebra to a unique connected, simply connected Lie group. This result extends to integrable (!) Lie algebroids as follows. 
\begin{defn}\cite{ruimarius}\label{defn: s-(simply )connected groupoid}
A groupoid is called \textit{$\source$-(simply) connected} (read: source-(simply) connected) if all source fibers are (simply) connected.
\end{defn}
Any groupoid $\groupoid$ has an open $\source$-connected subgroupoid $\group^0 \rra \base$. This is the content of the following lemma. 
\begin{lemm}\cite{ruimarius}\label{lemm: open s-connected subgroupoid}
Let $\groupoid$ be a groupoid and denote \[\group^0 \coloneqq \{g \in \group \mid g  \in \source^{-1}(\source(g))^0\} = \bigcup_{x \in \base} \source^{-1}(x)^0,\]
where $\source^{-1}(x)^0$ denotes the connected component of $\source^{-1}(x)$ containing $\id_x$. Then $\group^0 \subset \group$ is open, and the restriction to $\group^0$ defines an $\source$-connected Lie groupoid $\group^0 \rra \base$ that has the same Lie algebroid as $\groupoid$.
\end{lemm}
\begin{proof}
We first prove that $\group^0 \subset \group$ is open. Since $\source$ is a submersion, it induces a foliation $\mathcal{F}_\source$ of $\group$ whose leaves are the connected components of the source fibers, say of codimension $q$. Let $x \in \base$ and let $g \in \source^{-1}(x)^0$. Pick a foliation chart $(U,\Phi)$ around $\id_x$. Since $T=\identity(\base)$ is transverse to the foliation, the subset $T_{\textnormal{leaf}}$ from Proposition \ref{prop: the collection of leaves intersecting an open subset} is open. By definition of $\group^0$, $T_{\textnormal{leaf}} \subset \group^0$. The element $g$ lies in this subset, and by Proposition \ref{prop: the collection of leaves intersecting an open subset}, this set is open. To see that $\group^0 \rra \base$ is a groupoid, it remains to show that $\mult|_{(\group^0)^{(2)}}$ and $\inv|_{\group^0}$ map into $\group^0$. Indeed, for all $g \in \group^0$, the diffeomorphism $R_g: \source^{-1}(\target(g)) \ra \source^{-1}(\source(g))$ maps $\source^{-1}(\target(g))^0$ bijectively onto $\source^{-1}(\source(g))^0$. The first claim follows immediately, and the second claim follows by the observation that $g^{-1}$ maps to $\id_{\source(g)}$, so that $g^{-1} \in \source^{-1}(\target(g))^0$. The source fibers of $\group^0$ are given by $\source^{-1}(x)^0$, so $\group^0 \rra \base$ is an $\source$-connected Lie groupoid. The last statement follows by construction of the Lie algebroid of a Lie groupoid (alternatively, it follows from Proposition \ref{prop: Lie groupoid morphism integration source-simply connected}, which can only be discussed later).
\end{proof}
From this we can deduce the following partial analogue of Lie's third theorem.
\begin{prop}\cite{ruimarius}\label{prop: s-(simply )connected groupoid of integrable algebroid}
Let $\algebroid$ be an integrable Lie algebroid, integrating to, say, $\groupoid$. Then there is a unique (up to isomorphism) $\source$-connected, $\source$-simply connected groupoid integrating $\algebroid$. Moreover, it is locally diffeomorphic to $\group^0 \rra \base$.
\end{prop}
\begin{proof}
Since $\group^0 \rra \base$ is also an integration of $\algebroid$, we may assume $\groupoid$ is $\source$-connected. 

The foliation $\mathcal{F}_\source$, by source fibers of $\group$, induces a set-theoretical groupoid
\[\Mon(\group,\mathcal{F}_\source) \rra \group.\]
Apart from Hausdorffness assumptions, this groupoid is a Lie groupoid. Now consider the canonical right principal $\group$-action on $\group$ (with moment map $\target$). Then, by inspecting the definition, it is obvious that the $\group$-action maps a leaf of $\mathcal{F}_\source$ to another leaf of $\mathcal{F}_\source$. Hence, the $\group$-action 
\[\group \tensor[_{\source}]{\times}{_{\target_{\Mon}}} \Mon(\group,\mathcal{F}_{\source}) \ra \Mon(\group,\mathcal{F}_\source) \textnormal{ given by } (g,[\gamma]) \mapsto [R_g \circ \gamma]\]
is well-defined. This is again a right principal $\group$-action, so we can form the quotient groupoid $\widetilde{\group} \coloneqq \Mon(\group,\mathcal{F}_\source)/\group$ over $\base$ (with our conventions on Hausdorffness, this really is a Lie groupoid). Now, any monodromy groupoid is $\source$-connected and $\source$-simply connected, and $\widetilde{\group} \rra \base$ inherits these properties. The Lie groupoid morphism
%Since, for all $x \in \base$, $\source^{-1}(x)$ is connected, we have an associated universal cover of $\source_{\group^0}^{-1}(x)$, which can be realised as the homotopy classes of paths (relative to base-point) starting at $\id_x$. These universal covers assemble (as a disjoint union) into a groupoid $\widetilde{\group} \rra \base$ with $\source_{\widetilde{\group}}([\gamma]) = \gamma(0)$, $\target_{\widetilde{\group}}([\gamma]) = \gamma(1)$ and $\mult_{\widetilde{\group}}([\gamma],[\gamma']) = [(R_{\gamma'(1)} \circ \gamma) \star \gamma']$, where $R_g$ denotes the right-multiplication of $g \in \group$ in $\group$. Another way of defining $\widetilde{\group}$, which makes its smooth structure apparent, is as follows: since $\source$ is a submersion, it induces a foliation $\mathcal{F}$ on $\group$ whose leaves are the source fibers. Then we can realise $\mathcal{G}$ as the inverse image of $\base$ along the source map of the monodromy groupoid $\Mon(\group,\mathcal{F}) \rra \group$ (note: this is not a groupoid using our conventions because $\group$ is not necessarily Hausdorff, but all other groupoid axioms are satisfied, and the base of $\mathcal{G}$ is $\base$ which is Hausdorff nonetheless). The obvious projection 
\[\Mon(\group,\mathcal{F}_\source) \ra \group \textnormal{, given by } [\gamma] \mapsto \gamma(1) \cdot \gamma(0)^{-1},\]
with base map $\target$, is $\group$-invariant, and so it descends to a Lie groupoid morphism $\widetilde{\group} \ra \group$. It is surjective (since $\group$ is $\source$-connected), and it is readily verified, by using the charts constructed for the monodromy groupoid, that it is a local diffeomorphism. The proof that the Lie algebroids are equal, and the uniqueness statement, are a consequence of Proposition \ref{prop: Lie groupoid morphism integration source-simply connected} (which can only be discussed later).
\end{proof}

We end this section with two remarks. 

\begin{rema}\cite{meinrenken}\label{rema: other definitions of Lie algebroid of a Lie groupoid}
Notice that, by looking at left-invariant vector fields instead, we could just as well have defined the above Lie algebroid as $\algebr \coloneqq \ker d\target|_{\id(\base)}$, with bracket $\bracket{\cdot,\cdot}$ induced by the left-invariant vector fields, and anchor map $\anchor \coloneqq d\source|_{\algebr}$. Both Lie algebroids are canonically isomorphic. There is even another way to define a Lie algebroid from a Lie groupoid by defining the Lie algebroid as the normal bundle of $\id(X)$ in $\group$. Normal bundles are important in the theory of blow-ups, and this point of view is very useful, but more on this treatment later. 
\end{rema}
\begin{rema}\label{rema: invariant forms}
Consider $\Omega^\bullet(\algebr)$ of the Lie algebroid $\algebroid$ of a Lie groupoid $\groupoid$. An appropriate notion of right-invariant forms on $\groupoid$ yields a bijection between these forms and $\omega^\bullet(\algebr)$. Indeed, we call a form $\omega \in \Omega^k(\ker d\source)$ \textit{right-invariant} if, for all $(g,h) \in \group^{(2)}$ and $\xi_1,\dots,\xi_k \in \ker d\source(g)$, we have 
\[\omega(gh)(dR_h(g)\xi_1,\dots,dR_h(g)\xi_k)=\omega(g)(\xi_1,\dots,\xi_k).\]
Then, similar to before, the map $\omega \mapsto (x \mapsto \omega(\id_x))$ defines a bijection $\Omega^\bullet_{R\textnormal{-inv}} \xra{\sim} (\group)\Omega^\bullet(\algebr)$.
\end{rema}

\subsection{Examples of Lie algebroids}\label{sec: Examples of Lie algebroids}
Just as Lie groupoids are intimately related to Lie groups (e.g. Lie groups are given by Lie groupoids over a point, principal bundles are given by transitive Lie groupoids and Lie group actions are action Lie groupoids), Lie algebroids are intimately related to Lie algebras in a similar way. We also show, whenever we can, what Lie groupoids the following Lie algebroids integrate to. We start with the Lie algebroid analogues of the Lie groupoids defined in Section \ref{sec: Examples of Lie groupoids}.
\begin{exam}\cite{ruimarius}\label{exam: trivial Lie algebroid}
The tangent bundle $TM$ of a manifold $M$, with anchor map $\anchor=\id_{TM}$ and Lie bracket the usual Lie bracket $[\cdot,\cdot]$ of vector fields, is a Lie algebroid.
\end{exam}
\begin{rema}\cite{ruimarius}\label{rema: trivial Lie algebroid and pair groupoid}
The tangent bundle $TM$ integrates to the pair groupoid $M \times M \rra M$. Indeed, $\source=\pr_2$ and $\target=\pr_1$, so we see that
\[d\source|_{\Delta(M)}: T(M \times M)|_{\Delta M} \ra TM \textnormal{, is given by } (x,x,\xi_1,\xi_2) \mapsto (x,\xi_2),\]
which shows that $\ker d\source|_{\Delta(M)} \cong TM$, and, under this identification, we see that $d\target|_{TM}=\id_{TM}$. The Lie bracket of $A$ is induced by the right-invariant vector fields of $M \times M \rra M$. From the proof of Proposition \ref{prop: sections of Lie algebroid of a Lie groupoid are in bijective correspondence with right invariant vector fields}, it is readily verified that the induced Lie bracket, in terms of derivations, coincides with the usual Lie bracket on $\mathfrak{X}(M)=\Gamma(A)$.
\end{rema}
By viewing $M$ as a principal $\{e\}$-bundle, we can generalise this trivial example to the following class of examples.
\begin{exam}[Atiyah Lie algebroids]\cite{ruimarius}\label{exam: Atiyah Lie algebroids}
To any principal $G$-bundle $\pi: P \ra M$ we can associate the \textit{Atiyah Lie algebroid} 
\[A(P): TP/G \ra M,\]
($G$ acts on $TP$ via $(p,\xi) \cdot g \coloneqq (pg,dR_g(p)\xi)$) with anchor map $\anchor: TP/G \ra TM$, which is the map induced by the $G$-equivariant map $d\pi$ (with trivial $G$-action on $TM$). The Lie bracket is defined as the bracket of
\[\mathfrak{X}_G(P) \coloneqq \{\sigma \in \mathfrak{X}(P) \mid \sigma(pg)=\sigma(p) \textnormal{ for all } g \in G\},\]
(which is a Lie subalgebra of $\mathfrak{X}(P)$), under the bijective correspondence with sections of $A(P)$, which are given by the induced sections from $\mathfrak{X}_G(P)$.
\end{exam}
\begin{rema}[Lie algebroid of Gauge groupoid is Atiyah algebroid]\label{rema: Lie(Gauge groupoid) = Atiyah algebroid}
Let $P \ra M$ be a principal $G$-bundle. We will show that the Lie algebroid $A$ of $P \otimes_G P$ is the Atiyah Lie algebroid $A(P)$. To see this, notice that the source map $\source: (P \times P)/G \ra M$ is induced by the projection map $\pr_2: P \times P \ra P$. Since $d\pr_2=\pr_2: TP \times TP \ra TP$, $d\source$ is the induced map 
\[(TP \times TP)/TG \ra TP/TG \cong T(P/G) \cong TM,\]
with $TG$-action on $TP$ given by the differential of the action map (and $TG$ acts diagonally on $TP$). Notice that $(TP \times TP)/TG|_{\identity(M)}$ is given by $\{[p,\xi_1,p,\xi_2] \in (TP \times TP)/TG\}$, so
\[A = \ker d\source|_{\identity(M)} = \{[p,\xi,p,0] \in (TP \times TP)/TG\}.\]
Denote by $R: P \times G \ra P$ the action of $G$ on $P$. Then, for all $(p,\xi) \in TP$ and $(g,\eta) \in TG$, $dR(p,g)(\xi,\eta) = dR_g(p)\xi + d\hat{p}(g)\eta$ (where $\hat{p}: G \ra P$ is the map $h \mapsto ph$). In particular,
\[(p,\xi,p,0) \cdot (g,\eta) = (pg,dR_g(p)\xi + d\hat{p}(g)\eta,pg,d\hat{p}(g)\eta),\]
which shows that the smooth vector bundle morphism $TP \times TP \ra TP$ given by $(p_1,\xi_1,p_2,\xi_2) \mapsto (p_1,\xi_1 - \xi_2)$ descends to a vector bundle isomorphism $\varphi: A \ra A(P)$ (it has an obvious inverse). The base map of $\varphi$ is the isomorphism  $TP/TG \cong T(P/G)$ given by
\[[p,\xi] \mapsto ([p],d\pi(p)\xi).\]
To see that the map commutes with the anchor maps, notice that on $A$ the anchor map is given by 
\[[p,\xi,p,0] \mapsto [p,\xi] \in TP/TG \cong TM\]
and on $A(P)$ the anchor map is given by 
\[[p,\xi] \mapsto ([p],d\pi(p)\xi) \in T(P/G) \cong TM,\]
where $\pi: P \ra P/G \cong M$ is the projection map, which shows that, indeed, $\varphi$ commutes with the anchor maps (i.e. is an anchored vector bundle isomorphism). It remains to show that $\varphi$ intertwines the Lie brackets. However, this is obvious once we recognise that the Lie bracket on $A(P)$ is induced by the bijective correspondence between $\Gamma(A(P))$ and the Lie subalgebra $\mathfrak{X}_G(P) \subset \mathfrak{X}(P)$, and, similarly, the Lie bracket on $A$ is induced by the bijective correspondence between $\Gamma(A)$ and the Lie subalgebra $\mathfrak{X}_{R\textnormal{-inv}}((P \times P)/G) \subset \mathfrak{X}((P \times P)/G)$. Under these bijective correspondences, the vector bundle isomorphism $\varphi$ maps a section $\sigma=(\sigma^1,\sigma^2) \in \mathfrak{X}_{R\textnormal{-inv}}((P \times P)/G)$ to the section $\sigma^1-\sigma^2 \in \mathfrak{X}_G(P)$ (given by $p \mapsto \sigma^1([p,p]) - \sigma^2([p,p])$), and so 
\[[\sigma,\tau] = ([\sigma^1,\tau^1],[\sigma^2,\tau^2]) \textnormal{ maps, under $\varphi$, to } [\sigma^1 - \sigma^2, \tau^1 - \tau^2],\]
where we used that $\varphi$ induces a bijective map $\mathfrak{X}_{R\textnormal{-inv}}(P \times P/G) \ra \mathfrak{X}_G(P)$, so that, for all $\widetilde\sigma \in \mathfrak{X}_{R\textnormal{-inv}}(P \times P/G)$, we have $\widetilde\sigma=(\widetilde\sigma^1-\widetilde\sigma^2,0)$ in $\mathfrak{X}_{R\textnormal{-inv}}(P \times P/G)$. This shows that, indeed, $A(P)$ can be seen as the Lie algebroid of $P \otimes_G P$.
\end{rema}
Recall that principal bundles are in bijective correspondence with transitive groupoids (even via an equivalence of categories; see Remark \ref{rema: transitive groupoids are in one to one correspondece with Gauge groupoids}). In the setting of Lie algebroids there is also an appropriate notion of transitivity, and Atiyah Lie algebroids belong to this class of Lie algebroids.
\begin{defn}\cite{ruimarius}\label{defn: transitive Lie algebroid}
A Lie algebroid is called \textit{transitive} if the anchor map is surjective.
\end{defn}
Later, in Remark \ref{rema: transitive Lie algebroids and Atiyah algebroids} we will give, explicitly, a class of transitive Lie algebroids which do not integrate to a transitive Lie groupoid; in fact, these Lie algebroids do not integrate at all to a Lie groupoid. 
\begin{rema}\cite{ruimarius}\label{rema: Atiyah Lie algebroids are transitive}
Let $\pi: P \ra M$ be a principal $G$-bundle. Then $A(P)$ is transitive, because the anchor map $\anchor_{A(P)}$ is the induced map from the surjective $G$-equivariant map $d\pi: TP \ra TM$ (with trivial $G$-action on $TM$). 
\end{rema}
For the infinitesimal counterpart of general linear groupoids, we introduce the following definition.
\begin{defn}\cite{meinrenken}\label{defn: derivations of a vector bundle}
Let $E \ra M$ be a vector bundle. Then a \textit{derivation} of $E$ is an $\mathbb{R}$-linear map $D: \Gamma(E) \ra \Gamma(E)$ together with a (uniquely determined) vector field $\sigma_D \in \mathfrak{X}(M)$ such that the Leibniz-type rule
\[D(f\alpha) = fD(\alpha) + \sigma_D(f)\alpha, \textnormal{ for all } f \in C^\infty(M) \textnormal{ and } \alpha \in \Gamma(E),\]
is satisfied.
\end{defn}
\begin{exam}[General linear algebroids]\cite{meinrenken}\label{exam: general linear algebroid}
Let $E \ra M$ be a vector bundle. Then there is an associated \textit{general linear algebroid} $\mathfrak{gl}(E) \ra M$, where $\mathfrak{gl}(E)$ consists of derivations $D$ of $E$ (which comes with a uniquely determined vector field $\sigma_D \in \mathfrak{X}(M)$), with anchor map 
\[\anchor: \mathfrak{gl}(E) \ra TM \textnormal{ given by } D \mapsto \sigma_D.\]
The Lie bracket is given by 
\[[D,D']_{\mathfrak{gl}(E)} \coloneqq DD' - D'D\]
which is readily verified to be a derivation again (note: $\sigma_{[D,D']_{\mathfrak{gl}(E)}} = [\sigma_D,\sigma_{D'}]$).
\end{exam}
\begin{rema}\cite{meinrenken}\label{rema: general linear algebroid is atiyah algebroid frame bundle}
Alternatively, $\mathfrak{gl}(E)$ of a vector bundle $E \ra M$ is the Atiyah algebroid associated to the principal $\textnormal{GL}(k,\mathbb{R})$-bundle $\textnormal{Fr } E$, i.e. the Lie algebroid of the Lie groupoid $\textnormal{GL}(E)$. To establish this correspondence, notice that a section $\alpha$ of the Lie algebroid $\algebr$ of a groupoid $\group$ (in this case $A(\textnormal{Fr } E)$) induces a flow 
\[\phi_\alpha^t(x) \coloneqq \phi^t_{\sigma_\alpha}(\id_x) \in \group \textnormal{ } (x \in M)\]
(see the proof of Proposition \ref{prop: sections of Lie algebroid of a Lie groupoid are in bijective correspondence with right invariant vector fields} for the definition of $\sigma_\alpha \in \mathfrak{X}_{R\textnormal{-inv}}(\group)$), which, in this case, associates to $\alpha \in \Gamma(\algebr)$, where $\algebr\coloneqq A(\textnormal{Fr } E)$, the linear maps 
\[\phi_\alpha^t(x): E_x \ra E_{\phi_{\anchor_\algebr(\alpha)}^t(x)},\]
of course, only whenever defined. If we pick $t$ small enough, we can make sure this map exists and is invertible (with inverse $\phi_\alpha^{-t}(x)$), so we obtain a derivation $D_\alpha: \Gamma(E) \ra \Gamma(E)$ by setting
\[D_\alpha(s) \coloneqq \left.\frac{d}{dt}\right\vert_{t=0} (\phi^t_\alpha)^*s = (x \mapsto \left.\frac{d}{dt}\right\vert_{t=0} \phi^{-t}_\alpha(x) \circ s \circ \phi^t_{\anchor_\algebr(\alpha)}(x)),\]
with $\sigma_{D_\alpha} = \anchor_\algebr(\alpha)$, by an application of the chain rule. 

Conversely, if we start with a derivation $D$ of $E$, then the formula $D(s) = \left.\tfrac{d}{dt}\right\vert_{t=0} (\phi^t_D)^*s$ defines an ODE from which we obtain automorphisms $\phi^t_D$ of $E$ (for small enough $t$). Therefore, we can set
\[\alpha_D(x) \coloneqq \left.\frac{d}{dt}\right\vert_{t=0} \phi^t_D(x),\] 
and this defines a section of $\algebr$. It is readily verified that $\alpha_{D_{\alpha'}}=\alpha'$ for all $\alpha' \in \Gamma(\algebr)$, and that $D_{\alpha_{D'}}$ for all derivations $D'$ of $E$.

To show that the Lie brackets coincide under the correspondence, notice that, by the local nature, we can assume $E$ is trivial, and in this case it follows immediately by definition of the Lie brackets in question.
\end{rema}
\begin{exam}[Action algebroids]\cite{ruimarius}\label{exam: action Lie algebroid}
Suppose we have a Lie algebra $\mathfrak{g}$ that acts on a manifold $M$ by a Lie algebra morphism $\rho: \mathfrak{g} \ra \mathfrak{X}(M)$. The \textit{action Lie algebroid} $\mathfrak{g} \ltimes M$ associated to $\rho$ is the trivial vector bundle $\mathfrak{g} \times M \ra M$ with anchor map given by the action map:
\[\anchor: \mathfrak{g} \times M \ra TM \textnormal{ given by } (\xi,x) \mapsto (x,\rho(\xi)(x)).\]
To define a bracket $[\cdot,\cdot]_{\mathfrak{g} \times M}$, note that $\Gamma(\mathfrak{g} \times M) \cong C^\infty(M,\mathfrak{g})$, so we can set 
\[[\alpha,\beta]_{\mathfrak{g} \times M} \coloneqq [\alpha,\beta]_{\mathfrak{g}} + \anchor(\alpha)(\beta) - \anchor(\beta)(\alpha) \textnormal{ for all } \alpha,\beta \in  C^\infty(M,\mathfrak{g}),\]
where, for $\alpha,\beta \in C^\infty(M,\mathfrak{g})$ and $x \in M$, $\anchor(\alpha)(\beta)(x) \coloneqq d\beta(x)\anchor(\alpha(x),x)$, and the expression should be interpreted as first evaluating the sections before doing the operations. Since $[\cdot,\cdot]_{\mathfrak{g}}$ and $\rho$ are $\mathbb{R}$-bilinear, and, for all $\alpha,\beta \in C^\infty(M,\mathfrak{g})$ and $x \in M$,
\begin{align*}
    \anchor(\alpha)(f\beta)(x) &= d(f\beta)(x)\anchor(\alpha(x),x)\\
    &=f(x)d\beta(x)\anchor(\alpha(x),x)+df(x)\anchor(\alpha(x),x)\beta(x) \\
    &= \left(f\anchor(\alpha)(\beta) + \anchor(\alpha)(f)\beta\right)(x),
\end{align*}
we see that the Leibniz rule holds. The skew-symmetry and the Jacobi identity are immediate consequences of the definition. Notice that for constant sections $\alpha,\beta \in C^\infty(M,\mathfrak{g})$ we have $[\alpha,\beta]_{\mathfrak{g} \times M} = [\alpha,\beta]_{\mathfrak{g}}$. By the Leibniz rule, $[\cdot,\cdot]_{\mathfrak{g} \times M}$ is uniquely determined by this condition. 
\end{exam}
\begin{rema}\cite{ruimarius}\label{rema: Lie algebroid of action groupoid is action algebroid}
Let a Lie group $G$ act on a manifold $M$. We will show that the Lie algebroid $A$ of the corresponding action groupoid $G \ltimes M$ can be identified with $\mathfrak{g} \ltimes M$, where $\mathfrak{g}$ acts on $M$ via the induced infinitesimal action:
\[\rho(\xi)(x) = \left.\frac{d}{dt}\right\vert_{t=0} \exp_G(t\xi) \cdot x.\]
First of all, to avoid confusion, this infinitesimal action is an anti-homomorphism if we regard $\mathfrak{g}$ as the left-invariant vector fields on $G$, but it is a homomorphism if we regard $\mathfrak{g}$ as the right-invariant vector fields on $G$ (which corresponds to our treatment of the Lie algebroid of a Lie groupoid). To see that the Lie algebroid of $G \ltimes M$ equals $\mathfrak{g} \ltimes M$ with respect to $\rho$, notice that
\[\ker d\source|_{\identity(M)} = \ker d\pr_2(e,\cdot), \textnormal{ where } \pr_2: G \times M \ra M \textnormal{ is the projection},\]
so, since for all $x \in M$ we have that $d\pr_2(e,x): \mathfrak{g} \times TM \ra TM$ is the projection, we see that $\ker d\pr_2(e,\cdot) = \mathfrak{g} \times M$. The anchor is given by $d\target|_A = \rho|_{\mathfrak{g} \times M}$, and for constant sections $\alpha,\beta \in C^\infty(M,\mathfrak{g})$ we have
\[\bracket{\alpha,\beta}(x) = [X_\alpha,X_\beta]_{\mathfrak{X}_{R-\textnormal{inv}}(G \times M)}(e,x) = [\alpha(x),\beta(x)]_{\mathfrak{g}},\]
which proves that the bracket is equal to the bracket of $\mathfrak{g} \ltimes M$. So, indeed, the Lie algebroid of $G \ltimes M$ is $\mathfrak{g} \ltimes M$.
\end{rema}

We will define a Lie algebroid morphism $A \ra B$ between two Lie algebroids $\algebroid$ and $B \ra \subbase$ by requiring that the graph of this morphism is a Lie subalgebroid of $A \times B \ra \base \times \subbase$. The latter is always, naturally, a Lie algebroid.

\begin{exam}[Product algebroids]\cite{meinrenken}\label{exam: product of Lie algebroids}
Let $A_1 \ra M_1$ and $A_2 \ra M_2$ be two Lie algebroids. The \textit{product Lie algebroid} of $A_1$ and $A_2$ is the
vector bundle $A_1 \times A_2 \ra M_1 \times M_2$ with anchor map given by 
\[\anchor = \anchor_{A_1} \times \anchor_{A_2}: A_1 \times A_2 \ra TM_1 \times TM_2 \cong T(M_1 \times M_2).\] 
For convenience, we fix $a\in\{1,2\}$. To define the Lie bracket, let $\{\alpha^i_a\}$ be a local frame for $A_a$, and denote by $\pr_a: M_1 \times M_2 \ra M_a$ the projection. Then the vector bundle $A_1 \times A_2$ equals $\pr_1^*A_1 \oplus \pr_2^*A_2$, where $\oplus$ denotes the Whitney sum (over the base $M_1 \times M_2$). With this notation, $\{\gamma^k\} = \{\gamma^{i}_1,\gamma^{j}_2\} \coloneqq \{\pr_1^*\alpha^i_1,\pr_2^*\alpha^j_2\}$ is a local frame for $A_1 \times A_2$. We now define
\[[\gamma^{p}_a,\gamma^{q}_a]_{A_1 \times A_2} \coloneqq \pr_a^*[\alpha^p_a,\alpha^q_a]_{A_a}, \textnormal{ and } [\gamma^{p}_a,\gamma^{m}_b]_{A_1 \times A_2} \coloneqq 0 \textnormal{ where } b \in \{1,2\} \setminus \{a\},\]
so that, in general, for $u = \textstyle\sum_{k}u_k\gamma^k$ and $v = \textstyle\sum_\ell v_\ell\gamma^\ell$, where $u_k,v_\ell \in C^\infty(M_1 \times M_2)$, we have
\begin{align*} 
    [u,v]_{A_1 \times A_2} = \sum_{k,\ell}u_kv_\ell[\gamma^k,\gamma^\ell]_{A_1 \times A_2} + u_k\anchor(\gamma^k)(v_\ell)\gamma^\ell - v_\ell\anchor(\gamma^\ell)(u_k)\gamma^k.
\end{align*}
%\begin{align*}\label{eq: bracket for product algebroid}
%    [u,v]_{A_1 \times A_2} = &\bigodot_{1ik,2j\ell} u^1_iv^1_k[\gamma^i_1,\gamma^k_1]_{A_1 \times A_2} + u_i^1\anchor(\gamma^i_1)(v^1_k)\gamma^k_1 - v_k^1\anchor(\gamma^k_1)(u^1_i)\gamma^i_1 \\
%    + &\bigodot_{i\ell,jk} u_i^1\anchor(\gamma^i_1)(v^2_\ell)\gamma_2^\ell - v_\ell^2\anchor(\gamma^\ell_2)(u_i^1)\gamma^i_1.     
%\end{align*}
%so that, in general, for $u = \textstyle\sum_{i,j}u_{ij}\gamma^{ij}$ and $v = \textstyle\sum_{k,\ell}v_{k\ell}\gamma^{k\ell}$, where $u_{ij},v_{k\ell} \in C^\infty(M_1 \times M_2)$, we have
%\begin{equation}
%    [u,v]_{A_1 \times A_1} = \sum_{i,j,k,\ell} u_{ij}v_{k\ell}[\gamma^{ij},\gamma^{k\ell}]_{A_1 \times A_2} + u_{ij}(\anchor_{A_1 \times A_2}(v_{k\ell})(\gamma^{ij}))\gamma^{k\ell} - v_{k\ell}(\anchor_{A_1 \times A_2}(u_{ij})(\gamma^{k\ell}))\gamma^{ij},
%\end{equation}
by the Leibniz rule. It is immediate that the Lie bracket is skew-symmetric. To show that the Jacobi identity holds, it suffices to prove the identity in case all sections are part of the local frame $\{\gamma^i_1,\gamma^j_2\}$ (see Remark \ref{rema: defining a Lie algebroid locally out of an anchored vector bundle}). Since $[\gamma^i_a,\gamma^j_b]_{A_1 \times A_2}=0$, where $b \in \{1,2\} \setminus \{a\}$, by definition, it suffices to prove that   
\[\bigodot_{p,q,r} [[\gamma^{p}_a,\gamma^{q}_a]_{A_1 \times A_2},\gamma^{r}_a]_{A_1 \times A_2}=0.\] 
We will show that the structure equation of the Jacobi identity (see Remark \ref{rema: Lie algbroid local nature}) is satisfied. We write
\[[\alpha^w_a,\alpha^z_a]_{A_a} = \sum_tc^{wz}_{at}\alpha^t_a \textnormal{, and } \anchor_{A_a}(\alpha^w_a) = \sum_tb^{wt}_a\frac{\partial}{\partial x_a^t}.\]
%and also write
%\[\sum_{r,s} u^{wz}_{ars}\gamma^{rs} \coloneqq \pr_a^*[\alpha^w_a,\alpha^z_a]_{A_a}, \textnormal{ and } v^w_a = \sum_t v^{wt}_a\frac{\partial}{\partial x^t_a}\]
%so that
%\[[\gamma^{ij},\gamma^{k\ell}]_{A_1 \times A_2} = \sum_{r,s} u_{rs}^{ijk\ell}\gamma^{rs} \coloneqq \sum_{r,s}(u^{ik}_{1rs}+u^{j\ell}_{2rs})\gamma^{rs}, \textnormal{ and } \anchor_{A_1 \times A_2}(\gamma^{ij}) = \sum_tv^{it}_1 \frac{\partial}{\partial x_1^t} + v^{jt}_2 \frac{\partial}{\partial x_2^t}.\]
By using that 
\[\pr_a^*(\sum_t f_t \alpha^t_a) = \sum_t \pr_a^*(f_t) \pr_a^*(\alpha^t_a), \textnormal{ where } f_t \in C^\infty(M_a),\]
we see that 
\[[\gamma^w_a,\gamma^z_a] = \sum_t \pr_a^*(c^{wz}_{at})\gamma^t_a, \textnormal{ and } \anchor(\gamma^w_a) = \sum_t \pr_a^*(b_a^{wt})\frac{\partial}{\partial x_a^t}\]
(where we also wrote $\textstyle\frac{\partial}{\partial x_a^t}$ for the section $\pr_a^*(\textstyle\frac{\partial}{\partial x_a^t})$ of $\pr_1^*TM_1 \oplus \pr_2^*TM_2 = T(M_1 \times M_2)$). Since we have the structure equations (see Remark \ref{rema: Lie algbroid local nature})
\begin{equation*}\label{eq: Jacobiator of product algebroid, seperately}
    \bigodot_{y,w,z}\sum_s c_{as}^{wz}c_{at}^{ys} + b^{ys}_a \frac{\partial c_{at}^{wz}}{\partial x_a^s}=0 \textnormal{ for all } 1 \le t \le \textnormal{rank}(A_a),
\end{equation*}
and using that the pullback of smooth functions commutes with taking derivatives, we see that the Jacobiator $\textstyle\bigodot_{p,q,r} [[\gamma^{p}_a,\gamma^{q}_a]_{A_1 \times A_2},\gamma^{r}_a]_{A_1 \times A_2}$ is equal to zero, as claimed.

Lastly, we have a canonical map $\Gamma(A_1) \oplus \Gamma(A_2) \ra \Gamma(A_1 \times A_2)$ given by $(\alpha_1,\alpha_2) \mapsto \pr_1^*\alpha_1 + \pr_2^*\alpha_2$, which is a morphism of Lie algebras with the Lie algebroid structure defined here. Moreover, this condition uniquely determine the Lie algebroid structure on $A_1 \times A_2$, since the bracket has to be defined the way we did%; in other words, the image of the map $\Gamma(A_1) \oplus \Gamma(A_2) \ra \Gamma(A_1 \times A_2)$ generates $\Gamma(A_1 \times A_2)$ as a $C^\infty(M_1 \times M_2)$-module
.
\end{exam}

A different approach to defining a morphism of Lie algebroids is by instead using the pullback Lie algebroid (compare with Example \ref{exam: pull-back of a groupoid}). However, we will use the Lie algebroid structure of the product to define this Lie algebroid, and we will not use the pullback Lie algebroid to define morphisms of Lie algebroids. Instead, we will wait to introduce this class of examples until after the discussion about morphisms of Lie algebroids.

As in the groupoid case, we can extend the notion of an action of a Lie algebra to an action of a Lie algebroid.
\begin{defn}\cite{meinrenken}\label{defn: action of Lie algebroid}
Let $\algebroid$ be an algebroid and let $\mu: E \ra \base$ be a smooth map (note: $E$ is assumed to be Hausdorff). A \textit{Lie algebroid action of $\algebroid$ on $E$}, with moment map $\mu$, is given by a morphism of Lie algebras 
\[\rho: \Gamma(\algebr) \ra \mathfrak{X}(E)\]
such that, for all $\alpha \in \Gamma(\algebr)$ and $f \in C^\infty(\base)$, $\rho(f\alpha) = \mu^*f\rho(\alpha)$, and the induced map $\mu^*\algebr \ra TE$ (given by $\alpha(\mu(e)) \mapsto \rho(\alpha)(e)$) is a smooth map. If, in addition, $\mu: E \ra \base$ is a vector bundle, and $\rho$ maps into the linear vector fields of $\mathfrak{X}(E)$ (linear meaning homogeneous of degree $0$ as in \ref{defn: homogeneous functions}), then it is called a \textit{representation of $\algebr$}.
\end{defn}
\begin{rema}\cite{meinrenken}\label{rema: representation of Lie algebroid}
As in the Lie groupoid case, one can show that a representation of a Lie algebroid is the same as a morphism of Lie algebroids $\algebr \ra \mathfrak{gl}(E)$.
\end{rema}
We can now generalise the class of examples in Example \ref{exam: action Lie algebroid}.
\begin{exam}\cite{meinrenken}\label{exam: action Lie algebroid}
Let $\algebroid$ be an algebroid acting on a manifold $E$ with moment map $\mu: E \ra \base$. We can associate the Lie algebroid
\[\algebr \ltimes E \coloneqq \mu^*\algebr (= \algebr \tensor[_{\pi}]{\times}{_{\mu}} E) \ra E\]
to it. To see that it is a Lie algebroid, notice that, since the map
\[\mu^*\algebr \ra TE\]
is smooth, it is a morphism of vector bundles; we let this map be the anchor map. Moreover, since $\Gamma(\mu^*\algebr) \cong C^\infty(E) \otimes_{C^\infty(\base)} \Gamma(\algebr)$, we can extend the Lie bracket of $\Gamma(\algebr)$ (uniquely) to a Lie bracket on $\Gamma(\mu^*\algebr)$ using the Leibniz rule.
\end{exam}
\begin{rema}\cite{meinrenken}\label{rema: Lie algebroid of action Lie groupoid 2}
Let $\groupoid$ be a Lie groupoid and let $\algebroid$ be its Lie algebroid. If $E$ is a (Hausdorff) manifold such that $\groupoid$ acts on $E$, with moment map $\mu: E \ra \base$, one can show, similar to the way we did in Remark \ref{rema: Lie algebroid of action groupoid is action algebroid}, that we can describe an infinitesimal Lie algebroid action of $\algebr$ on $E$, such that the Lie algebroid of $\group \ltimes E$ is $\algebr \ltimes E$ (with respect to the same moment map), but we will not go into this. However, equivalently, an action of $\algebr$ on $E$, with moment map $\mu$, is the same as a Lie algebroid structure on $\mu^*\algebr$ such that the map $\pr_1: \mu^*\algebr = \algebr \tensor[_{\pi}]{\times}{_{\mu}} E \ra \algebr$ is a morphism of Lie algebroids. Once we discussed morphisms of Lie algebroids, it is not hard to see this way that the Lie algebroid of $\group \ltimes E$ is $\algebr \ltimes E$.
\end{rema}
Lastly, foliations admit a particularly nice description in terms of Lie algebroids.
\begin{exam}\cite{meinrenken}\label{exam: foliations as Lie algebroids}
Let $\mathcal{F}$ be a foliation on a manifold $M$. Then, precisely because $T\mathcal{F} \subset TM$ is an involutive distribution, the vector bundle
\[T\mathcal{F} \ra M,\]
equipped with the usual Lie bracket of vector fields, is a Lie algebroid with anchor map the inclusion map $T\mathcal{F} \hookrightarrow TM$.
\end{exam}
\begin{rema}\cite{IekeMariuscyclic}\label{rema: foliation integration as Lie algebroids}
It is not hard to see, for a foliated manifold $(M,\mathcal{F})$, that the Lie algebroid of the Monodromy and Holonomy groupoids are given by the Lie algebroid $T\mathcal{F} \ra M$. Later, we will see that, if $\groupoid$ is a groupoid integrating $T\mathcal{F}$, then we obtain morphisms of Lie groupoids
\[\Mon(M,\mathcal{F}) \rightarrowdbl \group \rightarrowdbl \Hol(M,\mathcal{F}).\]
These morphisms are surjective local diffeomorphisms, and the composition of these maps equals the map $F$ from Remark \ref{rema: groupoid morphism Mon -> Hol}. In particular, if $\Mon(M,\mathcal{F})$ and $\Hol(M,\mathcal{F})$ are both non-Hausdorff, then the Lie algebroid $T\mathcal{F} \ra M$ does not admit a Hausdorff integration (e.g. the Reeb foliation; see also \cite{hoyohausdorff}).
\end{rema}
\subsection{Morphisms of Lie algebroids}\label{sec: Lie subalgebroids and morphisms}
We start this section by introducing one of the simplest type of morphisms of Lie algebroids: Lie subalgebroids (the inclusion of which is a morphism of Lie algebroids). 
\subsubsection{Lie subalgebroids}
\begin{defn}\cite{meinrenken}\label{defn: Lie subalgebroid}
A \textit{Lie subalgebroid} of a Lie algebroid $\algebroid$ is an anchored subbundle $\subalgebroid$ (i.e. $\anchor_\algebr(\subalgebr) \subset T\subbase$), and
\[\Gamma(\algebr,\subalgebr) \coloneqq \{\alpha \in \Gamma(\algebr) \mid \alpha|_\subbase \in \Gamma(\subalgebr)\}\]
is a Lie subalgebra of $\Gamma(\algebr)$.
\end{defn}
\begin{exam}\cite{meinrenken}\label{exam: Lie subalgebroid foliation}
Let $M$ be a manifold. Notice that Lie subalgebroids of $TM$, over $M$, are precisely given by involutive distributions, i.e. foliations. In fact, from this we see that a foliation is a Lie algebroid $\algebroid$ which has an injective anchor map. Indeed, if $\anchor: \algebr \ra T\base$ is injective, then $\anchor(\algebr) \subset T\base$ is a distribution, and this distribution is involutive by Proposition \ref{prop: anchor induces Lie algebra morphism on sections}.
\end{exam}
\begin{prop}\cite{meinrenken}\label{prop: Lie subalgebroid}
A Lie subalgebroid $\subalgebroid$ of a Lie algebroid $\algebroid$ admits a unique Lie algebroid structure such that the restriction map $\Gamma(\algebr,\subalgebr) \ra \Gamma(\subalgebr)$ is a Lie algebra morphism. If $\subbase \subset \base$ is closed (and embedded), then this map is surjective.
\end{prop}
\begin{proof}
We have to check that $\Gamma(\subalgebr)$ inherits a Lie bracket such that $\anchor_\subalgebr \coloneqq \anchor_\algebr|_\subalgebr$ together with this Lie bracket satisfy the Leibniz rule. It suffices to prove that $\Gamma(\algebr,0_\subbase)$ is an ideal of $\Gamma(\algebr,\subalgebr)$, since then $\Gamma(\subalgebr) \cong \Gamma(\algebr,\subalgebr)/\Gamma(\algebr,0_\subbase)$ inherits a Lie bracket from $\Gamma(\algebr,\subalgebr)$ (which is a Lie subalgebra of $\Gamma(\algebr)$ by assumption). Observe that
\[\Gamma(\algebr,0_\subbase) = \{\alpha \in \Gamma(\algebr) \mid \alpha|_\subbase = 0\} = \{\alpha \in \Gamma(\algebr) \mid \alpha = \textstyle\sum_i f_i\alpha_i \textnormal{ where } f_i|_\subbase = 0 \textnormal{ and } \alpha_i \in \Gamma(\algebr)\},\]
hence for all $\beta \in \Gamma(\algebr,\subalgebr)$ and $\alpha = \textstyle\sum_if_i\alpha_i \in \Gamma(\algebr,0_\subbase)$, we have 
\[\bracket{\beta,\textstyle\sum_if_i\alpha_i} = \textstyle\sum_i\bracket{\beta,f_i\alpha_i} = \textstyle\sum_if_i\bracket{\beta,\alpha_i} + \anchor(\beta)(f_i)\alpha_i\]
which is an element of $\Gamma(\algebr,0_\subbase)$, because $f_i|_\subbase=0$ and $\anchor(\beta)(f_i)|_\subbase=0$ since $\anchor(\beta)|_\subbase \in \mathfrak{X}(\subbase)$. The last statement follows from the fact that, if $\subbase \subset \base$ is closed and embedded, then we can extend a section of $\subalgebr \ra \subbase$ to a section of $\algebroid$.
\end{proof}
\begin{rema}\cite{meinrenken}\label{rema: Lie subalgebroid sections}
Let $\subalgebroid$ be a Lie subalgebroid of a Lie algebroid $\algebroid$. Then the image of the restriction map $\Gamma(\algebr,\subalgebr) \ra \Gamma(\subalgebr)$ generates $\Gamma(\subalgebr)$ as a $C^\infty(\subbase)$-module. This follows by noticing that a local frame on $\subalgebr$ can be realised as the restriction of a local frame on $\algebr$. 
\end{rema}
We define pairs of Lie algebroids (and pairs of vector bundles) as follows.
\begin{defn}\label{defn: pairs of Lie algebroids}
A pair of Lie algebroids is a tuple $(\algebroid,\subalgebroid)$ such that $\subalgebroid$ is a subalgebroid of $\algebroid$ for which $\subalgebr \subset \algebr$ and $\subbase \subset \base$ are closed embedded submanifolds. Similarly, a pair of (anchored) vector bundles is a tuple $(E \ra M, F \ra N)$ such that $F \ra N$ is a(n) (anchored) subbundle of $E \ra M$ for which $F \subset E$ and $N \subset M$ are closed embedded submanifolds. If $E \ra M$ is of rank $k$, and $F \ra N$ is of rank $\ell$, then we will say that the \textit{rank} of $(E \ra M, F \ra N)$ is $(k,\ell)$.
\end{defn}

The following lemma is very useful to prove that (anchored) vector subbundles are Lie subalgebroids.
\begin{lemm}\cite{meinrenken}\label{lemm: anchored subbundle is Lie subalgebroid iff commutator}
Let $\algebroid$ be a Lie algebroid and let $\subalgebroid$ be an anchored vector subbundle (so $\anchor_\algebr(\subalgebr) \subset T\subbase$). If the image of a subset $R \subset \Gamma(\algebr,\subalgebr)$ under the map $\Gamma(\algebr,\subalgebr) \ra \Gamma(\subalgebr)$ generates $\Gamma(\subalgebr)$ as $C^\infty(\subbase)$-module, then $\subalgebr$ is a Lie subalgebroid if and only if $\bracket{R,R} \subset \Gamma(\algebr,\subalgebr)$.
\end{lemm}
\begin{proof}
By the Leibniz rule, we may replace $R$ with the $C^\infty(\base)$-module generated by it. In particular, the image of $R$ under the restriction map $\Gamma(\algebr,\subalgebr) \ra \Gamma(\subalgebr)$ equals $\Gamma(\subalgebr)$. Now, if we have a section $\alpha \in \Gamma(\algebr,\subalgebr)$, then $\alpha|_\subbase \in \Gamma(\subalgebr)$, so we can find a section $\alpha_R \in R$ such that $\alpha_R|_\subbase=\alpha|_\subbase$. From the equation $\alpha = \alpha_R + (\alpha - \alpha_R)$, we see that
\[\Gamma(\algebr,\subalgebr) = R + \Gamma(\algebr,0_\subbase).\]
Since $\bracket{\Gamma(\algebr,\subalgebr),\Gamma(\algebr,0_\subbase)} \subset \Gamma(\algebr,0_\subbase)$ (see Proposition \ref{prop: Lie subalgebroid}), the result follows.
\end{proof}
\subsubsection{Morphisms of Lie algebroids}
We are now ready to define what a Lie algebroid morphism is.
\begin{defn}\cite{meinrenken}\label{defn: morphism of Lie algebroids}
Let $\algebr \xra{\pi_\algebr} \base$ and $\subalgebr \xra{\pi_\subalgebr} \subbase$ be Lie algebroids. A morphism of Lie algebroids is a bundle map $\varphi: \algebr \ra \subalgebr$ (with base map $f$) such that the graph 
\[\gr\varphi \xra{\pi_\algebr \times \pi_\subalgebr} \gr f,\]
is a Lie subalgebroid of $\algebr \times \subalgebr$.
\end{defn}
\begin{rema}\cite{meinrenken}\label{rema: Lie algebroid morphism commutes with anchors}
Notice that, for a bundle map $\varphi: \algebr \ra \subalgebr$ between Lie algebroids (with base map $f$), the diagram
\begin{center}
\begin{tikzcd}
    \algebr \ar[r,"\varphi"] \ar[d,"\anchor_\algebr"] & \subalgebr \ar[d,"\anchor_\subalgebr"] \\
    T\base \ar[r,"df"] & T\subbase
\end{tikzcd}
\end{center}
commutes if and only if $\anchor_{\algebr \times \subalgebr}(\gr \varphi) = (\anchor_\algebr \times \anchor_\subalgebr)(\gr \varphi) \subset \gr df = T(\gr f)$.
\end{rema}
As mentioned before, a Lie algebroid morphism over the same base takes the following form.
\begin{prop}\cite{meinrenken}\label{prop: morphism of Lie algebroids over the same base}
Let $\algebroid$ and $\subalgebr \ra \base$ be two Lie algebroids over the same base. A morphism of vector bundles $\varphi: \algebr \ra \subalgebr$, with base map $\id_\base$, is a morphism of Lie algebroids if and only if $\varphi$ is a Lie algebra morphism on the level of sections and $\anchor_\algebr = \anchor_\subalgebr \circ \varphi$. 
\end{prop}
\begin{proof}
If $\varphi$ is a morphism of Lie algebroids, then by Remark \ref{rema: Lie algebroid morphism commutes with anchors} we have $\anchor_\algebr = \anchor_\subalgebr \circ \varphi$ if and only if $\anchor_{\algebr \times \subalgebr}(\gr \varphi) \subset T(\gr \id_\base)$. It remains to prove that, subject to the condition that $\anchor_{\algebr \times \subalgebr}(\gr \varphi) \subset T(\gr \id_\base)$, we have that
\[\Gamma(\algebr \times \subalgebr,\gr \varphi) \textnormal{ is a Lie subalgebra of } \Gamma(\algebr \times \subalgebr)\]
is equivalent to requiring that $\varphi$ is a Lie algebra morphism on the level of sections. This follows by taking $R = \{(\alpha,\varphi(\alpha)) \in \Gamma(\algebr \times \subalgebr, \gr \varphi) \mid \alpha \in \Gamma(\algebr)\}$ in the lemma below.
\end{proof}
Later in this section we will also discuss ``differentiating'' Lie groupoid morphisms and ``integrating'' Lie algebroid morphisms.

\subsubsection{Clean intersection and Lie algebroids}
In general, an intersection of two Lie subalgebroids is not a Lie subalgebroid again. However, a sufficient condition on the intersection to ensure that this will be the case, is to require for it to be a clean intersection. Before we prove this statement, we will need two intermediate results regarding vector bundles. First of all, we have the following elementary statement about clean intersecting submanifolds% (for a proof, see \cite{Hormander})
.
\begin{lemm}\cite{Hormander}\label{lemm: clean intersecting submanifolds linear subspace}
Let $M$ be a manifold of dimension $n$, and let $N_1$ and $N_2$ be embedded submanifolds of $M$ which intersect cleanly. Then, for all $y \in N_1 \cap N_2$, we can find a chart $(U,\varphi)$ of $M$ around $y$ that is adapted to $N_1$, $N_2$, and $N_1 \cap N_2$. In particular, we can find a chart $(U,\varphi)$ in which $\varphi(U \cap N_1)$, $\varphi(U \cap N_2)$ and $\varphi(U \cap N_1 \cap N_2)$ are vector subspaces of $\mathbb{R}^n$.
\end{lemm}
%\begin{proof}
%First of all, choose local coordinates $(U,\varphi)=(U,x_1,\dots,x_p,y_1,\dots,y_q,z_1,\dots,z_\ell)$ of $M$ such that $\varphi(y)=0$, and $\{y_1=\dots=y_q=z_1=\dots=z_\ell=0\}$ corresponds to $N_1 \cap N_2$ and $\{z_1=\dots=z_\ell=0\}$ corresponds to $N_1$. By shrinking $U$, can describe $N_2$ as the zero set of, say, a submersion $f=(f_1,\dots,f_k): U \ra \mathbb{R}^k$. Notice that now $f \circ \varphi^{-1}(x_1,\dots,x_p,0)=0$, so we can write
%\[f_m \circ \varphi^{-1} = \sum_i f_{mi}y_i + \sum_j f_{mj}z_j; \quad df_m \circ d\varphi^{-1} = \sum_i f_{mi}dy_i + \sum_j f_{mj}dz_j\]
%Let $y \in N_1 \cap N_2$. Then, pick a basis $\{a_i\}$ of $T_y(N_1 \cap N_2) = T_yN_1 \cap T_yN_2$. We can extend this basis to a basis $\{a_i,b_j\}$ of $T_yN_1$, and to a basis $\{a_i,c_k\}$ of $T_yN_2$, and then extend the basis $\{a_i,b_j,c_k\}$ of $T_yN_1 + T_yN_2$ to a basis $\{a_i,b_j,c_k,d_\ell\}$ of $T_yM$. This gives us an isomorphism $T_yM \ra \mathbb{R}^n$, so by the inverse function theorem we find a diffeomorphism $\varphi: U \ra O$, where $U \subset M$ and $O \subset \mathbb{R}^n$ are open subsets. This gives us a chart $(U,\varphi)$ on $M$, and by construction it has the desired properties. The last statement holds by making sure $\varphi$ maps onto $\mathbb{R}^n$ (and mapping $y$ to $0$), e.g. by shrinking $U$ so that $O$ is an open ball, and then post-compose the chart with a scaling diffeomorphism $O \ra \mathbb{R}^n$. 
%\end{proof}
To ensure that extensions of sections of a vector bundle exist, we will restrict to the case of closed embedded submanifolds.
\begin{lemm}\cite{meinrenken}\label{lemm: clean intersection of vector bundles}
Let $F_1 \ra N_1$ and $F_2 \ra N_2$ be subbundles of a vector bundle $E \ra M$. If $F_1$ and $F_2$ intersect cleanly (as manifolds), then $F_1 \cap F_2 \ra N_1 \cap N_2$ is a vector subbundle of $E \ra M$. If $N_1$ and $N_2$ are closed (embedded) submanifolds of $M$, then the restriction map
\[\Gamma(E,F_1) \cap \Gamma(E,F_2) \ra \Gamma(F_1 \cap F_2)\]
is surjective.
\end{lemm}
\begin{proof}
First of all, that $F_1 \cap F_2$ is a vector bundle over $N_1 \cap N_2$ follows from the more general statement Corollary \ref{coro: vector subbundle is invariant submanifold} which is proven much later. To show that $N_1$ intersects cleanly with $N_2$, it remains to prove that $TN_1 \cap TN_2 \cong T(N_1 \cap N_2)$. However, since $TF_1 \cap TF_2 = T(F_1 \cap F_2)$ by assumption, this follows from the sequence of vector bundle isomorphisms
\begin{align*}
    (F_1 \cap F_2) \oplus (TN_1 \cap TN_2) &\cong (F_1 \oplus TN_1) \cap (F_2 \oplus TN_2) \\
    &\cong TF_1|_{N_1} \cap TF_2|_{N_2} \\
    &\cong T(F_1 \cap F_2)|_{N_1 \cap N_2} \cong (F_1 \cap F_2) \oplus T(N_1 \cap N_2),
\end{align*}
where we repeatedly used that for a vector bundle $E' \ra M'$ we have a canonical isomorphism $TE'|_{M'} \cong E' \oplus TM'$. We will prove the last statement by using the following claim: whenever we have a vector bundle $E' \ra M'$, and $N_1'$ and $N_2'$ are clean intersecting closed embedded submanifolds of $M'$, then we can extend sections $\sigma_a' \in \Gamma(E'|_{N_a'})$ (for $a=1,2$), which agree on $N_1' \cap N_2'$, to a section $\sigma' \in \Gamma(E')$ extending both $\sigma_1'$ and $\sigma_2'$. We will first prove the statement using the claim: let $\sigma \in \Gamma(F_1 \cap F_2)$. Since $N_1$ and $N_2$ intersect cleanly, they intersect cleanly in $N_a$ (for $a=1,2$). Since $N_1 \cap N_2 \subset N_a$ is closed, we can find extensions $\widetilde{\sigma}_a$ of $\sigma$ to $F_a$ (for $a=1,2$) such that $\widetilde{\sigma}_a|_{N_1 \cap N_2}=\sigma$. Now, again because $N_1$ intersects cleanly with $N_2$, and $N_a \subset M$ is closed, there exists an extension $\widetilde{\sigma}$ to $E$ such that $\widetilde{\sigma}|_{N_a} = \widetilde{\sigma}_a$ (for $a=1,2$). Then $\widetilde{\sigma}$ maps to $\sigma$ under the map $\Gamma(E,F_1) \cap \Gamma(E,F_2) \ra \Gamma(F_1 \cap F_2)$. It remains to prove the claim (we drop the apostrophes). It suffices to prove the claim for local sections $\sigma_a \in \Gamma(U, E|_{N_a})$ (for $a=1,2$), where $(U,\varphi)$ is a chart of $M$ such that $N_1$, $N_2$ and $N_1 \cap N_2$ are locally given by vector subspaces of $\mathbb{R}^n$ (see Lemma \ref{lemm: clean intersecting submanifolds linear subspace}). So, to ease the notation, we replace $M$ with $\varphi(U)$ and $N_a$ with $\varphi(U \cap N_a)$ (for $a=1,2$). The sections are now locally of the form
\[\sigma_a = \sum_j f_{a,j} \cdot s^j|_{N_a}, \textnormal{ where } f_{a,j} \in C^{\infty}(N_a)\]
such that $f_{a,1}|_{N_1 \cap N_2} = f_{a,2}|_{N_1 \cap N_2}$, where $s^1,\dots,s^r$ is a local frame of $E$. Since $N_1$ and $N_2$ are vector subspaces, $N_1 + N_2 \subset M$ is also a vector subspace. The sections $\sigma_1$ and $\sigma_2$ agree on $N_1 \cap N_2$, so the section 
\[\sigma_{3} \coloneqq \sigma_1 + \sigma_2 \in \Gamma(E|_{N_1 + N_2})\]
is well-defined, and upon restriction to $N_a$ (for $a=1,2$), $\sigma_{3}$ becomes the section $\sigma_a$. Since $N_1 + N_2$ is, in particular, a closed embedded submanifold of $M$, we can extend the section $\sigma_{3}$ further to a section of $\Gamma(E)$. This proves the statement.
\end{proof}

%\begin{theo}\label{theo: invariant submanifold under scalar multiplication are vector subbundles}
%Let $\pi: E \ra M$ be a vector bundle and let $F \subset E$ be an embedded submanifold. Then $F \ra N$ (with $N \coloneqq \pi(F)$) is a subbundle if and only if it is invariant under scalar multiplication; that is, if $(x,\xi) \in F$, then $(x,\lambda\xi) \in F$ for all $\lambda \in \mathbb{R}$.
%\end{theo}

\begin{prop}\cite{meinrenken}\label{prop: clean intersection of Lie algebroids}
Let $(\algebroid, \subalgebr_1 \ra \subbase_1)$ and $(\algebroid, \subalgebr_2 \ra \subbase_2)$ be pairs of Lie algebroids (see Definition \ref{defn: pairs of Lie algebroids}). If $\subalgebr_1$ and $\subalgebr_2$ intersect cleanly (as manifolds), then $\subalgebr_1 \cap \subalgebr_2 \ra \subbase_1 \cap \subbase_2$ inherits a Lie algebroid structure such that $\subalgebr_1 \cap \subalgebr_2 \subset \subalgebr_1$ and $\subalgebr_1 \cap \subalgebr_2 \subset \subalgebr_2$ are subalgebroids. Moreover, the restriction map 
\[\Gamma(\algebr,\subalgebr_1) \cap \Gamma(\algebr,\subalgebr_2) \ra \Gamma(\subalgebr_1 \cap \subalgebr_2)\]
is surjective.
\end{prop}
\begin{proof}
To see that $\subbase_1$ and $\subbase_2$ automatically intersect cleanly, and that $\subalgebr_1 \cap \subalgebr_2 \ra \subbase_1 \cap \subbase_2$ is a vector subbundle of both $\subalgebr_1 \ra \subbase_1$ and $\subalgebr_2 \ra \subbase_2$, are consequences of Lemma \ref{lemm: clean intersection of vector bundles}. Since
\[\anchor(\subalgebr_1 \cap \subalgebr_2) \subset \anchor(\subalgebr_1) \cap \anchor(\subalgebr_2) \subset T\subbase_1 \cap T\subbase_2 = T(\subbase_1 \cap \subbase_2),\]
and $\Gamma(\algebr,\subalgebr_1) \cap \Gamma(\algebr,\subalgebr_2) \subset \Gamma(\algebr, \subalgebr_1 \cap \subalgebr_2)$ is a Lie subalgebra of $\Gamma(\algebr)$, it suffices to prove the last statement by Lemma \ref{lemm: anchored subbundle is Lie subalgebroid iff commutator}. This follows directly from Lemma \ref{lemm: clean intersection of vector bundles}.
\end{proof}

As a consequence, we obtain the following result.
\begin{coro}\cite{meinrenken}\label{coro: inverse image of algebroid}
Let $\varphi: \algebr \ra \algebrtwo$ be a morphism of Lie algebroids, with base map $f: \base \ra \basetwo$, and let $(\algebrtwoid,\subalgebrtwoid)$ be a pair of Lie algebroids. Then $\varphi^{-1}(\subalgebrtwo) \ra \varphi^{-1}(\subbasetwo)$ is a closed Lie subalgebroid of $\algebroid$.
\end{coro}
\begin{proof}
The result follows by the identification $\algebr \cong \gr \varphi$ as Lie algebroids. Indeed, under this identification, $\varphi^{-1}(\subalgebrtwo)$ maps bijectively onto $\gr \varphi \cap (\algebr \times \subalgebrtwo)$, which is a clean intersection of Lie subalgebroids of $\algebr \times \algebrtwo$. 
\end{proof}
Similar statements hold for Lie groupoids (e.g. the intersection of clean intersecting Lie subgroupoids is again a Lie subgroupoid).
\subsubsection{Pullback Lie algebroids}
As mentioned in Section \ref{sec: Examples of Lie algebroids}, we will discuss the class of examples arising as a pullback. We will use Corollary \ref{coro: inverse image of algebroid}.
\begin{exam}[Pullback algebroids]\cite{meinrenken}\label{exam: pullback of a Lie algebroid}
Let $\algebroid$ be a Lie algebroid and let $f: M \ra \base$ be a smooth map such that the maps
\[\anchor_{\algebr}: \algebr \ra T\base \textnormal{ and } df: TM \ra T\base\]
intersect cleanly (see Definition \ref{defn: clean intersetion}). By this last assumption, we have an associated \textit{pullback Lie algebroid} 
\[f^!\algebr \coloneqq TM \tensor[_{df}]{\times}{_{\anchor_{\algebr}}} \algebr \ra M,\]
with anchor map $\anchor$ defined as the map from the pullback square:
\begin{center}
\begin{tikzcd}
    f^!\algebr \ar[r,"\anchor"] \ar[d] & TM \ar[d,"df"] \\
    \algebr \ar[r,"\anchor_{\algebr}"] & T\base.
\end{tikzcd}
\end{center}
The Lie bracket is the one induced by the bracket of $TM \times \algebr$ (see Example \ref{exam: product of Lie algebroids}), by using Corollary \ref{coro: inverse image of algebroid} applied to the fact that $f^!\algebr$ is the pre-image of $\gr df \cong T(\gr f)$ under the map
\[\id_{TM} \times \anchor_{\algebr}: TM \times \algebr \ra TM \times T\base \cong T(M \times \base)\]
(note: we used here that $M$ can be identified with the submanifold $\gr f \subset M \times \base$). Observe that, from the above diagram, we obtain an anchored vector bundle morphism $f^!\algebr \ra \algebr$. It is the restriction of the map $\pr_2: TM \times \algebr \ra \algebr$, and therefore it is a morphism of Lie algebroids.
\end{exam}
We now discuss the integrability of the above class of examples (including the integrability of product Lie algebroids).
\begin{rema}\cite{meinrenken}\label{rema: Lie algebroid integration of product and pullback}
Let $\algebr_1 \ra M_1$ and $\algebr_2 \ra M_2$ be two Lie algebroids that integrate to groupoids $\group_1 \rra M_1$ and $\group_2 \rra M_2$, respectively. Then the product Lie algebroid $\algebr_1 \times \algebr_2 \rra M_1 \times M_2$ integrates to $\group_1 \times \group_2 \rra M_1 \times M_2$. That $\ker d(\source_1 \times \source_2)|_{\identity(M_1 \times M_2)}$ is naturally isomorphic to $\algebr_1 \times \algebr_2$ as anchored vector bundles follows from the commutativity of the diagrams
\begin{center}
\begin{tikzcd}
    T(\group_1 \times \group_2) \ar[r,"\sim"] \ar[d,"d(\source_1 \times \source_2)"] & T\group_1 \times T\group_2 \ar[d,"d\source_1 \times d\source_2"] & \phantom{A} \ar[d,phantom,"\text{ and }"] & T(\group_1 \times \group_2) \ar[r,"\sim"] \ar[d,"d(\target_1 \times \target_2)"] & T\group_1 \times T\group_2 \ar[d,"d\target_1 \times d\target_2"] \\
    T(M_1 \times M_2) \ar[r,"\sim"] & TM_1 \times TM_2 & \phantom{A} & T(M_1 \times M_2) \ar[r,"\sim"] & TM_1 \times TM_2.
\end{tikzcd}
\end{center} 
To check that the Lie brackets are the same, recall that the Lie bracket is uniquely determined by the property that the canonical map $\Gamma(\algebr_1) \oplus \Gamma(\algebr_2) \ra \Gamma(\algebr_1 \times \algebr_2)$ is a Lie algebra morphism. But, by identifying $\Gamma(A_a) \cong \mathfrak{X}_{R\textnormal{-inv}}(\group_a)$ ($a=1,2$) and $\Gamma(A_1 \times A_2) \cong  \mathfrak{X}_{R\textnormal{-inv}}(\group_1 \times \group_2)$, this map is the restriction of the canonical map $\mathfrak{X}(\group_1) \oplus \mathfrak{X}(\group_2) \ra \mathfrak{X}(\group_1 \times \group_2)$ which is readily verified to be a morphism of Lie algebras. %The inclusion $\group_a \hookrightarrow \group_1 \times \group_2$ gives rise to a Lie algebra morphism $\mathfrak{X}_{R\textnormal{-inv}}(\group_a) \ra \mathfrak{X}_{R\textnormal{-inv}}(\group_1 \times \group_2)$ (for $a=1,2$), which, under the Lie algebra isomorphisms $\Gamma(A_a) \cong \mathfrak{X}_{R\textnormal{-inv}}(\group_a)$ and $\Gamma(A_1 \times A_2) \cong  \mathfrak{X}_{R\textnormal{-inv}}(\group_1 \times \group_2)$ (note: this is how we defined the Lie algebra structures on $\Gamma(A_a)$ and $\Gamma(A_1 \times A_2)$), is precisely the canonical map $\Gamma(\algebr_a) \ra \Gamma(\algebr_1 \times \algebr_2)$. 

Now let $\algebroid$ be a Lie algebroid, and suppose $f: M \ra \base$ is a smooth map such that $\anchor_\algebr$ and $df$ have a clean intersection. If $\algebroid$ integrates to $\groupoid$, and $f$ has a clean intersection with $(\target,\source)$, then $f^!\algebr$ integrates to $f^!\group$. Indeed, by using the identification $M \cong \gr f$, $f^!\group$ is a subgroupoid of the product groupoid of $M \times M \rra M$ (the pair groupoid) with $\group \rra \base$. As anchored vector bundles, it is now readily verified that $f^!\algebr$ is isomorphic to the Lie algebroid belonging to $f^!\group$. Since the Lie bracket of the Lie algebroid of the product groupoid of $M \times M$ with $\group$ coincides with the Lie bracket of the product Lie algebroid $TM \times \algebr$, and $f^!\algebr$ inherits this bracket, it is clear that $f^!\algebr$ really is the Lie algebroid of $f^!\group$.
\end{rema}
We now state the following universal property of pullback Lie algebroids.
\begin{prop}\cite{waldron2015lie}\label{prop: universal property of pullbacks}
Let $\algebroid$ be a Lie algebroid, and let $f: M \ra \base$ be a smooth map such that the maps $\anchor_\algebr$ and $df$ intersect cleanly. Then $f^!\algebr$ is the unique (up to unique isomorphism) Lie algebroid such that whenever $\subalgebr \ra N$ is a Lie algebroid together with a morphism $\varphi: \subalgebr \ra \algebr$ over a base map $f \circ g$, where $g$ is a smooth map $N \ra M$, then there is a unique Lie algebroid morphism $\widetilde\varphi: \subalgebr \ra f^!\algebr$, with base map $g$, such that the post-composition with the map $f^!\algebr \ra \algebr$ equals $\varphi$.
\end{prop}
We state the following useful result as a consequence of this property.
\begin{prop}\cite{meinrenken}\label{prop: (fg)!=g!f!}
Let $\algebroid$ be a Lie algebroid, and let $f: M \ra \base$ and $g: N \ra M$ be smooth maps such that $\anchor_\algebr$ and $df$ have clean intersection, and $\anchor_{f^!\algebr}$ and $dg$ have clean intersection. Then $\anchor_{\algebr}$ and $d(f \circ g)$ have clean intersection, and $(f \circ g)^!\algebr \cong g^!f^!\algebr$.
\end{prop}
\begin{proof}
The clean intersection statement is readily verified, and the last statement follows by the universal property of pullback Lie algebroids.
\end{proof}

\subsubsection{Differentiation and integration of morphisms}
We will now prove that a Lie groupoid morphism induces a morphism of Lie algebroids between the induced Lie algebroids. We will prove this in three steps.
\begin{lemm}\cite{meinrenken}\label{lemm: differentiaion of morphisms}
Let $\groupoid$ be a groupoid, and let $\subgroupoid$ be a subgroupoid of $\groupoid$. Denote by $\algebroid$ the Lie algebroid of $\groupoid$ and by $\subalgebroid$ the Lie algebroid of $\subgroupoid$. Then the inclusion $\iota: \subgroup \ra \group$ induces an inclusion of Lie algebroids $\textnormal{Lie}(\iota): B \ra A$.
\end{lemm}
\begin{proof}
It is clear that $B = \ker d\source_\subgroup|_{\identity(\subbase)}$ can be seen as a vector subbundle of $A = \ker d\source|_{\identity(\base)}$ (note: $\source_\subgroup = \source|_{\subgroup}$). Observe that $d\target$ commutes with $d\iota$ (by the chain rule), so it is clear that $\anchor(B) \subset T\subbase$. It remains to check that $\Gamma(A,B) \subset \Gamma(A)$ is a Lie subalgebra. Under the Lie algebra isomorphism $\Gamma(A) \cong \mathfrak{X}_{R\textnormal{-inv}}(\group)$, $\Gamma(A,B)$ corresponds to the right-invariant vector fields of $\group$ which are tangent to $\subgroup$. The Lie bracket of two such vector fields is again tangent to $\subgroup$, which shows that $\Gamma(A,B) \subset \Gamma(A)$ is a Lie subalgebra.
\end{proof}
\begin{lemm}\cite{meinrenken}\label{lemm: graph of Lie groupoid morphism}
Let $\groupoid$ and $\subgroupoid$ be two groupoids and let $F: \group \ra \subgroup$ be a map together with a base map $f: \base \ra \subbase$. Then $(F,f)$ is a morphism of Lie groupoids if and only if $\gr F \rra \gr f$ is a subgroupoid of $\group \times \subgroup$.
\end{lemm}
\begin{proof}
The graph of a map is a smooth submanifold of a product manifold if and only if the map is smooth, so $\gr F$ is a smooth manifold if and only if $F: \group \ra \subgroup$ is a smooth map. It remains to show that a smooth map $F: \group \ra \subgroup$ with base map $f$ is a morphism of Lie groupoids if and only if $\gr F \ra \gr f$ is a subcategory of $\group \times \subgroup$. This follows immediately by spelling out what it means for $\gr F \ra \gr f$ to be a subgroupoid of $\group \times \subgroup$.
\end{proof}
\begin{prop}\cite{meinrenken}\label{prop: differentiation of morphisms}
Let $\groupoid$ and $\subgroupoid$ be two groupoids and let $F: \group \ra \subgroup$ be a morphism of Lie groupoids with base map $f: \base \ra \subbase$. Then
\[dF: T\group \ra T\subgroup\]
restricts to a morphism $\textnormal{Lie}(F)$ of Lie algebroids $\algebr \ra \subalgebr$.
\end{prop}
\begin{proof}
First of all, $dF$ restricts to a map $A \ra B$. Indeed, if $\xi \in \ker d\source_\group(\id_x)$, where $x \in \base$, then
\[d\source_\subgroup(F(\id_x))dF(\id_x)\xi = d(\source_\subgroup \circ F)(\id_x)\xi = d(f \circ \source_\group)(\id_x)\xi = 0.\]
Since $\ker d\source_\group \subset T\group$ and $\ker d\source_\subgroup \subset T\subgroup$ are subbundles (and so are $A$ and $B$, respectively), and $dF$ is a morphism of vector bundles, the map $\textnormal{Lie}(F)$ is a morphism of vector bundles. Since $dF$ commutes with $d\target$, $\textnormal{Lie}(F)$ is a morphism of anchored vector bundles. It remains to prove that the graph of $\textnormal{Lie}(F)$ is a Lie subalgebroid of $A \times B$. This follows from the two preceding lemmas: by Lemma \ref{lemm: graph of Lie groupoid morphism}, $\gr F \rra \gr f$ is a Lie subgroupoid of $\group \times \subgroup \rra \base \times \subbase$. It is clear that $\gr \textnormal{Lie}(F) \ra \gr f$ equals the Lie algebroid of $\gr F \rra \gr f$ as subalgebroids of $\algebr \times \subalgebr$, so by Lemma \ref{lemm: differentiaion of morphisms} the result follows.
\end{proof}
The appropriate notion of integrability of morphisms of Lie algebroids can be described as follows. 
\begin{defn}\cite{ruimarius}\label{defn: integrability of morphisms}
Let $\algebroid$ and $\subalgebroid$ be Lie algebroids that integrate to Lie groupoids $\groupoid$ and $\subgroupoid$, respectively. A morphism of Lie groupoids $F: \group \ra \subgroup$ is said to \textit{integrate} a morphism of Lie algebroids $\varphi: \algebr \ra \subalgebr$ if the morphism of Lie algebroids $\textnormal{Lie}(F)$ equals to $\varphi$. 
\end{defn}
We will now finish the proof of Proposition \ref{prop: s-(simply )connected groupoid of integrable algebroid} by proving a general, and useful, result regarding integrability of morphisms of Lie algebroids. We will use the theory of principal groupoid bundles (see Section \ref{sec: Morita equivalence}), and we will use the notion of a partial $\mathcal{F}$-connection with respect to a foliation $\mathcal{F}$ of the base space of a principal bundle.
\begin{defn}\cite{moerdijk_mrcun_2003}\label{defn: partial F-connection}
Let $\groupoid$ be a Lie groupoid, let $\pi: P \ra M$ be a principal $\group$-bundle with multiplication $\nu$ and moment map $\mu$, and let $\mathcal{F}$ be a foliation of $M$. A subbundle $\mathcal{H}$ of $T\pi^*\mathcal{F}$, where $\pi^*\mathcal{F}$ is the pullback foliation along $\pi$, is called a \textit{partial $\mathcal{F}$-connection} on $P$ if 
\begin{align*}
    d\mu(\mathcal{H})=0 \textnormal{ and } &\mathcal{H}_{p \cdot g} = d\nu_g(p)\mathcal{H}_p \textnormal{ for all $(g,p) \in \group \tensor[_{\target}]{\times}{_{\mu}} P$, and } \\
    T\pi^*\mathcal{F} &= \mathcal{V} \oplus \mathcal{H},
\end{align*} 
where $\mathcal{V} \coloneqq \ker d\pi$.
The connection is called \textit{flat} if $\mathcal{H}$ is integrable.
\end{defn}
Before we prove the statement regarding integrability of Lie algebroid morphisms, we prove the following lemma.
\begin{lemm}\cite{moerdijk_mrcun_2003}\label{lemm: covering projections leaves of connection}
Let $\groupoid$ be a Lie groupoid, let $\pi: P \ra M$ be a principal $\group$-bundle with moment map $\mu$, let $\mathcal{F}$ be a foliation of $M$, and let $\mathcal{H}$ be a $\mathcal{F}$-partial connection on $P$. If $\mathcal{H}$ is flat, then each leaf $\widetilde L$ of $\mathcal{H}$ comes with a map $\widetilde L \ra L$, where $L$ is a leaf of $\mathcal{F}$, and this map is a covering projection
\end{lemm}
\begin{proof}
Since $\mathcal{H}$ is a subbundle of $T\pi^*\mathcal{F}$, the restriction to a leaf $\widetilde{L}$ of the map $\pi: P \ra M$ becomes a local diffeomorphism $\widetilde{L} \ra L$. By the condition that $d\mu(\mathcal{H})=0$, we see that $\widetilde{L}$ is contained in some fiber of $\mu$, say $\widetilde{L} \subset \mu^{-1}(x)$. Then
\[\widetilde{L}_x \coloneqq \{g \in \group_x \mid \widetilde{L}g = \widetilde{L}\}\]
is well-defined and defines a (Lie) group (with the discrete topology). The action $\widetilde{L}_x$ on $\widetilde{L}$ is obviously free, and, for all $y \in \widetilde{L}$, we can find an open set $U \ni y$ of $\widetilde{L}$ such that $gU \cap U = \emptyset$ for all $\id_x \neq g \in \widetilde{L}_x$. Since a fiber of $\pi$ is naturally isomorphic to an isotropy group, we see that $\widetilde{L} \ra L$ is precisely this quotient projection (up to identification of $L$ with $\widetilde{L}/\widetilde{L}_x$). This proves the statement. 
\end{proof}
\begin{prop}\cite{moerdijk_mrcun_2003}\label{prop: Lie groupoid morphism integration source-simply connected}
Let $\groupoid$ and $\grouptwoid$ be two Lie groupoids with Lie algebroids $\algebroid$ and $C \ra \basetwo$, respectively. Assume that $\group$ is $\source$-connected and $\source$-simply connected. If $\varphi: \algebr \ra C$ is a morphism of Lie algebroids, then there exists a unique morphism of Lie groupoids $F: \group \ra \grouptwo$ that integrates $\varphi$.
\end{prop}
\begin{proof}
Denote by $f: \base \ra \basetwo$ the base map of $\varphi$ and let $P \coloneqq \group \tensor[_{f \circ \target_{\group}}]{\times}{_{\target_{\grouptwo}}} \grouptwo$ be the canonical principal $\grouptwo$-bundle over $\group$ with moment map $\mu=\target_\grouptwo \circ \pr_2$ and multiplication map 
\[\nu: P \tensor[_{\mu}]{\times}{_{\target_{\grouptwo}}} \grouptwo \ra P \textnormal{ given by } (g,h,k) \mapsto (g,\mult_\grouptwo(h,k)).\]
The idea is now as follows: out of $\varphi$, we will construct a flat $\mathcal{F}_{\source_\group}$-connection of $P$. Then, the foliation that this connection induces on $P$, induces covering projections of the leaves of the foliation to the leaves of $\mathcal{F}_{\source_\group}$ (by Lemma \ref{lemm: covering projections leaves of connection}). These projection maps are connected covering spaces, and by the assumption that $\group$ is $\source$-simply connected, these projections are diffeomorphisms. The inverse of these maps then assemble into a smooth map $\group \ra P$, which we can post-compose with $\pr_2: P \ra \grouptwo$ to obtain the desired map.

The $\mathcal{F}_{\source_\group}$-partial connection $\mathcal{H}$ on $P$ is defined as
\[\mathcal{H}_{(g,k)} \coloneqq \{\left(dR_g(\id_{\target_\group(g)})\xi, dR_k(\id_{\target_\group(g)})\varphi_{\target_\group(g)}(\xi)\right) \mid \xi \in \algebr_{\target_\group(g)}\}.\]
To see that this defines a $\mathcal{F}_{\source_\group}$-partial connection, observe first that
\[d\mu(g,h)(dR_g(\id_{\target_\group(g)})\xi,dR_k(\id_{\target_\group(g)})\varphi_{\target_\group(g)}(\xi)) = d(\target_\grouptwo \circ R_k)(\id_{\target_\group(g)})\varphi_{\target_\group(g)}\xi=0,\]
because $\target_\grouptwo \circ R_k$ is a constant map. Similarly, $\mathcal{H}_{(g,hk)} = d\nu_k(g,h)\mathcal{H}_{(g,h)}$ follows by an application of the chain rule. The leaves of $\pi^*\mathcal{F}$ are given by $\source^{-1}(x) \tensor[_{f \circ \target_\group}]{\times}{_{\target_{\grouptwo}}} \grouptwo$. Since $R_g$ and $R_k$ (in resp. $\group$ and $\grouptwo$) are diffeomorphisms, we see that the last requirement for $\mathcal{H}$ to be a $\mathcal{F}_{\source_\group}$-partial connection is also satisfied. The subbundle $\mathcal{H}$ is involutive, hence integrable, precisely because $\varphi$ is a Lie algebroid morphism, so the connection is flat. By Lemma \ref{lemm: covering projections leaves of connection}, we obtain covering projections $\widetilde{L} \ra L$ for each source fiber $L$, which are diffeomorphisms because $\group$ is $\source$-simply connected. In particular, if we restrict our attention to the leaves $\widetilde{L}$ corresponding to the points $(\id_x,\id_{\varphi(x)})$, the inverse maps $F_L: L \ra \widetilde{L}$ assemble into a map $F: \group \ra P$. To see that this map is smooth, note that, when restricted to $\identity(\base)$, it is the map $\id_x \mapsto (\id_x,\id_{\varphi(x)})$. The extension of this map using holonomy of paths is the map $\group \ra P$. More precisely, $\identity(\base)$ is a transversal for the foliation $\mathcal{F}_{\source_\group}$. So, we can define the map $F: \group \ra P$ as follows: if $g \in L = \source^{-1}(x)$, then, via a path to $\id_x$ (note: $\group$ is $\source$-connected), we obtain a holonomy diffeomorphism from a transversal $S \ni g$ to a transversal $S'$ lying in $\identity(\base)$. This diffeomorphism can be realised as a smooth map $H(1,\cdot)$, where $H: [0,1] \times S \ra \group$ is a smooth map with the following properties:
\[H(0,\cdot) = \iota_S \textnormal{, } H(1,\cdot)=\iota_{S'} \circ \hol^{S',S}(\gamma) \textnormal{ and for all } s \in S, \textnormal{ } \gamma_s \coloneqq H(\cdot,s) \textnormal{ is a path in the leaf } L_s,\]
with $\iota_S: S \hookrightarrow M$ (resp. $\iota_{S'}: S' \hookrightarrow M$) the inclusion maps. Since $S'$ maps, under $F|_{\identity(\base)}$, to a transversal $T'$, we can use the path $F_L \circ \gamma$ to obtain again such a holonomy diffeomorphism mapping onto a transversal $T \ni F_L(g)$. There is a unique such transversal which is given by the collection of endpoints of $F_{L_s} \circ \gamma_s$, so then the composition of maps $S \ra S' \ra T' \ra T$ is given by the restriction of the map $F$. This shows that the map $F$ is smooth.
By construction, the map $\pr_2 \circ F$ is the required morphism of groupoids.
\end{proof}
We end the section with showing that the bijective correspondence between Lie algebroids and vector bundles equipped with a differential can be extended to an isomorphism of categories. We first prove the statement in case the Lie algebroid morphism is an inclusion of Lie algebroids.
\begin{lemm}\cite{meinrenken}\label{lemm: morphism of diff vector bundle is equivalent to morphism of Lie algebroids}
A vector subbundle $\subalgebroid$ of a Lie algebroid $\algebroid$ is a Lie algebroid if and only if $\Omega^\bullet(\subalgebr)$ inherits a (necessarily unique) differential $d_\subalgebr$ such that $\Omega^\bullet(\algebr) \ra \Omega^\bullet(\subalgebr)$ is a cochain map. 
\end{lemm}
\begin{proof}
Denote by $\iota$ the inclusion $\subalgebr \hookrightarrow \algebr$. Suppose first that $\subalgebr$ is a subalgebroid of $\algebr$. Then $\subalgebr$ comes with the de Rham-differential $d_\subalgebr$. By Remark \ref{rema: Lie subalgebroid sections} we only have to show that, for all $\omega \in \Omega^\bullet(\algebr)$, the relation $d_\subalgebr \circ \iota^*(\omega) = \iota^* \circ d_\algebr(\omega)$ holds when applied to the restriction to $\subbase$ of sections $\alpha^1,\dots,\alpha^k$ in $\Gamma(\algebr,\subalgebr)$. Now, since $\subalgebr$ is a subalgebroid of $\algebr$, we have $\bracket{\alpha^i,\alpha^j}|_{\subbase} = [\alpha^i|_{\subalgebr},\alpha^j|_{\subalgebr}]_\subalgebr$, and $\anchor_\algebr(\alpha^i) \in \mathfrak{X}(\base)$ is $\iota$-related to $\anchor_\subalgebr(\alpha^i|_{\subbase})$. The result now follows by using the definition of the de Rham-differentials $d_\algebr$ and $d_\subalgebr$.

Conversely, assume $\subalgebr$ inherits a differential $d_\subalgebr$ such that $\iota^*$ is a cochain map. Since $\subalgebr$ now is equipped with a differential, it has a Lie algebroid structure for which $d_\subalgebr$ is the de Rham-differential. Let $\alpha \in \Gamma(\algebr,\subalgebr)$. Then, clearly, $\iota^* \circ \iota_\alpha = \iota_{\alpha|_{\subbase}} \circ \iota^*$ (where $\iota_\alpha$ and $\iota_{\alpha|_{\subbase}}$ denote interior derivatives). We now use Cartan Calculus identities: by Lemma \ref{lemm: Cartan's magic formula}, we have that 
\[\iota^* \circ \mathcal{L}_{\anchor_\algebr}(\alpha) = \mathcal{L}_{\anchor_\subalgebr}(\alpha|_{\subbase}) \circ \iota^*.\] 
It follows that $\subalgebr$ is an anchored subbundle of $\algebr$. Indeed, if $f \in C^\infty(\base)$ vanishes along $\subbase$, then, for all $\alpha \in \Gamma(\algebr,\subalgebr)$, we have
\[\iota^*\anchor_\algebr(\alpha) = \iota^*(\iota_\alpha \circ d_\algebr f) = \iota_{\alpha|_{\subbase}} \circ d_\subalgebr \circ \iota^*f = 0,\]
where we used Lemma \ref{lemm: Cartan's magic formula} in the first equality. This shows that $\anchor_\algebr(\alpha) \in \mathfrak{X}(\base,\subbase)$, so $\anchor_\algebr: \algebr \ra T\base$ maps $\subalgebr$ into $T\subbase$, i.e. $\subalgebr$ is an anchored subbundle of $\algebr$.
Lastly, if $\alpha^1,\alpha^2 \in \Gamma(\algebr,\subalgebr)$, we have, by Lemma \ref{lemm: interior product of bracket}, that 
\[\iota^* \circ \iota_{\bracket{\alpha^1,\alpha^2}} = \iota_{[\alpha^1|_{\subbase},\alpha^2|_{\subbase}]_\subalgebr} \circ \iota^*.\]
In particular, if $\omega \in \Omega^\bullet(\algebr)$ vanishes along $\subbase$, then $\iota_{\bracket{\alpha^1,\alpha^2}}\omega$ vanishes along $\subbase$ as well. From this we see that $\bracket{\alpha^1,\alpha^2} \in \Gamma(\algebr,\subalgebr)$, so $\Gamma(\algebr,\subalgebr) \subset \Gamma(\algebr)$ is a Lie subalgebra.
\end{proof}

\begin{prop}\cite{meinrenken}\label{prop: morphism of diff vector bundle is equivalent to morphism of Lie algebroids}
Let $\varphi: \algebr \ra \subalgebr$ be a morphism of vector bundles between Lie algebroids $\algebroid$ and $\subalgebroid$. Then $\varphi$ is a Lie algebroid morphism if and only if the map $\varphi^*: \Omega^\bullet(\subalgebr) \ra \Omega^\bullet(\algebr)$ is a cochain map.
\end{prop}
\begin{proof}
Recall first that we can construct the graded differential $C^\infty(\subbase)$-algebra 
\[\Omega^\bullet(\algebr) \otimes \Omega^\bullet(\subalgebr) = \Omega^\bullet(\algebr) \otimes_{C^\infty(\subbase)} \Omega^\bullet(\subalgebr) \coloneqq \bigoplus_{i+j=n} \Omega^i(\algebr) \otimes_{C^\infty(\subbase)} \Omega^j(\subalgebr)\]
equipped with the differential
\[d_{\otimes}(\omega_\algebr \otimes \omega_\subalgebr) \coloneqq d_{\algebr}(\omega_\algebr) \otimes \omega_\subalgebr + (-1)^{|\omega_\algebr|} \omega_{\algebr} \otimes d_\subalgebr(\omega_\subalgebr), \textnormal{ where } \omega_\algebr \in \Omega^\bullet(\algebr)_{\textnormal{hom}}, \omega_\subalgebr \in \Omega^\bullet(\subalgebr)_{\textnormal{hom}}.\]
By a direct computation, it is readily verified that the map
\[\Omega^\bullet(\algebr) \otimes \Omega^\bullet(\subalgebr) \xra{\sim} \Omega^\bullet(\algebr \times \subalgebr) \textnormal{ given by } \omega_\algebr \otimes \omega_\subalgebr \mapsto \pr_1^*\omega_\algebr \wedge \pr_2^*\omega_\subalgebr\]
(we use here the identification $\algebr \times \subalgebr \cong \pr_1^*\algebr \oplus \pr_2^*\subalgebr$) is an isomorphism of cochain complexes. Now, the map
\[\iota^*: \Omega^\bullet(\algebr) \otimes \Omega^\bullet(\subalgebr) \ra \Omega^\bullet(\algebr)\]
that the vector bundle morphism
\[\iota: \algebr \xra{\Delta_\algebr} \algebr \times \algebr \xra{\id_\algebr \times \varphi} \algebr \times \subalgebr\]
induces, is the multiplication map $\omega_\algebr \otimes \omega_\subalgebr \mapsto \omega_\algebr \cdot \varphi^*(\omega_\subalgebr)$, because the diagonal morphism $\Delta_\algebr$ induces the multiplication map on $\Omega^\bullet(\algebr)$. The map $\iota$ is the inclusion of the graph $\gr \varphi \cong \algebr$ into $\algebr \times \subalgebr$, so it is a subalgebroid (i.e. $\varphi$ is a morphism of Lie algebroids) if and only if $\iota^*$ is a cochain map. The latter statement means that
\[d_\algebr(\omega_\algebr \cdot \varphi^*(\omega_\subalgebr)) = d_\algebr(\omega_\algebr) \cdot \varphi^*(\omega_\subalgebr) + (-1)^{|\omega_\algebr|} \omega_\algebr \cdot \varphi^*(d_\subalgebr(\omega_\subalgebr)),\]
which holds if and only if it holds for $\omega_\algebr=1$, i.e. when $\varphi$ is a cochain map. This proves the statement.
\end{proof}
\subsection{Isotropy Lie algebras and orbits}\label{sec: Isotropy algebras and orbits}
The Lie algebroid analogue of the isotropy Lie groups are isotropy Lie algebras. The Lie algebroid also induces a partition of the base manifold into orbits. This partition carries the structure of a so-called \textit{singular foliation}, which should be compared to foliations as defined in Section \ref{sec: Foliations and the monodromy and holonomy groupoids}. We will go into the theory of singular foliations in the next section. We fix a Lie algebroid $\algebroid$.
\begin{defn}\cite{ruimarius}\label{defn: Isotropy Lie algebra}
Let $x \in \base$. Then the \textit{isotropy Lie algebra at $x$} is the Lie algebra 
\[\mathfrak{g}_x = \mathfrak{g}_x(\algebr) \coloneqq \ker \anchor(x),\]
with Lie bracket given by
$[\alpha(x),\beta(x)]_{\mathfrak{g}_x} \coloneqq [\alpha,\beta]_\algebr(x)$, where $\alpha,\beta \in \Gamma(\algebr)$ with $\alpha(x),\beta(x) \in \mathfrak{g}_x$.
\end{defn}
\begin{rema}\cite{ruimarius}\label{rema: Isotropy Lie algebra}
Notice that the Lie bracket on $\mathfrak{g}_x$ is well-defined: by the Leibniz-rule, 
\begin{equation}\label{eq: Leibniz rule for isotropy}
    [\alpha,f\beta]_\algebr(x) = f(x)[\alpha,\beta]_\algebr(x) \textnormal{ for all } f \in C^\infty(\base).
\end{equation} %\textnormal{ and } \alpha,\beta \in \Gamma(\algebr) \textnormal{ with } \alpha(x),\beta(x) \in \mathfrak{g}_x.
So, if $\alpha,\beta,\sigma,\tau \in \Gamma(\algebr)$ with $\alpha(x)=\sigma(x),\beta(x)=\tau(x) \in \mathfrak{g}_x$, then 
\[\bracket{\alpha,\beta}-\bracket{\sigma,\tau}=\bracket{\alpha,\beta-\tau} - \bracket{\tau,\alpha-\sigma},\]
(we use here the same trick as in Proposition \ref{prop: Lie bracket on A defines a Lie bracket on sheaf of A}). This reduces the problem to showing that, if $\beta(x)=0$, then $\bracket{\alpha,\beta}(x)=0$. From Proposition \ref{prop: Lie bracket on A defines a Lie bracket on sheaf of A}, it follows that we may assume that $\alpha$ and $\beta$ are local sections. If $\{\alpha^i\}$ is a local frame of $\algebr$, and $\beta=\textstyle\sum_i f_i\alpha^i$, then $f_i(x)=0$ for all $i$, so from \eqref{eq: Leibniz rule for isotropy} we see that $\bracket{\alpha,\beta}(x)=0$, as required.

To see that $[\cdot,\cdot]_{\mathfrak{g}_x}$ defines a Lie bracket, notice that \eqref{eq: Leibniz rule for isotropy} also shows that the Lie bracket is $\mathbb{R}$-bilinear. The skew-symmetry and the Jacobi-identity are immediate consequences of $\bracket{\cdot,\cdot}$ being a Lie bracket on $\Gamma(\algebr)$. This shows that, as claimed, $\mathfrak{g}_x$, equipped with the bracket $[\cdot,\cdot]_{\mathfrak{g}_x}$, is a Lie algebra.
\end{rema}
In fact, if $\algebroid$ integrates to $\groupoid$, then the Lie algebra of $\group_x$ is isomorphic to $\mathfrak{g}_x$.

We will define orbits of a Lie algebroids by using the following notion.
\begin{defn}\cite{ruimarius}\label{defn: A-path}
A path $a: [0,1] \ra \algebr$ is called an \textit{\algebr}-path if there is a path $\gamma: [0,1] \ra \base$, such that $a(t) \in \algebr_{\gamma(t)}$, and $\anchor(a(t)) = \tfrac{d\gamma}{dt}(t)$ for all $t \in [0,1]$. In this case, $a$ is said to be an \textit{$\algebr$-path above $\gamma$}.
\end{defn}
\begin{rema}\cite{ruimarius}\label{rema: A-path}
If $a: [0,1] \ra \algebr$ is a path, then there is at most one path $\gamma: [0,1] \ra \base$ for which $a$ is an $\algebr$-path $\gamma$, namely $\gamma=\pi \circ a$, where $\pi: \algebr \ra \base$ is the bundle projection. Moreover, notice that the above definition already makes sense for anchored vector bundles, and that if $\varphi: E \ra F$ is a morphism of anchored vector bundles (with base map $g$), and $e: [0,1] \ra E$ is an $E$-path above a path $\gamma$ ($=\pi_E \circ e$), then $f \coloneqq \varphi \circ e: [0,1] \ra F$ is an $F$-path above the path $\delta \coloneqq \pi_F \circ f = g \circ e$. Indeed, this follows from the sequence of equalities
\begin{align*}
    \anchor_F(f(t)) &= \anchor_F(\varphi \circ e(t)) \\
    &= dg \circ \anchor_E(e(t)) \\
    &= dg \circ \frac{d\gamma}{dt}(t) = \frac{d\delta}{dt}(t),
\end{align*}
where we used in the second equality that $\varphi$ is a morphism of anchored vector bundles.
\end{rema}
We can now define an equivalence relation on $\base$: 
\begin{equation}\label{eq: equivalence relation of A-paths}
    x \sim y \textnormal{ if and only if there is an $\algebr$-path } a \textnormal{ such that } \pi \circ a \textnormal{ is a path from } x \textnormal{ to } y.
\end{equation} 
The equivalence classes of this equivalence relation are immersed submanifolds of $\base$, but we will wait with giving a proof. The reason is that we will use the so-called \textit{local splitting theorem for Lie algebroids}, and we will give a proof of this statement using the theory of deformation to the normal cones.
\begin{defn}\cite{ruimarius}\label{defn: orbits of a Lie algebroid}
Let $x \in \base$. Denote by $\orbit_x \subset \base$ the equivalence class of $x$ with respect to the equivalence relation defined in (\ref{eq: equivalence relation of A-paths}). Then $\orbit_x$ is called the \textit{orbit} of $x$.
\end{defn}
In case $\algebroid$ is integrable, we can also describe the orbits using Lie groupoids (see also Proposition \ref{prop: s-(simply )connected groupoid of integrable algebroid}).
\begin{prop}\cite{ruimarius}\label{prop: orbits of Lie groupoid are orbits of Lie algebroid}
Let $\groupoid$ be an $\source$-connected Lie groupoid integrating $\algebroid$. Then the orbits of $\groupoid$ coincide with the orbits of $\algebroid$.
\end{prop}
\begin{proof}
Let $x \in \base$. We will write $\orbit_x^\group$ for the orbit of $x$ with respect to $\groupoid$, and $\orbit_x^\algebr$ for the orbit of $x$ with respect to $\algebroid$. We have to prove that \[\orbit_x^\group=\orbit_x^\algebr.\]
We will first show that $\orbit_x^\group \subset \orbit_x^\algebr$; let $y \in \orbit_x^\group$. This means that there exists a $g \in \source^{-1}(x)$ such that $\target(g)=y$. Since $\source^{-1}(x)$ is connected, we can find a path $u: [0,1] \ra \source^{-1}(x)$ from $\id_x$ to $g$, which, by post-composing with $\target$, descends to a path $\gamma$ in $\orbit_x^\group$. In particular, 
\[a: [0,1] \ra \algebr \textnormal{ given by } t \mapsto dR_{u(t)^{-1}}(u(t))\frac{du}{dt}(t)\]
defines a path in $\algebr$. It is even an $\algebr$-path above $\gamma$, because $\target \circ R_{u(t)^{-1}}=\target$ for all $t \in [0,1]$, so $y \in \orbit_x^\algebr$. To show that $\orbit_x^\algebr \subset \orbit_x^\group$, we have to show that every $\algebr$-path $a: [0,1] \ra \algebr$, with base path $\gamma$ from $x$ to $y$, can be written as above. To do this, observe that we can find a time-dependent section $\alpha: [0,1] \times \base \ra \algebr$ such that $\alpha(t,\gamma(t))=a(t)$. The flow of the time-dependent vector field $\sigma_\alpha$, the corresponding time-dependent right-invariant vector field defined by $(t,g) \mapsto dR_g(\id_{\target(g)})\alpha(t,\gamma(t))$, i.e. the local diffeomorphisms $\varphi^{t,0}_{\sigma_\alpha}: \group \ra \group$ determined as the solution of
\[\frac{d}{dt} \varphi^{t,0}_{\sigma_\alpha}(p) = \sigma_\alpha(t,\varphi^{t,0}_{\sigma_\alpha}(p)); \quad \varphi^{0,0}_{\sigma_\alpha}(p)=p,\]
gives rise to the smooth map $u: [0,1] \ra \source^{-1}(x)$ given by $u(t) \coloneqq \varphi^{t,0}_{\sigma_\alpha}(x)$. That the diffeomorphism $u$ maps into $\source^{-1}(x)$ follows from the fact that $\sigma_\alpha(t,\cdot) \in \Gamma(\ker d\source)$ (and applying the chain rule to the above ODE, but with $\source \circ \varphi^{t,0}_{\sigma_\alpha}$ everywhere instead of $\varphi^{t,0}_{\sigma_\alpha}$).
That $u$ is defined on all of $[0,1]$ is because $\sigma_\alpha$ is right-invariant: this follows from the fact that, for any smooth path $g: [0,1] \ra \source^{-1}(x)$ with $\target \circ g = \gamma$ and $g(0)=\id_x$, we have that $u(t) \cdot g(t)^{-1}$ is also a solution of the above initial value problem. Using right-invariancy once again proves the statement:
\begin{align*}
    dR_{u(t)^{-1}}(u(t))\frac{du}{dt}(t) &= dR_{u(t)^{-1}}(u(t))\sigma_\alpha(t,u(t))\\
    &=\sigma_\alpha(t,R_{u(t)^{-1}}u(t))\\
    &=\sigma_\alpha(t,\id_x)=\alpha(t,\gamma(t))=a(t),
\end{align*}
so we see that, indeed, every $\algebr$-path is of the form $t \mapsto dR_{u(t)^{-1}}(u(t))\tfrac{du}{dt}(t)$ for some path $u: [0,1] \ra \source^{-1}(x)$.
%notice that the equation
%\[\frac{du}{dt}(t) = dR_{u(t)}(u(t)^{-1})a(t)\]
%defines an initial value problem (at least locally) with initial value $u(0)=\id_x$. Since such a solution is unique if it exists, we see that $a(t)=dR_{u(t)^{-1}}(u(t))\tfrac{du}{dt}(t)$, as required. 
This proves the statement.
\end{proof}
\begin{coro}\cite{ruimarius}\label{coro: orbits of groupoid^0 coincide with Lie algebroid orbits}
Let $\groupoid$ be a Lie groupoid and let $\algebroid$ be its Lie algebroid. Then the orbits of $\group^0 \rra \base$ coincide with the orbits of $\algebroid$.
\end{coro}
As already mentioned before (see the remark made before Remark \ref{rema: Atiyah Lie algebroids are transitive}), not all transitive Lie algebroids arise as an Atiyah Lie algebroid. From the above corollary, we can give the following standard example:
\begin{rema}\cite{meinrenken}\label{rema: transitive Lie algebroids and Atiyah algebroids}
Let $M$ be a connected, simply connected manifold equipped with a closed $2$-form $\omega$. We can exactly describe for which $\omega$ the Lie algebroid
\[\algebr \coloneqq TM \times \mathbb{R} \ra M; \quad \anchor = \pr_1,\]
with Lie bracket given by 
\[\bracket{(\sigma_1,f_1),(\sigma_1,f_2)} \coloneqq ([\sigma_1,\sigma_2], \omega(\sigma_1,\sigma_2) + \sigma_1(f_2) - \sigma_2(f_1)),\]
is integrable. To do this, assume it is integrable to a Lie groupoid $\group \rra M$. Then $\group^0 \rra M$ is also an integration, and this groupoid has the same orbits as $\algebr$ (see Proposition \ref{prop: orbits of Lie groupoid are orbits of Lie algebroid}), so $\group^0$ is transitive. We replace $\group$ with $\group^0$. Then, $\group \cong P \otimes_G P$ for some principal $G$-bundle $P \ra M$ (see Remark \ref{rema: transitive groupoids are in one to one correspondece with Gauge groupoids}). If we pick $x \in M$, then we can realise $G$ as the Lie group $\group_x$. Since $\mathfrak{g}_x$ is the Lie algebra of $G=\group_x$, and $\mathfrak{g}_x$ is obviously one-dimensional, $G$ will be a $1$-dimensional Lie group. Since $M$ is connected and simply connected, we may assume that the fibers of $P$ are connected, i.e. $G$ is connected. In particular, $G \cong \mathbb{R}/\pi_1(G)$ (so $G \cong \mathbb{R}$ or $G \cong S^1$), because the exponential map $\textnormal{exp}: \mathbb{R} \cong \mathfrak{g} \ra G$ is surjective. The splitting $TM \ra A$ given by $(x,v) \mapsto (x,v,0)$ gives rise to a (principal) connection $\nabla$ on $P$, and it is readily verified that $\omega$ is its curvature $2$-form. 

Now, recall that, as for connections on vector bundles, connections on principal bundles $Q \ra X$ give rise to a notion of parallel transport, and therefore also of holonomy: if $\gamma: [0,1] \ra X$ is a loop based at $x$, and $q$ lies in the fiber over $x$, then $\gamma$ lifts to a unique path $\widetilde\gamma: [0,1] \ra Q$ with $\widetilde{\gamma}(0)=q$. The element $\widetilde{\gamma}(1) \in Q$ still lands in the fiber over $x$, so for a unique $g \in G$, which we denote by $\hol(\gamma)$, we have $p \cdot \hol(\gamma) = \widetilde{\gamma}(1)$. In our case, we can calculate the holonomy very explicitly: a loop $\gamma$ in $M$ is contractible, so by viewing $\gamma$ as a map $S^1 \ra M$, it can be extended to a smooth map $f: D^2 \ra M$ (where $D^2$ is the disc of radius $1$ in $\mathbb{R}^2$). We claim that
\[\int_{D^2} f^*\omega \mod \pi_1(G) = \hol(\gamma).\]
Indeed, the pullback principal bundle $\gamma^*P$ over $S^1$ carries the connection $\gamma^*\nabla$, which has curvature $2$-form $\gamma^*\omega$. The loop $\id_{S^1}: S^1 \ra S^1$, i.e. the pullback of $\gamma: S^1 \ra M$ to $S^1$, has the same holonomy as $\gamma$ by checking the relation between these elements via the $G$-equivariant map $\gamma^*P \ra P$. Since $G$ is abelian, the holonomy of $\id_{S^1}$ is given by $\exp \int_{S^1} \gamma^*\nabla$, so the result follows by using that $\omega$ is the curvature $2$-form of $\nabla$ and an application of Stokes' theorem.

If $f': D^2 \ra M$ is another extension of $\gamma$, then $f$ and $f'$ agree on $\partial D^2$, so the difference $\int_{D^2} f^*\omega - \int_{D^2} (f')^*\omega$ is the integral 
\[\int_{S^2} \delta^*\omega \mod \pi_1(G),\]
where $\delta: S^2 \ra M$ is the map $f$ on the upper hemisphere, and $f'$ on the lower hemisphere (here, we identified $D^2$ with the upper and lower hemispheres of $S^2$, respectively). In particular, the subgroup $\Lambda \subset \mathbb{R}$ given as the image of the group homomorphism
\[\pi_2(M) \ra \mathbb{R} \textnormal{, given by } [\delta] \mapsto \int_{S^2} \delta^*\omega\]
has to lie in $\pi_1(G)$, which is a discrete group, so $\Lambda$ must be a discrete subgroup of $\mathbb{R}$ as well. An explicit counterexample is now given by taking $M=S^2 \times S^2$ and $\omega = \omega_{\textnormal{std}} \oplus \lambda \omega_{\textnormal{std}}$, where $\omega_{\textnormal{std}}$ is the standard area form on $S^2$ and $\lambda \in \mathbb{R} \setminus \mathbb{Q}$. The group $\Lambda$ is then $\mathbb{Z} + \lambda\mathbb{Z}$, and $\mathbb{R} = \overline{\mathbb{Z} + \lambda\mathbb{Z}}$.
\end{rema}

\subsection{Singular foliations}\label{sec: Singular foliations}
We have seen that a foliation of a manifold is a partition into immersed submanifolds of the same dimension satisfying a locally triviality condition (see II in Proposition \ref{prop: equivalent definitions of foliations}). We have also seen in Section \ref{sec: Isotropy groups and orbits} that Lie groupoids admit a partition into immersed submanifolds, namely the partition into its orbit spaces, and in the former section we saw that Lie algebroids admit a partition into immersed submanifolds as well (except for the proof that the orbits are immersed submanifolds). We will see in this section that the partition of the base manifold corresponding to a Lie algebroids carries the structure of a so-called \textit{singular foliation}. Recall that a Lie algebroid with injective anchor map can be seen as a foliation, and it is readily verified that the leaves of this foliation coincide with the orbits of the Lie algebroid. So, we can try to find an appropriate generalisation of a foliation to general Lie algebroids%, which will be the notion of a \textit{singular foliation}
. While every Lie algebroid gives rise to a singular foliation, not every singular foliation comes from a Lie algebroid \cite{androu2}. Moreover, a partition into immersed submanifolds, which can be given the structure of a singular foliation, in general admits many such singular foliation structures. This section is based on \cite{androu,lavau2018lie}.
\begin{term}\label{term: foliations to regular foliations}
We will still call the foliations defined in Section \ref{sec: Foliations and the monodromy and holonomy groupoids} \textit{foliations}, and we will call the singular foliations defined below \textit{singular foliations}. If we want to be precise, we will call the former type of foliation a \textit{regular foliation}.
\end{term}
We will define (singular) foliations by generalising regular foliations viewed as integrable distributions. However, we will take the slightly different approach of sheaves. 
\begin{rema}\cite{androu}\label{rema: Serre Swan and foliations}
Recall that the assignment [$E \ra M$ a vector bundle $\mapsto \Gamma(E)$] can be extended to an equivalence of categories between vector bundles over a manifold $M$ and projective, locally finitely generated $C^\infty(M)$-modules (the Serre-Swan theorem). Moreover, subbundles $F$ of $E$, over the same base as $E$, can be recovered via this equivalence as follows: for all $x \in M$, $F_x$ is naturally isomorphic to $\Gamma(F)/I_x$, where $I_x$ is the vanishing ideal of $x$ in $C^\infty(M)$, via the evaluation map 
\[\textnormal{ev}_x: \Gamma(F)/I_x \ra F_x \textnormal{ given by } f\mod I_x \mapsto f(x).\]
Moreover, if we replace the modules in this discussion with the modules obtained by only allowing sections that have compact support, then the same statements hold. Therefore, a regular foliation, in terms of modules, can be described as a projective, locally finitely generated, and involutive submodule $\mathcal{F}$ of $\mathfrak{X}_c(M)$, where $\mathfrak{X}_c(M) \subset\mathfrak{X}(M)$ consist of the vector fields of $M$ with compact support, such that the evaluation maps $\textnormal{ev}_x: \mathcal{F}/I_x \ra T_xM$ are injective (mapping onto the tangent space of a leaf). 
\end{rema}
We now define singular foliations as follows.
\begin{defn}\cite{androu}\label{defn: singular subalgebroid}
Let $\algebroid$ be a Lie algebroid. We call a locally finitely generated and involutive $C^\infty(\base)$-submodule $\mathcal{F}$ of $\Gamma_c(\algebr) \coloneqq \{\alpha \in \Gamma(\algebr) \mid \alpha \textnormal{ has compact support}\}$ a \textit{singular subalgebroid} of $\algebr$. In case $A=T\base$, such a submodule of $\mathfrak{X}_c(\base)$ is called a \textit{singular foliation} of $\base$, and $\base$ is called a \textit{singularly foliated manifold}.
\end{defn}
\begin{rema}\cite{androu}\label{rema: singular foliation of a Lie algebroid}
To be precise, locally finitely generated in this case means that we can find an open cover $\{U_i\}$ of $\base$ such that for each $U_i$ we can find finitely many sections in
\[\{\sigma \in \mathfrak{X}(U) \mid f\sigma \in \iota_i^*\mathcal{F} \textnormal{ for all } f \in C_c^\infty(U)\},\]
where $\iota_i: U_i \hookrightarrow \base$ is the inclusion, that generate $\iota_i^*\mathcal{F}$. Here, 
\[\iota_i^*\mathcal{F} \coloneqq \{\sigma|_U \in \mathfrak{X}(U) \mid \sigma \in \mathcal{F} \textnormal{ with supp } \sigma \subset U\}.\]
Observe that we require $\mathcal{F}$ to consist of only compactly supported vector fields, but the generators may not be compactly supported (in particular, generators of $\mathcal{F}$ are typically not elements of $\mathcal{F}$).
\end{rema}
\begin{exam}\cite{androu}\label{exam: Lie algebroid foliation}
It is readily verified that a Lie algebroid $\algebroid$ determines the singular foliation $\anchor(\Gamma_c(\algebr))$ of $\base$. 
\end{exam}
As for regular foliations, for a singularly foliated manifold $(M,\mathcal{F})$, the assumption that the submodule of $\mathfrak{X}_c(M)$ is involutive ensures that the manifold is partitioned into immersed submanifolds (which may have varying dimension now). In fact, for all $x \in M$, we have an exact sequence
\[0 \ra \mathfrak{g}_x \ra \mathcal{F}/I_x\mathcal{F} \xra{\textnormal{ev}_x} T_xL \ra 0,\]
where $L$ is a leaf, and
\[\mathfrak{g}_x \coloneqq \mathcal{F}(x)/I_x\mathcal{F},\]
with $\mathcal{F}(x) \coloneqq \{\sigma \in \mathcal{F} \mid \sigma(x)=0\}$. The vector space $\mathfrak{g}_x$ inherits a Lie bracket from $\mathfrak{X}_c(M)$ (This works in the same way as for Lie algebroids). It is clear that, in case $\mathcal{F}$ is the singular foliation of a Lie algebroid, $\mathfrak{g}_x$ is isomorphic to the isotropy Lie algebra at $x$.  

We can think about the leaves of a singular foliation as follows: if $\sigma \in \mathfrak{X}_c(M)$ is a compactly supported vector field, then its flow $\varphi_{\sigma}^t$ exists for all time, and therefore we obtain a diffeomorphism
\[\exp \sigma: M \ra M, \textnormal{ given by } x \mapsto \varphi_\sigma^1(x).\]
One can show that $\exp \sigma$ respects the singular foliation of a singularly foliated manifold $M$ (i.e. the singular foliation on the source induces the same singular foliation on the target; see \cite{alfonso}). Therefore, the subgroup of Diff $M$ generated by elements of the form $\exp \sigma$ has orbits equal to the leaves of the singular foliation. One can show that the leaves of the singular foliation induced by a Lie algebroid are precisely the orbits described in the former section.
\begin{exam}\label{exam: foliation of R^2}
As a simple example, consider the foliation of $\mathbb{R}^2 \setminus \{0\}$ by circles, i.e. the foliation determined by the submersion $\mathbb{R}^2 \setminus \{0\} \ra \mathbb{R}_{>0}$ given by $(x,y) \mapsto x^2 + y^2$. By ``adding $0$ back in as a leaf'' we obtain a singular foliation $\mathcal{F}$ of $\mathbb{R}^2$. Namely, we can take $\mathcal{F}$ to be the $C^\infty_c(\mathbb{R}^2)$-span of the vector field
\[\sigma(x,y) \coloneqq -y\frac{\partial}{\partial x} + x\frac{\partial}{\partial y}.\]
%For later reference, notice that this vector field is linear. 
\end{exam}
%To end this section, we will show that the leaves of the singular foliation induced by a Lie algebroid are precisely the orbits of the Lie algebroid as described in the former section.
%\begin{rema}\label{rema: orbits of singular foliation = algebroid orbits}
%Consider the singular foliation $\anchor(\Gamma_c(\algebr))$ induced by a Lie algebroid $\algebroid$. If $\{\alpha^i\}$ is a local frame for $\algebr$ on an open subset $U \subset \base$, then the collection of vector fields $\{\anchor(\alpha^i)\}$ (which is usually not a frame anymore) generate $\iota^*\mathcal{F}$, where $\iota: U \hookrightarrow \base$ is the inclusion. Now, if $a: [0,1] \ra \algebr$ is an $\algebr$-path, with base path $\gamma$ from $x$ to $y$, then we can find a time-dependent section $\alpha: [0,1] \times \base \ra \algebr$ such that $\alpha(t,\gamma(t))=a(t)$. Since 
%\[\anchor(a(t)) = \anchor(\alpha(t,\gamma(t))) = \frac{d\gamma}{dt}(t),\]
%we see that $\gamma$ equals the map $t \mapsto \varphi^{t,0}_{\anchor(\alpha)}(x)$. By the characterisation of leaves of a singular foliation described before, we see that $\gamma$ is contained in a leaf. Conversely, if $\alpha$ is a time-dependent section of $\algebr$ such that $\anchor(\alpha)$ is tangent to the leaves, then, the flow $\varphi^t$
%\end{rema}
\addtocontents{toc}{\protect\thispagestyle{myplain}}\newpage

\section{Deformation to the normal cone}\label{sec: deformation to the normal cone}
In this section we will make the first important step towards constructing the blow-up of a groupoid. The first approach we will take is by introducing the \textit{deformation to the normal cone} functor 
\[\DNC: \mathcal{C}^\infty_2 \ra \mathcal{C}^\infty.\]
Here, $\mathcal{C}^\infty$ is the category of (possibly non-Hausdorff) smooth manifolds and $\mathcal{C}^\infty_2$ is the following category of pairs of (possibly non-Hausdorff) smooth manifolds. The objects of $\mathcal{C}^\infty_2$ are pairs of smooth (possibly non-Hausdorff) manifolds $(\base,\subbase)$, where $\subbase \subset \base$ is a \textbf{closed embedded submanifold} (closed in the topological sense; they can be non-compact), and the morphisms between two pairs $(\base,\subbase)$ and $(\basetwo,\subbasetwo)$ are smooth maps $f: \base \ra \basetwo$ that map $\subbase$ into $\subbasetwo$; we denote such a map by 
\[f: (\base,\subbase) \ra (\basetwo,\subbasetwo).\] 
For some of the theory we will discuss, it is not necessary to assume the submanifolds are closed. In fact, for most of the theory discussed in this section, this assumption is unnecessary. However, blow-ups are only defined for closed embedded submanifolds, and it does not make much sense for more general submanifolds. For convenience, we will therefore restrict our attention to closed embedded submanifolds.

Roughly speaking, for a pair of smooth manifolds $(\base,\subbase)$, the deformation to the normal cone of $\subbase$ in $\base$ is obtained by replacing $\base \times \{0\}$ in $\base \times \mathbb{R}$ by the normal bundle $\normal(\base,\subbase)$ of $\subbase$ in $\base$. At least heuristically, this makes sense, because we can glue $\base \times \mathbb{R}^\times$ to $\normal(\base,\subbase)$ by using a tubular neighbourhood. The space $\DNC(\base,\subbase)$ carries a natural $\mathbb{R}^\times$-action, and the blow-up manifold $\blup(\base,\subbase)$ will be the quotient of $\DNC(\base,\subbase) \setminus (\subbase \times \mathbb{R})$ (where we identified the zero section $0_\subbase \subset \normal(\base,\subbase)$ with $\subbase$) by this canonical action of $\mathbb{R}^\times$. This way we effectively ``replace'' $\subbase$ with the projectivisation of the normal bundle of $\subbase$ in $\base$. If we have a pair of Lie groupoids $(\groupoid,\subgroupoid)$, we will see that, by removing a suitable ``degeneracy'' from $\DNC(\group,\subgroup)$ first, we automatically obtain a blow-up groupoid by using the functoriality of the deformation to the normal cone construction.
 
The construction on the level of manifolds is based on \cite{rouse2008schwartz}. The extension to Lie groupoids is based on \cite{2017arXiv170509588D}. 

\begin{term}\label{term: adapted chart}
Let $(\base,\subbase)$ be a pair of smooth manifolds. We will say it has \textit{dimension $(n,p)$} if $\base$ has dimension $n$ and $\subbase$ has dimension $p$ as manifolds (we will then write $q=n-p$ for the codimension of $\subbase$ in $\base$). Moreover, we will call a chart $(U,\varphi)$ on $\base$ such that for $V \coloneqq U \cap \subbase$ we have
\[\varphi(V) = \varphi(U) \cap (\mathbb{R}^p \times \{0\})\]
an \textit{adapted chart (adapted to $\subbase$)}. Observe that $\varphi$ can be seen as a morphism of pairs $(U,V) \ra (\mathbb{R}^n,\mathbb{R}^p \times \{0\})$. Lastly, instead of $\mathbb{R}^p \times \{0\} \subset \mathbb{R}^n$, we will often just write $\mathbb{R}^p \subset \mathbb{R}^n$.
\end{term}

\subsection{Short review of the normal bundle functor}\label{sec: short review of the normal bundle functor}
Whenever $(\base,\subbase)$ is a pair of smooth manifolds, we can consider the normal bundle of $\subbase$ in $\base$: we denote it by $\normal(\base,\subbase) \ra \subbase$
and it is defined as the cokernel of the inclusion $T\subbase \hookrightarrow T\base|_\subbase$, i.e.
\[\normal(\base,\subbase) = T\base|_\subbase/T\subbase \textnormal{ which has fibers } \normal_y(\base,\subbase) \coloneqq T_y\base/T_y\subbase.\]
Observe that if $f: (\base,\subbase) \ra (\basetwo,\subbasetwo)$ is a morphism of pairs (so $f$ is a smooth map $\base \ra \basetwo$ that maps $\subbase$ into $\subbasetwo$), then the bundle map $df: T\base \ra T\basetwo$ maps $T\subbase$ into $T\subbasetwo$. In particular, we obtain a vector bundle morphism $d_\normal f: \normal(\base,\subbase) \ra \normal(\basetwo,\subbasetwo)$. By the chain rule, we see that
\[\normal: \mathcal{C}^\infty_2 \ra \mathcal{C}^\infty, \textnormal{ given by } (\base,\subbase) \mapsto \normal(\base,\subbase); \quad f \mapsto d_\normal f,\]
is a functor, called the \textit{normal bundle functor}.
We will see that this functor lies at the basis of the deformation to the normal cone functor. 

\begin{rema}\label{rema: local coordinates normal bundle}
Let $(\base,\subbase)$ be a pair of smooth manifolds of dimension $(n,p)$ and let $(U,\varphi)$ be an adapted chart (adapted to $\subbase$). Observe that $(TU,d\varphi)$ is an adapted chart adapted to $TV$ (where $V \coloneqq U \cap \subbase$). For this reason, the induced map
\[d_\normal\varphi: \normal(U,V) \ra \normal(\mathbb{R}^n,\mathbb{R}^p) \cong \mathbb{R}^n\]
gives rise to the chart $(\normal(U,V),d_\normal\varphi)$ on $\normal(\base,\subbase)$ (it is even a vector bundle chart for $\normal(\base,\subbase) \ra \subbase$). More specifically, if we write $\varphi = (\varphi^1,\varphi^2) = (y^1,\dots,y^p,x^1,\dots,x^q)$, then for all $(y,\xi) \in \normal(U,V)$ (so $\xi \in \normal_y(\base,\subbase)$), we have
\[d_\normal\varphi(y)\xi = (\varphi^1(y),d\varphi^2(y)\xi) = (y_1(y),\dots,y_p(y),dy_1(y)\xi,\dots,dy_q(y)\xi).\]
To justify the notation $dy_i(y)\xi$, suppose that $\eta \in T_y\base$ with $\eta = \xi \mod T_y\subbase$, say, $\eta = \xi + v$ with $v \in T_y\subbase$. Then,
\[dy_i(y)\eta = dy_i(y)\xi + dy_i(y)v = dy_i(y)\xi \mod \normal_y(\base,\subbase),\] 
where we used that $y_i(V) = \{0\} \subset \mathbb{R}$ in the last equality, so that $dy_i(y)T_y\subbase = 0$. 
%We will also write $d_\normal\varphi = d\varphi^2|_X = (dy_1,\dots,dy_q)$.
\end{rema}
We will see that the normal bundle functor has many nice properties. In particular, applying the normal bundle functor to geometric structures, like principal bundles, Lie groupoids, vector bundles, and Lie algebroids, often yields another such structure. One can think of this procedure as ``linearising" the structure at a submanifold. As a first application, we will show here that the normal bundle functor applied to pairs of vector bundles yields another vector bundle, but we will need a result regarding maps of constant rank first. We will use the following terminology:
\begin{defn}\label{defn: constant rank pairs}
A smooth map of pairs $f: (\base,\subbase) \ra (\basetwo,\subbasetwo)$ is said to have \textit{constant rank $(k,\ell)$} if $f: \base \ra \basetwo$ has constant rank $k$, and $f|_{\subbase}: \subbase \ra \subbasetwo$ has constant rank $\ell$. Similarly, $f$ is said to be a submersion/immersion if $f: \base \ra \basetwo$ and $f|_\subbase: \subbase \ra \subbasetwo$ are both submersions/immersions.
\end{defn}
We will show that a submersion $f:(\base,\subbase) \ra (\basetwo,\subbasetwo)$ always induces the submersion $d_\normal f$, and an immersion $f:(\base,\subbase) \ra (\basetwo,\subbasetwo)$ always induces the immersion $d_\normal f$ under the extra assumption that $f^{-1}(\subbase)=\base$. To prove this, we will prove the following lemma about more general vector bundles first.
\begin{lemm}\label{lemm: vector bundle constant rank maps}
Let $\pi_E: E \ra M$ and $\pi_F: F \ra N$ be vector bundles, and let $\chi: E \ra F$ be a morphism of vector bundles with base map $f: M \ra N$. If $f$ has constant rank $k$, and $\chi$ has fiberwise constant rank $\ell$ (that is, for all $x \in M$, $\chi(x): E_x \ra F_{f(x)}$ has rank $\ell$), then $\chi$ has constant rank $k+\ell$.
\end{lemm}
\begin{proof}
Since $f: M \ra N$ has constant rank $k$, we can find a cover of $M$ by charts $(U,\varphi)$, and a cover of $N$ by charts $(W,\psi)$, such that $\psi \circ f \circ \varphi^{-1}$ takes the form $(x^1,\dots,x^n) \mapsto (x^1,\dots,x^k,0)$. Since $\im \chi(x)$ has constant rank $\ell$ for all $x \in U$, we can shrink $U$, if necessary, and assume that we have a local frame $\{s^i\}: U \ra E$ for which $\{s^i \mid \ell+1 \le i \le p\}: U \ra E|_U$ is a local frame for $\ker \chi$% (note: this shows that $\ker \chi \subset E$ is a vector subbundle over $M$)
. Then $\{\chi \circ s^i \mid 1 \le i \le \ell\}: U \ra F|_W$ is a local frame for $\im \chi$% (note: this shows that $\im \chi \subset F$ is a vector subbundle over $N$)
, and we can extend this collection of sections to a frame $\{t^j\}$ of $F|_W$. The vector bundle charts $\Phi: E|_U \ra \varphi(U) \times \mathbb{R}^p$ and $\Psi: F|_W \ra \psi(W) \times \mathbb{R}^q$ defined, via their inverses, by $\Phi^{-1}(\varphi(u),\xi) \coloneqq \textstyle\sum_i \xi_is^i(u)$ and $\Psi^{-1}(\psi(w),\xi) \coloneqq \textstyle\sum_j \xi_jt^j(w)$ define vector bundle charts for which 
\[\Psi \circ \chi \circ \Phi^{-1}(x^1,\dots,x^n,\xi_1,\dots,\xi_p) = (x^1,\dots,x^k,0,\xi_1,\dots,\xi_\ell,0).\]
This shows that $\chi$ has constant rank $k+\ell$.
\end{proof}
\begin{prop}\label{prop: normal bundle constant rank}
Let $f: (\base,\subbase) \ra (\basetwo,\subbasetwo)$ be a smooth map of pairs. If $f$ has constant rank $(k,\ell)$, and $d_\normal f$ has fiberwise constant rank $q$, then $d_\normal f$ has constant rank $\ell+q$. In particular, if $f$ is a submersion (as a map of pairs), then $d_\normal f$ is a submersion. If $f$ is an immersion (as a map of pairs), and $f^{-1}(\subbasetwo)=\subbase$, then $d_\normal f$ is an immersion.
\end{prop}
\begin{proof}
%Indeed, take charts $(U,\varphi)$ and $(W,\psi)$ around $x$ and $f(x)$, respectively, such that $f(U) \subset W$ and
%\[\psi \circ f \circ \varphi^{-1}: \varphi(U) \ra \psi(W) \textnormal{ is given by } (x^1,\dots,x^n) \mapsto (x^1,\dots,x^k,0).\]
%Then 
%\[d\psi \circ df \circ d\varphi^{-1} \textnormal{ is the map } (x^1,\dots,x^n,\xi^1,\dots,\xi^n) \mapsto (x^1,\dots,x^k,0,\xi^1,\dots,\xi^k,0).\]
%This shows that $df: T\base \ra T\subbase$ has constant rank $2k$, and a similar argument shows that $df|_{T\subbase}$ has constant rank $2\ell$. %Similarly, $d(df)$ has constant rank $(4k,4\ell)$. 
The first statement is an immediate consequence of Lemma \ref{lemm: vector bundle constant rank maps} once we recognise that $d_\normal f: \normal(\base,\subbase) \ra (\basetwo,\subbasetwo)$ is a vector bundle morphism over $f|_\subbase: \subbase \ra \subbasetwo$. For the second statement, it suffices to prove by the first statement that, for all $y \in \subbase$, $d_\normal f(y)$ has full rank. We distinguish the two cases \textbf{a)} $f$ is a submersion and \textbf{b)} $f$ is an immersion and $f^{-1}(\subbasetwo)=\subbase$. Say, the dimension of $(\base,\subbase)$ is $(r,s)$ and the dimension of $(\basetwo,\subbasetwo)$ is $(m,n)$.

\textbf{a)}
In this case, $f$ has constant rank $(m,n)$, so for all $x \in \base$, $df(x)$ has rank $(m,n)$. Pick a basis $a_1,\dots,a_{s-n}$ for $\ker df(y)|_{T_y\subbase}$ and extend it to a basis $a_1,\dots,a_{r-m}$ for $\ker df(y)$. Now pick elements $c_1 \coloneqq df(y)b_1,\dots,c_n \coloneqq df(y)b_n$ that form a basis for $\im df(y)|_{T_y\subbase}=T_{f(y)}\subbasetwo$. Then  $a_1,\dots,a_{s-n},b_1,\dots,b_n$ is a basis for $T_y\subbase$. In particular, if there exists a $s-n+1 \le j \le r-m$ such that
\[\lambda_j a_j = \sum_{i=1}^{s-n} \lambda_ia_i + \sum_{i'=1}^{n} \mu_{i'}b_{i'} \textnormal{ for some } \lambda_i,\mu_{i'} \in \mathbb{R},\]
then the $\mu_{i'}$ must be zero (apply this equation to $df(y)$), and so the $\lambda_i$ must be zero as well by construction of $a_j$. This shows that $a_1,\dots,a_{r-m},b_1,\dots,b_n$ is a linear independent collection, so, as before, we can extend it to a basis $a_1,\dots,a_{r-m},b_1,\dots,b_m$ of $T_y\base$ such that $c_1,\dots,c_n,c_{n+1} \coloneqq df(y)b_{n+1},\dots, df(y)b_m$ is a basis for $\im df(y)=T_y\base$. Now, $T_y\base/T_y\subbase$ has basis 
\[a_{s-n+1} \mod T_y\subbase,\dots,a_{r-m} \mod T_y\subbase,b_{n+1} \mod T_y\subbase,\dots,b_m \mod T_y\subbase,\] 
which maps under $d_\normal f(y)$ to 
\[c_{n+1} \mod T_{f(y)}\subbasetwo, \dots, c_m \mod T_{f(y)}\subbasetwo.\]
Since $c_1,\dots,c_n$ is a basis for $T_{f(y)}\subbasetwo$, we have that $c_{n+1},\dots,c_m \not\in T_{f(y)}\subbasetwo$. Therefore, the elements $c_{n+1} \mod T_{f(y)}\subbasetwo,\dots,c_m \mod T_{f(y)}\subbasetwo$ form a basis for $T_{f(y)}\basetwo/T_{f(y)}\subbasetwo$, so $d_\normal f(y)$ has rank $m-n$, as required.

\textbf{b)} Here, $f$ has constant rank $(r,s)$, so for all $x \in X$, $df(x)$ has rank $(r,s)$. Let $y \in Y$, and pick a basis $b_1,\dots,b_s$ for $T_y\subbase$. Then $c_1 \coloneqq df(y)b_1,\dots,c_s \coloneqq df(y)b_s$ is a linearly independent collection. We now extend $b_1,\dots,b_s$ to a basis for $T_y\base$, so then $c_1,\dots,c_s,c_{s+1}\coloneqq df(y)b_{s+1},\dots,c_r \coloneqq df(y)b_r$ forms a basis for $\im df(y)$. Since $c_1,\dots,c_s$ is a basis for $T_{f(y)}\subbasetwo$, we have $c_{s+1},\dots,c_r \not\in T_{f(y)}\subbasetwo$ by our assumption that $f^{-1}(\subbasetwo)=\subbase$ (otherwise, we can find a small transversal of $\subbase$ which maps into $\subbasetwo$). Therefore, the elements 
\[c_{s+1} \mod T_{f(y)}\subbasetwo,\dots, c_r \mod T_{f(y)}\subbasetwo\]
form a linear independent collection in $T_{f(y)}\basetwo/T_{f(y)}\subbasetwo$, so $d_\normal f$ has rank $r-s$, as required.
\end{proof}
\begin{rema}\label{rema: d_Nf immersion}
The extra condition $f^{-1}(\subbasetwo)=\subbase$ in the former statement (in the case that $f$ is an immersion), is necessary. A simple example to illustrate this is given by the map of pairs $\id_{\mathbb{R}}: (\mathbb{R},0) \ra (\mathbb{R},\mathbb{R})$. The normal derivative of this map is the map
\[T_0\mathbb{R} \cong \normal(\mathbb{R},0) \ra \normal(\mathbb{R},\mathbb{R}) \cong \mathbb{R} \textnormal{ given by } (0,\xi) \mapsto (0,\xi \mod T_0\mathbb{R}) = 0,\]
which, obviously, is not an immersion.
\end{rema}
Using the former statement, we know that, if $(E \xra{\pi} M, F \ra N)$ is a pair of vector bundles, then the map $d_\normal \pi: \normal(E,F) \ra \normal(M,N)$ is a submersion. We will now naturally equip $\normal(E,F) \ra \normal(M,N)$ with a vector bundle structure. One way to view this, is as being induced from the so-called \textit{secondary vector bundle structure} on $TE \ra TM$ (and $TF \ra TN$).
\begin{prop}\cite{KMS}\label{prop: normal bundle of vector bundle}
Let $(E \xra{\pi} M, F \ra N)$ be a pair of vector bundles of rank $(k,\ell)$. Then $TE \ra TM$ (and $TF \ra TN$) is naturally a vector bundle of rank $2k$ (and $2\ell$, respectively). Similarly, 
\[\normal(E,F) \xra{d_\normal\pi} \normal(M,N)\]
is a vector bundle of rank $k$, and we have a vector bundle morphism $TE|_F \ra \normal(E,F)$ over $TM|_N \ra \normal(M,N)$. 
\end{prop}
\begin{proof}
To define the vector bundle structure on $TE \ra TM$, denote by $+_E: E \tensor[_{\pi}]{\times}{_{\pi}} E \ra E$ and $\cdot_E: \mathbb{R} \times E \ra E$ the addition and scalar multiplication of $E$, respectively. Then we define the vector bundle structure on $TE \ra TM$ by setting 
\begin{align*}
    +_{TE} &\coloneqq d(+_E): T(E \tensor[_{\pi}]{\times}{_{\pi}} E) \cong TE \tensor[_{d\pi}]{\times}{_{d\pi}} TE \ra TE \textnormal{, and } \\
    \cdot_{TE} &\coloneqq d(\cdot_E)(\bullet,\star)(1,\diamond): \mathbb{R} \times TE \ra TE
\end{align*}
(it is readily verified that we have a canonical isomorphism $T(E \tensor[_{\pi}]{\times}{_{\pi}} E) \cong TE \tensor[_{d\pi}]{\times}{_{d\pi}} TE$, but see Lemma \ref{lemm: if clean intersection, then T respects fiber products} for a more general statement). The zero section of $TE \ra TM$ is given by $d0_M: TM \ra TE$. For $(x,\xi) \in TM$, we write $(TE)_{(x,\xi)}$ for the fiber of $TE$ at $(x,\xi)$. Observe that $TE \ra TM$ inherits local frames from $E$ as follows: let $\{s^j: U \ra E\}$ be a local frame. Then $ds^j: TU \ra TE$ maps into $d\pi^{-1}(TU)$ by the chain rule, and, obviously,
\[s^j: TU \xra{\pi} U \xra{s^j} E \xra{0_E} TE\]
also maps into $d\pi^{-1}(TU)$. Therefore, we only have to check, for all $(x,\xi) \in TU$, that the collection $\{(s^j(x),0),(s^j(x),ds^j(x)\xi)\}$ forms a basis for $(TE)_{(x,\xi)}$. To see this, suppose
\[\sum_j \mu_j \cdot_{TE} (s^j(x),0) +_{TE} \lambda_j \cdot_{TE} (s^j(x),ds^j(x)\xi) = 0 = (0_M(x),d0_M(x)\xi)\]
(note: $\textstyle\sum$ is with respect to $+_{TE}$), then, we see that $\textstyle\sum_j\mu_js^j + \lambda_js^j=0$, and by repeatedly using the chain rule, we see that $\lambda_1 s^1(x) + \cdots + \lambda_ks^k(x)=0$ also holds, so $\lambda_1=\cdots=\lambda_k=0$, and then also $\mu_1=\cdots=\mu_k=0$. This shows that $TE \ra TM$ is a vector bundle with $+_{TE}$ and $\cdot_{TE}$.

Now, to see that $\normal(E,F) \ra \normal(M,N)$ is a vector bundle, replace $d$, in the above discussion, with $d_\normal$ everywhere (again, it is readily verified that we have a canonical isomorphism $\normal(E \tensor[_{\pi}]{\times}{_{\pi}} E, F \tensor[_{\pi|_F}]{\times}{_{\pi|_F}} F) \cong \normal(E,F) \tensor[_{d_\normal \pi}]{\times}{_{d_\normal \pi}} \normal(E,F)$, but see Lemma \ref{prop: if clean intersection, then N respects fiber products} for a more general statement). For sections $s: U \ra E$, to make sense of $d_\normal s$ as a map $\normal(U,V) \ra \normal(E,F)$, where $V \coloneqq U \cap F$, $s$ has to restrict to a map $V \ra F$. So, to see that $\normal(E,F)$ inherits a local frame from $E$ (and $F$), pick a local frame $\{s^j: U \ra E\}$ for which $\{s^j|_V: V \ra F \mid 1 \le j \le \ell\}$ is a local frame for $F$. Then, $\{s^i, d_\normal s^j \mid \ell+1 \le i \le k, 1 \le j \le \ell\}$ is a local frame for $\normal(E,F) \ra \normal(M,N)$ using the same arguments as before. 
%let $(E|_U,\varphi)$ be a local trivialisation of $E$; that is, the map
%\[\varphi: \pi^{-1}(U) \eqqcolon E|_U \xra{\sim} U \times \mathbb{R}^k\]
%is smooth, satisfies $\pi|_{E|_U} = \pr_1 \circ \varphi$ and the maps $\varphi_x \coloneqq \pi^{-1}(x) \eqqcolon E_x \ra \mathbb{R}^k$ are linear isomorphisms. Then, also $d\pi|_{TE|_{TU}} = \pr_{TU} \circ d\varphi$ (where $TE|_{TU} \coloneqq d\pi^{-1}(TU)$), so $(TU,d\varphi)$ is a local trivialisation of $TE \ra TY$ (note that this means that for $(x,\xi) \in TU$ we equip $TE_{(x,\xi)} \coloneqq d\pi^{-1}(x,\xi)$ with the vector space structure induced by the bijection $d\varphi^{-1}(x,\xi,\cdot)|_{\mathbb{R}^{2k}}$). This proves that $TE \xra{d\pi} TY$ is a vector bundle; in particular, $TE|_F \ra TY|_X$ is a vector bundle (of rank $k + q$% since $TE|_F \ra F$ has rank $n+k$ and $TY|_X \ra X$ has rank $n$) which maps $TF$ into $TX$. The induced map 
%\[TE|_F/TF = \normal(E,F) \xra{d_\normal\pi} \normal(Y,X) = TY|_X/TX\]
%becomes a vector bundle of rank $k$, because a local trivialisation 
%\[(TE|_F)|_{TU|_{U \cap X}} \xra{\sim} TU|_{U \cap X} \times \mathbb{R}^k \times \mathbb{R}^q\]
%that restricts to a local trivialisation $TF|_{T(U \cap X)} \xra{\sim} T(U \cap X) \times \mathbb{R}^q \times \mathbb{R}^q$ induces a local trivialisation 
%\[\normal(E,F)|_{TU/T(U \cap X)} \xra{\sim} TU|_{U \cap X}/T(U \cap X) \times (\mathbb{R}^k/\mathbb{R}^q) \times \mathbb{R}^q \cong TU|_{U \cap X}/T(U \cap X) \times \mathbb{R}^k.\]

For the last statement, it suffices to prove that the maps $(TE|_F)_{(y,\xi)} \ra \normal(E,F)_{(y,\xi \mod T_yN)}$ are linear, but this is obvious from the above vector bundle constructions. This concludes the proof.
\end{proof}
\begin{rema}\cite{Gracia_Saz_2010}\label{rema: tangent and normal double vector bundle}
Notice that (using notation from above) we have commutative diagrams
\begin{center}
\begin{tikzcd}
TE \ar[r,"d\pi"] \ar[d] & TY \ar[d] \\
E \ar[r,"\pi"] & Y,
\end{tikzcd}\hspace{1cm}%
\begin{tikzcd}
\normal(E,F) \ar[r,"d_\normal\pi"] \ar[d] & \normal(Y,X) \ar[d]
\\
F \ar[r,"\pi|_F"] & X
\end{tikzcd}
\end{center}
such that every arrow involved is a vector bundle. Moreover, they satisfy certain compatibility conditions; for example, $d\pi$ and $d_\normal \pi$ are vector bundle morphisms over $\pi$ and $\pi|_F$, respectively. The ``vector bundle diagrams" above are so-called \textit{double vector bundles}. We will explain more about this notion later.
\end{rema}
As a special case, one can always construct, out of a vector bundle $E \ra M$, the vector bundle $\normal(E,0_M) \ra \normal(M,M) \cong M$. However, these two vector bundles are canonically isomorphic. It follows from the following simple (but useful) lemma.
\begin{lemm}\cite{meinrenken}\label{lemm: canonical iso TE|_M}
Let $\pi: E \ra M$ be a vector bundle. Then, there is a canonical isomorphism of vector bundles $E \oplus TM \xra{\sim} TE|_M$.
\end{lemm}
\begin{proof}
Observe that we have a short exact sequence of vector bundles (over the same base)
\begin{center}
    \begin{tikzcd}[row sep=small, column sep=small]
    0 \ar[r] & E \ar[r,hook,"0_E"] & TE|_M \ar[r,two heads,"d\pi|_M"] & TM \ar[r] & 0
    \end{tikzcd}
\end{center}
that has a canonical splitting $TM \ra TE|_M$ given by $d0_M$. The result follows, and the isomorphism $E \oplus TM \xra{\sim} TE|_M$ is given by $((x,\zeta),(x,\xi)) \mapsto 0_E(x,\zeta)+(x,d0_M(x)\xi)$.
\end{proof}
\begin{coro}\cite{meinrenken}\label{coro: canonical iso N(E,0_M)}
Let $E \ra M$ be a vector bundle. Then there is a canonical isomorphism of vector bundles $E \xra{\sim} \normal(E,0_M)$.
\end{coro}
\begin{proof}
First of all, by Lemma \ref{lemm: canonical iso TE|_M}, we have a canonical isomorphism $E \oplus TM \xra{\sim} TE|_M$. This isomorphism sends $(x,\xi) \in TM$ to $(x,d0_M(x)\xi) \in T0_M$ (note: we also wrote $0_M \subset E$). In particular, we see that the map 
\[E \oplus TM \xra{\sim} TE|_M \rightarrowdbl \normal(E,0_M)\]
has kernel $TM$, so the induced map $E \ra \normal(E,0_M)$ is an isomorphism of vector bundles, and it is given by $(x,\zeta) \mapsto 0_E(x,\zeta) \mod T_x0_M$.
\end{proof}
Lastly, with the normal bundle functor at hand, the notion of a tubular neighbourhood embedding can be described cleanly as follows:
\begin{defn}\cite{meinrenken}\label{defn: tubular neighbourhood}
Let $(\base,\subbase)$ be a pair of manifolds. A \textit{tubular neighbourhood} of $\subbase$ in $\base$ is a smooth map of pairs $\chi: (\normal(\base,\subbase),0_\subbase) \ra (\base,\subbase)$ such that $\chi: \normal(\base,\subbase) \ra \base$ is an embedding, and $d_\normal\chi = \id_{\normal(\base,\subbase)}$.
\end{defn}
\begin{rema}\label{rema: tubular neighbourhood}
Notice that if $\chi: (\normal(\base,\subbase),0_\subbase) \ra (\base,\subbase)$ is a tubular neighbourhood embedding, then, necessarily, $\chi|_{\subbase}=\id_\subbase$.
\end{rema}

\subsection{Deformation to the normal cone functor: objects}\label{sec: Deformation to the normal cone functor: objects}
Let $(\base,\subbase)$ be a pair of (possibly non-Hausdorff) smooth manifolds of dimension $(n,p)$. As a set, we define the deformation to the normal cone of $\subbase$ in $\base$ as
\[D^\base_\subbase \coloneqq (\normal(\base,\subbase) \times \{0\}) \sqcup (\base \times \mathbb{R}^\times).\]
The purpose of this section is to define a smooth structure on this set. Later, when this set is equipped with a smooth structure, we will denote the resulting manifold by $\DNC(\base,\subbase)$. To define the smooth structure, we will explicitly construct an atlas. The first step is to describe the smooth structure for 
\[\subbase = \mathbb{R}^p \times \{0\} = \mathbb{R}^p \subset \mathbb{R}^n = \mathbb{R}^{p+q} = \base,\]
in which case we will (temporarily) also write $D^n_p$ for $D^\base_\subbase$. Note that, as a set, we can canonically identify $\normal(\base,\subbase)$ with $\mathbb{R}^n$. So, there is a canonical bijection $D^n_p \ra \mathbb{R}^n \times \mathbb{R}$. We will transfer the (standard) smooth structure of $\mathbb{R}^n \times \mathbb{R}$ to $D^n_p$ using a different bijection, namely, via the map
\begin{equation*}
\Psi: \mathbb{R}^{p+q} \times \mathbb{R} \rightarrow D^n_p \textnormal{ given by } \left(y,\xi,t\right) \mapsto
    \begin{cases}
    \left(y,\xi,0\right) & \text{if}\ t=0 \\
    \left(y,t\xi,t\right) & \text{if}\ t\neq0
    \end{cases}
\end{equation*}
which has inverse
\begin{equation*}
\Psi^{-1}: D^n_p \rightarrow \mathbb{R}^{p+q} \times \mathbb{R} \textnormal{ given by } (y,\xi,t) \mapsto    
\begin{cases}
    (y,\xi,0) & \text{if}\ t=0 \\
    (y,t^{-1}\xi,t) & \text{if}\ t\neq0.
    \end{cases}
\end{equation*}
We denote $\DNC(\mathbb{R}^n,\mathbb{R}^p)$ for $D^n_p$ equipped with this smooth structure. Observe that if $U \subset \mathbb{R}^n$ is an open subset and $V \coloneqq U \cap \mathbb{R}^p$, then $\normal(U,V) \cong V \times (\mathbb{R}^n/\mathbb{R}^p) \cong V \times \mathbb{R}^q$ (again, canonically), so
\[D^U_V = (V \times \mathbb{R}^q \times \{0\}) \sqcup (U \times \mathbb{R}^\times) \subset \DNC(\mathbb{R}^n,\mathbb{R}^p).\]
Observe that $D^U_V$ is the open subset $\Psi(\Omega^U_V)$ of $\DNC(\mathbb{R}^n,\mathbb{R}^p)$, where
\[\Omega^U_V = \{\left(y,\xi,t\right) \in \mathbb{R}^{p+q} \times \mathbb{R} \mid (y,t\xi) \in U\},\]
which, indeed, is an open subset of $\mathbb{R}^n \times \mathbb{R}$ as it is the preimage of $U$ under the continuous map 
\[\mathbb{R}^{p+q} \times \mathbb{R} \ra \mathbb{R}^{p+q} \textnormal{ given by } \left(y,\xi,t\right) \mapsto (y,t\xi).\]
We write $\DNC^U_V$ for $D^U_V$ equipped with the smooth structure it inherits from $\DNC(\mathbb{R}^n,\mathbb{R}^p)$ as an open subset.

Now, in the general case where $(\base,\subbase)$ is a pair of smooth manifolds, we can cover $\base$ with adapted charts (that is, adapted whenever the domain of the chart has non-empty intersection with $\subbase$).
%charts $(U_\alpha,\varphi_\alpha)_{\alpha \in I}$, such that whenever $X \cap U_\alpha \neq \emptyset$, then $(U,\varphi)$ is an adapted chart. 
If $(U,\varphi)$ is an adapted chart, write $V \coloneqq U \cap \subbase$ and $\varphi = (\varphi^1,\varphi^2)$. 
By Remark \ref{rema: local coordinates normal bundle}, we see that the map
\begin{equation*}
\DNC(\varphi): D^U_V \ra \DNC^{\varphi(U)}_{\varphi(V)} \textnormal{ given by } z \mapsto
    \begin{cases}
    \left(\varphi^1(y),d\varphi^2(y)\xi,0\right) & \text{if}\ z = (y,\xi,0) \\
    \big(\varphi(x),t\big) & \text{if}\ z = (x,t),
    \end{cases}
\end{equation*}
is a bijection, so we can transfer the smooth structure constructed before for $\DNC^{\varphi(U)}_{\varphi(V)}$ to $D^{U}_{V}$. The latter set equipped with this smooth structure is denoted by $\DNC^{U}_{V}$. We claim that, if $(U_i,\varphi_i)_{i \in I}$ is a cover by compatible adapted charts (adapted to $\subbase$) of $\base$, then
\[\left(\DNC^{U_i}_{V_i}, \Phi_i: \DNC^{U_i}_{V_i} \ra \Omega^{U_i}_{V_{i}}\right), \textnormal{ where } \Phi_i = \Psi^{-1} \circ \DNC(\varphi_i); \quad \Omega^{U_i}_{V_{i}} \coloneqq \Omega^{\varphi_i(U_i)}_{\varphi_i(V_{i})},\]
constitute a cover by charts of $D^\base_\subbase$ (note: the chart $(\DNC^{U_i}_{V_i},\Phi_i)$ is defined even if $(U_i,\varphi_i)$ is a chart away from $\subbase$; $\Omega^U_V$, where $V \coloneqq U \cap \subbase = \emptyset$, is given as above, and $\Phi_i \coloneqq \Psi^{-1} \circ \varphi_i$). To see that the transition maps are smooth, observe first that if $i,j \in I$, then $\DNC^{U_i}_{V_i} \cap \DNC^{U_j}_{V_j}$ can be identified with $\DNC^{U_i \cap U_j}_{V_i \cap V_j}$, because for open subset $U,U' \subset \mathbb{R}^n$ we have
\[\Omega^U_V \cap \Omega^{U'}_{V'} = \Omega^{U \cap U'}_{V \cap V'} \textnormal{ } (\textnormal{where } V \coloneqq U \cap \mathbb{R}^p \textnormal{ and } V' \coloneqq U' \cap \mathbb{R}^p)\] 
from the definitions. %It follows that, if $k \in I$, then
%\begin{align*}
%    \DNC(\varphi_k)\left(\DNC^{U_k \cap U_i}_{V_k \cap V_i} \cap \DNC^{U_k \cap U_j}_{V_k \cap V_j}\right) &= \Psi\left(\Omega^{\varphi_k(U_k \cap U_i)}_{\varphi_k(V_k \cap V_i)}\right) \cap \Psi\left(\Omega^{\varphi_k(U_k \cap U_j)}_{\varphi_k(V_k \cap V_j)}\right) \\
%    &= \Psi\left(\Omega^{\varphi_k(U_k \cap U_i \cap U_j)}_{\varphi_k(V_k \cap V_i \cap V_j)}\right) = \DNC(\varphi_k)\left(\DNC^{U_k \cap U_i \cap U_j}_{U_k \cap U_i\cap U_j}\right)
%\end{align*} 
%which shows that indeed $\DNC^{U_i}_{V_i} \cap \DNC^{U_j}_{V_j}$ can be identified with $\DNC^{U_i \cap U_j}_{V_i \cap V_j}$. 
That the transition maps
%\[\Phi_\beta \circ \Phi_\alpha^{-1} = \Psi^{-1} \circ \widetilde{\varphi_\beta} \circ \widetilde{\varphi_\alpha}^{-1} \circ \Psi: \Phi_\alpha\left(\DNC^{U_\alpha}_{V_\alpha} \cap \DNC^{U_\beta}_{V_\beta}\right) \ra \Phi_\beta\left(\DNC^{U_\alpha}_{V_\alpha} \cap \DNC^{U_\beta}_{V_\beta}\right),\]
\[\Phi_j \circ \Phi_i^{-1} = \Psi^{-1} \circ \DNC(\varphi_j) \circ \DNC(\varphi_i)^{-1} \circ \Psi: \Phi_i\left(\DNC^{U_i \cap U_{j}}_{V_{i} \cap V_j}\right) \ra \Phi_j\left(\DNC^{U_i \cap U_j}_{V_i \cap V_j}\right)\]
%To see this, note that we have a diffeomorphism
%\begin{align*}
%    \chi: \normal(U \cap U',V \cap V') &\ra \normal(U,V) \cap \normal(U',V') \textnormal{ given by } (x,\xi + T_x(V \cap V')) \mapsto (x,\xi + T_xV \cap T_xV'),
%\end{align*}
%so the map 
%\begin{equation*}
%    \DNC^U_V \cap \DNC^{U'}_{V'} \ra \DNC^{U \cap U'}_{V \cap V'} \textnormal{ given by } 
%\begin{cases}
%    (x,\xi,0) \mapsto (\chi(x,\xi),0) \\
%    (y,t) \mapsto (y,t)
%\end{cases}
%\end{equation*}
%is a diffeomorphism.
are smooth can be deduced from the following lemma with (see Corollary \ref{coro: transition maps of DNC are smooth}):
\[\varphi_j \circ \varphi_i^{-1} = (\varphi_j^1\circ \varphi_i^{-1},\varphi_j^2\circ \varphi_i^{-1}) \eqqcolon (h_1,h_2) = h.\]

\begin{lemm}\cite{rouse2008schwartz}\label{lemm: transition maps of DNC are smooth}
Let $h = (h_1,h_2): U \ra U'$ be a diffeomorphism between two open subsets $U$ and $U'$ of $\mathbb{R}^n=\mathbb{R}^{p+q}$. We set
\[V \coloneqq U \cap \mathbb{R}^p \textnormal{ and } \Omega^U_V \coloneqq \{\left(y,\xi,t\right) \in \mathbb{R}^{p+q} \times \mathbb{R} \mid (y,t\xi) \in U\},\] 
and we define $V'$ and $\Omega^{U'}_{V'}$ similarly. If $h_2|_V = 0$, then the map 
\begin{equation*}
\widetilde{h}: \Omega^U_V \ra \Omega^{U'}_{V'} \textnormal{ given by } \left(y,\xi,t\right) \mapsto
    \begin{cases}
    \big(h_1(y,0),\tfrac{\partial h_2}{\partial \xi}(y,0) \xi, 0\big) & \text{if}\ t=0 \\
    \left(h_1(y,t\xi),t^{-1}h_2(y,t\xi),t\right) & \text{if}\ t\neq0,
    \end{cases}
\end{equation*}
is smooth.
\end{lemm}
\begin{proof}
Write $\widetilde{h} = (\widetilde{h}_1,\widetilde{h}_2,\widetilde{h}_3)$ (e.g. $\widetilde{h}_3\left(y,\xi,t\right) = t$). Then it is clear that $\widetilde{h}_1$ and $\widetilde{h}_3$ are smooth, so we only need to check that $\widetilde{h}_2$ is smooth. Doing this coordinate-wise, we can reduce to the case that $h$ is a map $U \ra \mathbb{R}$ (with $h(y,0) = 0$) and 
\begin{equation*}
\widetilde{h}: \Omega^U_V \ra \mathbb{R} \textnormal{ is given by } \left(y,\xi,t\right) \mapsto 
    \begin{cases}
    \tfrac{\partial h}{\partial \xi}(y,0) \xi & \text{if}\ t=0 \\
    t^{-1}h(y,t\xi) & \text{if}\ t\neq0.
    \end{cases}
\end{equation*}  
Using that $h|_V = 0$, we have that
\[h(y,\xi) = \int_0^1 \frac{d}{dt}h(y,t\xi) dt = \int_0^1 \frac{\partial h}{\partial \xi}(y,t\xi)\xi dt = \frac{\partial h}{\partial \xi}(y,0)\xi + r(y,\xi)\xi,\]
where $r(y,\xi) \coloneqq -\textstyle\frac{\partial h}{\partial \xi}(y,0) + \textstyle\int_0^1\frac{\partial h}{\partial \xi}(y,t\xi)dt$ is a smooth map with the property that $r|_V=0$.
Hence,
\[t^{-1}h(y,t\xi) = \frac{\partial h}{\partial \xi}(y,0)\xi + r(y,t\xi)\xi,\]
so using that $r|_V=0$, the result follows.
\end{proof}
\begin{coro}\cite{rouse2008schwartz}\label{coro: transition maps of DNC are smooth}
If $(U_i,\varphi_i)_{i \in I}$ is a collection of compatible adapted charts (adapted to $\subbase$), then the collection $\left(\DNC^{U_i}_{V_i}, \Phi_i: \DNC^{U_i}_{V_i} \ra \Omega^{U_i}_{V_{i}}\right)$ is a collection of compatible charts for $D^\base_\subbase$. In particular, this way we obtain a smooth structure on $D^\base_\subbase (=\DNC(\base,\subbase))$.
\end{coro}
\begin{proof}
As mentioned before, this follows from Lemma \ref{lemm: transition maps of DNC are smooth} applied to \[\varphi_j \circ \varphi_i^{-1} = (\varphi_j^1\circ \varphi_i^{-1},\varphi_j^2\circ \varphi_i^{-1}) \eqqcolon (h_1,h_2) = h.\]
Indeed, 
%recall that the inverse of $\Psi$ is given by 
%\begin{equation*}
%\Psi^{-1}: \DNC(\mathbb{R}^n,\mathbb{R}^p) \rightarrow \mathbb{R}^{p+q} \times \mathbb{R} \textnormal{ given by } \left(x,\xi,t\right) \mapsto
%    \begin{cases}
%    \left(x,\xi,0\right) & \text{if}\ t=0 \\
%    \left(x,t^{-1}\xi,t\right) & \text{if}\ t\neq0.
%    \end{cases}
%\end{equation*}
%Hence, 
the transition maps are given by 
\begin{align*}
    \Phi_j \circ \Phi_i^{-1}(y,\xi,t) &= \Psi^{-1} \circ \DNC(\varphi_j) \circ \DNC(\varphi_i)^{-1} \circ \Psi(y,\xi,t) \\
    &= 
    \begin{cases}
        \Psi^{-1}\left(\varphi^1_j \circ \varphi_i^{-1}(y),d_\normal\varphi_j\left(\varphi_i^{-1}(y)\right)d_\normal(\varphi_i^{-1})(y)\xi,0\right) \\
        \Psi^{-1}\left(\varphi_j \circ \varphi_i^{-1}(y,t\xi),t\right)
    \end{cases} \\
    &=
    \begin{cases}
        \Psi^{-1}\left(h_1(y),d_\normal(\varphi_j \circ \varphi_i^{-1})(y)\xi,0\right) \\
        \Psi^{-1}\left(h(y,t\xi),t\right)
    \end{cases} \\
    %&=
    %\begin{cases}
    %    \Psi^{-1}\big(h_1(x),\tfrac{\partial h_2}{\partial \xi}(x,0) \xi,0\big) \\
    %    \Psi^{-1}\big(h_1(x,t\xi),h_2(x,t\xi),t\big)
    %\end{cases} \\
    &=
    \begin{cases}
        \big(h_1(y),\tfrac{\partial h_2}{\partial \xi}(y,0) \xi,0\big) & \text{if}\  t=0 \\
        \big(h_1(y,t\xi),t^{-1}h_2(y,t\xi),t\big) & \text{if}\ t\neq0,
    \end{cases}
\end{align*}
which proves the statement by an application of Lemma \ref{lemm: transition maps of DNC are smooth}.
\end{proof}
We have now shown that to each pair of smooth manifolds $(\base,\subbase)$ we can associate the smooth manifold $\DNC(\base,\subbase)$. As already mentioned, we will refer to it as the deformation to the normal cone of $\subbase$ in $\base$.

\begin{term}\label{term: DNC_t}
We will write $\DNC(\base,\subbase)_0 = \normal(\base,\subbase) \times \{0\} \subset \DNC(\base,\subbase)$ and $\DNC(\base,\subbase)_t = \base \times \{t\} \subset \DNC(\base,\subbase)$ for $t\neq0$. We will often write $(y,\xi) \in \DNC(\base,\subbase)_0 \subset \DNC(\base,\subbase)$ instead of $(y,\xi,0)$. Moreover, for all $t \in \mathbb{R}$, we will call $\DNC(\base,\subbase)_t$ \textit{the slice at $t$} of $\DNC(\base,\subbase)$.
\end{term}

\begin{rema}\cite{2017arXiv170509588D}\label{rema: gauge action on DNC}
Observe that for a pair of smooth manifolds $(\base,\subbase)$, $\mathbb{R}^\times$ acts canonically on $\DNC(\base,\subbase)$ via the group action
\begin{equation*}
\lambda \cdot z =
\begin{cases}
    \left(x,\lambda^{-1}\xi,0\right) & \text{if}\ z=(y,\xi)  \\
    (x,\lambda t) & \text{if}\ z=(x,t).
\end{cases}
\end{equation*}
We claim that the action is smooth. To see this, let $(U,\varphi)$ be an adapted chart and let $(\DNC^U_V,\Phi)$ be the induced chart of $\DNC(\base,\subbase)$ (see the discussion before Lemma \ref{lemm: transition maps of DNC are smooth}). Then 
\begin{align*}
\Phi\left(\lambda \cdot \Phi^{-1}(y,\xi,t)\right) &= 
\begin{cases}
    \Phi((\varphi^1)^{-1}(y),\lambda^{-1}d_\normal\varphi^{-1}(y)\xi,0) \\
    \Phi(\varphi^{-1}(y,t\xi),\lambda t)
\end{cases} \\
&=
\begin{cases}
    (y,\lambda^{-1}\xi,0) & \text{if}\ t=0 \\
    (y,\lambda^{-1}\xi,\lambda t) & \text{if}\ t \neq 0,
\end{cases} 
\end{align*}
which clearly is a smooth map. This action will play a key role in the construction of the blow-up; we will refer to it simply as the $\mathbb{R}^\times$-action on $\DNC(\base,\subbase)$. 
\end{rema}
Lastly, we want to mention that often the deformation to the normal cone is constructed by ``gluing'' along a tubular neighbourhood. To make the word ``gluing'' more precise, recall that if 
\[\varphi: U_1 \xra{\sim} U_2\]
is a diffeomorphism between open subsets $U_1 \subset \basetwo_1$ and $U_2 \subset \basetwo_2$ of (possibly non-Hausdorff) manifolds $\basetwo_1$ and $\basetwo_2$, then we can form a new smooth (possibly non-Huasdorff) manifold
\[\basetwo_1 \cup_{\varphi} \basetwo_2 \coloneqq (\basetwo_1 \sqcup \basetwo_2)/\sim,\]
where $\sim$ is the equivalence relation which identifies elements related by $\varphi$; that is, $x \sim y$ if and only if $y=x$, $y=\varphi(x)$, or $\varphi(y)=x$. This equivalence relation is open, and using charts of $\basetwo_1,\basetwo_2$, we can describe compatible charts of $\basetwo_1 \cup_{\varphi} \basetwo_2$. Indeed, if $(W_1,\psi_1)$ and $(W_2,\psi_2)$ are charts of $\basetwo_1$ and $\basetwo_2$, respectively, such that $\psi_2 \circ \varphi = \psi_1$, %$\psi_2 \circ \varphi \circ \psi_1^{-1}: \psi_1(W_1 \cap \varphi^{-1}(W_2)) \ra \psi_2(W_2)$ is smooth, and therefore 
then one can show that the map
\[\psi_1 \cup_{\varphi} \psi_2: W_1 \cup_{\varphi} W_2 \ra \mathbb{R}^n \textnormal{ given by } x_i \mapsto \psi_i(x_i), \textnormal{ where } x_i \in W_i,\]
is well-defined, and that such maps define charts of $\basetwo_1 \cup_{\varphi} \basetwo_2$. Such manifolds have the property that, if $f_i: \basetwo_i \ra P$ ($i=1,2$) are smooth maps such that $f_2 \circ \varphi = f_1$, then the canonical map
\[f_1 \cup_{\varphi} f_2: \basetwo_1 \cup_{\varphi} \basetwo_2 \ra P \textnormal{ given by } x_i \mapsto f_i(x_i), \textnormal{ where } x_i \in \basetwo_i\]
is smooth. We will say that $\basetwo_1$ and $\basetwo_2$ are \textit{glued via $\varphi$}. 
\begin{rema}\cite{proper}\label{rema: tubular neighbourhood DNC}
Let $\chi: (\normal(\base,\subbase),0_\subbase) \ra (\base,\subbase)$ be a tubular neighbourhood of $\subbase$ in $\base$ onto $U \subset \base$. We will show here that $\DNC(\base,\subbase)$ can be obtained by gluing $\normal(\base,\subbase) \times \{0\}$ and $\base \times \mathbb{R}^\times$ using the tubular neighbourhood $\chi$. %Slightly more generally than in Definition \ref{defn: tubular neighbourhood}, we allow $\chi$ to be a diffeomorphism $W \ra U$, where $0_\subbase \subbase W \subset \normal(\base,\subbase)$ is open, and $\subbase \subset U \subset \base$ is open. 
To see how this works, %let
%\[\Omega \coloneqq \normal(\base,\subbase) \times \mathbb{R}\]
%(note: it is an open subset of $\normal(\base,\subbase) \times \mathbb{R}$), and 
define the following diffeomorphism:
\[\hat\chi: \normal(\base,\subbase) \times \mathbb{R}^\times \ra U \times \mathbb{R}^\times \textnormal{ given by } (y,\xi,t) \mapsto (\chi(y,t\xi),t).\]
Then we define
\[\DNC_{\textnormal{tub}}(\base,\subbase) \coloneqq (\normal(\base,\subbase) \times \mathbb{R}) \cup_{\hat\chi} (\base \times \mathbb{R}^\times)\]
(the subscript ``tub'' stands for tubular neighbourhood). To see that this manifold is diffeomorphic to $\DNC(\base,\subbase)$, notice that the maps
\begin{align*}
    f: \normal(\base,\subbase) \times \mathbb{R} &\ra \DNC(\base,\subbase) \textnormal{ given by } (y,\xi,t) \mapsto
\begin{cases}
    (y,\xi,0) & \text{if}\ t=0 \\
    (\chi(y,t\xi),t) & \text{if}\ t\neq0
\end{cases} \textnormal{ and } \\
    g: \base \times \mathbb{R}^\times &\hookrightarrow \DNC(\base,\subbase) \textnormal{ given by } (x,t) \mapsto (x,t)
\end{align*}
are smooth, and $g \circ \hat{\chi} = f$, so that we obtain a smooth map
\[\DNC_{\textnormal{tub}}(\base,\subbase) \ra \DNC(\base,\subbase)\]
which is readily verified to have a smooth inverse, namely the map induced by the canonical map
\[(\normal(\base,\subbase) \times \{0\}) \sqcup (\base \times \mathbb{R}^\times) \hookrightarrow (\normal(\base,\subbase) \times \mathbb{R}) \sqcup (\base \times \mathbb{R}^\times).\]
This shows the equivalence of the constructions $\DNC(\base,\subbase)$ and $\DNC_{\textnormal{tub}}(\base,\subbase)$ (and it also shows that $\DNC_{\textnormal{tub}}(\base,\subbase)$ is independent of the chosen tubular neighbourhood $\chi$).
\end{rema}

\subsection{Deformation to the normal cone functor: morphisms}\label{sec: Deformation to the normal cone functor: morphisms}
Let $(\base,\subbase)$ and $(\basetwo,\subbasetwo)$ be pairs of smooth manifolds and let $f: (\base,\subbase) \ra (\basetwo,\subbasetwo)$ be a morphism of pairs. We claim that the induced map
\begin{equation*}
\DNC(f): \DNC(\base,\subbase) \ra \DNC(\basetwo,\subbasetwo) \textnormal{ given by } z \mapsto
\begin{cases}
\left(f(y),d_\normal f(y)\xi,0\right) & \text{if}\ z=(y,\xi,0) \\ 
(f(x),t) & \text{if}\ z=(x,t)
\end{cases}
\end{equation*}
is smooth. Indeed, let $(U,\varphi)$ be an adapted chart of $\base$ and let $(W,\psi)$ be an adapted chart of $\basetwo$. Then the result follows by a straightforward application of Lemma \ref{lemm: transition maps of DNC are smooth}, now with $h= \psi \circ f \circ \varphi^{-1}$. From the definition, it is obvious that $\DNC(g \circ f) = \DNC(g) \circ \DNC(f)$, and $\DNC(\id_\base) = \id_{\DNC(\base,\subbase)}$ (where we view $\id_\base$ as a map $(\base,\base) \ra (\base,\base)$). Consequently, $\DNC: \mathcal{C}^\infty_2 \ra \mathcal{C}^\infty$ is a (covariant) functor. Moreover, notice that, if $U \subset \base$ is an open subset (with $V \coloneqq U \cap \subbase$) which maps via $f$ into $W \subset \basetwo$ (with $Z \coloneqq W \cap \subbasetwo$), the map $\DNC(f)$ restricts to the map 
\[\DNC(f|_U): \DNC(U,V) \ra \DNC(W,Z).\]
In particular, if $W$ is Hausdorff, $\DNC(f|_U)$ is the unique extension of the map $f|_U \times \id_{\mathbb{R}^\times}$. In turn, this shows that $\DNC(f)$ is the unique extension of the map $f \times \id_{\mathbb{R}^\times}$.
\begin{rema}\cite{2017arXiv170509588D}\label{rema: gauge action on DNC maps}
Clearly, in the above situation, the smooth map $\DNC(f)$ is $\mathbb{R}^\times$-equivariant with respect to the $\mathbb{R}^\times$-actions on $\DNC(\base,\subbase)$ and $\DNC(\basetwo,\subbasetwo)$ (see Remark \ref{rema: gauge action on DNC}).  
\end{rema}

Many of the statements made in Section \ref{sec: short review of the normal bundle functor} can be generalised to deformation to the normal cones. Later, we will list many properties of deformation to the normal cone spaces, and we will show that a pair of Lie groupoids and a pair of Lie algebroids gives rise to a deformation to the normal cone groupoid and deformation to the normal cone algebroid, respectively. First, we will digress, and introduce a new way of approaching vector bundles and discuss $\VB$-groupoids, double vector bundles, and $\VB$-algebroids. All of this theory naturally comes up when studying normal bundles and deformation to the normal cones, so it is good to be aware of the ideas behind this theory.

\subsection{Vector bundles in terms of homogeneous structures}\label{sec: different approach to vector bundles}
%A good understanding of the theory of vector bundles is obviously necessary for a good understanding of Lie algebroids. 
Fairly recently, in \cite{Grabowski_2012} and \cite{Grabowski_2009}, it was observed that there are other, simpler ways to define vector bundles. In many cases, this point of view can greatly simplify arguments, so we will quickly go through the philosophy behind these ideas and discuss some results in this direction that will come in handy later. 
%\begin{rema}\label{rema: homogeneous structure on a vector bundle}

A vector bundle $\pi: E \ra M$ comes, for all $\lambda \in \mathbb{R}$, with a scalar multiplication
\[h_\lambda: E \ra E \textnormal{ given by } (x,\xi) \mapsto (x,\lambda\xi).\]
Notice, however, that the maps $h_\lambda$ for $\lambda \in \mathbb{R}$ are completely determined by the maps $h_\lambda$ for $\lambda \in \mathbb{R}_{\ge 0}$. Of course, for all $\lambda,\mu \in \mathbb{R}_{\ge 0}$, the identity $h_\lambda \circ h_\mu = h_{\lambda\mu}$ holds, so we can view 
\[h: \mathbb{R}_{\ge 0} \times E \ra E \textnormal{ given by } h(\lambda,e) \coloneqq h_\lambda(e)\] 
as an action by the monoid $(\mathbb{R}_{\ge 0},\cdot)$. The action is smooth in the sense that we can extend it to a smooth function on an open set $U \subset \mathbb{R}$ containing $\mathbb{R}_{\ge 0}$ (and therefore to all of $\mathbb{R}$). The main idea we will exploit in this section, is that the above smooth action of $(\mathbb{R}_{\ge 0},\cdot)$ contains all of the structure that underlies the vector bundle. First of all, the map
\[h_0: E \ra E\]
maps onto $M$, which is a (closed) embedded submanifold of $E$, so it factors through a map $E \ra M$, which is of course the projection $\pi$. Second of all, for all $e \in E$, we obtain a smooth curve
\[e_h: \mathbb{R}_{\ge 0} \ra E \textnormal{ given by } t \mapsto h(t,e)\]
and these curves give rise to the smooth map
\[\mathcal{V}_h: E \ra TE \textnormal{ given by } e \mapsto (e_h(0),de_h(0)1),\]
called the \textit{vertical lift}. Since $h_\lambda|_M: M \ra M$ is the identity map for all $\lambda \in \mathbb{R}_{\ge 0}$, we see that $\mathcal{V}_h(x)=(x,0)$ for all $x \in M$. Moreover, it is not hard to see that $M$ is determined by this property, i.e. if $\mathcal{V}_h(e)=(e_h(0),0)$, where $e \in E$, then $e=e_h(0) \in M$. 

Recall the following construction:
\begin{rema}\label{rema: vertical bundle}
Let $\pi: E \ra M$ be a vector bundle. Since $\pi$ is a submersion, $d\pi: TE \ra TM$ has constant rank, so $\ker d\pi$ is a vector bundle, which is called the \textit{vertical bundle} $VE \ra E$ of $E$.
\end{rema}
The vertical lift $\mathcal{V}_h$ maps into $VE$ by the chain rule, and it is readily verified to be an isomorphism onto $VE|_M$. We will see that all of the above statements still hold if we start with a manifold $E$ together with a smooth (multiplicative) $\mathbb{R}_{\ge 0}$-monoid action for which the vertical lift $\mathcal{V}_h$ vanishes at $e \in E$ if and only $e=x \in M$. 
\begin{defn}\cite{Grabowski_2009}\label{defn: homogeneous structure}
A \textit{homogeneous structure} on a manifold $E$ is a smooth monoid action $h: \mathbb{R}_{\ge 0} \times E \ra E$ of $(\mathbb{R}_{\ge 0},\cdot)$ on $E$ (we denote $M \coloneqq h_0(E)$) such that the maps $e_h: \mathbb{R}_{\ge 0} \ra E$ given by $t \mapsto h(t,e) \eqqcolon h_t(e)$ are non-singular for all $e \not\in M$. The \textit{vertical lift} is the map $\mathcal{V}_h: E \ra TE$ given by $e \mapsto (e_h(0),de_h(0)1)$.
\end{defn}
The main result of this section is the following theorem.
\begin{theo}\cite{Grabowski_2009}\label{theo: homogeneneous structure}
Let $h: \mathbb{R}_{\ge 0} \times E \ra E$ be a homogeneous structure on a manifold $E$. Then there is a unique vector bundle structure on $E$ for which $h_\lambda$ ($\lambda \in \mathbb{R}_{\ge 0}$) is the scalar multiplication by $\lambda \in \mathbb{R}_{\ge 0}$. Moreover, a smooth map $\varphi: E \ra D$ between vector bundles is a morphism of vector bundles if and only if $\varphi \circ h_\lambda^E = h_\lambda^D \circ \varphi$ for all $\lambda \in \mathbb{R}_{\ge 0}$.
\end{theo}

The statement says that vector bundles are characterised by their scalar multiplication map. We will not prove the statement here, but, in a sentence, a sketch of the proof is as follows: if we have a homogeneous structure $h: \mathbb{R}_{\ge 0} \times E \ra E$ on a manifold $E$, then it turns out that the vertical lift determines a diffeomorphism onto its image. One can show that this image is a vector bundle (which will be the vertical bundle $VE|_M$), and then we can pull this vector bundle back to obtain a vector bundle structure on $E$. This vector bundle structure on $E$ is unique, because fiberwise the linear structure will be completely determined by $h$. 

We obtain a powerful corollary from the former statement. 
\begin{coro}\cite{Grabowski_2009}\label{coro: vector subbundle is invariant submanifold}
Let $\pi: E \ra M$ be a vector bundle. An embedded submanifold $F \subset E$ is a vector subbundle over $\pi(F)$ if and only if $F$ is invariant under scalar multiplication.
\end{coro}
\begin{proof}
It is obvious that a subbundle of $E$ is invariant under scalar multiplication. If $F \subset E$ is an embedded submanifold which is invariant under scalar multiplication, then the homogeneous structure $h: \mathbb{R}_{\ge 0} \times E \ra E$ at least restricts to a smooth $\mathbb{R}_{\ge 0}$-action on $F$, denoted by $h|_F$. Moreover,
\[\mathcal{V}_{h|_F} = \mathcal{V}_h|_F\]
simply because $e_h: \mathbb{R}_{\ge 0} \ra E$ maps into $F$ if $e \in F$. It follows that (1) $h|_F$ is a homogeneous structure on $F$, so that $F$ carries a unique vector bundle structure over $N \coloneqq h_0(F)=\pi(F)$ for which $(h|_F)_\lambda$ is the scalar multiplication by $\lambda \in \mathbb{R}_{\ge 0}$, and that (2) $F \cong VF|_N \subset VE|_M \cong E$ is a vector subbundle. The result now follows.
\end{proof}
To end the section, we state that a vector bundle comes with a unique graded algebra $C^\infty_{\textnormal{poly}}(E) \subset C^\infty(E)$ consisting of so-called homogeneous functions on $E$. 
\begin{defn}\cite{Grabowski_2009}\label{defn: homogeneous functions}
Let $E \ra M$ be a vector bundle. A homogeneous function of degree $k$ on $E$ is a smooth map $f \in C^\infty(E)$ such that $h_\lambda^*f \coloneqq f \circ h_k = \lambda^k \cdot f$. The subset of $C^\infty(E)$ consisting of homogeneous functions of degree $k$ is denoted by $C^\infty_{\textnormal{poly}}(E)^k$, and we denote $C^\infty_{\textnormal{poly}}(E) \coloneqq \textstyle\bigoplus_{k \ge 0} C^\infty_{\textnormal{poly}}(E)^k$.
\end{defn}
The subscript ``poly'' is used, because one can show that homogeneous functions on a vector bundle are in fact fiberwise polynomial (see \cite{Grabowski_2012}). Notice that the projection $\pi: E \ra M$ of a vector bundle gives rise to a map $C^\infty(M) \ra C^\infty_{\textnormal{poly}}(E)$ and it maps onto $C^\infty_{\textnormal{poly}}(E)^0$. 

We will discuss an algebraic view of the deformation to the normal cone construction at the end of this chapter. In view of this, it will be useful to observe the following fact. 
\begin{theo}\cite{leuther2012affine}\label{theo: smooth polynomial functions of vector bundle}
Let $E \ra M$ be a vector bundle. The graded algebra
\[C^\infty_{\textnormal{poly}}(E) = \bigoplus_{k \ge 0} C^\infty_{\textnormal{poly}}(E)^k \subset C^\infty(E)\]
determines the vector bundle $E$ uniquely. That is, two vector bundles $E$ and $D$ are isomorphic if and only if $C^\infty_{\textnormal{poly}}(E)$ and $C^\infty_{\textnormal{poly}}(D)$ are isomorphic as $\mathbb{R}$-algebras.
\end{theo}

\subsection{A discussion on $\VB$-groupoids and $\VB$-algebroids}\label{sec: VB groupoids and algebroids}
In this section we will briefly discuss $\VB$-groupoids and $\VB$-algebroids (read ``vector bundle" for $\VB$). The main reason for this is that every pair of Lie groupoids $(\groupoid,\subgroupoid)$ and every pair of Lie algebroids $(\algebroid, \subalgebroid)$ will induce commutative diagrams
\begin{center}
\begin{tikzcd}
\normal(\group,\subgroup) \ar[r,shift left]\ar[r,shift right] \ar[d] & \normal(\base,\subbase) \ar[d]\\
\subgroup \ar[r,shift left]\ar[r,shift right] & \subbase
\end{tikzcd} \textnormal{ }
\begin{tikzcd}
\normal(\algebr,\subalgebr) \ar[r] \ar[d] & \normal(\base,\subbase) \ar[d]\\
\subalgebr \ar[r] & \subbase,
\end{tikzcd} 
\end{center}
respectively, such that the vertical arrows are vector bundles, and the horizontal arrows are Lie groupoids/Lie algebroids. Moreover, there are several compatibility conditions satisfied. Below, we will make precise what conditions these are. 

\subsubsection{$\VB$-groupoids}\label{sec: VB groupoids}
We start with the notion of a $\VB$-groupoid. This section is based on \cite{Gracia_Saz_2017}. 
\begin{defn}\cite{Gracia_Saz_2017}\label{defn: VB-groupoid}
A \textit{$\VB$-groupoid} is a commutative diagram
\begin{center}
    \begin{tikzcd}
    \grouptwo \ar[r,shift left]\ar[r,shift right] \ar[d,"\pi_\grouptwo"] & E \ar[d,"\pi_E"]\\
    \group \ar[r,shift left]\ar[r,shift right] & \base,
    \end{tikzcd}
\end{center}
where horizontal arrows represent Lie groupoids and vertical arrows represent vector bundles, satisfying the following three conditions:
\begin{enumerate}
    \item $\source_\grouptwo$ and $\target_\grouptwo$, with base maps $\source_\group$ and $\target_\group$, respectively, are morphisms of vector bundles;
    \item $\pi_\grouptwo$, with base map $\pi_E$, is a morphism of Lie groupoids;
    \item if $(g_1,g_2),(g_3,g_4) \in \grouptwo^{(2)}$, such that $(g_1,g_3),(g_2,g_4) \in \grouptwo \tensor[_{\pi_\grouptwo}]{\times}{_{\pi_\grouptwo}} \grouptwo$, then
    \[\mult_\grouptwo(g_1 +_\grouptwo g_3, g_2 +_\grouptwo g_4) = \mult_\grouptwo(g_1,g_2) +_\grouptwo \mult_\grouptwo(g_3,g_4),\]
    called the \textit{interchange law}, where $+$ denotes the addition map of a vector bundle.
\end{enumerate}
If we require the vertical arrows to represent Lie algebroids, and, additionally, require $\source_\grouptwo$ and $\target_\grouptwo$ to be morphisms of Lie algebroids, then the $\VB$-groupoid is called an $\LA$-groupoid.
\end{defn}
\begin{exam}\cite{Gracia_Saz_2017}\label{exam: Tgroup is VB-groupoid}
For any Lie groupoid $\groupoid$, the tangent prolongation groupoid (see Example \ref{exam: tangent groupoid}), represented as follows:
\begin{center}
    \begin{tikzcd}
    T\group \ar[r,shift left]\ar[r,shift right] \ar[d] & T\base \ar[d]\\
    \group \ar[r,shift left]\ar[r,shift right] & \base,
    \end{tikzcd}
\end{center}
is a $\VB$-groupoid. Indeed, the first two conditions of Definition \ref{defn: VB-groupoid} are obvious, and the interchange law is satisfied, because in this case it means that, for all $(g,h) \in \group^{(2)}$, $d\mult_\group(g,h)$ is linear. By relaxing the Hausdorff assumptions we made for Lie algebroids, the maps $d\source_\group$ and $d\target_\group$ are morphisms of Lie algebroids, so then it is even an $\LA$-groupoid.
\end{exam}
For us, the important part will be that we can extract a natural vector bundle, called the \textit{core}, out of a $\VB$-groupoid that encodes a lot of information about the groupoid $\grouptwo \rra E$ (using notation from Definition \ref{defn: VB-groupoid}). In the remaining of this section we will be studying this particular vector bundle. Be aware that, in the case of $\LA$-groupoids, this vector bundle carries the structure of a Lie algebroid. We will see that this Lie algebroid for the tangent prolongation $\VB$-groupoid $T\group \rra T\base$ is precisely (or at least naturally isomorphic to) the Lie algebroid of $\groupoid$. We will, however, only discuss the vector bundle structure of the so-called core. Using the following lemma, we can show how to construct this vector bundle. We fix a $\VB$-groupoid using the notation from Definition \ref{defn: VB-groupoid}. 
\begin{lemm}\cite{surjectivesubm}\label{lemm: submersion of pullback in VB-groupoid}
The maps
\[\pi^R \coloneqq (\pi_\grouptwo,\source_\grouptwo): \grouptwo \ra \group \tensor[_{\source_\group}]{\times}{_{\pi_E}} E \textnormal{ and } \pi^L \coloneqq (\pi_\grouptwo,\target_\grouptwo): \grouptwo \ra \group \tensor[_{\target_\group}]{\times}{_{\pi_E}} E\]
are surjective submersions (note: $\group \tensor[_{\source_\group}]{\times}{_{\pi_E}} E = \source_\group^*E$ and $\group \tensor[_{\target_\group}]{\times}{_{\pi_E}} E = \target_\group^*E$).
\end{lemm}
\begin{proof}
Notice that, indeed, $\group \tensor[_{\source_\group}]{\times}{_{\pi_E}} E \textnormal{ and } \group \tensor[_{\target_\group}]{\times}{_{\pi_E}} E$ are smooth manifolds, because $\source_\group,\target_\group$ and $\pi_E$ are submersions. Observe that we only need to show that $\source_\grouptwo|_{\pi_\grouptwo^{-1}(g)}$ and $\target_\grouptwo|_{\pi_\grouptwo^{-1}(g)}$ are surjective (to show that the maps are surjectivity) and similarly for $d\source_\grouptwo$ and $d\target_\grouptwo$ (to show that the maps are submersive). However, since
\begin{center}
\begin{tikzcd}
    T\grouptwo \ar[r,shift left]\ar[r,shift right] \ar[d,"d\pi_\grouptwo"] & TE \ar[d,"d\pi_E"]\\
    T\group \ar[r,shift left]\ar[r,shift right] & T\base
\end{tikzcd}
\end{center}
is again a $\VB$-groupoid, simply by applications of the chain rule, it suffices to prove the statement only for $\source_\grouptwo$ and $\target_\grouptwo$ (by replacing $\source_\grouptwo$ and $\target_\grouptwo$ with $d\source_\grouptwo$ and $d\target_\grouptwo$). %Indeed, then the maps $\grouptwo_g \ra E_{\source_\group(g)}$ and $\grouptwo_g \ra E_{\target_\group(g)}$ that are induced by $(\pi_\grouptwo,\source_\grouptwo) \textnormal{ and } (\pi_\grouptwo,\target_\grouptwo)$, respectively, are surjective, and so  

Now, recall that there are canonical vector bundle isomorphisms $T\grouptwo|_\group \cong \grouptwo \oplus T\group$ and $TE|_\base \cong E \oplus T\base$ (see Lemma \ref{lemm: canonical iso TE|_M}). Since $\source_\grouptwo$ and $\target_\grouptwo$ are morphisms of vector bundles over $\source_\group$ and $\target_\group$, respectively, the differential of $\source_\grouptwo$ and $\target_\grouptwo$ (restricted to $\group$) can be decomposed as maps
\[d\source_\grouptwo|_\group,d\target_\grouptwo|_\group: \grouptwo \oplus T\group \ra E \oplus T\base\]
by using the above isomorphisms. By tracing the maps, we see that these maps are given by 
\[((g,\eta),(g,\xi)) \mapsto (\source_\grouptwo(g,\eta),d\source_\group(g)\xi), \textnormal{ and } ((g,\eta),(g,\xi)) \mapsto (\target_\grouptwo(g,\eta),d\target_\group(g)\xi),\] 
respectively (note: by $(g,\eta) \in \grouptwo$ we mean $\eta \in \grouptwo_g \eqqcolon \pi_\grouptwo^{-1}(g)$). These maps are also surjective, because $\source_\grouptwo$ and $\target_\grouptwo$ are submersions, so, in particular, $\source_\grouptwo|_{\pi_\grouptwo^{-1}(g)}$ and $\target_\grouptwo|_{\pi_\grouptwo^{-1}(g)}$ are surjective, as required. 
\end{proof}
Observe that $\source_\group^*E = \group \tensor[_{\source_\group}]{\times}{_{\pi_E}} E$ and $\target_\group^*E = \group \tensor[_{\target_\group}]{\times}{_{\pi_E}} E$ are vector bundles over $\group$, and that $\pi^R$ and $\pi^L$ are morphisms of vector bundles over $\id_\group$. In particular, the kernel $E^R \ra \group$ of $\pi^R$ and the kernel $E^L \ra \group$ of $\pi^L$ are vector bundles. In turn, we can define the following vector bundles.
\begin{defn}\cite{Gracia_Saz_2017}\label{defn: core of VB-groupoid}
The restriction $C^R \ra \base$ of $E^R$ to $\identity_\group(\base)$ is called the \textit{right-core}. Similarly, the restriction $C^L \ra \base$ of $E^L$ to $\identity_\group(\base)$ is called the \textit{left-core}.
\end{defn}
\begin{exam}\cite{Gracia_Saz_2017}\label{exam: tangent prolongation core}
Consider the tangent prolongation $\VB$-groupoid (see Example \ref{exam: Tgroup is VB-groupoid}). Then $\pi^R: T\group \ra \group \tensor[_{\source_\group}]{\times}{_{\pi_{T\base}}} T\base$ is given by $(g,\eta) \mapsto (g,(\source_\group(g),d\source_\group(g)\eta))$, so the kernel of this maps is precisely $\ker d\source$. Restricting it to $\identity_\group(\base)$ yields precisely the Lie algebroid of $\groupoid$. The vector bundle $C^L \ra \base$ is the Lie algebroid of $\groupoid$ treated as $\ker d\target_\group|_{\identity(\base)}$ via left-invariant vector fields (see Remark \ref{rema: other definitions of Lie algebroid of a Lie groupoid}).
\end{exam}
In the above example, the left-core and the right-core of a $\VB$-groupoid are naturally isomorphic. This is a consequence of the following more general statement.
\begin{prop}\cite{Gracia_Saz_2017}\label{prop: natural isomorphism cores}
The \textit{involution of the $\VB$-groupoid}
\[-\inv_\grouptwo: \grouptwo \xra{\sim} \grouptwo, \textnormal{ given by } g \mapsto -g^{-1}\]
is an isomorphism of vector bundles over $\inv_\group$. Moreover, it restricts to isomorphisms of vector bundles $C^R \xra{\sim} C^L$ and $C^L \xra{\sim} C^R$.
\end{prop}
\begin{proof}
Obviously, it suffices to prove the statement for $\inv_\grouptwo$ instead of $-\inv_\grouptwo$. To show that $\inv_\grouptwo$ is a morphism of vector bundles, i.e. that the induced maps $\grouptwo_g \ra \grouptwo_{g^{-1}}$ (here, $g \in \group$) are linear, notice that $(g,\eta_1) + (g,\eta_2)$ is sent to $\inv_\grouptwo((g,\eta_1) + (g,\eta_2))$. Then, by the interchange law, 
\begin{align*}
    \mult_\grouptwo((g,\eta_1) + (g,\eta_2),\inv_\grouptwo(g,\eta_1) + \inv_\grouptwo(g,\eta_2)) &= \mult_\grouptwo((g,\eta_1),\inv_\grouptwo(g,\eta_1)) \\ &+ \mult_\grouptwo((g,\eta_2),\inv_\grouptwo(g,\eta_2)) \\
    &= \id_{\grouptwo,\target_\grouptwo(g,\eta_1)} + \id_{\grouptwo,\target_\grouptwo(g,\eta_2)} \\
    &= \id_{\grouptwo,\target_\grouptwo(g,\eta_1+\eta_2)},
\end{align*}
where we used in the last equality that $\target_\grouptwo$ is a morphism of vector bundles. This proves that $\inv_\grouptwo((g,\eta_1) + (g,\eta_2)) = \inv_\grouptwo(g,\eta_1) + \inv_\grouptwo(g,\eta_2)$, so $\inv_\grouptwo$ is a vector bundle isomorphism. The second statement follows by the observation that
\[\inv_\grouptwo(c^r) = \id_{\target_{\grouptwo}(c^r)}-c^r \textnormal{ and } \inv_\grouptwo(c^l) = \id_{\source_{\grouptwo}(c^l)}-c^l, \textnormal{ where } c^r \in C^R \textnormal{ and } c^l \in C^L,\]
by the interchange law.
\end{proof}
By the above statement, we will use the following terminology.
\begin{term}\label{term: core of a VB-groupoid}
We will refer to either $C^R \ra \base$ or $C^L \ra \base$ as the \textit{core} and simply denote either one by $C \ra \base$. If we want, or need, to be explicit, then we will always specify if we mean $C^R$ or $C^L$.
\end{term}

\subsubsection{$\VB$-algebroids}\label{sec: DVB}
In this section we will mostly be looking at \textit{double vector bundles}, since most of the theory that we will use relies on the machinery that is provided by double vector bundles. At the end we will state the most important result for us regarding $\VB$-algebroids. This section is based on \cite{Gracia_Saz_2010} and \cite{Urbanski}. 
\begin{defn}\cite{Gracia_Saz_2010}\label{defn: double vector bundle}
A \textit{double vector bundle} (DVB) is a commutative diagram 
\begin{center}
\begin{tikzcd}
D \ar[r, "\pi_D^E"] \ar[d, "\pi_D^\algebr"] & E \ar[d, "\pi_E"]\\
\algebr \ar[r,"\pi_\algebr"] & \base,
\end{tikzcd}
\end{center}
where every arrow is a vector bundle, satisfying the following conditions
\begin{enumerate}
    \item $\pi_D^\algebr$ is a morphism of vector bundles over $\pi_E$, and $\pi_D^E$ is a morphism of vector bundles over $\pi_\algebr$;
    \item $+_D^\algebr: D \tensor[_{\pi_D^\algebr}]{\times}{_{\pi_D^\algebr}} D \ra D$ is a morphism of vector bundles over $+_E:  E \tensor[_{\pi_E}]{\times}{_{\pi_E}} E \ra E$.
\end{enumerate}
If, in addition, $D \ra E$ and $\algebr \ra \base$ are Lie algebroids, $\pi^\algebr_D$ is a morphism of Lie algebroids over $\pi_E$, and $+_D^\algebr$ is a morphism of Lie algebroids over $+_E$, then the DVB is called a $\VB$-algebroid. 
\end{defn}
\begin{term}\label{term: DVB}
In the situation from Definition \ref{defn: double vector bundle}, we will write $0_E^\base: \base \ra E$ and $0_\algebr^\base: \base \ra \algebr$ for the zero sections.
\end{term}
\begin{rema}\label{rema: ambiguity VB-algebroid}
In order for the definition of a $\VB$-algebroid to be well-defined, we should check the following: using the notation from the above definition, if the vector bundles $D \ra E$ and $\algebr \ra \base$ are Lie algebroids, and $\pi^\algebr_D$ is a morphism of Lie algebroids over $\pi_E$, then $D \tensor[_{\pi_D^\algebr}]{\times}{_{\pi_D^\algebr}} D \ra E \tensor[_{\pi_E}]{\times}{_{\pi_E}} E$ is again (naturally) a Lie algebroid. Indeed, this follows from the fact that 
\[\varphi = \pi^\algebr_D \times \pi^\algebr_D: D \times D \ra \algebr \times \algebr,\]
which is a Lie algebroid morphism by assumption, has a clean intersection with the diagonal $\Delta_\algebr$. Therefore, $\gr \varphi \cap (D \times D \times \Delta_\algebr) \cong D \tensor[_{\pi_D^\algebr}]{\times}{_{\pi_D^\algebr}} D \subset D \times D$ is a Lie subalgebroid by Proposition \ref{prop: clean intersection of Lie algebroids} (note: we use here our assumptions on Hausdorffness on the Lie algebroid $\algebr$ to ensure that the graph of $\varphi$ and $\Delta_\algebr$ are closed).
\end{rema}
\begin{exam}\cite{Gracia_Saz_2010}\label{exam: tangent DVB}
Let $(E \xra{\pi_E} \base,F \xra{\pi_F} \subbase)$ be a pair of vector bundles. Then,
\begin{center}
\begin{tikzcd}
TE \ar[r, "d\pi_E"] \ar[d, "\pi_{TE}"] & T\base \ar[d, "\pi_{T\base}"]\\
E \ar[r,"\pi_E"] & \base
\end{tikzcd} \textnormal{ and }
\begin{tikzcd}
\normal(E,F) \ar[r, "d_\normal\pi_E"] \ar[d, "\pi_{\normal(E,F)}"] & \normal(\base,\subbase) \ar[d, "\pi_{\normal(\base,\subbase)}"]\\
F \ar[r,"\pi_F"] & \subbase
\end{tikzcd}
\end{center}
are readily verified to be double vector bundles (see Proposition \ref{prop: normal bundle of vector bundle} for the vector bundle structures on $TE \ra T\base$ and $\normal(E,F) \ra \normal(M,N)$). We will refer to the left DVB as the \textit{tangent prolongation} of $E$, and to the right DVB as the \textit{normal double bundle} of $(E,F)$.
\end{exam}
For DVB's we also have an appropriate notion of a core vector bundle. To explain the construction, we use two intermediate results. For the rest of this section, we fix a DVB as in Definition \ref{defn: double vector bundle}. 
\begin{lemm}\cite{Gracia_Saz_2010}\label{lemm: DVB equivalence definitions}
The following properties are satisfied for the DVB: 
\begin{enumerate}
    \item if $(d_1,d_2) \in D \tensor[_{\pi^E_D}]{\times}{_{\pi^E_D}} D$, and $(d_3,d_4) \in D \tensor[_{\pi^\algebr_D}]{\times}{_{\pi^\algebr_D}} D$, then
    \[\pi_D^\algebr(d_1 +_D^E d_2) = \pi_D^\algebr(d_1) +_\algebr \pi_D^\algebr(d_2), \textnormal{ and } \pi_D^E(d_3 +_D^\algebr d_4) = \pi_D^E(d_3) +_E \pi_D^E(d_4);\]
    \item if $(d_1,d_2), (d_3,d_4) \in D \tensor[_{\pi^E_D}]{\times}{_{\pi^E_D}} D$, such that $(d_1,d_3), (d_2,d_4) \in D \tensor[_{\pi^\algebr_D}]{\times}{_{\pi^\algebr_D}} D$, then
    \[(d_1 +_D^E d_2) +_D^\algebr (d_3 +_D^E d_4) = (d_1 +_D^\algebr d_3) +_D^E (d_2 +_D^\algebr d_4).\]
\end{enumerate}
We will refer to the property $2.$ as the \textit{interchange law}.
\end{lemm}
\begin{proof}
The first statement is a reformulation that $\pi^\algebr_D$ and $\pi^E_D$ are morphisms of vector bundles. Indeed, to see that, if $1.$ is satisfied, these maps intertwine the scalar multiplication for all $\lambda \in \mathbb{R}$, notice first that this holds for all $\lambda \in \mathbb{N}$. It then also holds for all $\lambda \in \mathbb{Q}$ (e.g. $\tfrac{1}{n}\pi^\algebr_D(d) = \tfrac{1}{n}\textstyle\sum_{j=1}^n\pi^\algebr_D(\tfrac{1}{n}d) = \pi^\algebr_D(\tfrac{1}{n}d)$), and therefore for all $\lambda \in \mathbb{R}$ by continuity. Similarly, the second statement is a reformulation that $+^\algebr_D$ is a morphism of vector bundles (note: we see that this statement is equivalent to $+^E_D$ being a vector bundle morphism over $+_\algebr$). This concludes the proof.
\end{proof}
\begin{lemm}\label{lemm: submersion D to A plus E}
The map
\[(\pi^\algebr_D,\pi^E_D): D \ra \algebr \oplus E\]
is a surjective submersion.
\end{lemm}
\begin{proof}
The proof is analogous to the proof of Lemma \ref{lemm: submersion of pullback in VB-groupoid}. Since 
\begin{center}
\begin{tikzcd}
TD \ar[r, "d\pi_D^E"] \ar[d, "d\pi_D^\algebr"] & TE \ar[d, "d\pi_E"]\\
T\algebr \ar[r,"d\pi_\algebr"] & T\base
\end{tikzcd}
\end{center}
is again a DVB (by applications of the chain rule), we see that (by the symmetry of the problem) it suffices to prove that $\pi^\algebr_D|_{(\pi^E_D)^{-1}(e)}: (\pi^E_D)^{-1}(e) \ra A$ is surjective for all $e \in E$. Recall that there are canonical vector bundle isomorphisms $TD|_E \cong D \oplus TE$ and $T\algebr|_\base \cong \algebr \oplus T\base$ (see Lemma \ref{lemm: canonical iso TE|_M}). Since $\pi^\algebr_D$ is a morphism of vector bundles over $\pi_E$, the differential of the former map restricts to a morphism of vector bundles
\[D \oplus TE \ra \algebr \oplus T\base.\]
Tracing the maps, we see that this map is given by 
\[((e,\eta),(e,\xi)) \mapsto (\pi^\algebr_D(e,\eta),d\pi_E(e)\xi)\]
(note: by $(e,\eta) \in D$ we mean that $\eta \in D_e$). Since $\pi^\algebr_D$ is a submersion, this map is surjective, so, in particular, we see that $\pi^\algebr_D|_{(\pi^E_D)^{-1}(e)}$ is surjective, as required.  
\end{proof}
\begin{defn}\cite{Gracia_Saz_2010}\label{defn: core of a DVB}
The \textit{core} $\pi_C: C \ra \base$ is the vector bundle $\ker \pi^\algebr_D \cap \ker \pi^E_D = \ker(\pi^\algebr_D,\pi^E_D)$.
\end{defn}
\begin{rema}\label{rema: core of a DVB}
Observe that the vector space structures on $\pi^E_D: D \ra E$ and $\pi^\algebr_D$ agree on $\pi_C: C \ra \base$. To see this, notice that
\[0_C^\base: \base \ra C\]
is equal to $0_D^E \circ 0_E^\base$ and also to $0_D^\algebr \circ 0_\algebr^\base$:
\begin{align*}
(0_D^\algebr \circ 0_\algebr^\base) +_D^E (0_D^\algebr \circ 0_\algebr^\base) &= (0_D^\algebr \circ 0_\algebr^\base +_D^\algebr 0_D^\algebr \circ 0_\algebr^\base) +_D^E (0_D^\algebr \circ 0_\algebr^\base +^\algebr_D 0_D^\algebr \circ 0_\algebr^\base) \\
&= (0_D^\algebr \circ 0_\algebr^\base +_D^E 0_D^\algebr \circ 0_\algebr^\base) +_D^\algebr (0_D^\algebr \circ 0_\algebr^\base +_D^E 0_D^\algebr \circ 0_\algebr^\base),
\end{align*}
by the interchange law, so $0^\algebr_D \circ 0_\algebr^\base=(0_D^\algebr \circ 0_\algebr^\base) +_D^E (0_D^\algebr \circ 0_\algebr^\base)$, and therefore $0^\algebr_D \circ 0^\base_\algebr = 0^E_D \circ 0^\base_\algebr$. Now we can see that both vector bundle structures that $C$ inherits coincide: let $c,d \in C$. Then
\[c +^\algebr_D d = (c +^E_D 0_C^\base) +^\algebr_D (0_C^\base +^E_D d) = (c +^\algebr_D 0_C^\base) +^E_D (0_C^\base +^\algebr_D d) = c +^E_D d,\]
by the interchange law. By a similar argument as in the proof of Lemma \ref{lemm: DVB equivalence definitions} this also shows that the scalar multiplications of $\pi_D^E: D \ra E$ and $\pi_D^\algebr: D \ra \algebr$ on $C$ coincide.
\end{rema}
In the rest of this section, we will talk about the sections of $D \ra E$, and what we can say about them in the case of a $\VB$-algebroid. 
\begin{term}\label{term: sections of D}
We will denote by $\Gamma^E(D)$ the sections of $\pi^E_D: D \ra E$ and by $\Gamma^\algebr(D)$ the sections of $\pi^\algebr_D: D \ra \algebr$.
\end{term}
First, we will show that the sections of $D \ra E$ are determined by sections of $C \ra \base$, which will be called \textit{core sections}, and the so-called \textit{linear sections}. 
\begin{lemm}\cite{Gracia_Saz_2010}\label{lemm: core sections}
Let $\alpha \in \Gamma(C)$. The map
\[\alpha_C \coloneqq \iota_C \circ \alpha \circ \pi_E +^\algebr_D 0^E_D: E \ra D,\]
where $\iota_C: C \hookrightarrow D$ is the inclusion, is a section of $\pi^E_D: D \ra E$. 
\end{lemm}
\begin{proof}
We only have to check that $\pi^E_D \circ \alpha_C = \id_E$:
\begin{align*}
    \pi^E_D \circ (\iota_C \circ \alpha \circ \pi_E +^\algebr_D 0^E_D) &= \pi^E_D \circ \iota_C \circ \alpha \circ \pi_E +_E \pi^E_D \circ 0^E_D \\
    &= 0_E^\base \circ \pi_E +_E \id_E = \id_E
\end{align*}
where we used Lemma \ref{lemm: DVB equivalence definitions} in the first equality and that $C \subset \ker \pi^E_D$ in the second equality.
\end{proof}
\begin{defn}\cite{Gracia_Saz_2010}\label{defn: core sections}
The sections $\alpha_C \in \Gamma^E(D)$, for $\alpha \in \Gamma(C)$, are called \textit{core sections}. The sections $\beta \in \Gamma^E(D)$, which are also vector bundle morphisms $E \ra D$, are called \textit{linear sections}. The collection of core sections will be denoted by $\Gamma_C(D)$ (or $\Gamma_C^E(D)$), and the collection of linear sections will be denoted by $\Gamma_\ell(D)$ (or $\Gamma_\ell^E(D)$).
\end{defn}
\begin{exam}\cite{Gracia_Saz_2010}\label{exam: core of TE and normal bundle}
Recall the DVB's from Example \ref{exam: tangent DVB}. Recall that we can construct the so-called vertical bundle $VE \ra \base$ of $\pi_E: E \ra \base$, which is defined as $d\pi^{-1}(0^\base_{T\base})$. In terms of this vector bundle, we can describe the core of the tangent prolongation DVB of $E \ra \base$ as those elements in $VE \ra \base$ which are tangent to the zero section $0^E_{TE}$. We also see that the core sections $\Gamma_C(TE)$ are fiberwise constant, so they are the fiberwise constant and vertical vector fields. The linear sections $\Gamma_\ell(TE)$ are just the linear vector fields. Notice that the vertical bundle is in fact naturally isomorphic to $E \tensor[_{\pi_E}]{\times}{_{\pi_E}} E \ra \base$ via the isomorphism 
\[E \tensor[_{\pi_E}]{\times}{_{\pi_E}} E \ra VE \textnormal{ given by } (e_1,e_2) \mapsto (e_1,\left.\tfrac{d}{dt}\right\vert_{t=0}e_1+te_2),\]
(note: by the chain rule it indeed maps into $VE$), so we see that we can realise the core of the tangent prolongation as $E \ra \base$. In particular, we have $\Gamma(E) \cong \Gamma_C(TE)$.

Using the identification that the tangent prolongation has core $E \ra \base$, it is readily verified that the core of the double normal bundle is $E|_\subbase/F \ra \subbase$. Another way to see this, is by the sequence of vector bundle isomorphisms
\[\normal(E,F)|_\subbase = TE|_\subbase/TF|_\subbase \cong (E|_\subbase \oplus T\base|_\subbase)/(F \oplus T\subbase) \cong E|_\subbase/F \oplus \normal(\base,\subbase),\]
where we used Lemma \ref{lemm: canonical iso TE|_M} two times.
We now have $\Gamma_C(\normal(E,F)) \cong \Gamma(E|_\subbase/F)$ and the space of sections $\Gamma_\ell(\normal(E,F))$ consists of the linear vector fields (with respect to $d_\normal$ instead of $d$).
\end{exam}
It is now a good time to explain the local nature of DVB's. To do this, we will use the following result.
\begin{prop}\cite{Urbanski}\label{prop: ker of DVB as sum and pullback}
The vector bundle $\ker \pi^E_D$, with respect to the vector bundle structure on $\pi^E_D: D \ra E$, is canonically isomorphic to $\algebr \oplus C$. The vector bundle $\ker \pi^E_D$, with respect to the vector bundle structure on $\pi^\algebr_D: D \ra \algebr$, is canonically isomorphic to $\algebr \tensor[_{\pi_\algebr}]{\times}{_{\pi_\algebr \circ \pi_D^\algebr}} C$. Similar results hold for $\ker \pi^\algebr_D$. Moreover, if $d_1,d_2 \in D$ with $\pi^E_D(d_1)=\pi^E_D(d_2)=e$ and $\pi^\algebr_D(d_1)=\pi^\algebr_D(d_2)=a$, then there is a $c \in C$ such that 
\[d_1 = d_2 +_D^E (e,c) = d_2 +_D^\algebr (a,c),\]
with respect to the decompositions from above.
\end{prop}
\begin{proof}
Observe that the zero section $0^\algebr_D: \algebr \ra D$ is a morphism of vector bundles with respect to the vector bundle structure $\pi^E_D: D \ra E$ on $D$. In particular, $0^\algebr_D(\algebr)$ is a vector subbundle contained in $\ker \pi^E_D$, and we can identify it with $\algebr$. We will now define two maps
\[p_\algebr,p_C: \ker \pi^E_D \ra \ker \pi^E_D\]
whose images will provide the desired splitting of $\ker \pi^E_D$. Let $d \in \ker \pi^E_D$, and let $x \coloneqq \pi_E \circ \pi^E_D(d) = \pi_\algebr \circ \pi^\algebr_D(d)$. Then
\[\pi^E_D \circ 0^\algebr_D \circ \pi^\algebr_D(d) = 0^\base_E \circ \pi_\algebr \circ \pi^\algebr_D(d) = 0^\base_E(x) = \pi^E_D(d),\]
so $(d,0^\algebr_D \circ \pi^A_D(d)) \in D \tensor[_{\pi^E_D}]{\times}{_{\pi^E_D}} D$. We now define the maps $p_\algebr$ and $p_C$ by setting
\[p_\algebr(d) \coloneqq 0^\algebr_D(\pi^\algebr_D(d)); \quad p_C(d) \coloneqq d -^E_D 0^\algebr_D(\pi^\algebr_D(d)).\]
This is well-defined, precisely by the calculation above. It is clear that $p_\algebr + p_C = \id_{\ker \pi^E_D}$, and the relations $p_\algebr \circ p_\algebr=p_\algebr$, $p_\algebr \circ p_C = 0$ and $p_C \circ p_C = p_C$ are readily verified. Moreover,
\[\pi^\algebr_D(d-^E_D 0^\algebr_D(\pi^\algebr_D(d))) = \pi^\algebr_D(d) -^\algebr_D \pi^\algebr_D(0^\algebr_D(\pi^\algebr_D(d)) = 0^\base_\algebr(x),\]
where $x \coloneqq \pi_E \circ \pi^E_D(d) = \pi_\algebr \circ \pi^\algebr_D(d)$, so we see that $p_C$ maps into $C$. It is clear that $p_D$ maps into $\algebr$. The intersection of a fiber over $x \in X$ in $0^\algebr_D(\algebr)$ and a fiber over $x \in X$ in $C$ is $\{0^E_D \circ 0^\base_E(x)\}=\{0^\algebr_D \circ 0^\base_\algebr(x)\}$, so it now follows that $p_\algebr(\ker \pi^E_D) \oplus p_C(\ker \pi^E_D) = \algebr \oplus C$, hence $\ker \pi^E_D \cong A \oplus C$ with respect to the vector bundle structure $\pi^E_D: D \ra E$. 

To see that $\ker \pi^E_D \cong \algebr \tensor[_{\pi_\algebr}]{\times}{_{\pi_\algebr \circ \pi_C^\algebr}} C$ as vector bundles, notice that, by the former statement, we already know that this identity holds as manifolds. Now, the zero section $0^\algebr_D$ maps into $\ker \pi^E_D$, and under the above identification, this map becomes the map $a \mapsto (a,0)$. The map $+^\algebr_D$ restricts to a map $\ker \pi^E_D \tensor[_{\pi^\algebr_D}]{\times}{_{\pi^\algebr_D}} \ker \pi^E_D \ra \ker \pi^E_D$, and this map becomes a map
\[\algebr \oplus C \oplus C \ra \algebr \oplus C,\]
which is given by (note: $+^\algebr_D$ is a vector bundle morphism with respect to the vector bundle structure on $\pi^E_D: D \ra E$)
\begin{align*}
    (a,c,c') \mapsto +^\algebr_D(a,c,c') &= +^\algebr_D((a,0,0) +^E_D (0,c,c')) \\
    &= (a,0) +^E_D (0,c+c') = (a,c+c').
\end{align*}
This proves that, indeed, $\ker \pi^E_D \cong \algebr \tensor[_{\pi_\algebr}]{\times}{_{\pi_\algebr \circ \pi_C^\algebr}} C$ as vector bundles. By the symmetry of the problem, similar statements hold for $\ker \pi^\algebr_D$.

It remains to prove the last statement. By the assumptions on $d_1$ and $d_2$, we have $\pi^E_D(d_1-^\algebr_D d_2) = e -_E e = 0$ and $\pi^\algebr_D(d_1 -^E_D d_2)=0$. Therefore, we can find $c,c' \in C$ with
\[d_1-^E_Dd_2 = (e,c); \quad d_1-^\algebr_Dd_2 = (a,c').\]
It suffices to show now that $c=c'$. And indeed, 
\[(d_1,d_2) = (d_1 -^\algebr_D d_2, d_2 -^\algebr_D d_1) +^\algebr_D (d_2,d_1),\]
so, by applying $-^E_D$, we see that
\begin{align*}
    (e,c) = d_1 -^E_D d_2 &= ((d_1 -^\algebr_D d_2) -^E_D (d_2 -^\algebr_D d_1)) +^\algebr_D (d_2-^E_Dd_1) \\
    &= ((a,c')-^E_D(a,-c'))+^\algebr_D(e,-c),
\end{align*}
from which it follows that $2c=2c'$, and so $c=c'$, as required.
\end{proof}
We will see below that the above proposition shows that, locally, we can decompose $D$ as $A \oplus E \oplus C$. In fact, one can show that such a decomposition always exists globally (non-canonically) (see e.g. \cite{Grabowski_2009}), but we will not go into this. 
\begin{rema}\cite{Urbanski}\label{rema: local nature DVB}
Let $U \subset \base$ be an open subset on which we have local coordinates $(x^1,\dots,x^n)$, such that on $\algebr$ and $E$ we have vector bundle coordinates 
\[(x^1,\dots,x^n,a^1,\dots,a^{n_\algebr}) \textnormal{ and } (x^1,\dots,x^n,e^1,\dots,e^{n_E}),\] 
respectively. From the last statement of Proposition \ref{prop: ker of DVB as sum and pullback}, we see that we can find coordinates $(x^1,\dots,x^n,c^1,\dots,c^{n_C})$ such that 
\[(x^1,\dots,x^n,a^1,\dots,a^{n_\algebr},e^1,\dots,e^{n_E},c^1,\dots,c^{n_C})\]
are vector bundle coordinates on $D$. We also see that
\[c^j(d_1 +^E_D d_2) = c^j(d_1) + c^j(d_2); \quad c^j(d_1 +^E_D d_2) = c^j(d_1) + c^j(d_2),\]
and $c^j \circ 0^E_D = 0$; $c^j \circ 0^\algebr_D=0$, so $(x^1,\dots,x^n,c^1,\dots,c^{n_C})$ are vector bundle coordinates on $C$.

Now, let $\alpha_1,\dots,\alpha_{n_\algebr},\gamma_1,\dots,\gamma_{n_C}$ be the local frame of $\pi^E_D: D \ra E$ determined by the coordinates $(a^1,\dots,a^{n_\algebr},c^1,\dots,c^{n_C})$. For a linear combination of these sections to be a linear section, such a linear combination has to equal a section on $\pi_\algebr: \algebr \ra \base$ when restricted to $\base=\{e^j=0\}$ (note: this is w.r.t. the coordinates on $E$), so a section $\beta \in \Gamma(D)$ is linear if and only if it locally takes the form
\[\beta = \sum_i (f^i \alpha_i + \sum_j g^i_j e^j \gamma_i), \textnormal{ where } f_i,g^i_j \in C^\infty(\base).\]
It is clear that a section $\alpha \in \Gamma(D)$ is a core section if and only if it locally takes the form
\[\alpha = \sum_k h^k \gamma_k \textnormal{ where } h^k \in C^\infty(\base).\]
It is important to notice that this shows that the sections of $\pi_D^E: D \ra E$ are generated by the linear and the core sections. 
\end{rema}
Recall the definition of a $\VB$-algebroid (see Definition \ref{defn: double vector bundle}). Since the DVB's we come across will have the added structure of a $\VB$-algebroid, we end this section with stating the most important result for us (in fact, it is an equivalent definition of a $\VB$-algebroid) for which we need the concepts of linear and core sections. A proof can be found in \cite{Gracia_Saz_2010}.
\begin{theo}\cite{Gracia_Saz_2010}\label{theo: alternative defintion of VB-algebroid}
A $\VB$-algebroid (as in Definition \ref{defn: double vector bundle}) satisfies the following two properties:
\begin{enumerate}
    \item the anchor map $\anchor_D: D \ra TE$ is a morphism of vector bundles over $\anchor_\algebr: \algebr \ra T\base$;
    \item $[\Gamma_\ell(D),\Gamma_\ell(D)]_D \subset \Gamma_\ell(D)$; \quad $[\Gamma_\ell(D),\Gamma_C(D)]_D \subset \Gamma_C(D)$; \quad $[\Gamma_C(D),\Gamma_C(D)]_D = 0$.
\end{enumerate}
Conversely, if a DVB (as in Definition \ref{defn: double vector bundle}), for which $D \ra E$ is a Lie algebroid, satisfies the above properties, then the DVB is a $\VB$-algebroid (in particular, $\algebr \ra \base$ inherits a Lie algebroid structure).
\end{theo}
%The definition suggests that the orientation of the diagram can matter. Luckily, this is not the case.
%\begin{prop}\label{prop: double vector bundle equivalent definition}
%Suppose we have a DVB
%\begin{center}
%\begin{tikzcd}
%D \ar[r, "\pi_D^E"] \ar[d, "\pi_D^F"] & E \ar[d, "\pi_E"] & \phantom{A} \ar[d,phantom,"\text{ or }"] & D \ar[r, "\pi_D^F"] \ar[d, "\pi_D^E"] & F \ar[d,"\pi_F"]\\
%F \ar[r,"\pi_F"] & Y & \phantom{A} & E \ar[r, "\pi_E"] & Y,
%\end{tikzcd}
%\end{center}
%then the other diagram is a DVB as well.
%\end{prop}
%\begin{proof}
%By symmetry of the problem, it suffices to prove that the conditions 1 and 2 in Definition \ref{defn: double vector bundle} are equivalent to the following condition: 

%$(*)$ Let $d_1,d_2,d_3,d_4 \in D$ with $(d_1,d_2), (d_3,d_4) \in D \times_E D$ and $(d_1,d_3), (d_2,d_4) \in D \times_F D$. Then
%\begin{align*}
%    \pi_D^F(d_1 +_D^E d_2) = \pi_D^F(d_1) +_F \pi_D^F(d_2) &, \text{ } \pi_D^E(d_1 +_D^F d_3) = \pi_D^E(d_1) +_E \pi_D^E(d_3) \text{, and } \\
%    (d_1 +_D^E d_2) +_D^F (d_3 +_D^E d_4) &= (d_1 +_D^F d_3) +_D^E (d_2 +_D^F d_4).
%\end{align*}
%Suppose $(*)$ holds. Then for all $d_1,d_2,d_3,d_4 \in D$ as in $(*)$ we have
%\begin{align*}
%    \pi^E_D(d_1 +^F_D d_2) &= \pi^E_D(d_1) +_E \pi^E_D(d_2) \text{ and } \pi^E_D(\lambda \cdot^F_D d_1) = \lambda \cdot_E \pi^E_D(d_1), \text{ and}.
%\end{align*}
%where the second statement follows from the first: it holds for all $\lambda \in \mathbb{Z}$, hence for all $\lambda \in \mathbb{Q}$, and so for all $\lambda \in \mathbb{R}$.
%\end{proof}
\subsection{Properties related to the deformation to the normal cone and normal bundle functors}\label{sec: Properties related to the deformation to the normal cone and normal bundle functors}
Before we show that we can apply the normal bundle functor and deformation to the normal cone functor to pairs of Lie groupoids and pairs of Lie algebroids, and obtain a Lie groupoid and Lie algebroid, respectively, we will go into some properties of deformation to the normal cones. Naturally, it is useful to look at properties of normal bundles as well. Later, this will reduce the work significantly when introducing blow-ups. In this section we fix a pair of smooth manifolds $(\base,\subbase)$.
\begin{rema}\cite{meinrenken}\label{rema: identifying DNC(Y,Y) with Y times R}
Notice that we can identify $\DNC(\base,\base)$ with $\base \times \mathbb{R}$ by identifying $\normal(\base,\base) = T\base|_\base/T\base$ with $\base$. Explicitly, the canonical map
\begin{align*}
f: \DNC(\base,\base) \ra \base \times \mathbb{R} \textnormal{ given by } z \mapsto
\begin{cases}
    (x,0) & \text{if}\ z=(x,0,0) \\
    (x,t) & \text{if}\ z=(x,t)
\end{cases}
\end{align*}
is a diffeomorphism. Indeed, if $(\varphi,U)$ is a chart on $\base$, then of course it is adapted to $Y$. With respect to this chart, the above map locally becomes
\begin{align*}
(\varphi \times \id_{\mathbb{R}}) \circ f \circ (\DNC(\varphi^{-1}) \circ \Psi): \Omega^U_U \ra \varphi(U) \times \mathbb{R} \textnormal{ given by } (x,0,t) \mapsto (x,t),
\end{align*}
so $f$ is, indeed, a diffeomorphism.
\end{rema}
More generally, for every vector bundle $E \ra \base$, we can identify the deformation to the normal cone of $0_\base \subset E$ with $E \times \mathbb{R}$. We can phrase the proof of this statement as a consequence of Corollary \ref{coro: canonical iso N(E,0_M)}. More on this later (in this section).

The following simple map is of great importance to the theory of deformation to the normal cones.
\begin{rema}\cite{meinrenken}\label{rema: canonical submersion}
Consider the map
\[\hat{t} \coloneqq \DNC(\pt): \DNC(\base,\subbase) \ra \DNC(\pt,\pt) \cong \mathbb{R}\]
that is induced by the canonical map $(\base,\subbase) \ra (\pt,\pt)$ (also denoted by $\pt$). Then $\hat{t}$ is a surjective submersion whose fibers are given by $\DNC(\base,\subbase)_t$; the slices at $t$. In particular, this shows that $\DNC(\base,\subbase)_t \subset \DNC(\base,\subbase)$ are closed submanifolds and the inclusion $\base \times \mathbb{R}^\times \hookrightarrow \DNC(\base,\subbase)$ is an open embedding. Moreover, the inclusion $(\subbase,\subbase) \hookrightarrow (\base,\subbase)$ induces a closed embedding 
\[\DNC(\subbase,\subbase) \cong \subbase \times \mathbb{R} \ra \DNC(\base,\subbase),\] 
since its image is the pre-image of $\subbase \times \mathbb{R}$ under the map $\DNC(\id_\base): \DNC(\base,\subbase) \ra \DNC(\base,\base) \cong \base \times \mathbb{R}$, which is a submersion. 
\end{rema}

Recall the notion of clean intersection (see Definition \ref{defn: clean intersetion}). The appropriate notion of clean intersection for maps of pairs is the following.
\begin{defn}\label{defn; clean intersection of pairs}
Let $f_i: (\base_i,\subbase_i) \ra (\basetwo,\subbasetwo)$ (for $i=1,2$) be smooth maps of pairs. The maps $f_1$ and $f_2$ are said to have a \textit{clean intersection} if $f_1: \base_1 \ra \basetwo$ and $f_2: \base_2 \ra \basetwo$ have clean intersection, and $f_1|_{\subbase_1}: \subbase_1 \ra \subbasetwo$ and $f_2|_{\subbase_2}: \subbase_2 \ra \subbasetwo$ have a clean intersection.
\end{defn}
We will see below that the normal bundle and the deformation to the normal cone respect fiber products induced by clean intersections of maps of pairs. These results follow from the following lemma, which is almost immediate from the definition of clean intersection.
\begin{lemm}\label{lemm: if clean intersection, then T respects fiber products}
Let $f_i: \base_i \ra M$ ($i=1,2$) be smooth maps that have clean intersection. Then $T\base_1 \tensor[_{df_1}]{\times}{_{df_2}} T\base_2 \subset T\base_1 \times T\base_2$ is an embedded submanifold, and the canonical map
\begin{align*}
    T(\base_1 \tensor[_{f_1}]{\times}{_{f_2}} \base_2) &\ra T\base_1 \tensor[_{df_1}]{\times}{_{df_2}} T\base_2 \textnormal{ given by } \\
    (x_1,x_2,\xi_1,\xi_2) &\mapsto (x_1,\xi_1,x_2,\xi_2)
\end{align*}
is an isomorphism of vector bundles.
\end{lemm}
\begin{proof}
By assumption, we have for all $(x_1,x_2) \in \base_1 \tensor[_{f_1}]{\times}{_{f_2}} \base_2$ that
\[T_{(x_1,x_2)}(\base_1 \tensor[_{f_1}]{\times}{_{f_2}} \base_2) = d(f_1 \times f_2)(x_1,x_2)^{-1}(T_{(f(x_1),f(x_2))}\Delta_M).\]
Under the canonical isomorphism $T(\base_1 \times \base_2) \xra{\sim} T\base_1 \times T\base_2$, which is explicitly given by $(x_1,x_2,\xi_1,\xi_2) \mapsto (x_1,\xi_1,x_2,\xi_2)$, we see from the former observation that $T(\base_1 \tensor[_{f_1}]{\times}{_{f_2}} \base_2)$ is mapped to
\[\{(x_1,\xi_1,x_2,\xi_2) \in T\base_1 \times T\base_2 \mid df_1(x_1)\xi_1=df_2(x_2)\xi_2\} = T\base_1 \tensor[_{df_1}]{\times}{_{df_2}} T\base_2.\]
Since $T(\base_1 \tensor[_{f_1}]{\times}{_{f_2}} \base_2) \subset T(\base_1 \times \base_2)$ is an embedded submanifold, the result follows.
\end{proof}
\begin{prop}\label{prop: if clean intersection, then N respects fiber products}
Let $f_i: (\base_i,\subbase_i) \ra (\basetwo,\subbasetwo)$ ($i=1,2$) be smooth maps of pairs that have clean intersection. Then $\normal(\base_1,\subbase_1) \tensor[_{d_\normal f_1}]{\times}{_{d_\normal f_2}} \normal(\base_2,\subbase_2) \subset \normal(\base_1,\subbase_1) \times \normal(\base_2,\subbase_2)$ is an embedded submanifold, and the canonical map
\begin{align*}
    \normal(\base_1 \tensor[_{f_1}]{\times}{_{f_2}} \base_2, \subbase_1 \tensor[_{f_1}]{\times}{_{f_2}} \subbase_2) &\ra \normal(\base_1,\subbase_1) \tensor[_{d_\normal f_1}]{\times}{_{d_\normal f_2}} \normal(\base_2,\subbase_2) \textnormal{ given by } \\
    (y_1,y_2,\xi_1,\xi_2) &\mapsto (y_1,\xi_1,y_2,\xi_2)
\end{align*}
is an isomorphism of vector bundles.
\end{prop}
\begin{proof}
By the isomorphism $T(\base_1 \times \base_2)|_{\subbase_1 \times \subbase_2} \xra{\sim} T\base_1|_{\subbase_1} \times T\base_2|_{\subbase_2}$, which restricts to an isomorphism $T(\subbase_1 \times \subbase_2) \xra{\sim} T\subbase_1 \times T\subbase_2$, we obtain a vector bundle isomorphism
\[\normal(\base_1 \times \base_2,\subbase_1 \times \subbase_2) \xra{\sim} \normal(\base_1,\subbase_1) \times \normal(\base_2,\subbase_2) \textnormal{ given by } (y_1,y_2,\xi_1,\xi_2) \mapsto (y_1,\xi_1,y_2,\xi_2).\]
From Lemma \ref{lemm: if clean intersection, then T respects fiber products} (applied to $f_i$ and $f_i|_{\subbase_i}$) we see that this map, restricted to the embedded submanifold $\normal(\base_1 \tensor[_{f_1}]{\times}{_{f_2}} \base_2, \subbase_1 \tensor[_{f_1|_{\subbase_1}}]{\times}{_{f_2|_{\subbase_2}}} \subbase_2)$, maps onto $\normal(\base_1,\subbase_1) \tensor[_{d_\normal f_1}]{\times}{_{d_\normal f_2}} \normal(\base_2,\subbase_2)$, from which the result follows.
\end{proof}
\begin{prop}\label{prop: if clean intersection, then DNC respects fiber products}
Let $f_i: (\base_i,\subbase_i) \ra (\basetwo,\subbasetwo)$ ($i=1,2$) be smooth maps that have clean intersection. Then $\DNC(\base_1,\subbase_1) \tensor[_{\DNC(f_1)}]{\times}{_{\DNC(f_2)}} \DNC(\base_2,\subbase_2) \subset \DNC(\base_1,\subbase_1) \times \DNC(\base_2,\subbase_2)$ is an embedded submanifold, and the canonical map
\begin{align*}
    \DNC(\base_1 \tensor[_{f_1}]{\times}{_{f_2}} \base_2, \subbase_1 \tensor[_{f_1}]{\times}{_{f_2}} \subbase_2) &\ra \DNC(\base_1,\subbase_1) \tensor[_{\DNC(f_1)}]{\times}{_{\DNC(f_2)}} \DNC(\base_2,\subbase_2) \textnormal{ given by } \\ z &\mapsto
    \begin{cases}
        (y_1,\xi_1,y_2,\xi_2) & \text{if}\ z=(y_1,y_2,\xi_1,\xi_2,0) \\
        (y_1,t,y_2,t) & \text{if}\ z=(y_1,y_2,t)
    \end{cases}
\end{align*}
is a diffeomorphism.
\end{prop}
\begin{proof}
Recall that we have a canonical smooth submersion $\hat{t}_i: \DNC(\base_i,\subbase_i) \ra \mathbb{R}$ ($i=1,2$). In particular, $\hat{t}_1$ and $\hat{t}_2$ have clean intersection, so we see that $\DNC(\base_1,\subbase_1) \times_{\mathbb{R}} \DNC(\base_2,\subbase_2) \subset \DNC(\base_1,\subbase_1) \times \DNC(\base_2,\subbase_2)$ is an embedded submanifold. The map $(\DNC(\pr_1),\DNC(\pr_2))$%, which has constant rank,
maps into $\DNC(\base_1,\subbase_1) \times_{\mathbb{R}} \DNC(\base_2,\subbase_2)$, so we obtain a %diffeomorphism
smooth map
\begin{align*}
    \DNC(\base_1 \times \base_2, \subbase_1 \times \subbase_2) \xra{\sim} \DNC(\base_1,\subbase_1) \times_{\mathbb{R}} \DNC(\base_2,\subbase_2) \textnormal{ given by } \\ z \mapsto 
    \begin{cases}
        (y_1,\xi_1,y_2,\xi_2) & \text{if}\ z=(y_1,y_2,\xi_1,\xi_2,0) \\
        (y_1,t,y_2,t) & \text{if}\ z=(y_1,y_2,t).
    \end{cases}
\end{align*}
%by a dimension count. 
The inclusion $\iota_i: \base_i \ra \base_1 \times \base_2$ (for $i=1,2$) gives rise to a smooth map $\DNC(\iota_i): \DNC(\base_i,\subbase_i) \ra \DNC(\base_1 \times \base_2,\subbase_1 \times \subbase_2)$. Then we obtain a well-defined smooth map $(\DNC(\iota_1),\DNC(\iota_2)): \DNC(\base_1,\subbase_1) \times_{\mathbb{R}} \DNC(\base_2,\subbase_2) \ra \DNC(\base_1 \times \base_2,\subbase_1 \times \subbase_2)$ which is inverse to the above map.
Notice that, by definition of $\DNC(f_1)$ and $\DNC(f_2)$, $\DNC(\base_1,\subbase_1) \tensor[_{\DNC(f_1)}]{\times}{_{\DNC(f_2)}} \DNC(\base_2,\subbase_2) \subset \DNC(\base_1,\subbase_1) \times_{\mathbb{R}} \DNC(\base_2,\subbase_2)$, and that $\DNC(\base_1 \tensor[_{f_1}]{\times}{_{f_2}} \base_2, \subbase_1 \tensor[_{f_1|_{\subbase_1}}]{\times}{_{f_2|_{\subbase_2}}} \subbase_2)$ maps onto this space under $(\DNC(\pr_1),\DNC(\pr_2))$. This proves the statement.
\end{proof}

We will now show that deformation to the normal cones of pairs of vector bundles carry a vector bundle structure that extends the vector bundle structures of the ambient vector bundle (for each $t$-slice, where $t\neq0$) and the normal bundle (for the $0$-slice). Then, as a first example, we will see that $\DNC(E,0_\base) \ra \DNC(\base,\base) \cong \base \times \mathbb{R}$ is a vector bundle, and it is naturally isomorphic to the vector bundle $E \times \mathbb{R} \ra \base \times \mathbb{R}$. 

The analogue of Proposition \ref{prop: normal bundle constant rank} for deformation to the normal cones is as follows.
\begin{prop}\label{prop: deformation constant rank}
Let $f: (\base,\subbase) \ra (\basetwo,\subbasetwo)$ be a smooth map of pairs of constant rank $(k,\ell)$. If $\DNC(f)_0=d_\normal f: \normal(\base,\subbase) \ra \normal(\basetwo,\subbasetwo)$ is of constant rank $k-\ell$, then $\DNC(f): \DNC(\base,\subbase) \ra \DNC(\basetwo,\subbasetwo)$ has constant rank $k+1$. In particular, if $f$ is a submersion (as a map of pairs), then $\DNC(f)$ is a submersion. If $f$ is an immersion (as a map of pairs), and $f^{-1}(\subbasetwo)=\subbase$, then $\DNC(f)$ is an immersion.
\end{prop}
\begin{proof}
Notice that, by assumption (and Proposition \ref{prop: normal bundle constant rank}), the maps $\DNC(f)_t: \DNC(\base,\subbase)_t \ra \DNC(\basetwo,\subbasetwo)_t$ are of constant rank $k$ each. Now observe that
\[\hat{t}_{\DNC(\basetwo,\subbasetwo)} \circ \DNC(f) = \hat{t}_{\DNC(\base,\subbase)},\]
simply by definition of $\DNC(f)$, and note that the differential of this map at a point, and restricted to $T_t\mathbb{R}$, has rank $1$. Let $z \in \DNC(\base,\subbase)$. Then $d\DNC(f)(z)$ restricted to $T_t\mathbb{R}$ has rank $1$ by the former observation, and the restriction to $T_x\base$, if $z=(x,t)$, or $T_{(y,\xi)}\normal(\base,\subbase)$, if $z=(y,\xi,0)$, maps into $T_{f(x)}\basetwo$ or $T_{(f(y),d_\normal f(y)\xi)} \normal(\basetwo,\subbasetwo)$, respectively, and has rank $k$ in each case. We can now represent $d\DNC(f)(z)$ as a matrix $A + B$ for which $A$ has rank $k$ and it has last row and column equal to zero, and $B$ has rank $1$ and the only non-zero rows and columns are the last row and column. Therefore, the rank of $\DNC(f)(z)$ has rank $k+1$ by the lemma below. 

%and that $\hat{t}_{\DNC(\base,\subbase)}$ is a submersion, we see that it suffices to prove that $\DNC(f)$ has constant rank $k+1$. Indeed, locally we can write the map as $\widetilde{h}=(\widetilde{h}_1,\widetilde{h}_2,\widetilde{h}_3)$, as in Lemma \ref{lemm: transition maps of DNC are smooth}, with $h=\psi \circ f \circ \varphi^{-1}$ (with $(U,\varphi)$ and $(W,\psi)$ adapted charts on $\base$ and $\basetwo$, respectively), and $\im d(\widetilde{h}_1,\widetilde{h}_2)(y,\xi,t)$ has constant rank $k$
The last statement now follows directly from Proposition \ref{prop: normal bundle constant rank}.
\end{proof}
%The proof of the following lemma is based on \cite{10.2307/2687710}.
\begin{lemm}\cite{10.2307/2687710}\label{lemm: constant rank A+B}
Let $A$ and $B$ be $m \times n$-matrices and denote by $\col(T)$ and $\row(T)$ the span of the columns and rows of a matrix $T$, respectively. Then
\[\rank (A + B) \le \rank A + \rank B\]
and equality holds if and only if $\col(A) \cap \col(B) = 0$ and $\row(A) \cap \row(B) = 0$.
\end{lemm}
\begin{proof}
Notice that $\dim \col(T) = \rank T$ (and $\dim \row(T) = \rank T$) for a matrix $T$. The first statement now follows from the fact that $\col(A+B) \subset \col(A) + \col(B)$. Since $\dim (\col(A) + \col(B)) = \dim \col(A) + \dim \col(B)$ if and only if $\col(A) \cap \col(B)=0$, and similarly for $\row$ instead of $\col$, we see that equality implies $\col(A) \cap \col(B) = 0$ and $\row(A) \cap \row(B) = 0$.

Conversely, suppose that $\col(A) \cap \col(B) = 0$ and $\row(A) \cap \row(B) = 0$. Recall that the rank of a matrix is invariant under row and column operations. We can therefore rewrite $A$ as $PAQ$ (for invertible matrices $P,Q$), because $\rank(PAQ + PBQ) = \rank(P(A+B)Q)=\rank(A+B)$. So, assume 
\[A = \begin{pmatrix} I_r & 0 \\ 0 & 0\end{pmatrix}\]
(note: $\rank A  = r$). Now write $B = ST$, where $S$ is obtained by choosing $s \eqqcolon \rank B$ linear independent columns of $B$ (so $S$ has $s$ columns, and $T$ has $s$ rows). Then $\rank S = \rank T = s$, and it is clear that $\row(T) \subset \row(B)$, so $\rank T \le s$.  Write
\[S = \begin{pmatrix} S_1 \\ S_2\end{pmatrix}; \quad T = \begin{pmatrix} T_1 & T_2 \end{pmatrix},\]
where $S_1$ and $T_1$ are $r \times s$-matrices. Then $S_2$ and $T_2$ both have an invertible $s \times s$-matrix. Indeed, if $S_2v=0$, then $Sv= \begin{pmatrix}S_1v & 0 \end{pmatrix}^T \in \col(A) \cap \col(B) = 0$ by our choice of $A$, so $v=0$, because $S$ has independent columns. Therefore, $S_2$ has rank $s$, and so it has an invertible $s \times s$-matrix, say $S_s$. By considering $T_2^T$ and using that $\row(A) \cap \row(B)=0$, it follows by a similar argument that $T_2$ has an invertible $s \times s$-matrix $T_s$. This yields an $(r+s) \times (r+s)$-submatrix
\[\begin{pmatrix} S_1 \\ S_s\end{pmatrix}\begin{pmatrix} T_1 & T_s \end{pmatrix} = \begin{pmatrix} S_1T_1 & S_1T_s \\ S_sT_1 & S_sT_s \end{pmatrix}\]
of $B$, and so 
\[\begin{pmatrix} I_r+S_1T_1 & S_1T_s \\ S_sT_1 & S_sT_s \end{pmatrix} \xrsquigarrow{$\row_1 - S_1S_s^{-1}\cdot\row_2$} \begin{pmatrix} I_r & 0 \\ S_sT_1 & S_sT_s \end{pmatrix}\]
is an $(r+s) \times (r+s)$-submatrix of $A+B$ which is invertible. This shows that $\rank(A+B) \ge \rank A + \rank B$, so $\rank(A+B)=\rank(A)+\rank(B)$, as required.
\end{proof}
We obtain the following useful corollary
\begin{coro}\label{coro: DNC submanifold}
Let $\iota: (\basetwo,\subbasetwo) \ra (\base,\subbase)$ be an embedding (that is, $\iota: \basetwo \ra \base$ is an embedding, and $\iota|_{\subbasetwo}: \subbasetwo \ra \subbase$ is an embedding) such that $\iota^{-1}(\subbase)=\subbasetwo$. Then 
\[\DNC(\iota): \DNC(\basetwo,\subbasetwo) \ra \DNC(\base,\subbase)\]
is an embedding, which, for all maps of pairs $f: (\base,\subbase) \ra (P,Q)$, maps $\DNC(\iota)^{-1}(\DNC_f(\base,\subbase))$ bijectively onto $\DNC_{f \circ \iota}(\basetwo,\subbasetwo)$. Moreover, if $\iota$ is open (resp. closed), then $\DNC(\iota)$ is open (resp. closed).
\end{coro}
\begin{proof}
We may assume that $\basetwo \subset \base$. Since $\DNC(\iota)$ is an (injective) immersion by Proposition \ref{prop: deformation constant rank}, it suffices to prove, for the first statement, that $\DNC(\iota)$ is a topological embedding, i.e. that the inverse of $f \coloneqq \DNC(\iota): \DNC(\basetwo,\subbasetwo) \ra \DNC(\iota)(\DNC(\basetwo,\subbasetwo))$ is continuous. Let $(W,\psi)$ be an adapted chart of $\basetwo$ to $\subbasetwo$ (with $Z \coloneqq W \cap \subbasetwo$) around $z \in \DNC(\basetwo,\subbasetwo)$.
To show that $f^{-1}$ is continuous, we will show that there is an adapted chart $(U,\varphi)$ of $\base$ adapted to $\subbase$ (with $V \coloneqq U \cap \subbase$) around $f(z)$ such that $\DNC(U,V) \cap \DNC(\iota)(\DNC(\basetwo,\subbasetwo)) \subset f(\DNC(W,Z))$. So, by picking some adapted chart $(U,\varphi)$ (of $\base$ to $\subbase$) around $f(z)$, and possibly shrinking it, we may assume, by also possibly shrinking $W$, that $U \cap \basetwo = W$ and $V \cap \subbasetwo = Z$ (note: $\basetwo \subset \base$ by assumption). Then 
\[f(\DNC(W,Z)) = \DNC(\iota)(\DNC(W,Z)) = \DNC(\iota)(\DNC(U \cap \basetwo, V \cap \subbasetwo)).\]
Therefore, it suffices to prove that 
\begin{equation}\label{eq: open or closed embedding}
    \DNC(\iota)(\DNC(U \cap \basetwo, V \cap \subbasetwo)) = \DNC(U,V) \cap \DNC(\iota)(\DNC(\basetwo,\subbasetwo))
\end{equation} 
in $\DNC(\base,\subbase)$. This, however, follows by checking it on every $t$-slice of $\DNC(\base,\subbase)$ (note: $T\basetwo \subset T\base$ and $T\subbasetwo \subset T\subbase$, so we can naturally view $\normal(\basetwo,\subbasetwo) \subset \normal(\base,\subbase)$, and $\normal(U \cap \basetwo, V \cap \subbasetwo) = \normal(U,V) \cap \normal(\basetwo,\subbasetwo)$, because $T(U \cap \basetwo) = TU \cap T\basetwo$ and similarly $T(V \cap \subbasetwo) = TV \cap T\subbasetwo$).

Since $\DNC(f) \circ \DNC(\iota) = \DNC(f \circ \iota)$, we have
\begin{align*}
    \DNC(\iota)^{-1}(\DNC_f(\base,\subbase)) &= \DNC(\iota)^{-1}(\DNC(\base,\subbase) \setminus \DNC(f)^{-1}(Q \times \mathbb{R})) \\
    &= \DNC(\basetwo,\subbasetwo) \setminus (\DNC(\iota)^{-1}(\DNC(f)^{-1}(Q \times \mathbb{R}))) = \DNC_{f \circ \iota}(\basetwo,\subbasetwo).
\end{align*}
It remains to show that if $\iota$ is open (resp. closed), then $\DNC(\iota)$ is open (resp. closed). However, this follows immediately by (taking complements in) \eqref{eq: open or closed embedding}. This proves the statement.
%Notice for this that we can view all the $t$-slices $\DNC(\basetwo,\subbasetwo)_t \subset \DNC(\base,\subbase)$ as subspaces, since $\DNC(\basetwo,\subbasetwo)_t$ is either diffeomorphic to $\basetwo$ (if $t\neq0$), in which case $\DNC(\basetwo,\subbasetwo)_t \cong \basetwo \subset \base \cong \DNC(\base,\subbase)_t$, and $\DNC(\basetwo,\subbasetwo)_0$ is diffeomorphic to $\normal(\basetwo,\subbasetwo)$, in which case $\DNC(\basetwo,\subbasetwo)_0 \cong \normal(\basetwo,\subbasetwo) \subset \normal(\base,\subbase) \cong \DNC(\base,\subbase)_0$. Since $\sqcup_{t \in \mathbb{R}}\DNC(\basetwo,\subbasetwo)_t$ has a unique topology for which $\basetwo \times \mathbb{R}^\times$ $\DNC(\basetwo,\subbasetwo)_t$ are subspaces 
\end{proof}
Here is the statement about applying the deformation to the normal cone functor to pairs of vector bundles:
\begin{prop}\cite{meinrenken}\label{prop: DNC of subbundle of a vector bundle is vector bundle}
Let $(E \xra{\pi} \base,F \ra \subbase)$ be a pair of vector bundles of rank $(k,\ell)$. Then 
\[\DNC(E,F) \xra{\DNC(\pi)} \DNC(\base,\subbase)\] 
is a vector bundle of rank $k$, and the vector bundle structure is uniquely determined by the property that $E \times \mathbb{R}^\times \xra{\pi \times \id_{\mathbb{R}^\times}} \base \times \mathbb{R}^\times$ is a vector subbundle.
\end{prop}
\begin{proof}
By Proposition \ref{prop: deformation constant rank}, $\DNC(\pi)$ is a submersion, and it is given by
\begin{equation*}
z \mapsto
    \begin{cases}
    (\pi(y),d_\normal\pi(x)\xi,0) & \text{if}\ z = (y,\xi) \\
    (\pi(x),t) & \text{if}\ z = (x,t).
    \end{cases}
\end{equation*}
Every fiber of this map is contained in one of the submanifolds $\DNC(E,F)_t$, i.e. in a fiber of the vector bundle $E \times \mathbb{R}^\times \ra \base \times \mathbb{R}^\times$ or $\normal(E,F) \ra \normal(\base,\subbase)$, so the fibers naturally inherit a vector space structure. That is, the addition and scalar multiplication maps are given by 
\[\DNC(+_E): \DNC(E,F) \tensor[_{\DNC(\pi_E)}]{\times}{_{\DNC(\pi_E)}} \DNC(E,F) \ra \DNC(E,F) \textnormal{ and } \DNC(\cdot_E): \mathbb{R} \times \DNC(E,F) \ra \DNC(E,F)\] 
(see Proposition \ref{prop: if clean intersection, then DNC respects fiber products} to see that the domains of these maps are indeed correct). 

Now, if we have a local trivialisation $\Phi: E|_U \xra{\sim} U \times \mathbb{R}^k$ that restricts to a local trivialisation $F|_V \xra{\sim} V \times \mathbb{R}^q$ (where $V \coloneqq U \cap \subbase$), then we obtain a diffeomorphism
\[\Phi_\DNC: \DNC(E|_U, F|_V) \xra[\sim]{\DNC(\Phi)} \DNC(U \times \mathbb{R}^k, V \times \mathbb{R}^\ell) \xra{\sim} \DNC(U,V) \times \mathbb{R}^k,\]
where we used Proposition \ref{prop: if clean intersection, then DNC respects fiber products} in the latter map (and the identification $\DNC(\mathbb{R}^k, \mathbb{R}^\ell) \cong \mathbb{R}^k \times \mathbb{R}$ which holds by definition; see Section \ref{sec: Deformation to the normal cone functor: objects}).
%\[\DNC(U \times \mathbb{R}^k, V \times \mathbb{R}^q) \xra{\sim} \DNC(U,V) \]
%\[\left(\DNC(\pr_U),\pr_{\mathbb{R}^k} \circ \Psi^{-1} \circ \DNC(\pr_{\mathbb{R}^k})\right), \textnormal{ with } \Psi: \mathcal{D}(\mathbb{R}^k,\mathbb{R}^q) \xra{\sim} \mathbb{R}^q \times \mathbb{R}^{k-q} \times \mathbb{R} \textnormal{ (see Section \ref{sec: Deformation to the normal cone functor: objects}}).\]
%To see that this is a diffeomorphism, observe that it is a submersion (all maps involved are submersions) and
%\[\left(d_\normal\pr_{U},d_\normal\pr_{\mathbb{R}^k}\right): \normal(U \times \mathbb{R}^k, (U \cap X) \times \mathbb{R}^q) \ra \normal(U,U \cap X) \times \mathbb{R}^k\]
%is a submersion and a bijection (hence a diffeomorphism) with inverse given by 
%\[(x,\xi,v) \mapsto (x,\pr_{\mathbb{R}^q}(v),\xi,\pr_{\mathbb{R}^{n-q}}(v))\]
%so that the map
%\begin{equation*}
%(z,v) \mapsto 
%    \begin{cases}
%    (x,\pr_{\mathbb{R}^q}(v),\xi,\pr_{\mathbb{R}^{k-q}}(v),0) & \text{if}\ z=(x,\xi,0) \\
%    (y,\pr_{\mathbb{R}^q}(v),t\pr_{\mathbb{R}^{k-q}}(v),t) & \text{if}\ z=(y,t),
%    \end{cases}
%\end{equation*}
%is inverse to $\left(\DNC(\pr_U),\pr_{\mathbb{R}^k} \circ \Psi^{-1} \circ \DNC(\pr_{\mathbb{R}^k})\right)$. 
Notice that the latter map is explicitly given by $(\DNC(\pr_U),\pr_{\mathbb{R}^k} \circ \Psi^{-1} \circ \DNC(\pr_{\mathbb{R}^k}))$ (see the proof of Proposition \ref{prop: if clean intersection, then DNC respects fiber products}). It is now easy to see that this does indeed give a local trivialisation:
\begin{align*}
    \pr_{\DNC(U,V)} \circ \Phi_\DNC &= \pr_{\DNC(U,V)} \circ \left(\DNC(\pr_U) \circ \DNC(\Phi),\pr_{\mathbb{R}^k} \circ \Psi^{-1} \circ \DNC(\pr_{\mathbb{R}^k}) \circ \DNC(\Phi)\right) \\
    &= \DNC(\pr_U \circ \Phi) = \DNC(\pi|_{E|_U}) = \DNC(\pi)|_{\DNC(E|_U,F|_V)},
\end{align*}
and, fiberwise, $\Phi_\DNC$ is given by 
\begin{align*}
    \DNC(E,F)_{(x,t)} \ra \mathbb{R}^k, \textnormal{ } z \mapsto 
    \begin{cases}
    (\Phi^1(x)f,t) & \text{if}\ z=(f,t) \in F_x \times \{t\} \\
    (\Phi^1(x)e,\tfrac{1}{t}\Phi^2(x)e,t) & \text{if}\ z=(e,t) \in (E_x \setminus F_x) \times \{t\},
    \end{cases}
\end{align*} 
where we wrote $\Phi(x) = (\Phi^1(x),\Phi^2(x)): E_x \ra \mathbb{R}^\ell \times \mathbb{R}^{k-\ell}$, or $d_\normal \Phi_{(y,\xi)}: \normal(E,F)_{(y,\xi)} \ra \mathbb{R}^k$, which are indeed linear isomorphisms. 

The uniqueness statement holds, because the structure maps of $\DNC(E,F)$ (that is, addition, scalar multiplication, and the projection) restrict to the structure maps of $E \times \mathbb{R}^\times$. Since $E \times \mathbb{R}^\times \subset \DNC(E,F)$ is dense, there is at most one extension of the structure maps of $E \times \mathbb{R}^\times$ to $\DNC(E,F)$. This proves the statement.
\end{proof}
We will now describe the sections of the above vector bundle.
\begin{rema}\label{rema: local sections of DNC}
Let $(E \xra{\pi} \base, F \ra \subbase)$ be a pair of vector bundles.
First of all, we claim that, for a local section $\alpha: U \ra E$, the section 
\[t\cdot (\alpha \times 0): U \times \mathbb{R}^\times \ra E \times \mathbb{R}^\times \textnormal{ given by } (u,t') \mapsto (t' \cdot \alpha(u),t'),\]
viewing $E \times \mathbb{R}^\times$ as a subbundle of $\DNC(E,F)$, is smooth, and extends to a smooth local section $\hat{\alpha}: \DNC(U,V) \ra \DNC(E,F)$. Indeed, let $(U,\varphi)$ be an adapted chart on $\base$, and let $(E|_U,\Phi)$ be a vector bundle chart over $\varphi$. Then, in local coordinates, the section is given (on $\Omega^U_V \times \mathbb{R}^\times$; see Section \ref{sec: Deformation to the normal cone functor: objects} for the definition of $\Omega^U_V$) by
\begin{align*}
    %(\Psi^{-1} \circ \DNC(\varphi)) &\times \id_{\mathbb{R}^k} \circ \Phi_\DNC \circ [t \cdot (\alpha \times 0)] \circ \DNC(\varphi^{-1}) \circ \Psi(y,\xi,t) \\
    %&= (\Psi^{-1} \circ \DNC(\varphi)) \times \id_{\mathbb{R}^k} \circ \Phi_\DNC \circ [t \cdot (\alpha \times 0)] \circ \DNC(\varphi^{-1})(y,t\xi,t) \\
    %&=(\Psi^{-1} \circ \DNC(\varphi)) \times \id_{\mathbb{R}^k} \circ \Phi_\DNC \circ [t \cdot (\alpha \times 0)](\varphi^{-1}(y,t\xi),t)\\
    %&=(\Psi^{-1} \circ \DNC(\varphi)) \times \id_{\mathbb{R}^k} \circ (\varphi^{-1}(y,t\xi),t, \Phi_\DNC(\varphi^{-1}(y,t\xi))t \cdot \alpha(\varphi^{-1}(y,t\xi))) \\
    %&= 
    (y,\xi,t) \mapsto (y,\xi,t, t \cdot \Phi_\DNC(x)\alpha(x)) = \begin{cases}
    (y,\xi,t, t \cdot \Phi^1(x)\alpha(x)) & \text{if}\ \alpha(x) \in F_x \\
    (y,\xi,t, t \cdot \Phi^1(x)\alpha(x),\Phi^2(x)\alpha(x)) & \text{if}\ \alpha(x) \in E_x \setminus F_x
    \end{cases}
\end{align*}
where we wrote $x \coloneqq \varphi^{-1}(y,t\xi)$. This map does indeed extend to $\DNC(U,V)_0$, and we see that it is the map $\alpha|_V: V \ra E|_V \ra E|_V/F$ on $\normal(U,V)$ (it is a core section of $\normal(E|_U,F|_V) \ra \normal(U,V)$). Similarly, we see that the section $\alpha \times 0: U \times \mathbb{R}^\times \ra E \times \mathbb{R}^\times$ extends to $\DNC(U,V)_0$ if and only if $\alpha|_V \in \Gamma(F)$, and then it is just the map $\DNC(\alpha)$. 

Now, if $\{s^j: U \ra E\}$ is a local frame of $E$ that restricts to a local frame $\{s^j|_V: V \ra F \mid 1 \le j \le \ell\}$ of $F$, then $\{\hat{s^i},\DNC(s^j) \mid \ell+1 \le i \le k, 1 \le j \le \ell\}$ is a local frame of $\DNC(E,F)$. It is clear that these sections are linearly independent over $E \times \mathbb{R}^\times \subset \DNC(E,F)$. Over $\normal(E,F) \subset \DNC(E,F)$, the collection becomes $\{s^i|_V,d_\normal s^j \mid \ell+1 \le i \le k, 1 \le j \le \ell\}$, and so the claim follows by the proof of Proposition \ref{prop: normal bundle of vector bundle}.
\end{rema}
Notice that from this remark we can, in particular, describe the algebra of smooth functions on $\DNC(\base,\subbase)$.
\begin{rema}\cite{Bischoff_2020}\label{rema: smooth functions DNC}
Notice that we came across three types of smooth functions on $\DNC(\base,\subbase)$. If $f \in C^\infty(\base)$ vanishes along $\subbase$, then we obtain a map $\DNC(f): \DNC(\base,\subbase) \ra \DNC(\mathbb{R},0)$. If we look at the definition of this map, we see that we can split $\DNC(f)$ up into two smooth functions of $\DNC(\base,\subbase)$, namely $\hat{\pr}_1 \circ \DNC(f)$ and $\hat{t} \circ \DNC(f)=\hat{t}$, where 
\[\hat{\pr}_1 \coloneqq \DNC(\mathbb{R},0) \cong \mathbb{R}^2 \xra{\pr_1} \mathbb{R}; \quad \hat{t} = \hat{\pr}_2 \coloneqq \DNC(\mathbb{R},0) \cong \mathbb{R}^2 \xra{\pr_2} \mathbb{R}\]
(notice that $\hat{\pr}_2$ is indeed equal to the map $\hat{t} = \DNC(\pt): \DNC(\base,\subbase) \ra \DNC(\pt,\pt) \cong \mathbb{R}$ as we defined it before). 
By abuse of notation, we will often write $\DNC(f)$ for the smooth map $\hat{\pr}_1 \circ \DNC(f): \DNC(\base,\subbase) \ra \mathbb{R}$. If $f \in C^\infty(\base)$, then by Remark \ref{rema: local sections of DNC} (applied to $E=\base \times \mathbb{R}$, $F = \subbase \times \mathbb{R}$), we also obtain smooth maps $\hat{f}: \DNC(\base,\subbase) \ra \mathbb{R}$. Alternatively,
\[\hat{f} = \DNC(\base,\subbase) \xra{\DNC(\id_\base)} \DNC(\base,\base) \cong \base \times \mathbb{R} \xra{\pr_1} \base \xra{f} \mathbb{R}.\]
Summarising, we have the following three types of smooth functions on $\DNC(\base,\subbase)$:
\begin{align*}
    \hat{f}_0(z) &\coloneqq \begin{cases} f_0(y) & \text{if}\ z=(y,\xi) \\ f_0(x) & \text{if}\ z=(x,t)\end{cases}; \\
    \DNC(f_1)(z) &\coloneqq \begin{cases}d_\normal f_1(y)\xi & \text{if}\ z=(y,\xi) \\ \tfrac{1}{t}f_1(x) & \text{if}\ z=(x,t) \end{cases}; \\
    \hat{t}(z) &\coloneqq \begin{cases} 0 & \text{if}\ z=(y,\xi) \\ t & \text{if}\ z=(x,t) \end{cases},
\end{align*}
where $f_0,f_1 \in C^\infty(\base)$ with $f_1$ vanishing along $\subbase$.
\end{rema}
\begin{prop}\cite{meinrenken}\label{prop: identify DNC(E,0_Y) with E times R}
Let $E \xra{\pi} \base$ be a vector bundle. Then we have a canonical isomorphism of vector bundles $\DNC(E,0_\base) \cong E \times \mathbb{R}$.
\end{prop}
\begin{proof}
By Proposition \ref{coro: canonical iso N(E,0_M)}, we can canonically identify $\normal(E,0_\base)$ with $E$. Now, the map $\id_E: (E,0_\base) \ra (E,E)$ induces a smooth map $\DNC(\id_E): \DNC(E,0_\base) \ra \DNC(E,E) \cong E \times \mathbb{R}$. However, notice that the restriction to $\normal(E,0_\base)$ is the canonical identification of $\normal(E,0_Y)$ with $E$, and since this map is an isomorphism of vector bundles (see Corollary \ref{coro: canonical iso N(E,0_M)}), $\DNC(\id_E)$ is an isomorphism of vector bundles, as required.
\end{proof}
\begin{rema}\label{rema: normal of normal in DNC}
Similar to the result above, we have $\normal(\DNC(\base,\subbase),\normal(\base,\subbase)) \cong \normal(\base,\subbase) \times \mathbb{R}$. Indeed, %if $(U,y,x)$ is an adapted chart of $\base$ (with $V \coloneqq U \cap \subbase$), then $\DNC(\base,\subbase)$ has local coordinates $(y,\hat{x},t)$. Since $\normal(\base,\subbase)$ in these coordinates corresponds to $t=0$, we see that  $(y,dx,dt)$ are local coordinates for $\normal(\DNC(\base,\subbase),\normal(\base,\subbase))$. From this we can see that 
$d_\normal\hat{t}: \normal(\DNC(\base,\subbase),\normal(\base,\subbase)) \ra \normal(\mathbb{R},0) \cong \mathbb{R}$ is fiberwise given by isomorphisms, so
\[\normal(\DNC(\base,\subbase),\normal(\base,\subbase)) \ra \normal(\base,\subbase) \times \mathbb{R} \textnormal{ given by } (y,\xi,\eta) \mapsto (y,\xi,d_\normal\hat{t}(y,\xi)\eta)\]
is an isomorphism.
\end{rema}

In fact, the functoriality of the deformation to the normal cone construction (and of the normal bundle construction) extend to functorial properties in the vector bundle setting.
\begin{prop}\label{prop: DNC morphism of vector bundles}
Let $\varphi: (E,F) \ra (\algebrtwo,\subalgebrtwo)$ be a morphism of pairs of vector bundles (that is, $\varphi: E \ra \algebrtwo$ is a morphism of vector bundles, and so is $\varphi|_F: F \ra \subalgebrtwo$). Then $\DNC(\varphi): \DNC(E,F) \ra \DNC(\algebrtwo,\subalgebrtwo)$.
\end{prop}
\begin{proof}
By the vector bundle structures of $\DNC(E,F)$ and $\DNC(\algebrtwo,\subalgebrtwo)$, and by definition of $\DNC(\varphi)$, it suffices to prove that $d_\normal \varphi: \normal(E,F) \ra \normal(\algebrtwo,\subalgebrtwo)$ is a morphism of vector bundles. But this follows by the observation that $d\varphi: TE \ra T\algebrtwo$ is a morphism of vector bundles restricting to a morphism of vector bundles $TF \ra T\subalgebrtwo$.
\end{proof}

Consider the smooth submersion $\hat{t}: \DNC(\base,\subbase) \ra \mathbb{R}$ (see Remark \ref{rema: canonical submersion}). We claim that the vector bundle $\DNC(T\base,T\subbase)$ can be identified with the kernel of $d\hat{t}$ (see Proposition \ref{prop: DNC of tangent spaces}). To prove this, we will construct an embedding $\DNC(T\base,T\subbase) \ra T\DNC(\base,\subbase)$ and then show that it maps onto $d\hat{t}$. However, to make sense of this, we will begin by showing that there is a canonical isomorphism between $\normal(T\base,T\subbase)$ and $T\normal(\base,\subbase)$ as vector bundles over $T\subbase$, for which we will need the following lemma.
\begin{lemm}\cite{KMS}\label{lemm: canonical flip}
Let $\base$ be a smooth manifold of dimension $n$. Then there is a canonical vector bundle isomorphism
\begin{center}
\begin{tikzcd}
TT\base \ar[r,"\kappa_\base","\sim"'] \ar{d}{\pi_{TT\base}} & TT\base \ar{d}{d\pi_{T\base}} \\
T\base \ar{r}{=} & T\base
\end{tikzcd}
\end{center}
which is called the \textit{canonical flip}.
\end{lemm}
\begin{proof}
Recall from Example \ref{exam: tangent DVB} that a vector bundle $\pi_E: E \ra \base$ (of rank $k$) induces the tangent prolongation DVB 
\begin{center}
\begin{tikzcd}
TE \ar[r,"\pi_{TE}"] \ar{d}{d\pi_E} & E \ar{d}{\pi_{E}} \\
T\base \ar{r}{\pi_{T\base}} & \base
\end{tikzcd}
\end{center}
In local coordinates: let $\{(U_\alpha,\varphi_\alpha)\}_{\alpha \in I}$ be an atlas of $\base$ and let $\{(E|_{U_\alpha},\psi_\alpha')\}_{\alpha \in I}$ be a cover of local trivialisations of $E$ such that 
\[\{(E_\alpha,\psi_\alpha)\}_{\alpha \in I} \coloneqq \{(E|_{U_\alpha},(\varphi_\alpha \times \id_V) \circ \psi_\alpha')\}_{\alpha \in I}\]
is an atlas of $E$; we denoted $V \coloneqq \mathbb{R}^k$ for notational purposes. Then,
\[\{(T(E|_{U_\alpha}),d\psi_\alpha)\}_{\alpha \in I}\]
constitutes an atlas for $TE$ with transition maps
\begin{align*}
    d\psi_\alpha \circ d\psi_\beta^{-1}&: \varphi_{\beta}(U_{\alpha\beta}) \times V \times \mathbb{R}^n \times V \ra \varphi_{\alpha}(U_{\alpha\beta}) \times V \times \mathbb{R}^n \times V \textnormal{ given by } \\
    (x,v,\xi,w) &\mapsto \left(\varphi_{\alpha\beta}(x),
    \psi_{\alpha\beta}'(\varphi_{\beta}^{-1}(x))v,
    d\varphi_{\alpha\beta}(x)\xi,
    [d(\psi_{\alpha\beta}' \circ \varphi_\beta^{-1})(x)\xi]v + \psi_{\alpha\beta}'(\varphi_\beta^{-1}(x))w\right),
\end{align*}
where $U_{\alpha\beta} = U_{\alpha} \cap U_{\beta}$, $\varphi_{\alpha\beta}=\varphi_\alpha \circ \varphi_\beta^{-1}$ and $\psi_{\alpha}' \circ (\psi_{\beta}')^{-1}=(\varphi_{\alpha\beta},\psi_{\alpha\beta}')$ with $\psi_{\alpha\beta}': U_{\alpha\beta} \ra \textnormal{GL}(V)$.

Now, put $E=T\base$. Then from this discussion we obtain a double vector bundle
\begin{center}
\begin{tikzcd}
    TT\base \ar{r}{d\pi_\base} \ar{d}{\pi_{T\base}} & T\base \ar{d}{\pi_\base} \\
    T\base \ar{r}{\pi_\base} & \base.
\end{tikzcd}
\end{center}
By inspecting the above transition maps (note: $\psi_\alpha'=d\varphi_\alpha$ using the notation from above), we see that the map
\[\kappa_\base: TT\base \ra TT\base, \textnormal{ given by } \left.\frac{d}{dt}\right\vert_{t=0}\left.\frac{d}{ds}\right\vert_{s=0}\gamma(t,s) \mapsto \left.\frac{d}{ds}\right\vert_{s=0}\left.\frac{d}{dt}\right\vert_{t=0}\gamma(t,s),\] 
locally, with respect to a chart $(U_\alpha,\varphi_\alpha)$ on $\base$, is described by the formula
\[dd\varphi_\alpha \circ \kappa_\base \circ dd\varphi_\alpha^{-1} : TT\varphi_\alpha(U_\alpha) \ra TT\varphi_\alpha(U_\alpha) \textnormal{ given by } (x,\xi,\eta,\zeta) \ra (x,\eta,\xi,\zeta).\]
This proves the statement.
\end{proof}
\begin{rema}\cite{KMS}\label{rema: canonical flip is involution}
Observe that $\kappa_\base: TT\base \ra TT\base$ is inverse to itself, i.e. an involution.
\end{rema}
\begin{lemm}\cite{meinrenken}\label{lemm: N of tangent spaces}
Let $(\base,\subbase)$ be a pair of smooth manifolds. Then $\normal(T\base,T\subbase) \cong T\normal(\base,\subbase)$ as vector bundles over $T\subbase$ by a canonical diffeomorphism.
\end{lemm}
\begin{proof}
Consider the tangent bundle $T\base \xra{\pi_{T\base}} \base$ restricted to $\subbase$:
\[\pi^\subbase_{T\base}: T\base|_\subbase \ra \subbase.\]
Then $d\pi^\subbase_{T\base}: T(T\base|_\subbase) \ra T\base$ is a vector bundle. Also, the vector bundle $TT\base \xra{\pi_{TT\base}} T\base$ induces the vector bundle $\pi^{T\subbase}_{TT\base}: T(T\base)|_{T\subbase} \ra T\subbase$. From the definition of the canonical flip $\kappa_\base$ (see Lemma \ref{lemm: canonical flip}), we see that it restricts to an isomorphism between the vector bundle bundles $T(T\base)|_{T\subbase} \xra{\pi^{T\subbase}_{TT\base}} T\subbase$ and $T(T\base|_\subbase) \xra{d\pi^{\subbase}_{T\base}} T\subbase$, so they are canonically isomorphic; we denote this bundle isomorphism $T(T\base)|_{T\subbase} \ra T(T\base|_\subbase)$ by $\kappa^\subbase_\base$. Now, the differential of the bundle morphism $\tau_{\normal(Y,X)}: T\base|_\subbase \ra \normal(\base,\subbase)$ (over $\subbase$) induces a morphism of vector bundles $T(T\base)|_{T\subbase} \ra T\normal(\base,\subbase)$ (over $T\subbase$), by factoring through $\kappa^\subbase_\base$. This morphism has kernel $TT\subbase$, hence we obtain the desired (canonical) isomorphism 
\begin{center}
\begin{tikzcd}
d\tau_{\normal(\base,\subbase)} \circ \kappa^\subbase_\base: T(T\base)|_{T\subbase}/TT\subbase = \normal(T\base,T\subbase) \ar[r,"\sim"] & T\normal(\base,\subbase).
\end{tikzcd}
\end{center} 
This concludes the proof.
\end{proof}
\begin{prop}\cite{meinrenken}\label{prop: DNC of tangent spaces}
Consider the smooth submersion $\hat{t}: \DNC(\base,\subbase) \ra \mathbb{R}$ (see Remark \ref{rema: canonical submersion}). Then $\ker d\hat{t} \cong \DNC(T\base,T\subbase)$ as vector bundles by a canonical diffeomorphism.
\end{prop}
\begin{proof}
%The vector bundle inclusion 
%\begin{center}
%\begin{tikzcd}
%TY|_X \ar[d,"\pi^Y_X"] \ar[r,hook,"\iota^Y"] & TY \ar[d,"\pi^Y"] \\
%X \ar[r,hook] & Y
%\end{tikzcd}
%\end{center}
%induces the vector bundle inclusion
%\begin{center}
%\begin{tikzcd}
%T(TY|_X) \ar[d,"d\pi^Y_X"] \ar[r,hook,"d\iota^Y"] & T(TY) \ar[d,"d\pi^Y"] \\
%TX \ar[r,hook] & TY.
%\end{tikzcd}
%\end{center}
%Clearly, $d\pi^Y=\pi^{TY}$ and since, as manifolds, $T(TY|_X)$ embeds into $T(TY)$ via $d\iota^Y$ and it maps into the embedded submanifold $T(TY)|_{TX}$ of $T(TY)$, we obtain a vector bundle morphism
%\begin{center}
%\begin{tikzcd}
%T(TY|_X) \ar[d,"d\pi^Y_X"] \ar[r,hook,"d\iota^Y"] & T(TY)|_{TX} \ar[d,"\pi^{TY}_{TX}"] \\
%TX \ar[r,"="] & TX,
%\end{tikzcd}
%\end{center}
%which is an isomorphism, because the vector bundles have the same rank (namely $2n$, where $n\coloneqq\dim Y$). 
Notice that the map
\begin{equation*}
\DNC(T\base,T\subbase) \ra T\DNC(\base,\subbase) \textnormal{ given by } z \mapsto
    \begin{cases}
    d\tau_{\normal(\base,\subbase)}\circ \kappa^\subbase_\base(z) & \text{if}\ z \in \normal(T\base,T\subbase) \\
    (x,t,v,0) & \text{if}\ z=(x,v,t),
    \end{cases}
\end{equation*}
(we used the same notation as in Lemma \ref{lemm: N of tangent spaces}) is injective and maps $\DNC(T\base,T\subbase)_t$ onto $T\DNC(\base,\subbase)_t$. 
%Since $\hat{t}$ is a submersion with fibers $\DNC(\base,\subbase)_t$, the map maps bijectively onto $\ker d\hat{t}$. 
It is even a morphism of vector bundles over the identity: indeed, the restriction to $T\base \times \mathbb{R}^\times$ is the map 
\[T\base \times \mathbb{R}^\times \xra{\id_{T\base} \times 0_{\mathbb{R}^\times}} T\base \times T\mathbb{R}^\times \cong T(\base \times \mathbb{R}^\times) = T\DNC(\base,\subbase)|_{\base \times \mathbb{R}^\times},\] 
and the restriction to $\normal(T\base,T\subbase)$ is the isomorphism of vector bundles $\normal(T\base,T\subbase) \xra{\sim} T\normal(\base,\subbase)$. The map $\DNC(T\base,T\subbase) \ra T\DNC(\base,\subbase)$ clearly maps into $\ker d\hat{t}$, so it maps bijectively onto $\ker d\hat{t}$ by a dimension count. As manifolds, $\ker d\hat{t}$ is an embedded submanifold of $T\DNC(\base,\subbase)$ (because $d\hat{t}$ is a submersion), so this proves the statement.
\end{proof}
We can wonder now in what way the $\DNC$-functor commutes with the differential $d$. Let $f: (\base,\subbase) \ra (\basetwo,\subbasetwo)$ be a map of pairs. Then $\DNC(df)$ is a map $\DNC(T\base,T\subbase) \ra \DNC(T\basetwo,T\subbasetwo)$, and $d\DNC(f)$ is a map $T\DNC(\base,\subbase) \ra T\DNC(\basetwo,\subbasetwo)$, so the only thing we can expect to happen is that we have an equality of maps $d\DNC(f)|_{\DNC(T\base,T\subbase)}=\DNC(df)$. 
This is indeed the case:
\begin{prop}\label{prop: DNC of differentials}
Let $f: (\base,\subbase) \ra (\basetwo,\subbasetwo)$ be a map of pairs. Then 
\[d\DNC(f)|_{\DNC(T\base,T\subbase)} = \DNC(df).\]\vspace{-\baselineskip}
\end{prop}
\begin{proof}
To see this, note that, by definition of $\DNC(f)$, we have $\DNC(f)(\DNC(\base,\subbase)_t) \subset \DNC(\basetwo,\subbasetwo)_t$. Since, for all $t \in \mathbb{R}$, we have $T\DNC(\base,\subbase)_t = \DNC(T\base,T\subbase)_t$ (and $T\DNC(\basetwo,\subbasetwo)_t = \DNC(T\basetwo,T\subbasetwo)_t$; see Proposition \ref{prop: DNC of tangent spaces}), $d\DNC(f)|_{\DNC(T\base,T\subbase)}$ is a map $\DNC(T\base,T\subbase) \ra \DNC(T\basetwo,T\subbasetwo)$.
In particular, since
\begin{align*}
d\DNC(f)|_{\DNC(T\base,T\subbase)}: \DNC(T\base,T\subbase) &\ra \DNC(T\basetwo,T\subbasetwo) \textnormal{ is given by } \\
z &\mapsto
\begin{cases}
   (f(y),df(y)v,d(d_\normal f)(f(y),df(y)v)\xi,0) & \text{if}\ z=(y,v,\xi) \\
   (f(x),df(x)v,t) & \text{if}\ z=(x,v,t),
\end{cases}
\end{align*}
we see that $d\DNC(f)|_{\DNC(T\base,T\subbase)}$ and $\DNC(df)$ agree on the dense open subset $T\base \times \mathbb{R}^\times \subset \DNC(T\base,T\subbase)$. Therefore, there is a unique extension of this map to all of $\DNC(T\base,T\subbase)$ (if it exists), and so the maps have to be equal.
\end{proof}
\begin{rema}\label{rema: DNC of differentials}
%The proof of the above statement only holds in case $\basetwo$ is Hausdorff. 
As a consequence of the above proposition, we see that $d(d_\normal f) = d_\normal(df)$ under the identification of $T\normal(\base,\subbase)$ with $\normal(T\base,T\subbase)$. Alternatively, we could have proven this identity directly, and the statement would follow from the expression for $d\DNC(f)|_{\DNC(T\base,T\subbase)}$ (see the proof above). To do this, consider the diagram %This yields a proof also for non-Hausdorff $\basetwo$; to do this, consider the diagram
\begin{center}
\begin{tikzcd}
T(T\base|_\subbase) \ar[rrr, bend left, "d(df|_\subbase)"] \ar[ddr,twoheadrightarrow, bend right, "d\tau_{\normal(\base,\subbase)}"'] \ar[r,"\kappa^\subbase_\base","\sim"'] & T(T\base)|_{T\subbase} \ar[d,twoheadrightarrow,"\tau_{\normal(T\base,T\subbase)}"] \ar[r,"d(df)|_{T\subbase}"] & T(T\basetwo)|_{T\subbasetwo} \ar[d,twoheadrightarrow,"\tau_{\normal(T\basetwo,T\subbasetwo)}"] & T(T\basetwo|_{\subbase}) \ar[l,"\kappa^{\subbasetwo}_{\basetwo}"', "\sim"] \ar[ddl,twoheadrightarrow, bend left, "d\tau_{\normal(\basetwo,\subbasetwo)}"] \\
& \normal(T\base,T\subbase) \ar[d,"\sim"] \ar[r,"d_\normal(df)"] & \normal(T\basetwo,T\subbasetwo) \ar[d,"\sim"]& \\
& T\normal(\base,\subbase) \ar[r, "d(d_\normal f)"] & T\normal(\basetwo,\subbasetwo), &
\end{tikzcd}
\end{center}
so since the outer pentagon commutes:
\[d\tau_{\normal(\basetwo,\subbasetwo)} \circ d(df|_\subbase) = d(\tau_{\normal(\basetwo,\subbasetwo)} \circ df|_\subbase) = d(d_\normal f \circ \tau_{\normal(\base,\subbase)}) = d(d_\normal f) \circ d\tau_{\normal(\base,\basetwo)}\]
we see that the outer square in the middle commutes and hence also the lower square.
\end{rema}
%Recall that the the total space of the Lie algebroid of a Lie groupoid $\groupoid$ can be viewed as the normal bundle $\normal(\group,\base)$ (see Remark \ref{rema: other definitions of Lie algebroid of a Lie groupoid}). 
Since the blow-up construction can be viewed as a special case of a principal $\mathbb{R}^\times$-bundle, the following corollary is very helpful (e.g. recall Remark \ref{rema: other definitions of Lie algebroid of a Lie groupoid} that we can view the Lie algebroid of a Lie groupoid $\groupoid$ as the normal bundle $\normal(\group,\base) \ra \base$).

\begin{prop}\label{prop: induced G action on normal bundle}
Let $P \xra{\pi} P/G$ be a principal $G$-bundle. Suppose $(P,Q)$ is a pair of manifolds such that $Q$ is $G$-invariant. Then $\normal(P,Q) \ra \normal(P/G,Q/G)$ is naturally a principal $G$-bundle.
\end{prop}
\begin{proof}
The action map $\rho: G \times P \ra P$ is a map of pairs
\[\rho: (G \times P, G \times Q) \ra (P,Q)\]
by assumption. So, it descends to a map
\[G \times \normal(P,Q) \cong \normal(G \times P, G \times Q) \xra{d_\normal\rho} \normal(P,Q)\]
where we used Proposition \ref{prop: if clean intersection, then N respects fiber products} with $f$ and $g$ the canonical maps $(P,Q) \ra (\pt,\pt)$ and $(G,G) \ra (\pt,\pt)$. To see that it is a principal $G$-bundle, observe that the submersion
\[d_\normal\pi: \normal(P,Q) \ra \normal(P/G,Q/G)\]
is $G$-equivariant (with trivial $G$-action on $\normal(P/G,Q/G)$). Indeed,
\[d_\normal\pi(g \cdot q,dL_g(q)\xi)= d_\normal\pi \circ d_\normal\rho(g,q,\xi) = d_\normal(\pi \circ \rho)(g,q,\xi) = d_\normal\pi(q,\xi)\]
where in the last equality we used that $\pi$ is $G$-equivariant. Moreover, if $U \subset P/G$ is open and
\[\psi_U: \pi^{-1}(U) \ra G \times U\]
is a local trivialisation, then it is a map of pairs $(\pi^{-1}(U),Q \cap \pi^{-1}(U)) \ra (G \times U, G \times (U \cap Q/G))$. Hence, we obtain a diffeomorphism
\[d_\normal\psi_U: \normal(\pi^{-1}(U),Q \cap \pi^{-1}(U)) \xra{\sim} \normal(G \times U, G \times (U \cap Q/G)) \cong G \times \normal(U,U \cap Q/G).\]
Again by the chain rule, this is a local trivialisation. Since the open sets $\normal(U,U \cap Q/G)$ cover $\normal(P/G,Q/G)$, we conclude that, indeed, $\normal(P,Q) \ra \normal(P/G,Q/G) \cong \normal(P,Q)/G$ is a principal $G$-bundle.
\end{proof}

\subsection{Deformation to the normal cone groupoids and algebroids}\label{sec: Deformation to the normal cone groupoid}
In this section we will see that the deformation to the normal cone construction can be extended to constructions on Lie groupoids and Lie algebroids. We start with groupoids. 
\begin{prop}\cite{meinrenken}\label{prop: DNC groupoid defn}
Let $(\groupoid,\subgroupoid)$ be a pair of Lie groupoids (see Definition \ref{defn: pair of Lie subgroupoids}). Then there is a unique Lie groupoid structure on \[\DNCgroupoid\] for which $\group \times \mathbb{R}^\times \ra \base \times \mathbb{R}^\times$ is a Lie subgroupoid. Moreover, $\DNC(\group,\subgroup)^{(2)}$ can canonically be identified with $\DNC(\group^{(2)},\subgroup^{(2)})$, and the structure maps are obtained by applying the functor $\DNC$ to the structure maps of $\group$. 
\end{prop}
\begin{proof}
That $\DNC(\group,\subgroup)^{(2)}$ can be canonically identified $\DNC(\group^{(2)},\subgroup^{(2)})$ follows immediately from Proposition \ref{prop: if clean intersection, then DNC respects fiber products}. It follows that $\DNC(\mult)$ can be seen as a smooth map $\DNC(\group,\subgroup)^{(2)} \ra \DNC(\group,\subgroup)$. Checking that it defines a groupoid is now easy by using the functorial properties of the deformation to the normal cone: let $g \in \DNC(\group,\subgroup)$. Then 
\[\DNC(\target)(\DNC(\inv)(g)) = \DNC(\target \circ \inv)(g) = \DNC(\source)(g),\]
so $(g,\DNC(\inv)(g)) \in \DNC(\group,\subgroup)^{(2)}$ and, moreover,
\begin{align*}
    \DNC(\mult)(g,\DNC(\inv)(g)) &= \DNC(\mult) \circ \DNC(\id_\group \times \inv)(g,g) \\
    &= \DNC(\mult \circ (\id_\group \times \inv))(g,g) \\
    &= \DNC(\identity \circ \target \circ \pr_1)(g,g) = \DNC(\identity)(\DNC(\target)(g)).
\end{align*}\
It is a Lie groupoid because all the structure maps are smooth, and the source map $\DNC(\source)$ and the target map $\DNC(\target)$ are submersions by Proposition \ref{prop: deformation constant rank}. To see that it is the unique groupoid structure for which $\group \times \mathbb{R}^\times$ is a subgroupoid, notice that all of the structure maps restrict to the structure maps on this subgroupoid, and because $\group \times \mathbb{R}^\times \subset \DNC(\group,\subgroup)$ is dense, one can show that there is at most one extension of the structure maps of $\group \times \mathbb{R}^\times$ to $\DNC(\group,\subgroup)$.  This proves the statement.
\end{proof}
\begin{rema}\label{rema: normal bundle groupoid}
In particular, if $(\groupoid,\subgroupoid)$ is a pair of Lie groupoids, then 
\[\normal(\group,\subgroup) \xrra[d_\normal\target]{d_\normal\source} \normal(\base,\subbase)\]
inherits the structure of a Lie groupoid from $\DNC(\group,\subgroup)$. Alternatively, it follows by the same proof as above. The structure maps are given by applying the normal derivatives to all structure maps of $\group$, and $\normal(\group^{(2)},\subgroup^{(2)})$ can be identified with $\normal(\group,\subgroup)^{(2)}$ by Proposition \ref{prop: if clean intersection, then N respects fiber products}. In particular, we see that
\begin{center}
    \begin{tikzcd}
    \normal(\group,\subgroup) \ar[r,shift left] \ar[r,shift right] \ar[d] & \normal(\base,\subbase) \ar[d] \\
    \subgroup \ar[r, shift left] \ar[r, shift right] & \subbase
    \end{tikzcd}
\end{center}
is a $\VB$-groupoid.
\end{rema}
We now state the Lie algebroid analogue of the former statement.
\begin{prop}\cite{meinrenken}\label{prop: DNC of algebroid}
Let $(\algebr \xra{\pi} \base,\subalgebroid)$ be a pair of Lie algebroids. Then there is a unique Lie algebroid structure on 
\[\DNC(\algebr,\subalgebr) \xra{\DNC(\pi)} \DNC(\base,\subbase)\]
such that $\algebr \times \mathbb{R}^\times \ra \base \times \mathbb{R}^\times$ is a Lie subalgebroid.
\end{prop}
\begin{proof}
To see this, recall that we can embed $\DNC(T\base,T\subbase)$ into $T\DNC(\base,\subbase)$ as the kernel of the differential of the canonical map $\hat{t}: \DNC(\base,\subbase) \ra \mathbb{R}$. Therefore, we can take as anchor map the map
\[\anchor_{\DNC(\algebr,\subalgebr)} \coloneqq \DNC(\anchor_\algebr): \DNC(\algebr,\subalgebr) \ra \DNC(T\base,T\subbase) \hookrightarrow T\DNC(\base,\subbase).\]
Now, $\DNC(\pi_\algebr): \DNC(\algebr,\subalgebr) \ra \DNC(\base,\subbase)$ is an anchored vector bundle (see Example \ref{prop: DNC of subbundle of a vector bundle is vector bundle} to see that it is a vector bundle). 
Before we introduce a Lie bracket $[\cdot,\cdot]_{\DNC(\algebr,\subalgebr)}$, we remark that the sections $\Gamma(\DNC(A,B))$ are generated by $\Gamma(A)$ (which give rise to the sections $\hat{\alpha}$) and $\Gamma(A,B)$ (which give rise to the sections $\DNC(\alpha)$) as a $C^\infty(\DNC(\base,\subbase))$-module (see Remark \ref{rema: local sections of DNC}). 
%We will explain this, in detail, later. In short, a section $\sigma \in \Gamma(A)$ determines the section 
%\[\sigma \times 0 \in \Gamma(A \times \mathbb{R}^\times).\]
%By viewing $Y \times \mathbb{R}^\times \subset \DNC(Y,X)$ and $A \times \mathbb{R}^\times \subset \DNC(A,B)$, and writing $\sigma \times 0$ in local coordinates there, one sees that $\sigma \times 0$ extends to a section on all of $\DNC(A,B)$ if and only if $\sigma|_X \in \Gamma(B)$ (the extension is then equal to $\DNC(\sigma)$) and $t(\sigma \times 0)$ extends always to a section $\widetilde{\sigma}$ (the restriction to $\normal(Y,X)$ is then equal to $\sigma|_X$ viewed as a section of $A|_X/B$ inside $\normal(A,B)$). 
Now, in order to define a Lie bracket on $\DNC(A,B) \ra \DNC(\base,\subbase)$ it suffices to specify what the bracket should be for a pair of generating sections satisfying the Lie bracket conditions (see Remark \ref{rema: defining a Lie algebroid locally out of an anchored vector bundle}); then we extend the definition so that it satisfies the Leibniz rule. The extensions are
\begin{align*}
    &[\hat{\alpha}_1, \hat{\alpha}_2]_{\DNC(A,B)} \coloneqq t \cdot \widehat{[\alpha_1,\alpha_2]}_{\algebr}, \\
    &[\DNC(\beta),\hat{\alpha}]_{\DNC(A,B)} \coloneqq \widehat{[\beta,\alpha]}_{\algebr}, \\
    &[\DNC(\beta_1),\DNC(\beta_2)]_{\DNC(A,B)} \coloneqq \DNC([\beta_1,\beta_2]_{\algebr}).
\end{align*}
Since the conditions for the Lie bracket obviously hold on $\algebr \times \mathbb{R}^\times \subset \DNC(\algebr,\subalgebr)$, the relations hold on all of $\DNC(\algebr,\subalgebr)$ since $\algebr \times \mathbb{R}^\times \subset \DNC(\algebr,\subalgebr)$ is dense (so that there is at most one extension). This proves the statement.
\end{proof}
\begin{rema}\cite{meinrenken}\label{rema: N of algebroid}
From the above we obtain a Lie algebroid structure on $d_\normal\pi: \normal(\algebr,\subalgebr) \ra \normal(\base,\subbase)$. The anchor map is given by $d_\normal\anchor_\algebr: \normal(\algebr,\subalgebr) \ra \normal(T\base,T\subbase) \cong T\normal(\base,\subbase)$, and the Lie bracket is determined by 
\begin{align*}
    [\alpha_1|_\subbase \mod \subalgebr,\alpha_2|_\subbase \mod \subalgebr]_{\normal(\algebr,\subalgebr)} &= 0, \\
    [d_\normal \beta,\alpha|_\subbase \mod \subalgebr]_{\normal(\algebr,\subalgebr)} &= [\beta,\alpha]_{\algebr}|_\subbase \mod \subalgebr, \\
    [d_\normal \beta_1, d_\normal \beta_2]_{\normal(\algebr,\subalgebr)} &= d_\normal[\beta_1,\beta_2]_\algebr
\end{align*}
where we note that sections of the form $\alpha|_\subbase \mod \subalgebr$ are the core sections, and sections of the form $d_\normal \beta$ are the linear sections (which generate $\Gamma(\normal(\algebr,\subalgebr))$; see Section \ref{sec: DVB}). In particular, we see that the diagram
\begin{center}
    \begin{tikzcd}
    \normal(\algebr,\subalgebr) \ar[r,"d_\normal\pi"] \ar[d] & \normal(\base,\subbase) \ar[d] \\
    \subalgebr \ar[r] & \subbase
    \end{tikzcd}
\end{center}
is a $\VB$-algebroid (see Theorem \ref{theo: alternative defintion of VB-algebroid} and the remark before it). 
\end{rema}
%\begin{exam}\label{exam: algebroid of N of groupoid}
%Let $(\groupoid,\subgroupoid)$ be a pair of groupoids and denote by $\pi_\algebr: \algebroid$ and $\pi_\subalgebr: \subalgebroid$ the associated Lie algebroids respectively. The Lie algebroid $C \ra \normal(\base,\subbase)$ of the Lie groupoid $\normal(\group,\subgroup) \rra \normal(\base,\subbase)$ can canonically be identified with the Lie algebroid $\normal(\algebr,\subalgebr) \ra \normal(\base,\subbase)$, i.e. 
%\[\normal(\normal(\group,\subgroup),\normal(\base,\subbase)) \cong \normal(\normal(\group,\base),\normal(\subgroup,\subbase))\]
%as Lie algebroids. We have a map
%\[\normal(\group,\subgroup) \ra \normal(\]
%\end{exam}
With all the work that we've done in the former section, it is easy to show that the Lie algebroid of a deformation to the normal cone Lie groupoid coincides with the Lie algebroid structure defined above.

\begin{prop}\label{prop: algebroid of DNC of groupoid}
Let $(\groupoid,\subgroupoid)$ be a pair of groupoids and denote by $(\algebroid,\subalgebroid)$ the associated pair of Lie algebroids. The Lie algebroid of $\DNC(\group,\subgroup) \rra \DNC(\base,\subbase)$ is $\DNC(\algebr,\subalgebr) \ra \DNC(\base,\subbase)$.
\end{prop}
\begin{proof}
We will show that the Lie algebroid $C \ra \DNC(\base,\subbase)$ of $\DNC(\group,\subgroup) \rra \DNC(\base,\subbase)$ can naturally be identified with $\DNC(\algebr,\subalgebr) \xra{\DNC(\pi_\algebr)} \DNC(\base,\subbase)$. To see why, observe that it suffices to prove the claim: $\DNC(\algebr,\subalgebr)$ is the Lie algebroid of $\DNC(\group,\subgroup)$ as a manifold, and $d\DNC(\source)|_{\identity_{\DNC(\group,\subgroup)(\DNC(\base,\subbase)})}$ becomes the map $\DNC(\pi_\algebr)$ under this identification. Indeed, the Lie algebroid of $\DNC(\group,\subgroup)$ contains the Lie subgroupoid $\algebr \times \mathbb{R}^\times$, and $\DNC(\pi_\algebr): \DNC(\algebr,\subalgebr) \ra \DNC(\base,\subbase)$ has a unique Lie algebroid structure for which $\algebr \times \mathbb{R}^\times$ is a Lie subalgebroid. To prove the claim, note that the canonical map 
\[\hat{t}_\group: \DNC(\group,\subgroup) \ra \mathbb{R}\]
equals the map
\[\DNC(\group,\subgroup) \xra{\DNC(\source)} \DNC(\base,\subbase) \xra{\hat{t}_\base} \mathbb{R},\]
because $\DNC(\source)\left(\DNC(\group,\subgroup)_t\right) \subset \DNC(\base,\subbase)_t$. Hence, using the identification $\ker d\hat{t}_\group \cong \DNC(T\group,T\subgroup)$, we can view $C = \ker d\DNC(\source)|_{\DNC(\identity)(\DNC(\base,\subbase))} \subset \DNC(T\group,T\subgroup)$. The claim now follows by Proposition \ref{prop: DNC of differentials}.
\end{proof}
To end this section, we observe that the functorial behaviour of the deformation to the normal cone construction extends to the setting of Lie groupoids and Lie algebroids.
\begin{prop}\label{prop: groupoid morphism DNC}
Let $F: (\group,\subgroup) \ra (\grouptwo,\subgrouptwo)$ be a morphism of pairs of Lie groupoids between pairs of Lie groupoids $(\groupoid,\subgroupoid)$ and $(\grouptwoid,\subgrouptwoid)$ (that is, $F: \group \ra \grouptwo$ is a morphism of Lie groupoids and a map of pairs). Then
\[\DNC(F): \DNC(\group,\subgroup) \ra \DNC(\grouptwo,\subgrouptwo)\]
is a morphism of Lie groupoids.
\end{prop}
\begin{proof}
Denote the base map by $f: (\base,\subbase) \ra (\basetwo,\subbasetwo)$. The statement follows immediately from the functorial behaviour of the deformation to the normal cone construction. Indeed,
\begin{align*}
    \DNC(F)(\DNC(\mult_\group)(g,h)) &= \DNC(F \circ \mult_\group)(g,h) \\
    &= \DNC(\mult_\grouptwo \circ (F \times F))(g,h) = \DNC(\mult_\grouptwo)(\DNC(F)(g),\DNC(F)(h))
\end{align*}
(and $\DNC(F)(\DNC(\identity_\group)(g)) = \DNC(\identity_\grouptwo)(\DNC(f)(g))$), as required.
\end{proof}
The proof of the analogue statement for Lie algebroids is very similar to a proof we will give later for blow-ups, and is therefore postponed (see Proposition \ref{prop: DNC of subalgebroid is subalgebroid} and Corollary \ref{coro: LA morphism to DNC}).

Recall that we still did not prove that the Lie algebroid orbits are immersed submanifolds. We will do this by proving a much more powerful theorem: the splitting theorem for Lie algebroids. With the theory of deformation to the normal cones at hand, we can give an interesting proof. Before we prove it, however, we will use the notion of Euler-like vector fields. This notion rather naturally comes up in the next section, so we will delay the proof of the theorem until the next section. 

In the next section, we will look into the theory of deformation to the normal cones from an algebraic point of view. This sheds new light on the topic%, and it will also give a new way to look at blow-ups later on
.

\subsection{Deformation to the normal cone: an algebraic approach}\label{sec: Deformation to the normal cone: an algebraic approach}
Here we will describe the deformation to the normal cone construction algebraically. It is largely based on \cite{sadegh2018eulerlike}. Notice that it is most sensible to \textbf{require manifolds to be Hausdorff in this section}, since then the algebra of smooth functions encodes all of the information of a manifold. The first thing we will do, is be precise about what we mean when ``approaching smooth manifolds from an algebraic perspective''. One way to do this, is via the following three definitions. We fix a commutative algebra $R$ (associative and with unit) over $\mathbb{R}$.

\begin{defn}\cite{sadegh2018eulerlike}\label{defn: spectrum of comm algebra}
A character of $R$ is a non-zero algebra morphism $R \ra \mathbb{R}$. We call the set of characters the \textit{spectrum} of $R$ and denote it by $\spec R$. Moreover, we equip it with the coarsest topology such that, for all $r \in R$, the evaluation map
\[\hat{r}: \spec R \ra \mathbb{R} \textnormal{ given by } \chi \mapsto \chi(r)\]
is a continuous map.
\end{defn}
\begin{defn}\cite{sadegh2018eulerlike}\label{defn: smallest sheaf that contains nice global sections}
We denote by $\mathcal{F}_R$ the smallest subsheaf of $C^0(\spec R)$ (the sheaf of continuous real-valued maps on $\spec R$) that includes all global sections of the form 
\[g \circ (\hat{r}_1,\dots,\hat{r}_k) = g(\hat{r}_1,\dots,\hat{r}_k),\]
where $k>0$, $r_1,\dots,r_k \in R$ and $g \in C^\infty(\mathbb{R}^k)$ (the sheaf of smooth real-valued maps on $\mathbb{R}^k$).
\end{defn}
\begin{defn}\cite{sadegh2018eulerlike}\label{defn: smoothly generated sheaves}
Let $\mathcal{F}$ be a sheaf of continuous real-valued maps on a topological space $Y$ and let $U \subset Y$ be an open subset. The elements $f_1,\dots,f_k \in \mathcal{F}(U)$ \textit{smoothly generate} $\mathcal{F}(U)$ if whenever $f \in \mathcal{F}(U)$ we can find a smooth map $g \in C^\infty(\mathbb{R}^k)$ such that
\[f=g(f_1,\dots,f_k).\]
In case that $\mathcal{F}=\mathcal{F}_R$, we will say that elements $r_1,\dots,r_k \in R$ generate $\mathcal{F}(U)$ if the above condition holds for $\hat{r}_1|_U,\dots,\hat{r}_k|_U$.
\end{defn}
\begin{prop}\cite{sadegh2018eulerlike}\label{prop: when spectrum is a manifold}
The spectrum $\spec R$ of $R$ is a smooth manifold of dimension $n$ with $\mathcal{F}_R$ its sheaf of smooth functions if and only if for all elements of $\spec R$ there is an open subset $U \subset \spec R$ and $r_1,\dots,r_n \in R$ such that $r_1,\dots,r_n$ smoothly generate $\mathcal{F}_R(U)$ and
\[\varphi = (\hat{r}_1,\dots,\hat{r}_n): U \ra \mathbb{R}^n\]
is a homeomorphism onto its image.
\end{prop}
\begin{proof}
Suppose we can find for all $\chi \in \spec R$ an open subset $U \subset \spec R$ and $r_1,\dots,r_n \in R$ as stated above. To check whether $\spec R$ is a manifold, it suffices to check that the transition maps are smooth. So, let $\varphi: U \ra \mathbb{R}^n$ and $\varphi': U' \ra \mathbb{R}^n$ be charts as above, determined by $r_1,\dots,r_n$ and $r_1',\dots,r_n'$, respectively. Then $g=\varphi' \circ \varphi^{-1}: \varphi(U \cap U') \ra \mathbb{R}^n$ satisfies
\[\varphi' = g(\hat{r}_1,\dots,\hat{r}_n).\]
Moreover, $g$ is the unique map to $\mathbb{R}^n$, defined on the open subset $\varphi(U \cap U')$ of $\mathbb{R}^n$, that satisfies this property. Since $r_1,\dots,r_n$ smoothly generates $\mathcal{F}_R$, the map $g$ must be smooth.

Conversely, to see that if $\spec R$ is a smooth manifold with $\mathcal{F}_R$ its sheaf of smooth functions, suppose that $(U,\varphi)$ is a chart around $\chi \in \spec R$, with $\varphi=(x^1,\dots,x^n)$. Then, by definition of $\mathcal{F}_R$, if we shrink $U$, then for each $1 \le i \le n$, we can find a smooth function $g^i \in C^\infty(\mathbb{R}^{k_i})$ such that
\[x^i = g^i(\hat{r}_1^i,\dots,\hat{r}_{k_i}^i).\]
We can replace $g^i$ with a smooth function $h^i$ of $\mathbb{R}^{k}$, with $k=\textstyle\sum_i k_i$, by setting $h^i \coloneqq g^i \circ \pr_i$, where $\pr_i: \mathbb{R}^k \rightarrowdbl \mathbb{R}^{k_i}$ is the projection with respect to the decomposition $\mathbb{R}^{k_1} \times \dots \times \mathbb{R}^{k_n}$. Then
\[x^i = h^i(\hat{r}_1,\dots,\hat{r}_k) = g^i(\hat{r}_1^i,\dots,\hat{r}_{k_i}^i),\]
where $\hat{r}_{(i-1)k_{i-1}+j} \coloneqq \hat{r}_j^i$ (with $k_0=0$). By replacing the $\hat{r}_j$ with $\hat{r-\chi(r_j)}_j$ (this makes sense: $R$ is an $\mathbb{R}$-algebra, and $\mathbb{R}$ is a field, so we can view $\mathbb{R}$ as a subring of $R$) and changing $h^i$ by a translation, we may assume $\hat{r}_j(\chi)=0$ for all $1 \le j \le k$. Now, since $\varphi$ is a chart, we can shrink $U$ once again to ensure that, for all $1 \le j \le k$, we can write
\[\hat{r}_j = f^j \circ \varphi\]
for some $f^j \in C^\infty(\mathbb{R}^n)$. Assemble the maps $h^i$ and $f^j$ into smooth maps $h: \mathbb{R}^k \ra \mathbb{R}^n$ and $f: \mathbb{R}^n \ra \mathbb{R}^k$, respectively. Then $h \circ f: \mathbb{R}^n \ra \mathbb{R}^n$ is the identity map when restricted to $\varphi(U) \ni 0$. In particular, $h$ is a submersion, so we can find $n$ standard basis elements $e^{k_1},\dots,e^{k_n}$ of $\mathbb{R}^k$ such that when we post-compose the map 
\[\mathbb{R}^n \ra \mathbb{R}^n \textnormal{ given by } e^\ell \mapsto e^{k_\ell}\]
with $h$, it becomes a local diffeomorphism around $0$. The collection $r_{k_1},\dots,r_{k_n}$ has the desired properties.
\end{proof}

We are almost ready to define the deformation to the normal cone with this algebraic viewpoint, but this still requires some preparation. We fix here a pair of smooth manifolds $(\base,\subbase)$.
\begin{defn}\cite{Higson}\label{defn: vanishing to order}
Let $f \in C^\infty(\base)$. The smooth function $f$ is said to \textit{vanish to order $k$ along $\subbase$} if we can find for all $x \in \base$ a chart $(U,\varphi)$ around $x$ such that $f \circ \varphi^{-1}$ is a sum of functions whose terms are products of $k$ or more smooth functions on $\varphi(U)$ that vanish on $\varphi(U \cap \subbase)$. We denote the ideal of $C^\infty(\base)$ consisting of those $f \in C^\infty(\base)$ that vanish to order $k$ along $\subbase$ by $I_k(\base,\subbase)$.
\end{defn}
\begin{term}\label{term: vanishing to order}
If $f \in C^\infty(\base)$, then, for all $k \le 0$, we will say that, always, $f$ vanishes to order $k$ along $\subbase$. 
\end{term}
Observe that if $f \in C^\infty(\base)$ vanishes to order $k$, then it also vanishes to order $k-1$.
\begin{lemm}\cite{Higson}\label{lemm: vanishing to order of vector fields}
Let $\sigma \in \mathfrak{X}(\base)$ and let $f \in C^\infty(\base)$. If $f$ vanishes to order $k$ along $\subbase$, then $\sigma(f)$ vanishes to order $k-1$ along $\subbase$.
\end{lemm}
\begin{proof}
Write $\sigma$ locally as $\textstyle\sum_i b^i \frac{\partial}{\partial x^i}$. It suffices to prove the statement in case $f$ locally can be written as $f^1\cdots f^k$ for smooth functions $f^1,\dots,f^k$ defined on an opens set of $\mathbb{R}^n$ that vanish along $\mathbb{R}^p$. Then
\begin{align*}
    \sigma(f) & = \sum_i b^i \frac{\partial f}{\partial x^i} = \sum_i b^i \frac{\partial f^{1} \cdots f^{k}}{\partial x^i},
\end{align*}
and by using the Leibniz rule, we can write it as a sum in which each term is of the form
\[b^if^{1}\cdots\widehat{f}^{j}\cdots f^{k} \frac{\partial f^{j}}{\partial x^i}.\]
These terms are written as a product of $k-1$ smooth functions $b^if^1, f^2, \dots, f^{k-1}, f^k \tfrac{\partial f^j}{\partial x^i}$ that vanish along $\subbase$, so this concludes the proof.
\end{proof}
The deformation to the normal cone algebraically takes the following form:
\begin{defn}\cite{sadegh2018eulerlike}\label{defn: algebraic deformation to the normal cone}
The \textit{(algebraic) deformation to the normal cone} is the spectrum associated to the $\mathbb{R}$-algebra
\[R_\DNC(\base,\subbase) = 
\{\sum_{k \in \mathbb{Z}} f_k t^{-k} \mid f_k \textnormal{ vanishes to order } k \textnormal{ on } \subbase\} \subset C^\infty(\base)[t,t^{-1}]\]
(note: when we write $\textstyle\sum_{k \in \mathbb{Z}} f_k t^{-k}$, all but finitely many $f_k$ are understood to be zero).
\end{defn}
\begin{rema}\cite{sadegh2018eulerlike}\label{rema: algebraic deformation to the normal cone}
Notice that $R_\DNC(\base,\subbase)$ is indeed an $\mathbb{R}$-algebra: from the definition of order of vanishing, it is clear that whenever $f_k \in C^\infty(\base)$ vanishes to order $k$ along $\subbase$, and $f_\ell \in C^\infty(\base)$ vanishes to order $\ell$ along $\subbase$, then $f_kf_\ell$ vanishes to order $k+\ell$ along $\subbase$. The ring map $\mathbb{R} \ra R_\DNC(\base,\subbase)$ is given by $\lambda \mapsto \lambda t^0$. 

Alternatively, $R_\DNC(\base,\subbase)$ is the $\mathbb{R}$-algebra generated by $C^\infty(\base)$, $t$, and $t^{-1} \cdot I_1(\base,\subbase)$.
\end{rema}
To give some motivation for the definition, let $f \in C^\infty(\base)$ be a smooth map that vanishes to order $k$ along $\subbase$. Then the map $(x,t) \mapsto t^{-k}f(x)$ on $\base \times \mathbb{R}^\times \subset \DNC(\base,\subbase)$ extends to a smooth map
\[\DNC^k(f): \DNC(\base,\subbase) \ra \mathbb{R} \textnormal{ given by } z \mapsto \begin{cases} d_\normal^{(k)} f(y)\xi & \text{if}\ z=(y,\xi) \\ \tfrac{1}{t^k} f(x) & \text{if}\ z=(x,t)\end{cases}\]
precisely because $f$ vanishes to order $k$ along $\subbase$ (the map $d_\normal^{(k)} f: \normal(\base,\subbase) \ra \mathbb{R}$ is obtained by applying $d_\normal$ to $f$ $k$ times; this makes sense because $\normal(E,0_M) \cong E$ for a vector bundle $E \ra M$, and $f$ vanishes to order $k$ along $\subbase$), so we can construct smooth maps for $\DNC(\base,\subbase)$ out of $R_\DNC(\base,\subbase)$. Conversely, we will show that $R_\DNC(\base,\subbase)$ carries a smooth structure from which we naturally obtain the smooth maps from Remark \ref{rema: smooth functions DNC}. For example, the analogue of the submersion $\hat{t}: \DNC(\base,\subbase) \ra \mathbb{R}$ in this algebraic set-up is the element $t \in R_\DNC(\base,\subbase)$, which gives rise to the continuous map
\[\hat{t}: \spec R_\DNC(\base,\subbase) \ra \mathbb{R}\]
(this is why we used the $\hat{ }$-notation here as well).
\begin{rema}\cite{Higson}\label{rema: characters DNC}
Let $(x,s) \in \DNC(\base,\subbase)_s$ ($s\neq0$) and $(y,\xi) \in \DNC(\base,\subbase)_0$. Let $\sigma_\xi \in \mathfrak{X}(\base)$ such that $\sigma_\xi(y) \mod T_y\subbase = \xi$. Then these elements induce characters of the algebraic deformation to the normal cone:
\begin{align*}
    \chi_{(x,s)} &: R_\DNC(\base,\subbase) \ra \mathbb{R} \textnormal{ given by } \sum_k f_kt^{-k} \mapsto \sum_k f_k(x)s^{-k} \textnormal{ and } \\
    \chi_{(y,\xi)} &: R_\DNC(\base,\subbase) \ra \mathbb{R} \textnormal{ given by } \sum_k f_kt^{-k} \mapsto \sum_k \frac{1}{k!} \sigma_\xi^k(f_k)(y),
\end{align*}
where $\sigma^k_\xi$ denotes the $k$-th iterated derivative along $\sigma_\xi$ (note: the values $\sigma_\xi^k(f_k)(y)$ depend only on $\xi$). In fact, we will see that these are all characters of $R_\DNC(\base,\subbase)$. Since they are all distinct, it will follow that we can identify $\spec R_\DNC(\base,\subbase)$ with $\DNC(\base,\subbase)$.
\end{rema}
Notice that
\[\hat{t}(\chi_{(x,s)}) = s, \textnormal{ and } \hat{t}(\chi_{(y,\xi)}) = 0.\]
To make the above statement about $\hat{t}$ more precise, we will show that it satisfies the following properties.
\begin{prop}\cite{sadegh2018eulerlike}\label{prop: characters DNC}
Let $\lambda \in \mathbb{R}$. Then 
\[\hat{t}^{-1}(\lambda) \cong \spec R_\DNC^\lambda(\base,\subbase),\]
where $R_\DNC^\lambda(\base,\subbase) = R_\DNC(\base,\subbase)/(t-\lambda)$. Moreover, if $\lambda \neq 0$, then we have 
\[R_\DNC^\lambda(\base,\subbase) \cong C^\infty(\base) \textnormal{ and } R_\DNC^0(\base,\subbase) \cong \bigoplus_{k \ge 0} I_k(\base,\subbase)/I_{k+1}(\base,\subbase),\]
where $\textstyle\bigoplus_{k \ge 0} I_k(\base,\subbase)/I_{k+1}(\base,\subbase)$ is viewed as a graded $\mathbb{R}$-algebra.
\end{prop}
\begin{proof}
To show that $\hat{t}^{-1}(\lambda) \cong \spec R^\lambda_\DNC(\base,\subbase)$, notice that the natural map
\[R_\DNC(\base,\subbase) \ra R_\DNC^\lambda(\base,\subbase) \textnormal{ given by } r \mapsto r \mod (t-\lambda)\]
induces a continuous map 
\[\spec R_\DNC^\lambda(\base,\subbase) \ra \spec R_\DNC(\base,\subbase) \textnormal{ given by } \chi \mapsto [r \mapsto \chi(r \mod (t-\lambda))].\]
Obviously, this map maps onto $\hat{t}^{-1}(\lambda)$. It has an obvious inverse:
\[\hat{t}^{-1}(\lambda) \ra \spec R_\DNC^\lambda(\base,\subbase) \textnormal{ given by } \chi \mapsto [r \mod (t-\lambda) \mapsto \chi(r)],\]
which is well-defined, because $\chi(t)=\lambda$ for all $\chi \in \hat{t}^{-1}(\lambda)$, and it is not hard to see that this map is continuous as well. 

To see that $R^\lambda_\DNC(\base,\subbase) \cong C^\infty(\base)$, we rewrite an element $\textstyle\sum_kf_kt^{-k}=0$ as
\[\sum_kf_kt^{-k} = (t-\lambda)\sum_k(\sum_\ell f_{k-\ell}\lambda^\ell)t^{-k-1}.\]
In particular, we see that the evaluation map
\[R_\DNC(\base,\subbase) \ra C^\infty(\base) \textnormal{ given by } \sum_k f_kt^{-k} \mapsto \sum_kf_k\lambda^{-k}\]
has kernel $(t-\lambda)$, so, indeed, $R_\DNC^\lambda(\base,\subbase) \cong C^\infty(\base)$. Lastly, the natural map
\[R_\DNC(\base,\subbase) \ra \bigoplus_{k \ge 0} I_k(\base,\subbase)/I_{k+1}(\base,\subbase) \textnormal{ given by } \sum_k f_kt^{-k} \mapsto (f_k \mod I_{k+1})_k\]
has kernel $(t)$, so $R_\DNC^0(\base,\subbase) \cong \textstyle\bigoplus_{k \ge 0} I_k(\base,\subbase)/I_{k+1}(\base,\subbase)$. This concludes the proof.
\end{proof}
Notice that we can identify $I_0(\base,\subbase)/I_1(\base,\subbase)$ with $C^\infty(\subbase)$, so we can view the degree zero part of the graded algebra $R_\DNC^0(\base,\subbase)$ as $C^\infty(\subbase)$. Moreover, for a map $\chi: R_\DNC^0(\base,\subbase)  \ra \mathbb{R}$ to be a character, it has to restrict to a character $C^\infty(\subbase) \ra \mathbb{R}$, and such a character is given by the evaluation of a smooth function (note: the Hausdorffness assumption is important here). So, if we define 
\[I^0(\base,\subbase)_y \subset R^0_\DNC(\base,\subbase), \textnormal{ with } y \in \subbase,\] 
to be the ideal in $R^0_\DNC(\base,\subbase)$ generated by the vanishing ideal of $y$ in $C^\infty(\subbase)$, then we see that $\chi$ factors through $R_\DNC^0(\base,\subbase)_y \coloneqq R_\DNC^0(\base,\subbase)/I^0(\base,\subbase)_y$. 

Before we move on, recall Theorem \ref{theo: smooth polynomial functions of vector bundle}, which relates vector bundles with their graded $\mathbb{R}$-algebra consisting of fiberwise polynomial functions $C^\infty_{\textnormal{poly}}(E)$.
\begin{prop}\cite{sadegh2018eulerlike}\label{prop: characters DNC 2}
Let $y \in \subbase$. The $\mathbb{R}$-algebra $R_\DNC^0(\base,\subbase)_y$ is naturally isomorphic to the $\mathbb{R}$-algebra of real-valued polynomials $P_\normal(\base,\subbase)_y$ on $T_y\base/T_y\subbase$.
\end{prop}
\begin{proof}
We identify $R_\DNC^0(\base,\subbase)$ with $\textstyle\bigoplus_{k \ge 0} I_k(\base,\subbase)/I_{k+1}(\base,\subbase)$, and define the map
\[R_\DNC^0(\base,\subbase) \ra P_\normal(\base,\subbase)_y \textnormal{ by } (f_k \mod I_{k+1})_k \mapsto (\xi \mapsto \sum_k \frac{1}{k!} \sigma_\xi^k(f_k)(y))\]
(where, as before, $\sigma_\xi \in \mathfrak{X}(\base)$ with $\sigma_\xi(y) \mod T_y\subbase = \xi$). Observe that this map is well-defined by Lemma \ref{lemm: vanishing to order of vector fields}. The restricted map to degree $1$ has the vanishing ideal of $y$ in $C^\infty(\subbase)$ as its kernel. Since the above map is completely determined by this restricted map, and $I^0(\base,\subbase)_y$ is the ideal generated by the aforementioned vanishing ideal, the map has kernel $I^0(\base,\subbase)_y$, so we obtain the desired isomorphism.
\end{proof}
\begin{rema}\cite{sadegh2018eulerlike}\label{rema: smooth structure algebraic DNC}
From the above two propositions (and the remark before Proposition \ref{prop: characters DNC 2}) we can conclude that the fibers $\hat{t}^{-1}(\lambda)$ of $\hat{t}$ can be identified with $\base$ if $\lambda \neq 0$, and with $\normal(\base,\subbase)$ if $\lambda=0$.
Therefore, $\spec R_\DNC(\base,\subbase) = \DNC(\base,\subbase)$; at least set-theoretically. Since $R_\DNC(\base,\subbase)$ is generated by the elements $t$,$f_0t^0$, and $f_1t^{-1}$, where $f_0 \in C^\infty(\base)$, and $f_1 \in C^\infty(\base)$ vanishes along $\subbase$, the topology on $\spec R_\DNC(\base,\subbase)$ agrees with the topology on $\DNC(\base,\subbase)$. Indeed, these maps induce the continuous maps $\hat{t}$, $\hat{f}_0$, and $\widehat{t^{-1}f_1} = \DNC(f_1)$, respectively (see Remark \ref{rema: smooth functions DNC}). Moreover, if $\varphi=(y^1,\dots,y^p,x^1,\dots,x^q)$ are coordinates on an open subset $U$ of $\base$ (with $V \coloneqq U \cap \subbase=\{x^j=0\}$), then the collection
\[y^1,\dots,y^p,t^{-1}x^1,\dots,t^{-1}x^q,t \in R_\DNC(\base,\subbase)\]
satisfies the hypotheses from Proposition \ref{prop: when spectrum is a manifold}; indeed, $y^1,\dots,y^p,t^{-1}x^1,\dots,t^{-1}x^q,t$ smoothly generates $\mathcal{F}_{R_\DNC(\base,\subbase)}(\DNC(U,V))$ by a version of Taylor's theorem, and 
\[(\hat{y}^1,\dots,\hat{y}^p,\widehat{t^{-1}x^q},\dots,\widehat{t^{-1}x^q},\hat{t}): \DNC(U,V) \ra \mathbb{R}^{n+1}\] 
is a homeomorphism onto the open set $\Omega^U_V = \{(y,\xi,t) \in \mathbb{R}^{n+1} \mid (y,t\xi) \in \varphi(U)\}$. Of course, the smooth structure defined here on $\spec R_\DNC(\base,\subbase)$ agrees with the smooth structure defined on $\DNC(\base,\subbase)$.
\end{rema}

\subsection{Euler-like vector fields and the splitting theorem for Lie algebroids}
To end this chapter, we will give a proof of the \textit{splitting theorem} for Lie algebroids. We will prove this by using the deformation to the normal cone construction, and we will state the theorem in terms of so-called \textit{Euler-like vector fields}. The theory of flows of vector fields will be important in the rest of the section, so it is good to be explicit about our conventions.
\begin{term}\label{term: Lie derivative}
If $\sigma \in \mathfrak{X}(\base)$ is a vector field, then we denote by $\phi^t_\sigma$ its flow, and the Lie derivative is given as $\Lie(\sigma) = \left.\tfrac{d}{dt}\right\vert_{t=0} (\phi^t_\sigma)^*$.
\end{term}
\begin{defn}\cite{meinrenken}\label{defn: Euler vector field}
The \textit{Euler vector field} of a vector bundle $E \ra M$ is the vector field $\mathcal{E}_E \in \mathfrak{X}(E) = \Gamma(TE)$ given, as a derivation, by 
\[\mathcal{E}_E(f) \coloneqq \left.\frac{d}{dt}\right\vert_{t=0} (e^t)^*f,\]
where $f \in C^\infty(E)$.
\end{defn}
Notice that the Euler vector field of a vector bundle $E \ra M$ is complete, and also tangent to the fibers of $E$. Moreover, if $f \in C^\infty(M)$, and it is linear (that is, locally, we can always write $f$ as a linear function), then $\mathcal{E}_E(f)=f$. The Euler vector field of $\mathbb{R}^k$ (i.e. fiberwise on $E$) is
\begin{equation}\label{eq: fiberwise Euler-vf}
    \mathcal{E}_{\mathbb{R}^k} = \sum_i x^i \frac{\partial}{\partial x^i}.
\end{equation}
\begin{defn}\cite{meinrenken}\label{defn: Euler-like}
A vector field $\sigma \in \mathfrak{X}(\base)$ is said to be \textit{Euler-like along $\subbase$} if it is complete, 
tangent to $\subbase$, and
\[d_\normal\sigma: \normal(\base,\subbase) \ra \normal(T\base,T\subbase) \cong T\normal(\base,\subbase)\]
is the Euler vector field $\mathcal{E}_{\normal(\base,\subbase)}$ of $\normal(\base,\subbase) \ra \subbase$.
\end{defn}
Notice that an Euler-like vector field $\sigma \in \mathfrak{X}(\base)$ along $\subbase$ is not only tangent to $\subbase$, but it has to vanish along $\subbase$, because $\mathcal{E}_{\normal(\base,\subbase)}$ is tangent to the fibers of $\normal(\base,\subbase)$. It follows that (equivalently)
%One can check that the property $\mathcal{E}_{\normal(\base,\subbase)}(f)-f=0$ for all linear functions $f \in C^\infty(\subbase)$ implies
\begin{equation}\label{eq: Euler-like}
    \sigma(f_1) - f_1 \in I_{2}(\base,\subbase) \textnormal{ for all } f_1 \in I_{1}(\base,\subbase)
\end{equation}
by applying the normal derivative to $\sigma(f_1) - f_1$ (see Remark \ref{rema: characters DNC}).
\begin{rema}\cite{meinrenken}\label{rema: local coordinates Euler-like vector field}
Euler-like vector fields $\sigma \in \mathfrak{X}(\base)$ are characterised locally as follows. Let $\varphi=\{y^i,x^j\}$ be adapted coordinates on an open subset $U \subset \base$ (so $V \coloneqq U \cap \subbase = \{x^j=0\}$). From \eqref{eq: fiberwise Euler-vf} we see that, locally, 
\[d_\normal \sigma = \sum_i x^i \frac{\partial}{\partial x^i},\]
so 
\[\sigma = \sum_i (x^i+b^i)\frac{\partial}{\partial x^i} + \sum_j a^j \frac{\partial}{\partial y^j}, \textnormal{ with } a^j,b^i \in C^\infty(\base),\]
where $a^j \circ \varphi^{-1}(x,0)=0$ and $b^i \circ \varphi^{-1}(x,0)$ vanishes to order $2$. Moreover, it is readily verified that if we can write a complete vector field $\sigma \in \mathfrak{X}(\base)$ like this locally, then $\sigma$ is an Euler-like vector field along $\subbase$.
\end{rema}
The deformation to the normal cone construction comes into play via the following characterisation of Euler-like vector fields.
\begin{prop}\cite{Bischoff_2020}\label{prop: Euler-like and deformation to the normal cone}
A vector field $\sigma \in \mathfrak{X}(\base)$ is Euler-like if and only if the vector field
\[W_\sigma \coloneqq \frac{1}{t}\sigma + \frac{\partial}{\partial t} \in \mathfrak{X}(\base \times \mathbb{R}^\times)\]
extends to a smooth vector field on $\DNC(\base,\subbase)$.
\end{prop}
\begin{proof}
Denote by $\mu_\lambda: \DNC(\base,\subbase) \ra \DNC(\base,\subbase)$ the multiplication by $\lambda \in \mathbb{R}^\times$ on $\DNC(\base,\subbase)$ with respect to the canonical $\mathbb{R}^\times$-action. We have a smooth map 
\begin{equation}\label{eq: multiplication DNC as R-action}
    \mu: \DNC(\base,\subbase) \times \mathbb{R} \ra \DNC(\base,\subbase) \textnormal{ given by } (z,s) \mapsto \mu_{e^s}(z)
\end{equation} 
which is the flow of a vector field; namely, $t \tfrac{\partial}{\partial t}$ on $\base \times \mathbb{R}^\times \subset \DNC(\base,\subbase)$ and $-\mathcal{E}_{\normal(\base,\subbase)}$ on $\normal(\base,\subbase) \subset \DNC(\base,\subbase)$ (see Definition \ref{defn: Euler vector field}). Moreover, if $\sigma \in \mathfrak{X}(\base)$, then we have the vector field $\DNC(\sigma)$ of $\DNC(\base,\subbase)$, which is given by $\sigma \times 0$ (from now on just denoted by $\sigma$) on $\base \times \mathbb{R}^\times \subset \DNC(\base,\subbase)$ and $d_\normal\sigma$ on $\normal(\base,\subbase) \subset \DNC(\base,\subbase)$. Now, the vector field $tW_\sigma = \sigma + t \tfrac{\partial}{\partial t}$ on $\base \times \mathbb{R}^\times \subset \DNC(\base,\subbase)$ extends to a smooth vector field on $\DNC(\base,\subbase)$, and it vanishes along $\normal(\base,\subbase)$ if and only if $\sigma$ is Euler-like. Lastly, $W_\sigma$ extends to a smooth vector field if and only if $tW_\sigma$ vanishes along $\DNC(\base,\subbase)_0 \cong \normal(\base,\subbase)$, so the result follows.
\end{proof}
The statement says, very roughly speaking, that Euler-like vector fields are those vector fields that can be ``deformed'' on $\DNC(\base,\subbase)$ into the Euler vector field of $\normal(\base,\subbase)$. In fact, by first examining some basic properties of the vector field $W_\sigma$ from above, one can show the following result that leads to many normal form results.
\begin{theo}\cite{Bischoff_2020}\label{theo: Euler-like vector fields and tubular neighbourhood embeddings}
To each Euler-like vector field $\sigma$ on $\base$ there is a unique tubular neighbourhood embedding $\chi=\chi_\sigma$ such that $\mathcal{E}_{\normal(\base,\subbase)}$ is $\chi$-related to $\sigma$.
\end{theo}
We fix an Euler-like vector field $\sigma \in \mathfrak{X}(\base)$. The properties of $W_\sigma$ (from Proposition \ref{prop: Euler-like and deformation to the normal cone}) that we need are the following. 
\begin{lemm}\cite{Bischoff_2020}\label{lemm: properties of W_sigma 1}
For all $\tau \in \mathfrak{X}(\base)$, $\tau + [\sigma,\tau] \in \mathfrak{X}(\base)$ is tangent to $\subbase$ and
\[[W_\sigma,\hat{\tau}] = \DNC(\tau + [\sigma,\tau]) \in \mathfrak{X}(\DNC(\base,\subbase))\]
(where $\hat{\tau}$ is the extension of $t(\tau \times 0) \in \mathfrak{X}(\base \times \mathbb{R}^\times)$ to $\DNC(T\base,T\subbase) \subset T\DNC(\base,\subbase)$; see Remark \ref{rema: local sections of DNC} and Proposition \ref{prop: DNC of tangent spaces}). If $\tau$ is in addition tangent to $\subbase$, then,
\[\frac{1}{t}\DNC([\sigma,\tau])|_{\base \times \mathbb{R}^\times} \in \mathfrak{X}(\base \times \mathbb{R}^\times)\]
extends to a smooth vector field on $\DNC(\base,\subbase)$ and is equal to $[W_\sigma,\DNC(\tau)]$. 
\end{lemm}
\begin{proof}
To prove the first statement, we will show that $(\tau + [\sigma,\tau])(f) = \tau(f) + [\sigma,\tau](f)$ vanishes along $\subbase$ for all $f \in C^\infty(\base)$ that vanish along $\subbase$. Indeed,
\[\tau(f) + [\sigma,\tau](f) = \tau(f) + \sigma(\tau(f)) - \tau(\sigma(f)) = \tau(f - \sigma(f)) + \sigma(f)\tau(f),\]
and both terms on the right-hand side vanish, because $\sigma$ vanishes along $\subbase$ and $f-\sigma(f)$ vanishes to order $2$ along $\subbase$ (see \eqref{eq: Euler-like}). 
The last statements follow by checking the identities on $\base \times \mathbb{R}^\times \subset \DNC(\base,\subbase)$, and using that $\base \times \mathbb{R}^\times \subset \DNC(\base,\subbase)$ is dense.
\end{proof}
\begin{lemm}\cite{Bischoff_2020}\label{lemm: properties of W_sigma 2}
For all $(y,\xi) \in \normal(\base,\subbase)$, the integral curve $\phi^t_W(y,\xi)$ of $W\coloneqq W_\sigma$ is defined for all $t \in \mathbb{R}$.
\end{lemm}
\begin{proof}
Observe first that $W\coloneqq W_\sigma$ is $\hat{t}$-related to $\tfrac{\partial}{\partial t} \in \mathfrak{\mathbb{R}}$. Therefore, at least for small $t>0$ (and $t<0$), $\phi^t_W$ will map into $\base \times \mathbb{R}_{>0}$ (resp. $\base \times \mathbb{R}_{<0}$). Now, by definition, $W$ is given by $\tfrac{1}{t}\sigma + \tfrac{\partial}{\partial t}$ on $\base \times \mathbb{R}^\times \subset \DNC(\base,\subbase)$, and we can describe its integral curve e.g. for $(x,s) \in \base \times \mathbb{R}_{>0} \subset \DNC(\base,\subbase)$ explicitly:
\[\phi^{t'}_W(x,s) = (\phi^{-\log(1-\frac{t'}{s})}_\sigma(x),s+t')\]
for all $t' \in \mathbb{R}$ such that $\tfrac{t'}{s} < 1$ (since $\sigma$ is complete). Indeed, 
\begin{align*}
    \left.\frac{d}{dt'}\right\vert_{t'=0}\phi^{-\log(1-\tfrac{t'}{s})}_\sigma(x) &= \left.\frac{d}{dt}\right\vert_{t'=0}\phi^{t'}_\sigma(x) \cdot (-\left.\frac{d}{dt'}\right\vert_{t'=0} \log(1-\frac{t'}{s})) \\
    &= \frac{1}{s} \cdot \left.\frac{d}{dt'}\right\vert_{t'=0}\phi^{t'}_\sigma(x) = \left.\frac{d}{dt'}\right\vert_{t'=0}\phi^t_{\frac{1}{s}\sigma}(x).
\end{align*}
But for $t'<0$ the condition $\tfrac{t'}{s} < 1$ is always satisfied (since $s>0$). Similarly, we can write down the flow for $(x,s) \in \base \times \mathbb{R}_{<0} \subset \DNC(\base,\subbase)$ and see that it exists for all $t'>0$, so, since the integral curve $\phi^t_W(y,\xi)$ lands for small $t>0$ (and $t<0$) into $\base \times \mathbb{R}^\times$, and then it exists for all time $t'>t$ (resp. $t'<t$), the result follows.
\end{proof}
\begin{proof}[Proof of Theorem \ref{theo: Euler-like vector fields and tubular neighbourhood embeddings}]
By Lemma \ref{lemm: properties of W_sigma 2}, the domain of $\phi_W^t$ ($t \in \mathbb{R}$), where $W\coloneqq W_\sigma$, is an open neighbourhood $U^t$ of $\normal(\base,\subbase) \subset \DNC(\base,\subbase)$. By the proof of the lemma, we see that, for all $t \in \mathbb{R}$, we obtain an embedding
\[\chi_t: \normal(\base,\subbase) \hookrightarrow U^t \xra{\phi^t_W} \DNC(\base,\subbase)_t \cong \base.\]
Moreover, $[W,\DNC(\sigma)]=0$ by Lemma \ref{lemm: properties of W_sigma 1}, so, in particular, the embedding $\chi\coloneqq\chi_1$ intertwines the vector fields $\mathcal{E}_{\normal(\base,\subbase)} = d_\normal \sigma = \DNC(\sigma)|_{\normal(\base,\subbase)}$ and $\sigma=\DNC(\sigma)|_{\DNC(\base,\subbase)_1}$. Since $W_\sigma$ equals $\tfrac{\partial}{\partial t}$ on $\subbase \times \mathbb{R} \subset \DNC(\base,\subbase)$, we see that $\chi$ is the identity map when restricted to $\subbase$. Now we use that $[W,\hat{\tau}] = \DNC(\tau + [\sigma.\tau])$ for all $\tau \in \mathfrak{X}(\base)$ by Lemma \ref{lemm: properties of W_sigma 1}, so that the map $\phi^t_W$ intertwines the vector field $\hat{\tau}$ with itself up to vector fields tangent to $\subbase \times \mathbb{R} \subset \DNC(\base,\subbase)$. Since $\hat{\tau}|_{\normal(\base,\subbase)}= \tau|_\subbase \mod T\subbase$, and $\hat{\tau}|_{\DNC(\base,\subbase)_t} = t\tau$, it follows that
\[d_\normal \chi_t: \normal(\base,\subbase) \cong \normal(\normal(\base,\subbase),\subbase) \ra  \normal(\DNC(\base,\subbase)_t,\subbase) \cong \normal(\base,\subbase)\]
is equal to the map $t \cdot \id_{\normal(\base,\subbase)}$, so $\chi=\chi_1$ is a tubular neighbourhood embedding for which $\mathcal{E}_{\normal(\base,\subbase)}$ is $\chi$-related to $\sigma$. 

To prove uniqueness, notice that if $\chi': \normal(\base,\subbase) \ra \base$ is another tubular neighbourhood embedding for which $\mathcal{E}_{\normal(\base,\subbase)}$ is $\chi'$-related to $\sigma$, then we obtain a map
\[\DNC(\chi'): \normal(\base,\subbase) \times \mathbb{R} \cong \DNC(\normal(\base,\subbase),\subbase) \ra \DNC(\base,\subbase)\]
(see Proposition \ref{prop: identify DNC(E,0_Y) with E times R}). Since the maps $\DNC(\chi')_t: \normal(\base,\subbase) \ra \DNC(\base,\subbase)_t \cong \base$ have normal derivative $t \cdot \id_{\normal(\base,\subbase)}$ (because $d_\normal(\chi') = \id_{\normal(\base,\subbase)}$ and Proposition \ref{prop: DNC of differentials}), and since $\mathcal{E}_{\normal(\base,\subbase)}$ is $\chi'$-related to $\sigma$, the above map is the flow of the vector field $W=W_\sigma$, so the maps are equal to the maps $\chi_t$ from above. Therefore, $\chi'=\chi$, which proves uniqueness.
%To prove uniqueness, suppose that $\chi: \base \ra \normal(\base,\subbase)$ is a tubular neighbourhood embedding for which the Euler vector field $\mathcal{E}=\mathcal{E}_{\normal(\base,\subbase)}$ is $\chi$-related to the Euler-like vector field $\sigma$. We claim that we can express $\chi$ in terms of the flow of $\sigma$. Of course, from this the result follows. To prove the claim, recall that $\mathcal{E}$ has time $s$-flow $\phi^s_{\mathcal{E}}$ given by $(y,\xi) \mapsto \mu_{e^{-s}}(y,\xi)$ (see \eqref{eq: multiplication DNC as R-action}); in particular, $\kappa_s \coloneqq \mu_s|_{\normal(\base,\subbase)}: \normal(\base,\subbase) \ra \normal(\base,\subbase)$ equals $\phi^{-\log s}_\tau$ for $s>0$. Now observe that 
%\[\lim_{s \ra 0} \kappa_s: \normal(\base,\subbase) \ra \normal(\base,\subbase)\]
%exists and equals the retraction of $\normal(\base,\subbase)$ onto $\subbase$.
\end{proof}
We are now ready to state the splitting theorem.
\begin{theo}\cite{meinrenken}[Splitting theorem for Lie algebroids]\label{theo: splitting theorem Lie algebroids}
Let $\algebroid$ be a Lie algebroid and let $\iota: \subbase \hookrightarrow \base$ be a transversal (that is, an embedded submanifold transverse to the anchor map). Let $\alpha \in \Gamma(\algebr)$ with $\alpha|_\subbase=0$, and the additional property that $\sigma \coloneqq \anchor(\alpha) \in \mathfrak{X}(\base)$ is Euler-like along $\subbase$ (which always exist). Then we obtain a Lie algebroid isomorphism 
\[\tau_{\normal(\base,\subbase)}^!\iota^!\algebr \xra{\sim} \algebr|_U\]
(here, $\tau_{\normal(\base,\subbase)}$ is the projection $\normal(\base,\subbase) \ra \subbase$) over the tubular neighbourhood embedding $\chi_\sigma: \normal(\base,\subbase) \xra{\sim} U \subset \base$ determined by $\sigma$.
\end{theo}
We will use two lemmas for the proof. 
\begin{lemm}\cite{meinrenken}\label{lemm: Euler like section of A}
Let $\algebroid$ be a Lie algebroid such that $\subbase \subset \base$ is a transversal. Then there is a section $\alpha \in \Gamma(\algebr)$ such that $\anchor(\alpha)$ is Euler-like along $\subbase$.
\end{lemm}
\begin{proof}
Since $\subbase$ is a transversal, we obtain a surjective morphism of vector bundles
\[\algebr|_\subbase \rightarrowdbl \normal(\base,\subbase).\]
Now let $\{y^i,x^j\}$ be adapted coordinates on an open subset $U$ of $\base$ (so $V \coloneqq U \cap \subbase = \{x^j=0\}$). Then, by surjectivity of the above map, we can find sections $\alpha_i \in \Gamma(A|_U)$ such that
\[\anchor(\alpha_i) = \frac{\partial}{\partial y^i} \mod T\subbase.\]
Now, along an open cover of adapted charts, patch the sections $\textstyle\sum_i y^i\alpha_i$ together via a partition of unity type argument. By Remark \ref{rema: local coordinates Euler-like vector field}, this section will map to an Euler-like vector field, except that it might not be complete. However, we can overcome this problem by picking open neighbourhoods $\widetilde{V} \subset \widetilde{U}$ of $\subbase$ for which we have a bump function supported on $\widetilde{U}$ and equal to $1$ on $\widetilde{V}$, and then multiply the vector field with this bump function.
\end{proof}
For the following lemma it is useful to remark the following first.
\begin{rema}\cite{meinrenken}\label{rema: normal(A,i^!A) = p^!i^!A}
Notice that, by checking the universal property of the pullback (see Proposition \ref{prop: universal property of pullbacks}), we have an isomorphism $\normal(\algebr,\iota^!\algebr) \cong \tau_{\normal(\base,\subbase)}^!\iota^!\algebr$ of Lie algebroids over $\normal(\base,\subbase)$ (note: we viewed $\iota^!\algebr$ as a subalgebroid of $\algebr$).
\end{rema}
\begin{lemm}\cite{meinrenken}\label{lemm: Euler vector field of normal(A,i^!A)}
There is a section $\epsilon \in \Gamma(\normal(\algebr,i^!\algebr))$ for which $\anchor_{\normal(\algebr,i^!\algebr)}(\epsilon) \in \mathfrak{X}(\normal(\base,\subbase))$ is the Euler vector field $\mathcal{E}_{\normal(\base,\subbase)}$. Moreover, there exists a section $\alpha \in \Gamma(\algebr,\iota^!\algebr)$ such that $d_\normal\alpha = \epsilon$; such a section will be called an \textit{Euler-like section} of $\algebr$.
\end{lemm}
\begin{proof}
The first statement follows by Remark \ref{rema: normal(A,i^!A) = p^!i^!A} and by observing that the Lie algebroid $T\normal(\base,\subbase) \times \iota^!\algebr$ carries the section $\mathcal{E}_{\normal(\base,\subbase)} \times 0$ which restricts to a section of $\tau_{\normal(\base,\subbase)}^!\iota^!\algebr$. Since the anchor map of a pullback Lie algebroid is given by the projection map from the pullback diagram, we see that $\anchor_{\tau_{\normal(\base,\subbase)}^!\iota^!\algebr}(\mathcal{E}_{\normal(\base,\subbase)} \times 0)=\mathcal{E}_{\normal(\base,\subbase)}$, so this proves the first statement. The second statement follows by noticing that any section $\alpha$ obtained from Lemma \ref{lemm: Euler like section of A} has the property that $d_\normal\alpha=\epsilon$.
\end{proof}
%Lastly, we observe the following fact.
%\begin{rema}\label{rema: linear vector fields on vector field}
%Let $E \ra M$ be a vector bundle. Then an $\mathbb{R}$-linear map 
%\[\psi: \Gamma(E) \ra \Gamma(E)\]
%for which $\tau(f) \coloneqq \psi(f\beta) - f\psi(\beta)$ defines a vector field $\tau \in \mathfrak{X}(M)$ determines a linear vector field $\widetilde{\tau}$ in $\mathfrak{X}(E)$ (linear meaning homogeneous of degree $1$).
%\end{rema}
The proof is based on \cite{Bischoff_2020} and \cite{meinrenken}.
\begin{proof}[Proof of Theorem \ref{theo: splitting theorem Lie algebroids}]
%Let $\alpha \in \Gamma(\algebr,\iota^!\algebr)$ be an Euler-like section. Then $\DNC(\alpha) \in \Gamma(\DNC(\algebr,\iota^!\algebr))$, and the derivation $\textnormal{ad } \DNC(\alpha) \coloneqq \bracket{\DNC(\alpha),\cdot}: \Gamma(\DNC(\algebr,\iota^!\algebr)) \ra \Gamma(\DNC(\algebr,\iota^!\algebr))$ defines a vector field $\sigma$ on $\algebr$ by 
%\[\sigma(f)(a) \coloneqq \textnormal{ad } \alpha(fa) - f\textnormal{ad } \alpha.\]
%Notice that this vector field restricts to a vector field on $\iota^!\algebr$, so we obtain a vector field $\DNC(\sigma) \in \mathfrak{X}(\DNC(\base,\subbase))$. Since
%\[\anchor \circ \textnormal{ad } \alpha = [\anchor(\alpha),\cdot] \circ \anchor,\]
%we see that the flow of $\sigma$ intertwines the anchor map, and 
The map $\hat{\pr}_1: \DNC(\base,\subbase) \ra \base$ (given by $\DNC(\base,\subbase) \xra{\DNC(\id_\base)} \DNC(\base,\base) \cong \base \times \mathbb{R} \xra{\pr_1} \base$) is transverse to the anchor map, because $\subbase$ is a transversal, so we have a Lie algebroid $\hat{\pr}_1^!\algebr \ra \DNC(\base,\subbase)$. Similarly as in the proof of Proposition \ref{prop: Euler-like and deformation to the normal cone}, if we take an Euler-like section $\alpha$ of $\algebr$ (and $\sigma \coloneqq \anchor_\algebr(\alpha)$), then 
\[\frac{1}{t}\alpha + \frac{\partial}{\partial t} \in \Gamma(\algebr \times T\mathbb{R}^\times)\]
extends to a smooth section $\widetilde{W}_\alpha$ of $\hat{\pr}_1^!\algebr$, and it is readily verified to be a lift of the vector field $W_\sigma \in \mathfrak{X}(\DNC(\base,\subbase))$. Now,
\[\ad \widetilde{W}_\alpha \coloneqq [\widetilde{W}_\alpha,\cdot]_{\hat{\pr}_1^!\algebr}\]
is a derivation of $\hat{\pr}_1^!\algebr$ by the Leibniz rule (see Definition \ref{defn: derivations of a vector bundle}), and so the ODE
\[\ad \widetilde{W}_\alpha(\beta) = \left.\frac{d}{dt}\right\vert_{t=0} (\phi^t_{\ad \widetilde{W}_\alpha})^*\beta\] 
gives rise to vector bundle automorphisms $\phi^t_{\ad \widetilde{W}_\alpha}$ (see also Remark \ref{rema: general linear algebroid is atiyah algebroid frame bundle}). Since $\ad \widetilde{W}_\alpha$ satisfies
\[\anchor_{\hat{\pr}_1^!\algebr}(\ad \widetilde{W}_\alpha(\beta)) = \anchor_{\hat{\pr}_1^!\algebr}([\widetilde{W}_\alpha,\beta]_{\hat{\pr}_1^!\algebr})=[\anchor_{\hat{\pr}_1^!\algebr}(\widetilde{W}_\alpha),\anchor_{\hat{\pr}_1^!\algebr}(\beta)] = \ad \anchor_{\hat{\pr}_1^!\algebr}(\widetilde{W}_\alpha) \circ \anchor_{\hat{\pr}_1^!\algebr}(\beta),\]
and it intertwines the Lie bracket by the Jacobi identity, the maps $\phi^t_{\ad \widetilde{W}_\alpha}$ are even Lie algebroid automorphisms, over $\phi_{W_\sigma}^t$. Again, as for $W_\sigma$, one can check that these automorphisms exist for all time $t$ (see Lemma \ref{lemm: properties of W_sigma 2}; alternatively, it follows directly by the existence of the flow of $W_\sigma$ for all time). Since $\widetilde{W}_\alpha$ lifts the vector field $W_\sigma$, and $\phi_{W_\sigma}^1$ is a tubular neighbourhood embedding onto the open set $U \subset \base$, we obtain an isomorphism of Lie algebroids 
\[\iota_0^!\hat{\pr}_1^!(\algebr) \xra{\sim} \iota_{1}^!\hat{\pr}_1^!(\algebr|_U)\] such that the post-composition with the canonical projection $\iota_{1}^!\hat{\pr}_1(\algebr|_U) \ra \hat{\pr}_1(\algebr|_U)$ equals the map $\iota_0^!\hat{\pr}_1(\algebr) \ra \hat{\pr}_1(\algebr)$ post-composed with $\phi^1_{\widetilde{W}_\alpha}$ (here, $\iota_{s}$ is the inclusion of $\DNC(\base,\subbase)_s$ into $\DNC(\base,\subbase)$). Since $\hat{\pr}_1 \circ \iota_0 = \iota \circ \tau_{\normal(\base,\subbase)}$ (here, $\iota: \subbase \hookrightarrow \base$ is the inclusion), and $\hat{\pr}_1 \circ \iota_1 = \id_\base$, the result follows by the relation $(f \circ g)^!B \cong g^!f^!B$ for Lie algebroids $B$ (see Proposition \ref{prop: (fg)!=g!f!}).
\end{proof}
Lastly, we state the local splitting theorem for Lie algebroids as a corollary of the above statement. 
\begin{coro}\cite{rui}\label{coro: local splitting theorem LA}
Let $\algebroid$ be a Lie algebroid and let $x \in \base$. If $\iota: \subbase \hookrightarrow \base$ is an embedded submanifold such that $x \in \subbase$, and 
\[T_x\subbase \oplus \anchor(x)\algebr_x = T_x\base,\]
then, by possibly shrinking $\subbase$, there is an open subset $U \ni x$ of $\base$ and an isomorphism of Lie algebroids
\[\iota^!\algebr \times TL_x \cong A|_U\]
for a connected embedded submanifold $L_x \subset \base$ with $T_xL_x = \anchor(x)\algebr_x$. Moreover, in $\algebr|_U$, $L_x = \orbit_x$, so, in $\algebr$, $\orbit_x$ is the unique maximal connected immersed submanifold $L \subset \base$ such that $\anchor(x)\algebr_x=T_xL$ for all $x \in L$.
\end{coro}
%\begin{proof}
%We only have to show that $L_x = \orbit_x$ in $\algebr|_U \cong \iota^!\algebr \times TL_x$. This goes as follows: if $y,z \in L_x$ $a: [0,1] \ra \algebr$ be an $\algebr$-path above $\gamma: [0,1] \ra \base$ such that $\gamma(0)=x$ and $\gamma(1)=y$. Since $T_{\gamma(t)}L_{\gamma(t)} = \anchor(y)\algebr_y$ we see that
%\[\frac{d\gamma}{dt}(t) = \anchor(a(t)) \in \anchor(y)\algebr_{\gamma(t)} = T_{\gamma(t)}L_{\gamma(t)},\]
%so $\gamma$ is tangent to $L_{\gamma(t)}$ for all $t \in [0,1]$.
%\end{proof}
\addtocontents{toc}{\protect\thispagestyle{myplain}}\newpage

\section{Lie groupoid and Lie algebroid blow-up}\label{sec: Blow-up of a pair of Lie groupoids}
We will see in this chapter that to each pair of smooth manifolds $(\base,\subbase)$ we can associate a blow-up manifold $\blup(\base,\subbase)$. We will do this in various ways, and we will carefully check the necessary details. Then, we show that to a pair of groupoids $(\groupoid,\subgroupoid)$ (see Definition \ref{defn: pair of Lie subgroupoids}) we can (naturally) associate a blow-up groupoid $\blup_{\source,\target}(\group,\subgroup) \rra \blup(\base,\subbase)$, and to a pair of algebroids $(\algebroid,\subalgebroid)$ (see Definition \ref{defn: pairs of Lie algebroids}), we can associate a blow-up algebroid $\blup_\pi(\algebr,\subalgebr) \ra \blup(\base,\subbase)$ (the notation will be explained later). 

Blow-ups in algebraic geometry are a well-studied topic, and one way to construct blow-ups of smooth manifolds is to mimic the usual construction for varieties and/or schemes.

\subsection{Blow-ups in algebraic geometry}\label{sec: Blow-up in the algebraic setting}
Before we go to the blow-up construction for manifolds, we will digress to explain the blow-up construction for varieties and schemes. If you are not familiar with the theory of algebraic varieties, replace $k$ everywhere with $\mathbb{C}$, and replace ``variety'' everywhere with the zero locus 
\[Z(f) \coloneqq \{z \in \mathbb{C}^n \mid f(z)=0\}\] 
of polynomials $f \in \mathbb{C}[x_1,\dots,x_n]$ (it is equipped with the Zariski topology; the inherited topology from $\mathbb{C}^n \supset Z(f)$ equipped with the coarsest topology for which all $Z(g) \subset \mathbb{C}^n$, $g \in \mathbb{C}[x_1,\dots,x_n]$ are closed). If you are not familiar with the theory of schemes, feel free to skip the part about the blow-up of schemes. The constructions are based on the blow-up construction explained in \cite{hartshorne}. This is the standard way to approach blow-ups, and it does not use deformation to the normal cones (which is also a well-known construction in algebraic geometry). Moreover, deformation to the normal cones are often constructed out of blow-ups, and we will see, in our differential geometric context, that we can indeed reconstruct deformation to the normal cones from the blow-up construction (see Proposition \ref{prop: DNC out of blup}). We will, however, approach blow-ups by using the theory of deformation to the normal cones instead (and one can also do so in the algebraic context). Still, we will also explain, and compare, the more ``direct'' blow-up construction that is similar to the one discussed in this section (see Remark \ref{rema: other blow-up constructions}). 

\subsubsection{Blow-up of a point in an algebraic variety}
First, we define the simplest, and perhaps most important (we will motivate this later in the manifold context), example of a blow-up in the variety context: the blow-up of a point in a variety. Let $k$ be an algebraically closed field and denote by $Z$ the functor from subsets of $k[x_1,\dots,x_n]$ (ordered by inclusion) to subsets of $k^n$ that associates to each $S \subset k[x_1,\dots,x_n]$ the zero locus of $S$ in $k^n$, i.e.
\[Z(S) = \{x \in k^n \mid f(x)=0 \textnormal{ for all } f \in S\}.\]
Moreover, denote affine $n$-space by $\mathbb{A}^n$, i.e. $\mathbb{A}^n = k^n$ equipped with the Zariski topology (i.e. the closed sets are given by the 
zero loci as above), and denote projective $n$-space $(\mathbb{A}^n \setminus \{0\})/k^\times$, where $k^\times$ acts diagonally on $\mathbb{A}^n \setminus \{0\}$, by $\mathbb{P}^n$. 
\begin{defn}\cite{hartshorne}\label{defn: blow-up 0 in A^n}
The \textit{blow-up of a point} (which we take to be the origin) \textit{in $\mathbb{A}^n$} is
\[\blup_{\textnormal{var}}(\mathbb{A}^n,0) \coloneqq \{\left((x_1,\dots,x_n),[y_1:\dots:y_n]\right) \mid x_iy_j = y_ix_j \textnormal{ for all } 1 \le i,j \le n\} \subset \mathbb{A}^n \times \mathbb{P}^{n-1},\]
where $\mathbb{A}^n \times \mathbb{P}^{n-1}$ is the product of $\mathbb{A}^n$ with $\mathbb{P}^{n-1}$ in the category of varieties.
\end{defn}
Observe that $\blup_{\textnormal{var}}(\mathbb{A}^n,0)$ is a closed subset of $\mathbb{A}^n \times \mathbb{P}^{n-1}$ and comes with a morphism of varieties
\[\bldown_{\textnormal{var}} \coloneqq \pr_1|_{\blup_{\textnormal{var}}(\mathbb{A}^n,0)}: \blup_{\textnormal{var}}(\mathbb{A}^n,0) \mapsto \mathbb{A}^n,\]
called the \textit{blow-down map}. It is readily verified that it restricts to an isomorphism
\[\blup_{\textnormal{var}}(\mathbb{A}^n,0) \setminus \bldown_{\textnormal{var}}^{-1}(0) \cong \mathbb{A}^n \setminus \{0\};\] 
the inverse is given by sending $x=(x_1,\dots,x_n) \in \mathbb{A}^n \setminus \{0\}$ to $(x,[x_1 : \cdots : x_n])$. Moreover, notice that $\bldown_{\textnormal{var}}^{-1}(0)=\{0\} \times \mathbb{P}^{n-1} \subset \blup_{\textnormal{var}}(\mathbb{A}^n,0)$ which is often called the \textit{exceptional divisor}; one can check that $\blup_{\textnormal{var}}(\mathbb{A}^n,0)$ is irreducible and of dimension $n$, so the exceptional divisor is indeed a divisor. Sometimes it is called the \textit{singular locus}, and we will adopt this terminology. We see that in $\blup_{\textnormal{var}}(\mathbb{A}^n,0)$ one keeps track, not only of the element of $\mathbb{A}^n$, but also of a line through the origin on which this element lies, so, effectively, we ``replaced the origin with the lines through the origin''. 
\begin{defn}\cite{hartshorne}\label{defn: blow-up 0 in affine variety}
If $Y \subset \mathbb{A}^n$ is closed (and $0 \in Y$), then the \textit{blow-up of a point} ($=0$) \textit{in $Y$} is
\[\blup_{\textnormal{var}}(Y,0) \coloneqq \overline{\bldown_{\textnormal{var}}^{-1}(Y\setminus\{0\})}.\]
By blowing up locally, we can extend the construction to all algebraic varieties. 
\end{defn}
Notice that for $Y \subset \mathbb{A}^n$ closed ($0 \in Y$), the blow-down map restricts to a blow-down map $p_{\textnormal{var}}: \blup_{\textnormal{var}}(Y,0) \ra Y$, which further restricts to an isomorphism $\blup_{\textnormal{var}}(Y,0) \setminus p_{\textnormal{var}}^{-1}(0) \cong Y \setminus \{0\}$. For more general varieties $X$, this map is defined locally on an open affine subset $U$ around the point $x$ we blew up, and we can extend it as the identity map on $X \setminus U \cong \blup_{\textnormal{var}}(X,x) \setminus \blup_{\textnormal{var}}(U,x)$, and then the resulting blow-down map $p_{\textnormal{var}}: \blup_{\textnormal{var}}(X,x) \ra X$ again restricts to an isomorphism $\blup_{\textnormal{var}}(X,x) \setminus p_{\textnormal{var}}^{-1}(x) \cong X \setminus \{x\}$.

%If we replace $Y$ with a variety $X$, then we can define the map $p_{\textnormal{var}}: \blup_{\textnormal{var}}(X,0) \ra X$ first on an open affine set $U$ around $0 \in X$ (we called $0$ the element in $X$ which corresponds to $0$ in $\mathbb{A}^n \supset Y \cong U$, where $Y$ is closed), and then extend the map $p_{\textnormal{var}}: \blup_{\textnormal{var}}(U,0) \ra U$ as the identity map
There is one obvious application of blow-ups: we can ``resolve singularities''. As a simple example, consider $X \coloneqq Z(f) \subset \mathbb{A}^2$, where $f \in k[x,y]$, that has self-intersections in $0 \in \mathbb{A}^2$. %If we write $f=\textstyle\sum_{k \ge 0} f_k$ (finitely many $f_k$ are non-zero), where the $f_k$ are homogeneous polynomials of degree $k$, then one can check that (1) the $f_k$ split into linear polynomials (that is, polynomials that are homogeneous of degree $1$), and (2) for $\ell$ the lowest $k \ge 0$ such that $f_k \neq 0$, the tangent lines of $f$ at $0$ are precisely the zero-loci of these linear polynomials corresponding to $f_\ell$. 
Now, if the tangent lines of $f$ at $0$ are all different, then we can blow up the point $0$ in $X$, and the blow-up $\widetilde{X} \coloneqq \blup_{\textnormal{var}}(X,0)$ will have no self-intersections anymore. This is because one can show that $\widetilde{X} \cap p_{\textnormal{var}}^{-1}(0)$ consists precisely of those points $(0,y) \in \mathbb{A}^2 \times \mathbb{P}^1$ for which the line that $y$ corresponds to in $\mathbb{A}^2$ equals a tangent line. 

To explain the blow-up construction in this algebraic setting more generally, we move to the theory of schemes.
\subsubsection{Blow-up of schemes}
We will now explain the more general construction of blow-ups of (Noetherian) schemes. Feel free to skip this part if you are not familiar with the theory of schemes.

Recall that a graded ring $S$ is a ring $\textstyle\bigoplus_{j\ge0}S_j$, such that whenever $i,j \ge 0$, then $S_i\cdot S_j \in S_{i+j}$. Elements of $S_j$ are called \textit{homogeneous elements} and ideals $\mathfrak{a} \subset S$ satisfying $\mathfrak{a} = \textstyle\bigoplus_{j\ge0} (\mathfrak{a} \cap S_j)$ are called \textit{homogeneous ideals}. It is readily verified that the sum, product, intersection and radical of such homogeneous ideals is again homogeneous and that a homogeneous ideal $\primeid$ is prime if and only if whenever $f,g$ are homogeneous with $fg \in \primeid$, then $f \in \primeid$ or $g \in \primeid$. We denote by $\Proj S$ the set of homogeneous prime ideals of $S$, and we denote the ideal $\textstyle\bigoplus_{j > 0} S_j$ of $S$ by $S_+$.
\begin{rema}\cite{hartshorne}\label{rema: proj of a graded ring}
Let $S$ be a graded ring. Then $\Proj S = \{\textnormal{homogeneous prime ideals of } S\}$ is equipped with the topology whose closed sets are of the form 
\[V(\mathfrak{a}) \coloneqq \{\primeid \in \Proj S \mid \mathfrak{a} \subset \primeid\},\]
where $\mathfrak{a}$ is a homogeneous ideal of $S$. It is readily verified that for homogeneous ideals $\mathfrak{a},\mathfrak{b}$ and $\{\mathfrak{a}_i\}_{i \in I}$ of $S$ we have $V(\mathfrak{a}) \cup V(\mathfrak{b}) = V(\mathfrak{a} \cdot \mathfrak{b})$ and $\textstyle\bigcap_{i \in I}V(\mathfrak{a}_i) = V(\textstyle\sum_{i \in I} \mathfrak{a}_i)$, so the above sets really do form the closed sets of a topology on $\Proj S$. 
\end{rema}
We will now define a sheaf of rings on $\Proj S$. Let $\primeid \in \Proj S$. Then
\[S_{(\primeid)}=T^{-1}S\]
is the localisation of $S$ with respect to the multiplicative closed subset 
\[T \coloneqq \{a \in S \setminus \primeid \mid a \textnormal{ is homogeneous}\}\]
(i.e. we invert the elements in $T$).
Let $U \subset \Proj S$ be open and consider a map $s: U \ra \textstyle\prod_{\primeid \in U} S_{(\primeid)}$. Then we define $\OX_{\Proj S}(U)$ as follows: $s \in \OX_{\Proj S}(U)$ if and only if for all $\primeid \in U$ there is an open subset $V \ni \primeid$ of $U$ together with homogeneous elements $a,f \in S$ such that whenever $\mathfrak{q} \in V$ we have $f \not\in \mathfrak{q}$ and $s(\mathfrak{q})=\tfrac{a}{f} \in S_{(\mathfrak{q})}$.
Similar to the ``Spec''-construction, this defines a sheaf of rings on $\Proj S$ with stalks $\OX_{\primeid}$ canonically isomorphic to $S_{(\primeid)}$. 
\begin{term}\label{term: open set and localisation with respect to f in S_+}
Let $S$ be a graded ring and let $f \in S_+ = \textstyle\bigoplus_{j>0}S_j$ be a homogeneous element. Then we write 
\[D_+f \coloneqq \Proj S \setminus V(f) = \{\primeid \in \Proj S \mid f \not\in \primeid\},\]
and we write $S_{(f)}$ for the subring of homogeneous elements of degree $0$ in $S_f$ (where $S_f$ is the localisation of $S$ with respect to the multiplicative closed subset $\{f^n\}_{n \ge 0}$ of $S$).
\end{term}
We state the following proposition without proof.
\begin{prop}\cite{hartshorne}\label{prop: Proj S is a scheme}
Let $S$ be a graded ring and let $f \in S_+ = \textstyle\bigoplus_{j>0}S_j$ be a homogeneous element. Then there is a canonical isomorphism of ringed spaces
\[(D_+f,\OX_{\Proj S}|_{D_+f}) \ra \Spec S_{(f)}.\]
In particular, $\Proj S$ is a scheme.
\end{prop}
\begin{defn}\cite{hartshorne}\label{defn: projective scheme}
Let $R$ be a ring. Then \textit{projective $n$-space over $R$} is the scheme $\mathbb{P}^n_R \coloneqq \Proj R[x_0,\dots,x_n]$.
\end{defn}
\begin{rema}\cite{hartshorne}\label{rema: projective scheme vs projective variety}
Let $X$ be a variety over an algebraically closed field $k$. Let $Y \subset X$ be closed and consider the set
\[t(Y) = \{Z \subset Y \mid Z \textnormal{ is non-empty and irreducible}\}.\]
Then $t(Y_1 \cup Y_2)= t(Y_1) \cup t(Y_2)$ and $t(\textstyle\bigcap_{i \in I} Y_i) = \textstyle\bigcap_{i \in I} t(Y_i)$, so the subsets $t(Y)$ of $t(X)$, where $Y \subset X$ is closed, are the closed sets of a topology on $t(X)$. Define a map
\[F: X \ra t(X) \textnormal{ by } x \mapsto \overline{\{x\}}.\]
Then one can show that $(t(X),F_*(\OX_X))$ is a scheme over $k$ (which means it is equipped with a morphism of schemes to $\Spec k$). Moreover, if $f: X_1 \ra X_2$ is a continuous map, then 
\[t(f): t(X_1) \ra t(X_2) \textnormal{ given by } Y \mapsto \overline{Y}\]
is also continuous. A careful assessment shows that this defines an equivalence of categories from varieties over $k$ to reduced schemes of finite type over $k$ (see e.g. \cite{haiman}). 

The point we want to make, is that the scheme $\mathbb{P}^n_k$ corresponds to the scheme $t(\mathbb{P}^n)$, where $\mathbb{P}^n$ is projective space (over $k$) as a variety. Indeed, let $X$ consist of the closed points of $\mathbb{P}^n_k$. Then $X=\mathbb{P}^n$. But $t(X)$ and $\mathbb{P}^n_k$ have the same underlying topological space, so since $t(X)$ is reduced (and so is $\mathbb{P}^n_k$) we see that $t(X)$ and $\mathbb{P}^n_k$ are equal as schemes (since there is a unique reduced scheme associated to a closed subspace).
\end{rema}
\begin{term}\label{term: noetherian schemes and quasi-coherent sheaves}
For the rest of this section, we will assume that $X$ is a Noetherian scheme (i.e. $X$ can be covered by finitely many affine schemes $\Spec A$ with $A$ Noetherian). Moreover, $\mathcal{F}=\textstyle\bigoplus_{j\ge0} \mathcal{F}_j$ is a quasi-coherent sheaf of graded $\OX_X$-algebras with $\mathcal{F}_0=\OX_X$ and $\mathcal{F}_1$ coherent. Lastly, $\mathcal{F}$ is locally generated by $\mathcal{F}_1$ as an $\OX_X$-algebra.
\end{term}
Let $\Spec R \cong U \subset X$ be open and affine. Then $\mathcal{F}(U) = \Gamma(U,\mathcal{F}|_U)$ is a graded $R$-algebra, so we can consider $\Proj \mathcal{F}(U)$. Notice that the map $R \ra \mathcal{F}(U)$ induces a morphism of schemes
\[\pi_U: \Proj \mathcal{F}(U) \ra U \cong \Spec R.\]
Since $\mathcal{F}$ is quasi-coherent, if $f \in R$, then $\pi_U^{-1}(U_f) \cong \Proj \mathcal{F}(U_f)$ (where $U_f \cong \Spec R_f$, and $R_f$ is the localisation of the multiplicative closed set $\{f^n\}_{n \ge 0}$). So, if $V$ is another affine open set of $X$, then $\pi_U^{-1}(U \cap V) \cong \pi_V^{-1}(U \cap V)$. In particular, we see that we can glue the schemes $\Proj \mathcal{F}(U)$ together to obtain a scheme $\PROJ \mathcal{F}$ together with a morphism of schemes
\[\pi_{\mathcal{F}}: \PROJ \mathcal{F} \ra X\]
that has the property that whenever $U \subset X$ is open and affine, then $\pi^{-1}(U) \cong \Proj \mathcal{F}(U)$.
\begin{defn}\cite{hartshorne}\label{defn: generalised Proj}
Let $Y \subset X$ be a closed subscheme corresponding to a coherent sheaf of ideals $\mathcal{Y}$ on $X$; set $\mathcal{F}=\textstyle\bigoplus_{j\ge0}\mathcal{Y}^j$ (with $\mathcal{Y}^0=\OX_X$). Then
\[\blup_{\textnormal{sch}}(X,Y) = \PROJ \mathcal{F}\]
is called the blow-up of $X$ along $Y$. The blow-down map is the morphism $\pi_{\mathcal{F}}$.
\end{defn}
One can mimic the construction so that it becomes a blow-up construction in the differential geometric setting. However, we will see that we can go a different route via the deformation to the normal cone construction already discussed. As we will see, this approach has numerous advantages, and we will explain this construction first. 

\subsection{Blow-ups of manifolds}\label{sec: Blow-up of a pair of manifolds}
In the former section we have seen that, roughly speaking, the blow-up of a point $x$ (at least in an algebraic variety) effectively replaces the point with the lines through it, i.e. the projectivisation of the tangent space at $x$. If we start with a pair of smooth manifolds $(\base,\subbase)$ (from now on it will be important that $\subbase \subset \base$ is closed), then, again roughly speaking, the blow-up of $\subbase$ in $\base$ replaces $\subbase$ in $\base$ with the \textit{projectivisation of the normal bundle $\normal(\base,\subbase)$}; we will explain this terminology later in this section. As we will see in this section, and what might already be (intuitively) clear from the former observation, is that the blow-up of $\subbase$ in $\base$ can be realised as the quotient of the canonical $\mathbb{R}^\times$-action on $\DNC(\base,\subbase) \setminus (\subbase \times \mathbb{R})$. Because we will have to work with quotient manifolds (by Lie group actions) a lot, and we will have to show that these quotient manifolds are again manifolds, it is useful to recall some facts in this direction. The statements are stronger than presented here (as can be seen from the proofs), but we state them in the form we need (based on \cite{kaliszewski2015properness}). %Since we will also have to work with non-Hausdorff manifolds, we start with the following observation.
\begin{lemm}\label{lemm: proper iff G_C,D}
Let a Lie group $G$ act on a possibly non-Hausdorff manifold $P$. Then the action is proper if and only if for all compact sets $C,D \subset P$ the set
\begin{equation}\label{eq: G_C,D}
    G_{C,D} \coloneqq \{g \in G \mid Cg \cap D \neq \emptyset\}
\end{equation}
is compact.
\end{lemm}
\begin{proof}
The following standard proof does not use Hausdorffness: suppose that the action is proper. Then this means that $f: P \times G \ra P \times P$ given by $(p,g) \mapsto (p,pg)$ is a proper map. Let $C,D \subset P$ be compact sets. Then $C \times D$ is compact in $P \times P$. So, $G_{C,D} = \pr_2(f^{-1}(C \times D))$ is compact. 

Conversely, assume $G_{C,D}$ is compact for all compact subsets $C,D \subset P$. If $\widetilde{C} \subset P \times P$ is compact, then $C \coloneqq \pr_1(\widetilde{C})$ and $D \coloneqq \pr_2(\widetilde{C})$ are compact, and $\widetilde{C} \subset C \times D$. Since $f^{-1}(\widetilde{C})$ is closed and contained in $f^{-1}(C \times D)$, we only have to show that $f^{-1}(C \times D)$ is compact. But it is closed and contained in $C \times G_{C,D}$, which is compact. It follows that $f$ is proper.
\end{proof}
\begin{coro}\label{coro: proper iff G_K,K}
Let a Lie group $G$ act on a (possibly non-Hausdorff) manifold $P$. Then the action is proper if and only if for all compact sets $K \subset P$ the set $G_{K,K}$ is compact.
\end{coro}
\begin{proof}
If the action is proper, then the statement follows from Lemma \ref{lemm: proper iff G_C,D}. The converse statement follows from the proof of Lemma \ref{lemm: proper iff G_C,D}, but where we replace both $C \coloneqq \pr_1(\widetilde{C})$ and $D \coloneqq \pr_2(\widetilde{C})$ with $K \coloneqq C \cup D$.
\end{proof}
\begin{coro}\cite{kaliszewski2015properness}\label{coro: proper iff relatively compact}
Let a Lie group $G$ act on a (possibly non-Hausdorff) manifold $P$. Then the action is proper if and only if for all $p,q \in P$ there are compact neighbourhoods $C \ni p$ $D \ni q$ (i.e. $C$ and $D$ are compact, and there are open sets $O_C \ni p$ and $O_D \ni q$ such that $O_C \subset C$ and $O_D \subset D$) such that $G_{C,D}$ (see \eqref{eq: G_C,D}) is compact.
\end{coro}
\begin{proof}
If the action is proper, then the statement follows from Lemma \ref{lemm: proper iff G_C,D}. To prove the converse statement, we will show that $G_{K,K}$ is compact for all compact subset $K \subset P$ (see Corollary \ref{coro: proper iff G_K,K}). By second countability, we may prove that every sequence $(g_n)$ in $G_{K,K}$ has a convergent subsequence. By definition of $G_{K,K}$, we can find $p_n \in K$ such that $\{p_ng_n\} \subset K$. Since $K$ is assumed to be compact, we can pass to subsequences and assume that $p_n \ra p$ and $p_ng_n \ra q$ for some $p,q \in K$. We now pick compact neighbourhoods $C$ and $D$ of $p$ and $q$, respectively. The sequences $(p_n)$ and $(p_ng_n)$ eventually land in $C$ and $D$, respectively, so $(g_n)$ eventually lands in $G_{C,D}$, which is compact. This shows that $(g_n)$ has a convergent subsequence, as claimed.
\end{proof}
Notice, especially in the proof of Corollary \ref{coro: proper iff relatively compact} that we did not use uniqueness of limits, so that it does indeed work for non-Hausdorff manifolds $P$. If $P$ is Hausdorff, then we can give another characterisation of properness in terms of local properness. 
\begin{rema}\cite{kaliszewski2015properness}\label{rema: locally proper}
Recall that a continuous map $f: \base \ra \basetwo$ between topological spaces is called locally proper if there exists an open cover $\{U_i\}$ of $\base$ such that $f|_{U_i}: U_i \ra \basetwo$ is proper (compare with Lemma \ref{lemm: proper maps} below). In case of the action by a Lie group, the action is called locally proper if there exists an open cover $\{U_i\}$ of $P$ by $G$-invariant open subsets such that the restrictions of $f: P \times G \ra P \times P$, given by $(p,g) \mapsto (p,pg)$, to subsets of the form $U_i \times G$ are proper.
\end{rema}
The correct analogue for locally proper actions of the equivalent characterisations given for proper actions is the following (compare with Corollary \ref{coro: proper iff relatively compact}).
\begin{lemm}\cite{kaliszewski2015properness}\label{lemm: locally proper iff relatively compact}
Let a Lie group $G$ act on a (possibly non-Hausdorff) manifold $P$. Then the action is locally proper if and only if for all $p \in P$ there is a compact neighbourhood $K \ni p$ such that $G_{K,K}$ is compact. 
\end{lemm}
\begin{proof}
If the action is locally proper, then we can find an open cover $\{U_i\}$ by $G$-invariant open subsets such that the action is proper when restricted to one of the $U_i$. The result follows by Corollary \ref{coro: proper iff relatively compact} (by shrinking $C$ such that $C \subset D$, and then taking $K \coloneqq D$). 

Conversely, let $p \in P$ and pick a compact neighbourhood $K \ni p$ such that $G_{K,K}$ is compact. By possibly shrinking $K$, we can find an open set $p \in O \subset K$, such that $\overline{O} = K$ is compact, and we will show that the action restricted to $U \coloneqq O \cdot G$ is proper (by possibly shrinking $K$). We will use Corollary \ref{coro: proper iff relatively compact}: let $u,v \in U$, i.e. $u=xg$ and $v=yh$ for some $x,y \in O$ and $g,h \in G$. Then we can find compact neighbourhoods $C$ and $D$ of $u$ and $v$, respectively, such that $G_{C,D}$ is compact, and we may assume (by shrinking $K$,$C$ and $D$) that $C=K \cdot g$ and $D = K \cdot h$. But then
\begin{align*}
    g^{-1}G_{K,K}h &= \{g' \in G \mid Kgg'h^{-1} \cap K \neq \emptyset\} \\
    &= \{g' \in G \mid Kgg' \cap Kh \neq \emptyset\} = G_{C,D}
\end{align*}
is a compact set, so $G_{K,K}$ is too. This proves the statement.
\end{proof}
We now give another characterisation of proper actions in case the manifold $P$ is Hausdorff.
\begin{prop}\cite{kaliszewski2015properness}\label{prop: proper iff locally proper + Hausdorff}
Let a Lie group $G$ act on a Hausdorff manifold $P$. Then the action is proper if and only if it is locally proper and $P/G$ is Hausdorff.  
\end{prop}
\begin{proof}
First assume the action is proper. Obviously the action is then locally proper as well. By second-countability, we may prove, for Hausdorffness of $P/G$, that every sequence $([p_n])=([p_n])_{n \in\mathbb{N}}$ in $P/G$ has a unique limit. So, assume $[p_n] \ra [p]$ and $[p_n] \ra [q]$. Let $\pi: P \ra P/G$ be the quotient map, and let $U \ni p$ be an open subset of $P$. Then $\pi(U) \subset P/G$ is open and contains $[p]$, so the sequence $([p_n])$ eventually lies in $\pi(U)$, i.e. we can find $g_n \in G$ such that $(p_ng_n)$ eventually lies in $U$. If $U$ is relatively compact (i.e. $\overline{U}$ is compact), we see that $(p_ng_n)$ has a subsequence converging to some element in $P$, which has to be of the form $pg$ (because $[p_ng_n] \ra [p]$). So, by relabelling and repeating the argument, we find sequences $(p_n)$ and $(p_ng_n)$ such that $p_n \ra p$ and $p_ng_n \ra q$. We now choose compact neighbourhoods $C,D \subset P$ around $p$ and $q$, respectively, so, in particular, $G_{C,D}$ is compact (see Corollary \ref{coro: proper iff relatively compact}). Then the sequences $(p_n)$ and $(p_ng_n)$ eventually lie in the open sets $O_C \subset C$ and $O_D \subset D$ of $P$, respectively, so the sequence $(g_n)$ lies eventually in $G_{C,D}$. Therefore, we can pass to a subsequence and assume $(g_n)$ converges to an element $g \in G$. But then $p_ng_n \ra pg$ and $p_ng_n \ra q$, so $q=pg$ by Hausdorffness of $P$. This shows that $[p]=[q]$, so $P/G$ is Hausdorff, as claimed. 

Conversely, suppose the action is locally proper and let $P/G$ be Hausdorff. Let $p,q \in P$. By Corollary \ref{coro: proper iff relatively compact}, we want to show that we can find compact neighbourhoods $C$ and $D$ of $p$ and $q$, respectively, such that $G_{C,D}$ is compact. By Lemma \ref{lemm: locally proper iff relatively compact}, we can find a compact neighbourhood $C \ni p$ such that $G_{C,C}$ is compact. We will show that for any compact neighbourhood $D \ni q$, $G_{C,D}$ is compact as well. By second countability, it suffices to show that any sequence $(g_n)$ in $G_{C,D}$ has a convergent subsequence. By definition of $G_{C,D}$, we can find $c_n \in C$ such that $\{c_ng_n\} \subset D$.  Since $C$ and $D$ are compact, we can pass to subsequences and assume that $c_n \ra c$ and $c_ng_n \ra d$ for some $c \in C$ and $d \in D$. Since $[c_n] \ra [c]$ and $[c_n] = [c_ng_n] \ra [d]$, and $P/G$ is Hausdorff, we see that $[c]=[d]$, i.e. there is a $g \in G$ such that $cg=d$. Therefore, $c_ng_n \ra cg$, and so $c_ng_ng^{-1} \ra c$. In particular, the sequence $(c_ng_ng^{-1})$ eventually lies in $C$, so the sequence $(g_ng^{-1})$ eventually lies in $G_{C,C}$. Therefore, it has a convergent subsequence converging to, say, $h \in G$, and so $(g_n)$ has a convergent subsequence to $hg$. This shows that $G_{C,D}$ is compact, as required.
\end{proof}

Using the above proposition, it is not hard to see that we obtain the following result.
\begin{prop}\label{prop: locally proper}
Let a Lie group $G$ act on a possibly non-Hausdorff smooth manifold $P$. If the action of $G$ on $P$ is free and locally proper, then $P/G$ is a smooth possibly non-Hausdorff manifold and $\pi: P \ra P/G$ is a (possibly non-Hausdorff) principal $G$-bundle. Moreover, $\pi^*: C^\infty(P/G) \ra C^\infty(P)$ maps isomorphically onto $C^\infty_G(P) \coloneqq \{f \in C^\infty(P) \mid f(pg)=f(p) \textnormal{ for all } g \in G\}$.
\end{prop}

The reason the above result holds is simply by checking that a usual proof of the above result, with ``locally proper'' replaced with ``proper'', uses properness only to show that $P/G$ is Hausdorff. For the rest of the proof, however, it suffices to use local properness. 

Before we move on to blow-ups, we give some useful facts about more general proper maps.
\begin{lemm}\cite{schultz}\label{lemm: proper maps}
Let $f: \base \ra \basetwo$ $g: \basetwo \ra K$ be continuous maps between locally compact topological spaces. Then
\begin{enumerate}
    \item $f$ is proper if and only if there is an open cover $\{W_i\}_{i \in I}$ of $\basetwo$ such that $f|_{f^{-1}(W_i)}$ is proper;
    \item if $\base$ and $\basetwo$ are Hausdorff, and $\basetwo$ is locally compact, then $f$ is proper if and only if $f$ is closed and all fibers $f^{-1}(x)$ ($x \in \basetwo$) are compact;
    \item if $f$ and $g$ are proper, then so is $g \circ f$. Conversely, if $g \circ f$ is proper, then $f$ is proper, and if $f$ is, in addition, surjective, then $g$ is proper as well.
\end{enumerate}\vspace{-\baselineskip}
\end{lemm}
To explain the terminology ``projectivisation of the normal bundle'' (see the beginning of this section), observe the following:
\begin{lemm}\label{lemm: free and proper scalar multiplication action on E}
Let $E \ra M$ be a vector bundle. Then the $\mathbb{R}^\times$-action on $E \setminus 0_M$ given by scalar multiplication is free and (locally) proper (if $E$ is possibly non-Hausdorff).
\end{lemm}
\begin{proof}
The freeness of the action is clear, since we removed $0_M$ from $E$, and we can check it fiberwise: $\lambda\xi=\xi$, where $\lambda \in \mathbb{R}^\times$ and $\xi \in E_x \setminus \{0\}$ ($x \in M$), then $\lambda=1$.

We will check now that the action is proper in the Hausdorff case. For the non-Hausdorff case (and proving local properness instead of properness), notice that we may replace $E$ with $U \times \mathbb{R}^k$, where $U$ is the domain of a chart of $M$. Since such an open subset is Hausdorff, it suffices to prove that the action is proper in case $E$ (and hence $M$) is Hausdorff.
To check properness in the Hausdorff setting, we will prove, whenever $C,D \subset E$ are compact subsets, that 
\[(\mathbb{R}^\times)_{C,D} = \{\lambda \in \mathbb{R}^\times \mid \lambda C \cap D \neq 0\}\]
is compact. We will do this by showing that every sequence $(\lambda_n)$ in $(\mathbb{R}^\times)_{C,D}$ admits a convergent subsequence. Indeed, for such a sequence, choose a sequence $(c_n)$ in $C$ such that $\{\lambda_nc_n\} \subset D$.
By passing to a subsequence, we may assume that $c_n \ra c$ and $\lambda_nc_n \ra d$ for some $c \in C$ and $d \in D$. Observe that $(\lambda_n)$ is bounded; otherwise $\lambda_n^{-1} \ra 0$, and so $c_n = \lambda_n^{-1} \cdot (\lambda_nc_n) \ra 0 \cdot d \in 0_M$, but this is not the case by assumption. Therefore, we can pass to a convergent subsequence (in $\mathbb{R}$). If $\lambda_n \ra 0$, then $\lambda_nc_n \ra 0 \cdot c \in 0_M$, which is not the case, so $(\lambda_n)$ converges to an element $\lambda$ in $\mathbb{R}^\times$. Since $\lambda_nc_n \ra \lambda c$, we see that $\lambda c = d$ (by the Hausdorffness assumption). This shows that $\lambda \in \mathbb{R}^\times$, so $(\mathbb{R}^\times)_{C,D}$ is compact, as required.
\end{proof}
With the former result at hand, we can define the following.
\begin{defn}\label{defn: projectivisation of the normal bundle}
Let $E \ra M$ be a vector bundle of rank $k+1$ ($k\ge0$). The \textit{projectivisation of $E$} is the fiber bundle with typical fiber $\mathbb{RP}^k$
\[\mathbb{P}(E) \ra M,\]
with $\mathbb{P}(E) \coloneqq (E \setminus 0_M)/\mathbb{R}^\times$, where $\mathbb{R}^\times$ acts on $E \setminus 0_M$ by scalar multiplication. Moreover, we denote
$\mathbb{P}(\base,\subbase) \coloneqq \mathbb{P}(\normal(\base,\subbase))$ (warning: this is not standard notation). If $E \ra M$ is a vector bundle of rank $0$, then $\mathbb{P}(E) \coloneqq E$.
\end{defn}
We will now define blow-ups of pairs of manifolds. Fix a pair of manifolds $(\base,\subbase)$ for the rest of this section. Recall that $\DNC(\subbase,\subbase) \cong \subbase \times \mathbb{R}$ can be embedded into $\DNC(\base,\subbase)$ (as a closed embedded submanifold; see Remark \ref{rema: canonical submersion}). %In fact, with the point of view that the blow-up replaces $\subbase$ in $\base$ with the projectivisation of the normal bundle $\normal(\base,\subbase)$, it might not come as a surprise that we can realise the blow-up by taking the quotient of $\DNC(\base,\subbase) \setminus (\subbase \times \mathbb{R})$ by the canonical $\mathbb{R}^\times$-action (see Remark \ref{rema: gauge action on DNC}).

\begin{prop}
Consider the $\mathbb{R}^\times$-action on $\DNC(\base,\subbase)$ (see Remark \ref{rema: gauge action on DNC}). When restricted to $\DNC(\base,\subbase) \setminus \subbase \times \mathbb{R}$, this action is free and (locally) proper (if $\base$ is possibly non-Hausdorff).
\end{prop}
\begin{proof}
The action is easily seen to be free:
\begin{equation*}
\begin{cases}
    \lambda \cdot \left(y,\xi,0\right)= \left(y,\xi,0\right) \textnormal{ iff } \left(y,\lambda^{-1}\xi,0\right) = \left(y,\xi,0\right) & \textnormal{ then from } \xi \neq 0, \textnormal{ we see } \lambda = 1 \\
    \lambda \cdot (x,t) = (x,t) \hspace{1cm} \textnormal{ iff } (x,t) = (x,\lambda t) & \textnormal{ then from } t \neq 0, \textnormal{ we see } \lambda = 1 .
\end{cases}
\end{equation*}

We will show that, in the possibly non-Hausdorff case, the action is locally proper. Then, properness in the Hausdorff case follows by showing that $(\DNC(\base,\subbase) \setminus \subbase \times \mathbb{R})/\mathbb{R}^\times$ is Hausdorff (see Proposition \ref{prop: proper iff locally proper + Hausdorff}), which we will do at the end.

Observe that if we have a chart $(\DNC(U,V),\Phi)$ of $\DNC(\base,\subbase)$ induced by an adapted chart $(U,\varphi)$ (where $V \coloneqq U \cap X$), then the open subset $\DNC(U,V)$ is $\mathbb{R}^\times$-invariant. Therefore, we can restrict to the case that $\DNC(\base,\subbase) = \DNC(\mathbb{R}^{p+q},\mathbb{R}^p)$. Moreover, the map
\begin{equation*}
\Psi^{-1}: \DNC(\mathbb{R}^{p+q},\mathbb{R}^p) \rightarrow \mathbb{R}^{p+q} \times \mathbb{R} \textnormal{ given by } \left(x,\xi,t\right) \mapsto
    \begin{cases}
    \left(x,\xi,0\right) & \text{if}\ t=0 \\
    \left(x,t^{-1}\xi,t\right) & \text{if}\ t\neq0,
    \end{cases}
\end{equation*}
is by construction a diffeomorphism and maps $X \times \mathbb{R}$ to $\mathbb{R}^p \times \mathbb{R}$. Define an $\mathbb{R}^\times$-action on $\mathbb{R}^{p+q} \times \mathbb{R}$ by 
\[\lambda \cdot (x,\xi,t) = (x,\lambda^{-1}\xi,\lambda t).\]
Then this action is smooth and we see that $\Psi^{-1}$ is an equivariant map with respect to the given $\mathbb{R}^\times$-actions. Hence, it suffices to prove that the $\mathbb{R}^\times$-action on $\mathbb{R}^n \times \mathbb{R} \setminus (\mathbb{R}^{p} \times \mathbb{R})$ is proper. 
To show properness in this case, we can prove that whenever $C,D \subset 
\mathbb{R}^n\times\mathbb{R} \setminus \left(\mathbb{R}^p \times \mathbb{R}\right)$ are compact sets, then the set
\[(\mathbb{R}^\times)_{C,D} = \{\lambda \in \mathbb{R}^\times \mid \lambda C \cap D \neq \emptyset\}\] 
is compact in $\mathbb{R}^\times$. 
It suffices to prove that if $(\lambda_n) \subset (\mathbb{R}^\times)_{C,D}$ is a sequence, then it admits a convergent subsequence. Pick a sequence $(c_n)\coloneqq(x_n,\xi_n,t_n)
$ such that $\{\lambda_nc_n\}\eqqcolon\{(x_n,\eta_n,s_n)\} 
\subset D$. By passing to subsequences, we may assume $c_n \ra (x,\xi,t) \in C$ and $\lambda_nc_n \ra (x,\eta,s) \in D$. 
Observe that $(\lambda_n)$ is bounded; if not, then we can pass to a subsequence such that $\lambda_{n}^{-1} \ra 0$. But then
\[(x_{n},\lambda_{n}^{-1}\xi_{n},\lambda_{n}t_{n}) = \lambda_{n} \cdot c_{n} = (x_{n},\eta_{n},s_{n}) \ra (x,\eta,s)\]
and $\lambda_{n}^{-1}\xi_{n} \ra 0$, so $\eta = 0$, which we assumed not to be true. So, by passing to a subsequence, we can assume $(\lambda_n)$ converges in $\mathbb{R}$. If $\lambda = 0$, then by a similar argument, we would have that $\xi=0$, which is not the case. We now have
\[\lambda_{n} \cdot c_{n} = (x_{n},\lambda_{n}^{-1}\xi_{n},\lambda_{n}t_{n}) \ra (x,\lambda^{-1}\xi,\lambda t) = \lambda \cdot c \textnormal{ and } \lambda_{n} \cdot c_n \ra d,\]
so that $\lambda \cdot c = d$. We see that $\lambda \in (\mathbb{R}^\times)_{C,D}$, so this shows that $(\mathbb{R}^\times)_{C,D}$ is compact.

To show that $(\DNC(\base,\subbase) \setminus \subbase \times \mathbb{R})/\mathbb{R}^\times$ is Hausdorff (in case $\base$ is Hausdorff), observe that the open subsets $\DNC(U,V)$ (where $U$ is adapted to $\subbase$ and $V \coloneqq U \cap \subbase$) induce open subsets of $(\DNC(\base,\subbase) \setminus \subbase \times \mathbb{R})/\mathbb{R}^\times$, and since they are $\mathbb{R}^\times$-invariant, it is readily verified that we can separate distinct elements of this space by such open subsets. This proves the statement.
\end{proof}
\begin{defn}\cite{2017arXiv170509588D}\label{defn: blow-up manifold}
We call the smooth (Hausdorff if $\base$ is Hausdorff) manifold
\[\blup(\base,\subbase) \coloneqq \left(\DNC(\base,\subbase) \setminus \subbase \times \mathbb{R}\right) / \mathbb{R}^\times\]
the \textit{blow-up of $\subbase$ in $\base$}. 
\end{defn}

%\begin{rema}\label{rema: blup(Y,Y) is the empty set}
%Observe that $\blup(Y,Y) = (\DNC(Y,Y) \setminus Y \times \mathbb{R}) / \mathbb{R}^\times = \emptyset$.
%\end{rema}

\begin{rema}\label{rema: smooth map gives rise to smooth map of Blup}
By the functoriality of the deformation to the normal cone construction, we can immediately conclude the following functoriality statement for blow-ups: if $(\basetwo,\subbasetwo)$ is another pair of smooth manifolds, then $\subbasetwo \times \mathbb{R} \subset \DNC(\basetwo,\subbasetwo)$ is closed. In particular, if $f: (\base,\subbase) \ra (\basetwo,\subbasetwo)$ is a map of pairs, then we can restrict $\DNC(f): \DNC(\base,\subbase) \ra \DNC(\basetwo,\subbasetwo)$ to the open, $\mathbb{R}^\times$-invariant subset
\[\DNC_f(\base,\subbase) \coloneqq \DNC(\base,\subbase) \setminus \DNC(f)^{-1}(\subbasetwo \times \mathbb{R}) \subset \DNC(\base,\subbase)\]
so that we obtain a smooth map 
\[\blup(f): \blup_f(\base,\subbase) \ra \blup(\basetwo,\subbasetwo),\] 
where $\blup_f(\base,\subbase)$ is the open subset $\DNC_f(\base,\subbase) / \mathbb{R}^\times$ of $\blup(\base,\subbase)$ (we use here that $\DNC(f)$ intertwines the $\mathbb{R}^\times$-action; see Remark \ref{rema: gauge action on DNC maps}). %If $f$ is a submersion, then $\blup(f)$ is also a submersion, as $\DNC(f)|_{\DNC_f(\base,\subbase)}$ is the restriction of a submersion to an open subset. Moreover,
Since $\DNC(\id_\base) = \id_{\DNC(\base,\subbase)}$, where we view $\id_\base$ as a map $(\base,\subbase) \ra (\base,\subbase)$, we have $\blup(\id_\base) = \id_{\blup(\base,\subbase)}$.
\end{rema}
\begin{term}\label{term: subscript f for blup and DNC}
Let $f: (\base,\subbase) \ra (\basetwo,\subbasetwo)$ be a map of pairs. Then we denote
\begin{align*}
    \DNC_f(\base,\subbase) \coloneqq \DNC(\base,\subbase) \setminus \DNC(f)^{-1}(\subbasetwo \times \mathbb{R})&; \quad \blup_f(\base,\subbase) \coloneqq \DNC_f(\base,\subbase)/\mathbb{R}^\times; \\
    \normal_f(\base,\subbase) \coloneqq \normal(\base,\subbase) \setminus d_\normal f^{-1}(0_\subbasetwo)&; \quad \mathbb{P}_f(\base,\subbase) \coloneqq \normal_f(\base,\subbase)/\mathbb{R}^\times.
\end{align*}
Moreover, if $f_i: (\base,\subbase) \ra (\basetwo_i,\subbasetwo_i)$ are maps of pairs ($i=1,2$), then we write $\DNC_{f_1,f_2}(\base,\subbase) \coloneqq \DNC_{f_1}(\base,\subbase) \cap \DNC_{f_2}(\base,\subbase)$ (and similarly if we replace ``$\DNC$'' with ``$\blup$'', ``$\normal$'' or ``$\mathbb{P}$'').
\end{term}
\begin{rema}\cite{2017arXiv170509588D}\label{rema: notation f_1,f_2}
If $f_i: (\base,\subbase) \ra (\basetwo_i,\subbasetwo_i)$ are maps of pairs ($i=1,2$), then, explicitly, 
\begin{equation*}
    \DNC_{f_1,f_2}(\base,\subbase) = \DNC(\base,\subbase) \setminus (\DNC(f_1)^{-1}(\subbasetwo_1 \times \mathbb{R}) \cup \DNC(f_2)^{-1}(\subbasetwo_2 \times \mathbb{R})),
\end{equation*}
and
\begin{equation*}
    \normal_{f_1,f_2}(\base,\subbase) = \normal(\base,\subbase) \setminus (d_\normal f_1^{-1}(0_{\subbasetwo_1}) \cup d_\normal f_2^{-1}(0_{\subbasetwo_2})).
\end{equation*}\vspace{-\baselineskip}
\end{rema}
As in Section \ref{sec: Blow-up in the algebraic setting}, we have a natural blow-down map associated to the blow-up of a pair of manifolds. This map relates the blow-up to the original manifold, and is therefore very important.
\begin{rema}\label{rema: blowdown map}
Recall that we have a smooth map $\hat{\pr}_1: \DNC(\base,\subbase) \ra \base$ given by 
\[\DNC(\base,\subbase) \xra{\DNC(\id_\base)} \DNC(\base,\base) \cong \base \times \mathbb{R} \xra{\pr_1} \base.\]
The map is a composition of $\mathbb{R}^\times$-equivariant maps (with trivial action on $\base$), so we get a canonical \textit{blowdown-map} (compare with the discussion below Definition \ref{defn: blow-up 0 in A^n}) given explicitly by
\begin{equation*}
\bldown: \blup(\base,\subbase) \ra \base, \textnormal{ given by } z \mapsto
\begin{cases}
y & \text{if}\ z = [y,\xi] \\
x & \text{if}\ z = [x,1], \\
\end{cases}
\end{equation*}
(we use here that $d_\normal\id_\base: \normal(\base,\subbase) \ra \normal(\base,\base) \xra{\sim} \base$ is the map $(y,\xi) \mapsto (y,0) \mapsto y$). 
Note that if $x,x' \in \base\setminus\subbase$ are distinct, then $(x,\lambda t) \neq (x',t)$ for any $t,\lambda \in \mathbb{R}^\times$, so $\bldown$ restricts to a bijection between $\bldown^{-1}(\base\setminus\subbase)$ and $\base \setminus \subbase$. It is even a diffeomorphism: recall that we can view $\base \times \mathbb{R}^\times$ as an open submanifold of $\DNC(\base,\subbase)$. In particular,
\[\iota: (\base \setminus \subbase) \times \mathbb{R}^\times \hookrightarrow \DNC(\base,\subbase) \setminus \subbase \times \mathbb{R},\]
is an open embedding, and this map is obviously $\mathbb{R}^\times$-equivariant, so it descends to an embedding $\base \setminus \subbase \ra \blup(\base,\subbase)$.
%Therefore, we can embed $\base\setminus\subbase$ into $\blup(\base,\subbase)$ (note: this embedding is open). 
The restriction of $\bldown$ to 
\[\bldown^{-1}(\subbase) = \mathbb{P}(\base,\subbase) = (\normal(\base,\subbase) \setminus 0_\subbase)/\mathbb{R}^\times \subset \blup(\base,\subbase)\] 
is the (projection map of the) projectivisation of the normal bundle $\mathbb{P}(\base,\subbase)=\mathbb{P}(\normal(\base,\subbase))$. In particular, $\bldown|_{\bldown^{-1}(\subbase)}: \bldown^{-1}(\subbase) \ra \subbase$ is a submersion. %, which is the locally trivial fiber bundle \[\mathbb{P}(\normal(\base,\subbase)) \ra X, \textnormal{ where } \mathbb{P}(\normal(\base,\subbase))_x = \normal(\base,\subbase)_x/\mathbb{R}^\times,\] where the $\mathbb{R}^\times$-action is the usual action that identifies the lines through the origin of a vector space.
Roughly speaking, we can say that $\blup(\base,\subbase)$ is obtained by ``replacing $\subbase$ with the projectivisation of the normal bundle of $\subbase$ in $\base$''. Lastly, notice that if $(U,\varphi)$ (with $V \coloneqq U \cap \subbase$) is an adapted chart of $\base$, then we can realise $\bldown^{-1}(U)$ as $\blup(U,V)$ (more about this later). 
\end{rema}
From the former two remarks, we can easily conclude the following property of blow-ups.
\begin{lemm}\cite{2017arXiv170509588D,10.36045/bbms/1292334057}\label{lemm: induced global smooth map on blow up}
Let $(\base,\subbase)$ and $(\basetwo,\subbasetwo)$ be pairs of manifolds. If $f: (\base,\subbase) \ra (\basetwo,\subbasetwo)$ is a map of pairs, then there is a unique smooth map $\blup(f): \blup_f(\base,\subbase) \ra \blup(\basetwo,\subbasetwo)$ such that the diagram
\begin{center}
\begin{tikzcd}
    \blup_f(\base,\subbase) \ar{r}{\blup(f)} \ar{d}{\bldown} & \blup(\basetwo,\subbasetwo) \ar{d}{\bldown} \\
    \base \ar{r}{f} & \basetwo
\end{tikzcd}
\end{center}
commutes. Moreover, if $f^{-1}(\subbasetwo)=\subbase$ and $d_\normal f$ is fiberwise injective, then $\blup_f(\base,\subbase)=\blup(\base,\subbase)$. %If $\basetwo$ is Hausdorff, then $\blup(f)$ is the unique smooth map $\blup(\base,\subbase) \ra \blup(\basetwo,\subbasetwo)$ making the diagram commute.
\end{lemm}
\begin{proof}
The first statement (apart from the uniqueness statement) follows from Remark \ref{rema: smooth map gives rise to smooth map of Blup} and by tracing the definitions to check that the diagram commutes:
\begin{align*}
    \bldown \circ \blup(f)(z) = 
    \begin{cases}
    \bldown([f(y),d_\normal f(y)\xi,0]) \\
    \bldown([f(x),t])
    \end{cases}
    = 
    \begin{cases}
    f(y) & \text{if}\ z=[y,\xi,0] \\
    f(x) & \text{if}\ z=[x,1].
    \end{cases}
\end{align*}
Since
\[\DNC(f)^{-1}(\subbasetwo \times \mathbb{R}) = \{(y,\xi) \in \normal(\base,\subbase) \mid d_\normal f(y)\xi=0\} \cup \{(x,t) \in \base \times \mathbb{R}^\times \mid f(x) \in \subbasetwo\}\]
the induced map $\blup(f)$ is a smooth map $\blup(\base,\subbase) \ra \blup(\basetwo,\subbasetwo)$ if $f^{-1}(\subbasetwo)=\subbase$ and $d_\normal f$ is fiberwise injective.
The uniqueness statement can be proven as follows: the map $\blup(f)$ restricted to $\base \setminus \subbase \subset \blup(\base,\subbase)$ has to be equal to $f|_{\base \setminus \subbase}$ by commutativity of the diagram. Since $\base \setminus \subbase$ is dense as a subspace of $\blup(\base,\subbase)$ (the complement is the codimension $1$ submanifold $\mathbb{P}(\base,\subbase) \subset \blup(\base,\subbase)$), there is at most one extension of the map $f|_{\base \setminus \subbase}$, seen as a map from $\base \setminus f^{-1}(\subbasetwo) \subset \blup_f(\base,\subbase)$ to $\blup(\basetwo,\subbasetwo)$, to a map $\blup_f(\base,\subbase) \ra \blup(\basetwo,\subbasetwo)$ if $\basetwo$ is Hausdorff. If $\basetwo$ is not Hausdorff, then observe that, if $U \subset \base$ is open (with $V \coloneqq U \cap \subbase$), then the inclusion $\iota: (U,V) \hookrightarrow (\base,\subbase)$ satisfies $\iota^{-1}(\subbase)=V$ and $d_\normal\iota: \normal(U,V) \ra \normal(\base,\subbase)$ is fiberwise injective (even fiberwise an isomorphism). Therefore, it induces a smooth map $\blup(U,V) \ra \blup(\base,\subbase)$, and one can check that it is an open embedding onto $\bldown^{-1}(U) \subset \blup(\base,\subbase)$. It is now readily verified that, if $f$ restricts to a map $f|_U: (U,V) \ra (W,Z)$, $W$ is Hausdorff, then $\blup(f)$ restricts to the map
\[\blup(f|_U): \blup_{f|_U}(U,V) \ra \blup(W,Z),\] 
which is uniquely determined by the property that the diagram
\begin{center}
\begin{tikzcd}
    \blup_{f|_U}(U,V) \ar{r}{\blup(f|_U)} \ar{d}{\bldown} & \blup(W,Z) \ar{d}{\bldown} \\
    U \ar{r}{f|_U} & W
\end{tikzcd}
\end{center}
commutes. Therefore, $\blup(f)$ is uniquely determined by the mentioned property also when $\basetwo$ is not Hausdorff.
This proves the statement.
\end{proof}
\begin{rema}\label{rema: clean intersection normal derivative}
Notice that a map of pairs $f: (\base,\subbase) \ra (\basetwo,\subbasetwo)$ satisfies the properties $f^{-1}(\subbasetwo)=\subbase$ and $d_\normal f$ is fiberwise injective if and only if $f: \base \ra \basetwo$ has a clean intersection with $\subbasetwo \subset \basetwo$. Indeed, this follows from the fact that $d_\normal f$ being fiberwise injective means that
\[T_yf^{-1}(\subbasetwo) = df(y)^{-1}T_{f(y)}\subbasetwo\]
for all $y \in f^{-1}(\subbasetwo)$.
\end{rema}
From this we even obtain a universal property of blow-ups. We will, however, postpone this result until after we gave a local description of blow-ups. From the lemma above we obtain a ``blow-up functor'' if we restrict attention to maps of pairs satisfying the hypothesis given in that lemma.
\begin{rema}\label{rema: blup is functor}
Let $f$ be as in Lemma \ref{lemm: induced global smooth map on blow up} and let $g: (\basetwo,\subbasetwo) \ra (P,Q)$ also satisfy $g^{-1}(Q)=\subbasetwo$ and $d_\normal g$ is fiberwise injective. Then
\[\blup(g \circ f) = \blup(g) \circ \blup(f),\]
since $\DNC(g \circ f) = \DNC(g) \circ \DNC(f)$. If we restrict our category of pairs of smooth manifolds by only considering morphisms $f: (\base,\subbase) \ra (\basetwo,\subbasetwo)$ between pairs such that $f^{-1}(\subbasetwo)=\subbase$ and $d_\normal f$ is fiberwise injective, then we obtain a functor $\blup$ from this category to the category of smooth manifolds. 
\end{rema}
We will now give the local description of blow-ups.
\begin{rema}\label{rema: local coordinates for blowup}
To extract local coordinates for $\blup(\base,\subbase)$ from the local coordinates of the manifold $\DNC(\base,\subbase)$, recall that if $(U,\varphi)$ is an adapted chart of $\base$, then we obtain a chart of the form $(\DNC(U,V), \Phi: \DNC(U,V) \ra \Omega^U_V)$ of $\DNC(\base,\subbase)$ (see the discussion before Lemma \ref{lemm: transition maps of DNC are smooth}). If we write $\varphi=(\varphi^1,\varphi^2) = (y^1,\dots,y^p,x^1,\dots,x^q)$, recall that
\begin{equation*}
    \Phi(z) = (y^1,\dots,y^p,\hat{x}^1,\dots,\hat{x}^q,t)(z) =
\begin{cases}
(\varphi^1(y), d\varphi^2(y)\xi,0) & \text{if}\ z=(y,\xi) \\
(\varphi^1(x),t^{-1}\varphi^2(x),t) & \text{if}\ z=(x,t)
\end{cases}
%\begin{cases}
%(y^1(y),\dots,y^p(y),dx^1(y)\xi,\dots,dx^q(y)\xi,0) & \text{if}\ z=(y,\xi) \\
%(y^1(y),\dots,y^p(y),t^{-1}x^1(x),\dots,t^{-1}x^q(x),t) & \text{if}\ z=(x,t).
%\end{cases}
\end{equation*}
(see Remark \ref{rema: smooth functions DNC}). 
We claim that we get $q$ induced charts $(U_i,\Phi_i)$ on $\blup(\base,\subbase)$ by setting
\[U_i \coloneqq (\DNC(U,V) \setminus \Phi^{-1}(\{(y,\xi,t) \in \Omega^U_V \mid \xi_i = 0\}))/\mathbb{R}^\times\]
and
\begin{align*}
    \Phi_i \coloneqq (y^1,\dots,y^p,\tilde{x}^1_i,\dots,\tilde{x}^q_i)
\end{align*}
where 
\begin{align*}
\tilde{x}^j_i(z) \coloneqq 
\begin{cases}
    \frac{dx^j(y)\xi}{dx^i(y)\xi} & \text{if}\ z=[y,\xi]\\
    \frac{x^j(x)}{x^i(x)} & \text{if}\ z=[x,1]
\end{cases} \textnormal{ if } j \neq i, \textnormal{ and } \tilde{x}^i_i \coloneqq \begin{cases}
    0 & \text{if}\ z=[y,\xi]\\
    x^i(x) & \text{if}\ z=[x,1],
\end{cases} 
\end{align*}
that is,
\begin{equation}\label{eq: local coordinates blup}
    \Phi_i(z) \coloneqq 
\begin{cases}
(y^1(y),\dots,y^p(y),\tfrac{dx^1(y)\xi}{dx^i(y)\xi},\dots,\tfrac{dx^{i-1}(y)\xi}{dx^i(y)\xi},0,\tfrac{dx^{i+1}(y)\xi}{dx^i(y)\xi},\dots,\tfrac{dx^q(y)\xi}{dx^i(y)\xi}) & \text{if}\ z=[y,\xi] \\ 
(y^1(x),\dots,y^p(x),\tfrac{x^1(x)}{x^i(x)},\dots,\tfrac{x^{i-1}(x)}{x^i(x)},x^i(x),\tfrac{x^{i+1}(x)}{x^i(x)},\dots,\tfrac{x^q(x)}{x^i(x)}) & \text{if}\ z=[x,1].
\end{cases}
\end{equation}
%First let us restrict our attention to the case where 
%\[\subbase = \mathbb{R}^p \times \{0\} \subset \mathbb{R}^{p+q} = \base.\]
%Then, with respect to the $\mathbb{R}^\times$-action
Notice that we have an $\mathbb{R}^\times$-action on $\mathbb{R}^{p+q} \times \mathbb{R}$ given by
\[\lambda \cdot (y,\xi,t) = (y,\lambda^{-1}\xi,\lambda t),\]
so that $\Psi^{-1}: \DNC(\mathbb{R}^n,\mathbb{R}^p) \ra \mathbb{R}^{p+q} \times \mathbb{R}$ becomes an $\mathbb{R}^\times$-equivariant map (note: the open sets $\Omega^U_V \subset \mathbb{R}^{p+q} \times \mathbb{R}$ are $\mathbb{R}^\times$-invariant). Since $\Phi = \Psi^{-1} \circ \DNC(\varphi)$, we see that $\Phi$ is also $\mathbb{R}^\times$-equivariant (see Remark \ref{rema: gauge action on DNC maps}), and it maps $V \times \mathbb{R}$ onto $\Omega^U_V \cap (\mathbb{R}^p \times \mathbb{R})$. Therefore, we get an induced diffeomorphism
\[\Phi_{\blup}: (\DNC(U,V) \setminus V \times \mathbb{R})/\mathbb{R}^\times \ra (\Omega^U_V \setminus \mathbb{R}^p \times \mathbb{R})/\mathbb{R}^\times.\]
To ease the notation, set $\base' \coloneqq \Omega^U_V$, $\subbase' \coloneqq \Omega^U_V \cap (\mathbb{R}^p \times \mathbb{R}) \subset \base'$, denote by $\pi'$ the quotient map $\base' \setminus \subbase' \ra (\base' \setminus \subbase')/\mathbb{R}^\times$, and let, for all $1 \le i \le q$, $V_i$ be the open subset $\{(y,\xi,t) \in \base'\setminus \subbase' \mid \xi_i \neq 0\}$. First, observe that $V_i$ is $\mathbb{R}^\times$-invariant. Now, define, for all $1 \le i \le q$, the map
\[\Psi_i: \pi'(V_i) \ra \mathbb{R}^{p+q} \textnormal{ by } [y,\xi,t] \mapsto (y,\xi_{i1},\dots,\xi_{ii-1},t\xi_i,\xi_{ii+1},\dots,\xi_{iq}),\]
where we denoted $\xi_{ik} = \tfrac{\xi_k}{\xi_i}$ (note: these maps are well-defined). The image of this map is the open subset
\[\Omega_i \coloneqq \{(y,\xi) \in \mathbb{R}^{p+q} \mid (y,\xi_i\xi_1,\dots,\xi_i\xi_{i-1},\xi_i,\xi_i\xi_{i+1},\dots,\xi_i\xi_q) \in \varphi(U)\} \subset \mathbb{R}^{p+q}.\]
To see this, notice that if $[y,\xi,t] \in V_i$, then $(y,\xi,t) \in X'=\Omega^U_V$, i.e. $(y,t\xi) \in \varphi(U)$. So,
\[(y,t\xi_i \cdot \xi_{i1},\dots,t\xi_i \cdot \xi_{ii-1},t\xi_i,t\xi_i \cdot \xi_{ii+1},\dots,t\xi_i \cdot \xi_{iq}) = (y,t\xi) \in \varphi(U),\]
from which we see that the image of $\Psi_i$ lies in $\Omega_i$. Conversely, if $(y,\xi) \in \varphi(U)$, then, by definition of $\base'=\Omega^U_V$ and $\subbase'$, $(y,\xi_1,\dots,\xi_{i-1},1,\xi_{i+1},\dots,\xi_q,\xi_i) \in \base' \setminus \subbase'$. Since we have 
\[\Psi_i([y,\xi_1,\dots,\xi_{i-1},1,\xi_{i+1},\dots,\xi_q,\xi_i])=(y,\xi),\] 
this shows that indeed the image of $\Psi_i$ is $\Omega_i$.
We also see from this what the candidate is for the inverse of $\Psi_i$:
\[\Omega_i \ra \pi'(V_i) \textnormal{ given by } (y,\xi) \mapsto [y,\xi_1,\dots,\xi_{i-1},1,\xi_{i+1},\dots,\xi_q,\xi_i]\]
(we will already denote this map by $\Psi_i^{-1}$). It is, indeed, inverse to $\Psi_i$:
\begin{align*}
    \Psi_i([y,\xi_1,\dots,\xi_{i-1},1,\xi_{i+1},\dots,\xi_q,\xi_i]) &= (y,\xi_1,\dots,\xi_{i-1},\xi_i,\xi_{i+1},\dots,\xi_q) \textnormal{ and } \\
    \Psi_i^{-1}(y,\xi_{i1},\dots,\xi_{ii-1},t\xi_i,\xi_{ii+1},\dots,\xi_{iq}) &= [y,\xi_{i1},\dots,\xi_{ii-1},1,\xi_{ii+1},\dots,\xi_{iq},t\xi_i] \\
    &= [\xi_i \cdot (y,\xi_1,\dots,\xi_{i-1},\xi_i,\xi_{i+1},\dots,\xi_{q},t)] \\
    &= [y,\xi_1,\dots,\xi_{i-1},\xi_i,\xi_{i+1},\dots,\xi_q,t].
\end{align*}
Now, the maps $\Psi_i$ and $\Psi_i^{-1}$ are smooth, because the maps $V_i \ra \Omega_i$ and $\Omega_i \ra V_i$ (defined by the same formulas) are smooth. This proves that, for all $1 \le i \le q$, $(\pi'(V_i),\Psi_i)$ is a chart of $(\base'\setminus \subbase')/\mathbb{R}^\times$. Now, define 
\[\Phi_i \coloneqq \Psi_i \circ \Phi_{\blup}: \Phi_{\blup}^{-1}(\pi'(V_i)) \ra \Omega_i.\]
Then $\Phi_{\blup}^{-1}(\pi'(V_i)) = U_i$, %because
%\[U_i = \{[x,\xi,0] \in (\DNC(U,V) \setminus (V \times \mathbb{R}))/\mathbb{R}^\times \mid (d\varphi^2(x)\xi)_i \neq 0\} \cup \{[y,t] \in (\DNC(U,V) \setminus (V \times \mathbb{R}))/\mathbb{R}^\times \mid \varphi^2(y)_i \neq 0\},\]
and $\Phi_i$ is, as claimed, the map \eqref{eq: local coordinates blup}, with inverse given by 
\begin{equation}\label{eq: inverse local coordinates blup}
\Phi_i^{-1}(y,\xi) \coloneqq
\begin{cases}
    [\varphi^{-1}(y),d_\normal\varphi^{-1}(y)(\xi_1,\dots,\xi_{i-1},1,\xi_{i+1},\dots,\xi_q)] & \text{if}\ \xi_i=0 \\
    [\varphi^{-1}(y,\xi_i\xi_1,\dots,\xi_i\xi_{i-1},\xi_i,\xi_i\xi_{i+1},\dots,\xi_i\xi_q),1] & \text{if}\ \xi_i \neq 0.
\end{cases}
\end{equation}
The transition maps are the maps
\[\Phi_{ij} = \Psi_i \circ \Psi_j^{-1},\]
which are smooth because the maps $\Psi_i$ are diffeomorphisms.
\end{rema}
As promised, we will prove that blow-ups have a universal property. 
\begin{prop}[Universal property of blow-ups]\label{prop: universal property blow-up}
The blow-up of $\subbase$ in $\base$, seen as a pair of manifolds $(\blup(\base,\subbase),\mathbb{P}(\base,\subbase))$, is the unique (up to unique diffeomorphism) pair of manifolds $(\widetilde{\base},\widetilde{\subbase})$ that comes with a map of pairs 
\[\widetilde{\bldown}: (\widetilde{\base},\widetilde{\subbase}) \ra (\base,\subbase)\] 
such that (a) $\widetilde{\subbase} \subset \widetilde{\base}$ is of codimension $1$, (b) $\widetilde{\bldown}^{-1}(\subbase) = \widetilde{\subbase}$, and (c) $d_\normal\widetilde{\bldown}$ is fiberwise injective, together with the following (universal) property: whenever $f: (\basetwo,\subbasetwo) \ra (\base,\subbase)$ is a map of pairs such that (1) $\subbasetwo \subset \basetwo$ is of codimension $1$, (2) $f^{-1}(\subbase)=\subbasetwo$, and (3) $d_\normal f$ is fiberwise injective, then there is a unique map of pairs $\widetilde{f}: (\basetwo,\subbasetwo) \ra (\widetilde{\base},\widetilde{\subbase})$ such that the diagram
\begin{center}
\begin{tikzcd}
     & (\widetilde{\base},\widetilde{\subbase}) \ar{d}{\widetilde{\bldown}} \\
    (\basetwo,\subbasetwo) \ar{r}{f} \ar[ru, dotted, "\widetilde{f}"] & (\base,\subbase)
\end{tikzcd}
\end{center}
commutes.
\end{prop}
\begin{proof}
To show that $\blup(\base,\subbase)$ satisfies (a) (b) and (c), observe that, by Remark \ref{rema: blowdown map}, we only have to show that the map $d_\normal\bldown: \normal(\blup(\base,\subbase),\mathbb{P}(\base,\subbase)) \ra \normal(\base,\subbase)$ is fiberwise injective. This will be proven last. To see that $\blup(\base,\subbase)$ has the mentioned universal property, we have to show that
\[\blup(\basetwo,\subbasetwo) \cong \basetwo\] 
if $\subbasetwo \subset \basetwo$ is of codimension $1$ (in fact, the blow-down map is a diffeomorphism). %In principal, the latter statement just follows from the fact that $\mathbb{P}(\basetwo,\subbasetwo) \cong \basetwo$ (because $\normal(\basetwo,\subbasetwo)$ has rank $1$), so that one can show that the blow-down map is a diffeomorphism in this case. However, f
To see this, notice that, from Remark \ref{rema: local coordinates for blowup}; \eqref{eq: inverse local coordinates blup}), we see that $\bldown$ can locally be written as (where $1 \le i \le q$, and we use notation as in Remark \ref{rema: local coordinates for blowup})
\begin{align*}
    \varphi \circ p \circ \Phi_i^{-1}(y,\xi) &= 
    \begin{cases}
    \varphi \circ p[\varphi^{-1}(y),d_\normal\varphi^{-1}(y)(\xi_1,\dots,\xi_{i-1},1,\xi_{i+1},\dots,\xi_q)] & \text{if}\ \xi_i=0 \\
    \varphi \circ p [\varphi^{-1}(y,\xi_i\xi_1,\dots,\xi_i\xi_{i-1},\xi_i,\xi_i\xi_{i+1},\dots,\xi_i\xi_q),1] & \text{if}\ \xi_i \neq 0.
\end{cases} \\
    &= (y,\xi_i\xi_1,\dots,\xi_i\xi_{i-1},\xi_i,\xi_i\xi_{i+1},\dots,\xi_i\xi_q).
\end{align*}
The latter statement now follows immediately, since we see that the blow-down map is a diffeomorphism if $\subbasetwo \subset \basetwo$ is of codimension $1$. That $(\blup(\base,\subbase),\mathbb{P}(\base,\subbase))$ is unique up to unique diffeomorphism follows now from Lemma \ref{lemm: induced global smooth map on blow up}. It remains to show that $d_\normal\bldown$ is fiberwise injective. Recall for this that $\Phi_i^{-1}$ is a map $\Omega_i \ra U_i$ (see Remark \ref{rema: local coordinates for blowup} for the notation), and $\Omega_i \cap \{(y,\xi) \in \mathbb{R}^{p+q} \mid \xi_i=0\}$ maps onto $U_i \cap \mathbb{P}(\base,\subbase)$. From the above local form of $p$, we see that we can write $d_\normal p$ locally as the map
\[(y,\xi,\eta) \mapsto (y,\xi_1\eta, \dots, \eta, \dots, \xi_q\eta),\]
which is fiberwise injective, as required. 
\end{proof}

Here is a simple example illustrating what blowing up a point does (see also Remark \ref{rema: connected sum}).
\begin{exam}\label{exam: blup of point in sphere}
We will show that blowing up a point in the two-sphere $S^2 \subset \mathbb{R}^3$ results in $\mathbb{RP}^2$. One can easily extend the argument to show that blowing up a point in the $n$-sphere $S^n$ yields $\mathbb{RP}^n$. 

We denote by $+1,-1 \in S^2$ the north ($=(0,0,1)$) and south ($=(0,0,-1)$) pole, respectively. Consider the standard charts given by stereographic projection:
\[(U_\pm,\varphi_\pm) \textnormal{ where } U_\pm \coloneqq S^2 \setminus \pm 1 \textnormal{ and } (x^1_{\pm},x^2_{\pm}) \coloneqq \varphi_\pm(x_0,x_1,x_2) \coloneqq \frac{1}{1 \mp x_2}(x_0,x_1).\]
We will show that $\blup(S^2,+1) \cong \mathbb{RP}^2$. Indeed, observe first that, with respect to the embedding $S^2 \subset \mathbb{R}^3$, we have 
\[d\varphi_{\pm}(x_0,x_1,x_2) = \begin{pmatrix} \frac{1}{1 \mp x_2} & 0 & \pm \frac{x_0}{(1\mp x_2)^2} \\
0 & \frac{1}{1 \mp x_2} & \pm \frac{x_0}{(1\mp x_2)^2} \end{pmatrix} \textnormal{, and } T_{+1}S^2 = \{(\xi_0,\xi_1,\xi_2) \in \mathbb{R}^3 \mid \xi_2=0\}.\]
We see that, by Remark \ref{rema: local coordinates for blowup}, $\varphi_-$ induces the local coordinates on $\blup(S^2,+1)$ given by
\begin{align*}
&z \mapsto
\begin{cases}
    (0,\tfrac{dx_-^2(+1)\xi}{dx_-^1(+1)\xi}) \\
    (x_-^1(x),\tfrac{x_-^2(x)}{x_-^1(x)})
\end{cases} = 
\begin{cases}
    (0,\tfrac{\xi_1}{\xi_0}) & \text{if}\ z=[+1,\xi] \\
    (\tfrac{x_0}{1+x_2},\tfrac{x_1}{x_0}) & \text{if}\ z=[x,1]
\end{cases}
\textnormal{ and }\\
&z \mapsto
\begin{cases}
    (\tfrac{dx^1_-(+1)\xi}{dx_-^2(+1)\xi},0) \\
    (\tfrac{x_-^1(x)}{x_-^2(x)},x_-^2(x))
\end{cases} = 
\begin{cases}
    (\tfrac{\xi_0}{\xi_1},0) & \text{if}\ z=[+1,\xi] \\
    (\tfrac{x_0}{x_1},\tfrac{x_1}{1+x_2}) & \text{if}\ z=[x,1].
\end{cases}
\end{align*}    
Recall that these maps map onto $\Omega_1 = \{(a,b) \in \mathbb{R}^2 \mid (a,ab) \in \mathbb{R}^2\} = \mathbb{R}^2 = \Omega_2$. The local coordinates induced by $\varphi_+$ on $\blup(S^2,+1)$ are given by 
\begin{align*}
&[x,1] \mapsto (x_+^1(x),\frac{x_+^2(x)}{x^1(x)}) = (\frac{x_0}{1-x_2},\frac{x_1}{x_0})
\textnormal{ and } \\
&[x,1] \mapsto (\frac{x_+^1(x)}{x_+^2(x)},x_+^2(x)) = (\frac{x_0}{x_1},\frac{x_1}{1-x_2})
\end{align*}    
(since $+1 \not\in U_+$). Consider now the map 
%\begin{align*}
%f: \DNC(S^2,+1) \setminus (\{+1\} \times \mathbb{R}) \ra \mathbb{RP}^2 \textnormal{ given by } z \mapsto
%\begin{cases}
%    [\xi_0:\xi_1:0] & \text{if}\ z=(+1,\xi_0,\xi_1,0) \\
%    [x_0:x_1:1-x_2] & \text{if}\ z=(x_0,x_1,x_2,t).
%\end{cases}
%\end{align*}
%Then this map is $\mathbb{R}^*$-equivariant (with trivial action on $\mathbb{RP}^2$) and so we get an induced map
\begin{align*}
    f: \blup(S^2,+1) \ra \mathbb{RP}^2 \textnormal{ given by } z \mapsto
\begin{cases}
    [\xi_0:\xi_1:0] & \text{if}\ z=[+1,\xi] \\
    [x_0:x_1:1-x_2] & \text{if}\ z=[x,1].
\end{cases}
\end{align*}
The map $f$ is readily verified to be well-defined, and bijective, with inverse given by 
\begin{align*}
    z \mapsto
\begin{cases}
    [+1,a_0,a_1] & \text{if}\ z=[a_0:a_1:0] \\
    [\frac{2a_0}{a_0^2+a_1^2+1},\frac{2a_1}{a_0^2+a_1^2+1},\frac{a_0^2+a_1^2-1}{a_0^2+a_1^2+1},1] & \text{if}\ z=[a_0:a_1:1].
\end{cases}
\end{align*}
To see that $f$ is a diffeomorphism, we can check it locally using the charts: in order of presentation, these maps are given by
\begin{equation}\label{eq: local expressions f sphere projective space}
    (a,b) \mapsto (a(b^2+1),b); \quad (a,b) \mapsto (a,b(a^2+1)); \quad (a,b) \mapsto (a,ab); \quad (a,b) \mapsto (ab,b).
\end{equation}
For example, the first map is calculated as follows:
\begin{align*}
\begin{cases}
    (0,\tfrac{\xi_1}{\xi_0}) \\
    (\tfrac{x_0}{1+x_2},\tfrac{x_1}{x_0})
\end{cases} \mapsto 
\begin{cases}
    [+1,\xi] \\
    [x,1]
\end{cases} \mapsto 
\begin{cases}
    [\xi_0:\xi_1:0] \\
    [x_0:x_1:1-x_2]
\end{cases} \mapsto 
\begin{cases}
    (0,\tfrac{\xi_1}{\xi_0}) \\
    (\tfrac{1-x_2}{x_0},\tfrac{x_1}{x_0}),
\end{cases}
\end{align*}
so using that 
\[\frac{1-x_2}{x_0} = \frac{1-x_2^2}{(1+x_2)x_0} = \frac{1-x_2^2-x_0^2}{(1+x_2)x_0} + \frac{x_0}{1+x_2} = \frac{x_0}{(1+x_2)} \cdot \frac{x_1^2}{x_0^2} + \frac{x_0}{1+x_2},\]
we see that the map is given by $(a,b) \mapsto (a(b^2+1),b)$, as claimed. The above local expressions for $f$ are smooth, and similarly one can check that the inverse is smooth. %Alternatively, notice that the first two local expressions of $f$ (see \eqref{eq: local expressions f sphere projective space}) are diffeomorphisms $\mathbb{R}^2 =\Omega_1=\Omega_2 \ra \mathbb{R}^2$. This shows that the restriction of $f$ to $\blup(S^2,+1)\setminus\{-1\}$ (note: $-1 \in S^2 \setminus \{+1\} \subset \blup(S^2,+1)$) is a diffeomorphism, and $f(-1)=[0:0:1]$ is not in the domain of the chart of $\mathbb{RP}^2$ that we used for the local expressions. 
This proves that $f$ is a diffeomorphism between $\blup(S^2,+1)$ and $\mathbb{RP}^2$.
\end{exam}
\begin{rema}\label{rema: connected sum}
More generally, if $\base$ is a manifold, one can show that $\blup(\base,\textnormal{pt})\cong \mathbb{RP}^n \# \base$, where $\#$ denotes the connected sum operation (see e.g. \cite{Lee} for the connected sum operation). 
\end{rema}
Later, in Section \ref{sec: Blow-up of a pair of manifolds: a different approach (part 1)}, we will see that there are other approaches to blow-ups. In view of this, let us describe the blow-up of a point; at least locally. Compare it with Definition \ref{defn: blow-up 0 in A^n}.
\begin{prop}\label{prop: blup 0 in R^n is same as our construction}
There is a canonical embedding 
\[\blup(\mathbb{R}^n,0) \hookrightarrow \mathbb{R}^n \times \mathbb{RP}^{n-1} \textnormal{ given by } z \mapsto
\begin{cases}
    (0,[\xi_1:\cdots:\xi_n]) & \text{if}\ z=[0,\xi] \\
    (x,[x_1:\cdots:x_n]) & \text{if}\ z=[x,1]
\end{cases}\]
onto $\widetilde{\mathbb{R}^n} \coloneqq \{((x_1,\dots,x_n),[y_1:\cdots:y_n]) \mid x_iy_j=y_ix_j \textnormal{ for all } 1 \le i,j \le n\} \subset \mathbb{R}^n \times \mathbb{RP}^{n-1}$, which is a closed embedded submanifold of codimension $n-1$.
\end{prop}
\begin{proof}
Denote the above map by $f$. Notice that $f$ is injective, and it indeed maps onto $\widetilde{\mathbb{R}^n}$: if $x \in \mathbb{R}^n \setminus \{0\}$, say $x_i \neq 0$, then there is only one $[y_1:\cdots:y_n] \in \mathbb{RP}^{n-1}$ with the property that $x_iy_j=y_ix_j$ for all $1 \le i,j \le n$, namely 
\[[y_1:\cdots:y_n] \coloneqq [\frac{x_1}{x_i}:\cdots:\frac{x_{i-1}}{x_i}:1:\frac{x_{i+1}}{x_i}:\cdots:\frac{x_n}{x_i}] = [x_1:\cdots:x_n],\]
and we also see that $f([x,1])=(x,[y_1:\cdots:y_n])$. If $x=0 \in \mathbb{R}^n$, then all $y=[y_1:\cdots:y_n] \in \mathbb{RP}^{n-1}$ satisfy $x_iy_j=y_ix_j$ for all $1 \le i,j \le n$, and $f([0,y])=(0,y)$.
We will now show that $f$ is an injective immersion. Denote by $(U_i^{\mathbb{RP}},\varphi_i)$ the standard charts of $\mathbb{RP}^{n-1}$, and consider the charts $(U_i,\Phi_i)$ of $\blup(\mathbb{R}^n,0)$ (note: $\Phi_i$ is a map $U_i \ra \Omega_i=\mathbb{R}^n$; see Remark \ref{rema: local coordinates for blowup} for the notation). We can then write the above map locally as
\begin{align*}
    (\id_{\mathbb{R}^n} \times \varphi_i) \circ f &\circ \Phi_i^{-1}: \mathbb{R}^n \ra \mathbb{R}^n \times \mathbb{R}^{n-1} \textnormal{ given by } \\
    x &\mapsto 
%\begin{cases}
%    [x_1,\dots,x_{i-1},1,x_{i+1},\dots,x_n] \\
%    [x_ix_1,\dots,x_ix_{i-1},x_i,x_ix_{i+1},\dots,x_ix_n]
%\end{cases} \mapsto 
%\begin{cases}
%    [x_1:\dots:x_{i-1}:1:x_{i+1}:\dots,x_n] \\
%    (x_ix_1,\dots,x_ix_{i-1},x_i,x_ix_{i+1},\dots,x_ix_n,[x_ix_1:\dots:x_ix_{i-1}:x_i:x_ix_{i+1}:\dots:x_ix_n])
%\end{cases} \mapsto
\begin{cases}
    (0,x_1,\dots,x_{i-1},x_{i+1},\dots,x_n) & \text{if}\ x_i=0 \\
    (x_ix_1,\dots,x_ix_{i-1},x_i,x_ix_{i+1},\dots,x_ix_n,x_1,\dots,x_{i-1},x_{i+1},\dots,x_n) & \text{if}\ x_i \neq 0,
\end{cases}
\end{align*}
so it is a smooth map. Moreover, by post-composing this map with the diffeomorphism
\begin{align*}
    \mathbb{R}^n \times \mathbb{R}^{n-1} &\xra{\sim} \mathbb{R}^n \times \mathbb{R}^{n-1} \textnormal{ given by } \\
    (x,y) &\mapsto (x_1-x_iy_1,\dots,x_{i-1}-x_iy_{i-1},x_i,x_{i+1}-x_iy_{i+1},\dots,x_n-x_iy_n,y),
\end{align*}
we see that $f$ is an (injective) immersion. It remains to show that $\widetilde{\mathbb{R}^n}$ is a closed embedded submanifold of codimension $n-1$ (it follows that $f$ is an embedding). Indeed, $\widetilde{\mathbb{R}^n} \cap (\mathbb{R}^n \times U_i)$ is given by the zero set of the map
\[\mathbb{R}^n \times U_i \ra \mathbb{R}^{n-1} \textnormal{ given by } (x,y) \mapsto (x_i\frac{y_1}{y_i}-x_1,\dots,\widehat{x_i\frac{y_i}{y_i} - x_i},\dots,x_i\frac{y_n}{y_i}-x_n),\]
which is a smooth submersion, since, locally, the map becomes
\begin{align*}
    \mathbb{R}^n \times \mathbb{R}^{n-1} &\ra \mathbb{R}^{n-1} \textnormal{ given by } \\
    (x,y) &\mapsto (x,[y^1:\cdots:1:\cdots:y^n]) \\
    &\mapsto (x_iy_1-x_1,\dots,x_iy_{i-1}-x_{i-1},x_iy_{i+1}-x_{i+1},\dots,x_iy_n-x_n)
\end{align*}
so we see that it becomes a projection map if we pre-compose the map with the diffeomorphism $\mathbb{R}^n \times \mathbb{R}^{n-1} \xra{\sim} \mathbb{R}^n \times \mathbb{R}^{n-1}$ from above. This proves the statement.
\end{proof}
\subsection{Properties related to the blow-up functor}\label{sec: properties related to the blow-up functor}
Here, we will go through some more properties that blow-ups have. Most of these results will follow from the results presented in Section \ref{sec: Properties related to the deformation to the normal cone and normal bundle functors}. In this section, we fix a pair of smooth manifolds $(\base,\subbase)$.
\begin{prop}\cite{10.36045/bbms/1292334057}\label{prop: blowdown is proper}
The blow-down map 
\[\bldown: \blup(\base,\subbase) \ra \base\]
(see Remark \ref{rema: blowdown map}) is proper.
\end{prop}
\begin{proof}
%By Lemma \ref{lemm: proper maps}, and using that $\bldown^{-1}(U)=\blup(U,V)$ for adapted charts $(U,\varphi)$ (with $V \coloneqq U \cap \subbase$), we may reduce to the local case; that is, if $(\base,\subbase)=(\mathbb{R}^n,\mathbb{R}^p)$. Now, 
Let $C \subset \base$ be compact and let $\{\widetilde U_i\}$ be an open cover of $p^{-1}(C)$. Pick an open cover $\{U_j\}$ of $C$ consisting of the domains of adapted charts $(U_j,\varphi_j)$ (with $V_j \coloneqq U_j \cap \subbase$) such that $p^{-1}(U_j) \subset \widetilde U_{i_j}$ for some $i_j$. Since $C$ is compact, there is a finite subcover $\{U_{j_k}\}$ of $\{U_j\}$. Notice that now 
\[\{\bldown^{-1}(U_{j_k})\} = \{\blup(U_{j_k},V_{j_k})\}\]
is a finite open cover of $p^{-1}(C)$, so $\{\widetilde U_{i_{j_k}}\}$ yields the desired finite subcover of $\{\widetilde U_j\}$.
%Again, using Lemma \ref{lemm: proper maps}, we only have to show that the map is closed: the fibers are diffeomorphic to $\mathbb{RP}^{q-1}$ (with $q \coloneqq n-p$). To see that this map is closed
\end{proof}
\begin{rema}\cite{2017arXiv170509588D}\label{rema: blow-up of trivials}
Since $\DNC(\base,\base) \cong \base \times \mathbb{R}$, we see that
\[\blup(\base,\base) = (\DNC(\base,\base) \setminus \base \times \mathbb{R})/\mathbb{R}^\times = \emptyset.\]
Also, recall from the proof of Proposition \ref{prop: universal property blow-up}, that if $\subbase \subset \base$ is of codimension $1$, then 
\[\blup(\base,\subbase) \cong \base.\]
So, we see that blow-ups only yield something non-trivial if $\subbase \subset \base$ is of codimension $\ge 2$.
\end{rema}
The blow-up respects fiber products, but we have to keep in mind Remark \ref{rema: smooth map gives rise to smooth map of Blup}. Before we give the result, it might be useful to give the analogue of Proposition \ref{prop: deformation constant rank} for blow-ups first.
\begin{prop}\label{prop: analogue constant rank when applying blup}
Let $f: (\base,\subbase) \ra (\basetwo,\subbasetwo)$ be a map of pairs. %of constant rank (see Definition \ref{defn: constant rank pairs}). 
If $\DNC(f)|_{\DNC_f(\base,\subbase)}$ is of constant rank $k+1$ (e.g. $\DNC(f)$ is of constant rank $k+1$), then $\blup(f): \blup_f(\base,\subbase) \ra \blup(\basetwo,\subbasetwo)$ is of constant rank $k$. In particular, if $f$ is a submersion (as a map of pairs), then $\blup(f)$ is a submersion, and if $f$ is an immersion, and $f^{-1}(\subbasetwo) = \subbase$, then $\blup(f)$ is an immersion.
\end{prop}
\begin{proof}
That $\blup(f)$ is of constant rank if $\DNC(f)|_{\DNC_f(\base,\subbase)}$ is follows from the lemma below. The other statement is then a direct consequence of Proposition \ref{prop: deformation constant rank}.
\end{proof}
\begin{lemm}\label{lemm: constant rank maps quotient vector space}
Let $f: V \ra W$ be a linear map and let $V' \subset V$ and $W' \subset W$ be vector subspaces with $\ell \coloneqq \dim V' = \dim W'$. If $f$ has rank $k$, and $f(V')=W'$, then $\overline{f}: V/V' \ra W/W'$ has rank $k-\ell$.
\end{lemm}
\begin{proof}
We have
\begin{align*}
    \dim \ker \overline{f} + \dim \im \overline{f} &= \dim V/V' \\
    &= \dim \ker f + \dim \im f - \dim V' \\
    &= \dim \ker f + \dim \im f - \dim W' = \dim \ker f + \dim \im \overline{f},
\end{align*}
so $\dim \ker \overline{f} = \dim \ker f$. Therefore, from the second equality we see that $\dim \im \overline{f} = k - \ell$, which proves the statement.
\end{proof}
From Corollary \ref{coro: DNC submanifold}, we immediately obtain the following useful fact:
\begin{coro}\label{coro: blup submanifold}
Let $(\basetwo,\subbasetwo)$ be a pair of manifolds together with an embedding $\iota: (\basetwo,\subbasetwo) \hookrightarrow (\base,\subbase)$ (that is, $\iota: \basetwo \ra \base$ is an embedding, and $\iota|_{\subbasetwo}: \subbasetwo \ra \subbase$ is an embedding) such that $\iota^{-1}(\subbase)=\subbasetwo$. Then we obtain an embedding
\[\blup(\iota): \blup(\basetwo,\subbasetwo) \ra \blup(\base,\subbase)\]
which, for any map of pairs $f: (\base,\subbase) \ra (P,Q)$, satisfies
\[\blup(\iota)^{-1}(\blup_f(\base,\subbase)) = \blup_{f \circ \iota}(\basetwo,\subbasetwo).\]
Moreover, if $\iota$ is open (resp. closed), then $\blup(\iota)$ is open (resp. closed).
\end{coro}
Recall the notion of clean intersection of pairs (see Definition \ref{defn; clean intersection of pairs}). The statement that blow-ups respect fiber products is as follows:
\begin{prop}\label{prop: blup respects fiber products}
Let $f_i: (\base_i,\subbase_i) \ra (\basetwo,\subbasetwo)$ ($i=1,2$) be smooth maps of pairs that have clean intersection. Then $\blup_{f_1}(\base_1,\subbase_1) \tensor[_{\blup(f_1)}]{\times}{_{\blup(f_2)}} \blup_{f_2}(\base_2,\subbase_2) \subset \blup(\base_1,\subbase_1) \times \blup(\base_2,\subbase_2)$ is an embedded submanifold, and 
the canonical map
\begin{align*}
    \blup_{f_1 \times f_2}(\base_1 \tensor[_{f_1}]{\times}{_{f_2}} \base_2, \subbase_1 \tensor[_{f_1}]{\times}{_{f_2}} \subbase_2) &\ra \blup_{f_1}(\base_1,\subbase_1) \tensor[_{\blup(f_1)}]{\times}{_{\blup(f_2)}} \blup_{f_2}(\base_2,\subbase_2) \\
    z &\mapsto
    \begin{cases}
    ([y_1,\xi_1],[y_2,\xi_2]) & \text{if}\ z=[y_1,y_2,\xi_1,\xi_2] \\
    ([x_1,1],[x_2,1]) & \text{if}\ z=[x_1,x_2,1]%,
    \end{cases}
\end{align*}
%where we denoted $\blup_{f_1,f_2}(\base_1 \tensor[_{f_1}]{\times}{_{f_2}} \base_2, \subbase_1 \tensor[_{f_1}]{\times}{_{f_2}} \subbase_2) \coloneqq \blup_{(f_1 \circ \pr_1,f_2 \circ \pr_2)}(\base_1 \tensor[_{f_1}]{\times}{_{f_2}} \base_2, \subbase_1 \tensor[_{f_1}]{\times}{_{f_2}} \subbase_2)$, 
is a diffeomorphism.
\end{prop}
\begin{proof}
From Proposition \ref{prop: if clean intersection, then DNC respects fiber products}, we know that $\DNC(\base_1,\subbase_1) \tensor[_{\DNC(f_1)}]{\times}{_{\DNC(f_2)}} \DNC(\base_2,\subbase_2) \subset \DNC(\base_1,\subbase_1) \times \DNC(\base_2, \subbase_2)$ is an embedded submanifold, and we have a canonical, and $\mathbb{R}^\times$-equivariant, diffeomorphism $\DNC(\base_1 \tensor[_{f_1}]{\times}{_{f_2}} \base_2, \subbase_1 \tensor[_{f_1}]{\times}{_{f_2}} \subbase_2) \xra{\sim} \DNC(\base_1,\subbase_1) \tensor[_{\DNC(f_1)}]{\times}{_{\DNC(f_2)}} \DNC(\base_2,\subbase_2)$ (given by the same formula as above). Since
\[\DNC_{f_1}(\base_1,\subbase_1) \tensor[_{\DNC(f_1)}]{\times}{_{\DNC(f_2)}} \DNC_{f_2}(\base_2,\subbase_2) \subset \DNC(\base_1,\subbase_1) \tensor[_{\DNC(f_1)}]{\times}{_{\DNC(f_2)}} \DNC(\base_2,\subbase_2)\]
is an open subset, it is an embedded submanifold of $\DNC(\base_1,\subbase_1) \times \DNC(\base_2,\subbase_2)$ as well. Now, the embedding
\[\DNC_{f_1}(\base_1,\subbase_1) \tensor[_{\DNC(f_1)}]{\times}{_{\DNC(f_2)}} \DNC_{f_2}(\base_2,\subbase_2) \hookrightarrow \DNC(\base_1,\subbase_1) \times \DNC(\base_2,\subbase_2)\]
descends to an embedding (see Proposition \ref{prop: analogue constant rank when applying blup})
\[\blup_{f_1}(\base_1,\subbase_1) \tensor[_{\blup(f_1)}]{\times}{_{\blup(f_2)}} \blup_{f_2}(\base_2,\subbase_2) \hookrightarrow \blup(\base_1,\subbase_1) \times \blup(\base_2,\subbase_2).\]
For the last statement, observe that
\[\DNC(\base_1 \tensor[_{f_1}]{\times}{_{f_2}} \base_2,\subbase_1 \tensor[_{f_1}]{\times}{_{f_2}} \subbase_2) \cap \DNC(f_1 \times f_2)^{-1}(\subbasetwo \times \subbasetwo \times \mathbb{R}) %\subset \DNC(\base_1 \tensor[_{f_1}]{\times}{_{f_2}} \base_2, \subbase_1 \tensor[_{f_1}]{\times}{_{f_2}} \subbase_2)
\]
maps, under the diffeomorphism from Proposition \ref{prop: if clean intersection, then DNC respects fiber products}, onto 
\[\DNC(f_1)^{-1}(\subbasetwo) \tensor[_{\DNC(f_1)}]{\times}{_{\DNC(f_2)}} \DNC(\base_2,\subbase_2) \cup \DNC(\base_1,\subbase_1) \tensor[_{\DNC(f_1)}]{\times}{_{\DNC(f_2)}} \DNC(f_2)^{-1}(\subbasetwo), %\subset \DNC(\base_1,\subbase_1) \tensor[_{\DNC(f_1)}]{\times}{_{\DNC(f_2)}} \DNC(\base_2,\subbase_2)
\]
which equals $\DNC(f_1)^{-1}(\subbasetwo) \tensor[_{\DNC(f_1)}]{\times}{_{\DNC(f_2)}} \DNC(f_2)^{-1}(\subbasetwo)$. Therefore, we obtain a diffeomorphism
\[\DNC_{f_1 \times f_2}(\base_1 \tensor[_{f_1}]{\times}{_{f_2}} \base_2, \subbase_1 \tensor[_{f_1}]{\times}{_{f_2}} \subbase_2) \xra{\sim} \DNC_{f_1}(\base_1,\subbase_1) \tensor[_{\DNC(f_1)}]{\times}{_{\DNC(f_2)}} \DNC_{f_2}(\base_2,\subbase_2),\]
(note: we used here that the latter manifold is an embedded submanifold of $\DNC(\base_1,\subbase_1) \times \DNC(\base_2,\subbase_2)$) which descends to the desired isomorphism. 
\end{proof}

The following fact, also a consequence of Proposition \ref{prop: if clean intersection, then DNC respects fiber products}, will be used many times.
\begin{prop}\label{prop: blup(X x M, Y x M) = blup(X,Y) x M}
Let $\basetwo$ be a smooth manifold. Then the canonical map
\[\blup(\base \times \basetwo,\subbase \times \basetwo) \ra \blup(\base,\subbase) \times \basetwo \textnormal{ given by } z \mapsto 
\begin{cases}
    ([y,\xi],p) & \text{if}\ z=[y,p,\xi,0] \\
    ([x,1],p) & \text{if}\ z=[x,p,1]
\end{cases}\]
is a diffeomorphism.
\end{prop}
\begin{proof}
Notice that we have a diffeomorphism
\[\DNC(\base \times \basetwo, \subbase \times \basetwo) \xra{\sim} \DNC(\base,\subbase) \times \basetwo \textnormal{ given by } z \mapsto
\begin{cases}
(y,\xi,p) & \text{if}\ z=(y,p,\xi,0) \\
(x,t,p) & \text{if}\ z=(x,p,t)
\end{cases}\]
(where we used $\DNC(M,M) \cong M \times \mathbb{R}$) by (the proof of) Proposition \ref{prop: if clean intersection, then DNC respects fiber products}. By considering the diagonal $\mathbb{R}^\times$-action on $\DNC(\base,\subbase) \times \basetwo$ (with trivial action on $\basetwo$), it is readily verified that the above map is $\mathbb{R}^\times$-equivariant. Moreover, the map maps $\subbase \times \basetwo \times \mathbb{R}$ onto $\subbase \times \mathbb{R} \times \basetwo$, so, by restricting to $\DNC(\base \times \basetwo, \subbase \times \basetwo) \setminus \subbase \times \basetwo \times \mathbb{R}$, and passing to quotient spaces, we obtain the desired diffeomorphism.
\end{proof}
\begin{rema}\label{rema: blow-up of a point is most important}
The above proposition shows that the blow-up of a point is most important to the theory of blow-ups in the following sense: the blow-up construction is local in nature and relies on the definition of 
\[\blup(\mathbb{R}^n,\mathbb{R}^p),\]
which, by Proposition \ref{prop: blup(X x M, Y x M) = blup(X,Y) x M}, is canonically isomorphic to $\blup(\mathbb{R}^q,0) \times \mathbb{R}^p$ (where $q=n-p$). This shows that the blow-up construction is completely determined by the property stated in Proposition \ref{prop: blup(X x M, Y x M) = blup(X,Y) x M} and the construction of blowing up a point. 
\end{rema}
Lastly, we will show that the deformation to the normal cone construction can be recovered from the blow-up construction.
\begin{prop}\cite{2017arXiv170509588D}\label{prop: DNC out of blup}
The canonical map
\[\DNC(\base,\subbase) \ra \blup(\base \times \mathbb{R}, \subbase \times \{0\}) \textnormal{ given by } z \mapsto \begin{cases}
    [y,\xi,1] & \text{if}\ z=(y,\xi) \\
    [x,t,1] & \text{if}\ z=(x,t)
\end{cases}\]
is an open embedding onto $\blup(\base \times \mathbb{R}, \subbase \times \{0\}) \setminus \blup(\base \times \{0\}, \subbase \times \{0\})$.
\end{prop}
\begin{proof}
%Notice that, indeed, the smooth map
%\[f: \DNC(\base,\subbase) \xra{\DNC(\id_\base)} \DNC(\base,\base) \cong \base \times \mathbb{R}\]
%is a map of pairs
%\[(\DNC(\base,\subbase),\normal(\base,\subbase)) \ra (\base \times \mathbb{R}, \subbase \times \{0\}).\]
%Moreover, $f^{-1}(\subbase \times \{0\}) = \normal(\base,\subbase)$, and $d_\normal f$ is fiberwise injective, which follows immediately by writing $f$ locally in an adapted chart $(U,\varphi)$ (with $V \coloneqq U \cap \subbase$):
%\[\Omega^U_V \ra \varphi(U) \times \mathbb{R} \textnormal{ given by } (y,\xi,t) \mapsto 
%\begin{cases}
%    (y,0,0) & \text{if}\ t=0 \\
%    (y,t\xi,t) & \text{if}\ t\neq0.
%\end{cases} = (y,t\xi,t).\]
%By applying the $\blup$ functor, we obtain the required map. 
Notice that the map is injective. In local coordinates, the map is readily verified to be an immersion, and observe that it maps $\base \times \mathbb{R}^\times$ onto $\base \times \mathbb{R}^\times \subset \base \times \mathbb{R} \setminus \subbase \times \{0\} \subset \blup(\base \times \mathbb{R},\subbase \times \{0\})$, and it maps $\normal(\base,\subbase)$ onto the set
\[\{[y,\xi,v] \in \mathbb{P}(\base \times \mathbb{R}, \subbase \times \{0\}) \mid v \neq 0\} = \mathbb{P}(\base \times \mathbb{R}, \subbase \times \{0\}) \setminus \mathbb{P}(\base \times \{0\}, \subbase \times \{0\}).\] 
We see that the map maps bijectively onto $\blup(\base \times \mathbb{R}, \subbase \times \{0\}) \setminus \blup(\base \times \{0\}, \subbase \times \{0\})$. By Corollary \ref{coro: blup submanifold}, the result follows. 
\end{proof}
\subsection{Blow-up groupoids and algebroids}\label{sec: Blow-up groupoids and algebroids}
In this section we will show that we can naturally blow up Lie groupoids and Lie algebroids. However, in general, we will see that we can not just blow up the total and base spaces and expect the structure maps to lift to the blow-up spaces. The obstruction is explained in Remark \ref{rema: smooth map gives rise to smooth map of Blup}: if $f: (\base,\subbase) \ra (\basetwo,\subbasetwo)$ is a map of pairs, then we only obtain a smooth map $\blup_f(\base,\subbase) \ra \blup(\basetwo,\subbasetwo)$, where
\[\blup_f(\base,\subbase) \coloneqq \DNC_f(\base,\subbase)/\mathbb{R}^\times; \quad \DNC_f(\base,\subbase) \coloneqq \DNC(\base,\subbase) \setminus \DNC(f)^{-1}(\subbasetwo \times \mathbb{R}).\]
However, we will see that, if $(\groupoid,\subgroupoid)$ is a pair of Lie groupoids (see Definition \ref{defn: pair of Lie subgroupoids}), then we can turn 
\[\blup_{\source,\target}(\group,
\subgroup) \coloneqq \blup_{(\target,\source)}(\group,
\subgroup) \rra \blup(\base,\subbase)\] 
into a Lie groupoid, and if $(\algebroid,\subalgebroid)$ is a pair of Lie algebroids (see Definition \ref{defn: pairs of Lie algebroids}), then we can turn 
\[\blup_\pi(\algebr,\subalgebr) \ra \base\] 
(where $\pi$ is the projection map) into a Lie algebroid. 
\begin{rema}\label{rema: notation blup groupoid and algebroid}
Notice that, if $\grouptwoid$ is a Lie groupoid, and $U \subset \basetwo$ is open, then
\[\grouptwo|_{U} \coloneqq \source^{-1}(U) \cap \target^{-1}(U)\]
is again a Lie groupoid (it also works if $U$ is so-called \textit{saturated} instead of open, i.e. a collection of orbits; more on that later), and it is called the \textit{restriction groupoid of $U$}. Therefore, we can also write
\[\blup_{\source,\target}(\group,\subgroup) = \DNC(\group,\subgroup)|_{\DNC(\base,\subbase) \setminus \subbase \times \mathbb{R}}/\mathbb{R}^\times; \quad \blup_{\pi}(\algebr,\subalgebr) = \DNC(\algebr,\subalgebr)|_{\DNC(\base,\subbase) \setminus \subbase \times \mathbb{R}}/\mathbb{R}^\times\]
(see Terminology \ref{term: subscript f for blup and DNC}).
\end{rema}
\subsubsection{Lie groupoid blow-up}
We start with the groupoid case.
\begin{theo}\cite{2017arXiv170509588D}\label{theo: groupoid blow-up}
Let $(\groupoid,\subgroupoid)$ be a pair of Lie groupoids. Then
\[\blup_{\source,\target}(\group,\subgroup) \rra \blup(\base,\subbase)\]
has a unique Lie groupoid structure for which $\group|_{\base \setminus \subbase} \subset \blup_{\source,\target}(\group,\subgroup)$ is a Lie subgroupoid. Moreover, $\blup_{\source,\target}(\group,\subgroup)^{(2)}$ can canonically be identified with $\blup_{\source \times \target}(\group^{(2)},\subgroup^{(2)})$, and the structure maps are obtained by applying the functor $\blup$ to the structure maps of $\group$ .
\end{theo}
\begin{proof}
The structure maps of the blow-up groupoid are simply given by applying the $\blup$-functor (see Lemma \ref{lemm: induced global smooth map on blow up}):
\[(\blup(\source),\blup(\target),\blup(\mult),\blup(\inv),\blup(\identity)).\]
To show that this is well-defined, we do have to show that the above maps restrict to maps with the right (co)domain; for example, we have to show that $\blup(\inv): \blup_\inv(\group,\subgroup) \ra \blup(\group,\subgroup)$ restricts to a map $\blup_{\source,\target}(\group,\subgroup) \ra \blup_{\source,\target}(\group,\subgroup)$. For $\blup(\source)$ and $\blup(\target)$ this is clear. For the other structure maps, we will show that (1)
\begin{align*}
    \blup_{\source,\target}(\group,\subgroup)^{(2)} &\cong \blup_{\source \times \target}(\group^{(2)},\subgroup^{(2)}) \subset \blup_\mult(\group^{(2)},\subgroup^{(2)}); \\
    \blup_{\source,\target}(\group,\subgroup) &\subset \blup_{\inv}(\group,\subgroup); \\
    \blup(\base,\subbase)&=\blup_{\identity}(\base,\subbase),
\end{align*}
where in the first inclusion we used that 
\[\blup_{\source,\target}(\group,\subgroup)^{(2)} \subset \blup_{\source}(\group,\subgroup) \tensor[_{\blup(\source)}]{\times}{_{\blup(\target)}} \blup_{\target}(\group,\subgroup) \cong \blup_{\source \times \target}(\group^{(2)},\subgroup^{(2)})\] 
(see Proposition \ref{prop: blup respects fiber products}), and (2) that $\mult$, $\inv$, and $\identity$ all map into $\blup_{\source,\target}(\group,\subgroup)$ when restricted to $\blup_{\source,\target}(\group,\subgroup)^{(2)}$, $\blup_{\source,\target}(\group,\subgroup)$, and $\blup(\base,\subbase)$, respectively. By definition, we can instead prove these statements by replacing ``$\blup$'' everywhere with ``$\DNC$''. 

To see that $\DNC_{\source,\target}(\group,\subgroup)^{(2)} \cong \DNC_{\source \times \target}(\group^{(2)},\subgroup^{(2)})$, we have to show that, under the diffeomorphism $\DNC(\group,\subgroup)^{(2)} \cong \DNC(\group^{(2)},\subgroup^{(2)})$, 
\[\DNC(\source)^{-1}(\subbase \times \mathbb{R}) \cup \DNC(\target)^{-1}(\subbase \times \mathbb{R}) \cong \DNC(\source \times \target)^{-1}(\subbase \times \subbase \times \mathbb{R}).\]
First, let $(\widetilde g_1,\widetilde g_2) \in \DNC(\group,\subgroup)^{(2)}$. Then this means that $\DNC(\source)(\widetilde g_1)=\DNC(\target)(\widetilde g_2)$, so, by definition of $\DNC(\source)$ and $\DNC(\target)$, we have
\[(\widetilde g_1,\widetilde g_2)=
\begin{cases} 
(h_1,\xi_1,h_2,\xi_2) & \text{with}\ (\source(h_1),d_\normal \source(h_1)\xi_1)=(\target(h_2),d_\normal \target(h_2)\xi_2), \textnormal{ or } \\ (g_1,t,g_2,t) & \text{with}\ \source(g_1)=\target(g_2). \end{cases}\]
We see that $(\widetilde g_1,\widetilde g_2) \in \DNC_{\source,\target}(\group,\subgroup)^{(2)}$ if and only if 
\[\begin{cases} 
(h_1,h_2,\xi_1,\xi_2) & \text{with}\ (\source(h_1),d_\normal \source(h_1)\xi_1)=(\target(h_2),d_\normal \target(h_2)\xi_2), \textnormal{ or } \\ (g_1,g_2,t) & \text{with}\ \source(g_1)=\target(g_2). \end{cases} \in \DNC_{\source \times \target}(\group^{(2)},\subgroup^{(2)}).\]
Moreover, if $(\widetilde g_1,\widetilde g_2) \in \DNC_{\source,\target}(\group,\subgroup)^{(2)}$, then
\[\DNC(\mult)(z) = 
\begin{cases} 
(\mult(h_1,h_2),d_\normal \mult(h_1,h_2)(\xi_1,\xi_2)) & \text{if}\ z=(h_1,h_2,\xi_1,\xi_2) \\
(\mult(g_1,g_2),t) & \text{if}\ z=(g_1,g_2,t)
\end{cases}\]
does not lie in $\subgroup \times \mathbb{R}$, since $\source \circ \mult(g_1,g_2) = \source(g_2) \not\in \subbase$ and $\target \circ \mult(g_1,g_2) = \target(g_1) \not\in \subbase$ (and similarly $d_\normal\source \circ d_\normal\mult(h_1,h_2)(\xi_1,\xi_2) = d_\normal\source(h_2)\xi_2 \not\in 0_\subbase$ and $d_\normal\target \circ d_\normal\mult(h_1,h_2)(\xi_1,\xi_2) = d_\normal\target(h_1)\xi_1 \not\in 0_\subbase$). This shows that $\blup_{\source,\target}(\group,\subgroup)^{(2)} \cong \blup_{\source \times \target}(\group^{(2)},\subgroup^{(2)}) \subset \blup_{\mult}(\group^{(2)},\subgroup^{(2)})$, but also that $\blup(\mult)$ restricts to a map
\[\blup_{\source,\target}(\group,\subgroup)^{(2)} \ra \blup_{\source,\target}(\group,\subgroup).\]
The inclusion $\blup_{\source,\target}(\group,\subgroup) \subset \blup_{\inv}(\group,\subgroup)$ follows similarly, but from $\source \circ \inv = \target$ and $\target \circ \inv = \source$, and it also follows from this that $\blup(\inv)$ restricts to a map $\blup_{\source,\target}(\group,\subgroup) \ra \blup_{\source,\target}(\group,\subgroup)$. Lastly, the equality $\blup(\base,\subbase) = \blup_{\identity}(\base,\subbase)$ follows from the fact that $\DNC(\identity)(z) \in \subgroup \times \mathbb{R}$ if and only if $z \in \subbase \times \mathbb{R} \subset \DNC(\base,\subbase)$, and from $\source \circ \identity = \id_\base$ and $\target \circ \identity = \id_\base$, we see that $\blup(\identity)$ is a map $\blup(\base,\subbase) \ra \blup_{\source,\target}(\group,\subgroup)$.

Since $\blup(\source)$ and $\blup(\target)$ are submersions (see Proposition \ref{prop: analogue constant rank when applying blup}) and all other structure maps are smooth (see Lemma \ref{lemm: induced global smooth map on blow up}), it remains to check that we obtain a groupoid this way. %To see that for all $(g,h) \in \blup_{\source,\target}(\group,\subgroup)^{(2)}$, 
%\[\blup(\mult)(g,h) \in \blup_{\source,\target}(\group,\subgroup),\]
%let $g_\DNC,h_\DNC \in U_{\source,\target}$ such that $\pi_\group(g_\DNC) = g$ \& $\pi_\group(h_\DNC) = h$ (where $\pi_\group: U_{\source,\target} \ra U_{\source,\target}/\mathbb{R}^\times = \blup_{\source,\target}(\group,\subgroup)$ is the canonical projection). Then since $\blup(\source)(g) = \blup(\target)(h)$ we have
%\[\pi_M \circ \DNC(\source)(g_\DNC) = \blup(\source) \circ \pi_\group(g_\DNC) = \blup(\target) \circ \pi_\group(h_\DNC) = \pi_M \circ \DNC(\target)(h_\DNC),\]
%(where $\pi_M: M \setminus (S \times \mathbb{R}) \ra \blup(M,S)$ is the canonical projection) by construction. It follows that there is $\lambda \in \mathbb{R}^\times$ for which $\DNC(\source)(g_\DNC) = \lambda \cdot \DNC(\target)(h_\DNC) = \DNC(\target)(\lambda \cdot h_\DNC)$, so that
%\[\blup(\mult)(g,h) = \pi_\group \circ \DNC(\mult)(g_\DNC,\lambda \cdot h_\DNC).\]
%Since
%\begin{align*}
%    \DNC(\source) \circ \DNC(\mult)(g_\DNC,\lambda \cdot h_\DNC) &= \DNC(\source)(\lambda \cdot h_\DNC) = \lambda \cdot \DNC(\source)(h_\DNC) \quad \& \\
%    \DNC(\target) \circ \DNC(\mult)(g_\DNC,\lambda \cdot h_\DNC) &= \DNC(\target)(g_\DNC),
%\end{align*}
%we see that $\DNC(\mult)(g_\DNC,\lambda \cdot h_\DNC) \in U_{\source,\target}$, so that indeed $\blup(\mult)(g,h) \in \blup_{\source,\target}(\group,\subgroup)$. It remains 
But for this we only have to check that, for all $g \in \blup_{\source,\target}(\group,\subgroup)$, $\blup(\inv)(g)$ is the inverse of $g$. This simply follows by $\mathbb{R}^\times$-equivariance of the maps $\DNC(\mult), \DNC(\id_\group,\inv)$ \& $\DNC(\identity)$ and the expression
\[\DNC(\mult) \circ \DNC(\id_\group,\inv) = \DNC(\identity)\]
which therefore descends to the expression 
\[\blup(\mult) \circ \blup(\id_\group,\inv) = \blup(\mult) \circ (\id_\group,\blup(\inv)) = \blup(\identity).\]
We have now shown that indeed we naturally obtain a blow-up groupoid from a pair of Lie groupoids $(\group,\subgroup)$. Moreover, notice that 
\[(\group \setminus \subgroup) \cap 
\blup_{\source,\target}(\group,\subgroup) = \group \setminus (\source^{-1}(\subbase) \cup \target^{-1}(\subbase)) = \group|_{\base \setminus \subbase}\]
is closed under multiplication (by definition of $\blup(\mult)$), so it defines a subgroupoid (over $\base \setminus \subbase$). The uniqueness statement is a consequence of Lemma \ref{lemm: induced global smooth map on blow up}, since each of the structure maps of $\group|_{\base \setminus \subbase}$ extends to at most one smooth map of $\blup_{\source,\target}(\group,\subgroup)$. 
This proves the statement.
\end{proof}
Notice that $\blup_{\source,\target}(\group,\subgroup)$ contains the subgroupoid
\[\mathbb{P}_{\source,\target}(\group,\subgroup) \rra \mathbb{P}(\base,\subbase),\]
where 
\[\mathbb{P}_{\source,\target}(\group,\subgroup) = \normal(\group,\subgroup)|_{\normal(\base,\subbase) \setminus 0_\subbase}/\mathbb{R}^\times.\]
Observe that $\blup(\source)$ and $\blup(\target)$ are maps of pairs
\[(\blup_{\source,\target}(\group,\subgroup),\mathbb{P}_{\source,\target}(\group,\subgroup)) \ra (\blup(\base,\subbase),\mathbb{P}(\base,\subbase))\]
such that $\blup(\source)^{-1}(\mathbb{P}(\base,\subbase)) = \mathbb{P}_{\source,\target}(\group,\subgroup) = \blup(\target)^{-1}(\mathbb{P}(\base,\subbase))$ (i.e. $\mathbb{P}(\base,\subbase) \subset \blup(\base,\subbase)$ is a so-called saturated submanifold). Indeed, this is by definition of $\DNC(\source)$ and $\DNC(\target)$, which are maps of pairs 
\[(\DNC(\group,\subgroup), \normal(\group,\subgroup)) \ra (\DNC(\base,\subbase),\normal(\base,\subbase)),\]
for which, similarly, $\DNC(\source)^{-1}(\normal(\base,\subbase)) = \normal(\group,\subgroup) = \DNC(\target)^{-1}(\normal(\base,\subbase))$.
From this we obtain the following result.
\begin{prop}\label{prop: isotropy groups and orbits of blow-up}
Let $(\groupoid,\subgroupoid)$ be a pair of groupoids. Then $(\blup_{\source,\target}(\group,\subgroup),\mathbb{P}_{\source,\target}(\group,\subgroup))$ is a pair of groupoids. Moreover, for all $z \in \blup(\base,\subbase)$, 
\begin{enumerate}
    \item the isotropy group $\blup_{\source,\target}(\group,\subgroup)_z$ equals the isotropy group $(\group|_{\base \setminus \subbase})_z$ if $z \in \base \setminus \subbase$, and it equals the isotropy group $\mathbb{P}_{\source,\target}(\group,\subgroup)_z$ if $z \in \mathbb{P}(\base,\subbase)$;
    \item the orbit $\orbit_z$ equals the orbit $\orbit_z \subset \base \setminus \subbase$ of $\group|_{\base \setminus \subbase}$ if $z \in \base \setminus \subbase$, and it equals the orbit $\orbit_z \subset \mathbb{P}(\base,\subbase)$ of $\mathbb{P}_{\source,\target}(\group,\subgroup)$ if $z \in \mathbb{P}(\base,\subbase)$.
\end{enumerate}\vspace{-\baselineskip}
\end{prop}
%\begin{term}\label{term: subscrupt for G setminus H}
%If $f: (\base,\subbase) \ra (\basetwo,\subbasetwo)$ is a map of pairs, then we will write $(\base \setminus \subbase)_{f}$ for $\base \setminus \subbase \cap \blup_f(\base,\subbase) = \base \setminus f^{-1}(\subbasetwo)$.
%\end{term}
%\begin{rema}\label{rema: Hausdorffness uniqueness blup groupoid}
%If $\group$ is Hausdorff, then the above Lie groupoid structure on $\blup_{\source,\target}(\group,\subgroup) \rra \blup(\base,\subbase)$ is unique with the property that $\group|_{\base \setminus \subbase} \subset \blup_{\source,\target}(\group,\subgroup)$ is a subgroupoid, since then all the structure maps of $\group|_{\base \setminus \subbase}$ have at most one extension to all of $\blup_{\source,\target}(\group,\subgroup)$ (since, as a subspace, $\group|_{\base \setminus \subbase} \subset \blup_{\source,\target}(\group,\subgroup)$ is dense).
%\end{rema}
It might be instructive to present a simple example.
\begin{exam}\label{exam: blow-up point in R^2}
Consider the pair groupoid $\mathbb{R} \times \mathbb{R} \rra \mathbb{R}$. Notice that the blow-up construction will only yield something interesting if we blow-up a point in $\mathbb{R}$, and a codimension $2$ submanifold in $\mathbb{R} \times \mathbb{R}$. For example, consider the closed subgroupoid $\{(0,0)\} \subset \mathbb{R} \times \mathbb{R}$ (over $\{0\} \subset \mathbb{R}$). First of all, we would like to understand 
\[\blup_{\source,\target}(\mathbb{R} \times \mathbb{R}, (0,0))\]
as a manifold. For this, notice that
\[\DNC(\mathbb{R} \times \mathbb{R}, (0,0)) \cong \mathbb{R}^3; \quad \DNC(\mathbb{R},0) \cong \mathbb{R}^2\]
and that, under these diffeomorphisms, $\DNC(\source)$ and $\DNC(\target)$ become the projection maps $\mathbb{R}^3 \ra \mathbb{R}^2$ onto the last two factors, and first and last factor, respectively. In particular, we see that
\[\DNC(\source)^{-1}(\{0\} \times \mathbb{R}) \cup \DNC(\target)^{-1}(\{0\} \times \mathbb{R}) = (\mathbb{R} \times \{0\} \times \mathbb{R}) \cup (\{0\} \times \mathbb{R} \times \mathbb{R}),\]
so that, explicitly, we have a diffeomorphism
\[\DNC_{\source,\target}(\mathbb{R} \times \mathbb{R}, (0,0)) \xra{\sim} \mathbb{R}^\times \times \mathbb{R}^\times \times \mathbb{R} \textnormal{ given by } z \mapsto
\begin{cases}
    (\xi_1,\xi_2,0) & \text{if}\ z=((0,0),\xi_1,\xi_2) \\
    (t^{-1}a,t^{-1}b,t) & \text{if}\ z=(a,b,t)
\end{cases}\]
Recall that, under this diffeomorphism, the $\mathbb{R}^\times$-action of $\DNC_{\source,\target}(\mathbb{R} \times \mathbb{R}, (0,0))$ transfers to $\mathbb{R}^\times \times \mathbb{R}^\times \times \mathbb{R}$ and is given by the $\mathbb{R}^\times$-action 
\[\lambda \cdot (a,b,t) \coloneqq (\lambda^{-1}a,\lambda^{-1}b,\lambda t)\]
(see e.g. Remark \ref{rema: local coordinates for blowup}),
so we see that 
\[\blup_{\source,\target}(\mathbb{R} \times \mathbb{R},(0,0)) \cong (\mathbb{R}^\times \times \mathbb{R}^\times \times \mathbb{R})/\mathbb{R}^\times \cong \mathbb{R}^\times \times \mathbb{R},\]
where the last diffeomorphism is given by e.g. $[\lambda,\mu,t] \mapsto (\tfrac{\lambda}{\mu},\mu \cdot t)$ (with inverse map given by $(\lambda,a) \mapsto [\lambda \cdot a,a,1]$). By tracing the diffeomorphisms, we see that, under this diffeomorphism, $\blup(\source)$ and $\blup(\target)$ become the maps
\[\mathbb{R}^\times \times \mathbb{R} \ra \mathbb{R} \textnormal{ given by } (\lambda,a) \mapsto a \textnormal{ and } \mathbb{R}^\times \times \mathbb{R} \ra \mathbb{R} \textnormal{ given by } (\lambda,a) \mapsto \lambda \cdot a,\]
respectively. It is now readily verified (by similarly inspecting the other structure maps as well) that $\blup_{\source,\target}(\mathbb{R} \times \mathbb{R}, (0,0)) \rra \blup(\mathbb{R},0)$ is isomorphic to the action groupoid $\mathbb{R}^\times \ltimes \mathbb{R}$. 

Two simple to verify, but interesting things happen: the isotropy Lie groups of the groupoid $\group \coloneqq \blup_{\source,\target}(\mathbb{R} \times \mathbb{R},(0,0))$ are $\group_x = \{x\}$ if $x \in \mathbb{R} \setminus \{0\}$, and $\group_0 = \mathbb{R}^\times$, and, moreover, the orbits are $\mathbb{R}^\times$ and $\{0\}$. Since the isotropy groups of pair groupoids are trivial, and pair groupoids are transitive, we see that these spaces already change quite a lot in simple cases. %We also see from this that blow-ups of pairs of proper groupoids $(\groupoid, \subgroupoid)$ (meaning the maps $(\source_\group,\target_\group)$ and $(\source_\subgroup,\target_{\subgroup})$ are proper maps) are not necessarily proper, and blow-ups of pairs of transitive groupoids $(\groupoid, \subgroupoid)$ (meaning $\group$ and $\subgroup$ are transitive groupoids) are not transitive.
\end{exam}
The above example is an example of a blow-up of a pair of transitive groupoids $(\group,\subgroup)$ (meaning $\group$ and $\subgroup$ are transitive groupoids) which is not transitive. However, from the following result we see that this is (almost) always the case.
\begin{rema}\label{rema: orbits G|X setminus Y}
Let $z \in \base \setminus \subbase$. Notice that the isotropy Lie group $\group_z$ of the Lie groupoid $\group|_{\base \setminus \subbase}$ is in fact equal (or canonically isomorphic) to the isotropy Lie group $\group_z$ of $\group$. Indeed, since $\source^{-1}(z) \subset \group \setminus \source^{-1}(\subbase)$ and $\target^{-1}(z) \subset \group \setminus \target^{-1}(\subbase)$, it follows that 
\[\group_z = \source^{-1}(z) \cap \target^{-1}(z) \subset \group \setminus (\source^{-1}(\subbase) \cup \target^{-1}(\subbase)) = \group|_{\base \setminus \subbase}.\]
The orbit of $z$ in $\group|_{\base \setminus \subbase}$ is given by $\orbit_z \setminus \subbase$. To see this, observe that 
\[\source_{\group|_{\base \setminus \subbase}}^{-1}(z) = \source^{-1}(z) \setminus \target^{-1}(\subbase).\]
Since $\target: \source^{-1}(z) \ra \orbit_z$ is surjective, it follows that 
\[\orbit_z \setminus \subbase = \target(\source^{-1}(z)) \setminus \target(\target^{-1}(\subbase)) \subset \target(\source^{-1}(z) \setminus \target^{-1}(\subbase)).\]
The other inclusion follows simply from the fact that if $g \in \source^{-1}(z) \setminus \target^{-1}(\subbase)$, then $\target(g) \in \orbit_z$, but $\target(g) \not\in \subbase$ by assumption.
\end{rema}
%In fact, we even have the following characterisation of the orbits of a blow-up groupoid.
%\begin{coro}\label{coro: orbits of blow-up}
%Let $z \in \blup(\base,\subbase)$. Then the orbit of $z$ in $\blup_{\source,\target}(\group,\subgroup) \rra \blup(\base,\subbase)$ is equal to $\overoverline{\orbit_{\bldown(z)}} \cap (\base \setminus \subbase)$ if $z \in \base \setminus\ \subbase$, and $\overoverline{\orbit_{\bldown(z)}} \cap \mathbb{P}(\base,\subbase)$ if $z \in \mathbb{P}(\base,\subbase)$.
%\end{coro}
%\begin{proof}
%If $z \in \base \setminus \subbase$, the statement follows by Remark \ref{rema: orbits G|X setminus Y} (note: $\overoverline{\orbit_{\bldown(z)}} \cap (\base \setminus \subbase) = \orbit_{\bldown(z)} \setminus \subbase$). 
%\end{proof}
With an eye on the Lie algebroid blow-up construction, we also present the following more general example (in comparison with Example \ref{exam: blow-up point in R^2}).
\begin{exam}\label{exam: blup pair groupoids}
Here, we will derive some properties of the blow-up of the pair groupoid:
\[\blup_{\source,\target}(\base \times \base, \subbase \times \subbase) \rra \blup(\base,\subbase)\]
In this discussion, we will also see that
\[\blup_{\pi_{T\base}}(T\base,T\subbase) \rra \blup(\base,\subbase)\]
is its Lie algebroid (we have not discussed the Lie algebroid structure on it, but we will see that, at least as manifolds, we can view it as the Lie algebroid).

Recall that $\DNC(T\base,T\subbase) \ra \DNC(\base,\subbase)$ is the Lie algebroid of $\DNC(\base \times \base, \subbase 
\times \subbase) \rra \DNC(\base,\subbase)$ (see Proposition \ref{prop: algebroid of DNC of groupoid}). Observe that, by (the proof of) Proposition \ref{prop: if clean intersection, then DNC respects fiber products}, we have a canonical diffeomorphism
\[(\DNC(\pr_1),\DNC(\pr_2)): \DNC(\base \times \base,\subbase \times \subbase) \ra \DNC(\base,\subbase) \times_{\mathbb{R}} \DNC(\base,\subbase)\]
which extends to an isomorphism of groupoids, where we view $\DNC(\base,\subbase) \times_{\mathbb{R}} \DNC(\base,\subbase)$ as a subgroupoid of the pair groupoid. This is readily verified simply by inspecting the definition of $\DNC(\pr_1)$ and $\DNC(\pr_2)$. Moreover, this morphism intertwines the $\mathbb{R}^\times$-actions (where we take the diagonal action on $\DNC(\base,\subbase) \times_{\mathbb{R}} \DNC(\base,\subbase)$). Under this identification, the open subgroupoid 
\[\DNC_{\source,\target}(\base \times \base,\subbase \times \subbase) \subset \DNC(\base \times \base,\subbase \times \subbase)\]
maps onto the open subgroupoid
\[P \times_{\mathbb{R}} P \subset \DNC(\base,\subbase) \times_{\mathbb{R}} \DNC(\base,\subbase)\]
where we wrote $P \coloneqq \DNC(\base,\subbase) \setminus \subbase \times \mathbb{R}$. Therefore, we obtain an isomorphism of groupoids
\[\blup_{\source,\target}(\base \times \base,\subbase \times \subbase) \ra (P \times_{\mathbb{R}} P)/\mathbb{R}^\times;\]
that is, we can view $\blup_{\source,\target}(\base \times \base,\subbase \times \subbase)$ as a closed subgroupoid of the Gauge groupoid 
\[P \otimes_{\mathbb{R}^\times} P = (P \times P)/\mathbb{R}^\times \rra P/\mathbb{R}^\times (= \blup(\base,\subbase)).\]
Now, by Remark \ref{rema: Lie algebroid integration of product and pullback}, the Lie algebroid 
\[\DNC(T\base,T\subbase)|_{\DNC(\base,\subbase) \setminus \subbase \times \mathbb{R}} = \DNC_{\pi_{T\base}}(T\base,T\subbase)\] is the Lie algebroid of $\DNC_{\source,\target}(\base \times \base,\subbase \times \subbase)$. By using that the induced Lie algebroid of a Lie groupoid $\grouptwoid$ can be identified with $\normal(\grouptwo,\basetwo) \ra \basetwo$, so that, as Lie algebroids, 
\[\DNC_{\pi_{T\base}}(T\base,T\subbase) \cong \normal(\DNC_{\source,\target}(\base \times \base, \subbase \times \subbase),\DNC(\base,\subbase) \setminus \subbase \times \mathbb{R}),\]
and using Proposition \ref{prop: induced G action on normal bundle}, so that
\begin{align*}
    \blup_{\pi_{T\subbase}}(T\base,T\subbase) &= \DNC_{\pi_{T\base}}(T\base,T\subbase)/\mathbb{R}^\times \\
    &\cong \normal(\DNC_{\source,\target}(\base \times \base, \subbase \times \subbase),\DNC(\base,\subbase) \setminus \subbase \times \mathbb{R})/\mathbb{R}^\times \\
    &\cong \normal(\DNC_{\source,\target}(\base \times \base, \subbase \times \subbase)/\mathbb{R}^\times,(\DNC(\base,\subbase) \setminus \subbase \times \mathbb{R})/\mathbb{R}^\times) \\
    &= \normal(\blup_{\source,\target}(\base \times \base, \subbase \times \subbase), \blup(\base,\subbase)) 
\end{align*}
we see that, as claimed, we can view the Lie algebroid of $\blup_{\source,\target}(\base \times \base,\subbase \times \subbase)$ as $\blup_{\pi_{T\base}}(T\base,T\subbase)$. Moreover, we can also conclude from the above discussion that it is isomorphic to the Lie algebroid of $(P \times_{\mathbb{R}} P)/\mathbb{R}^\times$, 
which is 
\[(\ker d\hat{t}|_{P})/\mathbb{R}^\times \subset TP/\mathbb{R}^\times\] 
(see also Proposition \ref{prop: DNC of tangent spaces}). Indeed, if we view the Lie algebroid of $P \times P$ as $\ker d\pr_2|_{P}$, then we see that $(p,\xi_1,p,\xi_2) \in TP \times_{T\mathbb{R}} TP \cong T(P \times_{\mathbb{R}} P)$ if and only if $\xi_2=0$, but also note that $(p,\xi_1,p,0) \in TP \times_{T\mathbb{R}} TP$ if and only if $d\hat{t}(p)\xi_1=0$. 
%To explicitly define the anchor map of this Lie algebroid, note that it is given by the induced map on Lie algebroids of the morphism of groupoids
%\[(P \times_{\mathbb{R}} P)/\mathbb{R}^\times \hookrightarrow (P \times P)/\mathbb{R}^\times \xra{(\target,\source)} P/\mathbb{R}^\times \times P/\mathbb{R}^\times \textnormal{ given by } [p,q] \mapsto ([p],[q]).\]

We will now describe the isotropy groups and orbits of $\blup_{\source,\target}(\base \times \base, \subbase \times \subbase) \cong (P \times_{\mathbb{R}} P)/\mathbb{R}^\times$. For $\group \coloneqq (P \times_{\mathbb{R}} P)/\mathbb{R}^\times$, observe that, if $[x,1] \in \base \setminus \subbase \subset \blup(\base,\subbase)$, then
\[[x_1,t_1,x_2,t_2] \in \group_{[x,1]} = \pr_2^{-1}([x,1]) \cap \pr_1^{-1}([x,1]) \cap (P \times_{\mathbb{R}} P)/\mathbb{R}^\times\]
if and only if $x=x_1=x_2$ and $t_1=t_2$, and, if $[y,\xi] \in \mathbb{P}(\base,\subbase) \subset \blup(\base,\subbase)$, then
\[[y_1,\xi_1,y_2,\xi_2] \in \group_{[y,\xi]}\]
if and only if $y=y_1=y_2$ and $\xi_1,\xi_2 \in \mathbb{R}^\times \cdot \xi$. This shows that $\group$ has the following isotropy groups:
\[\group_{[x,1]} = \{[x,1,x,1]\} \textnormal{ } ([x,1] \in \base \setminus \subbase); \quad \group_{[y,\xi]} = \{[y,\xi,y,\xi'] \mid \xi' \in \mathbb{R}^\times\} \cong \mathbb{R}^\times \textnormal{ } ([y,\xi] \in \mathbb{P}(\base,\subbase))\]
(compare this with Example \ref{exam: blow-up point in R^2}). The orbits can be computed similarly: if $[x,1] \in \base \setminus \subbase$, then
\[[x_1,t_1,x_2,t_2] \in \pr_2^{-1}([x,1])\]
if and only if $x=x_2$ and $t_1=t_2$, 
and if $[y,\xi] \in \mathbb{P}(\base,\subbase)$, then
\[[y_1,\xi_1,y_2,\xi_2] \in \pr_2^{-1}([y,\xi])\]
if and only if $y=y_2$ and $\xi_2 \in \mathbb{R}^\times \cdot \xi$. This shows that the orbits are 
\[\orbit_{[x,1]} = \pr_1(\pr_2^{-1}([x,1])) = \base \setminus \subbase; \quad \orbit_{[y,\xi]} = \mathbb{P}(\base,\subbase)\]
(also compare this with Example \ref{exam: blow-up point in R^2}).
\end{exam}

Observe that, by functoriality of the blow-up construction, we obtain the following statement:
\begin{prop}\label{prop: blow-up functoriality groupoids}
Let $F: (\group,\subgroup) \ra (\grouptwo,\subgrouptwo)$ be a morphism of pairs of Lie groupoids between pairs of Lie groupoids $(\groupoid,\subgroupoid)$ and $(\grouptwoid,\subgrouptwoid)$. If the base map $f$ of $F$ satisfies $f^{-1}(\subbasetwo)=\subbase$ and $d_\normal f$ is fiberwise injective, then $\blup(F)$ restricts to a morphism of Lie groupoids $\blup_{\source_\group,\target_\group}(\group,\subgroup) \ra \blup_{\source_\grouptwo,\target_\grouptwo}(\grouptwo,\subgrouptwo)$.
\end{prop}
\begin{proof}
We only have to show that $\blup(F)$ restricts to a map $\blup_{\source,\target}(\group,\subgroup) \ra \blup_{\source,\target}(\grouptwo,\subgrouptwo)$; the rest of the statement follows by functoriality (see e.g. Proposition \ref{prop: groupoid morphism DNC}). But this follows from the identities
\[\DNC(\source_\grouptwo) \circ \DNC(F) = \DNC(\source_\grouptwo \circ F) = \DNC(f) \circ \DNC(\source_\group); \quad \DNC(\target_\grouptwo) \circ \DNC(F) = \DNC(\target_\grouptwo \circ F) = \DNC(f) \circ \DNC(\target_\group).\]
Indeed, it proves that, if $\DNC(\source_{\group})(z) \not\in \subbase \times \mathbb{R}$, then $\DNC(f) \circ \DNC(\source_{\group})(z) \not\in \subbase \times \mathbb{R}$ (by our assumption on $f$), so $\DNC(F)(z) \not\in \subgrouptwo \times \mathbb{R}$, otherwise $\DNC(\source_{\grouptwo}) \circ \DNC(F)(z) = \DNC(f) \circ \DNC(\source_{\group})(z) \in \subbasetwo \times \mathbb{R}$. Similarly, this holds if we replace ``$\source$'' everywhere with ``$\target$''. Therefore, $\blup_{\source_{\group},\target_{\group}}(\group,\subgroup) \subset \blup_F(\group,\subgroup)$, but the same argument shows that $\blup(F)$ restricted to $\blup_{\source_{\group},\target_{\group}}(\group,\subgroup)$ maps into $\blup_{\source_\grouptwo,\target_\grouptwo}(\grouptwo,\subgrouptwo)$.
This proves the statement.
\end{proof}
\begin{rema}\label{rema: blow-up subgroupoids}
In particular, if $\iota: (\grouptwo,\subgrouptwo) \hookrightarrow (\group,\subgroup)$ is a pair of Lie subgroupoids such that $\subgrouptwo = \grouptwo \cap \subgroup$, then, from Corollary \ref{coro: blup submanifold}, we see that $\blup(\iota): \blup_{\source_\grouptwo,\target_\grouptwo}(\grouptwo,\subgrouptwo) \ra \blup_{\source_\group,\target_\group}(\group,\subgroup)$ is a Lie subgroupoid (and it is open/closed if $\iota$ is open/closed).
\end{rema}

%\subsubsection{Lie algebroid blow-up}\label{sec: sub blup lie algebroid}
We will now describe the blow-up of Lie algebroids. The way we will approach this, is by describing the blow-up construction for vector bundles first, and extend it to the blow-up construction for Lie algebroids afterwards.

\subsubsection{Vector bundle blow-up}\label{sec: subsub blup vector bundle}
We fix here a pair of vector bundles $(E \xra{\pi} \base,F \xra{\pi|_F} \subbase)$ (see Definition \ref{defn: pairs of Lie algebroids}).

\begin{theo}\cite{2017arXiv170509588D}\label{theo: vector bundle blow-up}
We can equip
\[\blup(\pi): \blup_\pi(E,F) \ra \blup(\base,\subbase)\]
with a unique vector bundle structure such that $E|_{\base \setminus \subbase} \subset \blup_\pi(E,F)$ is a vector subbundle. 
\end{theo}
\begin{proof}
%Write $\normal_\pi(E,F) \coloneqq \normal(E,F)|_{\normal(\base,\subbase) \setminus 0_\subbase}$. Moreover, write
%\[\mathbb{P}(\base,\subbase) \coloneqq \mathbb{P}(\normal(\base,\subbase)) \subset \blup(\base,\subbase), \textnormal{ and } \mathbb{P}_\pi(E,F) \coloneqq \mathbb{P}(\normal_\pi(E,F)) \subset \blup_\pi(E,F).\] 
We will first show that the fibers of
\[\blup(\pi)|_{\mathbb{P}_\pi(E,F)}: \mathbb{P}_\pi(E,F) \ra \mathbb{P}(\base,\subbase)\] 
naturally inherit vector space structures from the fibers of $\normal_\pi(E,F) \ra \normal(\base,\subbase) \setminus 0_\subbase$. Here,
\[\normal_\pi(E,F) = \normal(E,F)|_{\normal(\base,\subbase) \setminus 0_\subbase}\] 
with respect to the vector bundle structure described in Proposition \ref{prop: normal bundle of vector bundle} (see also Example \ref{exam: tangent DVB}), and 
\[\mathbb{P}_\pi(E,F) = \normal_\pi(E,F)/\mathbb{R}^\times.\]
Then, we will explain how we can extract vector bundle charts for $\blup_\pi(E,F)$ from $E \ra \base$.

Let $[y,\xi] \in \mathbb{P}(\base,\subbase)$. To see that the fiber $\mathbb{P}_\pi(E,F)_{[y,\xi]}$ over $[y,\xi]$ %of $\mathbb{P}_\pi(E,F)$ 
inherits a vector space structure from $\normal_\pi(E,F)_{(y,\xi)}$, notice that, automatically, 
\[[y,\xi,\eta] \in \mathbb{P}_\pi(E,F) \textnormal{ whenever } \eta \in \normal_\pi(E,F)_{(y,\xi)},\] 
because $d_\normal\pi(y,\xi,\eta)=(y,\xi) \not\in 0_Y \subset \normal(\base,\subbase)$. This shows that we have a natural map
\[\normal_\pi(E,F)_{(y,\xi)} \ra \mathbb{P}_\pi(E,F)_{[y,\xi]} \textnormal{ given by } (y,\xi,\eta) \mapsto [y,\xi,\eta].\]
It has an obvious inverse; to see that this inverse is well-defined, notice that if $[y,\xi,\eta]=[y,\xi,\zeta]$, then, since $\xi\neq0$, $\eta=\zeta$ by definition of the $\mathbb{R}^\times$-action on $\normal_\pi(E,F)$. So, along the above bijection, we can transfer the vector space structure on $\normal_\pi(E,F)_{(y,\xi)}$ to $\mathbb{P}_\pi(E,F)_{[y,\xi]}$. 
Observe now that 
\[(E \setminus F)_\pi = E \setminus F \cap \blup_\pi(E,F) = E \setminus F \cap E|_{\base \setminus \subbase} = E|_{\base \setminus \subbase},\]
because $E|_{\base \setminus \subbase} \subset E \setminus F$.
So, the fibers of $x \in \base \setminus Y$ is just the vector space $E_x$ viewed as a subspace of $E|_{\base \setminus \subbase} \subset \blup_\pi(E,F)$. We will now show that the blow-up coordinates of $\blup(E,F)$ become vector bundle charts upon restriction to $\blup_\pi(E,F)$. It is then clear that $E|_{\base \setminus \subbase} \subset \blup_{\pi}(E,F)$ is a vector subbundle, and since it is dense, this property uniquely determines the vector bundle structure of $\blup_\pi(E,F)$. 

Let 
\[(y,x): U \ra \mathbb{R}^{p+q}\] 
be a chart of $\base$, such that $V \coloneqq \subbase \cap U = \{v \in U \mid x(v)=0\}$, and let 
\[\psi=(y,x,f,e): E|_U \ra \mathbb{R}^{p+q} \times \mathbb{R}^{k+\ell}\]
be a vector bundle chart of $E$, such that $F \cap (E|_U) = F|_{V} = \{(v,\nu) \in E|_{V} \mid (x(v),e(\nu))=0\}$. Then,
\[d_\normal\psi=(y,dx,f,de): \normal(E|_U,F|_{V}) \ra \mathbb{R}^{p+q} \times \mathbb{R}^{k+\ell}\]
is a vector bundle chart of $\normal(E,F)$ with respect to both the vector bundle structure $\normal(E,F) \ra \normal(\base,\subbase)$, and the vector bundle structure $\normal(E,F) \ra F$ (see also Example \ref{exam: core of TE and normal bundle} and Remark \ref{rema: local nature DVB})%(namely, $(y,b,dx,da): \normal(E|_U,F|_{V}) \ra \mathbb{R}^{p+k} \times \mathbb{R}^{q+\ell}$)
. For all $1 \le r \le q$, write the induced chart of $\blup(\base,\subbase)$ by 
\[(y,\widetilde{x}_r): U_r \ra \mathbb{R}^n,\] 
that is,
\begin{equation}\label{eq: blup U_r}
    U_r = \{[v,\xi] \in \mathbb{P}(\base,\subbase) \mid v \in V, dx^r(v)\xi \neq 0\} \cup \{[u,1] \in \blup(\base,\subbase) \mid u \in U \setminus V, x^r(u) \neq 0\};
\end{equation}
\begin{equation}\label{eq: blup x_r^j}
\widetilde{x}_r^j(z) =
\begin{cases}
    \frac{dx^j(v)\xi}{dx^r(v)\xi} & \text{if}\ z=[v,\xi] \\
    \frac{x^j(u)}{x^r(u)} & \text{if}\ z=[u,1]
\end{cases} \textnormal{ if } j\neq r, \textnormal{ and } \widetilde{x}_r^r(z) =
\begin{cases}
    0 & \text{if}\ z=[v,\xi] \\
    x^r(u) & \text{if}\ z=[u,1]
\end{cases} 
\end{equation}
Similarly, on $\blup(E,F)$ we have, for all $1 \le r \le q+\ell$, an induced chart given by 
\[(y,\widetilde{x}_r,f,\widetilde{e}_r): V_r \ra \mathbb{R}^{p+q} \times \mathbb{R}^{k+\ell}.\]
But, on $\blup_{\pi}(E,F)$, the $V_r$ with $1 \le r \le q$ already constitute a cover. In fact, we have
\begin{equation}\label{eq: V_r constitute a cover}
    \blup_\pi(E,F) \cap U = \bigcup_{1 \le r \le q} \blup(E,F) \cap V_r.
\end{equation}
Indeed, let $z \in \DNC(E,F)$% with $\DNC(\pi)(z) \not\in \subbase \times \mathbb{R}$
. If 
\[z=(v,\xi,\eta) \in \normal(E,F), \textnormal{ where } \eta \in \normal(E,F)_{(v,\xi)},\] 
then $\DNC(\pi)(z) \not\in \subbase \times \mathbb{R}$ if and only if $\xi \neq 0$, i.e. $dx^r(v)\xi \neq 0$ for some $1 \le r \le q$. If 
\[z=(u,\upsilon,1), \textnormal{ where } \upsilon \in E_u,\] 
then $\DNC(\pi)(z) \not\in \subbase \times \mathbb{R}$ if and only if $u \not\in \subbase$, i.e. $x^r(u) \neq 0$ for some $1 \le r \le q$. This proves that \eqref{eq: V_r constitute a cover} holds.

Now, we will show that, for all $1 \le r \le q$, the induced chart 
\[(y,\widetilde{x}_r,f,\widetilde{e}_r): V_r \ra \mathbb{R}^{p+q} \times \mathbb{R}^{k+\ell}\] 
are vector bundle coordinates on $\blup_{\pi}(E,F)$. Recall that here  
\begin{equation}\label{eq: blup V_r}
    V_r = \{[v,\xi,\eta] \in \mathbb{P}_\pi(E,F) \mid v \in V, dx^r(v)\xi \neq 0\}
    \cup \{[u,\upsilon,1] \in \blup_\pi(E,F) \mid u \in U \setminus V, x^r(u) \neq 0\};
\end{equation}
%\begin{align*}
%\widetilde{x}_r^j(z) &=
%\begin{cases}
%    \frac{dx^j(v)\xi}{dx^r(v)\xi} & \text{if}\ z=[v,\xi,\eta] \\
%    \frac{x^j(u)}{x^r(u)} & \text{if}\ z=[u,\upsilon,1]
%\end{cases} \textnormal{ if } 1 \le j \le q, j \neq r, \textnormal{ } \widetilde{x}_r^r(z) &=
%\begin{cases}
%    0 & \text{if}\ z=[v,\xi,\eta] \\
%    x^r(u) & \text{if}\ z=[u,\upsilon,1]
%\end{cases}, \textnormal{ and }
%\end{align*}
\begin{equation}\label{eq: blup a_r^j}
\widetilde{e}_r^j(z) =
\begin{cases}
    \frac{de^j(v,\xi)\eta}{dx^r(v)\xi} & \text{if}\ z=[v,\xi,\eta] \\
    \frac{e^j(u,\upsilon)}{x^r(u)} & \text{if}\ z=[u,\upsilon,1]
\end{cases} \textnormal{ if } 1 \le j \le k
\end{equation}
(see also \eqref{eq: blup x_r^j}). Notice that if $\eta \in \normal_{(v,\xi)}(E,F)$, then we can write (recall: $d_\normal\psi=(y,dx,f,de)$)
\[\eta=(\eta^1,\eta^2)\] 
relative to the decomposition $\normal_{(v,\xi)}(E,F) = d_\normal\psi(v,\xi)^{-1}(\mathbb{R}^k) \oplus d_\normal\psi(v,\xi)^{-1}(\mathbb{R}^\ell)$.
We need to check that, for all $u \in U \setminus V$ and $(v,\xi) \in \normal(\base,\subbase)$ with $\xi \neq 0$, the maps 
\[(f,\widetilde{e}_r)([u,\cdot,t]): \blup_{\pi}(E,F)_{[u,1]} \ra \mathbb{R}^{k+\ell} \textnormal{ and } (f,\widetilde{e}_r)([v,\xi,\cdot]): \blup_{\pi}(E,F)_{[v,\xi]} \ra \mathbb{R}^{k+\ell}\] 
are linear. This, however, just follows from a simple computation: let $\upsilon,\upsilon' \in E_u$. Then,
\begin{align*}
    (f,\widetilde{e}_r)([u,\upsilon + \upsilon',1]) &= (f(u,\upsilon + \upsilon'),\frac{e^1(u,\upsilon + \upsilon')}{x^r(u)},\dots,\frac{e^k(u,\upsilon + \upsilon')}{x^r(u)}) \\ 
    &= (f(u,\upsilon),\frac{e^1(u,\upsilon)}{x^r(u)},\dots,\frac{e^k(u,\upsilon)}{x^r(u)}) + (f(u,\upsilon'),\frac{e^1(u,\upsilon')}{x^r(u)},\dots,\frac{e^k(u,\upsilon')}{x^r(u)}) \\
    &= (f,\widetilde{e}_r)([u,\upsilon,1]) + (f,\widetilde{e}_r)([u,\upsilon',1]),
\end{align*}
where we used that $(f,e)(u,\cdot): E_u \ra \mathbb{R}^{k+\ell}$ is linear by assumption. Similarly, for $\upsilon \in E_u$ and $\mu \in \mathbb{R}$, we have that $(f,\widetilde{e}_r)([u,\mu \cdot \upsilon,1]) = \mu\cdot(f,\widetilde{e}_r)([u,\upsilon,1])$. Now let $\eta,\zeta \in \normal_{(v,\xi)}(E,F)$. Then,
\begin{align*}
    (f,\widetilde{e}_r)([v,\xi,\eta + \zeta]) & = (f(v,(\eta+\zeta)^1),\frac{de^1(v,\xi)(\eta+\zeta)}{dx^r(v)\xi},\dots,\frac{de^k(v,\xi)(\eta+\zeta)}{dx^r(v)\xi}) \\
    &= (f(v,\eta^1),\frac{de^1(v,\xi)\eta}{dx^r(v)\xi},\dots,\frac{de^k(v,\xi)\eta}{dx^r(v)\xi}) \\
    &+ (f(v,\zeta^1),\frac{de^1(v,\xi)\zeta}{dx^r(v)\xi},\dots,\frac{de^k(v,\xi)\zeta}{dx^r(v)\xi}) \\
    &= (f,\widetilde{e}_r)([v,\xi,\eta]) + (f,\widetilde{e}_r)([v,\xi,\zeta]).
\end{align*}
This shows that our construction does indeed define a vector bundle structure on $\blup_{\pi}(E,F) \ra \blup(\base,\subbase)$. 
\end{proof}
From the proof we get a few immediate consequences.
\begin{coro}\label{coro: blow-up of vector bundles projectivisation}
We have a vector bundle 
\[\mathbb{P}_\pi(E,F) \ra \mathbb{P}(\base,\subbase),\]
for which $\normal_\pi(E,F) \rightarrowdbl \mathbb{P}_\pi(E,F)$ is a morphism of vector bundles. Moreover, we can view it as a vector subbundle of $\blup_\pi(E,F)$.
\end{coro}
\begin{coro}\label{coro: morphism of vector bundles DNC -> Blup}
The quotient map $\DNC(E,F) \setminus F \times \mathbb{R} \rightarrowdbl \blup(E,F)$ restricts to a morphism of vector bundles
\[\DNC_\pi(E,F) \ra \blup_\pi(E,F)\]
over the quotient map $\DNC(\base,\subbase) \setminus \subbase \times \mathbb{R} \rightarrowdbl \blup(\base,\subbase)$.
\end{coro}
\begin{coro}\label{coro: blow-down map}
The blow-down map $\blup(E,F) \ra E$ restricts to a morphism of vector bundles
\[p: \blup_\pi(E,F) \ra E\]
which further restricts to an isomorphism between $E|_{\base \setminus \subbase} \subset \blup_\pi(E,F)$ and $E|_{\base \setminus \subbase}$.
\end{coro}
We now describe the space of sections of $\blup_\pi(E,F) \ra \blup(\base,\subbase)$.
\begin{prop}\label{prop: space of sections blup}
We have a linear map
\[\Gamma(E,F) \ra \Gamma(\blup_\pi(E,F)) \textnormal{ given by } \alpha \mapsto \blup(\alpha)\]
whose image generates $\Gamma(\blup_\pi(E,F))$ as a $C^\infty(\blup(\base,\subbase))$-module, and it is a bijection if $\subbase \subset \base$ is a codimension $1$ submanifold.
\end{prop}
\begin{proof}
First of all, a vector bundle chart $\psi=(y,x,f,e)$ of $E$ (adapted to $F$) gives rise to the local frame 
\[\{\psi^{-1}(\cdot,e^j_{\mathbb{R}}) \mid e^1_{\mathbb{R}},\dots,e^{k+\ell}_{\mathbb{R}} \textnormal{, the standard basis of } \mathbb{R}^{k+\ell}\}.\]
Write $\{\alpha^i\}$ for $\{\psi^{-1}(y,x,0,e^i_{\mathbb{R}}) \mid 1 \le i \le \ell\}$, and write $\{\beta^j\}$ for $\{\psi^{-1}(y,x,e^j_{\mathbb{R}},0) \mid 1 \le j \le k\}$. Then, by the same procedure, but with the vector bundle chart $(y,\widetilde{x}_r,f,\widetilde{e}_r)$ ($1 \le r \le q$; see the proof of Theorem \ref{theo: vector bundle blow-up}), we obtain a local frame 
\[\beta^1,\dots,\beta^k,\widetilde\alpha^1_r,\dots,\widetilde\alpha^\ell_r\]
defined on $V_r$ (see \eqref{eq: blup V_r}), where, for all $1 \le i \le \ell$, $\widetilde\alpha^i_r: U_r \ra V_r$ (see \eqref{eq: blup U_r}) is given by
\begin{align*}
z \mapsto
\begin{cases}
[v,\xi,d_\normal\psi(v,\xi)^{-1}(0,e^i_{\mathbb{R}}dx^r(v)\xi)] & \text{if}\ z=[v,\xi] \\
[\psi^{-1}(y(u),x(u),0,x^r(u) \cdot e^i_{\mathbb{R}}),1] & \text{if}\ z=[u,1],
\end{cases}
\end{align*}
where we used that, for $[v,\xi] \in U_r \cap \normal(\base,\subbase)$, $dx^r(v)\xi \in \mathbb{R}^\times$.
Moreover, upon restriction to $V$, the local section 
\[x^r\alpha^i=\psi^{-1}(y,x,0,x^re^i): U \ra E|_U\] 
becomes the zero section. Therefore, we get an induced local section $\DNC(U,V) \ra \DNC(E|_U,F|_{V})$, which is given by
\begin{align*}
    z \mapsto \DNC(x^r\alpha^i)(z) = 
\begin{cases}
    [v,\xi,d_\normal\psi(v,\xi)^{-1}(0,e^i_{\mathbb{R}}dx^r(v)\xi)] & \text{if}\ z=[v,\xi]  \\
    [\psi^{-1}(y(u),x(u),0,x^r(u)\cdot e^i_{\mathbb{R}}),t] & \text{if}\ z=[u,t].
\end{cases}
\end{align*}
Since this map is $\mathbb{R}^\times$-equivariant, it descend to the section $\blup(x^r\alpha^i)$ of $\blup_{\pi}(E,F) \ra \blup(\base,\subbase)$, but from the above formulas we see that this is just the section $\widetilde\alpha^i_r$. Since $\Gamma(E,F)$ has local frame $\{x^r\alpha^i,\beta^j\}$, we see that the image of
\[\Gamma(E,F) \ra \Gamma(\blup_\pi(E,F)) \textnormal{ given by } \alpha \mapsto \blup(\alpha)\]
generates $\Gamma(\blup_\pi(E,F))$ as a $C^\infty(\blup(\base,\subbase))$-module, and since the sections $d_\normal\psi^{-1}(v,\xi,0,e^i)$ and $\psi^{-1}(y,x,0,e^i)$ are linear (in the last component), the above map is linear. %To see that it is $C^\infty(\base)$-linear, notice that, indeed, we have a $C^\infty(\base)$-structure on $\Gamma(\blup_\pi(\algebr,\subalgebr))$: whenever $f \in C^\infty(\base)$ and $\widetilde\alpha \in \Gamma(\blup_\pi(E,F))$, we set 
%\begin{align*}
%   f \cdot \widetilde\alpha: \base &\ra \blup_\pi(\algebr,\subalgebr) \textnormal{ is the map } \\ z \mapsto
%&\begin{cases}
%    d_\normal f(y)\xi \cdot \widetilde\alpha([y,\xi]) & \text{if}\ z=[y,\xi] \\
%    f(x)\cdot \widetilde\alpha([x,1]) & \text{if}\ z=[x,1]
%\end{cases}
%\end{align*} 
%(note: $d_\normal f$ is a smooth map $\normal(\base,\subbase) \ra \normal(\mathbb{R},\mathbb{R}) \cong \mathbb{R}$). Then,
%\[\blup(f \cdot \widetilde\alpha) = f \cdot \blup(\widetilde\alpha).\]
%Again from the definition, $\Gamma(E,F) \rightarrowdbl \Gamma(\blup_\pi(E,F))$ is now readily verified to be $C^\infty(\base)$-linear.
Lastly, if $\subbase \subset \base$ is of codimension $1$, i.e. $q=1$ in our notation, we obtain a local frame 
\[\beta^1,\dots,\beta^k,\widetilde{\alpha}^1,\dots,\widetilde{\alpha}^\ell\]
for all of $\blup_\pi(E|_U,F|_V)$. On the other hand, since $\blup(\base,\subbase) \cong \base$, the blow-down map
\[\bldown: \blup_\pi(E,F) \ra E\]
induces a map on sections $\bldown: \Gamma(\blup_\pi(E,F)) \ra \Gamma(E)$. If $\widetilde\alpha \in \Gamma(\blup_\pi(E,F))$, then $\widetilde\alpha|_{\subbase} \in \Gamma(\mathbb{P}_\pi(E,F))$ by definition of the vector bundle structure on $\blup_\pi(E,F)$. It follows that
\[p(\widetilde{\alpha})|_{\subbase} \in \Gamma(F),\]
so the blow-down map, on the level of sections, maps into $\Gamma(E,F)$. Moreover, it is readily verified that this map is inverse to the map $\Gamma(E,F) \ra \Gamma(\blup_\pi(E,F))$. This concludes the proof.
\end{proof}
Lastly, we prove also the following two useful statements which are very important for the next section.
\begin{prop}\label{prop: embedding and blup vector bundles}
Let $(\algebrtwoid, \subalgebrtwoid)$ be a pair of vector bundles such that $\algebrtwo \subset E$ and $\subalgebrtwo \subset F$ are vector subbundles, and $\algebrtwo \cap F = \subalgebrtwo$. Then we can view $\blup_{\pi_{\algebrtwo}}(\algebrtwo,\subalgebrtwo)$ as a vector subbundle of $\blup_{\pi_E}(E,F)$.% Moreover, if $\algebrtwo \subset E$ is open (resp. closed), then $\blup_{\pi}(\algebrtwo,\subalgebrtwo) \subset \blup_{\pi}(E,F)$ is open (resp. closed). %Moreover, the map $\Gamma(E,F) \ra \Gamma(\blup_\pi(E,F))$ (from Proposition \ref{prop: space of sections blup}) restricts to a map 
%\begin{equation}\label{eq: restriction of Gamma(E,F) ra Gamma(blup(E,F))}
%    \{\alpha \in \Gamma(E,F) \mid \alpha|_{\basetwo} \in \Gamma(\algebrtwo,\subalgebrtwo)\} \ra \Gamma(\blup_{\pi}(E,F),\blup_{\pi}(\algebrtwo,\subalgebrtwo))
%\end{equation}
%whose image generates $\Gamma(\blup_{\pi}(E,F),\blup_{\pi}(\algebrtwo,\subalgebrtwo))$ as a $\blup(\base,\subbase)$-module if $(\base,\basetwo)$ and $(\subbase,\subbasetwo)$ are pairs of manifolds (that is, $\basetwo \subset \base$ and $\subbasetwo \subset \subbase$ are closed).
\end{prop}
\begin{proof}
Recall from Corollary \ref{coro: blup submanifold} that the embedding $\iota: (\algebrtwo,\subalgebrtwo) \hookrightarrow (E,F)$ induces the embedding $\blup(\iota): \blup(\algebrtwo,\subalgebrtwo) \hookrightarrow \blup(E,F)$,
so in particular $\blup_{\pi_{\algebrtwo}}(\algebrtwo,\subalgebrtwo)$ is a submanifold of $\blup_{\pi_E}(E,F)$. Since it is invariant under scalar multiplication, the %first
result follows from Corollary \ref{coro: vector subbundle is invariant submanifold}. %For the last statement, it is clear that the map $\Gamma(E,F) \ra \Gamma(\blup_\pi(E,F))$ from Proposition \ref{prop: space of sections blup} restricts to \eqref{eq: restriction of Gamma(E,F) ra Gamma(blup(E,F))}. To see that its image generates $\Gamma(\blup_{\pi}(E,F),\blup_{\pi}(\algebrtwo,\subalgebrtwo))$ as a $\blup(\base,\subbase)$-module if $\basetwo \subset \base$ and $\subbasetwo \subset \subbase$ are closed, so we can extend sections of $\subalgebrtwo$ to $F$, then to $\algebrtwo$, and then to $E$. So, if $\widetilde\alpha \in \Gamma(\blup_\pi(E,F),\blup_\pi(\algebrtwo,\subalgebrtwo))$, then we can first pick a lift $\alpha_\subalgebrtwo$ of $\widetilde\alpha|_{\blup(\basetwo,\subbasetwo)}$ (i.e. $\alpha_\subalgebrtwo \in \Gamma(\algebrtwo,\subalgebrtwo)$ with $\blup(\alpha)=\widetilde\alpha|_{\blup(\basetwo,\subbasetwo)}$, and extend it to a section
\end{proof}
\begin{prop}\label{prop: vector bundle morphism to blup}
Let $\varphi: (E,F) \ra (\algebrtwo,\subalgebrtwo)$ be a morphism of vector bundles between pairs of vector bundles $(E \ra \base,F \ra \subbase)$ and $(\algebrtwoid,\subalgebrtwoid)$. If the base map $f: (\base,\subbase) \ra (\basetwo,\subbasetwo)$ satisfies $f^{-1}(\subbasetwo)=\subbase$ and $d_\normal f$ is fiberwise linear, then $\blup(\varphi)$ restricts to a morphism of vector bundles $\blup_{\pi_E}(E,F) \ra \blup_{\pi_{\algebrtwo}}(\algebrtwo,\subalgebrtwo)$.
\end{prop}
\begin{proof}
That $\blup_{\pi_E}(E,F) \subset \blup_\varphi(E,F)$, and that the restriction of $\blup(\varphi): \blup_\varphi(E,F) \ra \blup(\algebrtwo,\subalgebrtwo)$ to $\blup_{\pi_E}(E,F)$ lands in $\blup_{\pi_{\algebrtwo}}(\algebrtwo,\subalgebrtwo)$, follow both from the identity 
\[\DNC(\pi_\algebrtwo) \circ \DNC(\varphi) = \DNC(\pi_\algebrtwo \circ \varphi) = \DNC(f) \circ \DNC(\pi_E)\]
(this works exactly like in the proof of Proposition \ref{prop: blow-up functoriality groupoids}).
The rest of the statement now follows from the fact that $\DNC(\varphi): \DNC(E,F) \ra \DNC(\algebrtwo,\subalgebrtwo)$ is a morphism of vector bundles (see Proposition \ref{prop: DNC morphism of vector bundles}), which, for the same reasons as for $\blup(\varphi)$, restricts to a map $\DNC(\varphi): \DNC_{\pi_E}(E,F) \ra \DNC_{\pi_\algebrtwo}(\algebrtwo,\subalgebrtwo)$ (note: $\DNC_{\pi_E}(E,F)=\DNC(E,F)|_{\DNC(\base,\subbase) \setminus \subbase \times \mathbb{R}}$ and $\DNC_{\pi_\algebrtwo}(\algebrtwo,\subalgebrtwo)=\DNC(\algebrtwo,\subalgebrtwo)|_{\DNC(\base,\subbase) \setminus \subbase \times \mathbb{R}}$, so this map is just given by restricting to $\DNC(\base,\subbase) \setminus \subbase \times \mathbb{R}$). This map is again a morphism of vector bundles, so, by construction of the vector bundle structures on $\blup_{\pi_E}(E,F)$ and $\blup_{\pi_{\algebrtwo}}(\algebrtwo,\subalgebrtwo)$, the induced map $\blup(\varphi)$ is fiberwise linear. This shows that $\blup(\varphi)$ is a morphism of vector bundles, as required. 
\end{proof}

\subsubsection{Blow-ups of Lie algebroids}\label{sec: subsub blow-up lie algebroid}
Before we move on to the blow-up of Lie algebroids, let us consider anchored vector bundles first. To obtain an anchor map for blow-ups of anchored vector bundles, we will start by blowing up tangent bundles.
So, fix a pair of manifolds $(\base,\subbase)$. The blow-up of $T\base \xra{\pi_{T\base}} \base$ along $T\subbase \ra \subbase$, yields, by Theorem \ref{theo: vector bundle blow-up}, a vector bundle
\[\blup(\pi_{T\base}): \blup_{\pi_{T\base}}(T\base,T\subbase) \ra \blup(\base,\subbase)\]
(see also Example \ref{exam: blup pair groupoids}).
\begin{prop}\label{prop: tanent bundle blow-up}
Denote by $q: \normal(\base,\subbase) \setminus 0_\subbase \rightarrowdbl \mathbb{P}(\base,\subbase)$ the quotient projection. We have a morphism of vector bundles
\begin{align*}
    \blup_{\pi_{T\base}}(T\base,T\subbase) &\ra T\blup(\base,\subbase) \textnormal{ given by } \\
    z &\mapsto
\begin{cases}
    ([y,\xi],dq(y,\xi)\eta) & \text{if}\ z=[y,\xi,\eta] \textnormal{ }(\eta \in T_{(y,\xi)}\normal(\base,\subbase)) \\
    a & \text{if}\ z=[a,1] \textnormal{ }(a \in T\base|_{\base \setminus \subbase}),
\end{cases}
\end{align*} 
which is an isomorphism when restricted to $T\base|_{\base \setminus \subbase}$, and restricts to a surjective vector bundle morphism $\mathbb{P}_{\pi_{T\base}}(T\base,T\subbase) \rightarrowdbl T\mathbb{P}(\base,\subbase)$ which is of constant rank with $1$-dimensional kernel.
\end{prop}
\begin{proof}
Recall that we have an embedding $\DNC(T\base,T\subbase) \hookrightarrow T\DNC(\base,\subbase)$ (see Proposition \ref{prop: DNC of tangent spaces}), which is a morphism of vector bundles over $\DNC(\base,\subbase)$ (see the proof of Proposition \ref{prop: DNC of tangent spaces}). Moreover, it restricts to a morphism of vector bundles
\[\DNC_{\pi_{T\base}}(T\base,T\subbase) \ra TP,\]
where $P \coloneqq \DNC(\base,\subbase) \setminus (Y \times \mathbb{R})$. Consider the induced $T\mathbb{R}^\times$-action on $TP$ given by 
\begin{align*}
    (\lambda,v) \cdot (z,\eta) \coloneqq (\lambda \cdot z, d\rho(\lambda,z)(v,\eta)); \quad (\lambda,v) \in T\mathbb{R}^\times, (z,\eta) \in TP
\end{align*}
where we denoted $\rho: \mathbb{R}^\times \times P \ra P$ for the action map. Then the zero section $\mathbb{R}^\times \hookrightarrow T\mathbb{R}^\times$ is a Lie group homomorphism, and, moreover, $T\blup(\base,\subbase) = TP/T\mathbb{R}^\times$. The above morphism of vector bundles is now $\mathbb{R}^\times$-equivariant (with $\mathbb{R}^\times$-action on $TP$ induced by the zero section $\mathbb{R}^\times \hookrightarrow T\mathbb{R}^\times$), so we obtain a morphism of vector bundles
\[\blup_{\pi_{T\base}}(T\base,T\subbase) \ra T\blup(\base,\subbase).\]
More explicitly, this map is an equality on $T\base|_{\base \setminus \subbase} \subset \blup_{\pi_{T\base}}(T\base,T\subbase)$, and for all $[y,\xi] \in \mathbb{P}(\base,\subbase)$, the induced map 
\[\normal_{\pi_{T\base}}(T\base,T\subbase)_{(y,\xi)} \cong  \mathbb{P}_{\pi_{T\base}}(T\base,T\subbase)_{[y,\xi]}  \cong \blup_{\pi_{T\base}}(T\base,T\subbase)_{[y,\xi]} \ra T_{[y,\xi]}\blup(\base,\subbase)\]
maps onto $T_{[y,\xi]}\mathbb{P}(\base,\subbase)$, because $\DNC_{\pi_{T\base}}(T\base,T\subbase) \ra TP$ maps $\normal_{\pi_{T\base}}(T\base,T\subbase)$ onto $T(\normal(\base,\subbase) \setminus 0_\subbase)$, and \[T(\normal(\base,\subbase) \setminus 0_\subbase)/T\mathbb{R}^\times \cong T\mathbb{P}(\base,\subbase)\]
(note: the isomorphism is induced by the differential of the quotient map $q: \normal(\base,\subbase) \setminus 0_\subbase \rightarrowdbl \mathbb{P}(\base,\subbase)$).
From this we see that $\blup_{\pi_{T\base}}(T\base,T\subbase) \ra T\blup(\base,\subbase)$ restricts to a surjective morphism of vector bundles
\[\mathbb{P}_{\pi_{T\base}}(T\base,T\subbase) \rightarrowdbl T\mathbb{P}(\base,\subbase).\]
To see why the last statement holds, we have to show that, fiberwise, the kernel of this map has dimension $1$ (since the base map is the identity map $\id_{\mathbb{P}(\base,\subbase)}$; see Lemma \ref{lemm: vector bundle constant rank maps}). But this follows from the fact that $dq$ has fiberwise a $1$-dimensional kernel given explicitly, for $dq(y,\xi): T_{(y,\xi)}(\normal(\base,\subbase) \setminus 0_\subbase) \ra T_{[y,\xi]}\mathbb{P}(\base,\subbase)$, by $d\rho(1,(y,\xi))T_1\mathbb{R}^\times$. Indeed, $q: \normal(\base,\subbase) \setminus 0_\subbase \ra \mathbb{P}(\base,\subbase)$ satisfies, for all $(y,\xi) \in \normal(\base,\subbase) \setminus 0_\subbase$, the identity $q \circ \rho(\cdot,(y,\xi)) = q(y,\xi)$, so it follows by an application of the chain rule. This proves the statement.
\end{proof}
\begin{rema}\label{rema: anchored vector bundle blow-up}
From this we see that if $(E \xra{\pi} \base,F \ra \subbase)$ is a pair of anchored vector bundles, then $\blup_\pi(E,F) \ra \blup(\base,\subbase)$ naturally inherits an anchor by setting
\[\anchor: \blup_\pi(E,F) \xra{\blup(\anchor_E)} \blup_{\pi_{T\base}}(T\base,T\subbase) \ra T\blup(\base,\subbase),\]
(see also Proposition \ref{prop: vector bundle morphism to blup}).
Moreover, note that the blow-down map $\bldown: \blup_\pi(E,F) \ra E$ is now an anchored vector bundle morphism, i.e. the diagram
\begin{center}
\begin{tikzcd}
    \blup_\pi(E,F) \ar[r,"\bldown_E"] \ar[d,"\anchor"] & E \ar[d,"\anchor_E"] \\
    T\blup(\base,\subbase) \ar[r,"d\bldown_\base"] & T\base
\end{tikzcd}
\end{center}
commutes, where $\bldown_E$ and $\bldown_\base$ are the blow-down maps. Indeed, it suffices to prove it in case $E=T\base$ (because $\blup(\anchor_E)$ commutes with the blow-down maps by Lemma \ref{lemm: induced global smooth map on blow up}), and the commutativity of the diagram is clear when restricting to $E|_{\base \setminus \subbase} \subset \blup_\pi(E,F)$, so the commutativity of the whole diagram follows because $E|_{\base \setminus \subbase} \subset \blup_\pi(E,F)$ is open and dense. 
%\begin{center}
%\begin{tikzcd}
%    \mathbb{P}_{\pi_{T\base}}(T\base,T\subbase) \ar[r,"\bldown_E|_{\mathbb{P}_\pi(T\base,T\subbase)}"] \ar[d,"\anchor"] & T\subbase \ar[d,"="] \\
%    T\mathbb{P}(\base,\subbase) \ar[r,"d\bldown_\base|_{T\mathbb{P}(\base,\subbase)}"] & T\subbase,
%\end{tikzcd}
%end{center}
%which follows from the commutativity of
%\begin{center}
%\begin{tikzcd}
%    \normal(T\base,T\subbase) \ar[r] \ar[d,"\sim"] & T\subbase \ar[d,"="] \\
%    T\normal(\base,\subbase) \ar[r] & T\subbase.
%\end{tikzcd}
%\end{center}

This shows that the blow-up of a pair of anchored vector bundles is again an anchored vector bundle, and that the blow-down map is then a morphism of anchored vector bundles.
\end{rema}
Now, fix a pair of Lie algebroids $(\algebr \xra{\pi} \base,\subalgebroid)$. Then, to show that the anchored vector bundle $\blup_\pi(\algebr,\subalgebr) \ra \blup(\base,\subbase)$ is a Lie algebroid, it remains to check that it inherits a Lie bracket that is compatible with the anchor map via the Leibniz rule. \begin{theo}\label{theo: Lie algebroid blow-up}
There is a unique Lie algebroid structure on 
\[\blup_\pi(\algebr,\subalgebr) \ra \blup(\base,\subbase)\]
such that $\algebr|_{\base \setminus \subbase} \subset \blup_\pi(\algebr,\subalgebr)$ is a Lie subalgebroid. Moreover, %we have 
%\[[\blup(\alpha),\blup(\beta)]_{\blup_\pi(\algebr,\subalgebr)}=\blup(\bracket{\alpha,\beta}) \textnormal{ for all } \alpha,\beta \in \Gamma(\algebr),\]
%and
we have a Lie algebra morphism
\[\Gamma(\algebr,\subalgebr) \ra \Gamma(\blup_\pi(\algebr,\subalgebr)) \textnormal{ given by } \alpha \mapsto \blup(\alpha)\]
whose image generates $\Gamma(\blup_\pi(E,F))$ as a $C^\infty(\blup(\base,\subbase))$-module, and it is an isomorphism if $\subbase \subset \base$ is a codimension $1$ submanifold.
\end{theo}
\begin{proof}
Let $(U,y,x)$ be a chart adapted to $\subbase$ such that $\{\beta^1,\dots,\beta^k,\alpha^1,\dots,\alpha^\ell\}$ is a local frame of $\algebr$ defined on $U$, and $\beta^1,\dots,\beta^k$ is a local frame of $\subalgebr$ defined on $V \coloneqq U \cap \subbase$.
Recall from the proof of Proposition \ref{prop: space of sections blup} that, for all $1 \le r \le q$, a local frame on $V_r \subset \blup_\pi(\algebr,\subalgebr)$ (see \eqref{eq: blup V_r}) is given by 
\[\beta^1,\dots,\beta^k,\widetilde{\alpha}^1_r,\dots,\widetilde{\alpha}^\ell_r,\]
where $\widetilde{\alpha}^i_r=\blup(x^r\alpha^i)$. Since
\begin{align*}
    \bracket{x^r\alpha^i,x^r\alpha^j} &= x^r\bracket{x^r\alpha^j,\alpha^j} + x^r\anchor(\alpha^i)(x^r)\alpha^j \\
    &= (x^r)^2\bracket{\alpha^i,\alpha^j} + x^r\anchor(\alpha^i)(x^r)\alpha^j - x^r\anchor(\alpha^j)(x^r)\alpha^i,
\end{align*}
which vanishes along $V$, we get an induced map which we take to be the Lie bracket:
\[[\widetilde{\alpha}^i_r,\widetilde{\alpha}^j_r]_{\blup_\pi(\algebr,\subalgebr)} \coloneqq \blup(\bracket{x^r\alpha^i,x^r\alpha^j}).\]
Similarly, we define
\[[\widetilde{\alpha}^i_r,\beta^j]_{\blup_\pi(\algebr,\subalgebr)} \coloneqq \blup(\bracket{x^r\alpha^i,\beta^j}), \textnormal{ and } [\beta^i,\beta^j]_{\blup_\pi(\algebr,\subalgebr)} \coloneqq \blup(\bracket{\beta^i,\beta^j}).\]
To see that these expressions are independent of chosen local frame follows from the fact that %the map $\Gamma(\algebr,\subalgebr) \ra \Gamma(\blup_\pi(\algebr,\subalgebr))$ is $C^\infty(\base)$-linear (see Theorem \ref{theo: vector bundle blow-up}), so
\begin{equation}\label{eq: blup is Lie algebra morphism}
    \blup(\bracket{\alpha,\beta}) = [\blup(\alpha),\blup(\beta)]_{\blup_\pi(\algebr,\subalgebr)}
\end{equation}
(since the expression holds locally). We now get a unique bracket on $\blup_\pi(\algebr,\subalgebr)$ satisfying the above equations by requiring the Leibniz identity to hold with respect to the anchor map
\[\blup_\pi(\algebr,\subalgebr) \xra{\blup(\anchor)} \blup_{\pi_{T\base}}(T\base,T\subbase) \ra T\blup(\base,\subbase)\] 
(see Remark \ref{rema: anchored vector bundle blow-up}). To see that this indeed defines a Lie bracket, %the map $\Gamma(\algebr,\subalgebr) \ra \Gamma(\blup_\pi(\algebr,\subalgebr))$ is $C^\infty(\base)$-linear (see Theorem \ref{theo: vector bundle blow-up}), so
%\begin{equation}\label{eq: blup is Lie algebra morphism}
%    \blup(\bracket{\alpha,\beta}) = [\blup(\alpha),\blup(\beta)]_{\blup_\pi(\algebr,\subalgebr)}
%\end{equation}
%holds by definition (since the expression holds locally). By 
recall Remark \ref{rema: defining a Lie algebroid locally out of an anchored vector bundle}. The skew-symmetry of $[\cdot,\cdot]_{\blup_\pi(\algebr,\subalgebr)}$ follows from the skew-symmetry of $\bracket{\cdot,\cdot}$, and the Jacobi identity also follows: let $\widetilde\alpha,\widetilde\beta,\widetilde\gamma \in \Gamma(\blup_\pi(\algebr,\subalgebr))$. By arguing locally (see Remark \ref{rema: defining a Lie algebroid locally out of an anchored vector bundle}), we can take $\widetilde\alpha$,$\widetilde\beta$, and $\widetilde\gamma$ such that there are lifts $\alpha,\beta,\gamma \in \Gamma(\algebr,\subalgebr)$, i.e.
\[\blup(\alpha)=\widetilde\alpha, \blup(\beta)=\widetilde\beta, \textnormal{ and } \blup(\gamma)=\widetilde\gamma.\]
By linearity of the map $\Gamma(\algebr,\subalgebr) \ra \Gamma(\blup_\pi(\algebr,\subalgebr))$, together with \eqref{eq: blup is Lie algebra morphism}, we have 
\begin{align*}
    0 = \blup(\bigodot_{\alpha,\beta,\gamma} \bracket{\bracket{\alpha,\beta},\gamma}) &= \bigodot_{\alpha,\beta,\gamma} \blup(\bracket{\bracket{\alpha,\beta},\gamma}) \\
    &= \bigodot_{\alpha,\beta,\gamma} [[\blup(\alpha),\blup(\beta)]_{\blup_\pi(\algebr,\subalgebr)},\blup(\gamma)]_{\blup_\pi(\algebr,\subalgebr)} \\
    &= \bigodot_{\alpha,\beta,\gamma} [[\widetilde\alpha,\widetilde\beta]_{\blup_\pi(\algebr,\subalgebr)},\widetilde\gamma]_{\blup_\pi(\algebr,\subalgebr)}.
\end{align*}
This shows that we obtain a Lie algebroid this way. The last statements follow immediately from Proposition \ref{prop: space of sections blup} and \eqref{eq: blup is Lie algebra morphism}. This concludes the proof.
%whenever $\widetilde\alpha,\widetilde\beta \in \Gamma(\blup_\pi(\algebr,\subalgebr))$, then for all expressions
%\[\bldown(\widetilde\alpha) = \sum_i f_i p^*\alpha^i\]
%To finish the proof, it remains to check the Jacobi-identity of $[\cdot,\cdot]_{\blup_\pi(\algebr,\subalgebr)}$.For all $1 \le j \le m$, we can write
%\[\anchor \gamma^j = \sum_i b^{ji} \frac{\partial}{\partial x^i}, \textnormal{ where } b^{ji} \in C^\infty(U),\]
%and, for all $1 \le j,j' \le m$, we can write
%\[\bracket{\gamma^j,\gamma^{j'}} = \sum_d c_d^{jj'}\gamma^d, \textnormal{ where } c_d^{jj'} \in C^\infty(U).\]
\end{proof}
Let us prove right away that this blow-up construction for Lie algebroids is the right analogue of the Lie groupoid blow-up:
\begin{prop}\label{prop: Lie groupoid blup Lie algebroid}
Let $(\groupoid,\subgroupoid)$ be a pair of Lie groupoids, and let $(\algebroid,\subalgebroid)$ be the corresponding pair of their Lie algebroids. Then 
\[\blup_{\source,\target}(\group,\subgroup) \rra \blup(\base,\subbase)\]
has Lie algebroid (canonically isomorphic to)
\[\blup_\pi(\algebr,\subalgebr) \ra \blup(\base,\subbase).\]\vspace{-\baselineskip}
\end{prop}
\begin{proof}
We only have to show that the total space of the Lie algebroid $C$ of $\blup_{\source,\target}(\group,\subgroup)$ is (canonically) diffeomorphic to $\blup_\pi(\algebr,\subalgebr)$, that $C$ contains $\algebr|_{\base \setminus \subbase}$ as a Lie subalgebroid, and that the diffeomorphism maps $\algebr|_{\base \setminus \subbase} \subset C$ identically to $\algebr|_{\base \setminus \subbase} \subset \blup_{\source,\target}(\group,\subgroup)$ (since the Lie algebroid structure on $\blup_\pi(\algebr,\subalgebr)$ is uniquely determined by the property that $\algebr|_{\base \setminus \subbase} \subset \blup_\pi(\algebr,\subalgebr)$ is a Lie subalgebroid). The diffeomorphism can be realised by the sequence of diffeomorphisms
\begin{align*}
    C &\cong \normal(\DNC_{\source,\target}(\group,\subgroup)/\mathbb{R}^\times,(\DNC(\base,\subbase) \setminus \subbase \times \mathbb{R})/\mathbb{R}^\times) \\
    &\cong \normal(\DNC_{\source,\target}(\group,\subgroup),(\DNC(\base,\subbase)\setminus \subbase \times \mathbb{R}))/\mathbb{R}^\times \\
    &\cong \DNC_\pi(\algebr,\subalgebr)/\mathbb{R}^\times \cong \blup_\pi(\algebr,\subalgebr),
\end{align*}
where we used in the first isomorphism that the Lie algebroid of a Lie groupoid $\grouptwoid$ is (canonically isomorphic to) the normal bundle $\normal(\grouptwo,\basetwo)$ (as Lie algebroids), in the second isomorphism we used Proposition \ref{prop: induced G action on normal bundle} (note: as an alternative proof, one can show that this is an isomorphism of Lie algebroids), and in the third isomorphism we used that the open subgroupoid $\DNC_{\source,\target}(\group,\subgroup) \rra \DNC(\base,\subbase) \setminus \subbase \times \mathbb{R}$ of $\DNC(\group,\subgroup)$ has Lie algebroid $\DNC_\pi(\algebr,\subalgebr) \ra \DNC(\base,\subbase)\setminus \subbase \times \mathbb{R}$. To see why this last assertion holds, observe that $\algebr|_{\base \setminus \subbase} \times \mathbb{R}^\times \subset \DNC_\pi(\algebr,\subalgebr)$ is the Lie algebroid of $\group|_{\base \setminus \subbase} \times \mathbb{R}^\times \subset \DNC_{\source,\target}(\group,\subgroup)$, i.e. that $\algebr|_{\base \setminus \subbase} \ra \base \setminus \subbase$ is the Lie algebroid of $\group|_{\base \setminus \subbase} \ra \base \setminus \subbase$; this follows simply from the fact that
\[\identity(\base) \cap (\source^{-1}(\subbase) \cup \target^{-1}(\subbase)) = \identity(\base \setminus \subbase).\]

That $C$ contains $\algebr|_{\base \setminus \subbase}$ as a Lie subalgebroid, follows from the fact that $\blup_{\source,\target}(\group,\subgroup)$ contains $\group|_{\base \setminus \subbase}$ as a subgroupoid, whose Lie algebroid is $\algebr|_{\base \setminus \subbase}$. By tracing the maps, it is clear that the diffeomorphism maps $\algebr|_{\base \setminus \subbase} \subset C$ identically to $\algebr|_{\base \setminus \subbase} \subset \blup_{\source,\target}(\group,\subgroup)$. This shows that, indeed, $\blup_\pi(\algebr,\subalgebr)$ is the Lie algebroid of $\blup_{\source,\target}(\group,\subgroup)$.
\end{proof}
For the rest of this section, we will mainly prove that the results of Section \ref{sec: subsub blow-up lie algebroid} also hold in the Lie algebroid setting. In the end we put a simple example (namely of the case $(\algebr,\subalgebr) = (T\base,T\subbase)$). The following statement follows by construction of $\blup_\pi(\algebr,\subalgebr)$.
\begin{coro}\label{coro: projective lie algebroid}
We have a Lie algebroid
\[\mathbb{P}_\pi(\algebr,\subalgebr) \ra \mathbb{P}(\base,\subbase),\]
for which $\normal_\pi(\algebr,\subalgebr) \rightarrowdbl \mathbb{P}_\pi(\algebr,\subalgebr)$ is a morphism of Lie algebroids. Moreover, we can view it as a subalgebroid of $\blup_\pi(\algebr,\subalgebr)$.
\end{coro}
For integrable Lie algebroids, we now obtain the analogue of Proposition \ref{prop: isotropy groups and orbits of blow-up} (using Proposition \ref{prop: orbits of Lie groupoid are orbits of Lie algebroid}), but this statement also holds for non-integrable Lie algebroids:
\begin{prop}\label{prop: isotropy groups and orbits of blow-up algebroid}
Let $(\algebroid,\subalgebroid)$ be a pair of algebroids. Then $(\blup_{\pi}(\algebr,\subalgebr),\mathbb{P}_{\pi}(\algebr,\subalgebr))$ is a pair of Lie algebroids. Moreover, for all $z \in \blup(\base,\subbase)$, 
\begin{enumerate}
    \item the isotropy algebra $\blup_{\pi}(\algebr,\subalgebr)_z$ equals the isotropy algebra $(\algebr|_{\base \setminus \subbase})_z$ if $z \in \base \setminus \subbase$, and it equals the isotropy algebra $\mathbb{P}_{\pi}(\algebr,\subalgebr)_z$ if $z \in \mathbb{P}(\base,\subbase)$;
    \item the orbit $\orbit_z$ equals the orbit $\orbit_z \subset \base \setminus \subbase$ of $\algebr|_{\base \setminus \subbase}$ if $z \in \base \setminus \subbase$, and it equals the orbit $\orbit_z \subset \mathbb{P}(\base,\subbase)$ of $\mathbb{P}_{\pi}(\algebr,\subalgebr)$ if $z \in \mathbb{P}(\base,\subbase)$.
\end{enumerate}
\end{prop}
\begin{proof}
As in Proposition \ref{prop: isotropy groups and orbits of blow-up}, this is just a consequence of the fact that $\mathbb{P}_{\pi}(\algebr,\subalgebr)$ is the restriction of $\blup_{\pi}(\algebr,\subalgebr)$ to $\mathbb{P}(\base,\subbase)$, and that $\mathbb{P}(\base,\subbase) \subset \blup(\base,\subbase)$ is saturated (i.e. a union of orbits). From Corollary \ref{coro: projective lie algebroid}, we see that $(\blup_{\pi}(\algebr,\subalgebr),\mathbb{P}_{\pi}(\algebr,\subalgebr))$ is a pair of Lie algebroids. Recall that the anchor map is a map of pairs
\[\anchor: (\blup_\pi(\algebr,\subalgebr),\mathbb{P}_\pi(\algebr,\subalgebr)) \ra (\blup_{\pi_{T\base}}(T\base,T\subbase),\mathbb{P}_{\pi_{T\base}}(T\base,T\subbase) \ra (T\blup(\base,\subbase),T\mathbb{P}(\base,\subbase))\]
(see Example \ref{exam: blup tangent bundle} and Remark \ref{rema: anchored vector bundle blow-up}),
so the anchor map restricts to the anchor maps of $\algebr|_{\base \setminus \subbase}$ and $\mathbb{P}_{\pi}(\algebr,\subalgebr)$. Now, $1.$ and $2.$ follow both from this fact. To see why $2.$ holds because of this, notice that if $z \in \base \setminus \subbase$, then $\anchor(z)$ is a map $\blup_\pi(\algebr,\subalgebr)_z \ra T_z\base$, and if $z \in \mathbb{P}(\base,\subbase)$, it is a map $\blup_\pi(\algebr,\subalgebr)_z \ra T_z\mathbb{P}(\base,\subbase)$; then $T_z\orbit_z = \im \anchor(z) \subset T_z\base$ (see Corollary \ref{coro: local splitting theorem LA}) if $z \in \base \setminus \subbase$ and $T_z\orbit_z \subset T_z\mathbb{P}(\base,\subbase)$ if $z \in \mathbb{P}(\base,\subbase)$, i.e. $\orbit_z$ is either already an orbit of $\algebr|_{\base \setminus \subbase}$ or already an orbit of $\mathbb{P}_{\pi}(\algebr,\subalgebr)$. This proves the statement.
\end{proof}
\begin{rema}\label{rema: orbits A|X setminus Y}
Let $z \in \base \setminus \subbase$. Notice that the isotropy Lie algebra $\mathfrak{g}_z$ of the Lie algebroid $\algebr|_{\base \setminus \subbase}$ is in fact equal (or canonically isomorphic) to the isotropy Lie algebra $\mathfrak{g}_z$ of $\algebr$. This holds simply because $\ker \anchor(z) \subset \algebr_z \subset \algebr|_{\base \setminus \subbase}$. The orbit of $z$ in $\algebr|_{\base \setminus \subbase}$ is given by the connected component of $\orbit_z \setminus \subbase$ containing $z$, because $\orbit_z \setminus \subbase \subset \orbit_z$ is open, so that $T_z(\orbit_z \setminus \subbase) = T_z\orbit_z$; see Corollary \ref{coro: local splitting theorem LA}.
\end{rema}
Before we move on, it is useful to first prove the following three results:
\begin{prop}\label{prop: blup of subalgebroid is subalgebroid}
Let $(\algebrtwo \ra \basetwo,\subalgebrtwo \ra \subbasetwo)$ be a pair of Lie algebroids, such that of $\algebrtwo \subset \algebr$ and $\subalgebrtwo \subset \subalgebr$ are Lie subalgebroids, and $\algebrtwo \cap \subalgebr = \subalgebrtwo$. Then we can view $\blup_{\pi}(\algebrtwo,\subalgebrtwo)$ as a Lie subalgebroid of $\blup_{\pi}(\algebr,\subalgebr)$. %If $\algebrtwo \subset \algebr$ is open (resp. closed), then $\blup_{\pi}(\algebrtwo,\subalgebrtwo) \subset \blup_\pi(\algebr,\subalgebr)$ is open (resp. closed).
\end{prop}
\begin{proof}
Recall from Proposition \ref{prop: embedding and blup vector bundles} that we can view $\blup_{\pi}(\algebrtwo,\subalgebrtwo)$ as a vector subbundle of $\blup_{\pi}(\algebr,\subalgebr)$. By functoriality (of the blow-up construction), it is immediate that it is even an anchored subbundle, so we only have to show that
\[\Gamma(\blup_{\pi}(\algebr,\subalgebr),\blup_{\pi}(\algebrtwo,\subalgebrtwo)) \subset \Gamma(\blup_{\pi}(\algebr,\subalgebr))\]
is a Lie subalgebra. First of all, observe that $\mathbb{P}_\pi(\algebrtwo,\subalgebrtwo) \subset \mathbb{P}_\pi(\algebr,\subalgebr)$ is a Lie subalgebroid. This follows, since $\normal(\algebrtwo,\subalgebrtwo) \subset \normal(\algebr,\subalgebr)$ is a subalgebroid, and therefore $\normal_\pi(\algebrtwo,\subalgebrtwo) \subset \normal_\pi(\algebr,\subalgebr)$ is a subalgebroid. Now, if $\widetilde\alpha,\widetilde\beta \in \Gamma(\blup_{\pi}(\algebr,\subalgebr),\blup_{\pi}(\algebrtwo,\subalgebrtwo))$, %Moreover, we can make sure that $\alpha|_\basetwo,\beta|_\basetwo \in \Gamma(\algebrtwo)$. Indeed,
then $\widetilde\alpha|_{\basetwo \setminus \subbasetwo}, \widetilde\beta|_{\basetwo \setminus \subbasetwo} \in \Gamma(\algebrtwo|_{\basetwo \setminus \subbasetwo})$ and $\widetilde\alpha|_{\mathbb{P}(\basetwo,\subbasetwo)}, \widetilde\beta|_{\mathbb{P}(\basetwo,\subbasetwo)} \in \Gamma(\mathbb{P}_\pi(\algebrtwo,\subalgebrtwo))$, simply because $\widetilde\alpha|_{\blup(\basetwo,\subbasetwo)},\widetilde\beta|_{\blup(\basetwo,\subbasetwo)} \in \Gamma(\blup_\pi(\algebrtwo,\subalgebrtwo))$ by assumption. Therefore,
\[[\widetilde\alpha,\widetilde\beta]_{\blup_\pi(\algebr,\subalgebr)} \in \Gamma(\blup_\pi(\algebr,\subalgebr),\algebrtwo|_{\basetwo \setminus \subbasetwo}) \cap \Gamma(\blup_\pi(\algebr,\subalgebr),\mathbb{P}_\pi(\algebrtwo,\subalgebrtwo)),\]
since $\algebrtwo|_{\basetwo \setminus \subbasetwo} \subset \algebr|_{\base \setminus \subbase} \subset \blup_\pi(\algebr,\subalgebr)$ and $\mathbb{P}_\pi(\algebrtwo,\subalgebrtwo) \subset \mathbb{P}_\pi(\algebr,\subalgebr) \subset \blup_\pi(\algebr,\subalgebr)$ are subalgebroids. By definition of the Lie algebroid structure on $\blup_\pi(\algebrtwo,\subalgebrtwo)$, we see that 
\[[\widetilde\alpha,\widetilde\beta]_{\blup_\pi(\algebr,\subalgebr)} \in \Gamma(\blup_\pi(\algebr,\subalgebr),\blup_\pi(\algebrtwo,\subalgebrtwo)).\] This proves the statement.
\end{proof}
A similar result holds if we replace $\blup_\pi$ everywhere with $\DNC$:
\begin{prop}\label{prop: DNC of subalgebroid is subalgebroid}
Let $(\algebrtwo \ra \basetwo,\subalgebrtwo \ra \subbasetwo)$ be a pair of Lie algebroids, such that of $\algebrtwo \subset \algebr$ and $\subalgebrtwo \subset \subalgebr$ are Lie subalgebroids, and $\algebrtwo \cap \subalgebr = \subalgebrtwo$. Then we can view $\DNC(\algebrtwo,\subalgebrtwo)$ as a Lie subalgebroid of $\DNC(\algebr,\subalgebr)$. If $\algebrtwo \subset \algebr$ is open (resp. closed), then $\DNC(\algebrtwo,\subalgebrtwo) \subset \DNC(\algebr,\subalgebr)$ is open (resp. closed).
\end{prop}
The proof is analogous to the proof of Proposition \ref{prop: blup of subalgebroid is subalgebroid}. If we replace $\DNC$ everywhere with $\DNC_\pi$ (which is the statement we will need), an alternative proof can be given by showing that the embedding 
\[\DNC_\pi(\algebr,\subalgebr) \ra \blup_{\pi \times \id_{\mathbb{R}}}(\algebr \times \mathbb{R}, \subalgebr \times \{0\}),\]
from Proposition \ref{prop: DNC out of blup}, is a Lie subalgebroid. The result then follows from Proposition \ref{prop: blup of subalgebroid is subalgebroid}.

For the rest of the statements in this section, the key is to use the following lemma. 
\begin{lemm}\label{lemm: blup(X,Y) x M = blup(X x M,Y x M) LA}
Let $\algebrtwoid$ be a Lie algebroid. The canonical diffeomorphism
\[\blup(\algebr \times \algebrtwo,\subalgebr \times \algebrtwo) \xra{\sim} \blup(\algebr,\subalgebr) \times \algebrtwo\]
(from Proposition \ref{prop: blup(X x M, Y x M) = blup(X,Y) x M})
restricts to an isomorphism of Lie algebroids
\[\varphi: \blup_{\pi_\algebr \times \pi_\algebrtwo}(\algebr \times \algebrtwo,\subalgebr \times \algebrtwo) \xra{\sim} \blup_{\pi_{\algebr}}(\algebr,\subalgebr) \times \algebrtwo.\]
\end{lemm}
\begin{proof}
That the restriction to $\blup_{\pi_{\algebr}}(\algebr,\subalgebr) \times \algebrtwo$ maps bijectively onto $\blup_{\pi_\algebr \times \pi_\algebrtwo}(\algebr \times \algebrtwo,\subalgebr \times \algebrtwo)$ follows from the fact that the isomorphism from Proposition \ref{prop: blup(X x M, Y x M) = blup(X,Y) x M} is induced by the map
\[g: \DNC(\algebr \times \algebrtwo, \subalgebr \times \algebrtwo) \xra{\sim} \DNC(\algebr, \subalgebr) \times \algebrtwo \textnormal{ given by } z \mapsto 
\begin{cases}
    (b,\xi,c) & \text{if}\ z=(b,c,\xi,0) \\
    (a,1,c) & \text{if}\ z=(a,c,t),
\end{cases}\]
and one can readily verify that $\DNC(\pi_\algebr) \circ \pr_1(g(z)) \in \subbase \times \mathbb{R}$ if and only if $\DNC(\pi_\algebr \times \pi_\algebrtwo)(z) \in \subbase \times \basetwo \times \mathbb{R}$% (it just follows from the fact that $\pi_\algebrtwo(\algebrtwo) = \basetwo$)
. Observe that, under the diffeomorphism $\varphi$, 
\[(\algebr \times \algebrtwo)|_{(\base \times \basetwo) \setminus (\subbase \times \basetwo)} = (\algebr \times \algebrtwo)|_{(\base \setminus \subbase) \times \basetwo}\]
maps bijectively onto $\algebr|_{\base \setminus \subbase} \times \algebrtwo$. Now the result follows: since $\blup_{\pi_\algebr \times \pi_\algebrtwo}(\algebr \times \algebrtwo, \subalgebr \times 
\algebrtwo)$ has a unique Lie algebroid structure for which $(\algebr \times \algebrtwo)|_{(\base \setminus \subbase) \times \basetwo}$ is a Lie subalgebroid, we see that the Lie algebroid structures of $\blup_{\pi_\algebr \times \pi_\algebrtwo}(\algebr \times \algebrtwo, \subalgebr \times 
\algebrtwo)$ and $\blup_{\pi_\algebr}(\algebr,\subalgebr) \times \algebrtwo$ have to agree (under the map $\varphi$). This proves the statement.
\end{proof}
\begin{rema}\label{rema: DNC(X,Y) x M = DNC(X x M,Y x M) LA}
Notice that a similar argument shows that we can replace ``$\blup$'' everywhere with ``$\DNC$'', and we even have an isomorphism of Lie algebroids $\DNC(\algebr \times \algebrtwo, \subalgebr \times \algebrtwo) \xra{\sim} \DNC(\algebr,\subalgebr) \times \algebrtwo$.
\end{rema}
We start with showing that the blow-down map for blow-ups of Lie algebroids becomes a morphism of Lie algebroids.
\begin{coro}\label{coro: blow-down map is a morphism of LA}
The blow-down map $\bldown: \blup_\pi(\algebr,\subalgebr) \ra \algebr$ is a morphism of Lie algebroids.
\end{coro}
\begin{proof}
We denote the blow-down map on the base manifold by $\bldown_\base$. 
%We will show that
%\[\gr p \xra{\blup(\pi) \times \pi} \gr p_\base\]
%is a Lie subalgebroid of $\blup_\pi(\algebr,\subalgebr) \times \algebr$. 
By Corollary \ref{coro: blow-down map}, we only have to check that
$\gr \bldown \subset \blup_\pi(\algebr,\subalgebr) \times \algebr$ is a Lie subalgebroid. By Lemma \ref{lemm: blup(X,Y) x M = blup(X x M,Y x M) LA}, there is a canonical isomorphism of Lie algebroids
\begin{equation}\label{eq: iso blups in blow-down map LA}
    \blup_\pi(\algebr,\subalgebr) \times \algebr \xra{\sim} \blup_{\pi \times \pi}(\algebr \times \algebr, \subalgebr \times \algebr).
\end{equation}
Notice that it restricts to an isomorphism of Lie algebroids
\[\mathbb{P}_\pi(\algebr,\subalgebr) \times \algebr \xra{\sim} \mathbb{P}_{\pi \times \pi}(\algebr \times \algebr, \subalgebr \times \algebr),\]
so, by explicitly writing down what elements of $\gr \bldown$ are mapped to under \eqref{eq: iso blups in blow-down map LA}, we see that, under \eqref{eq: iso blups in blow-down map LA}, $\gr \bldown \subset \blup_\pi(\algebr,\subalgebr) \times \algebr$ maps bijectively onto $\blup_{(\pi \times \pi)|_{\Delta_\algebr}}(\Delta_\algebr,\Delta_\subalgebr) \subset \blup_{\pi \times \pi}(\algebr \times \algebr,\subalgebr \times \algebr)$: 
\begin{align*}
(z,p(z)) = 
\begin{cases}
    ([b,\xi],b) \\
    ([a,1],a)
\end{cases} \mapsto 
\begin{cases}
    [b,b,\xi,0] & \text{if}\ z=[b,\xi] \\
    [a,a,1] & \text{if}\ z=[a,1].
\end{cases}
\end{align*}
We can now equivalently prove that
\[\blup_{(\pi \times \pi)|_{\Delta_\algebr}}(\Delta_\algebr,\Delta_\subalgebr) \subset \blup_{\pi \times \pi}(\algebr \times \algebr, \subalgebr \times \algebr)\] 
is a Lie subalgebroid. Since $\Delta_\algebr \cap (\subalgebr \times \algebr) = \Delta_\subalgebr$, the result follows from Proposition \ref{prop: blup of subalgebroid is subalgebroid}. 
\end{proof}
Alternatively, the above statement is a consequence of the fact that $\bldown: \blup_{\pi}(\algebr,\subalgebr) \ra \algebr$ is an isomorphism on a dense open subset of $\blup_{\pi}(\algebr,\subalgebr)$.

We now prove the following analogue of Proposition \ref{prop: vector bundle morphism to blup}.
\begin{coro}\label{coro: LA morphism to blup}
Let $\varphi: (\algebr,\subalgebr) \ra (\algebrtwo,\subalgebrtwo)$ be a morphism of Lie algebroids between pairs of Lie algebroids $(\algebroid,\subalgebroid)$ and $(\algebrtwoid,\subalgebrtwoid)$. If the base map $f: (\base,\subbase) \ra (\basetwo,\subbasetwo)$ of $\varphi$ satisfies $f^{-1}(\subbasetwo)=\subbase$ and $d_\normal f$ is fiberwise injective, then $\blup(\varphi)$ restricts to a morphism of Lie algebroids $\blup_{\pi_\algebr}(\algebr,\subalgebr) \ra \blup_{\pi_{\algebrtwo}}(\algebrtwo,\subalgebrtwo)$.
\end{coro}
\begin{proof}
By Proposition \ref{prop: vector bundle morphism to blup}, we only have to check that $\gr \blup(\varphi)$ is a Lie subalgebroid of $\blup_{\pi_\algebr}(\algebr, \subalgebr) \times \blup_{\pi_\algebrtwo}(\algebrtwo,\subalgebrtwo)$. First of all, by Lemma \ref{lemm: blup(X,Y) x M = blup(X x M,Y x M) LA}, there is an isomorphism of Lie algebroids
\begin{equation}\label{eq: iso blups in LA morphism to blup}
    \blup_{\pi_\algebr}(\algebr,\subalgebr) \times \blup_{\pi_{\algebrtwo}}(\algebrtwo,\subalgebrtwo) \xra{\sim} \blup_{\blup(\pi_\algebr) \times \pi_\algebrtwo}(\blup_{\pi_\algebr}(\algebr,\subalgebr) \times \algebrtwo, \blup_{\pi_\algebr}(\algebr,\subalgebr) \times \subalgebrtwo).
\end{equation}
Notice that it restricts to an isomorphism of Lie algebroids
\[\blup_{\pi_\algebr}(\algebr,\subalgebr) \times \mathbb{P}_{\pi_\algebrtwo}(\algebrtwo,\subalgebrtwo) \xra{\sim} \mathbb{P}_{\blup(\pi_\algebr) \times \pi_\algebrtwo}(\blup_{\pi_\algebr}(\algebr,\subalgebr) \times \algebrtwo,\blup_{\pi_\algebr}(\algebr,\subalgebr) \times \subalgebrtwo).\]
%Since $T\gr f = \gr df$, for smooth maps $f: P \ra Q$, it is readily verified that, under this isomorphism, $\gr d_\normal\varphi \subset \normal(\algebr,\subalgebr) \times \normal(\algebrtwo,\subalgebrtwo)$ maps onto $\normal(\gr \varphi,\gr \varphi|_{\subalgebr})$. 
Now, under \eqref{eq: iso blups in LA morphism to blup}, $\gr \varphi$ maps bijectively onto $\blup_{(\blup(\pi_\algebr) \times \pi_\algebrtwo)|_{\gr \varphi \circ \bldown}}(\gr \varphi \circ \bldown, \gr \varphi \circ \bldown|_{\mathbb{P}_{\pi_\algebr}(\algebr,\subalgebr)})$:
\[(z,\blup(\varphi)(z)) = 
\begin{cases}
([b,\xi],[\varphi(b),d_\normal\varphi(b)\xi]) \\
([a,1],[\varphi(a),1])
\end{cases} \mapsto 
\begin{cases}
[[b,\xi],\varphi(b),d_\normal\varphi(b)\xi] & \text{if}\ z=[b,\xi] \\
[[a,1],\varphi(a),1] & \text{if}\ z=[a,1],
\end{cases}\]
where we used that there is a canonical diffeomorphism
\[\mathbb{P}_{\pi_\algebr}(\algebr,\subalgebr) \cong \mathbb{P}_{(\blup(\pi_\algebr) \times \pi_\algebrtwo)|_{\gr \varphi \circ \bldown}}(\gr \varphi \circ \bldown, \gr \varphi \circ \bldown|_{\mathbb{P}_{\pi_\algebr}(\algebr,\subalgebr)})\]
by the sequence of canonical diffeomorphisms (note: the result is not a morphism of vector bundles)
\begin{align*}
    \normal_{\pi_\algebr}(\algebr,\subalgebr) &\cong (\normal_{\pi_\algebr}(\algebr,\subalgebr) \times \mathbb{R})/\mathbb{R}^\times \\
    &\cong \normal_{\DNC(\pi_\algebr)}(\DNC_{\pi_\algebr}(\algebr,\subalgebr),\normal_{\pi_\algebr}(\algebr,\subalgebr))/\mathbb{R}^\times \\
    &\cong \normal_{\blup(\pi_\algebr)}(\blup_{\pi_\algebr}(\algebr,\subalgebr),\mathbb{P}_{\pi_\algebr}(\algebr,\subalgebr)) \cong \normal_{(\blup(\pi_\algebr) \times \pi_\algebrtwo)|_{\gr \varphi \circ \bldown}}(\gr \varphi \circ \bldown, \gr \varphi \circ \bldown|_{\mathbb{P}_{\pi_\algebr}(\algebr,\subalgebr)}).
\end{align*}
%where $q$ denotes the quotient map $\DNC_{\pi_\algebr}(\algebr,\subalgebr) \rightarrowdbl \blup_{\pi_\algebr}(\algebr,\subalgebr)$. 
where we used the canonical diffeomorphism $\normal(\DNC_\pi(\algebr,\subalgebr),\normal_\pi(\algebr,\subalgebr)) \cong \normal_\pi(\algebr,\subalgebr) \times \mathbb{R}$ (see Remark \ref{rema: normal of normal in DNC}) in the second diffeomorphism, and we used Proposition \ref{prop: induced G action on normal bundle} in the third diffeomorphism. %To see that this last map is a diffeomorphism, notice that $q$ can locally be written as the map
%\[(y,\xi,t) \mapsto
%\begin{cases}
%    (y,\tfrac{\xi_1}{\xi_i},\dots,\tfrac{\xi_{i-1}}{\xi_i},0,\tfrac{\xi_{i+1}}{\xi_i},\dots,\tfrac{\xi_{q}}{\xi_i}) & \text{if}\ t=0 \\
%    (y,\tfrac{\xi_1}{\xi_i},\dots,\tfrac{\xi_{i-1}}{\xi_i},t\xi_i,\tfrac{\xi_{i+1}}{\xi_i},\dots,\tfrac{\xi_{q}}{\xi_i}) & \text{if}\ t \neq 0
%\end{cases} = (y,\frac{\xi_1}{\xi_i},\dots,\frac{\xi_{i-1}}{\xi_i},t\xi_i,\frac{\xi_{i+1}}{\xi_i},\dots,\frac{\xi_{q}}{\xi_i}),\]
%so, in particular, we see that $d_\normal q$ is fiberwise an isomorphism, and for all $\lambda\in \mathbb{R}^\times$, $d_\normal q(b,\lambda\xi)1 = d_\normal q(b,\xi)\lambda$. %(note: in the charts, $\normal_\pi(\algebr,\subalgebr) \subset \DNC_\pi(\algebr,\subalgebr)$ is given by $t=0$, and $\normal(\blup_\pi(\algebr,\subalgebr),\mathbb{P}_\pi(\algebr,\subalgebr))$ has coordinates given by $(y,\widetilde{x}_i^1,\dots,\widetilde{x}_i^{i-1},d\widetilde{x}_i^i,\widetilde{x}_i^{i+1},\dots,\widetilde{x}_i^q)$, where $(y,\widetilde{x}_i)$ are the coordinates of $\blup_\pi(\algebr,\subalgebr)$). 
We can now equivalently show that 
\[\blup_{\blup(\pi_\algebr) \times \pi_\algebrtwo}(\gr \varphi \circ \bldown, \gr \varphi \circ \bldown|_{\mathbb{P}(\algebr,\subalgebr)}) \subset \blup_{\blup(\pi_\algebr) \times \pi_\algebrtwo}(\blup_{\pi_\algebr}(\algebr,\subalgebr) \times \algebrtwo, \blup_{\pi_\algebr}(\algebr,\subalgebr) \times \subalgebrtwo)\] 
is a Lie subalgebroid. And indeed, since $\bldown: \blup_{\pi_\algebr}(\algebr,\subalgebr) \ra \algebr$ is a morphism of Lie algebroids (see Corollary \ref{coro: blow-down map is a morphism of LA}), so are $\varphi \circ \bldown$ and $\varphi \circ \bldown|_{\mathbb{P}_{\pi_\algebr}(\algebr,\subalgebr)}$. This shows that $\gr \varphi \circ \bldown \subset \blup_{\pi_\algebr}(\algebr,\subalgebr) \times \algebrtwo$ and $\gr \varphi \circ \bldown|_{\mathbb{P}_{\pi_\algebr}(\algebr,\subalgebr)} \subset \mathbb{P}_{\pi_\algebr}(\algebr,\subalgebr) \times \subalgebrtwo$ are closed Lie subalgebroids with 
\[\gr \varphi \circ \bldown \cap (\blup_{\pi_\algebr}(\algebr,\subalgebr) \times \subalgebrtwo) = \gr \varphi \circ \bldown|_{\mathbb{P}_{\pi_\algebr}(\algebr,\subalgebr)}.\]
By Proposition \ref{prop: blup of subalgebroid is subalgebroid}, the result follows. 
\end{proof}
Again, following a similar argument shows that we can replace $\blup_{\pi_\algebr}$ and $\blup_{\pi_\algebrtwo}$ everywhere with $\DNC$ (or $\DNC_{\pi_\algebr}$ and $\DNC_{\pi_\subalgebr}$, respectively):
\begin{coro}\label{coro: LA morphism to DNC}
Let $\varphi: (\algebr,\subalgebr) \ra (\algebrtwo,\subalgebrtwo)$ be a morphism of Lie algebroids between pairs of Lie algebroids $(\algebroid,\subalgebroid)$ and $(\algebrtwoid,\subalgebrtwoid)$. Then $\DNC(\varphi): \DNC(\algebr,\subalgebr) \ra \DNC(\algebrtwo,\subalgebrtwo)$ is a morphism of Lie algebroids.
\end{coro}
We can relate the Lie algebroids $\DNC_{\pi}(\algebr,\subalgebr)$ and $\blup_{\pi}(\algebr,\subalgebr)$ as follows:
\begin{coro}\label{coro: quotient map DNC to blup lie algebroid}
The quotient map $\DNC(\algebr,\subalgebr) \setminus \subalgebr \times \mathbb{R} \rightarrowdbl \blup(\algebr,\subalgebr)$ restricts to a morphism of Lie algebroids
\[\DNC_\pi(\algebr,\subalgebr) \rightarrowdbl \blup_\pi(\algebr,\subalgebr)\]
over the quotient map  $\DNC(\base,\subbase) \setminus \subbase \times \mathbb{R} \rightarrowdbl \blup(\base,\subbase)$.
\end{coro}
\begin{proof}
Denote the quotient map $\DNC_\pi(\algebr,\subalgebr) \rightarrowdbl \blup_\pi(\algebr,\subalgebr)$ by $q$. We have to show that 
\[\gr q \subset \DNC_\pi(\algebr,\subalgebr) \times \blup_\pi(\algebr,\subalgebr)\]
is a Lie subalgebroid. Similar to the proof of Corollary \ref{coro: blow-down map is a morphism of LA}, recall that if we have a pair of manifolds $(M,N)$, and $Z$ is a manifold, then we have a canonical diffeomorphism
\[\DNC(M,N) \times Z \xra{\sim} \DNC(M \times Z, N \times Z).\]
%(see Proposition \ref{prop: blup(X x M, Y x M) = blup(X,Y) x M}).
Again, it is readily verified that, in case $M$, $N$, and $Z$ are Lie algebroids, then it is even an isomorphism of Lie algebroids, and the restriction to $\DNC_{\pi_M}(M,N) \times Z$ is an isomorphism of Lie algebroids onto $\DNC_{\pi_M \times \pi_Z}(M \times Z, N \times Z)$. %The result still holds if we replace ``$\blup$'' everywhere with ``$\DNC$'', so we obtain an isomorphism of Lie algebroids
So, we obtain an isomorphism of Lie algebroids
%\[\DNC_\pi(\algebr,\subalgebr) \times \blup_\pi(\algebr,\subalgebr) \xra{\sim} \DNC_{\blup(\pi \times \pi)}(\blup_{\pi \times \pi}(\algebr \times \algebr, \subalgebr \times \algebr), \blup_{\pi \times \pi}(\algebr \times \subalgebr, \subalgebr \times \subalgebr))\]
\[\DNC_\pi(\algebr,\subalgebr) \times \blup_\pi(\algebr,\subalgebr) \xra{\sim} \DNC_{\blup(\pi) \times \pi}(\blup_\pi(\algebr,\subalgebr) \times \algebr, \blup_\pi(\algebr,\subalgebr) \times \subalgebr).\]
Observe now that $\gr q$ is mapped bijectively onto $\DNC_{(\blup(\pi) \times \pi)|_{\gr \bldown}}(\gr \bldown, \gr \bldown|_{\mathbb{P}_\pi(\algebr,\subalgebr)})$: 
\begin{align*}
    (z,q(z)) =
    \begin{cases}
        (b,\xi,[b,\xi]) \\
        (a,t,[a,1])
    \end{cases} &\mapsto 
    \begin{cases}
        ([b,\xi],b,\xi) & \text{if}\ z=(b,\xi) \\
        ([a,1],a,t) & \text{if}\ z = (a,t),
    \end{cases}
\end{align*}
where we used that there is a canonical diffeomorphism
\[\normal_\pi(\algebr,\subalgebr) \cong \normal_{\blup(\pi)}(\blup_\pi(\algebr,\subalgebr),\mathbb{P}_\pi(\algebr,\subalgebr)) \cong \normal_{(\blup(\pi) \times \pi)|_{\gr \bldown}}(\gr \bldown, \gr \bldown|_{\mathbb{P}_\pi(\algebr,\subalgebr)}) %\textnormal{ given by } (b,\xi) \mapsto ([b,\xi],d_\normal q(b,\xi)1)
\]
(see the proof of Corollary \ref{coro: LA morphism to blup}). This shows that we can equivalently prove that 
\[\DNC_{(\blup(\pi) \times \pi)|_{\gr \bldown}}(\gr \bldown, \gr \bldown|_{\mathbb{P}_\pi(\algebr,\subalgebr)}) \subset \DNC_{\blup(\pi) \times \pi}(\blup_\pi(\algebr,\subalgebr) \times \algebr,\blup_\pi(\algebr,\subalgebr) \times \subalgebr)\] 
is a Lie subalgebroid. Since $\gr \bldown \cap (\blup_\pi(\algebr,\subalgebr) \times \subalgebr) = \gr \bldown|_{\mathbb{P}_\pi(\algebr,\subalgebr)}$, the result follows from Corollary \ref{coro: blow-down map is a morphism of LA} and Proposition \ref{prop: DNC of subalgebroid is subalgebroid}.
\end{proof}
Lastly, we present a simple example of a blow-up of Lie algebroids. It is the analogue of Example \ref{exam: blup pair groupoids} in the Lie algebroid setting.
\begin{exam}\label{exam: blup tangent bundle}
Recall from Proposition \ref{prop: Lie groupoid blup Lie algebroid}  (or Example \ref{exam: blup pair groupoids}) that the Lie algebroid of $\blup_{\source,\target}(\base \times \base, \subbase \times \subbase) \rra \blup(\base,\subbase)$ is 
\[\blup_{\pi_{T\base}}(T\base,T\subbase) \ra \blup(\base,\subbase).\]
Also, recall from Example \ref{exam: blup pair groupoids} that, as Lie groupoids, $\blup_{\source,\target}(\base \times \base, \subbase \times \subbase)$ is isomorphic to $(P \times_{\mathbb{R}} P)/\mathbb{R}^\times \rra P/\mathbb{R}^\times = \blup(\base,\subbase)$, where $P \coloneqq \DNC(\base,\subbase) \setminus \subbase \times \mathbb{R}$, which has Lie algebroid 
\[(\ker d\hat{t}|_P)/\mathbb{R}^\times \ra \blup(\base,\subbase)\]
(see Example \ref{exam: blup pair groupoids}), so we can identify $\blup_{\pi_{T\base}}(T\base,T\subbase)$ with $(\ker d\hat{t}|_P)/\mathbb{R}^\times \subset TP/\mathbb{R}^\times$. \

Observe that our findings in Example \ref{exam: blup pair groupoids} (the blow-up of pair groupoids) agree with our findings for $\algebr \coloneqq \blup_{\pi_{T\base}}(T\base,T\subbase)$. Indeed, 
%We will now prove in the Lie algebroid setting that the isotropy Lie algebras and orbits of $\algebr \coloneqq \blup_{\pi_{T\base}}(T\base,T\subbase)$ agree with our findings in Example \ref{exam: blup pair groupoids}.% Recall from Proposition \ref{prop: tanent bundle blow-up} (see also Remark \ref{rema: anchored vector bundle blow-up}) that the anchor map
\[\anchor: \blup_{\pi_{T\base}}(T\base,T\subbase) \ra T\blup(\base,\subbase)\]
is an isomorphism when restricted to $T\base|_{\base \setminus \subbase} \subset \blup_{\pi_{T\base}}(T\base,T\subbase)$ (onto $T\base|_{\base \setminus \subbase} \subset T\blup(\base,\subbase)$), and, for all $[y,\xi] \in \mathbb{P}(\base,\subbase)$, $\anchor([y,\xi])$ is surjective with a $1$-dimensional kernel. Therefore,
\[\mathfrak{g}_z = 0 \textnormal{ if } z \in \base\setminus\subbase; \quad \mathfrak{g}_z \cong \mathbb{R} \textnormal{ if } z \in \mathbb{P}(\base,\subbase)\]
Similarly, by Proposition \ref{prop: isotropy groups and orbits of blow-up algebroid}, we see that the orbits are given by the connected components of $\base \setminus \subbase$ and the connected components with $\mathbb{P}(\base,\subbase)$; this is also in line with Example \ref{exam: blup pair groupoids}.
\end{exam}

\subsection{Blow-ups: different approaches}\label{sec: Blow-up of a pair of manifolds: a different approach (part 1)}
We will now describe different approaches of constructing a blow-up of manifolds, and later in this section also of constructing a blow-up of Lie groupoids and Lie algebroids. The difference lies in the fact that we do not ``factor'' through the deformation to the normal cone, but rather describe the blow-up directly. This direct approach to blow-ups has its advantages. For example, in the next approach to blow-up of manifolds, the local description is very easy to understand. On the other hand, describing the blow-up globally (e.g. the smooth structure), is a bit trickier. 

\subsubsection{Blow-ups of manifolds: a different approach}
This construction is based on \cite{10.36045/bbms/1292334057}.
Let $(\base,\subbase)$ be a pair of manifolds. We first restrict our attention to the case where
\[\subbase \coloneqq \mathbb{R}^p \subset \mathbb{R}^{p+q} = \mathbb{R}^n \eqqcolon \base.\]
Denote by $S^{q-1}$ the unit sphere in $\mathbb{R}^q$. Observe that on 
\[\mathbb{R}^{p} \times S^{q-1} \times \mathbb{R}\]
(it is useful to think of this space as describing polar coordinates on $\mathbb{R}^n$) we have a canonical $\mathbb{Z}_2$-action given by 
\[\tau \cdot (x,\theta,t) = (x,-\theta,-t).\]
The action is obviously free (and proper, because $\mathbb{Z}_2$ is finite), so we obtain a quotient manifold 
\[\blupp(\mathbb{R}^n,\mathbb{R}^p) \coloneqq (\mathbb{R}^{p} \times S^{q-1} \times \mathbb{R})/\mathbb{Z}_2\]
(the subscript ``pc" refers to polar coordinates). To obtain a blow-down map, observe that the map
\[\Psi_{\textnormal{pc}}: \mathbb{R}^{p} \times S^{q-1} \times \mathbb{R} \ra \mathbb{R}^{p+q} \textnormal{ given by } (x,\theta,t) \mapsto (x,t \cdot \theta)\]
is smooth and intertwines the $\mathbb{Z}_2$-actions (with trivial action on $\mathbb{R}^n)$. It therefore descends to a smooth map, which we will call the blow-down map:
\[\bldownn: \blupp(\mathbb{R}^n,\mathbb{R}^p) \ra \mathbb{R}^n.\]
Observe that if we restrict $\Psi_{\textnormal{pc}}$ to $\mathbb{R}^p \times S^{q-1} \times \mathbb{R}_{>0}$, then this map $\Psi_{\textnormal{pc}}^{> 0}$ is a diffeomorphism with inverse
\[(\Psi_{\textnormal{pc}}^{> 0})^{-1}: \mathbb{R}^{p+q} \setminus \mathbb{R}^p \ra \mathbb{R}^p \times S^{q-1} \times \mathbb{R}_{>0} \textnormal{ given by } (x,y) \mapsto (x, \frac{y}{\Vert y \Vert}, \Vert y \Vert).\]
This shows that the blow-down map $\bldownn$ restricts to a diffeomorphism
\[\bldownn^{-1}(\mathbb{R}^n \setminus \mathbb{R}^p) \ra \mathbb{R}^n \setminus \mathbb{R}^p.\]
Moreover, if we restrict $\Psi_{\textnormal{pc}}$ to $\mathbb{R}^p \times S^{q-1} \times \{0\}$, then this map maps onto $\mathbb{R}^p$, and the induced map
\[\mathbb{R}^p \times \mathbb{RP}^{q-1} \xra{\pr_1} \mathbb{R}^p\]
can be seen as the restriction of $\bldownn$ to $\bldownn^{-1}(\mathbb{R}^p)$. 

Using the lemma below (compare it with Lemma \ref{lemm: transition maps of DNC are smooth}), we will be able to lift transition maps of a pair of smooth manifolds $(\base,\subbase)$ to diffeomorphisms between the blow-ups of their respective domains. Then we will assemble the blow-ups of the (ranges of the) charts as a disjoint union. We will then define a relation on this space which identifies points that are related by the induced diffeomorphisms (that the transition maps define), and the resulting quotient space will be the blow-up $\blupp(\base,\subbase)$.
\begin{lemm}\cite{10.36045/bbms/1292334057}\label{lemm: transition maps of blupp are smooth}
Let $h=(h_1,h_2): \mathbb{R}^{p+q} \ra \mathbb{R}^{k+\ell}$ be a smooth map such that $h^{-1}(\mathbb{R}^k) = \mathbb{R}^p$. Moreover, assume that $dh_2(x,0)$ is injective for all $x \in \mathbb{R}^p$ (i.e. $d_\normal h$ is fiber-wise injective). Then the map
\begin{align*}
    \widetilde{h}: \mathbb{R}^p \times S^{q-1} \times \mathbb{R} &\ra \mathbb{R}^k \times S^{\ell-1} \times \mathbb{R} \textnormal{ given by } \\
    (x,\theta,t) &\mapsto
\begin{cases}
    (h_1(x,0),\frac{\tfrac{\partial h_2}{\partial \theta}(x,0)\theta}{\Vert \tfrac{\partial h_2}{\partial \theta}(x,0)\theta \Vert},0) & \text{if}\ t = 0 \\
    (h_1(x,t\theta), \tfrac{t}{\vert t \vert}\tfrac{h_2(x,t\theta)}{\Vert h_2(x,t\theta) \Vert}, \tfrac{t}{\vert t \vert}\Vert h_2(x,t\theta) \Vert) & \text{if}\ t \neq 0
\end{cases}
\end{align*}
is smooth.
\end{lemm}
\begin{proof}
If we write $\widetilde{h} = (\widetilde{h}_1,\widetilde{h}_2,\widetilde{h}_3)$, it is obvious that $\widetilde{h}_1$ is smooth, so we only have to verify that $\widetilde{h}_2$ and $\widetilde{h}_3$ are smooth. Using that $h_2(x,0) = 0$, we have that
\[h_2(x,\theta) = \int_0^1 \frac{d}{dt}h_2(x,t\theta) dt = \int_0^1 \frac{\partial h_2}{\partial \theta}(x,t\theta)\theta dt = \frac{\partial h_2}{\partial \theta}(x,0)\theta + r(x,\theta)\theta,\]
where $r(x,\theta) \coloneqq -\textstyle\frac{\partial h_2}{\partial \theta}(x,0) + \textstyle\int_0^1\frac{\partial h_2}{\partial \theta}(x,t\theta)dt$, is a smooth map with the property that $r(x,0)=0$.
Hence,
\[t^{-1}h_2(x,t\theta) = \frac{\partial h_2}{\partial \theta}(x,0)\theta + r(x,t\theta)\theta\]
(note: $r(x,0)=0$).
Observe that this shows that $\widetilde{h}_2$ is the map
\[(x,\theta,t) \mapsto \frac{t^2 \cdot (\frac{\partial h_2}{\partial \theta}(x,0)\theta + r(x,t\theta)\theta)}{|t^2| \cdot \Vert \frac{\partial h_2}{\partial \theta}(x,0)\theta + r(x,t\theta)\theta \Vert} = \frac{\frac{\partial h_2}{\partial \theta}(x,0)\theta + r(x,t\theta)\theta}{|\Vert \frac{\partial h_2}{\partial \theta}(x,0)\theta + r(x,t\theta)\theta \Vert}\]
and $\widetilde{h}_3$ is the map 
\[(x,\theta,t) \mapsto t \cdot \Vert \frac{\partial h_2}{\partial \theta}(x,0)\theta + r(x,t\theta)\theta \Vert.\]
Since the norm $\Vert \cdot \Vert$, seen as a map $\mathbb{R}^d \ra \mathbb{R}$ (where $d\ge0$), is a smooth function away from $0 \in \mathbb{R}^d$, we see that it suffices to prove that the smooth map 
\[\mathbb{R}^p \times S^{q-1} \times \mathbb{R} \ra \mathbb{R}^\ell \textnormal{ given by } (x,\theta,t) \mapsto \frac{\partial h_2}{\partial \theta}(x,0)\theta + r(x,t\theta)\theta\]
is nowhere vanishing. Let $(x,\theta,t) \in \mathbb{R}^p \times S^{q-1} \times \mathbb{R}$. If $t \neq 0$, then 
\[\frac{\partial h_2}{\partial \theta}(x,0)\theta + r(x,t\theta)\theta = t^{-1}h_2(x,t\theta) \neq 0,\]
because $h_2^{-1}(\mathbb{R}^k) = \mathbb{R}^p$. If $t=0$, then, since $r(x,0)=0$, wee see that 
\[\frac{\partial h_2}{\partial \theta}(x,0)\theta + r(x,t\theta)\theta = \frac{\partial h_2}{\partial \theta}(x,0)\theta \neq 0,\]
because $\frac{\partial h_2}{\partial \theta}(x,0)\theta$ is injective by assumption. This proves that both $\widetilde{h}_2$ and $\widetilde{h}_3$ are smooth maps, so this proves the statement. 
\end{proof}
Immediately we obtain the following two corollaries:
\begin{coro}\cite{10.36045/bbms/1292334057}\label{coro: induced smooth maps on blupp}
Let $h$ be as in Lemma \ref{lemm: transition maps of blupp are smooth}. Then $h$ induces a smooth map 
\[\blupp(h): \blupp(\mathbb{R}^n,\mathbb{R}^p) \ra \blupp(\mathbb{R}^{k + \ell}, \mathbb{R}^k \times \{0\}).\]
\end{coro}
\begin{proof}
This follows from the fact that the map $\widetilde{h}$ from Lemma \ref{lemm: transition maps of blupp are smooth} intertwines the $\mathbb{Z}_2$-actions.
\end{proof}
\begin{coro}\label{coro: functoriality of smooth maps blupp}
Let $h: \mathbb{R}^{p+q} \ra \mathbb{R}^{k + \ell}$ and $h': \mathbb{R}^{k + \ell} \ra \mathbb{R}^{m+n}$ satisfy the hypothesis of Lemma \ref{lemm: transition maps of blupp are smooth}. Then $\blupp(h' \circ h) = \blupp(h') \circ \blupp(h)$.
\end{coro}
\begin{proof}
This follows from the definitions of $\widetilde{h}$ and $\widetilde{h}'$.
\end{proof}
\begin{rema}\cite{10.36045/bbms/1292334057}\label{rema: blupp behaves well with intersections}
Before we move to the general case of blow-ups of pairs of smooth manifolds, observe the following: if $U,U' \subset \mathbb{R}^{p+q}$ are open subsets, then we can view $\blupp(U,U \cap \mathbb{R}^p) \subset \blupp(\mathbb{R}^n, \mathbb{R}^p)$ and $\blupp(U',U' \cap \mathbb{R}^p) \subset \blupp(\mathbb{R}^n, \mathbb{R}^p)$ as open subsets, and
\[\blupp(U,U \cap \mathbb{R}^p) \cap \blupp(U',U' \cap \mathbb{R}^p) = \blupp(U \cap U',(U \cap U') \cap \mathbb{R}^p)\]
as subsets of $\blupp(\mathbb{R}^n, \mathbb{R}^p)$. This fact is readily verified by tracing the definitions.
\end{rema}
Now, if $(\base,\subbase)$ is a pair of smooth manifolds, let $\{(U_\alpha,\varphi_\alpha)\}_{\alpha \in I}$ be an adapted atlas (to $\subbase$) on $\base$ such that whenever $\subbase \cap U_\alpha \neq \emptyset$, we have $\varphi_\alpha(V_\alpha) = \varphi_\alpha(U_\alpha) \cap \mathbb{R}^p$ (where $V_\alpha \coloneqq U_\alpha \cap \subbase$). Let $\alpha,\beta \in I$ and write $U_{\alpha\beta}$ (and $V_{\alpha\beta}$) for the open set $U_\alpha \cap U_\beta$ (resp. $V_{\alpha} \cap V_{\beta}$), and $\varphi_{\alpha\beta}$ for the transition map $\varphi_\alpha \circ \varphi_\beta^{-1}$. From Corollary \ref{coro: induced smooth maps on blupp} we obtain the map
\[\blupp(\varphi_\beta(U_{\alpha\beta}), \varphi_\beta(U_{\alpha\beta} \cap \subbase)) \xra{\blupp(\varphi_{\alpha\beta})} \blupp(\varphi_\alpha(U_{\alpha\beta}), \varphi_\alpha(U_{\alpha\beta} \cap \subbase)).\]
Define 
\[\blupp^{\sim}(\base,\subbase) = \bigsqcup_{\alpha \in I} \blupp(\varphi_\alpha(U_\alpha),\varphi_\alpha(V_\alpha))\]
and define an equivalence relation $\sim$ on it such that 
\begin{align*}
    z &\sim w \textnormal{ (where } z \in \blupp(\varphi_\alpha(U_\alpha),\varphi_\alpha(V_\alpha)) \textnormal{ and } w \in \blupp(\varphi_\beta(U_\beta),\varphi_\beta(V_\beta))\textnormal{)} \textnormal{ if and only if } \\
    z &\in \blupp(\varphi_\alpha(U_{\alpha\beta}),\varphi_\alpha(V_{\alpha\beta})), w \in \blupp(\varphi_\beta(U_{\alpha\beta}),\varphi_\beta(V_{\alpha\beta})), \textnormal{ and } \blupp(\varphi_{\alpha\beta})(w)=z.
\end{align*} 
By Corollary \ref{coro: functoriality of smooth maps blupp} and Remark \ref{rema: blupp behaves well with intersections}, this is indeed an equivalence relation and the resulting quotient space
\[eq: \blupp^{\sim}(\base,\subbase) \ra \blupp(\base,\subbase)\]
is (set-theoretically) the blow-up of $\subbase$ in $\base$. We define a smooth structure on it by requiring that $eq|_{\blupp(\varphi_\alpha(U_\alpha),\varphi_\alpha(V_\alpha))}$ is a diffeomorphism onto its image for all $\alpha \in I$.
\begin{lemm}\cite{10.36045/bbms/1292334057}\label{lemm: smooth structure on blupp is well-defined}
The smooth structure on $\blupp(\base,\subbase)$ is well-defined.
\end{lemm}
\begin{proof}
First of all, the maps $eq|_{\blupp(\varphi_\alpha(U_\alpha),\varphi_\alpha(V_\alpha))}$ are injective; this is obvious by definition of the equivalence relation $\sim$. It follows that the equivalence relation is open. Indeed, %the union of sets maps to the union of the images under any map, so 
to prove that $eq$ is an open map, we may restrict to the open sets $\blupp(\varphi_\alpha(U_\alpha),\varphi_\alpha(V_\alpha))$, but in this case it is obvious. The maps $eq|_{\blupp(\varphi_\alpha(U_\alpha),\varphi_\alpha(V_\alpha))}$ are now seen to be open (continuous) embeddings into $\blupp(\base,\subbase)$. The smooth structure on $\blupp(\base,\subbase)$ is defined by taking as open cover the subsets $eq(\blupp(\varphi_\alpha(U_\alpha),\varphi_\alpha(V_\alpha)))$ together with the maps $eq|_{\blupp(\varphi_\alpha(U_\alpha),\varphi_\alpha(V_\alpha))}^{-1}$ post-composed with a chart on $\blupp(\varphi_\alpha(U_\alpha),\varphi_\alpha(V_\alpha))$. Observe that the charts will be compatible using Remark \ref{rema: blupp behaves well with intersections} and Corollary \ref{coro: functoriality of smooth maps blupp}. This proves the statement.
\end{proof}
\begin{rema}\cite{10.36045/bbms/1292334057}\label{rema: if Y Hdorff, then blupp Hdorff}
Observe that if $\base$ is Hausdorff, then $\blupp(\base,\subbase)$ will be Hausdorff. To show this, observe that it suffices to prove that if $z,w \in \blupp(\base,\subbase)$, then we must have $z,w \in eq(\blupp(\varphi_\alpha(U_\alpha),\varphi_\alpha(V_\alpha)))$ for some $\alpha \in I$. Indeed, if, say, $z \in eq(\blupp(\varphi_\beta(U_\beta),\varphi_\beta(V_\beta)))$ and $w \in eq(\blupp(\varphi_\gamma(U_\gamma),\varphi_\gamma(U_\gamma \cap \subbase)))$, then, since $eq^{-1}$ is well-defined when restricted to any of these open sets, we can consider the map
\[eq(\blupp(\varphi_\beta(U_\beta),\varphi_\beta(V_\beta))) \xra{eq^{-1}} \blupp(\varphi_\beta(U_\beta),\varphi_\beta(V_\beta)) \xra{\bldownn} \varphi_\beta(U_\beta) \xra{\varphi_\beta^{-1}} U_\beta \hookrightarrow \base\]
and similarly the map $eq(\blupp(\varphi_\gamma(U_\gamma),\varphi_\gamma(U_\gamma \cap \subbase))) \ra \base$. If the image of $z$ under the first map equals the image of $w$ under the second map, then there exists a chart $(U_\alpha,\varphi_\alpha)$ containing the common image. In this case it is then obvious that $z,w \in eq(\blupp(\varphi_\alpha(U_\alpha),\varphi_\alpha(V_\alpha)))$. If the images of $z$ and $w$ under the respective maps are different, then we can find $\alpha',\beta' \in I$ such that $z \in U_{\alpha'}$, $w \in U_{\beta'}$, and $U_{\alpha'} \cap U_{\beta'} = \emptyset$. Then, we may assume $\varphi_{\alpha'}(U_{\alpha'}) \cap \varphi_{\beta'}(U_\beta') = \emptyset$, so we can consider the chart $(U_{\alpha'} \cup U_{\beta'},\varphi_{\alpha'} \cup \varphi_{\beta'}) \eqqcolon (U_{\alpha},\varphi_{\alpha})$ (where $\varphi_{\alpha'} \cup \varphi_{\beta'}: U_{\alpha'} \cup U_{\beta'} \ra \varphi_{\alpha'}(U_{\alpha'}) \cup \varphi_{\beta'}(U_\beta')$ is the map that sends $a$ to $\varphi_{\alpha'}(a)$ and $b \in U_{\beta'}$ to $\varphi_{\beta'}(b)$). It is clear now that $z,w \in eq(\blupp(\varphi_\alpha(U_\alpha),\varphi_\alpha(V_\alpha)))$. This shows that indeed $\blupp(\base,\subbase)$ is Hausdorff if $\base$ is Hausdorff.
\end{rema}
We will now prove that this blow-up construction really is equivalent to our earlier blow-up construction. To do this, recall that the blow-up has a universal property (see Proposition \ref{prop: universal property blow-up}). So, we only have to prove that the universal property holds for our new blow-up construction. This will follow from the following lemma (compare with Lemma \ref{lemm: induced global smooth map on blow up}).
\begin{lemm}\cite{10.36045/bbms/1292334057}\label{lemm: induced global smooth map on blupp}
Whenever $f: (\base,\subbase) \ra (\basetwo,\subbasetwo)$ is a map of pairs such that (1) $\subbasetwo \subset \basetwo$ is of codimension $1$, (2) $f^{-1}(\subbasetwo)=\subbase$, and (3) $d_\normal f$ is fiberwise injective, then there is a unique smooth map $\blupp(f): \blupp(\base,\subbase) \ra \blupp(\basetwo,\subbasetwo)$ such that the diagram
\begin{center}
\begin{tikzcd}
    \blupp(\base,\subbase) \ar{r}{\blupp(f)} \ar{d}{\bldownn} & \blupp(\basetwo,\subbasetwo) \ar{d}{\bldownn} \\
    \base \ar{r}{f} & \basetwo
\end{tikzcd}
\end{center}
commutes.
\end{lemm}
\begin{proof}
Let $x \in \base$ and let $(U,\varphi)$ and $(V,\psi)$ be adapted charts around $x \in \base$ and $f(x) \in \basetwo$ respectively. Then 
\[(\psi \circ f \circ \varphi^{-1})^{-1}(\psi(V \cap \subbasetwo)) = \varphi(U \cap f^{-1}(V \cap \subbasetwo)) = \varphi(U \cap f^{-1}(V) \cap \subbase)\] by the assumption that $f^{-1}(\subbasetwo)=\subbase$. Therefore, by Corollary \ref{coro: induced smooth maps on blupp}, we obtain the smooth map $\blupp(\psi \circ f \circ \varphi^{-1})$ such that the diagram 
\begin{center}
\begin{tikzcd}
    \blupp(\varphi(U \cap f^{-1}(V)), \varphi(U \cap f^{-1}(V) \cap \subbase)) \ar{r}{\blupp(\psi \circ f \circ \varphi^{-1})} \ar{d}{\bldownn} & \blupp(\psi(V), \psi(V \cap \subbasetwo)) \ar{d}{\bldownn} \\
    \varphi(U \cap f^{-1}(V)) \ar{r}{\psi \circ f \circ \varphi^{-1}} & \psi(V)
\end{tikzcd}
\end{center}
commutes. Since $\blupp(\psi \circ f \circ \varphi^{-1}) = \blupp(\psi) \circ \blupp(f) \circ \blupp(\varphi)^{-1}$ by Corollary \ref{coro: functoriality of smooth maps blupp}, we see that $\blupp(f)$ is a smooth map which has the desired property.
\end{proof}
\begin{rema}\cite{10.36045/bbms/1292334057}\label{rema: blupp if submanifold is hypersurface}
We will discuss here the universal property of $\blupp(\base,\subbase)$ that are similar to $\blup(\base,\subbase)$. First of all, using notation as in the discussion after Remark \ref{rema: blupp behaves well with intersections}, notice that the blow-down maps 
\[(\bldownn)_\alpha: \blupp(\varphi_\alpha(U_\alpha),\varphi_\alpha(V_\alpha)) \ra \varphi_\alpha(U_\alpha)\]
assemble into a smooth map $\bldownn: \blupp(\base,\subbase) \ra \base$. The restriction of $\blupp^\sim(\base,\subbase)$ to 
\[\bigsqcup_{\alpha \in I} (\bldownn)_\alpha^{-1}(\varphi_\alpha(V_\alpha))\]
gives rise to a closed embedded submanifold $\mathbb{P}_{\textnormal{pc}}(\base,\subbase) \subset \blupp(\base,\subbase)$ which is of codimension $1$, and one can easily check that $\bldownn^{-1}(\subbase)=\mathbb{P}_{\textnormal{pc}}(\base,\subbase)$. We also have the following property: if $\subbase$ is of codimension $1$, then $\blupp(\base,\subbase)$ can be canonically identified with $\base$. Indeed, this follows simply by the fact that $\mathbb{R}^p \times S^{q-1} \times \mathbb{R} = \mathbb{R}^p \times \{\pm 1\} \times \mathbb{R}$. Now, like in Proposition \ref{prop: universal property blow-up}, one can check (locally) that $d_\normal \bldownn$ is fiberwise injective. From Lemma \ref{lemm: induced global smooth map on blupp}, and using the same argument as in Proposition \ref{prop: universal property blow-up}, we see that $(\blupp(\base,\subbase),\mathbb{P}_{\textnormal{pc}}(\base,\subbase))$ satisfies the same universal property as $(\blup(\base,\subbase),\mathbb{P}(\base,\subbase))$. 
\end{rema}
From the former remark, we see that $(\blupp(\base,\subbase),\mathbb{P}_{\textnormal{pc}}(\base,\subbase))$ and $(\blup(\base,\subbase),\mathbb{P}(\base,\subbase))$ are isomorphic as pairs of manifolds by a unique isomorphism. More explicitly, we have:
\begin{prop}\label{prop: blup = blupp}
The functors $\blup$ and $\blupp$ (from the category of pairs of manifolds, together with maps of pairs $f: (\base,\subbase) \ra (\basetwo,\subbasetwo)$ for which $f^{-1}(\subbasetwo)=\subbase$ and $d_\normal f$ is fiberwise injective, to the category of manifolds) are inverse to each other.
\end{prop}
\begin{proof}
Consider the morphism of functors $\mathfrak{B}: \blup \ra \blupp$ given by
\begin{align*}
    \mathfrak{B}(\base,\subbase)&: \blup(\base,\subbase) \ra \blupp(\base,\subbase); \\
    z &\mapsto
\begin{cases}
    [\varphi^1(y),\frac{d_\normal \varphi(y)\xi}{\Vert d_\normal \varphi(y)\xi \Vert},0] & \text{if}\ z=[y,\xi] \\
    [\varphi^1(x),\frac{t}{\vert t \vert}\frac{\varphi^2(x)}{\Vert \varphi^2(x) \Vert},\frac{t}{\vert t \vert}\Vert \varphi^2(x) \Vert] & \text{if}\ z=[x,t]
\end{cases}
\end{align*}
where $(U,\varphi)$ is an adapted chart around $\bldown(z)$. To see that this map is well defined, suppose $(V,\psi)$ is another adapted chart around $\bldown(z)$. If $z=[y,\xi]$, then 
\begin{align*}
    [\varphi^1(y),\frac{d_\normal \varphi(y)\xi}{\Vert d_\normal \varphi(y)\xi \Vert}] &= [\psi^1 \circ \varphi^{-1}(\varphi^1(y)),\frac{d_\normal(\psi \circ \varphi^{-1})(\varphi^1(y))}{\Vert d_\normal(\psi \circ \varphi^{-1})(\varphi^1(y)) \Vert} \circ \frac{d_\normal \varphi(y)\xi}{\Vert d_\normal \varphi(y)\xi \Vert}]  \\
    &= [\psi^1(y),\frac{d_\normal \psi(y)\xi}{\Vert d_\normal \psi(y)\xi \Vert}]
\end{align*} 
and if $z=[x,t]$, then, similarly,
\begin{align*}
    [\varphi^1(x),\frac{t}{\vert t \vert}\frac{\varphi^2(x)}{\Vert \varphi^2(x) \Vert},\frac{t}{\vert t \vert}\Vert \varphi^2(x) \Vert] %&= [\psi^{1} \circ \varphi^{-1}(\varphi^2(y)),\frac{\psi^2 \circ \varphi^{-1}(\varphi^2(y))}{\Vert \psi^2 \circ \varphi^{-1}(\varphi^2(y)) \Vert}\frac{t}{\vert t \vert}\frac{\varphi^2(y)}{\Vert \varphi^2(y) \Vert},\Vert \psi^{2} \circ \varphi^{-1}(\varphi^2(y))\Vert\frac{t}{\vert t \vert}\Vert \varphi^2(y) \Vert] \\
    &= [\psi^1(x),\frac{t}{\vert t \vert}\frac{\psi^2(x)}{\Vert \psi^2(x) \Vert},\frac{t}{\vert t \vert}\Vert \psi^2(x) \Vert].
\end{align*}
%(note: $[\cdot]$ is with respect to the equivalence relation $\sim$ defined on $\blupp^{\sim}(\base,\subbase)$). 
It remains to show that the map $\mathfrak{B}(\base,\subbase)$ is smooth for all pairs $(\base,\subbase)$, and that the diagrams 
\begin{center}
\begin{tikzcd}
\blup(\base,\subbase) \ar{r}{\blup(f)} \ar{d}{\mathfrak{B}(\base,\subbase)} & \blup(\basetwo,\subbasetwo) \ar{d}{\mathfrak{B}(\base',\subbase')} \\
\blupp(\base,\subbase) \ar{r}{\blupp(f)} & \blupp(\basetwo,\subbasetwo)
\end{tikzcd}
\end{center}
commute for all maps of pairs $f: (\base,\subbase) \ra (\basetwo,\subbasetwo)$ for which $f^{-1}(\subbasetwo)=\subbase$ and $d_\normal f$ is fiberwise injective. All of these statements follow from the observation that the map $\mathfrak{B}(\base,\subbase)$ is obtained from Lemma \ref{lemm: induced global smooth map on blupp} applied to the map of pairs
\[\bldown: (\blup(\base,\subbase),\mathbb{P}(\base,\subbase)) \ra (\base,\subbase).\]
This morphism of functors has an inverse $\mathfrak{B}^{-1}: \blupp \ra \blup$ given by
\begin{align*}
    \mathfrak{B}(\base,\subbase)&: \blupp(\base,\subbase) \ra \blup(\base,\subbase); \\
    z=[x,\theta,t] &\mapsto
\begin{cases}
    [\varphi^{-1}(x),d_\normal\varphi^{-1}(x)\theta,0] & \text{if}\ t=0 \\
    [\varphi^{-1}(x,t \cdot \theta),t] & \text{if}\ t\neq0.
\end{cases}
\end{align*}
where $(U,\varphi)$ is an adapted chart around $\bldownn(z)$. Using the same argument, but now applying Lemma \ref{lemm: induced global smooth map on blow up} to $\bldownn$, the result follows. 
\end{proof}
We end this section with two remarks on other variations of blow-ups of manifolds. 
\begin{rema}\label{rema: other blow-up constructions}
Another way to approach blow-ups, is by mimicking the algebraic blow-up construction explained in Section \ref{sec: Blow-up in the algebraic setting}: since we know that blow-ups have to satisfy
\[\blup(\base \times \basetwo, \subbase \times \basetwo) \cong \blup(\base,\subbase) \times \basetwo,\]
we only have to describe, for the local description, the case
\[\subbase \coloneqq \{0\} \subset \mathbb{R}^n \eqqcolon \base.\]
To do this, observe that
\[\blup_{\textnormal{ag}}(\mathbb{R}^n,0) \coloneqq \{((x_1,\dots,x_n),[y_1:\cdots:y_n]) \mid x_iy_j=y_ix_j \textnormal{ for all } 1 \le i,j \le n\} \subset \mathbb{R}^n \times \mathbb{RP}^{n-1}\]
(the subscript ``ag'' refers to algebraic geometry) is a closed embedded submanifold of codimension $n-1$ (see Proposition \ref{prop: blup 0 in R^n is same as our construction}). 
We now define $\blup_{\textnormal{ag}}(\mathbb{R}^n,\mathbb{R}^p) \coloneqq \blup_{\textnormal{ag}}(\mathbb{R}^q,0) \times \mathbb{R}^p$. That is, $\blup_{\textnormal{ag}}(\mathbb{R}^n,\mathbb{R}^p)$ is equal to
\[\{((x_1,\dots,x_n),(y_{p+1}:\cdots:y_n)) \mid x_iy_j=y_ix_j \textnormal{ for all } p+1 \le i,j \le n\} \subset \mathbb{R}^n \times \mathbb{RP}^{n-p-1}.\]
Now observe that
\[S^{n-1} \times \mathbb{R} \ra \blup_{\textnormal{ag}}(\mathbb{R}^n,0) \textnormal{ given by } (\theta,t) \mapsto (t\theta,\theta)\]
(where we view $t\theta \in \mathbb{R}^n$ and $\theta \in \mathbb{RP}^{n-1}$) is $\mathbb{Z}_2$-equivariant, so it descends to a smooth map 
\[\blupp(\mathbb{R}^n,0) \ra \blup_{\textnormal{ag}}(\mathbb{R}^n,0),\]
which is readily verified to be a diffeomorphism. So, by tracing the steps we made in this section, and using the above diffeomorphism repeatedly, we can proceed to prove the analogues of the statements for ``$\blup_{\textnormal{ag}}$'' to obtain another way to approach blow-ups for general pairs of manifolds $(\base,\subbase)$.
\end{rema}
\begin{rema}\cite{proper}\label{rema: tubular neighbourhood blow-up}
In view of Remark \ref{rema: tubular neighbourhood DNC}, there is yet another way to perform the blow-up construction. Recall that we can glue manifolds $M_1$ and $M_2$ along a diffeomorphism $\varphi: U_1 \ra U_2$, defined between open subsets $U_1 \subset M_1$ and $U_2 \subset M_2$ to a new manifold $M_1 \cup_{\varphi} M_2$. Using a tubular neighbourhood $\chi: (\normal(\base,\subbase),0_\subbase) \ra (\base,\subbase)$, we can describe the blow-up using this gluing construction (compare with Remark \ref{rema: tubular neighbourhood DNC}). In fact, what we will do is mimic the constructions explained in this section, and apply them to the normal bundle instead of locally to $\mathbb{R}^n$. First observe that
\[\widetilde{\normal}(\base,\subbase) \coloneqq \{(y,\xi_1,[y,\xi_2]) \in \normal(\base,\subbase) \tensor[_{\pi_\subbase^\normal}]{\times}{_{\pi_\subbase^{\mathbb{P}}}} \mathbb{P}(\base,\subbase) \mid \xi_1=\lambda \xi_2 \textnormal{ for some } \lambda \in \mathbb{R}^\times\}\]
(where $\pi_\subbase^\normal: \normal(\base,\subbase) \ra \subbase$ and $\pi_\subbase^{\mathbb{P}}: \mathbb{P}(\base,\subbase) \ra \subbase$ are the projection maps)
is a smooth manifold (note: locally, it is precisely given by the blow-up ``$\blup_{\textnormal{ag}}(\base,\subbase)$'' from Remark \ref{rema: other blow-up constructions}). Alternatively, we can choose a metric on $\base$ and realise the sphere bundle
\[S(\base,\subbase) \coloneqq S(\normal(\base,\subbase)) = (\normal(\base,\subbase) \setminus 0_\subbase)/\mathbb{R}_{>0},\]
where $\mathbb{R}_{>0}$ acts on $\normal(\base,\subbase) \setminus 0_\subbase$ by scalar multiplication (this construction works for more general vector bundles; see Lemma \ref{lemm: free and proper scalar multiplication action on E}), as $\{(y,\xi) \in \normal(\base,\subbase) \mid \Vert\xi\Vert=1\}$. We have a free (and proper) $\mathbb{Z}_2$-action on $S(\base,\subbase) \times \mathbb{R}$ given by 
\[\tau \cdot (y,\xi,t) = (y,-\xi,-t).\] 
Then the map
\[S(\base,\subbase) \times \mathbb{R} \ra \normal(\base,\subbase) \tensor[_{\pi_\subbase^\normal}]{\times}{_{\pi_\subbase^{\mathbb{P}}}} \mathbb{P}(\base,\subbase) \textnormal{ given by } (y,\xi,t) \mapsto (y,t \cdot \xi, [y,\xi])\]
is smooth and intertwines the $\mathbb{Z}_2$-action (with trivial action on $\normal(\base,\subbase) \tensor[_{\pi_\subbase^\normal}]{\times}{_{\pi_\subbase^{\mathbb{P}}}} \mathbb{P}(\base,\subbase)$), so it descends to a smooth map
\[(S(\base,\subbase) \times \mathbb{R})/\mathbb{Z}_2 \ra \normal(\base,\subbase) \tensor[_{\pi_\subbase^\normal}]{\times}{_{\pi_\subbase^{\mathbb{P}}}} \mathbb{P}(\base,\subbase).\]
This map is an embedding, and it maps onto $\widetilde{\normal}(\base,\subbase)$ (note: locally, it is precisely given by the blow-up ``$\blupp(\base,\subbase)$'' explained in this section). Now, $\widetilde{\normal}(\base,\subbase)$ comes with a smooth map 
\[\widetilde{\normal}(\base,\subbase) \xra{\pr_1} \normal(\base,\subbase)\]
which is a diffeomorphism onto $\normal(\base,\subbase) \setminus 0_\subbase$ when restricted to $\widetilde{\normal}(\base,\subbase) \setminus (0_\subbase \tensor[_{\pi_\subbase^\normal}]{\times}{_{\pi_\subbase^{\mathbb{P}}}} \mathbb{P}(\base,\subbase))$. Therefore, we can glue $\widetilde{\normal}(\base,\subbase)$ and $\base \setminus \subbase$ via the map
\[\widetilde\chi \coloneqq \chi \circ \pr_1: \widetilde{\normal}(\base,\subbase) \setminus (0_\subbase \tensor[_{\pi_\subbase^\normal}]{\times}{_{\pi_\subbase^{\mathbb{P}}}} \mathbb{P}(\base,\subbase)) \ra \base \setminus \subbase\]
(which is a diffeomorphism onto an open subset of $\base \setminus \subbase$) to the smooth manifold
\[\blup_{\textnormal{tub}}(\base,\subbase) \coloneqq \widetilde\normal(\base,\subbase) \cup_{\widetilde\chi} \base \setminus \subbase.\]
We claim that this manifold is diffeomorphic to $\blup(\base,\subbase)$.To see this, notice that the blow-down map for $\blup_{\textnormal{tub}}(\base,\subbase)$ is given by 
\[\bldown_{\textnormal{tub}} \coloneqq \chi \circ \pr_1 \cup_{\widetilde\chi} \iota: \blup_{\textnormal{tub}}(\base,\subbase) \ra \base,\]
where $\iota: \base \setminus \subbase \hookrightarrow \base$ is the inclusion map (note: $\chi \circ \pr_1$ has domain $\widetilde\normal(\base,\subbase)$).
Notice now that $\widetilde\normal(\base,\subbase)$ contains the codimension $1$ submanifold $0_\subbase \tensor[_{\pi_\subbase^\normal}]{\times}{_{\pi_\subbase^{\mathbb{P}}}} \mathbb{P}(\base,\subbase)$, that $\bldown_{\textnormal{tub}}^{-1}(\subbase) = 0_\subbase \tensor[_{\pi_\subbase^\normal}]{\times}{_{\pi_\subbase^{\mathbb{P}}}} \mathbb{P}(\base,\subbase)$, and that $d_\normal\bldown_{\textnormal{tub}}$ is fiberwise injective (since $\chi \circ \pr_1$ is), so we obtain a smooth map
\[\blup(\bldown_{\textnormal{tub}}): \blup(\blup_{\textnormal{tub}}(\base,\subbase),0_\subbase \tensor[_{\pi_\subbase^\normal}]{\times}{_{\pi_\subbase^{\mathbb{P}}}} \mathbb{P}(\base,\subbase)) \cong \blup_{\textnormal{tub}}(\base,\subbase) \ra \blup(\base,\subbase).\]
Notice that it restricts to the canonical diffeomorphism $0_\subbase \tensor[_{\pi_\subbase^\normal}]{\times}{_{\pi_\subbase^{\mathbb{P}}}} \mathbb{P}(\base,\subbase) \cong \mathbb{P}(\base,\subbase)$. Moreover, the inclusion
\[\base \setminus \subbase \hookrightarrow \blup_{\textnormal{tub}}(\base,\subbase)\]
extends uniquely to a smooth map $\blup(\base,\subbase) \ra \blup_{\textnormal{tub}}(\base,\subbase)$ that restricts to the canonical diffeomorphism $\mathbb{P}(\base,\subbase) \cong 0_\subbase \tensor[_{\pi_\subbase^\normal}]{\times}{_{\pi_\subbase^{\mathbb{P}}}} \mathbb{P}(\base,\subbase)$. This shows that $\blup(\bldown_{\textnormal{tub}})$ is a diffeomorphism. %Lastly, we want to mention that we can, equivalently, start with a tubular neighbourhood $\chi$ which is a diffeomorphism from an open subset $W \supset 0_\subbase$ of $\normal(\base,\subbase)$ to an open subset $U \supset \subbase$ of $\base$. In that case one should replace $\widetilde\normal(\base,\subbase)$ with the open subset
%\[\widetilde W \coloneqq \widetilde\normal(\base,\subbase) \cap (W \tensor[_{\pi_\subbase^\normal}]{\times}{_{\pi_\subbase^{\mathbb{P}}}} \mathbb{P}(\base,\subbase)).\]
%of $\widetilde\normal(\base,\subbase)$ (see \cite{proper}).
\end{rema}

\subsubsection{Blow-ups of groupoids and algebroids: a different approach}\label{sec: Blow-up of a pair of manifolds: a different approach (part 2)}
In \cite{gualtieri2012symplectic}, a different blow-up construction is explained for Lie groupoids. There, they are only interested in blowing up Lie groupoids for which the base manifold does not change: this happens if and only if the submanifold of the base manifold is of codimension $1$. We will phrase the differences of our earlier blow-up construction for Lie groupoids (from \cite{2017arXiv170509588D}), and the one from \cite{gualtieri2012symplectic} by identifying two main differences. First of all, the blow-up construction used there, on the level of manifolds, is the construction from Section \ref{sec: Blow-up of a pair of manifolds: a different approach (part 1)} (see also Remark \ref{rema: other blow-up constructions}). This is one of the main differences. For the other main difference, note that the key to defining blow-ups for both Lie groupoids and Lie algebroids, lies in Remark \ref{rema: smooth map gives rise to smooth map of Blup}: it says that if $f: (\base, \subbase) \ra (\basetwo,\subbasetwo)$ is a smooth map of pairs, then we obtain canonically a smooth map $\blup(f): \blup_f(\base,\subbase) \ra \blup(\basetwo,\subbasetwo)$, where 
\[\blup_f(\base,\subbase) \coloneqq (\DNC(\base,\subbase) \setminus \DNC(f)^{-1}(\subbasetwo \times \mathbb{R}))/\mathbb{R}^\times.\] 
If we find the analogue of this Remark in the setting of Section \ref{sec: Blow-up of a pair of manifolds: a different approach (part 1)} (that is, by not having to factor through the deformation to the normal cone), then one can, again, define the structure maps this way, and proceed similarly as in the proofs of Theorem \ref{theo: groupoid blow-up}, Theorem \ref{theo: vector bundle blow-up}, and Theorem \ref{theo: Lie algebroid blow-up}. Notice that in \cite{gualtieri2012symplectic}, there is no geometric construction for the blow-up of Lie algebroids. They do, however, give the description on the level of sections (in the codimension $1$ case; it is called the \textit{lower elementary modification}). For a pair of Lie algebroids $(\algebroid,\subalgebroid)$, this space of sections is $\Gamma(\algebr,\subalgebr)$, and we have seen in Theorem \ref{theo: Lie algebroid blow-up}, that, indeed, our Lie algebroid construction has $\Gamma(\algebr,\subalgebr)$ as its space of sections (in the codimension $1$ case). 

We will now prove the analogue of Remark \ref{rema: smooth map gives rise to smooth map of Blup} in the setting of Section \ref{sec: Blow-up of a pair of manifolds: a different approach (part 1)}. We fix here a smooth map of pairs $f: (\base, \subbase) \ra (\basetwo,\subbasetwo)$.
\begin{prop}\label{prop: analogue functoriality of map of pairs blow-up}
We have a canonical smooth map 
\begin{align*}
    \blupp(f): \blup_{\textnormal{pc}, f}(\base,\subbase) \ra \blupp(\basetwo,\subbasetwo) 
\end{align*}
where 
\[\blup_{\textnormal{pc}, f} \coloneqq \blupp(\base,\subbase) \setminus \overline{\bldownn^{-1}(f^{-1}(\subbasetwo) \setminus \subbase)}.\] 
Moreover, under the canonical diffeomorphism $\blupp(\base,\subbase) \cong \blup(\base,\subbase)$ from Proposition \ref{prop: blup = blupp}, we have $\blup_{\textnormal{pc}, f} \cong \blup_f(\base,\subbase)$.
\end{prop}
\begin{proof}
Since we already proved the analogue of this statement for our earlier blow-up construction, it suffices to prove the second statement. We denote $Z \coloneqq \bldown^{-1}(f^{-1}(\subbasetwo) \setminus \subbase) \subset \blup(\base,\subbase)$. The statement follows if we prove that
\[\blup_f(\base,\subbase) = \blup(\base,\subbase) \setminus \overline{Z}.\]
Since $\bldown^{-1}(\subbase) = \mathbb{P}(\base,\subbase)$, and $f^{-1}(\subbasetwo) \subset \base$ is closed, we see that $(\base \setminus \subbase) \cap \overline{Z} = f^{-1}(\subbasetwo) \setminus \subbase$. Therefore, we only have to prove the following claim: 
\begin{equation}\label{eq: claim in proper transform}
    \textnormal{ if } z \coloneqq [y,\xi] \in \mathbb{P}(\base,\subbase), \textnormal{ then } z \textnormal{ is in the closure of } Z \textnormal{ if and only if } d_\normal f(y)\xi = 0.
\end{equation}
We fix an element $z \coloneqq [y,\xi] \in \mathbb{P}(\base,\subbase)$, and denote $q_\DNC: P \rightarrowdbl \blup(\base,\subbase)$ for the quotient projection, where $P \coloneqq \DNC(\base,\subbase) \setminus (\subbase \times \mathbb{R})$.

First, we will show that, if $d_\normal f(y)\xi \neq 0$, then $z$ is not in the closure of $Z$. Notice that $z$ is in the closure of $Z$ if and only if every $\mathbb{R}^\times$-invariant open subset $U \ni (y,\xi)$ of $P$ intersects 
\[Z_\DNC \coloneqq q_\DNC^{-1}(Z) = (f^{-1}(\subbasetwo) \setminus \subbase) \times \mathbb{R}^\times = f^{-1}(\subbasetwo) \times \mathbb{R}^\times \setminus (\subbase \times \mathbb{R}^\times) \subset P.\]
We will show that if $d_\normal f(y)\xi \neq 0$, then we can find an open invariant subset of $P$ that is disjoint from $Z_\DNC$. Indeed, $\DNC(f)$ is continuous, so we can find for all invariant open subsets $V$ around $\DNC(f)(y,\xi)$ an open subset $U \subset P$ around $(y,\xi)$ such that $\DNC(f)(U) \subset V$ (we use here that $\subbase \times \mathbb{R} \subset \DNC(\base,\subbase)$ is closed). By equivariance of $\DNC(f)$ and $V$, $\DNC(f)(\textstyle\cup_{\lambda \in \mathbb{R}^\times} \lambda \cdot U) \subset V$, so we can assume $U$ is invariant as well. By assumption, $\DNC(f)(y,\xi)=d_\normal f(y)\xi \not\in \subbasetwo \times \mathbb{R}$, so we may assume $V$ is disjoint from $\subbasetwo \times \mathbb{R}$ (because $\subbasetwo \times \mathbb{R}$ is closed in $\DNC(\basetwo,\subbasetwo)$). But now $\DNC(f)(U)$ is disjoint from $\subbasetwo \times \mathbb{R}$, so $U$ is disjoint from $Z_\DNC$.

Conversely, suppose $d_\normal f(y)\xi=0$. If $z$ is not in the closure of $Z$, then there is an invariant open subset $W \ni (y,\xi)$ of $P \coloneqq \DNC(\base,\subbase) \setminus (\subbase \times \mathbb{R})$ that is disjoint from $Z_\DNC$. By shrinking $W$ if necessary, we may assume $W = \DNC(U,V) \cap P$ for some adapted chart $(U,\varphi)$ of $\base$ (with $V \coloneqq U \cap \subbase$). Now, since $d_\normal f(y)\xi = 0$, $\xi \in T_y\base$ spans a one-dimensional subspace $\langle \xi \rangle$ of $T_y\base$ which maps into $T_{f(y)}\subbasetwo$, but it has a trivial intersection with $T_y\subbase$. We can integrate $\langle \xi \rangle$ to an embedded $1$-dimensional submanifold $S$ of $U$. By construction, we can shrink $S$ such that $f$ maps $S$ into $\subbasetwo$, but this makes it clear that $S \times \mathbb{R}^\times \subset W$ has a non-trivial intersection with $Z_\DNC$. This contradicts our assumption on $W$. This proves \eqref{eq: claim in proper transform}, so this concludes the proof.
\end{proof}
\begin{defn}\cite{gualtieri2012symplectic}\label{defn: proper transform}
Let $(\base,\subbase)$ be a pair of manifolds, and let $S \subset \base$. The closed set
\[\overoverline{S} \coloneqq \overline{\bldown^{-1}(S \setminus \subbase)} \subset \blup(\base,\subbase)\]
is called the \textit{proper transform} of $S$.
\end{defn}
With this notation, if $f: (\base,\subbase) \ra (\basetwo,\subbasetwo)$ is a map of pairs, then we have 
\[\blup_f(\base,\subbase)=\blup(\base,\subbase) \setminus \overoverline{f^{-1}(\subbasetwo)}; \quad \blup_{f,g}(\base,\subbase)=\blup(\base,\subbase) \setminus (\overoverline{f^{-1}(\subbasetwo)} \cup \overoverline{g^{-1}(\subbasetwo)}).\]
In Definition \ref{defn: blow-up 0 in affine variety}, which is the analogue of Definition \ref{defn: proper transform} for varieties (in case of blowing up a point), we already hinted upon thinking about this as another way to blow-up:
\begin{prop}\label{prop: proper transfrom = blup}
Let $S \subset \base$ be a closed embedded submanifold such that $\subbase \subset S$. Then
\[\blup(S,\subbase) = \overoverline{S}\]
viewed as subspaces of $\blup(\base,\subbase)$ (see also Corollary \ref{coro: blup submanifold}).
\end{prop}
\begin{proof}
Since $\bldown^{-1}(S \setminus \subbase) \subset \blup(S,\subbase)$, and $\blup(S,\subbase) \subset \blup(\base,\subbase)$ is closed (see Corollary \ref{coro: blup submanifold}), we see that $\overoverline{S} = \overline{\bldown^{-1}(S \setminus \subbase)} \subset \blup(S,\subbase)$. For the converse inclusion, observe that $\bldown^{-1}(S \setminus \subbase)$ is dense as a subspace of $\blup(S,\subbase)$, so if $[y,\xi] \in \mathbb{P}(S,\subbase) \subset \blup(S,\subbase)$, and $U \ni [y,\xi]$ is an open subset of $\blup(\base,\subbase)$, then $U \cap \blup(S,\subbase)$ intersects $\bldown^{-1}(S \setminus \subbase)$ non-trivially, which shows that $\blup(S,\subbase) \subset \overoverline{S}$. This proves the statement.
\end{proof}
%\begin{rema}\label{rema: immersed submanifold blow-up}

%\end{rema}
%Theoretically, this statement is quite handy, but in explicit examples it is not too useful in general: we would first have to understand the blow-up of $\subbase \subset \base$ (using notation from above), which is in general ``more complicated'' than the blow-up of $\subbase \subset S$ since the codimension of $\subbase \subset \base$ is higher than or equal to the codimension of $\subbase \subset S$ (think for example about blowing up a point in $\mathbb{R}^3$ versus blowing up a point in $S^2$). 

To illustrate the ideas presented in this section and the former section, let us repeat the simple example explained in Example \ref{exam: blow-up point in R^2}. 
\begin{exam}\label{exam: blow-up point in R^2 2}
We will reprove that the blow-up of $\{(0,0)\} \rra \{0\}$ in the pair groupoid $\mathbb{R} \times \mathbb{R} \rra \mathbb{R}$ yields $\mathbb{R}^\times \ltimes \mathbb{R}$. Indeed, recall that
\[\blupp(\mathbb{R}^2,0) = (S^1 \times \mathbb{R})/\mathbb{Z}_2,\]
and 
\[\bldownn: \blupp(\mathbb{R}^2,0) \ra \mathbb{R}^2 \textnormal{ is given by } [\theta,t] \mapsto t\theta.\]
From this we see that
\begin{align*}
    \bldownn^{-1}(\pr_2^{-1}(0) \setminus \{(0,0)\}) &= \{[\theta,t] \mid \theta=(\pm1,0), t\neq0\}; \\ 
    \bldownn^{-1}(\pr_1^{-1}(0) \setminus \{(0,0)\}) &= \{[\theta,t] \mid \theta=(0,\pm1), t\neq0\},
\end{align*}
and, in turn, this shows that
\[\overoverline{\pr_2^{-1}(0)} = \{[\theta,t] \mid \theta=(\pm1,0)\}; \quad \overoverline{\pr_1^{-1}(0)} = \{[\theta,t] \mid \theta=(0,\pm1)\}.\]
Therefore, 
\[\blup_{\textnormal{pc}, \source,\target}(\mathbb{R}^2,0) = ((S^1 \setminus \{\theta = (\theta_1,\theta_2) \in S^1 \mid \theta_1 \neq 0, \theta_2 \neq 0\}) \times \mathbb{R})/\mathbb{Z}_2.\]
Now, the smooth map
\[(S^1 \setminus \{\theta \in S^1 \mid \theta_1 \neq 0, \theta_2 \neq 0\}) \times \mathbb{R} \ra \mathbb{R}^\times \times \mathbb{R} \textnormal{ given by } (\theta,t) \mapsto (\frac{\theta_1}{\theta_2},t\theta_2)\]
is $\mathbb{Z}_2$-equivariant (with trivial action on $\mathbb{R}^\times \times \mathbb{R}$), so it descends to a smooth map
\begin{equation}\label{eq: diff to R^times times R}
    \blup_{\textnormal{pc},\source,\target}(\mathbb{R}^2,0) \ra \mathbb{R}^\times \times \mathbb{R},
\end{equation}
and this map has a smooth inverse, namely $(\lambda,t) \mapsto [\tfrac{1}{\sqrt{1+\lambda^{-2}}}, \tfrac{\lambda^{-1}}{\sqrt{1+\lambda^{-2}}}, t\tfrac{\lambda}{\sqrt{1+\lambda^{-2}}}+t\tfrac{\lambda^{-1}}{\sqrt{1+\lambda^{-2}}}]$; here is a direct computation:
\begin{align*}
    [\frac{1}{\sqrt{1+\lambda^{-2}}}&, \frac{\lambda^{-1}}{\sqrt{1+\lambda^{-2}}}, t\frac{\lambda}{\sqrt{1+\lambda^{-2}}}+t\frac{\lambda^{-1}}{\sqrt{1+\lambda^{-2}}}] \\
    &\mapsto (\frac{\frac{1}{\sqrt{1+\lambda^{-2}}}}{\frac{\lambda^{-1}}{\sqrt{1+\lambda^{-2}}}},(t\frac{\lambda}{\sqrt{1+\lambda^{-2}}}+t\frac{\lambda^{-1}}{\sqrt{1+\lambda^{-2}}}) \cdot \frac{\lambda^{-1}}{\sqrt{1+\lambda^{-2}}}) = (\lambda,t), \textnormal{ and } \\
    (\frac{\theta_1}{\theta_2},t\theta_2) 
    &\mapsto [\frac{1}{\sqrt{1+\frac{\theta_2^2}{\theta_1^2}}}, \frac{\frac{\theta_2}{\theta_1}}{\sqrt{1+\frac{\theta_2^2}{\theta_1^2}}}, t\theta_2\frac{\frac{\theta_1}{\theta_2}}{\sqrt{1+\frac{\theta_2^2}{\theta_1^2}}}+t\theta_2\frac{\frac{\theta_2}{\theta_1}}{\sqrt{1+\frac{\theta_2^{2}}{\theta_1^2}}}] \\
    &= [\frac{\textnormal{sign}(\theta_1)\theta_1}{\sqrt{\theta_1^2+\theta_2^2}}, \frac{\textnormal{sign}(\theta_1)\theta_2}{\sqrt{\theta_1^2+\theta_2^2}}, t\frac{\textnormal{sign}(\theta_1)\theta_1^2}{\sqrt{\theta_1^2+\theta_2^2}}+t\frac{\textnormal{sign}(\theta_1)\theta_2^2}{\sqrt{\theta_1^2+\theta_2^2}}] = [\theta,t].
\end{align*}
To see that the diffeomorphism \eqref{eq: diff to R^times times R} is a Lie groupoid isomorphism, we will only check %(as in Example \ref{exam: blow-up point in R^2}) 
that, under this diffeomorphism, $\blup(\source)$ and $\blup(\target)$ become the source map and the target map of $\mathbb{R}^\times \ltimes \mathbb{R}$, respectively (that this works for the other structure maps as well is readily verified by doing similar computations). Indeed,
\[\blup(\source)([\theta,t])=t\theta_2; \quad \blup(\target)([\theta,t])=t\theta_1\] (see Lemma \ref{lemm: transition maps of blupp are smooth}), and, indeed, 
\[\source_{\mathbb{R}^\times \ltimes \mathbb{R}}(\frac{\theta_1}{\theta_2},t\theta_2)=t\theta_2; \quad \target_{\mathbb{R}^\times \ltimes \mathbb{R}}(\frac{\theta_1}{\theta_2},t\theta_2) = t\theta_1,\]
and
\begin{align*}
    \blup(\source)([\frac{1}{\sqrt{1+\lambda^{-2}}}&, \frac{\lambda^{-1}}{\sqrt{1+\lambda^{-2}}}, t\frac{\lambda}{\sqrt{1+\lambda^{-2}}}+t\frac{\lambda^{-1}}{\sqrt{1+\lambda^{-2}}}]) \\
    &= (t\frac{1}{\sqrt{1+\lambda^{-2}}}+t\frac{\lambda^{-1}}{\sqrt{1+\lambda^{-2}}}) \cdot \frac{\lambda^{-1}}{\sqrt{1+\lambda^{-2}}} = t \\
    \blup(\target)([\frac{\lambda}{\sqrt{1+\lambda^{-2}}}&, \frac{\lambda^{-1}}{\sqrt{1+\lambda^{-2}}}, t\frac{\lambda}{\sqrt{1+\lambda^{-2}}}+t\frac{\lambda^{-1}}{\sqrt{1+\lambda^{-2}}}]) \\
    &= (t\frac{\lambda}{\sqrt{1+\lambda^{-2}}}+t\frac{\lambda^{-1}}{\sqrt{1+\lambda^{-2}}}) \cdot \frac{1}{\sqrt{1+\lambda^2}} = \lambda t.
\end{align*}
From this we see that $\blup_{\textnormal{pc},\source,\target}(\mathbb{R}^2,0) \cong \mathbb{R}^\times \ltimes \mathbb{R}$, as expected.
\end{exam}

\subsection{Examples of blow-ups of pairs of Lie groupoids and Lie algebroids}\label{sec: Examples of blow-ups of pairs of Lie groupoids and Lie algebroids}
We now turn to examples of blow-ups of Lie groupoids and Lie algebroids. The first type of examples comes from restricting to a particular class of submanifolds of the base manifold, called \textit{saturated submanifolds}.
\subsubsection{Blow-up along a saturated submanifold}\label{sec: blup of saturated submanifold}
We will start this section with discussing the blow-up of an action groupoid. Immediately after, we will discuss the blow-up of an action algebroid. After that, we will use these examples to show that a large class of blow-up groupoids and algebroids can be described, canonically, as action groupoids/action algebroids. As already mentioned, we will focus on the following specific class of submanifolds that we blow-up in the base manifold.
\begin{defn}\label{defn: saturated}
If $\groupoid$ is a groupoid, we call a submanifold $\subbase \subset \base$ for which
\[\source^{-1}(\subbase)=\target^{-1}(\subbase)\]
(i.e. $\subbase$ is a union of orbits) a \textit{saturated submanifold of $\base$ (with respect to $\group$)}. Similarly, if $\algebroid$ is an algebroid, we call a submanifold $\subbase \subset \base$ which is a union of orbits a \textit{saturated submanifold of $\base$ (with respect to $\algebr$)}.
\end{defn}
\begin{prop}\label{prop: blup of action groupoid (generalisation)}
Let $(\basetwo,\subbasetwo)$ be a pair of manifolds such that $\basetwo$ is equipped with a $\groupoid$ action with moment map $\mu: \basetwo \ra \base$, and assume that $\subbasetwo$ is saturated with respect to this action. Then
\[\blup_{\source,\target}(\group \ltimes \basetwo, \group \ltimes \subbasetwo) = \blup(\group \ltimes \basetwo, \group \ltimes \subbasetwo),\]
and we can describe a $\group$-action on $\blup(\basetwo,\subbasetwo)$, with moment map $\mu \circ \bldown: \blup(\basetwo,\subbasetwo) \ra \base$, such that
\[\blup(\group \ltimes \basetwo, \group \ltimes \subbasetwo) \cong \group \ltimes \blup(\basetwo,\subbasetwo)\]
as Lie groupoids.
\end{prop}
\begin{proof}
The first statement follows from the fact that $\subbasetwo \subset \basetwo$ is saturated, and that $\source_{\group \ltimes \basetwo} = \pr_2$, so that $d_\normal\source_{\group \ltimes \basetwo}$ is fiberwise a bijection, and therefore so is $d_\normal\target_{\group \ltimes \basetwo}$ (by the relation $d_\normal\source_{\group \ltimes \basetwo}=d_\normal\target_{\group \ltimes \basetwo} \circ d_\normal\iota_{\group \ltimes \basetwo}$). Now, the diffeomorphism
\begin{equation}\label{eq: blup(G times M,G times N) = G times blup(M,N)}
\blup(\group \times \basetwo, \group \times \subbasetwo) \xra{\sim} \group \times \blup(\basetwo,\subbasetwo) \textnormal{ given by } z \mapsto
\begin{cases}
    (g,[y,\xi])    & \text{if}\ z=[g,y,0,\xi] \\
    (g,[x,1])    & \text{if}\ z=[g,x,1]
\end{cases}
\end{equation}
from Proposition \ref{prop: blup(X x M, Y x M) = blup(X,Y) x M}, restricts to a diffeomorphism
\begin{equation}\label{eq: blup(G times M,G times N) = G times blup(M,N) 2}
    \blup(\group \tensor[_{\source}]{\times}{_{\mu}} \basetwo, \group \tensor[_{\source}]{\times}{_{\mu}} \subbasetwo) \xra{\sim} \group \tensor[_{\source}]{\times}{_{\mu \circ \bldown}} \blup(\basetwo,\subbasetwo).
\end{equation}
Indeed, since $\blup(\group \tensor[_{\source}]{\times}{_{\mu}} \basetwo, \group \tensor[_{\source}]{\times}{_{\mu}} \subbasetwo) \subset \blup(\group \times \basetwo, \group \times \subbasetwo)$ and $\group \tensor[_{\source}]{\times}{_{\mu \circ \bldown}} \blup(\basetwo,\subbasetwo) \subset \group \times \blup(\basetwo,\subbasetwo)$ are embedded submanifolds, we only have to show that the restriction of \eqref{eq: blup(G times M,G times N) = G times blup(M,N)} is a bijection, but this is readily verified. Moreover, if we denote $\nu: \group \tensor[_{\source}]{\times}{_{\mu}} \basetwo \ra \basetwo$ for the $\group$-action on $\base$, then
\[\blup_\nu(\group \tensor[_{\source}]{\times}{_{\mu}} \basetwo,\group \tensor[_{\source}]{\times}{_{\mu}} \subbasetwo) = \blup(\group \tensor[_{\source}]{\times}{_{\mu}} \basetwo,\group \tensor[_{\source}]{\times}{_{\mu}} \subbasetwo),\]
because $d_\normal\target_{\group \ltimes \basetwo} = \nu$ is fiberwise a bijection.
Therefore, we can define a multiplication on $\blup(\basetwo,\subbasetwo)$ by setting
\[\nu_{\blup(\basetwo,\subbasetwo)} \coloneqq \blup(\nu): \group \tensor[_{\source}]{\times}{_{\mu \circ \bldown}} \blup(\basetwo,\subbasetwo) \cong \blup(\group \tensor[_{\source}]{\times}{_{\mu}} \basetwo,\group \tensor[_{\source}]{\times}{_{\mu}} \subbasetwo) \ra \blup(\basetwo,\subbasetwo).\]
It remains to show that the diffeomorphism \eqref{eq: blup(G times M,G times N) = G times blup(M,N) 2} is a morphism of Lie groupoids. This is, however, precisely how we defined the action on $\blup(\basetwo,\subbasetwo)$: $\blup(\target) = \blup(\nu)$ becomes $\target_{\group \ltimes \blup(\basetwo,\subbasetwo)} = \nu_{\blup(\basetwo,\subbasetwo)}$ under the diffeomorphism \eqref{eq: blup(G times M,G times N) = G times blup(M,N) 2} by construction. That the map intertwines the other structure maps as well, is readily verified. This proves the statement.
\end{proof}
\begin{coro}\cite{10.36045/bbms/1292334057}\label{coro: action groupoid blow-up}
If a Lie group $G$ acts on $\base$, and $\subbase \subset \base$ is saturated with respect to the action, then 
\[\blup_{\source,\target}(G \times \base, G \times \subbase) = \blup(G \times \base, G \times \subbase),\]
and we can describe a $G$-action on $\blup(\base,\subbase)$ such that
\[\blup(G \ltimes \base, G \ltimes \subbase) \cong G \ltimes \blup(\base,\subbase)\]
as Lie groupoids.
\end{coro}
\begin{prop}\label{prop: blup of action algebroid (generalisation)}
Let $(\basetwo,\subbasetwo)$ be a pair of manifolds such that $\basetwo$ is equipped with an $\algebroid$ action with moment map $\mu: \basetwo \ra \base$, and assume that $\subbasetwo$ is saturated with respect to this action. Then 
\[\blup_\pi(\algebr \ltimes \basetwo, \algebr \ltimes \subbasetwo) = \blup(\algebr \ltimes \basetwo, \algebr \ltimes \subbasetwo),\]
and we can describe an $\algebr$-action on $\blup(\basetwo,\subbasetwo)$, with moment map $\mu \circ \bldown: \blup(\basetwo,\subbasetwo) \ra \base$, such that
\[\blup(\algebr \ltimes \basetwo, \algebr \ltimes \subbasetwo) \cong \algebr \ltimes \blup(\basetwo,\subbasetwo)\]
as Lie algebroids.
\end{prop}
\begin{proof}
A lot of arguments are completely analogous to the proof of Proposition \ref{prop: blup of action groupoid (generalisation)}. For example, the first statement follows from the fact that $\subbasetwo \subset \basetwo$ is saturated, and that $\pi_{\algebr \ltimes \basetwo}=\pr_2$, so that $d_\normal\pi_{\algebr \ltimes \basetwo}$ is fiberwise a bijection. Also, we obtain a diffeomorphism 
\[\varphi: \blup(\algebr \tensor[_{\pi}]{\times}{_{\mu}} \basetwo, \algebr \tensor[_{\pi}]{\times}{_{\mu}} \subbasetwo) \xra{\sim} \algebr \tensor[_{\pi}]{\times}{_{\mu \circ \bldown}} \blup(\basetwo,\subbasetwo)\]
in exactly the same way (by using Proposition \ref{prop: blup(X x M, Y x M) = blup(X,Y) x M}). Now, that $\algebr$ acts on $\basetwo$, with moment map $\mu: M \ra \base$, can be equivalently described by saying that $\mu^*\algebr=\algebr \tensor[_{\pi}]{\times}{_{\mu}} \basetwo$ carries a Lie algebroid structure such that the map $\pr_1: \mu^*\algebr \ra \algebr$ is a morphism of Lie algebroids. Since the blow-down map 
\[\bldown: \blup(\algebr \tensor[_{\pi}]{\times}{_{\mu}} \basetwo, \algebr \tensor[_{\pi}]{\times}{_{\mu}} \subbasetwo) \ra \algebr \tensor[_{\pi}]{\times}{_{\mu}} \basetwo\]
is a morphism of Lie algebroids, the map $\pr_1 \circ \bldown: \blup(\algebr \tensor[_{\pi}]{\times}{_{\mu}} \basetwo, \algebr \tensor[_{\pi}]{\times}{_{\mu}} \subbasetwo) \ra \algebr$ is too, and it is readily verified that this map is equal to $\pr_1 \circ \varphi$. This provides $\blup(M,N)$ with the desired $\algebr$-action, and $\varphi$ is an isomorphism of Lie algebroids by construction.
\end{proof}
\begin{coro}\label{coro: action algebroid blow-up}
If a Lie algebra $\mathfrak{g}$ acts on $\base$, and $\subbase \subset \base$ is saturated with respect to the action, then 
\[\blup_{\pi}(\mathfrak{g} \times \base, \mathfrak{g} \times \subbase) = \blup(\mathfrak{g} \times \base, \mathfrak{g} \times \subbase),\]
and we can describe a $\mathfrak{g}$-action on $\blup(\base,\subbase)$ such that
\[\blup(\mathfrak{g} \ltimes \base, \mathfrak{g} \ltimes \subbase) \cong \mathfrak{g} \ltimes \blup(\base,\subbase)\]
as Lie algebroids.
\end{coro}
\begin{rema}\label{rema: DNC instead of blup action groupoids and algebroids}
Notice that virtually the same proofs work if we replace ``$\blup$'' everywhere with ``$\DNC$'' in Propositions \ref{prop: blup of action groupoid (generalisation)} and \ref{prop: blup of action algebroid (generalisation)}.
\end{rema}
In the results above we were in the situation that $\subbase \subset \base$ is saturated. By observing that $\group$ is canonically isomorphic to the action groupoid $\group \times \base$ (with moment map $\source$ and multiplication map $\target$), we obtain that the blow-up of a groupoid/algebroid along the restriction groupoid/algebroid with respect to a saturated submanifold results in an action groupoid/algebroid. 
\begin{prop}\label{prop: blup along saturated}
Let $\groupoid$ be a Lie groupoid. If $\subbase \subset \base$ is saturated, then we can describe a (canonical) $\group$-action on $\blup(\base,\subbase)$, with moment map $\bldown$, such that
\[(\bldown,\blup(\source)): \blup(\group,\group|_{\subbase}) \xra{\sim} \group \ltimes \blup(\base,\subbase)\]
is an isomorphism of Lie groupoids. Similarly, if $\algebroid$ is a Lie algebroid, and $\subbase \subset \base$ is saturated, then we can describe a (canonical) $\algebr$-action on $\blup(\base,\subbase)$, with moment map $\bldown$, such that 
\[(\bldown,\blup(\pi)): \blup(\algebr,\algebr|_{\subbase}) \xra{\sim} \algebr \ltimes \blup(\base,\subbase)\]
is an isomorphism of Lie algebroids.
\end{prop}
\begin{rema}\label{rema: blup along saturated}
Since 
\begin{align*}
    (\group \tensor[_{\source}]{\times}{_{\bldown}} \blup(\base,\subbase))|_{\mathbb{P}(\base,\subbase)} &\cong \group|_{\subbase} \ltimes \mathbb{P}(\base,\subbase) \cong \mathbb{P}(\group,\group|_{\subbase}) \textnormal{ and } \\ 
    (\group \tensor[_{\source}]{\times}{_{\bldown}} \blup(\base,\subbase))|_{\base \setminus \subbase} &= \group|_{\base \setminus \subbase}, 
\end{align*}
%we see that $\group \tensor[_{\source}]{\times}{_{\bldown}} \blup(\base,\subbase)$ is diffeomorphic to $\blup(\group,\group|_{\subbase})$ via the map $(\bldown,\blup(\source))$, and 
we see that the action map $\nu: \group \tensor[_{\source}]{\times}{_{\bldown}} \blup(\base,\subbase) \ra \blup(\base,\subbase)$ is given by the map 
\[g \cdot z \coloneqq
\begin{cases}
    g \cdot [y,\xi]    & \text{if}\ z=[y,\xi] \\
    [x,1]    & \text{if}\ z=[x,1]
\end{cases} = 
\begin{cases}
    [\target(g),d_\normal\target(g)\eta]    & \text{if}\ z=[\source(g),d_\normal\source(g)\eta] \\
    [\target(g),1]    & \text{if}\ z=[\source(g),1]
\end{cases}\]
(note: if $z=[y,\xi] \in \mathbb{P}(\base,\subbase)$, then $(g,z) \in \group \tensor[_{\source}]{\times}{_{\bldown}} \blup(\base,\subbase)$ if and only if $g \in \group|_{\subbase}$%, and if $z=[x,1] \in \base \setminus \subbase \subset \blup(\base,\subbase)$, then $(g,z) \in \group \tensor[_{\source}]{\times}{_{\bldown}} \blup(\base,\subbase)$ if and only if $\source(g)=x$
\end{rema}%We start with the following observation.
We also include the following two statements which might be instructive at this point.
\begin{prop}\cite{proper}\label{prop: saturated submanifold iso normal bundle}
Let $\groupoid$ be a Lie groupoid, and let $\subbase \subset \base$ be a saturated embedded submanifold.
Then there is a $\group|_{\subbase}$-action on $\normal(\base,\subbase)$, with moment map $\pi_\subbase: \normal(\base,\subbase) \ra \subbase$, such that the map
\[(\pi_\subbase^\group,d_\normal\source): \normal(\group,\group|_{\subbase}) \xra{\sim} \group|_{\subbase} \ltimes \normal(\base,\subbase)\]
is an isomorphism of Lie groupoids. Moreover, $\mathbb{P}_{\source,\target}(\group,\group|_{\subbase})=\mathbb{P}(\group,\group|_{\subbase})$, and the above map descends to an isomorphism of Lie groupoids $\mathbb{P}(\group,\group|_{\subbase}) \xra{\sim} \group|_{\subbase} \ltimes \mathbb{P}(\base,\subbase)$.
\end{prop}
\begin{proof}
The $\group|_{\subbase}$-action on $\normal(\base,\subbase)$ is 
\begin{align*}
    \nu: \group|_{\subbase} \tensor[_{\source}]{\times}{_{\pi_\subbase}} \normal(\base,\subbase) &\ra \normal(\base,\subbase) \textnormal{ given by } (g,y,\xi) \mapsto (\target(g),d_\normal\target(g)\eta), \textnormal{ where } \\
    \eta &\in \normal_g(\group,\group|_{\subbase}) \textnormal{ such that } \xi = d_\normal\source(g)\eta.
\end{align*}
To see that this is well-defined, suppose that also $\xi = d_\normal\source(g)\zeta$. Since $\source$ and $\target$ are submersions, $\source^{-1}(\subbase)=\target^{-1}(\subbase)$ is an embedded submanifold with 
\[d\source(g)^{-1}(T_{\source(g)}\subbase) = T_g\source^{-1}(\subbase) = T_g\target^{-1}(\subbase) = d\target(g)^{-1}(T_{\target(g)}\subbase).\]
Since $d\source(g)(\eta-\zeta) \in T_{\source(g)}\subbase$, we see that $d\target(g)(\eta-\zeta) \in T_{\target(g)}\subbase$, so $d_\normal\target(g)\eta=d_\normal\target(g)\zeta$. 

That the given map is a morphism of groupoids is now readily verified. To see that it is even an isomorphism, notice that 
\begin{center}
\begin{tikzcd}
\normal(\group,\group|_{\subbase}) \ar[r] \ar[d,shift left] \ar[d,shift right] & \group|_{\subbase} \ar[d,shift left] \ar[d,shift right] \\
\normal(\base,\subbase) \ar[r] & \subbase
\end{tikzcd}
\end{center}
is a $\VB$-groupoid (see Section \ref{sec: VB groupoids}), and the given map is the map from Lemma \ref{lemm: submersion of pullback in VB-groupoid}, which is, in that lemma, proven to be a surjective submersion. Since the map is also injective in this case, it is a diffeomorphism, and therefore it is an isomorphism of Lie groupoids. The last statements follows, because the map restricts to an isomorphism
\[\normal_{\source,\target}(\group,\group|_{\subbase}) \xra{\sim} \group|_{\subbase} \ltimes (\normal(\base,\subbase) \setminus 0_\subbase)\]
as is readily verified, and it is obvious that the diffeomorphism $\mathbb{P}_{\source,\target}(\group,\group|_{\subbase}) \xra{\sim} \group|_{\subbase} \ltimes \mathbb{P}(\base,\subbase)$ it descends to is a morphism of Lie groupoids as well. Moreover, $\mathbb{P}_{\source,\target}(\group,\group|_{\subbase})=\mathbb{P}(\group,\group|_{\subbase})$% by Lemma \ref{lemm: blup_s,t=blup}
. This proves the statement.
\end{proof}
By using that, for a Lie algebroid $\algebroid$ and $\subbase \subset \base$ a saturated embedded submanifold,
\begin{center}
\begin{tikzcd}
\normal(\algebr,\algebr|_{\subbase}) \ar[r] \ar[d] & \algebr|_{\subbase} \ar[d] \\
\normal(\base,\subbase) \ar[r] & \subbase
\end{tikzcd}
\end{center}
is a $\VB$-algebroid, and using Lemma \ref{lemm: submersion D to A plus E} instead, we see that a similar results holds for a Lie algebroid $\algebroid$ and a saturated submanifold $\subbase \subset \base$ (note: an action of $\algebr|_{\subbase}$, with moment map $\pi_\subbase$, is a Lie algebroid structure on $\pi_\subbase^*(\algebr|_{\subbase})$ such that the map $\pr_1: \pi_\subbase^*(\algebr|_{\subbase}) \ra \algebr|_{\subbase}$ is a morphism of Lie algebroids).
\begin{prop}\label{prop: saturated submanifold iso normal bundle}
Let $\algebroid$ be a Lie algebroid, and let $\subbase \subset \base$ be a saturated embedded submanifold.
Then there is a $\algebr|_{\subbase}$-action on $\normal(\base,\subbase)$, with moment map $\pi_\subbase: \normal(\base,\subbase) \ra \subbase$, such that the map
\[(\pi_\subbase^\algebr,d_\normal\pi_\algebr): \normal(\algebr,\algebr|_{\subbase}) \xra{\sim} \algebr|_{\subbase} \ltimes \normal(\base,\subbase)\]
is an isomorphism of Lie algebroids. Moreover, $\mathbb{P}_{\pi}(\algebr,\algebr|_{\subbase})=\mathbb{P}(\algebr,\algebr|_{\subbase})$, and the above map descends to an isomorphism of Lie algebroids $\mathbb{P}(\algebr,\algebr|_{\subbase}) \xra{\sim} \algebr|_{\subbase} \ltimes \mathbb{P}(\base,\subbase)$.
\end{prop}

\subsubsection{Blow-up of a foliation}
As a next example, we will show that foliations behave well under the blow-up construction with respect to a leaf. We fix a pair of manifolds $(M,L)$ such that $M$ is equipped with a foliation $\mathcal{F}$ for which $L \subset M$ is a leaf, i.e. $L$ is a leaf of the foliation $\mathcal{F}$ that is a closed embedded submanifold (reminder: closed in the topological sense).
\begin{exam}[Blow-up of a foliated manifold]\label{exam: Blow-up along a leaf of a foliated manifold}
%As we have seen, this means that we have a Lie algebroid $T\mathcal{F} \ra M$ with injective anchor map $\anchor_{\mathcal{F}}: T\mathcal{F} \hookrightarrow TM$. 
We will show here that, for a closed embedded leaf $L \subset M$, $\blup(M,L)$ naturally inherits a foliation. %The easiest way to go about this is as follows: we can prove the statement by showing that $\blup(M,L)$ has a foliation atlas, so, by the local nature of the blow-up construction, we may assume $M=\mathbb{R}^n$ is trivially foliated, and $L=\mathbb{R}^p \times \{0\}$. Since
%\[\blup(M,L) = \blup(\mathbb{R}^n,\mathbb{R}^p) \cong \blup(\mathbb{R}^q,0) \times \mathbb{R}^p,\]
%we can even assume $L=0$. But then the statement is: if $M$ is foliated by points, then so is $\blup(M,L)$, which is of course true.
%
%Another approach is by using the Lie algebroid blow-up construction: 
The Lie algebroid $\pi: T\mathcal{F} \ra M$, whose anchor map $\iota: T\mathcal{F} \hookrightarrow TM$ is the inclusion, has the Lie subalgebroid $T\mathcal{F}|_{L} = TL$ (note: $\Gamma(T\mathcal{F},TL) = T\mathcal{F}$). It suffices now to show that the Lie algebroid
\[\blup(T\mathcal{F},TL) \ra \blup(M,L)\]
(recall Proposition \ref{prop: blup along saturated})
again has an injective anchor map. This anchor map is given by 
\[\blup(T\mathcal{F},TL) \xra{\blup(\iota)} \blup_{\pi_{TM}}(TM,TL) \xra{\widetilde\anchor} T\blup(M,L).\]
Since $T\mathcal{F} \cap TL = TL$, we see that, by Proposition \ref{prop: analogue constant rank when applying blup}, the map
\[\blup(\iota): \blup_{\pi}(T\mathcal{F},TL) \ra \blup_{\pi_{TM}}(TM,TL),\]
is an injective immersion. However, the map
\[\widetilde\anchor: \blup_{\pi_{TM}}(TM,TL) \ra T\blup(M,L)\]
has, over a point $[y,\xi] \in \mathbb{P}(M,L)$, a $1$-dimensional kernel (see Example \ref{exam: blup tangent bundle}). Therefore, we have to show that, fiberwise, the image of $\blup(\iota)|_{\mathbb{P}(M,L)}$ has trivial intersection with the kernel of $\widetilde\anchor|_{\mathbb{P}(M,L)}$. To show this, we can work in a vector bundle chart, so we can consider $\mathcal{F}$ to be a trivial foliation on $M=\mathbb{R}^n$, with $L=\mathbb{R}^p \times \{0\}$. The quotient map $q: \normal(M,L) \setminus 0_N \cong \mathbb{R}^{p} \times (\mathbb{R}^{q} \setminus \{0\}) \rightarrowdbl \mathbb{P}(M,L) \cong \mathbb{R}^p \times \mathbb{RP}^{q-1}$ can locally be expressed as the map
\[\varphi_i \circ q: \mathbb{R}^p \times (\mathbb{R}^q\setminus \{0\}) \ra \mathbb{R}^p \times \mathbb{R}^{q-1} \textnormal{ given by } (y,\xi) \mapsto (y,\frac{\xi_1}{\xi_i},\dots,\frac{\xi_{i-1}}{\xi_i},\frac{\xi_{i+1}}{\xi_i},\dots,\frac{\xi_q}{\xi_i}),\]
where $(U_i^{\mathbb{RP}},\varphi_i)$ is a standard chart of $\mathbb{RP}^{q-1}$. Therefore, we can write (locally), for $(y,\xi) \in \mathbb{R}^{n}$, 
\[dq(y,\xi)(\zeta,\eta)=(\zeta,\frac{\eta_1}{\xi_i} - \frac{\eta_i\xi_1}{\xi_i^2},\dots,\frac{\eta_{i-1}}{\xi_i} - \frac{\eta_i\xi_{i-1}}{\xi_i^2},\frac{\eta_{i+1}}{\xi_i} - \frac{\eta_{i+1}\xi_1}{\xi_i^2},\dots,\frac{\eta_q}{\xi_i} - \frac{\eta_i\xi_q}{\xi_i^2}).\]
Notice that this map has kernel given by
\[\{(0,\lambda \cdot \xi_1,\dots,%\lambda \cdot \xi_{i-1},\lambda \cdot \xi_i,\lambda \cdot \tfrac{\xi_{i+1}}{\xi_i},\dots,
\lambda \cdot \xi_q) \in \mathbb{R}^{p+q} \mid \lambda \in \mathbb{R}\}.\]
However, an element $(0,\lambda \cdot \xi_1,\dots,%\lambda \cdot \xi_{i-1},\lambda \cdot \xi_i,\lambda \cdot \tfrac{\xi_{i+1}}{\xi_i},\dots,
\lambda \cdot \xi_q) \in \mathbb{R}^{p+q}$ is tangent to $\mathbb{R}^p \times \{v\}$ for some $v \in \mathbb{R}^q$ if and only if $\lambda=0$. This shows that, indeed, fiberwise, the image of $\blup(\iota)|_{\mathbb{P}(M,L)}$ has trivial intersection with the kernel of $\widetilde\anchor|_{\mathbb{P}(M,L)}$ (note: $\widetilde\anchor([y,\xi]) = dq(y,\xi)$; see Proposition \ref{prop: tanent bundle blow-up}). 

Alternatively, we can use Proposition \ref{prop: blup along saturated}: $L \subset M$ is saturated, so 
\[\blup_{\pi}(T\mathcal{F},TL) = \blup(T\mathcal{F},TL),\]
and so that the statement can also be shown by using the commutativity of the following diagram:
\begin{center}
    \begin{tikzcd}
        \blup(T\mathcal{F},TL) \ar[r,"\anchor"] \ar[d,"\sim"] & T\blup(M,L) \ar[r,"d\bldown"] & TM \ar[d,"="] \\
        T\mathcal{F} \ltimes \blup(\base,\subbase)  \ar[r,"\pr_1"] & T\mathcal{F} \ar[r,"\iota",hook] & TM.
    \end{tikzcd}\vspace{-\baselineskip}
\end{center}
\end{exam}
\begin{term}
We denote the blow-up of the foliation $\mathcal{F}$ along $L$ by 
\[\blup(\mathcal{F},L).\] 
Similarly, the induced foliation on $\mathbb{P}(M,L)$ is written as $\mathbb{P}(\mathcal{F},L)$.
\end{term}
To describe the leaves of the blow-up foliation $\blup(\mathcal{F},L)$, recall Proposition \ref{prop: isotropy groups and orbits of blow-up algebroid}. That is, in $M \setminus L \subset \blup(M,L)$ the partition into leaves is equal to the partition into leaves of $M \setminus L \subset M$, because $\blup(T\mathcal{F},TL)|_{M \setminus L} \cong T\mathcal{F}|_{M \setminus L}$.
In $\mathbb{P}(M,L) \subset \blup(M,L)$, the leaves are described by the subbundle:
\[\mathbb{P}(T\mathcal{F},TL) \ra \mathbb{P}(M,L)\]
of the vector bundle $\blup(T\mathcal{F},TL)$.
To be more specific about the nature of this foliation, we will use that $\mathbb{P}(M,L) \ra L$ is a fiber bundle with typical fiber $\mathbb{RP}^{q-1}$. The foliation $\mathbb{P}(\mathcal{F},L)$ falls into the following class of foliations: 
\begin{defn}\cite{GeoTheoryFolio}\label{defn: transverse foliation of a fiber bundle}
Let $E \xra{\pi} \base$ be a fiber bundle with typical fiber $F$. A foliation $\mathcal{F}$ on $E$ is called \textit{transverse to the fibers} if $\codim \mathcal{F} = \dim F$ (i.e. $\dim \mathcal{F} = \dim \base$), and for all leaves $L$ of $\mathcal{F}$, if $e \in L$, then the fiber $E_x$ of $E$ over $x\coloneqq\pi(e)$ is transverse to $L$, and $L \xra{\pi|_L} \base$ is a covering projection.
\end{defn} 
Since the fibers of $\mathbb{P}(M,L) \ra L$ are compact, the following result ensures that the foliation $\mathbb{P}(\mathcal{F},L)$ is indeed transverse to the fibers.
\begin{prop}\cite{GeoTheoryFolio}\label{prop: compact fibers transverse foliation}
Let $E \xra{\pi} \base$ be a fiber bundle with typical fiber $F$, and let $\mathcal{F}$ be a foliation on $E$. If $F$ is compact, then $\mathcal{F}$ is transverse to the fibers if and only if for all leaves $L$ of $\mathcal{F}$, if $e \in L$, then the fiber $E_x$ of $E$ over $x\coloneqq\pi(e)$ is transverse to $L$ (i.e. the assumption in Definition \ref{defn: transverse foliation of a fiber bundle} that $L \xra{\pi|_L} \base$ is a covering projection is redundant in this case).
\end{prop}
\begin{proof}
Obviously, we only have to prove the ``if'' statement. The idea is as follows: each $x \in \base$ admits an open neighbourhood $U$ such that, for all $u \in U$, $E_u \cap L'$ consists of precisely one point for all leaves $L'$ of the foliation $\mathcal{F}|_{\pi^{-1}(U)}$, whose leaves are given by the connected components of the intersection of leaves in $\mathcal{F}$ with $\pi^{-1}(U)$. Then, for every leaf $L'$ of $\mathcal{F}|_{\pi^{-1}(U)}$, we will see that $\pi_{L'}: L' \ra U$ is a bijection, and therefore a diffeomorphism, which shows that $\pi|_L: L \ra \base$ is a covering projection (note: since $\dim L = \dim \base$, $\pi|_L: L \ra \base$ is a local diffeomorphism).

Since the fibers are compact, we can reduce the problem to showing that, for all $x \in \base$,
\begin{align}\label{eq: compact fibers claim}
&\textnormal{if } e \in E_x, \textnormal{ there is a foliation chart } (W,\varphi) \textnormal{ around } e \textnormal{ such that, for all } w \in \pi(W), \nonumber \\
&\textnormal{the intersection of } E_w \textnormal{ with a plaque of } W \textnormal{ consists of exactly one point}.
\end{align} 
To see why the statement follows from this claim, notice that, since $E_x$ is compact, we can cover $E_x$ by finitely many charts $(W_1,\varphi_1),\dots,(W_k,\varphi_k)$ with the mentioned property in \eqref{eq: compact fibers claim}. So, by choosing a connected open set
\[U \subset \bigcap_{i=1}^k \pi(W_k)\]
which contains $x$, if $L'$ is a leaf of $\mathcal{F}|_{\pi^{-1}(U)}$, it is contained in a plaque $L_1'$ of, say, $W_1$ (note: $\pi(L_1') = \pi(W_1) \supset U$ by \eqref{eq: compact fibers claim}). In particular, $\pi|_{L'}: L' \ra U$ is injective, and therefore also surjective, since 
\[\pi(L') = \pi(L_1' \cap \pi^{-1}(U)) = \pi(L_1') \cap \pi(\pi^{-1}(U)) = U.\] 
This shows that $\pi|_{L'}: L' \ra U$ is a bijective local diffeomorphism, i.e. it is a diffeomorphism.

It remains to prove \eqref{eq: compact fibers claim}. To do this, let $L$ be a leaf through $e \in E_x$. %Recall that we can view $L \subset E$ as an immersed submanifold. With the topology on $L$ that makes it into a smooth submanifold, it is readily verified that $e \in L$ is an isolated point in $E_x \cap L$ (by using that $L$ is transverse to $E_x$). Therefore, we 
We can find a connected open subset $V \ni e$ of $L$ such that, for all $v \in V$, $E_{\pi(v)} \cap V = \{v\}$. Then the map
\[\pi^{-1}(\pi(V)) \ra V \textnormal{ given by } z \mapsto v; \quad v \in E_{\pi(z)} \cap V\]
is well-defined, and it maps $V \subset \pi^{-1}(\pi(V))$ onto $V$. Moreover, by shrinking $V$, and arguing locally, the map is readily verified to be a submersion. Again, by shrinking $V$, we may assume that there is a foliation chart $(W',\varphi')$ such that $\overline{V}$ is contained in a plaque of $W'$, and so we may pick a foliation chart $(W,\varphi)$ such that $W \subset W' \cap \pi^{-1}(\pi(V))$ and $\varphi(V \cap W) = \varphi(W) \cap (\mathbb{R}^p \times \{0\})$ (where $p \coloneqq \dim L$). The chart $(W,\varphi)$ has the desired property from \eqref{eq: compact fibers claim}. This proves the statement.
\end{proof}
\begin{rema}\label{rema: P(F,L) is transverse to fibers}
To see that $\mathbb{P}(\mathcal{F},L)$ is transverse to the fibers, we have to prove, by the above result, that the leaves are transverse to the fibers. Since this is a local question, we may reduce to the case that $\mathcal{F}$ is the trivial foliation on $M=\mathbb{R}^n$, with $L=\mathbb{R}^p \times \{0\}$. But then $\mathbb{P}(M,L) \cong \mathbb{P}(\mathbb{R}^q,0) \times \mathbb{R}^p \cong \mathbb{RP}^{q-1} \times \mathbb{R}^p$, in which case it is obvious that the leaves are transverse to the fibers.
\end{rema}
The existence of a foliation on a fiber bundle, which is transverse to the fibers, is determined by whether or not the fiber bundle is isomorphic (as fiber bundles) to a fiber bundle with a discrete structure group. Explicitly, having a discrete structure group translates into the following definition:
\begin{defn}\cite{GeoTheoryFolio}\label{defn: discrete structure group}
Let $E \xra{\pi} \base$ be a fiber bundle with typical fiber $F$. If we can find a fiber bundle atlas $\{(U_i,\psi_i)\}$ such that the maps
\[U_{ij} \ra \textnormal{Diff } F \textnormal{ given by } x \mapsto \psi_{ij}(x,\cdot): F \xra{\sim} F\]
(here, $U_{ij} \coloneqq U_i \cap U_j$ and $\psi_{ij} \coloneqq \psi_i \circ \psi_j^{-1}: U_{ij} \times F \xra{\sim} U_{ij} \times F$) are locally constant, then $E$ is said to have a \textit{discrete structure group}. 
\end{defn}
\begin{rema}\cite{GeoTheoryFolio}\label{rema: transverse foliation of a fiber bundle}
If $E \xra{\pi} \base$ is a fiber bundle with typical fiber $F$, we see that $E$ admits a foliation whose plaques are given by the fibers of the map
\[\pi^{-1}(U_i) \xra[\psi_i]{\sim} U_i \times F \xra{\pr_2} F\]
(where $\{(U_i,\psi_i)\}$ is an atlas) if and only if $E$ has a discrete structure group (such that the condition from Definition \ref{defn: discrete structure group} is satisfied with the atlas $\{(U_i,\psi_i)\}$). Indeed, if we pick charts $(U_i,\varphi_i)$ of $\base$, and charts $(W_r,\chi_r)$ of $F$, then the maps
\[f \coloneqq \pr_2 \circ (\varphi_\ell \times \chi_s) \circ \psi_{ij} \circ (\varphi_k \times \chi_r)^{-1}: \varphi_k(U_{ijk\ell}) \times \chi_r(W_{rs}) \ra \chi_s(W_{rs})\]
satisfy $f(x,y)=f(y)$ if and only if the maps $x \mapsto \psi_{ij}(x,\cdot)$ are locally constant. It is readily verified that this foliation is transverse to the fibers. %All foliations transverse to the fibers arise this way; that is, if the foliation on $E$ is transverse to the fibers, then we can find a diffeomorphism $\varphi: E \ra E'$, with $E'$ a fiber bundle with discrete structure group, such that $\pi_{E'} \circ \varphi = \pi_E$ (i.e. $\varphi$ is an isomorphism of fiber bundles), but we will not discuss this (see e.g. \cite{GeoTheoryFolio}).
\end{rema}
We can also describe the foliation on $\mathbb{P}(M,L)$ by factoring through a foliation on $\normal(M,L)$. That the foliation $\mathcal{F}$ on $M$ induces a foliation on $\normal(M,L)$, we again use that the inclusion $\iota: T\mathcal{F} \hookrightarrow TM$ induces the injective anchor map
\[\normal(T\mathcal{F},TL) \xhookrightarrow{d_\normal\iota} \normal(TM,TL) \cong T\normal(M,L)\]
of the Lie algebroid $\normal(T\mathcal{F},TL) \ra \normal(M,L)$. This foliation is denoted by $\normal(\mathcal{F},L)$. %To see that this foliation is transverse to the fibers, we use that, as a fiber bundle, $\normal(M,L) \ra L$ has a discrete structure group. 
\begin{exam}\cite{GeoTheoryFolio}\label{exam: foliation of normal bundle}
We can apply the above idea to show that $\normal(\mathcal{F},L)$ is a foliation transverse to the leaves (and similarly for $\mathbb{P}(M,L) \ra L$). Indeed, if $\{(U_i,\varphi_i)\}=\{(U_i,y_i^1,\dots,y_i^p,x_i^1,\dots,x_i^q)\}$ is an adapted chart (to $L$), recall that then we have the vector bundle coordinates
\[\normal(U_i,V_i) \ra V_i \times \mathbb{R}^q \textnormal{ given by } (y,\xi) \mapsto (y,dx_i^1(y)\xi,\dots,dx_i^q(y)\xi)\]
(with $V_i \coloneqq U_i \cap L$). %So, the transition maps are the maps
%\begin{align*}
%   d_\normal \varphi_{ij}: \varphi_j(V_{ij}) \times \mathbb{R}^q &\ra \varphi_i(V_{ij}) \times \mathbb{R}^q \textnormal{ given by } \\
%    (y,\xi) &\mapsto (y_i^1 \circ \varphi_j^{-1}(y),\dots,y_i^p \circ \varphi_j^{-1}(y),d(x_i^1 \circ \varphi_j^{-1})(y)\xi,\dots,d(x_i^q \circ \varphi_j^{-1})(y)\xi).
%\end{align*}
So, if $\{(U_i,\varphi_i)\}$ is a foliation atlas (with $0 \in \varphi_i(U_i)$), %i.e. the maps
%\[y \mapsto (x_i^1 \circ \varphi_j^{-1}(y,\cdot),\dots,x_i^q \circ \varphi_j^{-1}(y,\cdot))\]
%are locally constant, then so are the maps
then showing that $\normal(M,L) \ra L$ has a transverse foliation structure comes down to the fact that the maps
\[y \mapsto (d(x_i^1 \circ \varphi_j^{-1})(y)\cdot,\dots,d(x_i^q \circ \varphi_j^{-1})(y)\cdot) = (d(x_i^1 \circ \varphi_j^{-1})(0)\cdot,\dots,d(x_i^q \circ \varphi_j^{-1})(0)\cdot)\]
are locally constant (see Remark \ref{rema: transverse foliation of a fiber bundle}).
This shows that $\normal(M,L) \ra L$ naturally inherits a foliation which is transverse to the fibers. We denote this foliation by $\normal(\mathcal{F},L)$. If $\iota: T\mathcal{F} \hookrightarrow TM$ is the inclusion, then the corresponding integrable subbundle $T\normal(\mathcal{F},L)$ of $T\normal(M,L)$ is $\normal(T\mathcal{F},TL)$.
We can similarly prove that $\mathbb{P}(M,L) \ra L$ inherits a transverse foliation structure this way, and the corresponding integrable subbundle $T\mathbb{P}(\mathcal{F},L)$ of $T\mathbb{P}(M,L)$ is $\mathbb{P}(T\mathcal{F},TL)$.
\end{exam} 
\begin{rema}\label{rema: quotient foliation}
The $\mathbb{R}^\times$-action on $\normal(M,L) \setminus 0_L$ leaves the foliation $\normal(\mathcal{F},L)|_{\normal(M,L) \setminus 0_L}$ invariant; that is, the diffeomorphisms 
\[\lambda: \normal(M,L) \setminus 0_L \xra{\sim} \normal(M,L) \setminus 0_L \textnormal{ given by } (y,\xi) \mapsto (y,\lambda\xi)\]
map leaves to leaves (we can reduce to the local case, where this property of the foliation is obvious). %Indeed, since a diffeomorphism $\lambda$ can be written as a composition $\lambda_k \circ \dots \circ \lambda_1$ for $\lambda_i \in (1-\epsilon,1+\epsilon)$ (for any $0 < \epsilon < 1$), we can reduce to the local case, where this property of the foliation is obvious. 
One can show now that the foliation descends to a foliation on the quotient $(\normal(M,L) \setminus 0_L)/\mathbb{R}^\times = \mathbb{P}(M,L)$, and it is precisely the foliation $\mathbb{P}(\mathcal{F},L)$. 
\end{rema}

Notice that the Lie groupoids
\begin{align*}
    \blup(\Mon(M,\mathcal{F}),\Mon(M,\mathcal{F})|_L) &\rra \blup(M,L) \textnormal{ and } \\
    \blup(\Hol(M,\mathcal{F}),\Hol(M,\mathcal{F})|_L) &\rra \blup(M,L)
\end{align*}
integrate the Lie algebroid $\blup(T\mathcal{F},TL) \ra \blup(M,L)$. We will show that the former groupoid is canonically isomorphic to the monodromy groupoid of the foliation $\blup(\mathcal{F},L)$.
We will use the following general fact about Lie groupoid-integrations of foliations.
\begin{lemm}\label{lemm: groupoid integration of foliation maps Hol and Mon}
Let $\group \rra M$ be an $\source$-connected groupoid, and assume it is an integration of the Lie algebroid $T\mathcal{F} \hookrightarrow TM$ of the foliated manifold $(M,\mathcal{F})$. Then there is a sequence of morphism of Lie groupoids, which are surjective local diffeomorphisms and integrate the identity map $\id_{T\mathcal{F}}$,
\[\Mon(M,\mathcal{F}) \xrightarrowdbl{F_\group^\Mon} \group \xrightarrowdbl{F_\group^\Hol} \Hol(M,\mathcal{F}),\]
such that the composition $F_\group^\Hol \circ F_\group^\Mon$ is the canonical map $F$ from Remark \ref{rema: groupoid morphism Mon -> Hol}. Moreover, the map $F_\group^\Mon$ is an isomorphism of groupoids if and only if $\group$ is, in addition, $\source$-simply connected.
\end{lemm}
\begin{proof}
Since $\Mon(M,\mathcal{F})$ is $\source$-connected and $\source$-simply connected (since a source fiber is the universal covering space of a leaf; see Proposition \ref{prop: (universal) cover of Mon and Hol}), the map $F_\group^\Mon$ is obtained by using Proposition \ref{prop: Lie groupoid morphism integration source-simply connected}, and the same proposition provides an inverse precisely when $\group$ is $\source$-simply connected. It remains to prove that there is a morphism of Lie groupoids $F_\group^\Hol$ fitting into the diagram 
\begin{center}
\begin{tikzcd}
\Mon(M,\mathcal{F}) \ar[rd,"F", twoheadrightarrow]\ar[rr,"F_\group^\Mon",twoheadrightarrow] & & \group \ar[ld,"F_\group^\Hol",twoheadrightarrow] \\
& \Hol(M,\mathcal{F}). &
\end{tikzcd}
\end{center}
But, by the surjectivity of these maps, there is only one candidate for the map $\group \rightarrowdbl \Hol(M,\mathcal{F})$, namely $F_\group^\Mon(g) \mapsto F(g)$. %To prove that this is well-defined, we describe the map $F_\group^\Mon$ more explicitly. 
This map is well-defined and surjective, but we will not prove it here; see e.g. \cite{IekeMariuscyclic}.
%Let $\gamma: [0,1] \ra M$ be a path contained in a leaf $L$, and let $\{U_i\}$ be a finite cover of $\gamma([0,1])$ such that $s_i: U_i \ra T_i$ are trivialising submersions onto transversals $T_i$ (see III in \ref{defn: foliation}; assume $\gamma([t_i,t_{i+1}]) \subset U_i$ for $0=t_0 < \dots < t_k = 0$.  
\end{proof}
\begin{prop}\label{prop: Monodromy and Holonomy blow-up}
We have a canonical Lie groupoid isomorphism
\begin{align*}
    \Mon(\blup(M,L),\blup(\mathcal{F},L)) &\xra{\sim} \blup(\Mon(M,\mathcal{F}),\Mon(M,\mathcal{F})|_L)  
    %\textnormal{ and } \\
    %\Hol(\blup(M,L),\blup(\mathcal{F},L)) &\xra{\sim} \blup(\Hol(M,\mathcal{F}),\Hol(M,\mathcal{F})|_L)
\end{align*}
%that restrict to isomorphisms 
which restricts to an isomorphism
\begin{align*}
    \Mon(\mathbb{P}(M,L),\mathbb{P}(\mathcal{F},L)) &\xra{\sim} \mathbb{P}(\Mon(M,\mathcal{F}),\Mon(M,\mathcal{F})|_L)  
    %\textnormal{ and } \\
    %\Hol(\mathbb{P}(M,L),\mathbb{P}(\mathcal{F},L)) &\xra{\sim} \mathbb{P}(\Hol(M,\mathcal{F}),\Hol(M,\mathcal{F})|_L)
\end{align*}
(recall Proposition \ref{prop: blup along saturated}).
\end{prop}
\begin{proof}
%We first prove the statement for the monodromy groupoid. 
Since the groupoid 
\[\blup(\Mon(M,\mathcal{F}),\Mon(M,\mathcal{F})|_L) \rra \blup(M,L)\] 
is also an integration of the Lie algebroid $T\blup(\mathcal{F},L) \ra \blup(M,L)$, we only have to prove, by Lemma \ref{lemm: groupoid integration of foliation maps Hol and Mon}, that $\blup(\Mon(M,\mathcal{F}),\Mon(M,\mathcal{F})|_L)$ is $\source$-connected and $\source$-simply connected. This follows from Proposition \ref{prop: blup along saturated}, because it says that 
\begin{equation}\label{eq: blup Mon saturated iso}
    \blup(\Mon(M,\mathcal{F}),\Mon(M,\mathcal{F})|_L) \cong \Mon(M,\mathcal{F}) \ltimes \blup(\base,\subbase)
\end{equation}
which is $\source$-connected and $\source$-simply connected because $\Mon(M,\mathcal{F})$ is.
%First of all, the subgroupoid 
%\[\Mon(M,\mathcal{F})|_{M \setminus L} \subset \blup_{\source,\target}(\Mon(M,\mathcal{F}),\Mon(M,\mathcal{F})|_L)\]
%is the monodromy groupoid of the foliation $\mathcal{F}|_{M \setminus L}$ on $M \setminus L$, which is $\source$-connected and $\source$-simply connected. By the groupoid structure of $\blup_{\source,\target}(\Mon(M,\mathcal{F}),\Mon(M,\mathcal{F})|_L)$, we only have to show now that the subgroupoid
%\[\mathbb{P}_{\source,\target}(\Mon(M,\mathcal{F}),\Mon(M,\mathcal{F})|_L) \rra \mathbb{P}(M,L)\]
%is $\source$-connected and $\source$-simply connected. This follows from Lemma \ref{lemm: saturated submanifold iso normal bundle}: we obtain an isomorphism of Lie groupoids
%\begin{equation}\label{eq: proj Mon saturated iso}
%    \mathbb{P}_{\source,\target}(\Mon(M,\mathcal{F}),\Mon(M,\mathcal{F})|_{L}) \xra{\sim} \Mon(M,\mathcal{F})|_{L} \ltimes \mathbb{P}(M,L)
%\end{equation}
%and the right-hand side is $\source$-connected and $\source$-simply connected because $\Mon(M,\mathcal{F})|_{L}$ is. 
%This proves that there are isomorphisms \eqref{eq: Mon blup} and \eqref{eq: Mon proj}.
This proves the statement.
\end{proof}
%Lastly, let us present a simple example of a blow-up of a singular foliation. 
%\begin{exam}\label{exam: singular folation blow-up}
%Recall the singular foliation of $\mathbb{R}^2$ presented in Example \ref{exam: foliation of R^2}. We will blow-up the origin in $\mathbb{R}^2$ and show that the foliation naturally lifts to a regular foliation on $\blup(\mathbb{R}^2,0)$. To do this, recall that the vector field which describes the singular foliation of $\mathbb{R}^2$ is given by $\sigma(x,y) \coloneqq -y\tfrac{\partial}{\partial x} + x\tfrac{\partial}{\partial y}$. Since this vector field is linear, its derivative at every point is equal to the map $(\xi_1,\xi_2) \mapsto (-\xi_2,\xi_1)$. From this we see that 
%\[\blup(\sigma): \blup(\mathbb{R}^2,0) \ra \blup(T\mathbb{R}^2,0) \ra T\blup(\mathbb{R}^2,0)\]
%is in fact the map 
%\end{exam}
\begin{rema}\label{rema: diagram Hol Mon blup}
Notice that the diagram
\begin{center}
\begin{tikzcd}
\Mon(\blup(M,L),\blup(\mathcal{F},L)) \ar[d,"F_{\blup(\mathcal{F},L)}",twoheadrightarrow] \ar[r,"\sim"] & \blup(\Mon(M,\mathcal{F}),\Mon(M,\mathcal{F})|_L) \ar[d,"\blup(F_{\mathcal{F}})",twoheadrightarrow] \\
\Hol(\blup(M,L),\blup(\mathcal{F},L)) \ar[r] & \blup(\Hol(M,\mathcal{F}),\Hol(M,\mathcal{F})|_L)
\end{tikzcd}
\end{center}
commutes.
\end{rema}
Similarly, we have the following statement.
\begin{prop}\label{prop: Monodromy and Holonomy normal}
We can describe a foliation $\DNC(\mathcal{F},L)$ on $\DNC(M,L)$, which induces, for all $t \in \mathbb{R}^\times$, the foliation $\mathcal{F}$ on $\DNC(M,L)_t \cong M$, and the foliation $\normal(\mathcal{F},L)$ on $\DNC(M,L)_0 \cong \normal(M,L)$. Moreover, we have a canonical Lie groupoid isomorphism
\begin{align*}
    \Mon(\DNC(M,L),\DNC(\mathcal{F},L)) &\xra{\sim} \DNC(\Mon(M,\mathcal{F}),\Mon(M,\mathcal{F})|_L)  
    %\textnormal{ and } \\
    %\Hol(\DNC(M,L),\DNC(\mathcal{F},L)) &\xra{\sim} \DNC(\Hol(M,\mathcal{F}),\Hol(M,\mathcal{F})|_L)
\end{align*}
that restricts to a Lie groupoid isomorphism
\begin{align*}
    \Mon(\normal(M,L),\normal(\mathcal{F},L)) &\xra{\sim} \normal(\Mon(M,\mathcal{F}),\Mon(M,\mathcal{F})|_L).  
    %\textnormal{ and } \\
    %\Hol(\normal(M,L),\normal(\mathcal{F},L)) &\xra{\sim} \normal(\Hol(M,\mathcal{F}),\Hol(M,\mathcal{F})|_L).
\end{align*}\vspace{-\baselineskip}
%Moreover, we also have an isomorphism $\normal(\Mon(M,\mathcal{F}),\Mon(M,\mathcal{F})|_L) \cong \Mon(M,\mathcal{F})|_L \ltimes \normal(M,L)$. %and $\normal(\Hol(M,\mathcal{F}),\Hol(M,\mathcal{F})|_L) \cong \Hol(M,\mathcal{F})|_L \ltimes \normal(M,L)$.
\end{prop}
\begin{proof}
If $\iota: T\mathcal{F} \hookrightarrow TM$ is the inclusion, then the foliation $\DNC(\mathcal{F},L)$ is described by the Lie algebroid
\[\DNC(T\mathcal{F},TL) \xhookrightarrow{\DNC(\iota)} \DNC(TM,TL) \hookrightarrow T\DNC(M,L).\]
%(for the last map, see Proposition \ref{prop: DNC of tangent spaces}). 
That the foliation restricts to the foliation $\mathcal{F}$ on $\DNC(M,L)_t \cong M$, and the foliation $\normal(\mathcal{F},L)$ on $\DNC(M,L)_0 \cong \normal(M,L)$, follows simply from the fact that $\DNC(T\mathcal{F},TL)$ has $\DNC(T\mathcal{F},TL)_t \cong T\mathcal{F}$ and $\DNC(T\mathcal{F},TL)_0 \cong \normal(T\mathcal{F},TL)$ as subalgebroids over the bases $\DNC(M,L)_t \cong M$ and $\DNC(M,L)_0 \cong \normal(M,L)$, respectively. 

The last statement is proven similarly as in Proposition \ref{prop: Monodromy and Holonomy blow-up}% by using (this time the first part of) Lemma \ref{lemm: saturated submanifold iso normal bundle}
. 
\end{proof}

%Let $D \subset TM$ be an integrable generalised distribution, and let $L \subset M$ be a leaf that is an embedded submanifold. 

\subsubsection{Other examples of blow-ups}
Here we will describe a few more examples of blow-ups of Lie groupoids and Lie algebroids. 
\begin{exam}[Principal bundles]\label{exam: blow-up of principal bundles}
Suppose $P \xra{\pi} \basetwo$ is a principal $G$-bundle, and let $Q \subset P$ be invariant with respect to the action of $G$; notice that then $Q \ra \pi(Q) \eqqcolon \subbasetwo$ is again a principal $G$-bundle. The action of $G$ on $\blup(P,Q)$ (see Corollary \ref{coro: action groupoid blow-up}) is again free and proper (e.g. because the blow-down map $\blup(P,Q) \ra P$ is proper). This yields a principal bundle
\[\blup(P,Q) \ra \blup(P,Q)/G.\]
Let us also discuss the Gauge groupoid $\blup(P,Q) \otimes_G \blup(P,Q)$ and the Atiyah algebroid $A(\blup(P,Q))$. % ask in what way $\blup(P,Q) \otimes_G \blup(P,Q)$ is related to $\blup_{\source,\target}(P \otimes_G P, Q \otimes_G Q)$.% First of all, since the projection map $(P \times P, Q \times Q) \ra ((P \times P)/G,(Q \times Q)/G)$ has a fiberwise injective normal derivative (since $Q \times Q$ equals the inverse image of $(Q \times Q)/G$), we obtain a map
We have an isomorphism
\[\blup_{\source,\target}((P \times P)/G, (Q \times Q)/G) \xra{\sim} \blup_{\source,\target}(P \times P, Q \times Q)/G,\]
%and it is readily verified to be $G$-invariant (with trivial action on $\blup((P \times P)/G, (Q \times Q)/G)$). 
and by Example \ref{exam: blup pair groupoids}, we obtain an isomorphism
\[\blup_{\source,\target}((P \times P)/G, (Q \times Q)/G) \xra{\sim} (\hat P \times_{\mathbb{R}} \hat P/\mathbb{R}^\times)/G,\]
where $\hat P = \DNC(P,Q) \setminus Q \times \mathbb{R}$. Similarly, we have an isomorphism
\[\blup_{\pi}(TP/G,TQ/G) \xra{\sim} (\ker d\hat{t}|_{\hat{P}}/\mathbb{R}^\times)/G\]
(see Example \ref{exam: blup tangent bundle}).
The map $(\hat{P} \times_{\mathbb{R}} \hat{P})/\mathbb{R}^\times \ra \blup(P,Q) \times \blup(P,Q)$ now induces a map 
\[\blup_{\source,\target}((P \times P)/G,(Q \times Q)/G) \ra (\blup(P,Q) \times \blup(P,Q))/G,\]
and the map $\ker d\hat{t}|_{\hat{P}}/\mathbb{R}^\times \ra T\blup(P.Q)$ induces a map $\blup_{\pi}(TP/G,TQ/G) \ra T\blup(P,Q)/G$. Lastly, notice that $\blup_{\source,\target}((P \times P)/G, (Q \times Q)/G)$ and $\blup_{\pi}(TP/G,TQ/G)$ are not transitive.
\end{exam}
%The next example (both the groupoid and algebroid cases) are very general, and can therefore be applied to many different situations. 
%Before we go into the next class of examples, observe that, 
A remark in between: if $\groupoid$ is a Lie groupoid such that $\identity(\base) \subset \group$ is closed (i.e. $\group$ is Hausdorff; see Proposition \ref{prop: G Hausdorff iff X closed}) with Lie algebroid $\algebroid$, then 
\[\mathbb{P}(\group,\identity(\base)) = (\normal(\group,\base) \setminus 0_\base)/\mathbb{R}^\times \cong (\algebr \setminus 0_\base)/\mathbb{R}^\times = \mathbb{P}(\algebr),\]
so we can think of $\blup(\group,\identity(\base))$ as ``replacing $\identity(\base) \subset \group$ with the projectivisation of the Lie algebroid of $\group$'' (that is, as a manifold; note: $\blup_{\source,\target}(\group,\identity(\base))= \emptyset$).
%\begin{exam}\label{exam: blup along unit subgroupoid}

%\end{exam}
We will examine now what we can say about the blow-up of pullback Lie groupoids/algebroids. 
\begin{prop}\label{prop: blup of pullback}
Let $(\groupoid,\subgroupoid)$ be a pair of Lie groupoids and let $f: (\basetwo,\subbasetwo) \ra (\base,\subbase)$ be a smooth map of pairs such that $f^{-1}(\subbase)=\subbasetwo$ and $d_\normal f$ is fiberwise injective. Moreover, assume that $f \times f$ (seen as a map of pairs) has a clean intersection with the map $(\target,\source)$ (seen as a map of pairs $(\group,\subgroup) \ra (\base \times \base, \subbase \times \subbase)$), and that $\blup(f) \times \blup(f)$ has a clean intersection with $(\blup(\target),\blup(\source))$. Then 
\[\blup_{\source^!,\target^!}(f^!\group,f^!\subgroup) \xra{\sim} \blup(f)^!\blup_{\source,\target}(\group,\subgroup)\]
as Lie groupoids (we denoted $\source^!$ and $\target^!$ for the source and target map of $f^!\group$, respectively).
\end{prop}
\begin{proof}
By assumption, the pull-back groupoids $f^!\group$, $f^!\subgroup$ and $\blup(f)^!\blup_{\source,\target}(\group,\subgroup)$ are Lie groupoids. Notice that, by viewing $f^!\group$ and $f^!\subgroup$ as subgroupoids of $\basetwo \times \group \times \basetwo$, we can view $(f^!\group,f^!\subgroup)$ as a pair of Lie groupoids (notice that $f^!\subgroup \subset f^!\group$ is indeed closed). From Proposition \ref{prop: blup respects fiber products}, we see that
\begin{align*}
    \blup_{\source^!,\target^!}(f^!\group,f^!\subgroup) &\xra{\sim} \blup(\basetwo,\subbasetwo) \tensor[_{\blup(f)}]{\times}{_{\blup(\source)}} \blup_{\source,\target}(\group,\subgroup) \tensor[_{\blup(\target)}]{\times}{_{\blup(f)}} \blup(\basetwo,\subbasetwo) \\
    &= \blup(f)^!\blup_{\source,\target}(\group,\subgroup),
\end{align*}
(observe that the source of this map is correct $\blup_{\source^!,\target^!}(f^!\group,f^!\subgroup)=\blup_{f \times (\target,\source) \times f}(f^!\group,f^!\subgroup)$),
and it is readily verified to be a morphism of Lie groupoids. This proves the statement.
\end{proof}
%\begin{rema}\label{rema: blup of pullback applied to bldown}
%In particular, we can apply the above proposition to $\bldown: (\blup(\base,\subbase),\mathbb{P}(\base,\subbase)) \ra (\base,\subbase)$. That is, if $\bldown \times \bldown$ has a clean intersection with $(\target,\source)$. Since $\bldown^{-1}(\subbase)=\mathbb{P}(\base,\subbase)$, $d_\normal\bldown$ is fiberwise injective, and $\blup(\bldown) \times \blup(\bldown)=\id_{\blup(\base,\subbase)} \times \id_{\blup(\base,\subbase)}$, it follows then from the above proposition that
%\[\blup_{\source^!,\target^!}(\bldown^!\group,\bldown^!\subgroup) \xra{\sim} \blup(\bldown)^!\blup_{\source,\target}(\group,\subgroup)=\blup_{\source,\target}(\group,\subgroup),\]
%where we again used that $\blup(\bldown)=\id_{\blup(\base,\subbase)}$ in the last equality. Notice that 
%\[(\bldown^!\group \rra \blup(\base,\subbase),\bldown^!\subgroup \rra \mathbb{P}(\base,\subbase))\]
%is a pair of groupoids $(\grouptwoid,\subgrouptwoid)$ such that $\subbasetwo \subset \basetwo$ is of codimension $1$. This therefore shows that, in case $\bldown \times \bldown$ has a clean intersection with $(\target,\source)$, the blow-up construction for Lie groupoids in the codimension $1$ case suffices to describe the blow-up construction for this pair of Lie groupoids.
%\end{rema}
Following a similar reasoning, we can prove a statement for pullbacks of Lie algebroids as well.
\begin{prop}\label{prop: blup of pullback}
Let $(\algebroid,\subalgebroid)$ be a pair of Lie algebroids and let $f: (\basetwo,\subbasetwo) \ra (\base,\subbase)$ be a smooth map of pairs such that $f^{-1}(\subbase)=\subbasetwo$ and $d_\normal f$ is fiberwise injective. Moreover, assume that $\anchor_{\algebr}$ (seen as a map of pairs) has a clean intersection with the map $df$ (seen as a map of pairs), and that $\anchor_{\blup_{\pi}(\algebr,\subalgebr)}$ has a clean intersection with $d\blup(f)$. Then 
\[\blup_{\pi^!}(f^!\algebr,f^!\subalgebr) \xra{\sim} \blup(f)^!\blup_{\pi}(\algebr,\subalgebr)\]
as Lie algebroids (we denoted $\pi^!$ for the projection map of $\pi^!\algebr$).
\end{prop}
\begin{proof}
By assumption, the pull-back algebroids $f^!\algebr$, $f^!\subalgebr$ and $\blup(f)^!\blup_{\pi}(\algebr,\subalgebr)$ are Lie algebroids. We can view $(f^!\algebr,f^!\subalgebr)$ as a pair of Lie algebroids (notice that $f^!\subalgebr \subset f^!\algebr$ is indeed closed). From Proposition \ref{prop: blup respects fiber products}, we see that
\begin{align*}
    \blup_{\pi^!}(f^!\algebr,f^!\subalgebr) &\xra{\sim} \blup_{\pi_{T\base}}(T\basetwo,T\subbasetwo) \tensor[_{\blup(df)}]{\times}{_{\blup(\anchor)}} \blup_{\pi}(\algebr,\subalgebr) \\
    &\xra{\sim} T\blup(\basetwo,\subbasetwo) \tensor[_{d\blup(f)}]{\times}{_{\anchor_{\blup_{\pi}(\algebr,\subalgebr)}}} \blup_{\pi}(\algebr,\subalgebr) \\
    &= \blup(f)^!\blup_{\pi}(\algebr,\subalgebr),
\end{align*}
where we used in the second isomorphism that the diagram
\begin{center}
\begin{tikzcd}
    \blup_{\pi_{T\basetwo}}(T\basetwo,T\subbasetwo) \ar[r,"\blup(df)"] \ar[d,"\widetilde{\anchor}"] & \blup_{\pi_{T\base}}(T\base,T\subbase) \ar[d,"\widetilde{\anchor}"] & \blup_{\pi}(\algebr,\subalgebr) \ar[l,"\blup(\anchor)"] \ar[d, "="]\\
    T\blup(\basetwo,\subbasetwo) \ar[r,"d\blup(f)"] & T\blup(\base,\subbase) & \blup_{\pi}(\algebr,\subalgebr) \ar[l,"\anchor_{\blup_{\pi}(\algebr,\subalgebr)}"]
\end{tikzcd}
\end{center}
commutes. 
%(observe that the source of this map is correct $\blup_{\source^!,\target^!}(f^!\group,f^!\subgroup)=\blup_{f \times (\target,\source) \times f}(f^!\group,f^!\subgroup)$) 
The above map is readily verified to be an isomorphism of Lie algebroids. This proves the statement.
\end{proof}

\subsection{Blow-up of differentiable stacks}
%Fix two pairs of groupoids $(\groupoid,\subgroupoid)$ and $(\grouptwoid,\subgrouptwoid)$. %The notion of Morita equivalence (see Section \ref{sec: Morita equivalence}) in this setting becomes the following.
%\begin{defn}
%Let $(E,F)$ be a pair of manifolds. If 
%\[\base \xla{\pi_\base} E \xra{\pi_\basetwo} \basetwo\] 
%is a Morita equivalence such that 
%\[\subbase \xla{\pi_\base|_{F}} F \xra{\pi_\basetwo|_{F}} \subbasetwo\] 
%is also a Morita equivalence, then $(\groupoid,\subgroupoid)$ and $(\grouptwoid,\subgrouptwoid)$ are called \textit{Morita equivalent}; we write $(\groupoid,\subgroupoid) \sim (\grouptwoid,\subgrouptwoid)$.
%\end{defn}
In the proofs presented in this section, it is useful to use a different characterisation of Morita equivalences. We need the following definition.\newpage
\begin{defn}\cite{gerbes}\label{defn: Morita morphism}
Let $\groupoid$ and $\grouptwoid$ be Lie groupoids. A morphism of Lie groupoids $\group \ra \grouptwo$ is called a \textit{Morita morphism} if the base map $f$ is a surjective submersion, and the map
\[\varphi \times (\target_{\group},\source_{\group}): \group \ra \grouptwo \tensor[_{(\target_{\grouptwo},\source_{\grouptwo})}]{\times}{_{f \times f}} (\base \times \base)\]
is a diffeomorphism. If $(\groupoid,\subgroupoid)$ and $(\grouptwoid,\subgrouptwoid)$ are pairs of Lie groupoids, then a Lie groupoid morphism $\varphi: (\group,\subgroup) \ra (\grouptwo,\subgrouptwo)$ of pairs is called a \textit{Morita morphism of pairs} if $\varphi: \group \ra \grouptwo$ and $\varphi|_{\subgroup}: \subgroup \ra \subgrouptwo$ are Morita morphisms. Additionally, the base map $f$ of $\varphi$ is required to satisfy $f^{-1}(\subbasetwo)=\subbase$.
\end{defn}
\begin{rema}\label{rema: base map Morita morphism}
%If $\varphi: (\group,\subgroup) \ra (\grouptwo,\subgrouptwo)$ is a Morita morphism of pairs between pairs of groupoids $(\groupoid,\subgroupoid)$ and $(\grouptwoid,\subgrouptwoid)$, $\varphi \times (\target_{\group},\source_{\group})$
%restricts to a diffeomorphism 
%\[f \times (\id_{\subbase},\id_{\subbase}): \subbase \xra{\sim} \subbasetwo \tensor[_{(\id_{\subbasetwo},\id_{\subbasetwo})}]{\times}{_{f \times f}} (\subbase \times \subbase).\] 
%But $\subbasetwo \tensor[_{(\id_{\subbasetwo},\id_{\subbasetwo})}]{\times}{_{f \times f}} (\subbase \times \subbase) \cong f^{-1}(\subbasetwo)$ (the diffeomorphism is simply given by $(n,y,y) \mapsto y$), so we see that we must have $f^{-1}(\subbasetwo)=\subbase$. Since, in addition, 
The reason we put the extra condition on the base map $f$ is to ensure that, if we apply the blow-up functor, then $\blup(f)$ is a map $\blup(\base,\subbase) \ra \blup(\basetwo,\subbasetwo)$. This condition suffices, because $f$ is a surjective submersion, so $f: \base \ra \basetwo$ has a clean intersection with $\subbasetwo$, and therefore $d_\normal f$ is also fiberwise injective (see Remark \ref{rema: clean intersection normal derivative}).
\end{rema}
\begin{prop}\cite{gerbes}\label{prop: Morita morphism is same as equivalence}
Two groupoids $\groupoid$ and $\grouptwoid$ are Morita equivalent if and only if there is a groupoid $\mathcal{P} \rra U$ and Morita morphisms $\mathcal{P} \ra \group$ and $\mathcal{P} \ra \grouptwo$.
\end{prop}
This suggest the following definition for pairs of Lie groupoids:
\begin{defn}\label{defn: pairs of Morita equivalences}
A \textit{Morita equivalence} between pairs of Lie groupoids $(\groupoid,\subgroupoid)$ and $(\grouptwoid,\subgrouptwoid)$ is given by a pair of Lie groupoids $(\mathcal{P} \rra U, \mathcal{Q} \rra V)$ together with Morita morphisms $(\mathcal{P},\mathcal{Q}) \ra (\group,\subgroup)$ and $(\mathcal{P},\mathcal{Q}) \ra (\grouptwo,\subgrouptwo)$. In this case, $(\groupoid,\subgroupoid)$ and $(\grouptwoid,\subgrouptwoid)$ are said to be \textit{Morita equivalent}; we write $(\groupoid,\subgroupoid) \sim (\grouptwoid,\subgrouptwoid)$.
\end{defn}
\begin{rema}\cite{2017arXiv170509588D}\label{rema: alternative morita equivalence}
Using Proposition \ref{prop: Morita morphism is same as equivalence}, we can describe another way to define Morita equivalence between pairs of groupoids. Namely, two pairs of groupoids $(\groupoid,\subgroupoid)$ and $(\grouptwoid,\subgrouptwoid)$ are Morita equivalent if and only if there is a pair of manifolds $(E,F)$ together with a Morita equivalence
\[\base \xla{\pi_\base} E \xra{\pi_\basetwo} \basetwo\] 
such that $F=\pi_\base^{-1}(\subbase)=\pi_\basetwo^{-1}(\subbasetwo)$, and
\[\subbase \xla{\pi_\base|_{F}} F \xra{\pi_\basetwo|_{F}} \subbasetwo\] 
is also a Morita equivalence.
\end{rema}
As in the case of the usual notion of Morita equivalence, Morita equivalence in the setting of pairs of groupoids is an equivalence relation.
By appropriately defining \textit{pairs of differentiable stacks}, we can construct the deformation to the normal cone and blow-up of differentiable stacks.
\begin{defn}\label{defn: pairs of differentiable stacks}
A \textit{pair of differentiable stacks} is a pair $(\mathfrak{X},\mathfrak{Y})$, where $\mathfrak{X}=[\groupoid]$ and $\mathfrak{Y}=[\grouptwoid]$ are differentiable stacks, together with a morphism $\mathfrak{Y} \ra \mathfrak{X}$ given by 
\[(\grouptwoid) \sim (\subgroupoid) \xra{\iota} (\groupoid),\]
where $\iota$ is an embedding of Lie groupoids.
\end{defn}
%Indeed, the condition that $\pi(p)=\pi(p')$ if and only if $p'= p \cdot g_{p'}$ for a unique $g_{p'} \in \group$ means that the map
%\[\pi^{-1}(p) \ra \target_{\group}^{-1}(\mu(p)) \textnormal{ given by } p' \mapsto g_{p'}\]
We will now show that, as claimed, the constructions of deformation to the normal cone, and of blow-up, are ``stable'' under Morita equivalence.
\begin{prop}\label{prop: DNC Morita equivalence}
Let $(\groupoid,\subgroupoid)$ and $(\grouptwoid,\subgrouptwoid))$ be Lie groupoids which are Morita equivalent. Then $\DNC(\groupoid,\subgroupoid)$ and $\DNC(\grouptwoid,\subgrouptwoid)$ are Morita equivalent.
\end{prop}
\begin{proof}
It suffices to show that, if $F: (\group,\subgroup) \ra (\grouptwo,\subgrouptwo)$ is a Morita morphism of pairs, then $\DNC(F)$ is a Morita morphism. And indeed, the base map of $\DNC(F)$ is $\DNC(f): \DNC(\base,\subbase) \ra \DNC(\basetwo,\subbasetwo)$, which is a surjective submersion by Proposition \ref{prop: deformation constant rank}. That the map 
\[\DNC(\varphi) \times (\DNC(\target_{\group}),\DNC(\source_{\group})): \DNC(\group,\subgroup) \ra \DNC(\grouptwo,\subgrouptwo) \tensor[_{(\DNC(\target_{\grouptwo}),\DNC(\source_{\grouptwo}))}]{\times}{_{\DNC(f) \times \DNC(f)}} \DNC(\base,\subbase) \times \DNC(\base,\subbase)\]
is a diffeomorphism goes as follows: notice that $(\DNC(\target_{\grouptwo})(g),\DNC(\source_{\grouptwo})(g))=(\DNC(f)(z),\DNC(f)(w))$ can only happen if $(z,w) \in \DNC(\base,\subbase) \times_{\mathbb{R}} \DNC(\base,\subbase)$ (because $(\DNC(\target_{\grouptwo})(g),\DNC(\source_{\grouptwo})(g))$ is an element of $\DNC(\basetwo,\subbasetwo) \times_{\mathbb{R}} \DNC(\basetwo,\subbasetwo)$). Therefore, by Proposition \ref{prop: if clean intersection, then DNC respects fiber products}, we can equivalently prove that
\[\DNC(\varphi) \times \DNC(\target_{\group},\source_{\group}): \DNC(\group,\subgroup) \ra \DNC(\grouptwo,\subgrouptwo) \tensor[_{\DNC(\target_{\grouptwo},\source_{\grouptwo})}]{\times}{_{\DNC(f \times f)}} \DNC(\base \times \base, \subbase \times \subbase)\]
is a diffeomorphism. The result follows by again using Proposition \ref{prop: if clean intersection, then DNC respects fiber products}, and the observation that the map $\DNC(\varphi \times (\target_{\group},\source_{\group}))$ is a diffeomorphism.
\end{proof}
\begin{prop}\label{prop: blow-up Morita equivalence}
Let $(\groupoid,\subgroupoid)$ and $(\grouptwoid,\subgrouptwoid))$ be Lie groupoids which are Morita equivalent. Then $\blup_{\source_\group,\target_\group}(\group,\subgroup)$ and $\blup_{\source_\grouptwo,\target_\grouptwo}(\grouptwo,\subgrouptwo)$ are Morita equivalent.
\end{prop}
\begin{proof}
Again, it suffices to prove that, if $F: (\group,\subgroup) \ra (\grouptwo,\subgrouptwo)$ is a Morita morphism of pairs, then $\blup(F)$ is a Morita morphism. Notice that, by Proposition \ref{prop: blow-up functoriality groupoids} and Remark \ref{rema: base map Morita morphism}, $\blup(F)$ is a morphism of Lie groupoids. The map $\blup(f): \blup(\base,\subbase) \ra \blup(\basetwo,\subbasetwo)$ is a surjective submersion by Proposition \ref{prop: analogue constant rank when applying blup}. To see in this case that the map $\blup(\varphi) \times (\blup(\target_{\group}),\blup(\source_{\group}))$ is a diffeomorphism, we use
Proposition \ref{prop: blup respects fiber products} and the fact that $\blup(\varphi \times (\target_{\group},\source_{\group}))$ is a diffeomorphism. This proves the statement.
%Let $\base \xla{\pi_\base} E \xra{\pi_\basetwo} \basetwo$ be a Morita equivalence that restricts to a Morita equivalence $\subbase \xla{\pi_\subbase} F \xra{\pi_\subbasetwo} \subbasetwo$. We will show that $\blup(\base,\subbase) \xla{\blup(\pi_\base)} \blup_{\pi_\base,\pi_\basetwo}(E,F) \xra{\blup(\pi_\basetwo)} \blup(\basetwo,\subbasetwo)$ is a Morita equivalence. We already know from Lemma \ref{lemm: DNC Morita equivalence} that $\DNC(\base,\subbase) \xla{\DNC(\pi_\base)} \DNC(E,F) \xra{\DNC(\pi_\basetwo)} \DNC(\basetwo,\subbasetwo)$ is a Morita equivalence. But then 
%\[\DNC(\base,\subbase) \setminus (\subbase \times \mathbb{R}) \xla{\DNC(\pi_\base)} \DNC_{\pi_\base,\pi_\basetwo}(E,F) \xra{\DNC(\pi_\basetwo)} \DNC(\basetwo,\subbasetwo) \setminus (\subbasetwo \times \mathbb{R})\]
%is a Morita equivalence as well. Taking the quotients with respect to the $\mathbb{R}^\times$-actions, the sequence descends to a sequence $\blup(\base,\subbase) \xla{\blup(\pi_\base)} \blup_{\pi_\base,\pi_\basetwo}(E,F) \xra{\blup(\pi_\basetwo)} \blup(\basetwo,\subbasetwo)$, which is now readily verified to be a Morita equivalence. %\textcolor{red}{Nog even preciezer zijn.}
\end{proof}
From this we see that the blow-up of a pair of differentiable stacks is well-defined and yields another differentiable stack. 

A special type of differentiable stack is a so-called orbifold. To describe, in this way, what an orbifold is, we need to know what a proper groupoid is.
\begin{defn}\cite{proper}\label{defn: proper groupoid}
A groupoid $\groupoid$ is called \textit{proper} if the map
\[(\target,\source): \group \ra \base \times \base\]
is a proper map.
\end{defn}
\begin{rema}\cite{proper}\label{rema: proper groupoid}
Observe that a Lie group $G$ (a groupoid over a point) is proper if and only if $G$ is compact, and an action groupoid $G \ltimes M$ is proper if and only if the action is proper. One can show that properness is a Morita invariant property.
\end{rema}
\begin{defn}\cite{Joao}\label{defn: orbifolds}
An \textit{orbifold} $\mathfrak{O}$ is a differentiable stack for which there exists a representative groupoid $\groupoid$ which is \'etale (i.e. $\dim \group = \dim \base$) and proper. We will call a groupoid $\groupoid$ an \textit{orbifold groupoid} if it is \'etale and proper.
\end{defn}
\begin{rema}\label{rema: orbifolds}
There are more concrete descriptions of orbifolds (see e.g. \cite{moerdijk_mrcun_2003}), but we will not go into this. The property of being \'etale is not Morita invariant. One can show that a groupoid is Morita equivalent to an \'etale one if and only if the groupoid is a so-called \textit{folaited groupoid}, which means that all the isotropy groups are discrete (this property is Morita invariant, and from this we also see that being \'etale is not a Morita invariant property; see \cite{moerdijk_mrcun_2003}).
\end{rema}
We now have the following result.
%Notice that if $\groupoid$ is an \'etale Lie groupoid, then, if $\subbase \subset \base$ is saturated, then $\group|_\subbase \rra \subbase$ is again an \'etale Lie groupoid. The reason is that, also in the non-\'etale case, $\group|_\subbase \subset \group$ has the same codimension as $\subbase \subset \base$. We now have the following. 
\begin{prop}\label{prop: orbifold blow-up}
Let $\groupoid$ be a Lie groupoid and let $\subbase \subset \base$ be saturated. If $\groupoid$ is (a) \'etale, or (b) proper, then 
\[\blup(\group,\group|_\subbase) \rra \blup(\base,\subbase)\]
is again (a) \'etale or (b) proper. In particular we can blow up orbifold groupoids, along this special class of sub-orbifold groupoids represented by Lie subgroupoids of the form $\group|_\subbase \rra \subbase$ (with $\subbase \subset \base$ saturated), and obtain another orbifold groupoid. 
\end{prop}
\begin{proof}
(a) Since $\subbase \neq \base$, we have $\dim \blup(\group,\group|_\subbase) = \dim \group$ and $\dim \blup(\base,\subbase) = \dim \base$.

(b) Consider the following commutative diagram:
\begin{center}
\begin{tikzcd}
    \blup(\group,\group|_{\subbase}) \ar{rr}{(\blup(\target),\blup(\source))} \ar{d}{\bldown} & & \blup(\base,\subbase) \times \blup(\base,\subbase) \ar{d}{\bldown \times \bldown} \\
    \group \ar{rr}{(\target,\source)} & & \base \times \base.
\end{tikzcd}
\end{center}
Since $(\target,\source)$ is proper by assumption, and the blow-down map of a manifold is proper (see Proposition \ref{prop: blowdown is proper}), we see, by the above commutative diagram, that the map $(\bldown \times \bldown) \circ (\blup(\target),\blup(\source))$ is proper. The result now follows from Lemma \ref{lemm: proper maps} (if $g \circ f$ is proper, then $f$ is proper).

Since $\group|_\subbase \subset \group$, also in the non-\'etale case, has the same codimension as $\subbase \subset \base$, we see that $\group|_\subbase \rra \subbase$ is again \'etale. Moreover, $(\target,\source): \group|_\subbase \ra \subbase \times \subbase$ is proper, because it equals the map 
\[\group|_\subbase \xhookrightarrow{\iota} \group \xra{(\target,\source)} \base \times \base\]
which is the composition of proper maps (note: $\iota$ is a proper embedding precisely because $\group|_\subbase \subset \group$ is a closed embedded submanifold). This proves the statement.
\end{proof}
\begin{rema}\label{rema: orbifold blow-up}
This yields a blow-up construction (at least in special cases) for orbifolds using the blow-up construction for differentiable stacks described in this section. For more details on suborbifolds, and how to describe them using the theory of Lie groupoids, see \cite{Joao}.
\end{rema}

\addtocontents{toc}{\protect\thispagestyle{myplain}}\newpage

\printbibliography

\end{document}